%% file: Gutzwiller_trace__onirban_islam.tex
\documentclass[a4paper,12pt]{book}
\usepackage[T1]{fontenc}
\usepackage[utf8]{inputenc}
\usepackage{tikz}
\usetikzlibrary{decorations.text,calc,arrows.meta}
\usepackage{tkz-euclide}
\usepackage{caption}
\usetikzlibrary{perspective}
\usepackage{slashed}
\usepackage{amsfonts, amssymb, amsmath, amsthm}
\usepackage{thm-restate}
\usepackage{nicefrac}
\usepackage{enumerate} 
\usepackage{enumitem}
\usepackage{xcolor}
\usepackage{a4wide}
\usepackage{bbding}
\usepackage{charter}
\usepackage{newtxmath} 
\usepackage{mathrsfs} 
\usepackage{calligra}
\usepackage{bbm}
\usepackage{hyperref}
\hypersetup{
    colorlinks,%
    citecolor=green!50!black,%
    filecolor=black,%
    linkcolor=blue,%
    urlcolor=cyan!50!blue
}
\usepackage[numbers,sort&compress,merge]{natbib}
\usepackage{graphicx}
\usepackage{booktabs}
\usepackage{subfig}
\usepackage{xparse} 
\usepackage[toc,page]{appendix}
\usepackage[final]{pdfpages}
\theoremstyle{plain}
\newtheorem{theorem}{Theorem}[section]
\newtheorem{proposition}[theorem]{Proposition}
\newtheorem{lemma}[theorem]{Lemma}
\newtheorem{corollary}[theorem]{Corollary}
\theoremstyle{definition}
\newtheorem{definition}[theorem]{Definition}
\newtheorem{example}[theorem]{Example}
\newtheorem{remark}[theorem]{Remark}

\newtheorem{assumption}[theorem]{Assumption}
\DeclareMathAlphabet{\mathsf}{OT1}{\sfdefault}{m}{n} 
\DeclareMathOperator{\Ad}{Ad} 

\DeclareMathOperator{\sgn}{sgn} 
\DeclareMathOperator{\codim}{codim}
\newcommand{\scalarProdOne}[2]{\left\langle {#1}, {#2} \right\rangle} 
\newcommand{\scalarProdTwo}[2]{\left\langle {#1} \left| {#2} \right. \right\rangle}

\newcommand{\jbracket}[1]{\left\langle {#1} \right\rangle}
\newcommand{\norm}[2]{  \left\| {#1} \right\|_{ \scriptscriptstyle{{#2}} }  }
\newcommand{\seminorm}[2]{\left| \! \left| \! \left| {#1} \right| \! \right| \! \right|_{{#2}}}

\newcommand{\sA}{\mathscr{A}}
\newcommand{\cA}{\mathcal{A}}

\newcommand{\fA}{\mathsf{A}}
\newcommand{\tA}{\mathtt{A}}
\newcommand{\kA}{\mathfrak{A}}
\newcommand{\fa}{\mathsf{a}}

\newcommand{\cB}{\mathcal{B}}

\newcommand{\B}{\mathbb{B}}
\newcommand{\fB}{\mathsf{B}}
\newcommand{\fb}{\mathsf{b}}

\newcommand{\tB}{\mathtt{B}}

\newcommand{\sC}{\mathscr{C}}
\newcommand{\cC}{\mathcal{C}}

\newcommand{\fc}{\mathsf{c}} 

\newcommand{\kC}{\mathfrak{C}}

\newcommand{\C}{\mathbb{C}}

\newcommand{\rc}{\mathrm{c}}

\newcommand{\sD}{\mathscr{D}}
\newcommand{\cD}{\mathcal{D}}

\newcommand{\rD}{\mathrm{D}}
\newcommand{\rd}{\mathrm{d}}
\newcommand{\fd}{\mathsf{d}}

\newcommand{\dM}{ d_{\ms M}}
\newcommand{\dN}{ d_{\ms N}}

\newcommand{\E}{\mathbb{E}}
\newcommand{\re}{\mathrm{e}}
\newcommand{\rE}{\mathrm{E}}
\newcommand{\cE}{\mathcal{E}}
\newcommand{\sE}{\mathscr{E}}
\newcommand{\fE}{\mathsf{E}}
\newcommand{\tE}{\mathtt{E}}

\newcommand{\cF}{\mathcal{F}}
\newcommand{\sF}{ \mathscr{F} }
\newcommand{\Fourier}{\mathcal{F}}
\newcommand{\ff}{\mathsf{f}}

\newcommand{\rF}{\mathrm{F}}

\newcommand{\fF}{\mathsf{F}}
\newcommand{\kF}{\mathfrak{F}}

\newcommand{\sG}{\mathscr{G}}

\newcommand{\fG}{\mathsf{G}}
\newcommand{\fg}{\mathsf{g}}

\newcommand{\sH}{\mathscr{H}}
\newcommand{\cH}{\mathcal{H}}

\newcommand{\fh}{\mathsf{h}}
\newcommand{\rH}{\mathrm{H}}

\newcommand{\cI}{\mathcal{I}}

\DeclareMathOperator{\ri}{i}
\newcommand{\fI}{\mathsf{I}}

\newcommand{\fJ}{\mathsf{J}}

\DeclareMathOperator{\rJ}{J}

\newcommand{\cK}{\mathcal{K}}

\newcommand{\fK}{\mathsf{K}}

\newcommand{\sL}{\mathscr{L}}
\newcommand{\cL}{\mathcal{L}}
\newcommand{\bbL}{\mathbb{L}} 
\newcommand{\bbl}{\mathbbm{l}}
\newcommand{\rL}{\mathrm{L}}
\DeclareMathOperator{\LieDeri}{L}
\newcommand{\fL}{\mathsf{L}}
\newcommand{\fl}{\mathsf{l}}

\newcommand{\kl}{\mathfrak{l}} 
\newcommand{\sM}{\mathscr{M}}

\newcommand{\cM}{\mathcal{M}}
\newcommand{\fM}{\mathsf{M}}
\newcommand{\fm}{\mathsf{m}}
\newcommand{\km}{\mathfrak{m}}
\newcommand{\kM}{\mathfrak{M}}
\newcommand{\rM}{\mathrm{M}}
\newcommand{\rrm}{\mathrm{m}}
\newcommand{\m}{\mathbbm{m}}
 
\newcommand{\M}{\mathbb{M}}

\newcommand{\rN}{\mathrm{N}}

\newcommand{\kn}{\mathfrak{n}}

\newcommand{\N}{\mathbb{N}}
\newcommand{\NO}{\mathbb{N}_{0}}
\newcommand{\sN}{\mathscr{N}}
\newcommand{\fN}{\mathsf{N}}

\newcommand{\nM}{n_{\ms M}}
\newcommand{\nN}{n_{\ms N}}

\newcommand{\fO}{\mathsf{O}}

\newcommand{\rP}{\mathrm{P}}

\newcommand{\bbP}{\mathbb{P}}

\newcommand{\sP}{\mathscr{P}}
\newcommand{\cP}{\mathcal{P}}

\newcommand{\kp}{\mathfrak{p}}

\newcommand{\fP}{\mathsf{P}}
\newcommand{\fp}{\mathsf{p}}

\newcommand{\sQ}{\mathscr{Q}}

\newcommand{\rQ}{ \mathrm{Q} }
\newcommand{\fQ}{\mathsf{Q}}

\newcommand{\fq}{\mathsf{q}}
\newcommand{\cR}{\mathcal{R}}
\newcommand{\sR}{\mathscr{R}}
\newcommand{\R}{\mathbb{R}}
\newcommand{\rR}{ \mathrm{R} }  
\newcommand{\fR}{\mathsf{R}} 
 
\newcommand{\kR}{\mathfrak{R}} 
\newcommand{\rr}{\mathrm{r}}

\newcommand{\Rd}{\mathbb{R}^{d}}
\newcommand{\Rn}{\mathbb{R}^{n}}

\newcommand{\dotRd}{\dot{\mathbb{R}}^{d}}
\newcommand{\dotRn}{\dot{\mathbb{R}}^{n}}
\newcommand{\sS}{\mathscr{S}}

\newcommand{\fS}{\mathsf{S}}

\newcommand{\cS}{\mathcal{S}}

\newcommand{\bbS}{\mathbb{S}}
\newcommand{\fs}{\mathsf{s}}

\newcommand{\ks}{\mathfrak{s}}

\newcommand{\cT}{\mathcal{T}}

\newcommand{\ft}{\mathsf{t}}
\newcommand{\bbt}{\mathbbm{t}}
\newcommand{\fT}{\mathsf{T}}

\newcommand{\tT}{\texttt{T}}
\newcommand{\T}{\mathbb{T}}

\newcommand{\dotCoTan}{\dot{\mathrm{T}}^{*}}
\newcommand{\dotCoTanRd}{\dot{\mathrm{T}}^{*} \! \mathbb{R}^{d}}
\newcommand{\dotCoTanRn}{\dot{\mathrm{T}}^{*} \! \mathbb{R}^{n}}

\newcommand{\dotCoTanM}{\dot{\mathrm{T}}^{*} \! M}
\newcommand{\dotCoTanN}{\dot{\mathrm{T}}^{*} \! N}
\newcommand{\dotCoTanMM}{\dot{\mathrm{T}}^{*} \! (M \times M)}
\newcommand{\dotCoTanMN}{\dot{\mathrm{T}}^{*} \! (M \times N)}
\newcommand{\dotCoTansM}{\dot{\mathrm{T}}^{*} \! \sM}
\newcommand{\sU}{\mathscr{U}}
\newcommand{\cU}{\mathcal{U}}
\newcommand{\kU}{\mathfrak{U}}
 
\newcommand{\fU}{\mathsf{U}}

\newcommand{\U}{\mathbb{U}}

\newcommand{\sV}{\mathscr{V}}
\newcommand{\cV}{\mathcal{V}}

\newcommand{\kV}{\mathfrak{V}}

\newcommand{\fv}{\mathsf{v}}

\newcommand{\cW}{\mathcal{W}}

\newcommand{\fW}{\mathsf{W}}
\newcommand{\fw}{\mathsf{w}}
\newcommand{\sX}{\mathscr{X}}

\newcommand{\x}{\mathbbm{x}}
\newcommand{\bx}{\boldsymbol{x}}

\newcommand{\cY}{\mathcal{Y}}

\newcommand{\by}{\boldsymbol{y}}
\newcommand{\Z}{\mathbb{Z}}

\newcommand{\bz}{\boldsymbol{z}}

\newcommand{\bxi}{\boldsymbol{\xi}}
\newcommand{\btheta}{\boldsymbol{\theta}} 
\newcommand{\bfeta}{\boldsymbol{\eta}}
\newcommand{\bzeta}{\boldsymbol{\zeta}}

\newcommand{\CN}{\mathbb{C}^{k}} 

\DeclareMathOperator{\grad}{grad}
\DeclareMathOperator{\Div}{div}

\DeclareMathOperator{\Hol}{Hol}

\newcommand{\spacetime}{(\mathscr{M}, \mathsf{g})}
\newcommand{\Cauchy}{(\varSigma, \mathsf{h})}
\newcommand{\SST}{(\mathscr{M}, \mathsf{g}, \varXi)}

\newcommand{\tangent}{\mathrm{T}}

\newcommand{\tansM}{\mathrm{T} \! \mathscr{M}}

\newcommand{\coTan}{\mathrm{T}^{*}}

\newcommand{\coTanM}{\mathrm{T}^{*} \! M}
\newcommand{\coTanN}{\mathrm{T}^{*} \! N}

\newcommand{\coTansM}{\mathrm{T}^{*} \! \mathscr{M}}
\newcommand{\coLightBun}{\dot{\mathrm{T}}_{0}^{*} \sM}

\newcommand{\coTanCauchy}{\mathrm{T}^{*} \! \varSigma}
\newcommand{\dotCoTanCauchy}{\dot{\mathrm{T}}^{*} \! \varSigma}

\newcommand{\halfDen}{\varOmega^{\nicefrac{1}{2}}}
\newcommand{\halfDenM}{\varOmega^{\nicefrac{1}{2}} \! M}
\newcommand{\halfDenN}{\varOmega^{\nicefrac{1}{2}} \! N}
\newcommand{\halfDenMM}{\varOmega^{\nicefrac{1}{2}} (M \times M)}
\newcommand{\halfDenMN}{\varOmega^{\nicefrac{1}{2}} (M \times N)}
\newcommand{\halfDenC}{\varOmega^{\nicefrac{1}{2}} C}

\newcommand{\Maslov}{\mathbb{M}}

\newcommand{\CComk}[1]{C_{\mathrm{c}}^{k} ({#1})}
\newcommand{\CInfinity}[1]{C^{\infty} ({#1})}
\newcommand{\CComInfinity}[1]{C_{\mathrm{c}}^{\infty} ({#1})}

\newcommand{\triv}{\textrm{triv}}
\newcommand{\nonTriv}{\textrm{ntriv}}

\newcommand{\PsiDO}[2]{\Psi \mathrm{DO}^{{#1}} ({#2})}
\newcommand{\FIO}[2]{\mathrm{FIO}^{{#1}} ({#2})}
\newcommand{\PsiDONU}[1]{\Psi \mathrm{DO}^{{#1}} (U',\matk)}
\newcommand{\PsiDON}[1]{\Psi \mathrm{DO}^{{#1}} (\mathbb{R}^{d}, \matk)}

\newcommand{\matN}{\C^{N \times N}}
\newcommand{\matk}{\C^{k \times k}}

\DeclareMathOperator{\PDO}{PDO}
\newcommand{\PDOU}{\mathrm{PDO}^{m} (U)}
\newcommand{\PsiDOU}{\Psi \mathrm{DO}^{m} (U)}

\newcommand{\PsiDOM}{\Psi \mathrm{DO}^{m} (M; \halfDen)} 
\newcommand{\PsiDOME}{\Psi \mathrm{DO}^{m} (M; \sE \otimes \halfDen)}
\newcommand{\PsiDOMEF}{\Psi \mathrm{DO}^{m} (M; \sE \otimes \halfDen, \sF \otimes \halfDen)}

\newcommand{\FIONM}{\mathrm{FIO}^{m} (\halfDenN \to N, \halfDenM \to M; C')}

\newcommand{\LagrangianDistM}{I^{m} (M, \varLambda; \halfDen)}
\newcommand{\LagrangianDistMN}{I^{m} \big( M \times N, C'; \halfDenMN \big)}

\newcommand{\LagrangianDistHomFE}{I^{m} \big( M \times N, C'; \Hom{\sF, \sE} \otimes \halfDenMN \big)}
\newcommand{\LagrangianDistGraphHomFE}{I^{m} \big( M \times N, \varGamma'; \Hom{\sF, \sE} \otimes \halfDenMN \big)}

\newcommand{\FIOFE}{\mathrm{FIO}^{m} (\sF \otimes \halfDenN \to N, \sE \otimes \halfDenM \to M; C')}

\newcommand{\secE}{C^{\infty} (M; \sE \otimes \halfDen)} 
\newcommand{\comSecE}{C_{\mathrm{c}}^{\infty} (M; \sE \otimes \halfDen)}
\newcommand{\secEStar}{C^{\infty} (M; \sE^{*} \otimes \halfDen)} 
\newcommand{\comSecEStar}{C_{\mathrm{c}}^{\infty} (M; \sE^{*} \otimes \halfDen)}
\newcommand{\dualSecE}{\mathcal{E}^{\prime} (M; \sE \otimes \halfDen)}
\newcommand{\dualComSecE}{\mathcal{D}^{\prime} (M; \sE \otimes \halfDen)}
\newcommand{\secME}{C^{\infty} (M; \sE \otimes \halfDenM)} 
\newcommand{\comSecME}{C_{\mathrm{c}}^{\infty} (M; \sE \otimes \halfDenM)}

\newcommand{\dualComSecME}{\mathcal{D}^{\prime} (M; \sE \otimes \halfDenM)} 
 
\newcommand{\comSecNF}{C_{\mathrm{c}}^{\infty} (N; \sF \otimes \halfDenN)}

\newcommand{\secMF}{C^{\infty} (M; \sF \otimes \halfDen)}

\newcommand{\secsME}{C^{\infty} (\mathscr{M; E})}
\newcommand{\comSecsME}{C_{\mathrm{c}}^{\infty} (\mathscr{M; E})}
\newcommand{\dualComSecsME}{\mathcal{D}^{\prime} (\mathscr{M; E})}
\newcommand{\dualSecsME}{\mathcal{E}^{\prime} (\mathscr{M; E})}
\newcommand{\secsMEStar}{C^{\infty} (\mathscr{M}; \mathscr{E}^{*})}
\newcommand{\comSecsMEStar}{C_{\mathrm{c}}^{\infty} (\mathscr{M}; \mathscr{E}^{*})}

\newcommand{\comSecsMEConjStar}{C_{\mathrm{c}}^{\infty} (\mathscr{M}; \bar{\mathscr{E}}^{*})}
\newcommand{\secECauchy}{C^{\infty} (\varSigma; \sE_{\varSigma})}

\newcommand{\comSecECauchy}{C_{\mathrm{c}}^{\infty} (\varSigma; \sE_{\varSigma})}

\newcommand{\connectionLC}{ \nabla^{\ms \mathrm{LC}} }
\newcommand{\connectionE}{ \nabla^{\ms \sE} } 
\newcommand{\connectionCoTanME}{ \nabla^{\ms \mathrm{LC} \otimes \sE} }

\newcommand{\connectionEndE}{ \nabla^{\ms \Hom{\sE, \sE}} }
\newcommand{\connectionEndPiE}{ \nabla^{\ms \pi^{*} \Hom{\sE, \sE}} }

\newcommand{\PFBSO}{( \mathscr{Q}, \mathrm{SO}_{0},  \mathscr{M})}
\newcommand{\PFBSpin}{( \mathscr{P}, \mathrm{Spin}_{0},  \mathscr{M})}

\newcommand{\symb}[1]{ \sigma_{\ms{#1}} }

\newcommand{\totSymb}[1]{\sigma^{\mathrm{tot}}_{\ms{#1}}}
\newcommand{\subSymb}[1]{ \sigma^{\mathrm{sub}}_{\ms{#1}} }
 
\newcommand{\sbracket}[1]{ \left\{ {#1} \right\} }

\newcommand{\equiBra}[1]{[{#1}]} 

\newcommand{\Dxbeta}{\mathrm{D}_{x}^{\beta}}

\newcommand{\parDeri}[2]{\frac{\partial {#2}}{\partial {#1}}}
\newcommand{\GreenOpAdv}{G^{\mathrm{adv}}}
\newcommand{\GreenOpRet}{G^{\mathrm{ret}}}
\newcommand{\GreenOpFeyn}{G^{\mathrm{Fyn}}}

\newcommand{\GreenOpAdvKernel}{\mathsf{G}^{\mathrm{adv}}}

\newcommand{\GreenKernelAdv}{\mathsf{G}^{\mathrm{adv}}}
\newcommand{\GreenKernelRet}{\mathsf{G}^{\mathrm{ret}}}
\newcommand{\GreenKernelFeyn}{\mathsf{G}^{\mathrm{Fyn}}}

\newcommand{\fundaSolAdv}{F^{\mathrm{adv}}}
\newcommand{\fundaSolRet}{F^{\mathrm{ret}}}
\newcommand{\fundaSolFeyn}{F^{\mathrm{Fyn}}}

\newcommand{\fundaSolAdvKernel}{\mathsf{F}^{\mathrm{adv}}}
\newcommand{\fundaSolRetKernel}{\mathsf{F}^{\mathrm{ret}}}

\newcommand{\parametrixLeft}{ E_{\mathrm{L}} } 
\newcommand{\parametrixRight}{ E_{\mathrm{R}} }
\newcommand{\parametrixAdv}{ E^{\mathrm{adv}} }
\newcommand{\parametrixRet}{ E^{\mathrm{ret}} }
\newcommand{\parametrixFeyn}{ E^{\mathrm{Fyn}} }
\newcommand{\parametrixAntiFeyn}{ E^{\mathrm{aFyn}} }

\newcommand{\parametrixKernelLeft}{ \mathsf{E}_{\mathrm{L}} } 
\newcommand{\parametrixKernelRight}{ \mathsf{E}_{\mathrm{R}} }
\newcommand{\parametrixKernelAdv}{\mathsf{E}^{\mathrm{adv}}}
\newcommand{\parametrixKernelRet}{\mathsf{E}^{\mathrm{ret}}}
\newcommand{\parametrixKernelFeyn}{\mathsf{E}^{\mathrm{Fyn}}} 
\newcommand{\parametrixKernelAntiFeyn}{\mathsf{E}^{\mathrm{aFyn}}}

\DeclareMathOperator{\clo}{clo} 
\DeclareMathOperator{\adv}{adv} 
\DeclareMathOperator{\ret}{ret} 
\DeclareMathOperator{\Feyn}{Fyn}
\DeclareMathOperator{\aFeyn}{aFyn}
\DeclareMathOperator{\Char}{Char}
\DeclareMathOperator{\ES}{ES} 
\DeclareMathOperator{\WF}{WF} 

\DeclareMathOperator{\WFPrime}{WF^{\prime}} 

\DeclareMathOperator{\spec}{Spec} 
\DeclareMathOperator{\supp}{supp} 
\DeclareMathOperator{\singsupp}{sing supp}

\DeclareMathOperator{\esssupp}{ess supp}
\DeclareMathOperator{\Tr}{Tr}
\DeclareMathOperator{\tr}{tr}
\DeclareMathOperator{\rk}{rk}

\DeclareMathOperator{\res}{res}

\newcommand{\LC}{\ms \mathrm{LC}}
\newcommand{\img}[1]{\mathrm{img} \left({#1} \right)}
\newcommand{\dVol}{\mathrm{d} \mathsf{v}}
\newcommand{\dVolg}{\mathrm{d} \mathsf{v}_{\mathsf{g}}}
\newcommand{\dVolh}{\mathrm{d} \mathsf{v}_{\mathsf{h}}}
\DeclareMathOperator{\vol}{vol} 

\newcommand{\Hom}[1]{\mathrm{Hom} ({#1})}
\DeclareMathOperator{\End}{End}
\DeclareMathOperator{\Aut}{Aut}
\DeclareMathOperator{\Iso}{Iso}

\DeclareMathOperator{\id}{id} 
\DeclareMathOperator{\tId}{\texttt{I}} 
\DeclareMathOperator{\I}{\mathbb{I}}    
\DeclareMathOperator{\one}{\mathbbm{1}} 
\DeclareMathOperator{\ric}{ric}
\DeclareMathOperator{\Ric}{Ric}

\DeclareMathOperator{\Hess}{Hess}
\DeclareMathOperator{\pr}{pr}
\DeclareMathOperator{\cst}{cst}

\newcommand{\measurexi}{\frac{\mathrm{d} \xi}{(2 \pi)^{d}}} 
\newcommand{\measuretheta}{\frac{\mathrm{d} \theta}{(2 \pi)^{d}}} 
\newcommand{\xtheta}{(x; \theta)}
\newcommand{\xxi}{(x, \xi)}
\newcommand{\yeta}{(y, \eta)}
 
\newcommand{\xytheta}{(x, y; \theta)}

\newcommand{\xxiyeta}{(x, \xi; y, \eta)}

\newcommand{\xthetaNot}{(x_{0}; \theta^{0})}
\newcommand{\xxiNot}{(x_{0}, \xi^{0})}

\newcommand{\xythetaNot}{(x_{0}, y_{0}; \theta^{0})}
\newcommand{\xxiyetaNot}{(x_{0}, \xi^{0}; y_{0}, \eta^{0})}

\newcommand{\rSpin}{\mathrm{Spin} \,}

\newcommand{\SpinO}{\mathrm{Spin}_{0} (1, d - 1)}

\newcommand{\SOO}{ \mathrm{SO}_{0} (1, d-1) }
\DeclareMathOperator{\rO}{O}
\DeclareMathOperator{\rSO}{SO}
\DeclareMathOperator{\rU}{U}
\newcommand{\CCR}{\mathscr{C \! C \! R}}
\newcommand{\CAR}{\mathscr{C \! A \! R}}
\newcommand{\ms}{\scriptscriptstyle}

\newcommand{\helix}[3] 
{%
  \draw[green!50!black,thick,#3] plot[domain=#1:#2,samples={(#2-#1)/3+1}]
    ({-\r*cos(\x)},{-\r*sin(\x)},\x*\h/540);
}
\begin{document}
    \frontmatter 
    \include{title}
    \include{dedication}

    \include{declaration} 
    \include{abstract}
    \include{acknowledgment}

    \tableofcontents
    \listoffigures
    \listoftables
     \mainmatter
     \include{introduction}
     \include{differential_operator}  
     \include{Feynman_propagator}

     \include{Gutzwiller_trace}

     \begin{appendix}
        \include{composition_density_canonical_relation}
        \include{symbol}

        \include{QFT}

     \end{appendix}
%

\end{document}

%% file: title.tex
\begin{titlepage}
\begin{center}
    \huge \textbf{A Gutzwiller Trace formula for Dirac Operators on a Stationary Spacetime}
\end{center}

\vspace*{3.0cm}
\begin{center}
	\includegraphics[width=7cm, height=2cm]{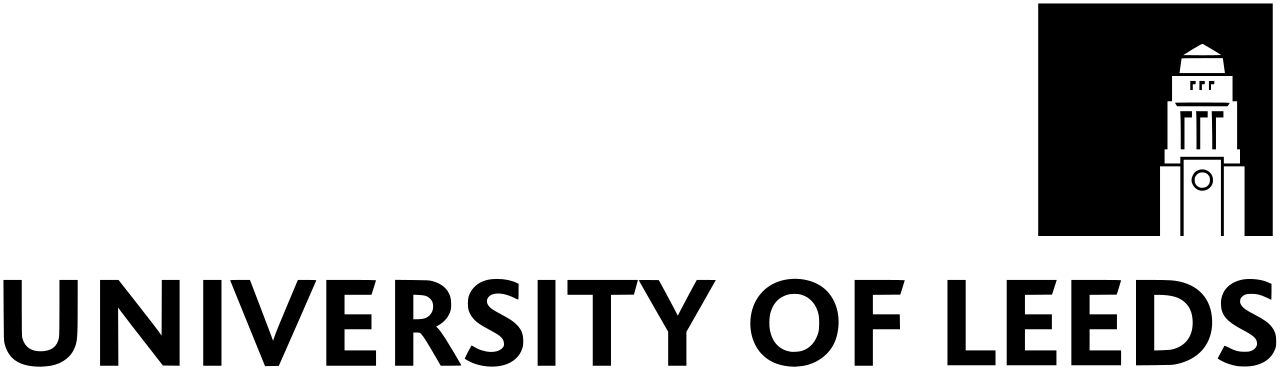}
\end{center}

\vspace*{2.5cm}
\begin{center}
    \textbf{Mohammad Onirban Islam} 
    \\ 
    Department of Pure Mathematics 
    \\ 
    School of Mathematics 
    \\ 
    University of Leeds 
    \\ 
    LS2 9JT, United Kingdom
    \\ 
    \vspace*{2cm}
    
    A thesis submitted for the degree of 
    \\ 
    Doctor of Philosophy 
    \\ 
    October 2022
\end{center}
\end{titlepage}

%% file: dedication.tex
\begin{center}
    \vspace*{\stretch{1}}
        \textrm{\calligra{... To all who have taught me}}
    \vspace*{\stretch{1}}
\end{center}

%% file: declaration.tex
\chapter*{Declaration}
To the best of my knowledge, I declare that the research work in this thesis is original and my own work, except where either cited or commonly known. 
The work presented in Chapter $3$ is identical to the joint paper with 
Alexander Strohmaier ``On microlocalization and the construction of Feynman propagators for normally hyperbolic operators'', arXiv:2012.09767 [math.AP] (accepted in \textit{Communications in Analysis and Geometry}) with a few added background materials. 
Chapter $4$ comprises a slightly expository version of the single-authored article 
``A Gutzwiller trace formula for Dirac operators on a stationary spacetime'', 
arXiv:2109.09219 [math.AP], published in \href{https://doi.org/10.1007/s12220-022-01084-x}{J Geom Anal $\mathbf{33}, 57 (2023)$}. 

%% file: abstract.tex
\chapter*{Abstract}
A Duistermaat-Guillemin-Gutzwiller trace formula for Dirac-type operators on a globally hyperbolic spatially compact stationary spacetime is achieved by generalising the recent construction by 
Strohmaier and Zelditch [Adv. Math. \textbf{376}, 107434 (2021)] 
to a vector bundle setting. 
We have analysed the spectrum of the Lie derivative with respect to a global timelike Killing vector field on the solution space of the Dirac equation and found that it consists of discrete real eigenvalues. 
The distributional trace of the time evolution operator has singularities at the periods of induced Killing flow on the manifold of lightlike geodesics. 
This gives rise to the Weyl law asymptotic at the vanishing period.    
A pivotal technical ingredient to prove these results is the Feynman propagator. 
In order to obtain a Fourier integral description of this propagator, we have generalised the classic work of 
Duistermaat and H\"{o}rmander [Acta Math. \textbf{128}, 183 (1972)] 
on distinguished parametrices for normally hyperbolic operators on a globally hyperbolic spacetime by propounding their microlocalisation theorem to a bundle setting. 
As a by-product of these analyses, another proof of the existence of Hadamard bisolutions for a normally hyperbolic operator (resp. Dirac-type operator) is reported. 

%% file: acknowledgment.tex
\chapter*{Acknowledgments}
I am indebted to my supervisor Alexander Strohmaier for providing me the opportunity to explore the ineffable interplay among Fourier integral operators, quantum field theories, and spectral geometry. 
Being a bachelor's in Electronics and Communication Engineering and a master's in Physics, it was a daunting task to step into the realm of pure mathematics unless I had received such an excellent supervision.  
Throughout the Ph.D. program, I have been supported by the Leeds International Doctoral Studentship, and I am grateful to everyone who made it possible for me. 
\newline 

I acknowledge Benjamin Sharp for his support in clarifying the Hessian of diagonal embedding and commenting on the presentation of the thesis, as my second supervisor.
\newline 

My gratitude to Yan-Long Fang, and Rakibur Rahman and Christian B\"{a}r for discussions on Fourier integral operators and Rarita-Schwinger operator, respectively. 
I am thankful to Ryan O'Loughlin and Ahasan Habib for their constant support to cope with the academic and non-academic life in Leeds. 
Special thanks to Joseph Elmes, Amrita Khan, Sayeed Shawon, Sadat Hasan, Mehboob Rehman Khan, Jowadul Kader and Sadia Rahman for providing mental support during my breaking points and helping me with the proof-reading.  
\textit{Grazie tante} Marta Fiori Carones for making the British weather colourful!
\newline 

I would like to use this opportunity to express my acknowledgment to all of my previous supervisors: 
Stefan Hollands, Ko Sanders, J\"{u}rgen Berges, Arshad Momen, ASM Moslehuddin,  and Zahid Hasan Mahmood. 
I have learned something more valuable than just academic research from them.  
Unfortunately, Professor Mahmood left us at the beginning of my Ph.D. but it will be the least to say that without him and his wife, Kaniz Fatema, I might be ended up being an engineer. 
\newline 

I am grateful to my parents, Ruksana Islam and ABM Rafiqul Islam. 
They let me pursue my ambition regardless of their initial objections, and tried their level best to provide educational facilities from my childhood to high school. 
Each family member has assisted my educational journey directly or indirectly and I pay my tribute to them. 
In particular, apology and love are the least I could have for my wife Rebeka Alam for her sacrifice.

%% file: introduction.tex
\chapter{Introduction}
\label{ch: introduction}
%
%
%
In 1966, Kac~\cite{Kac_AMM_1966} asked in a rather leisurely fashion: 
\begin{center}
    ``\textit{Can you hear the shape of a drum}?'' 
\end{center}
This is arguably the quintessential question in spectral geometry. 
The Gutzwiller trace formula is one of the paramount tools used to investigate such a beautiful interplay between geometry and the spectral world in the sense explained below. 
%
%
%
%
%
%
%
%
%
%
\section{Trace formulae}
\subsection{A prelude} 
In order to understand the notion of a trace formula precisely, let us consider the simplest non-trivial differential operator: the derivative 
\begin{equation}
    \rD := - \ri \frac{\rd}{\rd x}
\end{equation}
acting on smooth functions $C^{\infty} (\bbS)$ on a circle $\bbS$ of length $T$. 
Since $C^{\infty} (\bbS)$ is in a one-to-one correspondence with $T$-periodic functions on $\R$, one defines 
\begin{equation} \label{eq: def_L_two_circle}
	L^{2} (\bbS) := \Big\{ u \in C^{\infty} (\R) \,|\, \forall x \in \R : u (x + T) = u (x), \int_{0}^{T} \bar{u} (x)  \, u (x) \, \rd x < \infty \Big\}.  
\end{equation}
It is straightforward to compute the spectrum of $\rD$: 
\begin{equation} \label{eq: spec_spin_Dirac_trivial_circle}
    \spec \rD = \{ \lambda_{n} | n \in \Z \}, 
    \quad 
    \lambda_{n} := n \omega, 
    \quad 
    \omega := \frac{2 \pi}{T}
\end{equation}
corresponding to the orthonormal eigenfunctions $\phiup_{n} := \re^{\ri \lambda_{n} x} / \sqrt{T}$. 
For any Schwartz function $\rho \in \cS (\R)$ on $\R$ and for all $T > 0$, the \textbf{Poisson summation formula} reads  
(see e.g.~\cite[Thm. 1.2]{Uribe_Cuernavaca_1998})
\begin{equation} \label{eq: Poisson_summation_formula_derivative_circle}
    \sum_{n \in \Z} \rho (\lambda_{n}) = T \sum_{k \in \Z} (\Fourier^{-1} \rho) (kT), 
\end{equation}
where $\Fourier^{-1} (\rho)$ denotes the inverse Fourier transformation of $\rho$. 
One observes that the left-hand side of this identity entirely consists of integral multiple of the fundamental angular frequency $\omega$ on $\bbS$, whereas the right-hand side involves only the length $T$ of $\bbS$. 
In other words, \textit{the Poisson summation formula associates the spectral data with a geometric quantity}. 
As a natural suspicion, we may ask whether such a connection is something specific to the simple operator ($\rD$) and/or the geometry ($\bbS$) we have considered. 
However, this suspicion does not hold. 
For instance, the Poisson summation formula holds for a Laplacian on a torus (see e.g.~\cite[Thm. 1.3]{Uribe_Cuernavaca_1998}). 
Thus, one may wonder whether there are some deeper reasons behind such a bridge between the spectral data and the geometric world. 
\newline 

In order to gain a systematic understanding, let us consider the half-density bundle (see Section~\ref{sec: distribution_density_mf}) $\halfDen \varSigma \to \varSigma$ and a (positive-definite) hermitian vector bundle $(\sE_{\varSigma} \to \varSigma, \scalarProdTwo{\cdot}{\cdot})$, over a $d-1$-dimensional \textit{compact} manifold $\varSigma$ without boundary and let $P$ be an \textit{elliptic} (see Definition~\ref{def: char_PsiDO}), \textit{symmetric, first-order} pseudodifferential operator (see Definition~\ref{def: PsiDO}) on $\sE_{\varSigma} \otimes \halfDen \varSigma$.   
It is well-known that the spectrum of $P$ consists of discrete real eigenvalues $\lambda_{n}^{\pm}$ and is not necessarily semi-bounded; each $\lambda_{n}^{\pm}$ is isolated and has a finite multiplicity. 
Here, $\lambda_{n}^{\pm}$ are enumerated in increasing-order in the sense that 
$\ldots < \lambda_{2}^{-} < \lambda_{1}^{-} < 0 =: \lambda_{0} < \lambda_{1}^{+} < \lambda_{2}^{+} < \ldots$ 
and with multiplicities.  
In general, $P$ is not a trace-class operator but its suitable functions can be. 
For instance, the most important operators in geometric analysis, such as,  the Laplace-Beltrami operator (see~\eqref{eq: Laplacian_coordinate}), the d'Alembertian (see Example~\ref{exm: covariant_Klein_Gordon_op}), the Schr\"{o}dinger operator, and the Dirac operator (see Definition~\ref{def: Dirac_type_op}) are not trace-class yet their suitable functions are trace-class. 
We opt for an operator of the following form, originally due to 
Levitan~\cite{Levitan_1952, *Levitan_1955} 
and to 
Avakumovi\v{c}~\cite{Avakumovic_MathZ_1956}   
(see also~\cite{Hoermander_ActaMath_1968})   
\begin{equation} \label{eq: def_U_t}
	\forall t \in \R : U_{t} := \re^{- \ri t P}.  
\end{equation}
This is by no means a canonical choice; see Table~\ref{tab: choice_P_U_t} for some other possibilities. 
However, adopting~\eqref{eq: def_U_t} over the other options offers the following two advantages~\cite[p. 196]{Hoermander_ActaMath_1968} 
(see also~\cite[p. XXV]{Ivrii_I_Springer_2019}):  

%
%
%
\begin{table}[!h]
	\centering 
	\begin{tabular}{ccccc}
		\toprule
        $U_{t}$ & $t$ & $P$ & Partial Differential Equation (PDE) & Method  
		\\ \hline 
        $\frac{1}{t - P}$ & $\R \setminus \spec P$ & - & $(t - P) U_{t} = I$ & Resolvent
		\\ 
        $P^{-t}$ & $\C$ & - & - & $\zeta$-function 
		\\ 
        $\re^{- t P}$ & $\R_{+}$ &$\Delta$ & Heat eq. $(\partial_{t} + \Delta) u = 0$ & Hypoelliptic PDE  
		\\ \addlinespace[5pt]
		& $\R$ & $ \pm \sqrt{\Delta}$ & Half-wave eq. $(- \ri \partial_{t} \pm \sqrt{\Delta}) u = 0$ & Hyperbolic PDE 
		\\ 
        $\re^{- \ri t P}$ & $\R$ & $H$ & Schr\"{o}dinger eq. $(- \ri \partial_{t} - H) u = 0$ & Parabolic PDE
		\\ 
        & $\R$ & $H_{D}$ & Dirac eq. $D u = 0$ & Hyperbolic PDE
		\\ \bottomrule[1pt]
    \end{tabular}
    \caption[Different choices of $U_{t}$]{Different 
        choices of $U_{t}$ with corresponding domains of $t$ and feasible $P$ 
        (see e.g.~\cite{Hoermander_ActaMath_1968, Ivrii_I_Springer_2019}. 
        Here $\Delta, H, H_{D}, D$ are the Laplace-Beltrami operator, the non-relativistic Hamiltonian, the Dirac-Hamiltonian, and the Dirac-type operator, respectively.
    }
    \label{tab: choice_P_U_t}
\end{table} 
%
%
%

\begin{itemize}
    \item 
    $U_{t}$ is the solution operator to the Cauchy problem of the half-wave equation. 
    That is, given an initial data $k \in H^{1} (\varSigma; \sE_{\varSigma} \otimes \halfDen \varSigma)$, one has  
    \begin{equation} \label{eq: U_t_solution_op}
        (- \ri \partial_{t} + P) u = 0 \Rightarrow u = U_{t} (k),  
    \end{equation}
    where $H^{1} (\cdot)$ denotes the Sobolev space (see Section~\ref{sec: result_Feynman_propagator_NHOp}). 
    We can even consider $P$ of order $m \geq 2$ and replace $P$ in the preceding equation by $|P|^{1/m}$ without destroying its hyperbolic nature whilst $- \ri^{m} \partial_{t}^{m} + P$ is no more hyperbolic for $m > 2$. 
    This, in fact, shows our motivation to consider only first-order operators without loss of generality. 
    The primary reason in favour of the strictly hyperbolic operator $- \ri \partial_{t} + P$ compared to a hypoelliptic or a parabolic one is the \textit{finite propagation speed of singularities} of $U_{t}$. 
    \item 
    The Weyl eigenvalue counting function of $P$  
    \begin{equation} \label{eq: def_counting_function}
        \fN^{+} (\lambda) := \# \{ n \,|\, \lambda_{n}^{+} \leq \lambda^{+} \}, 
        \quad 
        \fN^{-} (\lambda) := \# \{ n \,|\, \lambda_{n}^{-} \geq \lambda^{-} \}
    \end{equation}
    is related to the distributional trace (as precised below) $\Tr U_{t}$ via the distributional Fourier transform 
    (see e.g.~\cite[Chap. 3]{Fang_PhD}): 
    \begin{equation} \label{eq: Weyl_counting_function_Tr_U_t}
        \frac{\rd \fN^{\pm}}{\rd \lambda} = \cF_{t \mapsto \lambda}^{-1} ( \Tr U_{t}).   
    \end{equation}
    Therefore, the information contained in $U_{t}$ can be proficiently transferred to the spectral function by means of the Tauberian arguments.  
\end{itemize}
%
%
%

Our choice comes with the caveat that $U_{t}$ is not a trace-class operator as, for instance, $U_{0}$ is the identity operator. 
Thus, $\Tr U_{t}$ makes sense only as a distribution on $\R$: 
\begin{subequations} \label{eq: def_Tr_U_t_distribution}
	\begin{eqnarray}
        && 
        \Tr U_{t} : \cS (\R) \to \C, ~ \rho \mapsto 
        \Tr U_{\rho} := \Tr \Big( \int_{\R} U_{t} \, \Fourier^{-1} (\rho) \, \rd t \Big),  
        \\ 
        && 
        U_{\rho} u := \Big( \int_{\R} U_{t} \, \Fourier^{-1} (\rho) \, \rd t \Big) u := \int_{\R} U_{t} \, \Fourier^{-1} (\rho) \, u \, \rd t, 
        \label{eq: def_smoothed_Cauchy_evolution_op}
        \\ 
        && 
        (\Fourier^{-1} \rho) (t) = \frac{1}{2 \pi} \int_{\R} \re^{\ri t \lambda} \rho (\lambda) \, \rd \lambda, 
	\end{eqnarray}
\end{subequations}
for any $u \in \CComInfinity{\varSigma; \sE_{\varSigma} \otimes \halfDen \varSigma}$, where the rightmost integral in~\eqref{eq: def_smoothed_Cauchy_evolution_op} is defined as a vector-valued Riemann integral.  
Furthermore, it is a $1$-parameter group of unitary operators (see Section~\ref{sec: spectral_theory_L}) on $L^{2} (\varSigma; \sE_{\varSigma} \otimes \halfDen \varSigma)$ and a bounded operator on $H^{s} (\varSigma; \sE_{\varSigma} \otimes \halfDen \varSigma)$, uniformly in $t \in \R$ and for any $s \in \R$ 
(see~\cite[Sec. 3]{Hoermander_ActaMath_1968} 
for the scalar, and, for 
instance,~\cite[\S $20$]{Shubin_Springer_2001} 
for an elaboration). 
\newline 

Therefore, it makes sense to talk about the distributional trace of $U_{t}$ on $L^{2} (\varSigma; \sE_{\varSigma} \otimes \halfDen \varSigma)$ and by the functional calculus $\Tr U_{t} = \sum_{n} \exp (- \ri t \lambda_{n}^{\pm})$. 
Equivalently, we can compute it deploying the Schwartz kernel $\fU_{t}$ of $U_{t}$. 
Putting both together, we arrive at the \textit{quintessential feature of  trace formulae} 
(see e.g.~\cite[Intro.]{Scott_OUP_2010})
\begin{equation} \label{eq: def_trace_formulae} 
    \sum_{n} \re^{- \ri t \lambda_{n}^{\pm}} = \Tr U_{t} = \int_{\varSigma} \tr \big( \fU_{t} (x, x) \big),   
\end{equation} 
where $\fU_{t} (x, x)$ is the restriction of the Schwartz kernel $\fU_{t} (x, y)$ of $U_{t}$ to the diagonal in $\varSigma \times \varSigma$ and $\tr$ is the endomorphism trace (see~\eqref{eq: def_endo_trace}). 
This unveils the underpinning of the bridge connecting the spectral data of $P$ to the geometric world where this operator acts on, as glimpsed in the Poisson summation formula~\eqref{eq: Poisson_summation_formula_derivative_circle} where $P = \rD$ and $\varSigma = \bbS$. 
\newline 

The geometry ($\bbS$) considered in~\eqref{eq: Poisson_summation_formula_derivative_circle} was simple enough to compute the spectrum and the Schwartz kernel explicitly, and hence one has an identity.  
But, this is an atypical situation in a generic compact manifold.  
Therefore, one can, at best, approximate $\fU_{t}$ employing the theory of 
hyperbolic\footnote{The 
    situation is opposite for the non-hyperbolic methods (resolvent, $\zeta$-function, heat-kernel) where the PDE part is simpler but the Tauberian part always fails to provide a decent remainder estimation. 
    In fact, the \textit{best} remainder estimations have been achieved employing hyperbolic methods~\cite{Hoermander_ActaMath_1968, Duistermaat_InventMath_1975} 
    (see also the monograph~\cite{Ivrii_I_Springer_2019} and references therein).
} 
differential equations, usually by constructing parametrices (see Definition~\ref{def: parametrix_PsiDO}), and then obtain the Tauberian part in a straightforward fashion. 
This leads to the notion of asymptotic trace formulae, where one investigates large eigenvalues only by shifting the test function $\rho (\cdot) \mapsto \rho_{\mu} (\cdot) := \rho (\cdot - \mu)$ by an amount $\mu$ and letting $\mu \to \infty$. 
In particular, if $\{ \phiup_{n} \}$ are orthonormal eigensections of $P$ on $\sE_{\varSigma}$, then 
(see Section~\ref{sec: adjoint} for the precise meaning of an anti-dual bundle $\bar{\sE}_{\varSigma}^{*} \to \varSigma$)
\begin{equation}
    \fU_{\rho} (x, y) = \int_{t \in \R} \sum_{n \in \N} \re^{- \ri t \lambda_{n}^{\pm}} (\bar{\phiup}_{n}^{*} \boxtimes \phiup_{n}) (x, y) \, (\Fourier^{-1} \rho) (t) \, \rd t,  
\end{equation}
where the sum converges in $C^{\infty} \big( \varSigma \times \varSigma; (\bar{\sE}_{\varSigma}^{*} \boxtimes \sE_{\varSigma}) \otimes \halfDen (\varSigma \times \varSigma) \big)$ 
(see e.g.~\cite[p. 133]{Grigis_CUP_1994} for the scalar version), 
and the \textit{asymptotic trace formula} reads  
\begin{equation} \label{eq: def_asymptotic_trace_formulae_distribution}
    \Tr (U_{\rho_{\mu}}) = \int_{\varSigma} \tr \big( \fU_{\rho_{\mu}} (x, x) \big), 
    \quad \mathrm{as} \quad \mu \to \pm \infty. 
\end{equation}

The field of asymptotic (semi-classical) trace formulae stems from the study of Green's operator for a Schr\"{o}dinger operator in the limit of vanishing Planck's constant by  
Gutzwill -er~\cite{Gutzwiller_JMP_1971} 
(see also~\cite[Ch. 17]{Gutzwiller_Springer_1990}).   
This seminal work is not entirely rigorous; see for instance, the 
expositions~\cite{Uribe_Cuernavaca_1998, Muratore-Ginanneschi_PR_2003} 
for a scrutinised discussion of his original idea.   
A number of mathematically diligent 
proofs~\cite{Verdiere_CompositioMathematica_1973, *Verdiere_CompositioMathematica_1973_II, Chazarain_InventMath_1974, Duistermaat_InventMath_1975, Chazarain_CPDE_1980, Albeverio_CMP_1982, Brummelhuis_CMP_1991, Meinrenken_ReptMathPhys_1992, Paul_JFA_1995, Combescure_CMP_1999} 
have been reported since then.  
Amongst these, 
Chazarain~\cite{Chazarain_InventMath_1974, Chazarain_CPDE_1980} (for Laplace-Beltrami operator on a closed Riemannian manifold)
and 
Duistermaat-Guillemin~\cite{Duistermaat_InventMath_1975} (for a scalar and positive $P$ with other assumptions as in the on going discussion)
deployed the global theory of Fourier integral operators to derive the complete singularity structure of the wave-trace. 
We refer to the  
monographs~\cite[Chap. XXIX]{Hoermander_Springer_2009},~\cite[Chap. 1]{Safarov_AMS_1997},~\cite[Chap. 11]{Guillemin_InternationalP_2013} 
for details.  
%
%
%
%
%
%
%
%
%
%
\subsection{Duistermaat-Guillemin-Gutzwiller trace formula}
\label{sec: Duistermaat_Guillemin_trace_formula}
To obtain a concrete understanding, let $\Delta$ be the Laplacian (in Geometers' convention) on a $d-1$-dimensional closed Riemannian manifold $(\varSigma, \fh)$.  
It is well-known that its spectrum $\spec \Delta$ comprises discrete eigenvalues $\lambda_{n}^{2}$ and that the counting function $\fN_{\Delta} (\lambda)$ satisfies the celebrated 
Weyl law~\cite{Hoermander_ActaMath_1968} 
\begin{equation} \label{eq: Weyl_law_Laplacian}
	\fN_{\Delta} (\lambda) = \left( \frac{\lambda}{2 \pi} \right)^{d-1} \vol (\bar{\B}^{d-1}) \, \vol (\varSigma), \quad \textrm{as} \quad \lambda \to \infty, 
\end{equation}
where $\vol (\bar{\B}^{d-1})$ is the volume of the closed unit ball $\bar{\mathbb{B}}^{d-1}$ in $\R^{d-1}$.  
As aforementioned, the counting measure $\rd \fN_{\Delta} / \rd \lambda$ is related to the distributional trace $\Tr U_{t}$ via the Fourier transform~\eqref{eq: Weyl_counting_function_Tr_U_t}, where 
\begin{equation} \label{eq: def_U_t_Laplacian}
    \forall t \in \R : U_{t} := \re^{- \ri t \sqrt{\Delta}}   
\end{equation}
is the solution map to the Cauchy problem of the half-wave operator $- \ri \partial_{t} + \sqrt{\Delta}$. 
Moreover, the singularity analysis of $\Tr U_{t}$ at $t=0$ gives rise to the half-wave trace invariants and its singular support is contained in the set of periods of periodic geodesics $\gamma$ on $\varSigma$. 
Furthermore, $\Tr U_{t}$ admits a singularity expansion around the non-zero periods $t = T \neq 0$ (under some technical assumptions related to the clean intersection (see Definition~\ref{def: clean_intersection}) condition between certain  manifolds)~\cite{Chazarain_InventMath_1974, Duistermaat_InventMath_1975} 
and its leading-order term can be expressed by means of the fundamental periods and eigenvalues of the Poincar\'{e} maps of $\gamma$'s~\cite[Thm. 4.5]{Duistermaat_InventMath_1975}.  
\newline 

In their seminal work, 
Duistermaat and Guillemin~\cite{Duistermaat_InventMath_1975} 
have actually considered $U_{t} := \re^{-\ri t \sqrt[m]{P}}$ where $P$ is a positive, elliptic, symmetric pseudodifferential operator of order $m>0$ on a compact boundaryless manifold $\varSigma$. 
\newline 

The Gutzwiller trace formula is structurally similar to that of Duistermaat-Guillemin but the asymptotic is in the limit of vanishing Planck's constant. 
To be precise, take $\varSigma = \R^{d-1}$ and $\sE$ a trivial $\R$-line bundle for simplicity. 
Let $H := \hbar^{2} \Delta + V$ be the non-relativistic Hamiltonian so that $U_{t, \hbar} := \re^{- \ri t H / \hbar}$, where $\hbar$ is some real parameter and $V$ is any appropriate smooth potential on $\R^{d-1}$. 
If $\hbar$ is small enough then the spectrum of $H$ on $L^{2} (\R^{d-1})$ is discrete and there exists an orthonormal basis $\{ \upphi_{n} \}$ for $L^{2} (\R^{d-1})$ of eigenfunctions of $H$. 
Let $\tE_{n, \hbar}$ be the associated eigenvalues of $H$. 
Then one seeks asymptotic expansion of $\sum_{n} \rho ( \nicefrac{\tE_{n, \hbar} - \tE}{\hbar})$ in the limit $\hbar \to 0$. 
Physically, $\hbar$ is identified with the reduced Planck's constant so that $\hbar \to 0$ can be reckoned as the quantum-to-classical transition and hence the name semi-classical trace formula 
(see e.g.~\cite[Sec. 2]{Uribe_Cuernavaca_1998}). 
In this thesis, however, we have only pursued along the tenet of the Duistermaat-Guillemin trace formula and set $\hbar = 1$ throughout our exploration. 
Note, $\re^{- \ri t \Delta}$ is \textit{not} a trace-class operator since, for instance, on a circle its singular support is the entire real axis~\cite[p. 46]{Duistermaat_InventMath_1975}.  
%
%
%
%
%
%
%
%
%
%
\subsection{Sandoval trace formula} 
\label{sec: Sandoval_trace_formula}
Sandoval has generalised the Duistermaat-Guillemin trace formula for a Dirac-type operator $\hat{D}$ on a hermitian vector bundle $(\sE_{\varSigma} \to \varSigma, \langle \cdot|\cdot \rangle)$ over a closed Riemannian manifold $(\varSigma, \fh)$ by studying 

\begin{equation} \label{eq: def_U_t_Riem_Dirac}
    \forall t \in \R : \hat{U}_{t} := \re^{- \ri t \hat{D}}.  
\end{equation}
Arguably, $\hat{D}$ is the most fundamental first-order differential operator in geometric analysis and it differs from a Laplace-Beltrami operator $\Delta$ on $\sE_{\varSigma}$ in that its spectrum $\spec \hat{D}$ is unbounded from both below and above whereas $\spec \Delta$ is bounded from the below (semi-bounded). 
Sandoval's investigation results in a straightforward generalisation of 
the Weyl law~\cite[$(5)$]{Sandoval_CPDE_1999}  
and Dirac-wave trace invariants~\cite[Thm. 2.2]{Sandoval_CPDE_1999}. 
But, genuine bundle features show up in terms of the holonomy group in the singularity analysis of $\Tr U_{T \neq 0} (\hat{D})$~\cite[Thm. 2.8]{Sandoval_CPDE_1999}. 
However, there is an ad-hoc choice regarding the connection used to deduce the holonomy group as elucidated below. 
\newline 

As a matter of fact, the Weitzenb\"{o}ck connection (see~\eqref{eq: Weitzenboeck_formula_Dirac_type}) $\hat{\nabla}$ determined by $\hat{D}$ and $\fh$ does not automatically induce a Clifford connection albeit such a connection always exists (see Remark~\ref{rem: Clifford_connection}) and the corresponding operator is usually called a \textit{compatible Dirac-type operator} $\check{D}$. 
They differ only by a smooth term $W := \hat{D} - \check{D} \in C^{\infty} (\varSigma; \End \sE)$. 
To compute the principal symbol $\symb{\hat{U}_{t}}$ of $\hat{U}_{t}$, Sandoval used the parallel transporter $\check{\cT}$ with respect to the Weitzenb\"{o}ck connection $\check{\nabla}$ of $\check{D}$ instead of $\hat{\cT}$ corresponding to $\hat{\nabla}$~\cite[Prop. 5.3]{Sandoval_CPDE_1999}.  
Moreover, she employed the trivialisation characterised by the vanishing subprincipal symbol $\subSymb{\check{D}^{2}}$ of $\check{D}^{2}$~\cite[Cor. 5.11]{Sandoval_CPDE_1999}. 
As a consequence, $\check{\cT}$ appearing in expression of $\symb{\hat{U}_{t}}$ is somewhat expedient and the splitting of $\hat{\cT}$ in terms of $\check{\cT}$ and the 
average value\footnote{See~\cite[$(16)$]{Sandoval_CPDE_1999}  
	for the precise expression. 
} 
of $W$ along the geodesic $[0, t] \ni s \mapsto \gamma (s) \in \mathbb{S}^{*} \varSigma$ is rather ad-hoc.
\newline 

The objective of our study is twofold. 
First, we have relaxed all the aforementioned ad-hoc considerations to compute $\symb{\hat{U}_{t}}$.   
Next, one notices that $\hat{U}_{t}$ solves the initial value problem of the differential equation 
\begin{equation} \label{eq: Dirac_wave}
	(- \ri \partial_{t} + \hat{D}) u = 0, \quad u (t_{0}) = k  
\end{equation}
on a $d$-dimensional product manifold  (in the context of general relativity, this is known as the \textit{ultrastatic spacetime};  
see Remark~\ref{rem: SST}~\ref{rem: hypersurface_orthogonality_SST})
\begin{equation} \label{eq: def_ultrastatic_spacetime}
	\spacetime := (\R \times \varSigma, \rd t^{2} - \fh)  
\end{equation}
admitting (at least) a global timelike Killing vector field $Z := \partial_{t}$ so that the external parameter $t$ in $\hat{U}_{t}$ can be considered as the canonical global time-coordinate on $\R$. 
But, $- \ri \partial_{t} + \hat{D}$ does not respect general relativistic covariance and so it seems more natural to instead work with the Lorentzian Dirac-type operator $D$. 
Thus the study is devoted to investigating a Lorentzian generalisation of Sandoval's trace formula.
%
%
%
%
%
%
%
%
%
%
\subsection{Riemannian trace formulae: A relativistic viewpoint}
\label{sec: DGTF_ultrastatic}
The Riemannian to Lorentzian transition is conceptually non-trivial because the solution map to the Cauchy problem of the Dirac equation 
\begin{equation} \label{eq: Dirac_eq}
	D u = 0, \quad u \upharpoonright \varSigma = k  
\end{equation} 
on a globally hyperbolic spacetime (see Definition~\ref{def: globally_hyperbolic_spacetime}) is no longer of the form~\eqref{eq: def_U_t_Riem_Dirac}, since the Dirac-Hamiltonian (the counterpart of $\hat{D}$) becomes time-dependent. 
In order to contemplate what should be looked for in a Lorentzian setting, we revisit the product scenario from a relativistic viewpoint. 
One reckons that the solutions of~\eqref{eq: Dirac_wave} are of the following form 
\begin{equation}
	\uppsi_{n}^{\pm} := \re^{- \ri t \lambda_{n}^{\pm}} \upphi_{n}^{\pm}, 
\end{equation}
where $\upphi_{n}^{\pm}$ are the eigensections of $\hat{D}$ corresponding to the eigenvalues $\lambda_{n}^{\pm}$. 
But 
\begin{equation} \label{eq: eigenvalue_KVF}
	Z \uppsi_{n}^{\pm} = - \ri \lambda_{n}^{\pm} \uppsi_{n}^{\pm}  
\end{equation}
and consequently, 
\begin{equation} \label{eq: trace_time_evolution_op_tr_Killing_flow}
    \Tr_{L^{2} (\varSigma; \sE)} \hat{U}_{t}  
	= \sum_{\lambda_{n}^{\pm} \in \spec \hat{D}} \re^{- \ri t \lambda_{n}^{\pm}} 
	= \Tr_{\ker (- \ri \partial_{t} + \hat{D})} \re^{tZ}. 
\end{equation}
In other words, the \textit{trace of the Dirac-wave group $\hat{U}_{t}$ on a closed Riemannian manifold is equivalent to the trace of the flow $\re^{tZ}$ induced by the global timelike Killing vector field $Z$ of an ultrastatic spacetime $\sM$, acting on the kernel $\ker (- \ri \partial_{t} + \hat{D})$ of the half-wave operator $- \ri \partial_{t} + \hat{D}$}. 
This is essentially a straightforward generalisation of the observation originally due to 
Strohmaier and Zelditch~\cite[Sec. 10.1]{Strohmaier_AdvMath_2021} 
in the context of a scalar wave operator.    
Thus, the deep results on Riemannian spectral geometry as mentioned in Section~\ref{sec: Duistermaat_Guillemin_trace_formula} and~\ref{sec: Sandoval_trace_formula}   
(and many more as available in the 
compendia~\cite{Safarov_AMS_1997, Hoermander_Springer_2009, Ivrii_BullMathSci_2016, Avetisyan_JST_2016, Li_JGP_2016}, 
for instance) 
are essentially achieved by investigating the solution operator of the half-wave equation on an ultrastatic spacetime with a compact Cauchy hypersurface (see Definition~\ref{def: Cauchy_hypersurface}) without boundary.    
\newline 

It is, therefore, evident that a global timelike Killing flow (see~\eqref{eq: def_Killing_flow}) $\varXi$ on a Lorentzian manifold plays a pivotal role and the most general class of spacetimes admitting such a flow is known as the stationary spacetime (see Definition~\ref{def: stationary_spacetime}) $\SST$.  
Physically speaking, these spacetimes admit a \textit{canonical flow of time} but unlike the ultrastatic case there is \textit{no} preferred time coordinate. 
They are interesting because a number of exact solutions of Einstein's equation (Schwarzschild and Kerr black holes) belong to this class. 
We demand the global hyperbolicity condition on $\SST$ to ensure the well-posed Cauchy problem for~\eqref{eq: Dirac_eq}.  
Therefore, the object of our study is the distributional trace $\Tr U_{t}$ (see~\eqref{eq: Tr_Cauchy_evolution_op_Dirac_Tr_Killing_flow}) of the \textit{time evolution operator} $U_{t}$ on a compact (to guarantee discrete eigenvalues) Cauchy hypersurface $\varSigma \subset \sM$ without boundary. 
In addition to the conceptual issues, our computational techniques (see Section~\ref{sec: proof_strategy} for details) are different from those used by Sandoval, primarily due to the fact that the Dirac Hamiltonian $H_{D}$ (see~\eqref{eq: def_Dirac_Hamiltonian_SSST}) on a stationary spacetime is \textit{not} of Dirac-type in contrast to $\hat{D}$.  
Furthermore, we have computed $\symb{U_{t}}$ (see Lemma~\ref{lem: symbol_Tr_U_0} and~\ref{lem: symbol_Tr_U_T}) using the Weitzenb\"{o}ck connection induced by $D$ in order to avoid the aforementioned expedient choice in Sandoval's work, as explained in the comment after the statement of the foremost Theorem~\ref{thm: trace_formula_t_T}. 
\newline 

As a quintessential relativistic operator, it furthermore allows to study the semi-classical limit of quantum field theoretic observables interacting with classical 
gravity\footnote{We 
	do not require Einstein field equation to be satisfied.
}. 
Whilst spin (see Example~\ref{ex: spin_Dirac_op})- or spin$^{\mathrm{c}}$-Dirac operators are desirable from a physics point of view, we prefer to work on more general Dirac-type operators because they do not enforce any topological restrictions on  
$\sM$\footnote{Global  
	existence of a Lorentz metric depends on the topology of $\sM$. 
	In particular, such a metric exists in all non-compact manifolds and compact manifolds with vanishing Euler characteristics.}
unlike the 
spin\footnote{A 
	spin-structure always exists locally but its global existence depends on some higher orientability property of the base manifold. 
	For instance, spin (resp. spin$^{\mathrm{c}}$) structure exists if and only if the second (resp. the third integral) Stiefel–Whitney class of $\sM$ 
    vanishes~\cite{Haefliger_1956, *Borel_AJM_1959, *Milnor_Enseign_1963, *Karoubi_1968}.}-Dirac 
operators. 
\newline 

We would like to close this subsection by mentioning that this is not the only literature on Lorentzian trace formulae rather the first step towards this direction has been taken up for the d'Alembertian operator (see Example~\ref{exm: covariant_Klein_Gordon_op}) on the same spacetime geometry as in this thesis by 
Strohmaier and Zelditch~\cite{Strohmaier_AdvMath_2021}.  
We are going to describe their work briefly in the next section. 
Therefore, the thesis (see Section~\ref{sec: result_Gutzwiller_trace} for the primary results) can be viewed as the Lorentzian generalisation of Sandoval's work propounding the framework of Strohmaier-Zelditch into a bundle setting. 
Recently, 
McCormick~\cite{McCormick} 
has also extended Strohmaier-Zelditch trace formula for the vector d'Alembertian on a spatially compact globally hyperbolic stationary spacetime equipped with the Kaluza-Klein metric by employing a few results of this investigation.  
%
%
%
%
%
%
%
%
%
\subsection{Strohmaier-Zelditch trace formula}
\label{sec: Strohmaier_Zelditch_trace_formula}
In the set-up of a globally hyperbolic stationary spacetime $\SST$ with a compact Cauchy hypersurface $(\varSigma, \fh)$, 
Strohmaier-Zelditch~\cite{Strohmaier_AdvMath_2021} 
have considered $\square$, the d'Alembert operator with a smooth potential invariant under the induced Killing flow $\varXi_{s}^{*}$, so that the spectral problem can be stated as  
\begin{equation} \label{eq: joint_eigenfunction_Lie_derivative_NHOp}
    \square \psiup_{n} = 0, \quad Z \psiup_{n} = \lambda_{n} \psiup_{n}, 
\end{equation}
where $\lambda_{n}$ are the eigenvalues of the (timelike) complete Killing vector field $Z$ (see Remark~\ref{rem: SST}~\ref{rem: complete_Killing_flow}) on $\sM$ and $\psiup_{n}$ are the joint eigenfunctions of $\square$ and $Z$.  
This eigenvalue problem is equivalent to that of an operator pencil of  
certain\footnote{See~\cite[$(29)$]{Strohmaier_AdvMath_2021}  
    for the precise form.
}  
form which entails that the eigenvalue problem cannot be interpreted for some operator on $\varSigma$. 
Furthermore, it reduces to the classical eigenvalue problem $(\Delta - \lambda^{2}) \upphi_{n} = 0$ only when the spacetime metric $\fg$ is of the product form~\eqref{eq: def_ultrastatic_spacetime}, where $\Delta$ is the Laplacian on $\varSigma$. 
\newline 

In the relativistic setting, the classical phase space is the conic symplectic manifold $\sN$ of lightlike (co)geodesics on the cotangent bundle $\coTansM$ of $\sM$ in contrast to the cotangent bundle $\coTanCauchy \to \varSigma$ over a Cauchy hypersurface $\varSigma$, as expected from the non-relativistic scenario (see Section~\ref{sec: classical_dyanmics} for details). 
The trace $\Tr U_{t}$ of the quantum dynamics, governed by the time-evolution operator $U_{t}$ of $\square$, is given by (cf.~\eqref{eq: trace_time_evolution_op_tr_Killing_flow})
\begin{equation}
    \Tr U_{t} = \Tr_{\ker \square} \re^{t Z}. 
\end{equation}
Then, tailoring 
Duistermaat-Guillemin~\cite{Duistermaat_InventMath_1975} framework in this Lorentzian set-up, they have achieved asymptotic expressions for $\Tr U_{t=0}$ and $\Tr U_{T \neq 0}$, and have subsequently obtained a Lorentzian Weyl law~\cite{Strohmaier_AdvMath_2021}  
(see also~\cite{Strohmaier_IndagMath_2021} 
and the 
review~\cite{Strohmaier_RMP_2021}). 
%
%
%
%
%
%
%
%
%
%
\section{Feynman propagators}
\label{sec: intro_Feynman_propagator}
On a globally hyperbolic spacetime, the Cauchy problem for any normally hyperbolic operator (see Definition~\ref{def: NHOp}) $\square$ is well-posed. 
This implies the existence of unique advanced and retarded fundamental solutions (propagators) (see Definition~\ref{def: fundamental_solution} and Remark~\ref{rem: exist_unique_advanced_retarded_Green_op_NHOp}).  
The properties of these fundamental solutions are extremely important for the understanding of classical wave propagation, such as the electromagnetic waves. 
They appear naturally because of causality: the retarded fundamental solution propagates to the future, whereas the advanced fundamental solution propagates to the past. 
In quantum field theory, the appearance of time-ordering (see Appendix~\ref{sec: timeordered_product}) and the enforcement of positivity of energy have led to the development of another type of fundamental solution ---  the Feynman propagator.  
It propagates \textit{positive energy solutions to the future} and \textit{negative energy solutions to the past}, and thus combines causality with the notion of positivity of energy 
(see Example~\ref{ex: Feynman_propagator_Minkowski}). 
The explicit formulae of Lorentz invariant Feynman propagators for normally hyperbolic operators in Minkowski spacetime is available in the standard textbooks of quantum field theory. 
As recalled in Example~\ref{ex: Feynman_propagator_Minkowski}, the usual construction employs the Fourier transform. 
However, on a generic globally hyperbolic spacetime there is \textit{neither a Fourier transform nor any reasonable notion of energy in the absence of a global timelike Killing vector field}. 
A priori, it is, therefore unclear what a Feynman propagator should be.
\newline 

In the theory of partial differential equations, the notion of a parametrix (see Definition~\ref{def: parametrix_PsiDO}) is often useful in the first stage on the  construction of a true fundamental solution. 
A parametrix, per se, is an inverse of the operator $\Box$ modulo smoothing operators. 
Parametrices are considered equivalent if they differ by smoothing operators.
It was a deep insight of 
Duistermaat and H\"{o}rmander~\cite{Duistermaat_ActaMath_1972}   
that there is a well-defined notion of Feynman parametrices and these parametrices are unique up to smoothing operators. 
In other words, they are \textit{unique as parametrices}. 
In fact, the notion of distinguished parametrices for any scalar pseudodifferential operator of real principal type (see Definition~\ref{def: real_principal_type_PsiDO_mf}) is given in their seminal article where they identified a geometric notion of pseudoconvexity which allows one to prove uniqueness of such parametrices. 
Roughly speaking, despite not not having any notion of energy due to the lack of a global timelike Killing vector field, there is still a microlocal notion of positivity of energy and a corresponding flow on the cotangent bundle. 
This microlocal notion can be used to characterise Feynman parametrices (see Definition~\ref{def: Feynman_parametrix}): they are distinguished by their wavefront sets rather than the support properties.
\newline

Feynman parametrices play an extremely important role in quantum field theory on curved spacetimes, and the theory has actually been developed to a certain extent first in the physics literature.
Canonical quantisation of linear fields can be done in two stages. 
In the first step (see Appendix~\ref{sec: field_algebra}), one constructs a field algebra from the space of solutions of the 
respective equation of motion
and next (see Appendix~\ref{sec: quantum_state}), some quantum state is required to construct a Hilbert space representation of this field algebra. 
In Minkowski spacetime these two steps are usually combined into one owing to the existence of the vacuum state (see Example~\ref{ex: Minkowski_vacuum}), while it is more fruitful to separate them in curved spacetimes. 
The first step can be done without any problems, exactly the same way as in Minkowski spacetime but the second step necessitates a notion of a reasonable state. 
It has been realised that a state compels certain conditions in order to perform the usual operations in 
perturbative quantum field theory~\cite{Brunetti_CMP_2000, Hollands_CMP_2001, Hollands_CMP_2002}. 
One of the identified conditions is a restriction of the type of singularity that one obtains from the state, the so-called Hadamard condition (Definition~\ref{def: Hadamard_state}). 
Although Duistermaat and H\"ormander were certainly aware of the developments in physics, it was realised only much later by 
Radzikowski~\cite{Radzikowski_CMP_1996} 
that the expectation values of the time-ordered products with respect to states satisfying the Hadamard condition are Feynman propagators. 
In fact, the construction of a Hadamard state is equivalent to the construction of a Feynman propagator that satisfies a certain (see Proposition~\ref{prop: positivity_Feynman_minus_adv_NHOp}) positivity property. 
The fact that this positivity property holds for parametrices was already shown by Duistermaat and H\"{o}rmander. 
This property for parametrices implies the existence of such Feynman propagators and hence of Hadamard states 
(see~\cite{Lewandowski_JMP_2022}). 
\newline 

Apart from its applications in physics, Feynman parametrices also play an enormously important role in mathematics.
For instance, they appear in the context of Lorentzian index theory as inverses modulo compact operators of the Dirac operator with Atyiah-Patodi-Singer boundary conditions~\cite{Baer_AJM_2019}, 
local index theorem~\cite{Baer},   
and asymptotically static spacetimes~\cite{Shen}.  
These concepts also arise in Vasy's treatment of asymptotically hyperbolic problems, see for example~\cite{GellRedman_CMP_2016}  
and references therein, where even a non-linear problem is discussed in this context.
Moreover, the notion of Feynman propagator also comes up naturally in the Lorentzian generalisation of the Duistermaat-Guillemin-Gutzwiller trace formula as observed first in~\cite{Strohmaier_AdvMath_2021}.  
\newline

So far, Duistermaat and H\"ormander's construction of distinguished parametrices for scalar pseudodifferential operators of real principal type has not been used much in the physics literature probably due to its great generality and the associated complex notation. 
Their main idea is, however, compelling and simple. 
The operator can be conjugated microlocally to a vector field and it is therefore sufficient to construct a parametrix  for the operator of differentiation on the real line. 
There are two distinguished fundamental solutions for the operator of differentiation. 
A choice of parametrix for each connected component of the characteristic set of the operator results in a distinguished parametrix. 
Thus, there are $2^N$ distinguished parametrices for such a pseudodifferential operator if $N$ is the number of connected components of its characteristic set. 
For the wave operator, this gives $4$ distinguished parametrices in dimensions $d \geq 3$ and $16$ distinghuished parametrices in $d=2$. 
\newline 

The aim of this investigation is threefold. 
First, we would like to revise and simplify the construction of 
Duistermaat-H\"ormander~\cite{Duistermaat_ActaMath_1972} 
in the special case of a normally hyperbolic operator. 
The second aim is to fill a gap in the literature: microlocalisation and the corresponding construction of distinguished parametrices in Duistermaat-H\"ormander is covered in the literature only for scalar operators.
Several constructions in index theory, in trace formulae, and also in physics require the existence and uniqueness of Feynman parametrices for operators acting on vector bundles. 
It is known that most of the constructions carry over to the case of any  normally hyperbolic operator, since its principal symbol is a scalar. 
There are, however, also important differences that appear on the level of subprincipal symbols (see Definition~\ref{def: subprincipal_symbol}). 
In this study, we would like to give a precise statement (see Theorems~\ref{thm: microlocalisation_P} and~\ref{thm: microlocalisation_NHOp})  of microlocalisation for these class of geometric operators. 
We then provide a detailed construction (Theorem~\ref{thm: exist_unique_Feynman_parametrix_NHOp}) of Feynman parametrices for vector bundles with complete proofs for their uniqueness and discuss the effect of curvature of the bundle connection. 
We show (see Proposition~\ref{prop: positivity_Feynman_minus_adv_NHOp}) that the construction can be carried out in such a way that the above mentioned positivity property holds, provided that the vector bundle has a hermitian inner product with respect to which the operator is formally selfadjoint. 
Subsequently, we promote the Feynman parametrices to the Feynman propagators (see Theorem~\ref{thm: existence_Feynman_propagator_NHOp}) by employing the well-posed Cauchy problem for a normally hyperbolic operator on a globally hyperbolic spacetime. 
Third and finally, for any Dirac-type operator, we give (see Theorem~\ref{thm: Hadamard_bisolution_Dirac_type_op}) a much more direct construction of Feynman propagators satisfying a positivity property.
We also discuss some consequences, including the usual propagation of singularity theorem (see Theorem~\ref{thm: propagation_singularity_Sobolev_WF}). 
These are well known to hold for very general operators on vector bundles. 
For instance, the propagation of polarisation sets has been proven by Dencker using 
microlocalisation in the matrix setting~\cite{Dencker_JFA_1982}. 
\newline 

As explained above (and see~\eqref{eq: def_Green_kernel_Feynman_timeordered_product}), the construction of Feynman propagators satisfying the positivity property is equivalent to the construction of Hadamard states.
There are several constructions of Hadamard states, even in the 
analytic category~\cite{Gerard_CMP_2019}. 
Amongst the methods to construct them, there are direct ones using singularity expansions employing the 
Lorentzian distance function 
(the so-called Hadamard expansion)~\cite{Fulling_CMP_1978, Brown_JMP_1984, Sahlmann_RMP_2001, Marecki_master, Hollands_RMP_2008, Dappiaggi_RMP_2009, Lewandowski_JMP_2022},  
spectral methods that rely on frequency splitting 
(deformation method) \cite{Fulling_AnnPhys_1981, Murro_AGAG_2021}, 
pseudodifferential methods~\cite{Junker_RMP_1996, Hollands_adiabatic_CMP_2001, Gerard_CMP_2014, Gerard_CMP_2015, Gerard, Gerard_bounded_geometry} 
(see also the monograph~\cite{Gerard_EMS_2019}), 
holography~\cite{Moretti_CMP_2008, Dappiaggi_JMP_2009, Dappiaggi_ATMP_2011, Gerard_AnalPDE_2016} 
(see also the exposition~\cite{Dappiaggi_Springer_2017}), 
and even 
global methods~\cite{GellRedman_CMP_2016, Vasy_AHP_2017, Derezinski_RMP_2018, Capoferri_JMAA_2020} 
(see also the review~\cite{Avetisyan_Math_2021}).  
Many of these constructions can and have been generalised to the bundle case. 
Each comes with their advantages and disadvantages.
%
%
%
%
%
%
%
%
%
%
\section{Organisation}
The thesis consists of three chapters and three appendices. 
Chapter~\ref{ch: FIO} essentially contains the background materials required for the thesis,  while our genuine contributions have been placed in Chapters~\ref{ch: Feynman_propagator} and~\ref{ch: Gutzwiller_trace}. 
Appendices~\ref{ch: canonical_relation_composition} and~\ref{ch: symbol} primarily supplement Chapter~\ref{ch: FIO} and substantially Chapter~\ref{ch: Gutzwiller_trace}, whereas Chapter~\ref{ch: Feynman_propagator} is complemented by Appendix~\ref{ch: Hadamard_state}.  
\newline 

Fourier integral operators are the prime technical tool used in this thesis. 
Specifically, the symbolic viewpoint has been utilised heavily. 
We have begun Chapter~\ref{ch: FIO} with the local theory to motivate the global formulation.  
To make the local-to-global transition smooth, first the naive theory of Fourier integral operators on Euclidean spaces has been given in a pedagogical fashion, and then an invariant local formalism has been presented. 
This characterisation paves the way to generalise the theory of Lagrangian distributions on manifolds and subsequently on vector bundles. 
Besides reviewing the standard formulae, we have proven bundle version of a few theorems usually available only for the scalar case. 
Although, these proofs follow the standard ones with minor modifications, they are not easy to find in the existing literature in a systematic way. 
\newline 

The subject matter of Chapter~\ref{ch: Feynman_propagator} is the Feynman propagators for any normally hyperbolic operator and any Dirac-type operator  on a globally hyperbolic spacetime. 
To keep the discourse self-contained, we start with explaining the necessary backgrounds on Lorentzian geometry and on these operators. 
Our primary findings are presented afterwards in Section~\ref{sec: result_Feynman_propagator_NHOp},~\ref{sec: result_Feynman_propagator_Dirac}, and~\ref{sec: Rarita_Schwinger_op}. 
In order to derive these results, we have introduced the notion of a $P$-compatible connection (see Definition~\ref{def: P_compatible_connection}) induced by the subprincipal symbol of a pseudodifferential operator $P$. 
This enables us to formulate (see Theorem~\ref{thm: microlocalisation_P} and Remark~\ref{rem: microlocalisation_real_principal_type}) the microlocalisation of a pseudodifferential operator of real-principal type on a vector bundle. 
Consequently, normally hyperbolic operators are microlocalised (see Theorem~\ref{thm: microlocalisation_NHOp}).  
\newline 

The trace of the time-flow on the kernel of any Dirac-type operator on a globally hyperbolic stationary spacetime has been investigated in Chapter~\ref{ch: Gutzwiller_trace}. 
This chapter is comprised adequate background on the geometry of stationary spacetimes followed by our primary findings in Section~\ref{sec: result_Gutzwiller_trace}. 
\newline 

Appendix~\ref{ch: canonical_relation_composition} summarises the composition of vector bundle-valued half-densities on homogeneous canonical relations and relevant materials on conic symplectic geometry. 
\newline 

A write up on vector bundle-valued polyhomogeneous symbol class on a canonical relation has been offered in Appendix~\ref{ch: symbol}. 
\newline 

Appendix~\ref{ch: Hadamard_state} contains a non-technical prelude to quantum field theory in curved spacetime in the algebraic approach and a technical introduction to algebraic quantum states. 
To relate this approach with the standard physics literature, we have used the Minkowski vacuum as an elucidating example. 
This also motivates the demand of imposing the Hadamard condition as a necessary criterion to select a physical state in curved spacetime.  
%
%
%
%
%
%
%
%
%
%
\section{Conventions} 
\label{sec: convention}
Throughout the thesis, a vector bundle $\sE \to M$ means a smooth $\C$-vector bundle over a Hausdorff, second countable $d \in \N$-dimensional topological space $M$ furnished with an equivalence class of smooth atlases. 
We use the notation $\dot{\sE}$ to symbolise the zero-section removed part of $\sE$. 
Occasionally, $\sE$ is endowed with a non-degenerate sesquilinear form $(\cdot|\cdot)$ which is assumed to be anti-linear in its first argument. 
By a hermitian form $\langle \cdot | \cdot \rangle$ on $\sE$, we mean a positive-definite $(\cdot|\cdot)$.
A Lorentzian manifold and its special case - a globally hyperbolic manifold, both are denoted by $\spacetime$ with metric signature $+ - \ldots -$,  and $d := \dim \sM \geq 2$.  
\newline 

Let $\halfDen \to M$ be the bundle of half-densities over $M$. 
We denote the vector spaces of smooth and compactly supported smooth half-densities on $\sE$ by $\secE$ and $\comSecE$, and endow them with the Fr\'{e}chet space and the inductive limit topologies, respectively. 
The vector spaces of (compactly supported) distributional half-densities ($\dualSecE := \big( \secEStar \big)'$ resp.) $\dualComSecE := \big( \comSecEStar \big)'$ on $\sE$ are defined by the topological duals of ($\secEStar$ resp.) $\comSecEStar$, which are equipped with the weak $*$-topologies induced by the topologies of ($\secEStar$ resp.) $\comSecEStar$, where $\sE^{*} \to M$ is the dual bundle of $\sE$    
(see e.g.~\cite[p. 307]{Guillemin_AMS_1977},~\cite[pp. 24-25, 146]{Guenther_AP_1988} for details of these spaces).  
On an oriented Lorentzian manifold $\spacetime$, we use the natural Lorentzian volume element to identify $\secsME = C^{\infty} (\sM; \sE \otimes \halfDen \sM)$ and $\dualComSecsME = \cD' (\sM; \sE \otimes \halfDen \sM)$. 
Additionally, one defines the space $C_{\mathrm{sc}}^{\infty} (\sM; \sE)$ of spatially compact smooth sections of $\sE$ as the set of all $u \in \secsME$ for which there exists a compact subset $K$ of $\sM$ such that $\supp{u} \subset J (K)$ where $J (K) := J^{+} (K) \cup J^{-} (K)$ and $J^{\pm} (K)$ are the causal future(past) of $K$; as a vector space $C_{\mathrm{sc}}^{\infty} (\sM; \sE) \subset \secsME$.   
Echoing this spirit $C^{\infty} (M)$ denotes the set of all $\C$-valued smooth functions on $M$. 
\newline 

The Fourier transform $\cF (u)$ of any $u \in L^{1} (\Rd, \rd x)$ is defined by    
\begin{equation} \label{eq: def_Fourier_transformation} 
	\mathcal{F} (u) := \int_{\Rd} \re^{ - \ri x \cdot \theta} u \, \rd x,   
	~\quad~\textrm{and}~\quad  
	u = \frac{1}{(2 \pi)^{d}} \int_{\Rd} \re^{ \ri x \cdot \theta} \cF (u) \, \rd \theta, 
\end{equation}
whenever $\cF (u) \in L^{1} \big( \Rd, \rd \theta/ (2 \pi)^{d} \big)$, and here $\cdot$ is either the Euclidean or the Minkowski inner product depending on the context. 
If $T : \comSecNF \to \dualComSecME $ is an integral operator then its Schwartz kernel is denoted by $\fT$. 
\newline 

The symbol $\sim$ stands for the asymptotic summation and only the polyhomogeneous symbol class $S^{m}$ is utilised in this thesis so that $\Psi \mathrm{DO}^{m}$ (resp. $I^{m}$) is the set of pseudodifferential operators (resp. Lagrangian distributions) having polyhomogeneous total symbol. 
We have used the notation $S^{m - [m']} := S^{m} / S^{m-m'}$ to represent the quotient space (and similarly $\Psi \mathrm{DO}^{m - [m']}, I^{m - [m']}$, etc.). 
We have succinctly written $u \equiv v$ to mean $u - v \in \secE$ for any $u, v \in \dualComSecE$. 
%
%
%

%% file: differential_operator.tex
\chapter{Fourier Integral Operators} 
\label{ch: FIO} 
\textsf{The standard theory of Lagrangian distributions on a vector bundle is recalled in this chapter. 
In order to motivate the global formalism, we begin with the local theory of partial differential operators and generalise in ascending order, i.e., pseudodifferential operators are followed by Fourier integral operators. 
The formulation of pseudodifferential (and Fourier integral) operators on a manifold has been used as an intermediate step to complement the bundle setting.}
%
%
%
%
%
%
%
%
%
%
\section{Local theory of Fourier integral operators}
\subsection{Differential operators}
\label{sec: PDO_Euclidean} 
Let $U$ be an open subset of a $d \in \N$-dimensional Euclidean space $\Rd$. 
A partial differential operator $L$ of order (at most) $m \in \R$ with smooth coefficients on $U$ can be expressed as 
\begin{equation} \label{eq: expression_PDO_Euclidean}
    L = \sum_{|\alpha| \leq m} f_{\alpha} \, \rD^{\alpha} : C^{\infty} (U) \to C^{\infty} (U), 
\end{equation}
for some $f_{\alpha} \in C^{\infty} (U)$ and here the multiindex notation (see Appendix~\ref{ch: symbol}) has been used. 
Employing Fourier transform $\Fourier (u)$ of $u \in \CComInfinity{U}$ one obtains    
\begin{eqnarray} 
    (L u) (x) 
    & = & 
    \int_{U} \fL (x, y) \, u (y) \, \rd y, 
    \label{eq: PDO_Euclidean_Fourier}
    \\ 
    \fL (x, y) 
    & := & 
    \int_{\Rd} \re^{\ri (x - y) \cdot \theta} \fl (x, \theta) \measuretheta, 
    \label{eq: def_expression_kernel_PDO_Euclidean}
    \\ 
    \fl (x, \theta) 
    & := & 
    \sum_{|\alpha| \leq m} f_{\alpha} \theta^{\alpha}, 
    \label{eq: def_coordinate_expression_total_symbol_PDO_Euclidean}
\end{eqnarray}
where $\fL$ and $\fl$ are called the \textbf{Schwartz kernel} and the \textbf{total symbol} of $L$, respectively. 
The leading order term 
\begin{equation} \label{eq: def_coordinate_expression_symbol_PDO_Euclidean}
	l (x, \theta) := \sum_{|\alpha| = m} f_{\alpha} \theta^{\alpha} 
\end{equation}
of the total symbol $\fl$ is called the \textbf{principal symbol} of $L$.  
Let us denote the set of all partial differential operators of order most $m \in \R$ on $U$ by $\PDOU$. 
\newline 

One observes that the total symbol  
\begin{equation} \tag{\ref{eq: def_total_symbol_PDO_Euclidean}'} \label{eq: tentative_def_total_symbol_PDO_Euclidean}
    \totSymb{} : \PDO^{m} (U) \to \bbP_{m} (U \times \Rd), \quad L \mapsto 
    \totSymb{L} := \fl
\end{equation}
is an isomorphism between $\PDOU$ and the set $\bbP_{m} (U \times \Rd)$ of all smooth functions on $U \times \Rd$ which are polynomials of degree at most $m$ in $\Rd$.  
This isomorphism tells that \textit{knowing the total symbol $\fl$ up to smooth terms is equivalent to the knowledge of the operator $L$ modulo smooth $\End \big( C^{\infty} (U) \big)$}, and hence, a sharper version of~\eqref{eq: tentative_def_total_symbol_PDO_Euclidean} is  
\begin{equation} \label{eq: def_total_symbol_PDO_Euclidean}
    \totSymb{} : \PDO^{m- [\infty]} (U) \to \bbP_{m - [\infty]} (\coTan U), \quad [L] \mapsto 
    \totSymb{[L]} := [\fl], 
\end{equation}
where $\coTan U$ is the cotangent bundle over $U$ and $\PDO^{m- [\infty]} (U)$ means $\PDO^{m} (U)$ modulo smooth operators on $U$. 
One also observes that changing next to the leading-order terms in~\eqref{eq: expression_PDO_Euclidean} does not change the principal symbol:  
\begin{equation} \label{eq: def_symbol_PDO_Euclidean}
    \symb{} : \PDO^{m- [1]} (U) \to \bbP_{m - [1]} (\coTan U), \quad [L] \mapsto 
    \symb{[L]} := [l], 
\end{equation}
where $\PDO^{m-[1]} := \PDO^{m} / \PDO^{m-1}$. 
Note that \textit{neither} $\fl$ \textit{nor} $\fL$ is well-behaved under the change of coordinates whilst $l$ is. 
More precisely, let us consider a diffeomorphism 
\begin{equation} \label{eq: diffeo_U_V}
	\kappaup : U \to V, ~ x \mapsto y := \kappaup (x)
\end{equation}
between two open subsets $U, V \subset \Rd$ whose cotangent lift is 
\begin{equation} \label{eq: coTan_diffeo_U_V}
	\coTan \kappaup : \coTan V \to \coTan U, ~ \yeta \mapsto 
    \xxi := \coTan \kappaup \, \yeta 
    := \left( \kappaup^{-1} (y), \coTan_{\kappaup^{-1} (y)} \kappaup \, (\eta) \right), 
\end{equation}
where $\coTan_{y} V \ni \eta \mapsto \coTan_{\kappaup^{-1} (y)} \kappaup \, (\eta) := \eta \circ \rd_{x} \kappaup \in \coTan_{x} U$  
and 
$\rd_{x} \kappaup : \tangent_{x} U \to \tangent_{\kappaup (x)} V$ 
is the tangent mapping of $\kappaup$. 
Hence one has the induced mapping  
\begin{equation}
	(\coTan \kappaup)^{*} : C^{\infty} (\coTan U) \to C^{\infty} (\coTan V), ~a \mapsto 
    \big( (\coTan \kappaup)^{*} a \big) \yeta := a \big( \coTan \kappaup \, \yeta \big).  
\end{equation}
Then, the chain rule  
\begin{equation}
	\rD_{x}^{\alpha} 
    \mapsto 
	\Big( \parDeri{x^{1}}{ \kappaup^{j_{1}} } \frac{\partial}{\ri \partial y^{j_{1}} } \Big)^{\alpha_{1}} 
	\ldots 
	\Big( \parDeri{x^{d}}{ \kappaup^{j_{d}} } \frac{\partial}{\ri \partial y^{j_{d}} } \Big)^{\alpha_{d}}  
\end{equation}
entails that the transformed total symbol $\fl \mapsto (\coTan \kappaup)^{*} \fl$ ends up with messy tangled terms involving higher derivatives of $\kappaup$, whilst the principal symbol transforms invariantly: 
\begin{equation}
	\big( (\coTan \kappaup)^{*} l \big) \yeta 
	= 
    \sum_{|\alpha| = m} \big( (\kappaup^{-1})^{*} f_{\alpha} \big)  (y) \,
	\Bigg( \Big( (\kappaup^{-1})^{*} \parDeri{x^{1}}{ \kappaup^{j_{1}} } \Big) (y) \, \eta_{j_{1}} \Bigg)^{\alpha_{1}} 
	\ldots 
	\Bigg( \Big( (\kappaup^{-1})^{*} \parDeri{x^{d}}{ \kappaup^{j_{d}} } \Big) (y) \, \eta_{j_{d}} \Bigg)^{\alpha_{d}} . 
\end{equation}
In other words, the following diagrams commute for any $L_{U} \in  \PDO^{m} (U), L_{V} \in \PDO^{m} (V)$: 
%
%
%
\begin{center}
	\begin{tikzpicture}
        \node (a) at (0, 0) {$C^{\infty} (V)$};
        \node (b) at (4, 0) {$C^{\infty} (V)$};
        \node (c) at (0, 2) {$C^{\infty} (U)$};
        \node (d) at (4, 2) {$C^{\infty} (U)$};
        \draw[->] (a) -- (b); 
        \draw[->] (c) -- (d); 
        \draw[->] (a) -- (c); 
        \draw[->] (b) -- (d); 
        \node[below] at (2, 0) {$L_{V}$}; 
        \node[above] at (2, 2) {$L_{U}$}; 
        \node[left] at (0, 1) {$\kappaup^{*}$}; 
        \node[right] at (4, 1) {$\kappaup^{*}$}; 
        \node (e) at (8, 0) {$\PDO^{m - [1]} (V)$}; 
        \node (f) at (8, 2) {$\PDO^{m - [1]} (U)$}; 
        \node (g) at (12, 0) {$S^{m - [1]} (\coTan V)$}; 
        \node (h) at (12, 2) {$S^{m - [1]} (\coTan U)$};
        \draw[->] (e) -- (g); 
        \draw[->] (f) -- (h); 
        \draw[->] (f) -- (e); 
        \draw[->] (h) -- (g);
        \node[left] at (8,1) {$\boldsymbol{\kappaup}$}; 
        \node[right] at (12, 1) {$(\coTan \kappaup)^{*}$}; 
        \node[above] at (10, 2) {$\symb{[L_{U}]}$}; 
        \node[below] at (10, 0) {$\symb{[L_{V}]}$}; 
    \end{tikzpicture}
    \captionof{figure}[Transformation of principal symbol of a differential operator]{Transformation 
        of a partial differential operator and its principal symbol under the diffeomorphism $\kappaup :=$ \eqref{eq: diffeo_U_V}.  
        In the right diagram, $\PDO^{m} (U) \ni L_{U} \mapsto L_{V} := \boldsymbol{\kappaup} (L_{U}) := (\kappaup^{*})^{-1} \circ L_{U} \circ \kappaup^{*} \in \PDO^{m} (V)$, and all other symbols are as defined in the ongoing discussion.
    }
    \label{fig: transformation_PDO_symbol}
\end{center}
%
%
%

In order to achieve an intrinsic characterisation of these notions, one observes  that~\cite{Hoermander_CPAM_1965} 
(see also~\cite{Weinstein_BullAMS_1976, *Weinstein_TransAMS_1978},~\cite[pp. 151-152]{Hoermander_Springer_2003}) 
\begin{equation}
	\forall \lambdaup \in \R_{+} : 
    \re^{- \ri \lambdaup x \cdot \xi} L (\re^{\ri \lambdaup x \cdot \xi}) 
    = 
    \sum_{|\alpha| \leq m} f_{\alpha} \xi^{\alpha} \lambdaup^{\alpha} 
\end{equation}
is a polynomial in $\lambdaup$ of degree $m$ for any $L \in \PDOU$. 
Conversely, any continuous linear operator satisfying the preceding property is a differential operator of order $m$ and can be expressed as~\eqref{eq: PDO_Euclidean_Fourier}. 
Evidently, hence 
\begin{equation} \label{eq: def_expression_symbol_PDO_Euclidean}
   	\totSymb{L} \xxi := \re^{- \ri x \cdot \xi} L (\re^{\ri x \cdot \xi}),  
   	\qquad 
   	\symb{L} \xxi := \lim_{\lambdaup \to \infty} \lambdaup^{-m} \re^{- \ri \lambdaup x \cdot \xi} L (\re^{\ri \lambdaup x \cdot \xi}). 
\end{equation}
%
%
%
%
%
%
%
%
%
%
\subsection{Pseudodifferential operators}
\label{sec: PsiDO_Euclidean} 
The essential notion of a pseudodifferential operator is to consider an operator of the form~\eqref{eq: PDO_Euclidean_Fourier} whose Schwartz kernel is also singular at the diagonal yet generalises a differential operator in the sense that it does not respect supports and that the total symbol~\eqref{eq: def_coordinate_expression_symbol_PDO_Euclidean} of a differential operator is replaced by a wider class of functions.  
To be precise, a \textbf{pseudodifferential operator $P$} on $U$ of \textbf{order} (at most) $m \in \R$ is a continuous linear map  
(see e.g.~\cite[p. 69]{Hoermander_Springer_2007}) 
\begin{equation} \label{eq: def_PsiDO_Euclidean}
	P : \CComInfinity{U} \to C^{\infty} (U),  ~ u \mapsto 
	(P u) (x) := \int_{U} \fP (x, y) \, u (y) \, \rd y  
\end{equation}
whose Schwartz kernel $\fP \in \cD' (U \times U)$ is of the form 
\begin{equation} \label{eq: def_naive_kernel_PsiDO_Euclidean}
	\fP (x, y) = \frac{1}{(2 \pi)^{d}} \int_{\Rd} \re^{\ri (x - y) \cdot \xi} \fp (x, y; \xi) \, \rd \xi, 
\end{equation}
for some element $\fp$ in the symbol class (see Definition~\ref{def: polyhomogeneous_symbol_Euclidean}) $S^{m} (U \times U \times \Rd) \supsetneq \bbP_{m} (U \times U \times \Rd)$. 
Note, the preceding expression is an \textit{oscillatory integral} and thus it must be understood in a distributional sense. 
We denote the set of all such operators by $\PsiDOU$. 
Recall, the set of \textbf{smoothing operators} $\PsiDO{-\infty}{U}$ is defined as the set of all linear continuous operators $\CComInfinity{U} \to C^{\infty} (U)$ whose Schwartz kernels are smooth. 
\newline 

Naturally, we are now tempted to have the notions of the total symbol and the principal symbol, for which, let us introduce a subset of $\PsiDOU$ which respects support. 

%
%
%
\begin{definition} \label{def: psPsiDO_Euclidean}
	Let $U \subset \Rd$ be an open set and $\PsiDOU$ the set of all pseudodifferential operators on $U$ of order (at most) $m \in \R$. 
    An element $P \in \PsiDOU$ is called \textbf{properly supported} if both projections $\supp \fP \to U, U$ from the support of its Schwartz kernel $\fP$ in $U \times U$ to $U$ are proper maps, i.e., for every compact set $K \subset U$ there is a compact set $\tilde{K} \subset U$ such that~\cite[p. 103]{Hoermander_ActaMath_1971}  
	(see also~\cite[Def. 18.1.21]{Hoermander_Springer_2007}) 
	\begin{equation}
		\supp u \subset K \Rightarrow \supp (Pu) \subset \tilde{K}, 
		\quad 
		u \upharpoonright \tilde{K} = 0 \Rightarrow (Pu) \upharpoonright K = 0. 
	\end{equation}
	Analogously, $\fp \in S^{m} (U \times U \times \Rd)$ is called properly supported if the projectors $\clo \big( \Pr (\supp \fp) \big)$ $\to U, U$ are proper, where $\Pr \fp$ projects $\fp$ onto $U \times U$ and $\clo$ denotes the closure. 
\end{definition}
%
%
%

\textit{A properly supported pseudodifferential operator can be defined with restricting supports}, i.e., $P : \CComInfinity{U} \to \CComInfinity{U}, P : \cE' (U) \to \cE' (U)$ and moreover $P : C^{\infty} (U) \to C^{\infty} (U), P : \cD' (U) \to \cD' (U)$. 
These operators are useful because 

%
%
%
\begin{proposition} \label{prop: PsiDO_psPsiDO_smooth}
	Every pseudodifferential operator of order $m$ can be written as a combination of an $m$-order properly supported pseudodifferential operator and a smoothing operator~\cite[p. 103]{Hoermander_ActaMath_1971}.  
\end{proposition}
%
%
%

The concept of the total symbol of a pseudodifferential operator is formalised as follows. 

%
%
%
\begin{theorem} \label{thm: characterisation_psPsiDO_Euclidean}
	As in the terminologies of Definition~\ref{def: psPsiDO_Euclidean}, any properly supported $P \in \PsiDO{m}{U}$ can be uniquely expressed as~\cite[Thm. 2.1.1]{Hoermander_ActaMath_1971} 
	\begin{equation} \label{eq: form_psPsiDO_Euclidean}
		(Pu) (x) = \int_{\Rd} \re^{\ri (x - y) \cdot \xi} \totSymb{P} \xxi \, u (y) \, \rd y \, \measurexi, 
	\end{equation}
	where the total symbol $\totSymb{P} \xxi := \re^{- \ri x \cdot \xi} P (\re^{\ri \bullet \cdot \xi})  (x)$ is given by means of the following asymptotic summation (see Definition~\ref{def: asymptotic_summation}) formula 
	\begin{equation} \label{eq: def_coordinate_expression_total_symbol_psPsiDO_Euclidean}
		\totSymb{P} \xxi \sim \sum_{\alpha \in \NO^{d}} \frac{1}{\alpha !} (\partial_{\xi}^{\alpha} \rD_{y}^{\alpha} \fp) (x, x; \xi),  
	\end{equation}
    where $\fp$ as in~\eqref{eq: def_naive_kernel_PsiDO_Euclidean} is an element of $ S^{m} (U \times U \times \Rd)$. 
	\newline 

	Conversely, given any proper $\fp \in S^{m} (U \times U \times \Rd)$, every operator of the form~\eqref{eq: form_psPsiDO_Euclidean} with total symbol~\eqref{eq: def_coordinate_expression_total_symbol_psPsiDO_Euclidean} is in the set of properly supported $\PsiDOU$. 
\end{theorem}
%
%
%

Since a generic pseudodifferential operator $P$ differs from its properly supported part $\tilde{P}$ only by a smoothing operator, one sets 
\begin{equation}
    \totSymb{} : \PsiDO{m - [\infty]}{U} \to S^{m - [\infty]} (\coTan U), ~ [P] \mapsto 
    \totSymb{[P]} := \totSymb{\tilde{P}}.  
\end{equation}
The preceding expression of $\totSymb{[P]}$ then exhibits that~\eqref{eq: def_coordinate_expression_total_symbol_PDO_Euclidean} holds for pseudodifferential operators in the asymptotic sense.  
Furthermore, this expression paves the way to define the principal symbol $\symb{\tilde{P}}$ of $\tilde{P}$, which is precisely the \textit{leading order homogeneous term} in~\eqref{eq: def_coordinate_expression_total_symbol_psPsiDO_Euclidean}. 
Clearly, it does not see smooth terms, and so $\symb{P}$ is identified with it. 

%
%
%
\begin{definition} \label{def: symbol_PsiDO_Euclidean} 
    As in the terminologies of Definition~\ref{def: psPsiDO_Euclidean}, let $S^{m} (\coTan U)$ be the space of symbols on the cotangent bundle $\coTan U$ over $U$.  
    Then, the \textbf{principal symbol} $\symb{P}$ of $P \in \PsiDOU$ is defined by the isomorphism~\cite[p. 110]{Hoermander_ActaMath_1971}  
	\begin{eqnarray} \label{eq: def_symbol_PsiDO_Euclidean}
        \symb{} & : & \PsiDO{m - [1]}{U} \to S^{m - [1]} (\coTan U), ~ [P] \mapsto 
        \nonumber \\ 
        && 
        \symb{[P]} := \textrm{homogeneous term of degree $m$ in~\eqref{eq: def_coordinate_expression_total_symbol_psPsiDO_Euclidean}} \mod S^{m-1} (\coTan U). \quad  
	\end{eqnarray} 
\end{definition}
%
%
%

In terms of transformation properties, $\totSymb{P}$ and $\symb{P}$ resemble those of a differential operator. 
More precisely, under the action of $\kappaup :=$~\eqref{eq: diffeo_U_V}, the total symbol of any properly supported $P \in \PsiDOU$ transforms like~\cite[$(2.1.14)$]{Hoermander_ActaMath_1971}
(see also~\cite[Thm. 18.1.17]{Hoermander_Springer_2007}) 
\begin{equation} \label{eq: total_symbol_PsiDO_Euclidean_coordinate_change}
	\big( (\coTan \kappaup)^{*} \totSymb{P} \big) \yeta 
	\sim 
	\sum_{|\alpha| \leq m} \frac{1}{\alpha !} \partial_{\xi}^{\alpha} \fp \xxi \,   \rD_{\tilde{x}}^{\alpha} \re^{\ri \rhoup (x, \tilde{x}) \cdot (\coTan_{x} \kappaup)^{-1} (\xi)} \big|_{\tilde{x} = x}, 
\end{equation}
where $\xxi := \big(\kappaup^{-1} (y), \coTan_{x} \kappaup \, \eta \big)$ and $\rhoup (x, \tilde{x}) := \kappaup (\tilde{x}) - \kappaup (x) - (\tilde{x} - x) \rd_{x} \kappaup$. 

%
%
%
\begin{remark} \label{rem: symbol_PsiDO_Euclidean_coordinate_change}
	The leading-homogeneous term $p$ in~\eqref{eq: total_symbol_PsiDO_Euclidean_coordinate_change} satisfies 
	(see e.g.~\cite[p. 83]{Hoermander_Springer_2007})
	\begin{equation} \label{eq: symbol_PsiDO_Euclidean_coordinate_change}
        \big( (\coTan \kappaup)^{*} p \big) \yeta = p (x, \coTan_{x} \kappaup \, \eta). 
	\end{equation}
    It is then evident that the principal symbol is an \textit{invariantly} defined \textit{homogeneous function} on the punctured cotangent bundle $\dotCoTan U$ in contrast to the total symbol.  
    In other words, Figure~\ref{fig: transformation_PDO_symbol} holds for pseudodifferential operators as well. 
\end{remark}
%
%
%

So far, we have specified a (properly supported) pseudodifferential operator by its Schwartz kernel~\eqref{eq: form_psPsiDO_Euclidean} which has been given as an oscillatory integral whose integrand is $\re^{\ri (x - y) \cdot \xi} \totSymb{P} \xxi$. 
But, neither the exponent nor the total symbol are invariant under a diffeomorphism, and so it is unforeseeable whether the Schwartz kernel has an invariant meaning. 
So, we ought to look for an intrinsic characterisation of such class of distributions. 
It is well-recognised that the singularity structure is a characteristic feature of a distribution. 
But, our purpose cannot be served by the 
singular support\footnote{Recall, 
    the \textbf{singular support} $\singsupp u$ of a distribution $u \in \cD' (U)$ is the set of all points in $U$ having no open neighbourhood to which the restriction of $u$ is a smooth function 
    (see e.g.~\cite[Def. 2.2.3]{Hoermander_Springer_2003}.} 
as it is \textit{not} diffeomorphism invariant. 
Thus, one introduces the \textbf{wavefront set}~\cite{Hoermander_Nice_1970},~\cite[$(2.5.2)$]{Hoermander_ActaMath_1971} 
\begin{equation} \label{eq: def_WF_Euclidean}
    \WF u := \bigcap_{\substack{P \in \PsiDO{0}{U} \\ Pu \in C^{\infty} (U)}} \Char P, 
    \quad 
    \Char P := \{ \xxiNot \in \dotCoTan U \,|\, \symb{P} \xxiNot = 0 \}. 
\end{equation}
of a distribution $u \in \cD' (U)$, where the intersection runs over all properly supported $P$ and $\Char P$ is called the \textbf{characteristic variety} of $P$. 
Evidently, a wavefront set transforms \textit{covariantly} under a diffeomorphism. 
Equivalently, $\dotCoTan U \ni \xxiNot \notin \WF{u}$ if and only if \textit{there exists a compactly supported smooth function $f$ on $U$ non-vanishing at $x_{0} \in U$ such that the Fourier transform of $\big( \Fourier (fu) \big) (\xi)$ is rapidly decreasing in a conic neighbourhood of} $\xi^{0}$~\cite[Prop. 2.5.5]{Hoermander_ActaMath_1971}. 
Note, wavefront set is a generalisation of the singular support because the latter is the base projection of the former~\cite[Thm. 2.5.3]{Hoermander_ActaMath_1971}. 
Speaking differently, $\singsupp u$ only supplies the \textit{location} of the singularities whilst $\WF u$ provides the \textit{cotangent vectors causing them, on top of their locations}.    
If $\fA \in \cD' (U \times V)$ is a bidistribution where $V$ is an open subset of any Euclidean space, then its \textbf{twisted wavefront set} is given by~\cite[Thm. 2.5.14]{Hoermander_ActaMath_1971} 
\begin{equation}
    \WF' \fA := \{ (x, \xi; y, - \eta) \in \dotCoTan U \times \dotCoTan V | \xxiyeta \in \WF \fA \}
\end{equation}
Details of wavefront set is available, for instance, 
in~\cite{Hoermander_Springer_2003, Strohmaier_Springer_2009, Brouder_JPA_2014}. 
\newline 

It follows that the wavefront set of the Schwartz kernel of a pseudodifferential operator~\cite[p. 124]{Hoermander_ActaMath_1971}  
(see also~\cite[Thm. 18.1.16, 18.1.26]{Hoermander_Springer_2007}):  
\begin{equation} \label{eq: WF_kernel_PsiDO}
    \WF \fP \subseteq \{ (x, \xi; x, - \xi) \in \dotCoTan U \times \dotCoTan U \} = (\varDelta U)^{\perp *}
\end{equation}
is contained in the conormal bundle (see Examples~\ref{exm: conormal_bundle_Lagrangian_submf} and~\ref{exm: diagonal_embedding_Lagrangian_submf})  $(\varDelta U)^{\perp *}$ of the diagonal embedding 
\begin{equation} \label{eq: def_diagonal_embedding_Euclidean}
	\varDelta : U \to U \times U, ~ x \mapsto \varDelta (x) := (x, x). 
\end{equation}
By projecting $\coTan U \times \coTan U \to \coTan U$ on the first factor, $(\varDelta U)^{\perp *}$ can be \textit{identified} with $\coTan U$ and 
\begin{equation} \label{eq: def_WF_PsiDO}
    \WF P = \ES P := \{ \xxi \in \dotCoTan U \,|\, (x, \xi; x, \xi) \in \WF' \fP \}  
\end{equation}
is called the \textbf{wavefront set}, also known as the \textbf{essential support} of a pseudodifferential operator $P$. 
This is the \textit{smallest conic set} such that $P$ is of \textit{order} $-\infty$ in $\dotCoTan U \setminus \WF P$. 
In particular, if $P$ is properly supported then $\dotCoTan U \setminus \WF P$ is the \textit{largest open conic set} where $\totSymb{P}$~\eqref{eq: def_coordinate_expression_total_symbol_psPsiDO_Euclidean} is rapidly decreasing~\cite[Prop. 2.5.8]{Hoermander_ActaMath_1971}. 
\newline 

Given a $P \in \PsiDOU$, there exists a unique formally dual operator $P^{*} \in \PsiDOU$ such that 
\begin{equation}
	P^{*} : \cE' (U) \to \cD' (U), ~ u \mapsto 
	(P^{*} u) (\phi) := u (P \phi) 
\end{equation}
for any $\phi \in \CComInfinity{U}$. 
Pseudodifferential operators are \textit{microlocal}, meaning that 
\begin{equation} \label{eq: PsiDO_Euclidean_microlocal}
	\WF (Pu) \subseteq \WF P \cap \WF u \Rightarrow \singsupp (Pu) \subseteq \singsupp u  
\end{equation}
for any $u \in \cE' (U)$. 
Elements of $\PsiDO{-\infty}{U}$ are called smoothing operators due to the fact that they map $\cE' (U) \to C^{\infty} (U)$. 
We remark that, by the Peetre theorem~\cite{Peetre_MathScand_1960} - \textit{every linear and local operator is a partial differential operator $L$}:  
\begin{equation}
	\forall u \in \cD' (U) : \supp (Lu) \subseteq \supp u.  
\end{equation}
Thus, pseudodifferential operators are pseudolocal (actually microlocal) generalisation of differential operators in the sense that they respect the singular support (actually wavefront set) instead of the support. 
\newline 

With the precise formulation~\eqref{eq: WF_kernel_PsiDO} of singularities of $\fP$, one furthermore notices that $(\varDelta U)^{\perp *}$ is generated by the non-degenerate phase function (see Definition~\ref{def: clean_phase_function}) $\varphi := (x - y) \cdot \xi$, as explained in Example~\ref{exm: diagonal_embedding_Lagrangian_submf}. 
In fact, it is the conormal bundle rather the oscillatory integral~\eqref{eq: def_naive_kernel_PsiDO_Euclidean} that characterises $\fP$ invariantly~\cite{Hoermander_ActaMath_1971} 
(for details, see e.g.~\cite{Duistermaat_CPAM_1974},~\cite[Sec. 2.3]{Duistermaat_Birkhaeuser_2011},~\cite[Sec. VIII.1]{Treves_Plenum_1980}). 

%
%
%
\begin{definition} \label{def: conormal_distribution_Euclidean}
    As in the terminologies of Definition~\ref{def: symbol_PsiDO_Euclidean}, let $(\varDelta U)^{\perp *}$ be the conormal bundle of the diagonal embedding~\eqref{eq: def_diagonal_embedding_Euclidean}. 
    Then the space $I^{m} \big( U \times U, (\varDelta U)^{\perp *} \big)$ of \textbf{conormal distributions} of \textbf{order} (at most) $m$ is defined as the set of all bidistributions $\fP \in \cD' (U \times U)$ on $U$ of the form~\cite[$(2.4.1)$]{Hoermander_ActaMath_1971}  
    \begin{equation}
        \fP (\phi \otimes u) = \int_{\Rd} \int_{U} \int_{U} \re^{\ri (x - y) \cdot \xi} \fp (x, y; \xi) \, \phi (x) \, u (y) \, \rd x \, \rd y \frac{\rd \xi}{(2 \pi)^{d}},  
    \end{equation}
    where $u, \phi \in \CComInfinity{U}$ and $\fp \in S^{m} (U \times U \times \Rd)$. 
\end{definition}
%
%
%

Then, $P \in \PsiDOU$ is invariantly defined as the integral operator~\eqref{eq: def_PsiDO_Euclidean} whose Schwartz kernel $\fP \in I^{m} \big( U \times U, (\varDelta U)^{\perp *} \big)$.  

%
%
%
\begin{example}
	Let $\Delta$ be the Laplacian on an open set $U \subset \Rd$. 
	Then $\sqrt{\Delta} \in \PsiDO{1}{U}$. 
\end{example}
%
%
%

%
%
%
%
%
%
%
\subsection{Fourier integral operators}
\label{sec: FIO_Euclidean} 
As observed in the last two sections, the partial differential operators are local whereas the pseudodifferential operators are microlocal. 
Fourier integral operators are a vast generalisation of the latter in the sense that it maps a function on a set to a distribution on another set by allowing a general phase function. 
Let $U \subset \R^{d_{\ms U}}$ and $V \subset \R^{d_{\ms V}}$ be open sets.  
Then, a \textbf{Fourier integral operator} of \textbf{order} (at most) $m \in \R$ is the continuous linear map~\cite[Sec. 1.4]{Hoermander_ActaMath_1971}  
(see also~\cite[Prop. 25.1.5']{Hoermander_Springer_2009} 
and, e.g.~\cite[Sec. VI.2]{Treves_Plenum_1980})
\begin{equation} \label{eq: def_FIO_Euclidean}
	A : \CComInfinity{V} \to \cD' (U), ~v \mapsto 
	(A v) (x) := \int_{V} \fA (x, y) \, v (y) \, \rd y  
\end{equation}
whose Schwartz kernel $\fA$ is given by the oscillatory integral of the form 
\begin{equation} \label{eq: def_naive_kernel_FIO_Euclidean}
    \fA (x, y) = (2 \pi)^{- (d_{\ms U} + d_{\ms V} + 2n - 2e) / 4} \int_{\Rn} \re^{\ri \varphi \xytheta} \fa \xytheta \, \rd \theta 
\end{equation} 
 where $\varphi$ is a clean phase function (see Definition~\ref{def: clean_phase_function}) with excess $e$ on $U \times V \times \dotRn$, $\fa \in S^{m + (d_{\ms U} + d_{\ms V} - 2n - 2e)/4} (U \times V \times \Rn)$ and $\rd \theta$ is the Lebesgue measure on $\Rn$. 
\newline 

Being an oscillatory integral,~\eqref{eq: def_naive_kernel_FIO_Euclidean} must be understood as a \textit{formal expression} that does not pointwise make sense rather it defines a bidistribution $\cD' (U \times V)$ in the sense that~\cite[Sec. 1.2]{Hoermander_ActaMath_1971}  
(see also, e.g.~\cite[Sec. 2.2-2.3]{Duistermaat_Birkhaeuser_2011},~\cite[Sec. 1.1-1.3]{Shubin_Springer_2001} for details) 
\begin{equation} \label{eq: def_naive_weak_FIO_Euclidean}
    (Av) (\phi) 
    := \fA (\phi \otimes v) 
    = (2 \pi)^{- (d_{\ms U} + d_{\ms V} + 2n - 2e) / 4} \int_{\Rn} \int_{V} \int_{U} \re^{\ri \varphi \xytheta} \fa \xytheta \, \phi (x) \, v (y) \, \rd x \, \rd y \, \rd \theta 
\end{equation} 
for any $\phi \in \CComInfinity{U}$. 
The precise meaning of this expression is as follows~\cite[Thm. 1.4.1]{Hoermander_ActaMath_1971}:  
\begin{enumerate}
    \item[(i)] 
	If $\varphi$ has no critical point as a function of $\xytheta$ then the oscillatory integral~\eqref{eq: def_naive_weak_FIO_Euclidean} exists and it is a continuous bilinear form for the $C_{\rc}^{k}$-topologies on $v$ and $\phi$ whenever   
    \begin{equation} \label{eq: FIO_1_4_3}
		m - k < - n.  
	\end{equation}
    When~\eqref{eq: FIO_1_4_3} is valid then one obtains a continuous linear map
    \begin{equation}
    	A : \CComk{V} \to \cD'^{k} (U)
    \end{equation}
    whose Schwartz kernel $\fA \in \cD'^{k} (U \times V)$ is given by the oscillatory integral 
	\begin{equation} \label{eq: FIO_I_1_4_4}
		\forall u \in \CComInfinity{U \times V}: 
        \fA (u)
		= 
		\int_{\Rn} \int_{V} \int_{U} \re^{\ri \varphi \xytheta} a \xytheta \, u (x, y) \, \rd x \, \rd y \, \rd \theta. 
	\end{equation}
	\item[(ii)] 
    If, for each fixed $x$, $\varphi$ has no critical point $(y_{0}; \theta^{0})$ then~\eqref{eq: def_naive_kernel_FIO_Euclidean} is defined as an oscillatory integral. 
    One obtains a continuous map
    \begin{equation} \label{eq: FIO_Euclidean_CCom_C}
    	A : C_{\rc}^{k} (V) \to C (U)
    \end{equation}
    when~\eqref{eq: FIO_1_4_3} is valid. 
    By differentiation under the integral sign, it follows that $A$ is also a continuous map from $C_{\rc}^{k} (V)$ to $C^{j} (U)$ provided that 
    \begin{equation} \label{eq: FIO_I_1_4_5}
		m + n + j < k. 
	\end{equation}
    \item[(iii)] 
    If, for each fixed $y$, $\varphi$ has no critical point $\xthetaNot$ then the adjoint of $A$ has the properties listed in (ii), so 
    \begin{equation}
    	A : \cE'^{j} (V) \to \cD'^{k} (U)
    \end{equation}
    is a continuous map when~\eqref{eq: FIO_I_1_4_5} is fulfilled. 
	In particular, $A$ defines a continuous map from $\cE' (V)$ to $\cD' (U)$. 
	\item[(iv)] 
	Let $\sR$ be the open set of all $(x, y) \in U \times V$ such that $\varphi \xytheta$ has no critical point as a function of $\theta$. 
	Then, the oscillatory integral 
	\begin{equation}
		\fA (x, y) = (2 \pi)^{- (d_{\ms U} + d_{\ms V} + 2n - 2e) / 4} \int_{\Rn} \re^{\ri \varphi \xytheta} a \xytheta \, \rd \theta 
	\end{equation}
	defines a function in $C^{\infty} (\sR)$ which is equal to the distribution~\eqref{eq: FIO_I_1_4_4} in $\sR$. 
    If $\sR = U \times V$, it follows that $A$ is an integral operator with the smooth Schwartz kernel $\fA \in C^{\infty} (U \times V)$, so 
    \begin{equation} \label{eq: smoothing_op_V_U}
        A : \cE' (V) \to C^{\infty} (U)
    \end{equation}
    is a continuous map, called the \textbf{smoothing operator} from $V$ to $U$.  
\end{enumerate}

As described for the kernel of a pseudodifferential operator, we would like to characterise the oscillatory integral $\fA$ in an invariant fashion starting with the observation that its wavefront set~\cite[Thm. 3.2.6]{Hoermander_ActaMath_1971}  
\begin{equation}
	\WF \fA \subseteq C', 
    \quad 
    C' := \big\{ (x_{0}, \rd_{x_{0}} \varphi; y_{0}, - \rd_{y_{0}} \varphi) \in \dotCoTan U \times \dotCoTan V | (\grad_{\theta} \varphi) \xythetaNot = 0 \big\}
\end{equation}
is contained in the homogeneous twisted canonical relation (see Definition~\ref{def: canonical_relation}) $C'$ from $\dotCoTan V$ to $\dotCoTan U$ which is closed in $\dotCoTan (U \times V)$. 
The Schwartz kernel of a Fourier integral operator is then intrinsically formalised as,  

%
%
%
\begin{definition} \label{def: Lagrangian_distribution_Euclidean}
	Let $U \subset \R^{d_{\ms U}}, V \subset \R^{d_{\ms V}}$ be open sets and $C \subset \dotCoTan U \times \dotCoTan V$ a homogeneous canonical relation which is closed in $\dotCoTan (U \times V)$. 
    Then the space $I^{m} (U \times V, C')$ of \textbf{Lagrangian distributions} on $U \times V$ of \textbf{order} (at most) $m \in \R$ is defined as the set of all bidistributions $\fA \in \cD' (U \times V)$ on $U \times V$ of the form~\eqref{eq: def_naive_weak_FIO_Euclidean}~\cite[Def. 3.2.2]{Hoermander_ActaMath_1971}. 
\end{definition}
%
%
%

Therefore, \textit{a Fourier integral operator $A$ associated with a closed homogeneous canonical relation $C \subset \dotCoTan U \times \dotCoTan V$ is defined by the integral operator~\eqref{eq: def_FIO_Euclidean} whose Schwartz kernel $\fA$ is an element of $I^{m} (U \times V, C')$}. 
We denote the set of all such operators by $\FIO{m}{V \to U, C'}$. 
Note, given an oscillatory integral representation of $\fA$ of the form~\eqref{eq: def_naive_kernel_FIO_Euclidean}, the degree $\mu$ of homogeneity of $\fa$ is related to the order $m$ of $A$ by 
(see e.g.~\cite[$(2.4.22)$]{Duistermaat_Birkhaeuser_2011})
\begin{equation} \label{eq: order_FIO_order_amplitude}
    m = \mu + \frac{n}{2} - \frac{d_{\ms U} + d_{\ms V}}{4}.   
\end{equation}

The apparently counter-intuitive $(2 \pi)^{- (d_{\ms U} + d_{\ms V} + 2n - 2e) / 4}$ factor in~\eqref{eq: def_naive_weak_FIO_Euclidean} has been chosen in order to match the $(2 \pi)^{-d}$ factor in~\eqref{eq: def_naive_kernel_PsiDO_Euclidean} when $d_{\ms U} = d_{\ms V} = d = n$ and $e = 0$. 
As expected, $I^{m} \big( U \times U, (\varDelta U)^{\perp *} \big)$ is a special case of $I^{m} (U \times V, C')$ when $C$ is given by the graph of the identity homogeneous symplectomorphism on $U$, globally generated by the non-degenerate phase function $\varphi := (x - y) \cdot \xi$ (as elucidated in Example~\ref{exm: diagonal_embedding_Lagrangian_submf}).  
Besides $C^{\infty} (U \times V) \subset I^{-\infty} (U \times V, C')$ 
(see e.g.~\cite[Prop. 3.2 $($p. 439$)$]{Treves_Plenum_1980}), 
i.e., the kernel of a smoothing operator~\eqref{eq: smoothing_op_V_U} can be thought as an element of $I^{- \infty} (U \times V, C')$.  
\newline 

In order to define the principal symbol of a Fourier integral operator, one cannot directly look for an isomorphism between $\FIO{m}{V \to U, C'}$ and the respective symbol class, as done for pseudodifferential operators (Definition~\ref{def: symbol_PsiDO_Euclidean}) because there is no analog of~\eqref{eq: def_coordinate_expression_total_symbol_psPsiDO_Euclidean} for the total symbol of a Fourier integral operator due to the fact that a canonical relation is more intricate than the conormal bundle $(\varDelta U)^{\perp *}$. 
Since a pseudodifferential operator is characterised by its Schwartz kernel and $(\varDelta U)^{\perp *} \cong \coTan U$, its principal symbol isomorphism map~\eqref{eq: def_symbol_PsiDO_Euclidean} is equivalent to~\cite[Thm. 2.4.2]{Hoermander_ActaMath_1971} 
\begin{equation}
    I^{m - [1]} \big( U \times U, (\varDelta U)^{\perp *} \big) \cong S \big( (\varDelta U)^{\perp *} \big). 
\end{equation}
This paves the right direction to address the principal symbol of a Fourier integral operator. 
In other words, one pursues an isomorphism between $I^{m} (U \times V, C')$  and respective symbol space on $C$. 
It turns out that one actually requires to consider half-density (see Section~\ref{sec: distribution_density_mf} and Remark~\ref{rem: subprincipal_symbol_PsiDO_mf_coordinate_change} for further motivation)-valued symbols on $C$ of order $m + (d_{\ms U} + d_{\ms V}) / 4$ in order to define an invariant principal symbol~\cite[Thm. 3.2.1]{Hoermander_ActaMath_1971}. 

%
%
%
\begin{definition} \label{def: symbol_FIO_Euclidean}
    As in the terminologies of Definition~\ref{def: Lagrangian_distribution_Euclidean}, suppose that $\halfDenC, \M \to C$ are the half-density bundle (see Appendix~\ref{sec: volume_canonical_relation}) and the Keller-Maslov bundle (see Definition~\ref{def: Keller_Maslov_bundle}) over $C$, respectively, and that $S^{m + (d_{\ms U} + d_{\ms V}) / 4 - [1]} (C; \M \otimes \halfDenC)$ is the space (see~\eqref{eq: symbol_canonical_relation_Maslov_density_function}) of $\Maslov \otimes \halfDenC$-valued symbols on $C$. 
    Then, the \textbf{principal symbol} of a Lagrangian distribution is defined by the isomorphism~\cite[$(7.6)$]{Duistermaat_InventMath_1975} 
    (see also~\cite[Prop. 25.1.5']{Hoermander_Springer_2009})   
	\begin{eqnarray} \label{eq: def_symbol_FIO_Euclidean}
        \symb{} \! &:& \! I^{m - [1]} (U \times V, C') \to S^{m + \frac{d_{\ms U} + d_{\ms V}}{4} - [1]} (C; \M \otimes \halfDenC), 
        \nonumber \\ 
        \left[ \fA \right] 
        & \mapsto & 
        \symb{[\fA]} \xxiyeta := \int_{\kC_{\xi, \eta}} \rd \theta'' a (x, y; \theta', \theta'') \m \otimes |\dVol_{\ms C}|^{\frac{1}{2}}  \! \!
        \mod S^{m + \frac{d_{\ms U} + d_{\ms V}}{4} - 1} (\cdot), \qquad  \;
	\end{eqnarray}
    where $I^{m-[1]} := I^{m} / I^{m-1}$, $a$ is the top-degree homogeneous term of $\fa$ in~\eqref{eq: def_naive_kernel_FIO_Euclidean}, $\m$ and $\sqrt{|\dVol_{\ms C}|}$ are sections of $\M$ and $\halfDenC$, respectively. 
    For each $\xxiyeta \in C$, $\varphi$ is a representative of the stable equivalence class (see Remark~\ref{rem: stable_equivalence}) of clean phase functions with excess $e$ (see Definition~\ref{def: clean_phase_function}) for $C$ and 
    \begin{equation} \label{eq: def_fibre_Lagrangian_fibration_symbol_FIO_Euclidean}
        \kC_{\xi, \eta} := \{ \xytheta \in (\grad_{\theta} \varphi)^{-1} (0) \,|\, \rd_{x} \varphi := \xi, \rd_{y} \varphi =: \eta \}
    \end{equation}
    is the $e$-dimensional fibre over the corresponding Lagrangian fibration (see~\eqref{eq: def_canonical_relation_fibration}). 
    Here, $\rd \theta''$ is the Lebesgue measure on $\dot{\R}^{e}$ and the variable $\theta''$ is defined by the 
    splitting\footnote{Such 
        splitting is always possible due to the Thom splitting (also known as the parametrised Morse) lemma  
        (see e.g.~\cite[App. C. 6]{Hoermander_Springer_2007},~\cite[p. 52]{Bates_AMS_1997}).
    } 
    $\dotRn \in \theta = (\theta', \theta'') \in \dot{\R}^{n - e} \times \dot{\R}^{e}$ such that the projection $\kC_{\xi, \eta} \ni (x, y; \theta', \theta'') \mapsto \theta'' \in \dot{\R}^{e}$ has a bijective differential so that, for a fixed $\theta'' = \cst$, $\varphi (x, y; \theta', \cst)$ is non-degenerate.   
\end{definition}
%
%
%

Explicit expressions of $\m$ (see Appendix~\ref{sec: Keller_Maslov_bundle}) and $\sqrt{|\dVol_{\ms C}|}$ (see Appendix~\ref{sec: volume_canonical_relation}) in terms of $\varphi$ entail~\cite[$(25.1.4)'$, p. 15]{Hoermander_Springer_2009} 
\begin{equation}
    \symb{\fA} \xxiyeta 
    =  
    \sqrt{|\rd \xi| |\rd \eta|}  \int_{\kC_{\xi, \eta}} a (x, y; \theta', \theta'') \, \dfrac{\re^{\nicefrac{\mathrm{i} \pi}{4} \sgn (\Hess_{x, y; \theta'} \varphi)}}{\sqrt{|\det (\Hess_{x, y; \theta'} \varphi)|}} \rd \theta'', 
\end{equation}
where $\Hess_{x, y; \theta'} \varphi$ is the Hessian matrix (see~\eqref{eq: def_Hessian}) of $\varphi (x, y; \theta', \cst)$ with respect to $(x, y; \theta') \in \kC_{\xi, \eta}$: 
\begin{eqnarray}
    \Hess_{x, y; \theta'} \varphi 
    = \left( 
   	\begin{array}{ccc}
        \dfrac{\partial^{2} \varphi}{ \partial x^{i} \partial x^{i^{\backprime}} } 
        & \dfrac{\partial^{2} \varphi}{\partial x^{i} \partial y^{j}} 
        & \dfrac{\partial^{2} \varphi}{\partial x^{i} \partial \theta'_{k}} 
   		\\ 
   		\dfrac{\partial^{2} \varphi}{\partial y^{j} \partial x^{i}} 
        & \dfrac{\partial^{2} \varphi}{ \partial y^{j} \partial y^{j^{\backprime}} } 
        & \dfrac{\partial^{2} \varphi}{\partial y^{j} \partial \theta'_{k}} 
   		\\ 
        \dfrac{\partial^{2} \varphi}{\partial \theta'_{k} \partial x^{i}} 
        & \dfrac{\partial^{2} \varphi}{\partial \theta'_{k} \partial y^{j}} 
        & \dfrac{\partial^{2} \varphi}{ \partial \theta'_{k} \partial \theta'_{k^{\backprime}} } 
   	\end{array}
    \right) 
\end{eqnarray} 
for $i, i^{\backprime} = 1, \ldots, d_{\ms U}; j, j^{\backprime} = 1, \ldots, d_{\ms V}; k, k^{\backprime} = 1, \ldots, n - e$. 
Here, $\rd \xi \, \rd \eta$ is the Lebesgue measure on $C$ at $\xxiyeta$.  
%
%
%

If $C$ is parametrised by a non-degenerate phase function $\varphiup$, i.e., $e = 0$, then the $\theta$-variable is no longer required to be split up and hence the expression of $\symb{\fA}$ simplifies to~\cite[Thm. 3.2.5]{Hoermander_ActaMath_1971} 
(see also~\cite[Prop. 25.1.5]{Hoermander_Springer_2009} and~\eqref{eq: volume_canonical_relation_nondegenerate_phase_function})
\begin{eqnarray} \label{eq: def_symbol_FIO_nondegenerate_Euclidean} 
	\symb{[\fA]} \xxiyeta 
    & := & 
    a \xytheta \, \dfrac{\re^{\nicefrac{\mathrm{i} \pi}{4} \sgn (\Hess  \varphiup)}}{\sqrt{|\det (\Hess \varphiup)|}} \sqrt{|\rd \xi| |\rd \eta|} \mod S^{m + \frac{d_{\ms U} + d_{\ms V}}{4} - 1} (\cdot), 
    \nonumber \\ 
    \Hess \varphiup 
    & = & 
    \left( 
    \begin{array}{ccc}
        \dfrac{\partial^{2} \varphiup}{ \partial x^{i} \partial x^{i^{\backprime}} } 
        & \dfrac{\partial^{2} \varphiup}{\partial x^{i} \partial y^{j}} 
        & \dfrac{\partial^{2} \varphiup}{\partial x^{i} \partial \theta_{k}} 
    	\\ 
    	\dfrac{\partial^{2} \varphiup}{\partial y^{j} \partial x^{i}} 
        & \dfrac{\partial^{2} \varphiup}{ \partial y^{j} \partial y^{j^{\backprime}} } 
        & \dfrac{\partial^{2} \varphiup}{\partial y^{j} \partial \theta_{k}} 
    	\\ 
        \dfrac{\partial^{2} \varphiup}{\partial \theta_{k} \partial x^{i}} 
        & \dfrac{\partial^{2} \varphiup}{\partial \theta_{k} \partial y^{j}} 
        & \dfrac{\partial^{2} \varphiup}{ \partial \theta_{k} \partial \theta_{k^{\backprime}} } 
    \end{array}
    \right).  
\end{eqnarray} 

%
%
%
\begin{remark} \label{rem: symbol_FIO_Euclidean_diffeo}
	The half-density $\symb{\fA}$ is invariant under any diffeomorphism apart from the Maslov factor $\m$ of absolute value $1$~\cite[$(25.1.14)$]{Hoermander_Springer_2009}. 
\end{remark}

%
%
%
\begin{example} \label{exm: pullback_Euclidean_FIO}
    Suppose that $U \subset \Rd$ and $V \subset \Rn$ are open sets and that $\kappaup \in \CInfinity{U, V}$. 
	Then, the definition of pullback 
	\begin{equation}
		\kappaup^{*} : \CComInfinity{V} \to \CInfinity{U}, ~ u \mapsto \kappaup^{*} u := u \circ \kappaup 
	\end{equation}
    entails 
    (see e.g.~\cite[p. 37]{Duistermaat_Birkhaeuser_2011},~\cite[Sec. 5.25]{Melrose_2007}) 
    \begin{eqnarray}
    	(\kappaup^{*} u) (x) 
    	= u \big( \kappaup (x) \big) 
        = \delta_{\kappaup (x)} (u) 
        = \frac{1}{(2 \pi)^{n}} \int_{\Rn} \int_{V} \re^{\ri ( \kappaup (x) - y ) \cdot \eta} u (y) \, \rd y \, \rd \eta.  
    \end{eqnarray}
    That is, $\kappaup^{*}$ is a Fourier integral operator whose Schwartz kernel is given by  
    \begin{equation}
    	\fK := (2 \pi)^{- (n-d + d+3n) / 4} \int_{\Rn} \re^{\ri \varphi} \rd \eta, 
        \quad 
        \varphi := \big( \kappaup (x) - y \big) \cdot \eta 
    \end{equation}
    where the phase function $\varphi$ is non-degenerate whose fibre-critical set (see~\eqref{eq: def_fibre_critical_set}) is  
    \begin{eqnarray}
        \sC := \big\{ (x_{0}, y_{0}; \eta^{0}) \in U \times V \times \dotRn \,|\, y_{0} = \kappaup (x_{0}) \big\}. 
    \end{eqnarray}
    Thereby, $\kappaup^{*}$ is of order (\eqref{eq: order_FIO_order_amplitude}) $\frac{n-d}{4}$ and is associated with the canonical relation (see~\eqref{eq: def_canonical_relation_fibration})  
    \begin{equation} \label{eq: canonical_relation_pullback_Euclidean}
        C := \{ (x, \xi; y, - \eta) \in \coTan U \times \dotCoTan V \,|\, y = \kappaup (x), \xi = \coTan_{x} \kappaup \, (\eta) \}, 
    \end{equation}
    where $\coTan_{x} \kappaup : \coTan_{\kappaup (x)} V \to \coTan_{x} U$ is the algebraic dual of the tangent mapping $\rd_{x} \kappaup : \tangent_{x} U \to \tangent_{\kappaup (x)} V$ of $\kappaup$. 
    Furthermore, its principal symbol is given by (cf.~\eqref{eq: def_symbol_FIO_nondegenerate_Euclidean}) 
    \begin{eqnarray} 
        \symb{\fK} \xxiyeta 
        := 
        \frac{\sqrt{|\rd \xi| |\rd \eta|}}{ (2 \pi)^{\frac{n-d}{4}} } \dfrac{\re^{\nicefrac{\mathrm{i} \pi}{4} \sgn (\Hess  \varphi)}}{ \sqrt{|\det (\Hess \varphi)|} }, 
        \quad 
        \Hess \varphi 
        =
        \left( 
        \begin{array}{ccc}
            \dfrac{\rd^{2} \kappaup^{j}}{ \rd x^{i} \rd x^{i^{\backprime}} } \eta_{j}
            & 0
            & \dfrac{\rd \kappaup^{j}}{\rd x^{i}} 
            \\ 
            0
            & 0
            & -1
            \\ 
            \dfrac{\rd \kappaup^{j}}{\rd x^{i}} 
            & -1  
            & 0 
        \end{array}
        \right), 
    \end{eqnarray} 
    where we have used $\kappaup (x) = \big( \kappaup^{1} (x), \ldots, \kappaup^{n} (x) \big)$. 
    \newpage  

    In addition, if $\kappaup$ is a diffeomorphism then we can consider its cotangent lift $\coTan \kappaup$, and thus $C' = \big\{ \big( \coTan \kappaup \, \yeta; \yeta \big) \big\}$ is the graph of the induced symplectomorphism 
    \begin{equation}
        \coTan \kappaup : \dotCoTan V \to \dotCoTan U, ~ \yeta \mapsto \varkappa 
        \yeta := \big( \kappaup^{-1} (y), \coTan_{x} \kappaup \, (\eta) \big). 
    \end{equation}
\end{example}
%
%
%

A careful look at $C$ in the preceding example reveals that $C$ is \textit{not} homogeneous because $\xi$ in~\eqref{eq: canonical_relation_pullback_Euclidean} is not necessarily non-zero even if $\eta$ is. 
As a consequence, $\Hess \varphi$ in $\symb{\fK}$ may vanish. 
This entails that $\kappaup^{*}$ is \textit{not} a \textit{homogeneous} Fourier integral operator. 
The assumption of homogeneity in this context ensures that a Fourier integral operator (and its dual) maps smooth functions to \textit{smooth functions} rather than  distributions 
(see e.g.~\cite[$(25.2.1)$]{Hoermander_Springer_2009}). 
We will discuss a plausible solution of this problem for the restriction map (a special case of the preceding example) in a bundle setting in Example~\ref{exm: restriction_op_FIO}. 
%
%
%
%
%
%
%
%
%
%
%
\section{Fourier Integral Operators on manifolds}
\label{sec: FIO_mf}
We would like to extend our preceding discussion on manifolds in this section. 
As pseudodifferential operators and their generalisation - Fourier integral operators are defined as integral operators by means of invariantly defined Schwartz kernels, it follows that those kernels must be integrable, at least in their second variable. 
However, there is no invariant way to integrate a function on a manifold. 
Thus, we are compelled to work with densities instead of functions and particularly, half-densities turn out to be convenient.
%
%
%
%
%
%
%
%
%
%
\subsection{Distributional half-densities}
\label{sec: distribution_density_mf}
\subsubsection{Densities}
We recall that an $s \in \C$-\textbf{density} on a $d$-dimensional real vector space $\sV$ is a map $\mu : \sV^{d} \to \C$ such that 
(see e.g.~\cite[Sec. 6.1-6.3]{Guillemin_InternationalP_2013},~\cite[pp. 374-377]{Treves_Plenum_1980},~\cite[pp. 428-432]{Lee_Springer_2013})  
\begin{equation}
	\mu (T v_{1}, \ldots, T v_{d}) = |\det T|^{s} \mu (v_{1}, \ldots, v_{d}) 
\end{equation}
for any linear transformation $T : \sV \to \sV$ and for all $(v_{1}, \ldots, v_{d}) \in \sV^{d}$. 
Such a quantity \textit{always exists} and the set of all such densities is denoted by $\varOmega^{s} (\sV)$. 
It is a $1$-dimensional vector space over $\C$. 
\newline 

The $s$-\textbf{density bundle} $\varOmega^{s} \! M \to M$ over a manifold $M$ is defined by 
\begin{equation}
    \forall s \in \C:~ \varOmega^{s} M := \bigsqcup_{x \in M} \varOmega^{s} (\tangent_{x} M), 
\end{equation}
which is a (smooth) \textit{trivialisable} $\C$-line bundle, i.e., for each $s$, $\varOmega^{s} (\tangent_{x} M) \cong \C$ but in a non-canonical fashion. 
A section of $\varOmega^{s} \! M$ are called an $s$-\textbf{density} on $M$; $s=1$-density is simply called a density. 
Let $\big( U, (x^{i}) \big)$ be a local chart around $x \in M$ and $\big( \tangent_{x} U, (\partial / \partial x^{i}) \big)$ (resp. $\big( \coTan_{x} U, (\rd x^{i}) \big)$) the induced natural chart for the tangent bundle $\tangent M \to M$ (resp. cotangent bundle $\coTanM \to M$). 
Then 
\begin{equation}
    |\rd x^{1} (x) \wedge \ldots \wedge \rd x^{d} (x)|^{s} \Big( \parDeri{x^{1}}{} (x) \wedge \ldots \wedge \parDeri{x^{d}}{} (x) \Big) = 1
\end{equation}
defines a canonical local section $|\rd x|^{s} := |\rd x^{1} \wedge \ldots \wedge \rd x^{d}|^{s} \in C^{\infty} (U; \varOmega^{s} M)$, where $\nicefrac{\partial}{\partial x^{1}} (x) \wedge \ldots \wedge \nicefrac{\partial}{\partial x^{d}} (x)$ is the generator of the line $\dot{\Lambda}^{d} \tangent_{x} M$ defined by local coordinates. 
Note, $|\rd x|^{s}$ defines a \textit{nowhere-vanishing} section of $\varOmega^{s} M$ by allowing $x$ to vary in the preceding equation. 
Any (compactly supported) smooth $s$-density $\mu$ on $M$ can locally be expressed as 
\begin{equation} \label{eq: halfdensity_local}
    \mu = f |\rd x|^{s} 
\end{equation}
for a unique $\C$-valued (compactly supported) smooth function $f = \mu (\partial/\partial x^{1} \wedge  \ldots \wedge \partial/\partial x^{d})$ on $U$.
A density $\mu$ on $M$ is said to be \textbf{positive} if its value at each point $x$ of $M$ has that property: $\mu (x) \, (X_{1}, \ldots, X_{d}) > 0$ for any linearly independent $d$-tuple of tangent vectors $(x; X_{1}, \ldots, X_{d})$ on $M$. 
Note, \textit{a smooth manifold $M$ always admits a smooth positive density on $M$} and the set of positive densities on $M$ is an open subset of $\varOmega \!M$ whose intersection with each fibre is convex. 

%
%
%
\begin{example} \label{exm: density}
	\begin{enumerate} [label=(\alph*)]
		\item \label{exm: density_Rd}
        Let $M := \Rd$ with standard coordinates $(x^{1}, \ldots, x^{d})$ and $f$ be a smooth function on $\Rd$. 
        A canonical trivialisation of the density bundle $\varOmega \Rd \to \Rd$ is provided by Lebesgue measure $\rd x$.  
        Then the most general smooth $s$-density on $\Rd$ is given by $\mu := f |\rd x|^{\alpha}$ and it is positive whenever $f$ is so. 
        \item 
        Let $M$ be an oriented manifold, i.e., it admits a volume form $\dVol$. Then, every smooth function $f$ on $M$ defines a smooth $s$-density $\mu := f |\dVol|^{s}$ on $M$, where $|\dVol|$ is a positive density, defined by $|\dVol|_{x} := |\dVol (x)|$ for each $x \in M$. 
        \item \label{exm: density_coTanM}
        Let $(\fM, \sigmaup)$ be a $2d$-dimensional symplectic manifold (see Section~\ref{sec: preliminary_conic_symplectic_geometry}) and $f \in C^{\infty} (\fM)$. 
        Such a manifold $\fM$ is always orientable because it admits a smooth canonical volume form, the \textbf{Liouville form} 
        \begin{equation} \label{eq: def_Liouville_form}
            \dVol :=(-1)^{d (d-1) / 2} \dfrac{\sigmaup^{d}}{d!}, 
            \quad 
            \sigmaup^{d} := \sigmaup \wedge \ldots \wedge \sigmaup 
        \end{equation} 
        and thus any smooth $s$-density on $\fM$ is given by $\mu := f |\dVol|^{s}$. 
        In particular, the cotangent bundle $\coTanM \to M$ over a $d$-dimensional manifold $M$ is a symplectic manifold and its Liouville form is given by (see Example~\ref{exm: cotangent_bundle_symplectic_mf})  
        \begin{equation}
            \dVol \upharpoonright U = \rd x^{1} \wedge \ldots \wedge \rd x^{d} \wedge \rd \xi_{1} \wedge \ldots \wedge \rd \xi_{d} = \rd x \, \rd \xi
        \end{equation}
        in the Darboux chart $\big( U, (x^{i}, \xi_{i}) \big)$.  
        Thus, any smooth $s$-density on $\coTanM$ is globally (resp. locally) given by $\mu = f |\dVol|^{s}$ (resp. $\mu = f (|\rd x| |\rd \xi|)^{s}$). 
        \item \label{exm: density_Lorentzian_mf} 
        Let $\spacetime$ be a $d$-dimensional Lorentzian manifold (see Section~\ref{sec: spacetime}). 
        At each $x \in \sM$, $\spacetime$ admits a \textit{canonical smooth volume element}, the \textbf{Lorentzian volume element} $\dVolg$, induced by the Lorentzian metric $\fg$ 
        (see e.g.~\cite[p. 195]{ONeill_Academic_1983}): 
        \begin{equation}
            \dVolg (x) (X_{1}, \ldots, X_{d}) := \sqrt{ |\det \big( \fg_{x} (X_{n}, e_{i}) \big)| }, 
        \end{equation}
        where $X_{n} \in \tangent_{x} \! \sM$ and $\{ e_{i} \}$ is an orthogonal basis for $\tangent_{x} \! \sM$. 
        In terms of a local oriented coordinate chart $\big( U, (x^{i}) \big)$ on $\sM$, this canonical volume element reads 
        \begin{equation} \label{eq: Lor_vol_form_coordinate}
            \dVolg (x) = \sqrt{| \det \fg |} \rd x^{1} \wedge \ldots \rd x^{d}. 
        \end{equation}
        Thus, every smooth $s$-density $\mu$ on $\sM$ is locally given by $\mu (x) = f (x) |\rd x|^{s}$ for a unique $f \in C^{\infty} (U)$.  
        However, $\dVolg$ \textit{might not exists globally} in contrast to $\mu$.  
        The (global) \textit{Lorentzian volume form exists if and only if $\sM$ is orientable} and then one identifies a smooth function on $\sM$ with a smooth $s$-density on $\sM$ employing the Lorentzian $s$-density $|\dVolg|^{s} := | \det \fg |^{\nicefrac{s}{2}} |\rd x|^{s}$:  
        \begin{equation}
            C^{\infty} (\sM) \ni f \mapsto f |\dVolg|^{s} \in C^{\infty} (\sM; \varOmega^{s} \! \sM). 
        \end{equation}
        \item 
        A semi-Riemannian hypersurface $(\varSigma, \iota)$ of an orientable semi-Riemannian manifold $\spacetime$ is orientable if and only if it admits a smooth unit normal vector field $\rN$. 
        In particular, if $c$ is a value of $f \in C^{\infty} (\sM, \R)$ then $\varSigma := f^{-1} (c)$ is a semi-Riemannian hypersurface of $\sM$ if and only if $\fg (\grad f, \grad f) \lessgtr 0$ on $\varSigma$. 
        In this case 
        \begin{equation}
            \rN := \frac{\grad f}{\norm{\grad f}{\fg}}
        \end{equation}
        is a (unique up to $\pm$) unit normal vector field along $\varSigma$. 
        Thus $\varSigma$ is orientable whose orientation is determined by $\rN$ and its volume form is given by 
        (see e.g.~\cite[Lem. 7.8, Prop. 4.17]{ONeill_Academic_1983},~\cite[Prop. 15.21]{Lee_Springer_2013}) 
        \begin{equation}
            \dVol_{\ms \varSigma} = \iota^{*} (\rN \lrcorner \dVol_{\ms \sM}), 
        \end{equation}
        where $\lrcorner$ is the interior multiplication. 
        \newline 

        To illustrate this, let us consider the standard $d$-sphere $\bbS^{d}$. 
        The inclusion $\iota : \bbS^{d} \hookrightarrow \R^{d+1}$ is a smooth embedding and the vector field $\rN := x^{i} \nicefrac{\partial}{\partial x^{i}}$ is nowhere tangent to $\bbS^{d}$ and thus $\rN$ induces a density on $\bbS^{d}$  via the preceding equation.  
        In particular, $\rd x = \rd x^{1} \wedge \rd x^{2} \wedge \rd x^{3}$ is a volume form of $\R^{3}$ and then the preceding equation entails that $|x^{1} \, \rd x^{2} \wedge \rd x^{3} + x^{2} \, \rd x^{3} \wedge \rd x^{1} + x^{3} \, \rd x^{1} \wedge \rd x^{2}|$ is the induced density on $\bbS^{2}$, where $(x^{1})^{1} + (x^{2})^{2} + (x^{3})^{2} = 1$ on $\bbS^{2}$ 
        (see e.g.~\cite[Exm. $15.22$, Ex. $16-9$]{Lee_Springer_2013}). 
	\end{enumerate}
\end{example}
%
%
%

In what follows, it is required to consider homogeneous (half-)densities on a punctured cotangent bundle and on its conic Lagrangian subsmanifolds. 
These are examples (see Examples~\ref{exm: cotangent_bundle_symplectic_mf} and~\ref{exm: diagonal_embedding_Lagrangian_submf}) of a conic symplectic manifold $(\cM, \fm_{\lambdaup})$ (see Definition~\ref{def: conic_mf} for a precise formulation). 
A half-density $\mu$ on $\cM$ is called \textbf{homogeneous} resp. \textbf{positively homogeneous} of degree $k \in \R$ if 
\begin{equation} \label{eq: def_homogeneous_density}
    \forall \lambdaup \in \dot{\R} \quad (\textrm{resp.} ~ \R_{+}) : 
    \fm_{\lambdaup}^{*} \mu = \lambdaup^{k} \mu. 
\end{equation}

%
%
%
\begin{example} \label{exm: order_homogeneous_halfdensity}
	\begin{enumerate}[label=(\alph*)]
		\item \label{exm: homogeneous_density_dotRd}
        The standard $s$-density (Example~\ref{exm: density}~\ref{exm: density_Rd}) $|\rd \theta|^{s}$ on $\dotRn$ is homogeneous of degree $sn$ 
        (see e.g.~\cite[p. 13]{Hoermander_Springer_2009},~\cite[$(2.1.13)$]{Duistermaat_Birkhaeuser_2011}).
        \item \label{exm: homogeneous_density_conormal_bundle}
        As in Example~\ref{exm: density}~\ref{exm: density_coTanM}, the conormal bundle (see Example~\ref{exm: conormal_bundle_Lagrangian_submf}) $(\varDelta M)^{\perp}$ of the diagonal in $M \times M$ is isomorphic to $\coTanM$. 
        Hence it has a natural half-density $|\rd x|^{1/2} |\rd \xi|^{1/2}$ of degree $\dim (M \times M) / 4 = d/2$ 
        (see e.g.~\cite[p. 104]{Hoermander_Springer_2007}).  
        \item \label{exm: homogeneous_density_canonical_relation} 
        The density $|\dVol_{\ms C}|$ on a canonical relation $C$ as constructed in Appendix~\ref{sec: volume_canonical_relation} and used in Definition~\ref{def: symbol_FIO_Euclidean} is of degree $n$ because the degree of homogeneity of the Lebesgue density (see~\eqref{eq: volume_fibre_critical_mf_local}) on the fibre-critical manifold is $\dM + \dN$ whereas that of the determinant in~\eqref{eq: volume_fibre_critical_mf_local} is $\dM + \dN - n$ 
        (see e.g.~\cite[p. $123$]{Grigis_CUP_1994}). 
	\end{enumerate}
\end{example}
%
%
%
%
%
%
%
%
%
%
\subsubsection{Distributional densities}
The pointwise product between a $s$-density and $1-s$ density is a $1$-density which permits us to define the following bilinear form 
(see e.g.~\cite[Sec. 4.2]{Scott_OUP_2010},~\cite[Sec. 1.1]{Duistermaat_Birkhaeuser_2011})  
\begin{equation}
    \scalarProdOne{\cdot}{\cdot} \, : C^{\infty} (M; \varOmega^{s} M) \times \CComInfinity{M; \varOmega^{1-s} M} \to \C, ~ (\mu, \phi) \mapsto \, 
	\scalarProdOne{\mu}{\phi} \, := \int_{M} \mu (x) \, \phi (x). 
\end{equation}
This form is continuous because 
$|\scalarProdOne{\mu}{\phi}| \leq \norm{\mu}{L^{1} (M; \sqrt{\varOmega})} \norm{\phi}{\infty}$. 
In this thesis only $s = 1/2$ will be used unless stated otherwise. 
Thus, we define the space $\cD' (M; \halfDenM) := \big( \CComInfinity{M; \halfDenM} \big)'$ of \textbf{distributional half-densities} on $M$ as the \textit{topological dual} of $C_{\rc}^{\infty} (M;$ $\halfDenM)$, i.e., the space of all $\C$-valued continuous $\C$-linear functionals on $\CComInfinity{M; \halfDenM}$. 
Throughout the thesis, $\cD' (M; \halfDenM)$ will be endowed with the weak $*$-topology. 
Note, one can define distributions $\cD' (M)$ on $M$ by patching up distributions $\cD' (U_{\alpha})$ on local charts $\{ (U_{\alpha}, \kappa_{\alpha}) \}_{\alpha}$ for $M$, where $\cD' (U_{\alpha})$ is obtained by pushing-forward the distributions $\cD' \big( \kappa_{\alpha} (U_{\alpha}); \varOmega \Rd \big)$. 
But, $\cD (M; \varOmega) \neq \cD' (M)$ 
(see e.g.~\cite[Sec. 6.3]{Hoermander_Springer_2003},~\cite[Sec. 4.2.2]{Scott_OUP_2010}). 
%
%
%
%
%
%
%
%
%
%
\subsection{Pseudodifferential operators}
\label{sec: PsiDO_mf}
To put forward the concept of a conormal distribution on a manifold, one can choose a local chart where the invariant Definition~\ref{def: conormal_distribution_Euclidean} applies and then, pullback the local pieces, glue together and prove chosen coordinate chart independence. 
However, it turns out that there is an intrinsic way to define this object which requires to introduce the notion of a Besov space as described below. 
\newline 

Let $1 \leq \kp \leq \infty$ and $m \in \R$. 
The \textbf{Besov space} ${}^{\kp}H^{m} (\Rd)$ is the Banach space of all temperate distributions $u \in \cS' (\Rd)$ on a Euclidean space $\Rd$ such that $\cF (u) \in L_{\mathrm{loc}}^{2} (\Rd)$ and the norm (for $\kp \neq \infty$) 
$^{\kp}\| u \|_{m} := \big( \sum_{i} \| (2 \pi)^{-d} \int_{\mathtt{U}_{i}} \re^{\ri x \cdot \xi} (\cF u) (\xi) \, \rd \xi \|_{m}^{\kp} \big)^{1/\kp}$ 
is finite; in case $\kp = \infty$ one sets 
$^{\infty}\| u \|_{m} := \big( \sup_{i} \| (2 \pi)^{-d} \int_{\mathtt{U}_{i}} \re^{\ri x \cdot \xi} (\cF u) (\xi) \, \rd \xi \|_{m}^{\kp} \big)^{1/\kp}$. 
Here, $\| \cdot \|_{s}$ denotes the Sobolev norm and $\mathtt{U}_{i} := \{ \xi \in \Rd ~|~ \mathtt{R}_{i-1} < |\xi| < \mathtt{R}_{i} \}, \mathtt{R}_{0} := 0, \mathtt{R}_{i} := 2^{i-1}$ for all $i > 0$ 
(see e.g.~\cite[Def. B.1.1]{Hoermander_Springer_2007}). 
If $U \subset \Rd$ is an open subset, then the \textbf{Besov space} ${}^{\kp} H_{\mathrm{loc}}^{s} (U)$ is defined as the set of all distributions $u \in \cD' (U)$ such that $\phi u \in {}^{\kp} H^{s} (\Rd)$ for every $\phi \in \CComInfinity{U}$. 
More generally, on a manifold $M$, the \textbf{Besov space} ${}^{\kp} H_{\mathrm{loc}}^{s} (M)$ is defined as the set of all distributions $u \in \cD' (M)$ such that $(\kappa^{-1})^{*} u \in {}^{\kp} H_{\mathrm{loc}}^{s} \big( \kappa (U) \big)$ for every smooth coordinate chart $(U, \kappa)$ for $M$. 
The topology is defined by the seminorms $u \mapsto {}^{\kp}\| \phi (\kappa^{-1})^{*} u \|_{s}$ where $\phi$ is an arbitrary element in $\CComInfinity{ \kappa (U) }$ 
(see e.g.~\cite[p. 475]{Hoermander_Springer_2007}). 

%
%
%
\begin{definition} \label{def: conormal_distribution_mf}
	Let $\halfDen \to M$ be the half-density bundle over a $d$-dimensional manifold $M$ and $S \subset M$ a closed submanifold. 
    Then the space $I^{m} (M, S; \halfDen)$ of \textbf{distributional half-densities} on $M$, \textbf{conormal} with respect to $S$ and of \textbf{degree} (at most) $m \in \R$, is defined as the set of all distributional half-densities $u \in \cD' (M; \halfDen)$ such that $L_{1} \ldots L_{N} u \in {}^{\infty} H_{\mathrm{loc}}^{- m - d / 4} (M; \halfDen)$ for all $N \in \NO$ and for all $L_{i} \in \PDO^{1} (M; \halfDen)$ with smooth coefficients tangential to $S$. 
    Here ${}^{\infty} H_{\mathrm{loc}}^{- m - d / 4} (M; \halfDen)$ denotes the Besov space and the topology is the weakest one which makes the maps $u \mapsto L_{1} \ldots L_{N} u \in {}^{\infty} H_{\mathrm{loc}}^{- m - d / 4} (M; \halfDen)$ continuous~\cite[Def. 18.2.6]{Hoermander_Springer_2007}. 
\end{definition}
%
%
%

\begin{definition} \label{def: PsiDO_mf}
	Let $\halfDen \to M$ be the bundle of half-densities over a manifold $M$ and $m \in \R$. 
    A \textbf{pseudodifferential operator} $P$ on $M$ of \textbf{order} (at most) $m$ is a continuous linear mapping  
    (see e.g~\cite[p. 100]{Hoermander_Springer_2007}) 
	\begin{equation}
		P : \CComInfinity{M; \halfDen} \to \CInfinity{M; \halfDen}, ~ u \mapsto 
		(P u) (x) := \int_{M} \fP (x, y) \, u (y)  
	\end{equation}
    whose Schwartz kernel $\fP$ is an element in the space of conormal distributions $I^{m} \big( M \times M, (\varDelta M)^{\perp *}; \halfDenMM \big)$ where $\varDelta : M \to M \times M$ is the diagonal embedding. 
    The set of all such operators is denoted by $\PsiDOM$. 
\end{definition}
%
%
%

The abstract definition implies that a pseudodifferential operator $P$ on a manifold $M$  is locally a collection $\{P_{\alpha}\}_{\alpha}$ of pseudodifferential operators  $P_{\alpha} \in \Psi \mathrm{DO}^{m} \big( \kappa_{\alpha} (U_{\alpha}) \big)$ on local charts $\{(U_{\alpha}, \kappa_{\alpha})\}_{\alpha}$ for $M$. 
To see this assertion, let $\sqrt{|\mu|}$ be an arbitrary but fixed half-density on $M$  and let $\hat{\rM}, \check{\rM} \in \End \big( \CComInfinity{U_{\alpha}} \big)$ be defined as multiplications by $\hat{\chi}, \check{\chi} \in \CComInfinity{U_{\alpha}}$, respectively. 
Then 
(see e.g.~\cite[p. 85]{Hoermander_Springer_2007},~\cite[Sec. 7.3]{vandenBan_2017}), 
\begin{eqnarray} \label{eq: localisation_PsiDO_mf}
    && P_{\alpha} := P_{\kappa_{\alpha}} := (\kappa_{\alpha}^{*})^{-1} \circ P \upharpoonright U_{\alpha} \circ \kappa_{\alpha}^{*} \in \Psi \mathrm{DO}^{m} \big( \kappa_{\alpha} (U_{\alpha}) \big), 
    \nonumber \\ 
    && P \upharpoonright U_{\alpha} := \hat{\rM} \circ P \circ \check{\rM} \in \PsiDO{m}{U_{\alpha}} 
\end{eqnarray}
is a local representative of $P$ with respect to the local charts and the chosen half-density. 
By Proposition~\ref{prop: PsiDO_psPsiDO_smooth}, the Schwartz kernel $\fP_{\alpha} (\bx, \by) \sqrt{|\mu_{\alpha} (\bx)|} \otimes \sqrt{|\mu_{\alpha} (\by)|}$ of $P_{\alpha}$ must be of the form 
\begin{equation}
    \fP_{\alpha} (\bx, \by) = \int_{\Rd} \re^{\ri (\bx - \by) \cdot \bxi} \totSymb{P_{\alpha}} (\bx, \bxi) \frac{\rd \bxi}{(2 \pi)^{d}} 
    \mod C^{\infty} \big( \kappa_{\alpha} (U_{\alpha}) \times \kappa_{\alpha} (U_{\alpha}) \big)   
\end{equation}
and vice-versa. 
Here, the total symbol $\totSymb{P_{\alpha}} (\bx, \bxi)$ is given by~\eqref{eq: def_coordinate_expression_total_symbol_psPsiDO_Euclidean} and $\mu_{\alpha} (\bx)$ is the Euclidean representative (cf.~\eqref{eq: halfdensity_local} and Example~\ref{exm: density}\ref{exm: density_Rd}) of the chosen density at $\bx := \kappa_{\alpha} (x) \in \kappa_{\alpha} (U_{\alpha})$.  
In general, it is necessary to have sufficiently many local patches $U_{\alpha}$ for $M$ so that the products $U_{\alpha} \times U_{\alpha}$ form an atlas for $M \times M$. 
But $\fP$ is off-diagonally smooth, and thus, it suffices to consider just an atlas $\{ (U_{\alpha}, \kappa_{\alpha}) \}_{\alpha}$ for $M$.  

%
%
%
\begin{remark} \label{rem: symbol_mf_coordinate_change}
    Let $\dotCoTanM \to M$ be the punctured cotangent bundle over a manifold $M$ and let $\{ (U_{\alpha}, \kappa_{\alpha}) \}_{\alpha}$ (resp. $\{ (\cU_{\alpha}, \varkappa_{\alpha}) \}_{\alpha}$) be an atlas (resp. a homogeneous symplectic atlas (see Theorem~\ref{thm: homogeneous_Darboux} and Example~\ref{exm: cotangent_bundle_symplectic_mf})) for $M$ (resp. $\dotCoTanM$) such that 
    \begin{equation} \label{eq: Darboux_chart_coTanM}
        \varkappa_{\alpha} : \cU_{\alpha} \to \kappa_{\alpha} (U_{\alpha}) \times \dotRd.     
    \end{equation} 
    The preceding map allows us to pullback $\totSymb{P_{\alpha}}$ to obtain $\totSymb{P \upharpoonright U_{\alpha}} \in S^{m} (\cU_{\alpha})$ (see Appendix~\ref{sec: symbol_mf}). 
    But, $\totSymb{P_{\alpha}}$ does not behave nicely under a coordinate transformation as observed in~\eqref{eq: total_symbol_PsiDO_Euclidean_coordinate_change} and thus it is impossible to patch up the local isomorphisms $\PsiDO{m}{U_{\alpha}} \cong S^{m} (\cU_{\alpha})$ to obtain a global one. 
    \newline 

    Analogously, $\symb{P \upharpoonright U_{\alpha}}$ is simply obtained by pulling back the leading order homogeneous term $\symb{P_{\alpha}}$ in $\totSymb{P_{\alpha}}$: 
    \begin{equation}
        \symb{} : \PsiDO{m}{U_{\alpha}} \to S^{m} (\cU_{\alpha}), ~P \upharpoonright U_{\alpha} \mapsto 
        \symb{P \upharpoonright U_{\alpha}} := \varkappa_{\alpha}^{*} \symb{P_{\alpha}}. 
    \end{equation}
    In contrast to $\totSymb{[P]}$, $\symb{[P]}$ can be promoted to the global isomorphism $\PsiDO{m}{M} \cong S^{m} (\coTanM)$ due to the fact that $\symb{P \upharpoonright U_{\alpha}}$ is an \textit{invariantly} defined \textit{homogeneous function} on $\cU_{\alpha}$ as Remarked in~\ref{rem: symbol_PsiDO_Euclidean_coordinate_change} 
    (see e.g.~\cite[pp. 85-86]{Hoermander_Springer_2007}). 
\end{remark}
%
%
%

\begin{definition} \label{def: symbol_PsiDO_mf}
    As in the terminologies of Definition~\ref{def: PsiDO_mf}, let $S^{m} (\coTanM)$ be the space of half-density-valued symbols (see Appendix~\ref{sec: symbol_mf}) on the cotangent bundle $\coTanM$ over $M$. 
    Then, the \textbf{principal symbol} $\symb{P}$ of $P \in \PsiDOM$ is defined by the isomorphism 
    (see e.g.~\cite[$(18.1.29)'$]{Hoermander_Springer_2007})
    \begin{equation} \label{eq: def_symbol_PsiDO_mf}
        \symb{} : \PsiDO{m - [1]}{M; \halfDen} \to S^{m - [1]} (\coTanM),  
	\end{equation}
    where $\Psi \mathrm{DO}^{m-[1]} := \Psi \mathrm{DO}^{m} / \Psi \mathrm{DO}^{m-1}$. 
    In a homogeneous symplectic chart (as in Remark~\ref{rem: symbol_mf_coordinate_change}) $\{ (\cU_{\alpha}, \varkappa_{\alpha}) \}_{\alpha}$ for $\dotCoTanM$, the preceding isomorphism is given by 
    \begin{equation}
        [P] \mapsto \symb{[P]} \xxi := \varkappa_{\alpha}^{*} (\symb{[P_{\alpha}]}) \sqrt{|\dVol_{\ms \coTanM} \xxi |} \mod S^{m-1} (\cdot), 
    \end{equation}
    where $\dVol_{\ms \coTanM}$ is the natural volume form (Example~\ref{exm: density}\ref{exm: density_coTanM}) on $\coTanM$ and $\symb{[P_{\alpha}]}$ is the principal symbol of a Euclidean representative $P_{\alpha} :=$~\eqref{eq: localisation_PsiDO_mf} of $P$. 
\end{definition}
%
%
%

At this point one may wonder whether the next-to-leading order term in $\totSymb{[P]}$ has some invariant meaning. 
Unfortunately, this naive expectation does not go well, rather the following combination works. 

%
%
%
\begin{definition} \label{def: subprincipal_symbol_mf}
    As in the terminologies of Definitions~\ref{def: PsiDO_mf} and~\ref{def: symbol_PsiDO_mf}, the \textbf{subprincipal symbol} $\subSymb{P}$ of $P \in \PsiDOM$ is defined by the isomorphism~\cite[$(5.2.8)$]{Duistermaat_ActaMath_1972} 
    (see also~\cite[Thm. 18.1.33]{Hoermander_Springer_2007}) 
    \begin{eqnarray} \label{eq: def_subprincipal_symbol_mf}
        \subSymb{} & : & \PsiDO{m - [2]}{M; \halfDen} \to S^{m - 1 - [1]} (\coTanM), 
        \nonumber \\ 
        && [P] \mapsto \subSymb{[P]} \xxi := p_{m-1} | \dVol_{\ms \coTanM} \xxi |^{\frac{1}{2}} + \frac{\ri}{2} \frac{\partial^{2} \symb{P}}{\partial x^{i} \partial \xi_{i}} \xxi \mod S^{m-2} (\cdot), \quad 
    \end{eqnarray}
    where $p_{m-1}$ is the next-to-leading order term in the expression (cf.~\eqref{eq: def_coordinate_expression_total_symbol_psPsiDO_Euclidean}) of the total symbol of $P$, pulled back to $\cU_{\alpha}$. 
\end{definition}
%
%
%

\begin{remark} \label{rem: subprincipal_symbol_PsiDO_mf_coordinate_change}
	$\subSymb{P \upharpoonright U}$ transforms as ~\cite[Prop. 5.2.1]{Duistermaat_ActaMath_1972} 
	(see also~\cite[$(18.1.33')$]{Hoermander_Springer_2007}):    
	\begin{equation}
        \big( (\coTan \kappa)^{*} \subSymb{P} \big) (\bx, \bxi) = \subSymb{P_{\kappa}} \big( \bx, \coTan_{x} \kappa \, \bxi \big) 
	\end{equation}
	under a coordinate change $U \ni x \mapsto \bx := \kappa (x) \in \kappa (U)$. 
    Thus it is an \textit{invariantly defined homogeneous function} on $\dotCoTanM$ (cf. Remark~\ref{rem: symbol_PsiDO_Euclidean_coordinate_change}). 
    \newline 

    We will compute subprincipal symbol of $P \in \PsiDO{m}{\sM}$ on a Lorentzian manifold $\spacetime$ in Chapter~\ref{ch: Feynman_propagator}.  
    Then, one has the \textit{canonical} Lorentzian volume form $\dVolg$ (Example~\ref{exm: density}~\ref{exm: density_Lorentzian_mf}) so that $\symb{P} = p \sqrt{|\dVolg|}$, and $\subSymb{P}$ is given by 
	(see e.g.~\cite[Rem. 2.1.10]{Safarov_AMS_1997})  
	\begin{equation} \label{eq: subprincipal_symbol_bundle_Lor_mf}
		\subSymb{P} \xxi = p_{m-1} \xxi + \frac{\ri}{2} \frac{\partial^{2} p}{\partial x^{i} \partial \xi_{i}} \xxi + \frac{\ri}{2} \Gamma_{ji}^{j} (x) \frac{\partial p}{\partial \xi_{i}} (\xi) 
		\mod S^{m - 2} (\dotCoTansM), 
	\end{equation} 
    where $p, p_{m-1}$ are functions and $\Gamma_{ji}^{j}$ are the contracted Christoffel symbols on $\sM$. 
    Therefore, subprincipal symbol is \textit{independent} of \textit{local coordinates} but \textit{depends} on the \textit{chosen half-density}. 
\end{remark}
%
%
%

\begin{center}
	\begin{tikzpicture}
        \node (a) at (0, 0) {$C_{\rc}^{\infty} \big( \kappa (U) \big)$};
        \node (b) at (4, 0) {$C^{\infty} \big( \kappa (U) \big)$};
        \node (c) at (0, 2) {$\CComInfinity{U}$};
        \node (d) at (4, 2) {$C^{\infty} (U)$};
        \draw[->] (a) -- (b); 
        \draw[->] (c) -- (d); 
        \draw[->] (a) -- (c); 
        \draw[->] (b) -- (d); 
        \node[below] at (2, 0) {$P_{\kappa}$}; 
        \node[above] at (2, 2) {$P$}; 
        \node[left] at (0, 1) {$\kappa^{*}$}; 
        \node[right] at (4, 1) {$\kappa^{*}$}; 
        \node (e) at (6.5, 0) {$S^{m} \big( \kappa (U) \big)$}; 
        \node (f) at (6.5, 2) {$S^{m} (\coTan U)$}; 
        \node (g) at (10, 0) {$\Psi \mathrm{DO}^{m} \big( \kappa (U) \big)$}; 
        \node (h) at (10, 2) {$\PsiDO{m}{U; \halfDen}$}; 
        \node (i) at (14, 0) {$S^{m-1} \big( \kappa (U) \big)$}; 
        \node (j) at (14, 2) {$S^{m-1} (\coTan U)$}; 
        \draw[->] (g) -- (e); 
        \draw[->] (h) -- (f); 
        \draw[->] (f) -- (e); 
        \draw[->] (h) -- (j); 
        \draw[->] (g) -- (i); 
        \draw[->] (j) -- (i); 
        \draw[->] (h) -- (g); 
        \node[right] at (6.5, 1) {$(\coTan \kappa)^{*}$}; 
        \node[right] at (10, 1) {$\eqref{eq: localisation_PsiDO_mf}$}; 
        \node[right] at (14, 1) {$(\coTan \kappa)^{*}$}; 
        \node[above] at (8, 2) {$\symb{P}$}; 
        \node[below] at (8.25, 0) {$\symb{P_{\kappa}}$}; 
        \node[above] at (12, 2) {$\subSymb{P}$}; 
        \node[below] at (12, 0) {$\subSymb{P_{\kappa}}$}; 
    \end{tikzpicture}
    \captionof{figure}[Transformation of (sub)principal symbol of a pseudodifferential operator]{Transformation  
        of a pseudodifferential operator, its principal symbol and subprincipal symbol under a coordinate chart $(U, \kappa)$. 
    }
\end{center}
%
%
%

We now introduce the following terminologies required to define a special class of pseudodifferential operators which is going to play a pivotal role in this thesis. 

%
%
%
\begin{definition} \label{def: bicharacteristic_strip}
	Let $M$ be a manifold and $H$ a Hamiltonian function on its cotangent bundle $\coTanM$. 
    A \textbf{bicharacteristic} (strip) of $H$ through a point $\xxiNot \in \coTanM$ is the integral curve 
    $\gamma : \R \supseteq \I \to \coTanM$ of the Hamiltonian vector field $X_{H}$ (see~\eqref{eq: def_HVF}) with initial data $\gamma (s_{0}) := \xxiNot$. 
    If $H \xxiNot = 0$ then $\gamma$ is called the \textbf{null bicharacteristic} (strip) of $H$. 
    The projection $c : \I \to M$ of $\gamma$ on $M$ is called the \textbf{bicharacteristic} (curve). 
\end{definition}
%
%
%

\begin{definition} \label{def: real_principal_type_PsiDO_mf}
	As in the terminologies of Definition~\ref{def: PsiDO_mf}, any properly supported $P \in \PsiDOM$ is of \textbf{real principal type} on a manifold $M$ if $P$ has a real homogeneous principal symbol $p$ of order $m$ and no complete null bicharacteristic (strip) of $p$ stays over a compact set in $M$~\cite[Def. 6.3.2]{Duistermaat_ActaMath_1972}.  
\end{definition}
%
%
%

\begin{remark} \label{rem: real_principal_type_PsiDO_mf_HVF_radial}
	The non-trapping property for the bicharacteristics of $p$ locally means that $\rd p$ is not proportional to the Liouville $1$-form (see Example~\ref{exm: cotangent_bundle_symplectic_mf}) of $\coTanM$ at any $\xxiNot \in \Char P$. 
    Equivalently, this can be rephrased as that $X_{p}$ and the radial vector field (see Definition~\ref{eq: def_radial_VF}) are linearly independent at every $\xxiNot \in \Char P$: 
    \begin{equation}
    	\forall \xxiNot \in \Char P: X_p \xxiNot \neq 0. 
    \end{equation}
    This property implies that $\Char P$ is a \textit{smooth, closed conic hypersurface of dimension $2d -1$}
    (see e.g.~\cite[p. 54]{Hoermander_Springer_2009},~\cite[Sec. VIII.7]{Treves_Plenum_1980}). 
\end{remark}
%
%
%
%
%
%
%
%
%
%
\subsection{Lagrangian distributions}  
\label{sec: Lagrangian_distribution_mf}
To define a Lagrangian distribution on $M$ one observes that $I^{m} (M, S; \halfDen)$ is the \textit{largest subspace} of the Besov space (introduced in Section~\ref{sec: PsiDO_mf}) ${}^{\infty} H_{\mathrm{loc}}^{-m - d/4} (M; \halfDen)$ which is \textit{invariant} under the action of all \textit{first-order partial differential operators tangent} to $S$. 
Moreover, any first-order pseudodifferential operator having vanishing principal symbol on the conormal bundle $S^{\perp *}$ of $S$, keeps $I^{m} (M, S; \halfDen)$ invariant~\cite[Thm. 18.2.12]{Hoermander_Springer_2007}. 
Since conormal bundles are a particular case of Lagrangian submanifolds (see Example~\ref{exm: conormal_bundle_Lagrangian_submf}), the last fact is used to define Lagrangian distributions. 

%
%
%
\begin{definition} \label{def: Lagrangian_distribution_mf} 
	Let $\halfDen, \dotCoTanM \to M$ be the bundle of half-densities and the punctured cotangent bundle over a $d$-dimensional manifold $M$ and $\varLambda \subset \dotCoTanM$ a closed conic Lagrangian submanifold (see Definition~\ref{def: Lagrangian_submf}).  
    Then the space $\LagrangianDistM$ of \textbf{Lagrangian distributional half-densities} on $M$, of \textbf{order} (at most) $m \in \R$, is defined as the set of all distributional half-densities $u \in \cD' (M; \halfDen)$ on $M$ such that $P_{1} \ldots P_{N} u \in {}^{\infty} H_{\mathrm{loc}}^{-m - d/4} (M; \halfDen)$ for all $N \in \NO$ and for all properly supported $P_{i} \in \PsiDO{1}{M; \halfDen}$ having vanishing principal symbols on $\varLambda$, where ${}^{\infty} H_{\mathrm{loc}}^{-m - d/4} (M; \halfDen)$ denotes the Besov space~\cite[Def. 25.1.1]{Hoermander_Springer_2009}.  
\end{definition}
%
%
%

In pursuance of utilising $\LagrangianDistM$ for the Schwartz kernel of a Fourier integral operator, $M$ has to be replaced by a product manifold $M \times N$ where $N$ is an arbitrary manifold and then the role of $\varLambda$ will be played by a canonical relation (see Definition~\ref{def: canonical_relation}). 

%
%
%
\begin{definition} \label{def: FIO_mf}
	Let $M, N$ be manifolds endowed with the bundles of half-densities $\halfDenM \to M, \halfDenN \to N$. 
	Suppose that $\dotCoTanM$ (resp. $\dotCoTanN$) is the punctured cotangent bundle of $M$ (resp. $N$) and that $C$ is a homogeneous canonical relation form $\dotCoTanN$ to $\dotCoTanM$ which is closed in $\dotCoTanMN$. 
	A linear continuous operator~\cite[Def. 25.2.1]{Hoermander_Springer_2009} 
	\begin{equation}
		A : \CComInfinity{N; \halfDenN} \to \cD' (M; \halfDenM),~ v \mapsto (A v) (x) 
		:= \int_{N} \fA (x, y) \, v (y)   
	\end{equation}
    whose Schwartz kernel $\fA$ belongs to the Lagrangian distribution $\LagrangianDistMN$, is called a \textbf{Fourier integral operator} of \textbf{order} (at most) $m \in \R$ from compactly supported half-densities $\CComInfinity{N; \halfDenN}$ on $N$ to distributional half-densities $\cD' (M; \halfDenM)$ on $M$ associated with the \textbf{canonical relation} $C$.  
	We will denote set of all such operators by $\FIONM$. 
\end{definition}
%
%
%

Similar to the pseudodifferential operators, Fourier integral operators can also be written as a collection of local Fourier integral operators on atlases for $M$ and for $N$.  
But their localisation is intricate because a generic canonical relation is more complicated than the conormal bundle $(\varDelta M)^{\perp *}$ and the requisite to localise in covectors together with their base points makes our previous localisation by cutoff functions inadequate. 
We localise an element of $\LagrangianDistMN$ by making use of 

%
%
%
\begin{lemma} \label{lem: localisation_Lagrangian_dist_mf}
	As in the terminologies of Definition~\ref{def: Lagrangian_distribution_mf}, if $u \in \LagrangianDistM$ then  
	\begin{equation}
        \WF u \subset \varLambda 
    \end{equation}
    and $Pu \in \LagrangianDistM$ if $P \in \PsiDO{0}{M; \halfDen}$. 
    Conversely, $u \in \LagrangianDistM$ if for every $\xxi \in \dotCoTanM$ one can find properly supported $Q \in \PsiDO{0}{M; \halfDen}$ non-characteristic at $\xxi$ such that $Qu \in \LagrangianDistM$~\cite[Lem. 25.1.2]{Hoermander_Springer_2009}.
\end{lemma}
%
%
%

Therefore, appropriate pseudodifferential operators play the same role for Lagrangian distributions that cutoff functions (cf.~\eqref{eq: localisation_PsiDO_mf}) do for the conormal bundle $(\varDelta M)^{\perp *}$ and the preceding lemma reduces the investigation of $\fA \in \LagrangianDistMN$ to the case where $\WF \fA$ is contained in a small closed conic neighbourhood of $\xxiyeta \in C'$ and $(x, y) \in \supp \fA$.  
More generally, given a homogeneous canonical relation $C$ of the open cone $\kU \times \kV \subset \dotCoTan M \times \dotCoTan N$, one says that $\fA \in \LagrangianDistMN$ at $\xxiyeta \in \kU \times \kV$ if there is an open conic neighbourhood $\cU \times \cV$ of $\xxiyeta$ contained in $\kU \times \kV$ such that $(P \otimes Q) \fA \in \LagrangianDistMN$ for all properly supported $P \in \PsiDO{0}{M; \halfDenM}, Q \in \PsiDO{0}{N; \halfDenN}$ having $\WF P \subset \clo (\cU), \WF Q \subset \clo (\cV)$. 
In fact, it suffices to know this for some such $P \otimes Q$ non-characteristic at $\xxiyeta$~\cite[p. 5]{Hoermander_Springer_2009}. 
\newline 

In order to get a local description of $\fA$, let $\{ (U_{\alpha}, \kappa_{\alpha}) \}_{\alpha}$ (resp. $ \{ (V_{\beta}, \rho_{\beta}) \}_{\beta}$) be an atlas for $M$ (resp. $N$). 
Since $\dotCoTanM$ and $\dotCoTanN$ are conic symplectic manifolds (see Definition~\ref{def: conic_mf} and Example~\ref{exm: cotangent_bundle_symplectic_mf}), there exists a homogeneous atlas $\{ (\cU_{\alpha}, \varkappa_{\alpha}) \}_{\alpha}$ (resp. $ \{ (\cV_{\beta}, \varrho_{\beta}) \}_{\beta}$) for $\dotCoTanM$ (resp. $\dotCoTanN$) such that 
\begin{equation} \label{eq: Darboux_atlas_coTanM_coTanN}
	\varkappa_{\alpha} : \cU_{\alpha} \to \kappa_{\alpha} (U_{\alpha}) \times \dot{\R}^{\dM} , 
	\qquad 
	\varrho_{\beta} : \cV_{\beta} \to \rho_{\beta} (V_{\beta}) \times \dot{\R}^{\dN},   
\end{equation}
by the homogeneous Darboux Theorem~\ref{thm: homogeneous_Darboux}. 
For each $\xxiyeta \in C'$, suppose that $[\varphi_{\alpha, \beta}]$ is an element of the stable equivalence class of phase functions for $C'$ with excess $e_{\alpha, \beta}$ (see Definition~\ref{def: clean_phase_function} and Remark~\ref{rem: stable_equivalence}), defined on an open conic neighbourhood of $\big( \bx := \kappa_{\alpha} (x), \by := \rho_{\beta}(y), \theta \big)$ in $\kappa_{\alpha} (U_{\alpha}) \times \rho_{\beta} (V_{\beta}) \times \dot{\R}^{n_{\alpha, \beta}}$ for some $n_{\alpha, \beta} \in \N$, and that 
\begin{equation} \label{eq: fibre_critical_mf_Euclidean}
    \sC_{\alpha, \beta} := \{ (\bx_{0}, \by_{0}; \theta^{0}) \in \kappa_{\alpha} (U_{\alpha}) \times \rho_{\beta} (V_{\beta}) \times \dot{\R}^{n_{\alpha, \beta}} | \grad_{\theta} \varphi_{\alpha, \beta} \, (\bx_{0}, \by_{0}; \theta^{0}) = 0 \}
\end{equation}
is the fibre-critical manifold of $\varphi_{\alpha, \beta}$. 
Then (see~\eqref{eq: def_canonical_relation_fibration}) 
\begin{equation} \label{eq: Euclidean_rep_canonical_relation}
    C'_{\alpha, \beta} := \{ (\bx, \bxi; \by, \bfeta) \in \varkappa_{\alpha} (\cU_{\alpha}) \times \varrho_{\beta} (\cV_{\beta}) | (\bx, \by; \theta) \in \sC_{\alpha, \beta}, \bxi = \rd_{\bx} \varphi_{\alpha, \beta}, \bfeta = \rd_{\by} \varphi_{\alpha, \beta}) \}
\end{equation}
is a Euclidean representative of $C'$ in a conic neighbourhood $\cU_{\alpha} \times \cV_{\beta}$ of $\xxiyeta \in C'$.  
Hence, locally 
\begin{equation}
    \fA = \sum_{\alpha, \beta} (\kappa_{\alpha} \times \rho_{\beta})_{*}^{-1} \, \fA_{\alpha, \beta}, 
\end{equation}
where module $C^{\infty} \big( \kappa_{\alpha} (U_{\alpha}) \times \rho_{\beta} (V_{\beta}) \big)$, $\fA_{\alpha, \beta}$ is an oscillatory integral of the form~\cite[Def. 3.2.2]{Hoermander_ActaMath_1971},~\cite[Prop. 25.1.5']{Hoermander_Springer_2009}:    
\begin{equation} \label{eq: Lagrangian_distribution_mf_local}
	\fA_{\alpha, \beta} (\bx, \by) \equiv (2 \pi)^{- (\dM + \dN + 2 n_{\alpha, \beta} - 2 e_{\alpha, \beta}) / 4} \int_{ \R^{n_{\alpha, \beta}} } \re^{\ri \varphi_{\alpha, \beta} (\bx, \by, \theta)} \fa_{\alpha, \beta} (\bx, \by; \theta) \, \rd \theta  
\end{equation}
and $\supp \fA_{\alpha, \beta}$ is locally finite. 
Here, $\fa_{\alpha, \beta} \in S^{m + (\dM + \dN - 2 n_{\alpha, \beta} - 2 e_{\alpha, \beta}) / 4} (\R^{\dM} \times \R^{\dN} \times \R^{n_{\alpha, \beta}})$ having support in the interior of a sufficiently small conic neighbourhood of $\sC_{\alpha, \beta}$ contained in the domain of definition of $\varphi_{\alpha, \beta}$, the symbol $\equiv$ means modulo smoothing kernels, and $\rd \theta$ is the Lebesgue measure at $\theta \in \dot{\R}^{n_{\alpha, \beta}}$. 
\newline 

Since the principal symbol of a scalar Lagrangian distribution has good transformation property under a diffeomorphism (Remark~\ref{rem: symbol_FIO_Euclidean_diffeo}), one defines $\symb{\fA}$ in a local chart to identify with $\symb{\fA_{\alpha, \beta}}$ by enjoying the independence with respect to the chosen coordinates. 

%
%
%
\begin{definition} \label{def: symbol_FIO_mf}
    As in the terminologies of Definition~\ref{def: Lagrangian_distribution_mf}, let $\halfDenC, \M \to C$ be the half-density bundle (see Appendix~\ref{sec: volume_canonical_relation}) and the Keller-Maslov bundle (see Definition~\ref{def: Keller_Maslov_bundle}) over $C$, respectively, and $S^{m + (\dM + \dN) / 4 - [1]} (C; \M \otimes \halfDenC)$ the space (see~\eqref{eq: symbol_canonical_relation_Maslov_density_function}) of $\Maslov \otimes \halfDenC$-valued symbols on $C$. 
    Then the \textbf{principal symbol} of a Lagrangian distribution is defined by the isomorphism~\cite[$(7.6)$]{Duistermaat_InventMath_1975} 
    (see also~\cite[Prop. 25.1.5']{Hoermander_Springer_2009})   
	\begin{equation}
		\symb{} : I^{m - [1]} \big( M \times N, C'; \halfDenMN \big) \to S^{m + \frac{\dM + \dN}{4} - [1]} (C; \M \otimes \halfDenC),   
	\end{equation}
    where $I^{m-[1]} := I^{m} / I^{m-1}$. 
    In a homogeneous symplectic chart $(\cU_{\alpha} \times \cV_{\beta}, \varkappa_{\alpha} \times \varrho_{\beta})$ for $C'$ 
    (see e.g.~\cite[Thm. 21.2.8]{Hoermander_Springer_2007}) 
    induced from those $\{ (\cU_{\alpha}, \varkappa_{\alpha}) \}_{\alpha}$ (resp. $ \{ (\cV_{\beta}, \varrho_{\beta}) \}_{\beta}$) for $\dotCoTanM$ (resp. $\dotCoTanN$) with ~\eqref{eq: Darboux_atlas_coTanM_coTanN}, the preceding isomorphism is given by 
    \begin{eqnarray}
        [\fA] \mapsto \symb{[\fA]} \xxiyeta 
        \! & \! := \! & \!
        (\varkappa_{\alpha} \times \varrho_{\beta})^{*} \left( \int_{\kC} \rd \theta'' a_{\alpha, \beta} (\bx, \by; \theta', \theta'') \right) \m \xxiyeta \otimes 
        \nonumber\\ 
        && 
        \sqrt{|\dVol_{\ms C} \xxiyeta|} 
        \mod S^{m + \frac{d_{\ms U} + d_{\ms V}}{4} - 1} (\cdot), 
    \end{eqnarray}
    where $a_{\alpha, \beta}$ is the top-order homogeneous term of $\fa_{\alpha, \beta}$ in~\eqref{eq: Lagrangian_distribution_mf_local} and all other symbols are as in Definition~\ref{def: symbol_FIO_Euclidean}. 
\end{definition}
%
%
%

\begin{example} \label{exm: restriction_op_FIO_mf}
    Let $(\varSigma, \iota_{\ms \varSigma})$ be an $n$-dimensional immersed submanifold of a $d$- dimensional manifold $M$ and $\halfDenM \to M$ the bundle of half-densities over $M$. 
    Then the pullback 
    \begin{equation}
        \iota_{\ms \varSigma}^{*} : \CComInfinity{M; \halfDenM} \to C^{\infty} (\varSigma; \halfDenM_{\varSigma}) 
    \end{equation}
    of the injective immersion $\iota_{\ms \varSigma} : \varSigma \hookrightarrow M$ is called the \textbf{restriction map}. 
    It is known that $\varSigma$ admits an adapted atlas, i.e., an atlas whose charts are induced from the charts of $M$. 
    In other words, for each $x' \in \varSigma$, there exist a chart $(V, \rho)$ for $\varSigma$ at $x'$ and a chart $(U, \kappa)$ for $M$ at $\iota_{\ms \varSigma} (x')$ such that $\iota_{\ms \varSigma} (V) \subset U$ and 
    \begin{equation}
        \kappa \circ \iota_{\ms \varSigma} \circ \rho^{-1} (\bx') = (0, \ldots, 0, \bx^{d - n + 1}, \ldots, \bx^{d}) 
    \end{equation}
    is a chart for $\varSigma$. 
    We read off, as a special case of Example~\ref{exm: pullback_Euclidean_FIO} that $\iota_{\ms \varSigma}^{*}$ is a Fourier integral operator whose Schwartz kernel is given by 
    \begin{eqnarray}
    	\fR_{\varSigma} (\bx', \by) 
        & := & 
        (2 \pi)^{- (d-n + n+3d) / 4} \int_{\Rd} \re^{\ri \varphi_{\ms \varSigma}} \rd \bfeta, 
        \nonumber \\ 
        \varphi_{\ms \varSigma} 
        & := & 
        (\bx^{i} - \by^{i}) \boldsymbol{\eta}_{i} - \by^{j} \bfeta_{j}, 
        \quad 
        i = \codim \varSigma + 1, \ldots, d, \; j = 1, \ldots, \codim \varSigma, \quad 
    \end{eqnarray}
    where $\codim \varSigma := \dim M - \dim \varSigma$ is the codimension of $\varSigma$. 
    The phase function $\varphi_{\ms \varSigma}$ is non-degenerate and thus 
    (see e.g.~\cite[$(2.4.4), (5.1.2)$]{Duistermaat_Birkhaeuser_2011})
    \begin{subequations}
    	\begin{eqnarray}
            \fR_{\varSigma} & \in & I^{\codim \varSigma / 4} \big( \varSigma \times M, \varLambda_{\varSigma}'; \halfDen (M_{\varSigma} \times M) \big), 
            \\ 
            \varLambda_{\varSigma}' & := & \{ (x', \xi'; x, \xi) \in \coTanM_{\varSigma} \times \dotCoTanM \,|\, x = \iota_{\ms \varSigma} (x'), \xi' = \xi \upharpoonright \tangent_{x'} M_{\varSigma} \}, 
            \label{eq: def_canonical_relation_restriction_op_mf}
            \\ 
            \symb{\fR_{\varSigma}} & := & 
            (2 \pi)^{- \codim \varSigma / 4} \one_{\varOmega^{\ms 1/2} (M_{\varSigma} \times M)} \sqrt{|\dVol_{\ms \varLambda_{\varSigma}}|} \otimes \bbl, 
    	\end{eqnarray}
    \end{subequations}
    where $|\dVol_{\ms \varLambda}|$ is the density on $\varLambda_{\varSigma}$ and $\bbl$ is a section of the Keller-Maslov bundle $\bbL \to \varLambda_{\varSigma}$. 
    The canonical symplectic form $\sigmaup$ on $\coTanM$ is given by (see  Example~\ref{exm: cotangent_bundle_symplectic_mf}) $\sigmaup = \rd x^{i} \wedge \rd \xi_{i} + \rd x^{j} \wedge \rd \xi_{j}$ in the adapted coordinates and so $|\dVol_{\ms \varLambda}| = |\rd x' \wedge \rd \xi|$ is the induced density on $\varLambda_{\varSigma}$. 
    The bundle $\bbL$ comprises global constant sections (cf. Example~\ref{exm: Maslov_bundle_conormal_bundle}) constructed from $\varphi_{\ms \varSigma}$  
    (see e.g.~\cite{Toth_GFA_2013}).
\end{example}
%
%
%

As observed in Example~\ref{exm: pullback_Euclidean_FIO}, scrutinising $\varLambda_{\varSigma}$ in the preceding example, one perceives that $\varLambda_{\varSigma}$ is \textit{not} homogeneous because elements $\big( y', 0; \iota_{\ms \varSigma} (y'), \eta \big)$ occur in $\varLambda_{\varSigma} :=$~\eqref{eq: def_canonical_relation_restriction_op_mf} whenever $\eta \in \varSigma_{\iota_{\varSigma} (y')}^{\perp *}$ where $\varSigma^{\perp *}$ is the conormal bundle (see Example~\ref{exm: conormal_bundle_Lagrangian_submf}). 
Therefore, $\iota_{\ms \varSigma}^{*}$ is \textit{not} a \textit{homogeneous} Fourier integral operator. 
We will present a solution to this issue for the bundle version of the restriction operator in Example~\ref{exm: restriction_op_FIO}. 
%
%
%
%
%
%
%
%
%
%
\section{Fourier Integral Operators on vector bundles}
We will now describe the foremost technical tool of this thesis equipped with the preparations done in the last two sections. 
Let us begin with some rudimentary facts about integral operators acting on vector bundles. 
%
%
%
%
%
%
%
%
%
%
\subsection{Integral operators} 
\label{sec: integral_op}
Let $\sE, \halfDenM \to M$ (resp. $\sF, \halfDenN \to N$) be a vector bundle and bundle of half-densities over a $\dM$ (resp. $\dN$)-dimensional manifold $M$ (resp. $N$). 
Recall, the set $\dualComSecME$ of \textbf{half-density valued distributions} on $\sE$ is the $\halfDenM$-valued $\C$-linear functionals on $\CComInfinity{M; \sE^{*} \otimes \halfDenM}$ where $\sE^{*} \to M$ is the dual bundle (see Section~\ref{sec: adjoint}) of $\sE$. 
If $\pr_{\ms M} : M \times N \to M, \pr_{\ms N} : M \times N \to N$ are the Cartesian projectors then 
\begin{equation} \label{eq: def_exterior_tensor_product}
	\sE \boxtimes \sF := \pr_{\ms M}^{*} \sE \otimes \pr_{\ms N}^{*} \sF \to M \times N. 
\end{equation}
is the \textbf{exterior tensor bundle} between $\sE$ and $\sF$.  
Suppose that $\{ (U_{\alpha}, \kappa_{\alpha}) \}_{\alpha}$ (resp. $\{ (V_{\beta}, \rho_{\beta}) \}_{\beta}$) is an atlas for $M$ (resp. $N$) such that $\sE$ (resp. $\sF$) admits a local trivialisation over $U_{\alpha}$ (resp. $V_{\beta}$) and that $\{ (U_{\alpha}, \chiup_{\alpha}) \}_{\alpha}$ (resp. $\{ (V_{\beta}, \tauup_{\beta}) \}_{\beta}$) is a bundle-atlas for $\sE$ (resp. $\sF^{*}$). 
Let $|\rd x|$ (resp. $|\rd y|$) be an arbitrary but fixed nowhere-vanishing global section (Section~\ref{sec: distribution_density_mf}) of $\halfDenM$ (resp. $\halfDenN$)
and let $(\rE_{r})$ (resp. $(\rF^{k})$) be a frame for $\sE$ (resp.$\sF$). 
Locally, $\phi \in \comSecME, \psi \in \CComInfinity{N; \sF \otimes \halfDenN}, \varUpsilon \in C_{\mathrm{c}}^{\infty} \big( M \times N; (\sE^{*} \boxtimes \sF) \otimes \halfDenMN \big), u \in \dualComSecME, v \in \cD' (N; \sF \otimes \halfDenN), \fT \in \cD' \big( M \times N, (\sE \boxtimes \sF^{*}) \otimes \halfDenMN \big)$ can be expressed as  
(see e.g.~\cite[Sec. 4.2]{Scott_OUP_2010},~\cite[pp. 148-150]{Guenther_AP_1988})  
\begin{subequations}
	\begin{eqnarray}
        \phi & = & \phi^{r} \rE_{r} \otimes \sqrt{|\rd x|}, \quad \phi^{r} \in \CComInfinity{U_{\alpha}}, 
        \\  
        \psi & = & \psi^{k} \rF_{k} \otimes \sqrt{|\rd y|}, \quad \psi^{k} \in \CComInfinity{V_{\beta}}, 
        \\ 
        \varUpsilon & = & \varUpsilon_{r}^{k} (\rE^{r} \boxtimes \rF_{k}) \otimes \sqrt{|\rd x| |\rd y|}, \quad \varUpsilon_{r}^{k} \in \CComInfinity{U_{\alpha} \times V_{\beta}}, 
        \\ 
        u & = & u^{r} \rE_{r} \otimes \sqrt{|\rd x|}, \quad  u^{r} \in \cD' (U_{\alpha}), 
        \\  
        v & = & v^{k} \rF_{k} \otimes \sqrt{|\rd y|}, \quad  v^{k} \in \cD' (V_{\beta}), 
        \\ 
        \fT & = & \fT^{r}_{k} (\rE_{r} \boxtimes \rF^{k}) \otimes \sqrt{|\rd x| |\rd y|}, \quad \fT^{r}_{k} \in \cD' (U_{\alpha} \times V_{\beta}), 
	\end{eqnarray}
\end{subequations}
where $r = 1, \ldots, \rk \sE$ and $k = 1, \ldots, \rk \sF$. 
Thus, one has the bijections depending on the choices of coordinate-charts and bundle-charts 
\begin{subequations}
	\begin{eqnarray}
		&& 
		C^{\infty} (U_{\alpha}; \sE_{\alpha} \otimes \halfDen_{\alpha}) 
        \cong 
        \big( C^{\infty} (U_{\alpha}; \halfDen_{\alpha}) \big)_{\rk \sE \times 1}
        \cong 
        \Big( C^{\infty} \big( \kappa_{\alpha} (U_{\alpha}) \big) \Big)_{\rk \sE \times 1}, 
        \\ 
        && 
        \cD' (U_{\alpha}; \sE_{\alpha} \otimes \halfDen_{\alpha}) 
        \cong 
        \big( \cD' (U_{\alpha}; \halfDen_{\alpha}) \big)_{\rk \sE \times 1}
        \cong 
        \Big( \cD' \big( \kappa_{\alpha} (U_{\alpha}) \big) \Big)_{\rk \sE \times 1}, 
        \\ 
        && 
        C^{\infty} (V_{\beta}; \sF_{\beta} \otimes \halfDen_{\beta}) 
        \cong 
        \big( C^{\infty} (V_{\beta}; \halfDen_{\beta}) \big)_{\rk \sF \times 1}
        \cong 
        \Big( C^{\infty} \big( \rho_{\beta} (V_{\beta}) \big) \Big)_{\rk \sF \times 1}, 
        \\ 
        && 
        \cD' (V_{\beta}; \sF_{\beta} \otimes \halfDen_{\beta}) 
        \cong 
        \big( \cD' (V_{\beta}; \halfDen_{\beta}) \big)_{\rk \sF \times 1}
        \cong 
        \Big( \cD' \big( \rho_{\beta} (V_{\beta}) \big) \Big)_{\rk \sE \times 1}. 
	\end{eqnarray}
\end{subequations}
According to the Schwartz kernel theorem 
(see e.g.~\cite[p. 93]{Hoermander_Springer_2009}), 
given a bidistribution $\fT \in \cD' \big( M \times N, (\sE \boxtimes \sF^{*}) \otimes \halfDenMN \big)$, one has a unique continuous linear operator  
\begin{equation}
	T : \comSecNF \to \dualComSecME, ~ v \mapsto (T v) (x) := \int_{N} \fT (x, y) \, v (y) 
\end{equation}
and vice-versa. 
One defines the restricted operator 
\begin{equation}
	T^{\alpha}_{\beta} := \big( T (\psi \upharpoonright V_{\beta}) \big) \upharpoonright U_{\alpha} 
	: \CComInfinity{V_{\beta}; \sF_{\beta} \otimes \halfDen_{\beta}} \to \cD' (U_{\alpha}; \sE_{\alpha} \otimes \halfDen_{\alpha}), 
\end{equation}
which can be identified with a $\rk \sE \times \rk \sF$-matrix of scalar operators 
(see e.g.~\cite[Sec. 4.3.1]{Scott_OUP_2010}) 
\begin{equation}
	T^{\alpha, r}_{\beta, k} := \chiup_{\alpha *} \circ T^{\alpha}_{\beta} \circ \tauup_{\beta}^{*}  
	: \CComInfinity{V_{\beta}; \halfDen_{\beta}} \to \cD' (U_{\alpha}; \halfDen_{\alpha}). 
\end{equation}
This can be further boiled down to the level of Euclidean spaces 
\begin{equation}
	T^{\alpha, r, i}_{\beta, k, j} := \kappa_{\alpha *} \circ T^{\alpha, r}_{\beta, k} \circ \rho_{\beta}^{*}  
	: C_{\mathrm{c}}^{\infty} \big( \rho_{\beta} (V_{\beta}) \big) \to \cD' \big( \kappa_{\alpha} (U_{\alpha}) \big), 
\end{equation}
as depicted in the commutative diagram below in Fig.~\ref{fig: local_FIO_F_E_N_M_R}. 
%
%
%
\begin{center}
		\begin{tikzpicture}
			\node (a) at (0, 0) {$C_{\mathrm{c}}^{\infty} \big( \rho_{\beta} (V_{\beta}) \big)$}; 
			\node (b) at (8, 0) {$\cD' \big( \kappa_{\alpha} (U_{\alpha}) \big)$};
			\node (c) at (0, 2) {$\CComInfinity{V_{\beta}; \halfDen_{\beta}}$};
			\node (d) at (8, 2) {$\cD' (U_{\alpha}; \halfDen_{\alpha})$};
			\node (e) at (0, 4) {$\CComInfinity{V_{\beta}; \sF_{\beta} \otimes \halfDen_{\beta}}$};
			\node (f) at (8, 4) {$\cD' (U_{\alpha}; \sE_{\alpha} \otimes \halfDen_{\alpha})$};
			\node[above] at (4, 4) {$T^{\alpha}_{\beta}$};
			\node[above] at (4, 2) {$T^{\alpha, r}_{\beta, k}$}; 
			\node[above] at (4, 0) {$T^{\alpha, r, i}_{\beta, k, j}$}; 
			\node[right] at (0, 3) {$\tauup_{\beta}^{*}$};
			\node[left] at (8, 3) {$\chiup_{\alpha *}$}; 
			\node[right] at (0, 1) {$\rho_{\beta}^{*}$};
			\node[left] at (8, 1) {$\kappa_{\alpha *}$};
			\draw[->] (a) -- (b);
			\draw[->] (c) -- (d); 
			\draw[->] (e) -- (f);
			\draw[->] (a) -- (c);
			\draw[->] (c) -- (e); 
			\draw[->] (f) -- (d);
			\draw[->] (d) -- (b);
		\end{tikzpicture}
		\captionof{figure}[Euclidean representative of an integral operator]{Euclidean representative of an integral operator acting on vector bundles.}
    \label{fig: local_FIO_F_E_N_M_R}
	\end{center}
%
%
%
    
Therefore, we can always locally express an integral operator acting on half-density-valued vector bundles as a system of scalar integral operators acting on functions on subsets of the respective Euclidean spaces:   
\begin{equation} \label{eq: integral_op_system_scalar_integral_op}
	T 
	\leftrightarrow 
	\big( T^{\alpha, r, i}_{\beta, k, j} \big), 
	\quad r = 1, \ldots, \rk \sE; k = 1, \ldots, \rk \sF; i = 1, \ldots, \dM; j = 1, \ldots, \dN. 
\end{equation}
%
%
%
%
%
%
%
%
%
%
\subsection{Pseudodifferential operators} 
\label{sec: PsiDO}
The notion of the Besov space ${}^{\infty} H_{\mathrm{loc}}^{- m - d / 4} (M; \halfDen)$ introduced in Section~\ref{sec: PsiDO_mf} generalises for a vector bundle $\sE \to M$, ${}^{\infty} H_{\mathrm{loc}}^{- m - d / 4} (M; \sE \otimes \halfDen)$, in a straightforward way. 

%
%
%
\begin{definition} \label{def: conormal_distribution}
	Let $\sE, \halfDen \to M$ be a vector bundle and the half-density bundle over a $d$-dimensional manifold $M$, respectively, and $S$ a closed submanifold of $M$. 
    Then the space $I^{m} (M, S; \sE \otimes \halfDen)$ of \textbf{half-density-valued distributions} of $\sE$, \textbf{conormal} with respect to $S$ and of \textbf{degree} (at most) $m \in \R$, is defined as the set of all half-density-valued distributions $u \in \cD' (M; \sE \otimes \halfDen)$ on $\sE$ such that $L_{1} \ldots L_{N} u$ belongs in the Besov space ${}^{\infty} H_{\mathrm{loc}}^{- m - d / 4} (M; \sE \otimes \halfDen)$ for all $N \in \NO$ and for all $L_{i} \in \PDO^{1} (M; \sE \otimes \halfDen)$ with smooth coefficients tangential to $S$. 
    Here the topology is the weakest one which makes the maps $u \mapsto L_{1} \ldots L_{N} u \in {}^{\infty} H_{\mathrm{loc}}^{- m - d / 4} (M; \sE \otimes \halfDen)$ continuous~\cite[Def. 18.2.6]{Hoermander_Springer_2007}. 
\end{definition}
%
%
%

Let $\sF \to M$ be a vector bundle over a manifold $M$. 
We make use of the natural vector bundle isomorphism 
\begin{equation}
    \Hom{\sE, \sF} \cong \sE^{*} \otimes \sF \to M 
\end{equation}
in the definition below. 

%
%
%
\begin{definition} \label{def: PsiDO}
	Let $\sE, \sF, \halfDen \to M$ be two vector bundles and the half-density bundle over a manifold $M$, respectively.  
    A \textbf{pseudodifferential operator} $P$ from $\sE$ to $\sF$ of \textbf{order} (at most) $m \in \R$ is a continuous linear map 
    (see e.g.~\cite[p. 100]{Hoermander_Springer_2007}) 
	\begin{equation*}
		P : \comSecE \to \secMF, ~u \mapsto 
		(P u) (x) := \int_{M} \fP (x, y) \, u (y)  
	\end{equation*}
    whose Schwartz kernel $\fP$ belongs to the space conormal distributions $I^{m} \big( M \times M, (\varDelta M)^{\perp *}; \mathrm{Hom} (\sE,$ $\sF) \otimes \halfDenMM \big)$ where $\varDelta: M \to M \times M$ is the diagonal embedding. 
	We denote the set of all such operators by $\PsiDOMEF$. 
\end{definition}
%
%
%

As inscribed in Section~\ref{sec: integral_op}, locally $P$ (resp. its Schwartz kernel $\fP$) is in a bijection with a matrix of properly supported scalar pseudodifferential operators $(P_{r}^{k})$ (resp. scalar kernels $(\fP_{r}^{k})$) on respective Euclidean spaces up to smoothing operators (resp. kernels); cf. Proposition~\ref{prop: PsiDO_psPsiDO_smooth} 
(for details, see, 
for instance~\cite[Def. 18.1.32]{Hoermander_Springer_2007},~\cite[Sec. 1.5.3]{Scott_OUP_2010}): 
\begin{eqnarray}
    && P \leftrightarrow (P^{k}_{r})_{\rk \sF \times \rk \sE} \mod (\varPsi^{k}_{r})_{\rk \sF \times \rk \sE}, 
    \quad \varPsi^{k}_{r} \in \PsiDO{-\infty}{\kappa_{\alpha} (U_{\alpha})}, 
    \nonumber \\  
    && \fP \leftrightarrow \big( \fP^{k}_{r} \big)_{\rk \sF \times \rk \sE} \mod \big( \mathsf{\Psi}^{k}_{r} \big)_{\rk \sF \times \rk \sE},  
    \quad \mathsf{\Psi}^{k}_{r} \in I^{-\infty} \big( \kappa_{\alpha} (U_{\alpha}) \times \kappa_{\alpha} (U_{\alpha}) \big)
\end{eqnarray}
where (cf. Theorem~\ref{thm: characterisation_psPsiDO_Euclidean}) 
\begin{equation} \label{eq: def_coordinate_expression_total_symbol_PsiDO}
   	\fP^{k}_{r} (\bx, \by) = \int_{\Rd} \frac{\rd \bxi}{(2 \pi)^{d}} \re^{\ri (\bx - \by) \cdot \bxi} \totSymb{P^{k}_{r}} (\bx, \bxi), 
   	\quad 
   	\totSymb{P^{k}_{r}} (\bx, \bxi) \sim \sum_{\alpha \in \NO^{d}} \frac{1}{\alpha !} (\partial_{\bxi}^{\alpha} \rD_{\by}^{\alpha} \fp^{k}_{r}) (\bx, \bx, \bxi) 
\end{equation}
for some proper $\fp^{k}_{r} \in S^{m} \big( \kappa_{\alpha} (U_{\alpha}) \times \kappa_{\alpha} (U_{\alpha}) \times \Rd \big)$, $\rd \bx$ is the Lebesgue measure at each $\bx := \kappa (x)$, and $\sim$ means the asymptotic summation (see Definition~\ref{def: asymptotic_summation}).  
Therefore, one sets 
\begin{equation}
	\totSymb{[P]} \xxi 
    \equiv \totSymb{\big[ \big( P^{k}_{r} \big)_{\rk \sF \times \rk \sE} \big]} (\bx, \bxi) 
    \equiv \big( \totSymb{[P^{k}_{r}]} (\bx, \bxi) \big)_{\rk \sF \times \rk \sE}  
\end{equation}
for the total symbol of $P$.  

%
%
%
\begin{definition} \label{def: symbol_PsiDO}
	As in the terminologies of Definitions~\ref{def: PsiDO} and~\ref{def: symbol_PsiDO_mf}, let $\pi : \coTanM \to M$ be the cotangent bundle over $M$ and $S^{m} \big( \coTanM, \Hom{\sE, \sF} \big)$ the symbol space of $\Hom{\sE, \sF}$-valued half-densities (see Appendix~\ref{sec: symbol}) on $\coTanM$. 
    Then the \textbf{principal symbol} is defined by the isomorphism 
    (see e.g.~\cite[p. 92]{Hoermander_Springer_2007})  
	\begin{equation}
		\symb{} : \PsiDO{m - [1]}{M; \sE \otimes \halfDen, \sF \otimes \halfDen} \to S^{m - [1]} \big( \coTanM, \Hom{\sE, \sF} \big),  
	\end{equation}
    where $\Psi \mathrm{DO}^{m-[1]} := \Psi \mathrm{DO}^{m} / \Psi \mathrm{DO}^{m-1}$. 
    The preceding map is locally given by 
    \begin{equation}
        [P] \mapsto \symb{[P]} \xxi := \big( \symb{[P^{k}_{r}]} \big)_{\rk \sF \times \rk \sE} \sqrt{| \dVol_{\ms \coTanM} \xxi |} \mod S^{m - [1]} (\cdot),  
    \end{equation}
    where $\symb{[P^{k}_{r}]}$ is the principal symbol of a Euclidean representative $P_{r}^{k}$ of $P$ in any homogeneous symplectic chart for $\dotCoTanM$ and any bundle-charts for $\sE, \sF$. 
\end{definition}
%
%
%

\begin{definition} \label{def: subprincipal_symbol}
	As in the terminologies of Definitions~\ref{def: PsiDO},~\ref{def: symbol_PsiDO},  and~\ref{def: subprincipal_symbol_mf}, the \textbf{subprincipal symbol} of $P \in \PsiDOMEF$ is defined by the following isomorphism~\cite[Thm. 18.1.33]{Hoermander_Springer_2007}) 
	\begin{eqnarray}
		\subSymb{} & : & \PsiDO{m - [2]}{M; \sE \otimes \halfDen, \sF \otimes \halfDen} \to S^{m - [2]} \big( \coTanM, \Hom{\sE, \sF} \big), 
		\nonumber \\ 
		&& 
		[P] \mapsto \subSymb{[P]} \xxi := \big( \subSymb{[P^{k}_{r}]} \big)_{\rk \sF \times \rk \sE}, 
	\end{eqnarray}
    where $\subSymb{[P^{k}_{r}]}$ is the subprincipal symbol of $P_{r}^{k}$. 
\end{definition}
%
%
%

\begin{remark} \label{rem: subprincipal_symbol_bundle_chart_change}
    The principal symbol is \textit{independent} of the choice of trivialisations of $\sE, \sF$ but the subprincipal symbol depends on the choice which has been explored in-depth in Section~\ref{sec: connection_form_subprincipal_symbol}.  
    Therefore, combining with Remark~\ref{rem: subprincipal_symbol_PsiDO_mf_coordinate_change} we summarise that the subprincipal symbol is \textit{independent} of \textit{local coordinates} but \textit{depends} on the \textit{chosen half-density} and \textit{bundle chart}. 
\end{remark}
%
%
%

In order to purvey the notion of real principal type pseudodifferential operators in the bundle setting, let $\sG \to M$ be a vector bundle.  
Then there is the natural bilinear mapping 
(see e.g.~\cite[Lem. 8.3.6]{vandenBan_2017},~\cite[p. 85]{Scott_OUP_2010})  
\begin{eqnarray}
    \cdot \circ \cdot & : & S^{m} \big( \coTanM, \Hom{\sE, \sG} \big) \times  S^{m'} \big( \coTanM, \Hom{\sG, \sF} \big) \to S^{m + m'} \big( \coTanM, \Hom{\sE, \sF} \big), 
    \nonumber \\ 
    && (p, p') \mapsto (p' \circ p) \, \xxi := p' \xxi \big( p \xxi \big).     
\end{eqnarray}
Besides, we note that $S^{0} \big( \coTanM, \Hom{\sE, \sF} \big)$ has a distinguished element $\tId$ given by the identity $\one_{\Hom{\sE, \sF}} = \one_{\Hom{\pi^{*} \sE, \pi^{*} \sF}}$ homomorphism $(\pi^{*} \sE)_{\xxi} \ni p \mapsto  p \in (\pi^{*} \sF)_{\xxi}$: 
\begin{equation}
    \coTan_{x} M \ni \xxi \mapsto \tId \xxi := \one_{\Hom{\sE, \sF}_{\xxi}} \in \Hom{\pi^{*} \sE, \pi^{*} \sF}, 
\end{equation}
where $\pi : \coTanM \to M$ is the cotangent bundle of $M$. 

%
%
%
\begin{definition} \label{def: real_principal_type_PsiDO}
    As in the terminologies of Definitions~\ref{def: PsiDO} and~\ref{def: symbol_PsiDO},  $P \in \PsiDOMEF$ having principal symbol $p$ is of \textbf{real principal type} at $\xxiNot \in \dotCoTanM$, if there exists a symbol $\tilde{p} \in S^{0} \big( \cU, (\pi^{*} \Hom{\sE, \sF})_{\cU} \big)$ such that~\cite[Def. 3.1]{Dencker_JFA_1982} 
	\begin{equation} \label{eq: def_real_principal_type_PsiDO}
        (\tilde{p} \circ p) \xxi = q \one_{\pi^{*} \Hom{\sE, \sF}_{\xxi}} 
	\end{equation} 
	for any $\xxi$ in an open conic neighbourhood $\cU$ of $\xxiNot$, where $q$ is the principal symbol of some scalar real principal type operator $Q \in \PsiDOM$. 
	We say that $P$ is of real principal type in $\dotCoTanM$ if it is so at every $\xxiNot \in \dotCoTanM$.
 \end{definition} 
%
%
%
 
The Euclidean formulation~\eqref{eq: def_WF_Euclidean} of the wavefront set generalises in bundle setting in a straightforward way: 
\begin{equation}
    \forall u \in \cD' (M; \sE \otimes \halfDen) : \WF u := \bigcap_{\substack{P \in \PsiDOME \\ Pu \in \secE}} \Char P
\end{equation}
except the characterisation of $\Char P$ used in~\eqref{eq: def_WF_Euclidean} is required to generalise as follows.  

%
%
%
\begin{definition} \label{def: char_PsiDO}
    As in the terminologies of Definitions~\ref{def: PsiDO} and~\ref{def: symbol_PsiDO}, let $P \in \Psi\mathrm{DO}^{m} (M;$ $\sE \otimes \halfDen, \sF \otimes \halfDen)$ whose principal symbol is denoted by $p$. 
    The operator $P$ is said to be \textbf{non-characteristic} at some point $\xxiNot \in \dotCoTanM$ if $p \circ \tilde{p} - \one \in S^{-1} \big( \cU, \Hom{\sE, \sF}_{\cU} \big)$ in a conic neighbourhood $\cU$ of $\xxiNot$ for some $\tilde{p} \in S^{m} \big( \cU, \Hom{\sE, \sF}_{\cU} \big)$. 
    The set of all \textbf{characteristic points} of $P$ is denoted by $\Char{P}$ 
    (see e.g.~\cite[Def. 18.1.25]{Hoermander_Springer_2007},~\cite[Sec. 8.4]{vandenBan_2017}).  
\end{definition}
%
%
%

This definition implies that the microlocal ellipticity condition is equivalent to 
\begin{equation}
	\forall \xxi \in \cU : \det \, p \xxi \neq 0. 
\end{equation}
It also follows that  
\begin{equation}
    \WF u = \bigcup_{r = 1, \ldots, \rk \sE} \WF (u^{r})
\end{equation}
where $u^{r} \in \cD' (U)$ is a component of $u$ with respect to a trivialisation of $\sE$ over $U \subset M$ for a chosen half-density. 
\newline 

To illustrate these concepts we present the wave operator as an example
(see~\cite[Prop. 4.2-4.4]{Radzikowski_CMP_1996} 
for the scalar case together with~\cite[Prop. 2.1]{Low_JMP_1989},~\cite[Sec. 2]{Low_NonlinearAnal_2001}).  

%
%
%
\begin{example} \label{exm: NHOp_real_principal_type}
	Let $\square$ be a normally hyperbolic operator (see Definition~\ref{def: NHOp}) on a vector bundle $\sE \to \sM$ over a globally hyperbolic spacetime (see Definition~\ref{def: globally_hyperbolic_spacetime}) $\spacetime$. 
    By definition, the principal symbol is given by the spacetime metric $\fg^{-1}$ on the cotangent bundle $\coTansM$ times the identity endomorphism. 
    Thus the condition for being a real principal type differential operator on $\sE$ is satisfied trivially because the scalar wave operator (see Example~\ref{exm: covariant_Klein_Gordon_op}) on $\spacetime$ is of real principal type. 
    The characteristic set of $\square$ is the (co)lightcone bundle (see Section~\ref{sec: causal_structure_Lorentzian_mf}) $\Char \square = \coLightBun \to \sM$ over $\sM$ and the integral curves of the vector field $X_{\fg}$ generated by the Hamiltonian $\fg^{-1}$ are the geodesic strips in the punctured cotangent bundle $\coTansM$ whose projections on $\sM$ are the geodesic curves. 
\end{example}
%
%
%
%
%
%
%
%
%
%
\subsection{Lagrangian distributions}
\label{sec: Lagrangian_distribution}
The Lagrangian distribution on a vector bundle, its principal symbol and their composition are some of the indispensable concepts involved in this thesis.  
Definition~\ref{def: Lagrangian_distribution_mf} generalises to half-density-valued sections of a vector bundle in a straightforward fashion yet presented for the sake of completeness.  

%
%
%
\begin{definition} \label{def: Lagrangian_distribution}
	Let $\sE, \halfDen, \dotCoTanM \to M$ be a vector bundle, the half-density bundle, and the punctured cotangent bundle over a $d$-dimensional manifold $M$, respectively, and $\varLambda \subset \dotCoTanM$ a closed conic Lagrangian submanifold (see Definition~\ref{def: Lagrangian_submf}).  
    Then the \textbf{space} $I^{m} (M, \varLambda; \sE \otimes \halfDen)$ of \textbf{half-density-valued Lagrangian distributional sections} of $\sE$, of \textbf{order} (at most) $m \in \R$, is defined as the set of all half-density-valued distributions $u \in \cD' (M; \sE \otimes \halfDen)$ such that $P_{1} \ldots P_{N} u \in {}^{\infty} H_{\mathrm{loc}}^{-m - d/4} (M; \sE \otimes \halfDen)$ for all $N \in \NO$ and for all properly supported $P_{i} \in \PsiDO{1}{M; \sE \otimes \halfDen}$ having vanishing principal symbols on $\varLambda$, where ${}^{\infty} H_{\mathrm{loc}}^{-m - d/4} (M; \sE \otimes \halfDen)$ denotes the Besov space (Section~\ref{sec: PsiDO})~\cite[Def. 25.1.1]{Hoermander_Springer_2009}.  
\end{definition}
%
%
%

Let $\sF, \halfDenN \to N$ be a vector bundle and the half-density bundle over a manifold $N$. 
On top of the replacements mentioned after Definition~\ref{def: Lagrangian_distribution_mf}, we are compelled to substitute $\sE \otimes \halfDenM \to M$ by $(\sE \boxtimes \sF^{*}) \otimes \halfDenMN \to M \times N$ in order to define a Fourier integral operator from $\sF \otimes \halfDenN$ to $\sE \otimes \halfDenM$. 
One makes use of the natural bundle isomorphism 
\begin{equation}
	\Hom{\sF, \sE} \cong \sE \boxtimes \sF^{*} \to M \times N,   
\end{equation}
where the homomorphism bundle $\Hom{\sF, \sE} \to M \times N$ is the vector bundle whose fibres $\Hom{\sF, \sE}_{(x, y)} = \Hom{\sF_{y}, \sE_{x}}$ over $(x, y) \in M \times N$ are $\C$-linear maps from the vector space $\sF_{y}$ to the vector space $\sE_{x}$. 

%
%
%
\begin{definition} \label{def: FIO}
	Let $\sE, \halfDenM, \dotCoTanM \to M$ (resp. $\sF, \halfDenN, \dotCoTanN \to N$) be a vector bundle, the half-density bundle, and the punctured cotangent bundle over a $\dM$ (resp. $\dN$)-dimensional manifold $M$ (resp. $N$) and $C$ a homogeneous canonical relation (see Definition~\ref{def: canonical_relation}) from $\dotCoTanN$ to $\dotCoTanM$ of $M$, which is closed in $\dotCoTanMN$.  
	A linear continuous operator~\cite[Def. 25.2.1]{Hoermander_Springer_2009}   
	\begin{equation}
		A : \comSecNF \to \dualComSecME, ~v \mapsto 
		(A v) (x) := \int_{N} \fA (x, y) \, v (y) 
	\end{equation}
    whose Schwartz kernel $\fA$ is an element of the Lagrangian distribution $I^{m} \big( M \times N, C'; \mathrm{Hom}$ $(\sF, \sE) \otimes \halfDenMN$ is called a \textbf{Fourier integral operator} from $\sF$ to $\sE$, of \textbf{order} (at most) $m \in \R$ associated with the \textbf{canonical relation} $C$. 
	We will denote the set of all such operators by $\FIOFE$.
\end{definition}
%
%
%

The localisation (Lemma~\ref{lem: localisation_Lagrangian_dist_mf}) of Lagrangian distributions on a manifold holds true in the bundle scenario as well~\cite[Lem. 25.1.2]{Hoermander_Springer_2009}. 
Following the demonstration in Section~\ref{sec: integral_op}, it is evident that $\fA \in \LagrangianDistHomFE$ can be written as a matrix of scalar Lagrangian distributions on respective Euclidean spaces   
\begin{equation}
    \fA \leftrightarrow \big( \fA^{r}_{k} \big)_{\rk \sE \times \rk \sF} \in \Big( I^{m} \big( \kappa_{\alpha} (U_{\alpha}) \times \rho_{\beta} (V_{\beta}); C'_{\varphi} \big) \Big)_{\rk \sE \times \rk \sF}
\end{equation}
where $C'_{\varphi} := C'_{\alpha, \beta} :=$~\eqref{eq: Euclidean_rep_canonical_relation} is given by means of a phase function $\varphi := \varphi_{\alpha, \beta}$ with excess $e := e_{\alpha, \beta}$ as elucidated in Section~\ref{sec: FIO_mf}. 
$\fA^{r}_{k}$ is of the form~\eqref{eq: Lagrangian_distribution_mf_local}~\cite[Prop. 25.1.5']{Hoermander_Springer_2009} 
\begin{equation} \label{eq: Lagrangian_distribution_local}
	\fA^{r}_{k} (\bx, \by) = (2 \pi)^{- (\dM + \dN + 2 n - 2 e) / 4} \int_{\Rn} \re^{\ri \varphi (\bx, \by; \theta)} \fa^{r}_{k} (\bx, \by; \theta) \, \rd \theta 
    \mod C^{\infty} \big( \kappa_{\alpha} (U_{\alpha}) \times \rho_{\beta} (V_{\beta}) \big) 
\end{equation}
where $\fa^{r}_{k} \in S^{m + (\dM + \dN - 2 n - 2 e) / 4} (\R^{\dM} \times \R^{\dN} \times \R^{n})$ having support in the interior of a sufficiently small conic neighbourhood of the fibre-critical manifold~\eqref{eq: fibre_critical_mf_Euclidean} of $\varphi$ contained in the domain of definition of $\varphi$. 

%
%
%
\begin{definition} \label{def: symbol_FIO}
    As in the terminologies of Definition~\ref{def: FIO}, let $\widetilde{\mathrm{Hom}} (\sF, \sE), \halfDenC, \M \to C$ be the pullback of $\Hom{\sF, \sE}$ to $\dotCoTanM \times \dotCoTanN$ restricted to $C$ (see~\eqref{eq: symbol_HomFE_canonical_relation}), the half-density bundle (see Appendix~\ref{sec: volume_canonical_relation}) and the Keller-Maslov bundle (see Definition~\ref{def: Keller_Maslov_bundle}) over $C$, respectively, and $S^{m + (\dM + \dN) / 4 - [1]} \big( C; \M \otimes \halfDenC \otimes \widetilde{\mathrm{Hom}} (\sF, \sE) \big)$ the space (see~\eqref{eq: symbol_canonical_relation_Maslov_density_function}) of $\Maslov \otimes \halfDenC \otimes \widetilde{\mathrm{Hom}} (\sF, \sE)$-valued symbols on $C$. 
    Then the \textbf{principal symbol} of a Lagrangian distribution is defined by the isomorphism~\cite[Prop. 25.1.5']{Hoermander_Springer_2009})     
	\begin{equation}
		\symb{} : I^{m - [1]} \big( M \times N, C'; \Hom{\sF, \sE} \otimes \halfDenMN \big) \to S^{m + \frac{\dM + \dN}{4} - [1]} \big( C; \M \otimes \halfDenC \otimes \widetilde{\mathrm{Hom}} (\sF, \sE) \big),   
	\end{equation}
    where $I^{m-[1]} := I^{m} / I^{m-1}$. 
    The preceding map is locally given by
    \begin{eqnarray}
        [\fA] \mapsto 
        \symb{[\fA]} \xxiyeta 
        & := & 
        \left( \int_{\kC} \rd \theta'' a^{r}_{k} (\bx, \by; \theta', \theta'') \right)_{\rk \sE \times \rk \sF} \m \xxiyeta \otimes  
        \nonumber \\ 
        && 
        \sqrt{|\dVol_{\ms C} \xxiyeta|} \mod S^{m + \frac{\dM + \dN}{4} - 1} (\cdot),   
    \end{eqnarray}
    where $a^{r}_{k}$ is the top-order homogeneous term of $\fa^{r}_{k}$ in~\eqref{eq: Lagrangian_distribution_local} and all other symbols are as in Definition~\ref{def: symbol_FIO_Euclidean}.  
\end{definition}
%
%
%
%
%
%
%
%
%
%
\subsection{Adjoint of a Fourier integral operator} 
\label{sec: adjoint}
To describe the adjoint of a Fourier integral operator, let us recall some standard notions from linear algebra. 
Let $\sU$ and $\sV$ be $\C$-vector spaces. 
Then the set $\Hom{\sU, \sV}$ of $\C$-linear maps $T : \sU \to \sV$ is a vector space and $\End{\sV} := \Hom{\sV, \sV}$ is an unital algebra over $\C$ under the composition of endomorphisms.  
If $\sV^{*} := \Hom{\sV, \C}$ is the algebraic dual of $\sV$ then the set of algebraic dual maps $T^{*} \in \Hom{\sV^{*}, \sU^{*}}$ is defined by $(T^{*} \psi) (u) = \psi (Tu)$ for all $\psi \in \sV^{*}, u \in \sU$.  
The conjugate of $\sV$ is the complex vector space $\bar{\sV} := \{ \bar{v} \,|\, v \in \sV \}$ with an antilinear isomorphism $\cC : \sV \ni v \mapsto \bar{v} \in \bar{\sV}$ where, for any given $v \in \sV$, 
$\cV^{*} \ni \psi \mapsto \bar{v} (\psi) := \overline{\psi (v)} \in \C$. 
In other words, $\sV$ is identical to $\bar{\sV} := \overline{\mathrm{Iso}} (\sV^{*}, \C)$ as a \textit{set but not as a vector space} because the latter is equipped with the scalar multiplication $\lambdaup \bullet v := \bar{\lambdaup} v$ in contrast to $\lambdaup \cdot v := \lambdaup v$ for the former for any $\lambdaup \in \C$. 
If $\sV$ is equipped with a bilinear form $\scalarProdOne{\cdot}{\cdot}$ then an intertwiner $\cB \in \Hom{\sV, \sV^{*}}$ is determined by $(\cB v) (u) := \scalarProdOne{v}{u}$ for any $u \in \sV$.
Analogously, for a given sesquilinear form $(\cdot|\cdot)$ on $\sV$, another intertwiner $\cA := \bar{\cB} \circ \cC : \sV \to \bar{\sV}^{*}$ is characterised by $(\cA v) (\bar{u}) := (u|v)$. 
One has the identifications $(\bar{\sV})^{*} = \overline{\sV^{*}} = \bar{\sV}^{*}$ and $\bar{\bar{\sV}} = \sV$. 
\newline 

Since $\sE$ is a vector bundle, the above description applies fibrewise $\sE_{x}$ and one has the (anti)dual $(\bar{\sE}^{*} := \sqcup_{x \in M} \bar{\sE}_{x}^{*} \to M) \sE^{*} := \sqcup_{x \in M} \sE_{x}^{*} \to M$ bundles over a manifold $M$  
(see e.g.~\cite[p. 22]{Guenther_AP_1988}). 
A Fourier integral operator (cf~\eqref{eq: FIO_Euclidean_CCom_C}) 
\begin{equation}
	A : \comSecNF \to \secME 
\end{equation}
associated with a canonical relation $C$, has a unique formal adjoint in the following sense.  
If $C^{-1}$ denotes the inverse relation of $C$ obtained by interchanging $\coTanM$ and $\coTanN$ then there exists a unique Fourier integral operator 
\begin{equation}
	\bar{A}^{*} : C_{\mathrm{c}}^{\infty} (M; \bar{\sE}^{*} \otimes \halfDenM) \to C^{\infty}  (N; \bar{\sF}^{*} \otimes \halfDenN)
\end{equation}
associated with $C^{-1}$ such that~\cite[Thm. 25.2.2]{Hoermander_Springer_2009} 
\begin{equation}
	\int_{N} (\bar{A}^{*} \phi) (y) \, v (y) := \int_{M} (Av) (x) \, \phi (x)  
	\Leftrightarrow 
    \bar{\fA}^{*} (v \boxtimes \phi) := \fA (\phi \boxtimes v) 
\end{equation}
for any $\phi \in C_{\mathrm{c}}^{\infty} (M; \bar{\sE}^{*} \otimes \halfDenM)$ and any $v \in C_{\mathrm{c}}^{\infty} (N; \sF \otimes \halfDenN)$. 
Its Schwartz kernel and principal symbol are given by 
\begin{subequations}
	\begin{eqnarray}
		&& 
        \bar{\fA}^{*} \in I^{m} \big( N \times M, C^{-1 \prime}; \Hom{\bar{\sE}^{*}, \bar{\sF}^{*}} \otimes \halfDen (N \times M) \big), 
        \\ 
        && 
        \symb{[\bar{\fA}^{*}]} = \mathsf{s}^{*} (\overline{\symb{[\fA]}}^{*}) \in S^{m + \frac{\dM + \dN}{4} - [1]} \big( C^{-1}; \halfDen C^{-1} \otimes \Maslov C^{-1} \otimes \widetilde{\mathrm{Hom}} (\bar{\sE}^{*}, \bar{\sF}^{*}) \big), \qquad 
	\end{eqnarray}
\end{subequations}
where $\mathsf{s} : N \times M \to M \times N$ is the interchanging (Schwartz) map and $\Maslov C^{-1} \to C^{-1}$ is the Keller-Maslov bundle over $C^{-1}$. 
\newline 

Locally $\bar{\fA}^{*} (y, x)$ is given by  
\begin{subequations}
	\begin{eqnarray} 
        && 
        \fA^{k}_{r} (\by, \bx) \equiv (2 \pi)^{- (\dM + \dN + 2 n - 2 e) / 4} \int_{\Rn} \re^{- \ri \varphi (\by, \bx; \btheta)} \bar{\fa}^{k}_{r} (\by, \bx; \btheta) \, \rd \btheta,  
        \\ 
        && 
        \symb{\fA^{k}_{r}} (\by, \bfeta; \bx, \bxi)  
        = 
        \sqrt{|\rd \bfeta| \, |\rd \bxi|} 
        \int_{\fs (\kC_{\xi, \eta})} \bar{\fa}^{k}_{r} (\by, \bx; \btheta', \btheta'') \, \dfrac{\re^{- \nicefrac{\mathrm{i} \pi}{4}  \sgn (\Hess_{\by, \bx; \btheta'} \varphi)}}{\sqrt{|\det (\Hess_{\by, \bx; \btheta'} \varphi)|}} \rd \btheta'' \qquad 
\end{eqnarray}
\end{subequations}
whenever $\fA$ is given by~\eqref{eq: Lagrangian_distribution_local} and $\fs (\kC_{\xi, \eta})$ is defined by interchanging $\bx$ and $\by$ in~\eqref{eq: def_fibre_Lagrangian_fibration_symbol_FIO_Euclidean}. 
%
%
%
%
%
%
%
%
%
%
\subsection{Algebra of Fourier integral operators}
\label{sec: algebra_FIO}
A necessary assumption for the product (composition) of two Fourier integral operators to be well-defined is that the first operator must be properly supported. 
Then the defined composition may, however, still fail to be a Fourier integral operator.
For instance, the composition of two canonical relations does not necessarily have the required properties to define another Fourier integral operator.  
To ensure a well-defined product that is again a Fourier integral operator, we are obliged to assume that their Schwartz kernels (Lagrangian distributions) are \textit{properly supported} and the composition of canonical relations is \textit{clean}, \textit{proper} (see Definition~\ref{def: clean_composition}), and \textit{connected}~\cite[Sec. 4.1]{Weinstein_Nice_1975},~\cite[Thm. 5.4]{Duistermaat_InventMath_1975}.   
Given a vector bundle $\sG \to O$ over a manifold $O$ besides $\sE \to M$ and $\sF \to N$, if   
$\fA \in I^{m} \big( M \times O, C'; \Hom{\sG, \sE} \otimes \halfDen (M \times O) \big)$ 
and 
$\fB \in I^{m'} \big( O \times N, \varLambda'; \Hom{\sF, \sG} \otimes \halfDen (O \times N) \big)$ 
with the required restrictions, then~\cite[Thm. 25.2.3]{Hoermander_Springer_2009} 
\begin{equation} \label{eq: product_Lagrangian_dist}
	\fA \fB := \fA \circ \fB \in I^{m + m' + e/2} \big( M \times N, (C \circ \varLambda)'; \Hom{\sF, \sE} \otimes \halfDenMN \big), 
\end{equation}
where $e$ is the excess (see Definition~\ref{def: clean_composition} and~\eqref{eq: def_excess}) of the clean composition $C \circ_{e} \varLambda$ and its principal symbol is given by 
\begin{equation} \label{eq: symbol_product_Lagrangian_dist}
    \symb{\fA \fB} = \symb{A} \diamond_{e} \symb{B}.  
\end{equation} 
The precise meaning of the right hand side of this equation demands a few new terminologies which are detailed in Appendix~\ref{sec: composition_halfdensity_Maslov_vector_bundle_canonical_relation}. 
\newline 

The space of Fourier integral operators is not an algebra unless the canonical relation is symmetric and transitive~\cite[Ex. 1 and comment on it in the succeeding page]{Guillemin_JFA_1993}. 
In case the composition of the canonical relation $C$ by itself is clean with excess zero, proper, and connected, the space of Fourier integral operators $\big( \mathrm{FIO}^\bullet (\sE \otimes \halfDen \to M; C'), \circ \big)$ is an \textit{associative graded algebra} over the field $\C$.
The algebra~\eqref{eq: symbol_product_Lagrangian_dist} of principal symbols constitute an algebra which is commutative only the scalar case but non-commutative in general. 
We remark that the product formula~\eqref{eq: symbol_product_Lagrangian_dist} becomes simpler (see~\eqref{eq: product_symbol_FIO_symplecto}) when both $C$ and $\varLambda$ are local canonical graphs and then the product is always defined as detailed later in Section~\ref{sec: FIO_symplecto}. 

%
%
%
\begin{example} \label{exm: restriction_op_FIO}
	Let $\sE, \halfDenM \to M$ be a vector bundle and the half-density bundle over a manifold $M$, respectively, and $\iota_{\ms \varSigma} : \varSigma \hookrightarrow M$ an immersed submanifold.  
    Given a fibrewise isomorphism $\hat{\iota}_{\ms \varSigma} : \sE_{ \varSigma_{x'}} \cong \sE_{\iota_{\ms \varSigma} (x')}$, the \textbf{restriction map}  
    \begin{equation}
        \iota_{\ms \varSigma}^{*} : \comSecME \to C^{\infty} (\varSigma; \sE_{\varSigma} \otimes \halfDenM_{\varSigma})
    \end{equation}
    is defined by the pullback of the morphism $(\iota_{\ms \varSigma}, \hat{\iota}_{\ms \varSigma})$. 
    This boils down to Example~\ref{exm: restriction_op_FIO_mf} after trivialising $\sE$ and thus its Schwartz kernel $\fR_{\varSigma}$ is a Lagrangian distribution of order $\codim \varSigma/4$ associated with the canonical relation $\varLambda_{\varSigma}'$ defined in~\eqref{eq: def_canonical_relation_restriction_op_mf}.   
    As commented after Example~\ref{exm: restriction_op_FIO_mf}, one discerns that $\varLambda_{\varSigma} \subset \coTanM_{\varSigma} \times \dotCoTanM$ is \textit{not} a homogeneous canonical relation because it is not guaranteed that restriction of all elements in $\dotCoTanM$ to $\varSigma$ is non-zero. 
    Geometrically speaking, the complication arises due to the elements of the conormal bundle (see Example~\ref{exm: conormal_bundle_Lagrangian_submf}) $\varSigma^{\perp *} \subset \coTanM$. 
    To circumvent this difficulty, one temporarily introduces a cutoff function $\chi$ on $\varSigma^{\perp *}$ and sets~\cite[Sec. 5.2]{Toth_GFA_2013}   
    \begin{equation}
        \fR_{\chi} (\bx', \by) := \int_{\Rd} \re^{\ri \varphi_{\ms \varSigma}} \big( 1 - \chi \yeta \big) \frac{\rd \bfeta}{(2 \pi)^{d}}, 
        \quad 
        \varphi_{\ms \varSigma} := (\iota_{\ms \varSigma} (\bx') - \by) \cdot \bfeta
    \end{equation}   
    so that such elements do \textit{not} occur in the support of $1 - \chi$. 
    Then the respective operator $\iota_{\chi}^{*}$ is a \textit{homogeneous} Fourier integral operator of order $\codim \varSigma/4$ associated with the homogeneous canonical relation $\varLambda_{\chi} \subset \dotCoTanM_{\varSigma} \times \dotCoTanM$. 
    Putting altogether: 
    \begin{subequations} \label{eq: restriction_kernel}
        \begin{eqnarray}
		&& 
        \fR_{\chi} \in I^{\codim \varSigma / 4} \big( \varSigma \times M, \varLambda_{\chi}'; \Hom{\sE, \sE_{\varSigma}} \otimes \halfDen (M_{\varSigma} \times M) \big), 
		\label{eq: restriction_kernel_Lagrangian_dist}
		\\ 
		&& 
        \varLambda_{\chi}' := \{ (x', \xi'; x, \xi) \in \dotCoTanM_{\varSigma} \times \dotCoTanM \,|\, x = \iota_{\ms \varSigma} (x'), \xi' = \xi \upharpoonright \tangent_{x'} M_{\varSigma} \},  
		\label{eq: def_canonical_relation_restriction} 
		\\ 
		&& 
        \symb{\fR_{\chi}} |_{\supp(1 - \chi)} (x', \xi'; x, -\xi) 
        := 
        \frac{ \sqrt{|\dVol_{\ms \varLambda} (x', \xi'; x, -\xi)|} }{(2 \pi)^{ \codim \varSigma / 4}} \one_{\Hom{\sE, \sE_{\varSigma}} \otimes \varOmega^{\ms 1/2}(M_{\varSigma} \times M) } \otimes \bbl (x', \xi'; x, -\xi),   
		\label{eq: symbol_restriction}
        \nonumber \\
        \end{eqnarray}
    \end{subequations} 
    where $|\dVol_{\ms \varLambda}|$ is the density on $\varLambda_{\chi}$ and $\bbl$ is a section of the Keller-Maslov bundle $\bbL \to \varLambda_{\chi}$ as in Example~\ref{exm: restriction_op_FIO_mf}. 
    \newline 

    Equivalently, one can say that $\iota_{\ms \varSigma}^{*}$ is not a homogeneous Fourier integral operator due to plausible occurrence of elements of $\varSigma^{\perp *}$ in its canonical relation and the cutoff above can be emulated by composing with a $P \in \PsiDO{0}{M; \sE \otimes \halfDenM}$ with $\ES{P} \cap \dot{\varSigma}^{\perp *} = \emptyset$. 
    The principal symbol of $\fR_{\varSigma}$ is then $(2 \pi)^{- \codim \varSigma / 4} \one |\dVol_{\ms \varLambda}|^{1/2} \otimes \bbl$ over $\WF{P}$ 
    (see e.g.~\cite[Lem. 8.3]{Strohmaier_AdvMath_2021} for the scalar version). 
    \newline 

    Any distribution on $\sE$ can be restricted to $\sE_{\varSigma}$ when its wavefront set is disjoint with $\varSigma^{\perp *}$. 
    In particular, if $\Lambda$ is such that the wavefront set of every $u \in I^{m} (M, \Lambda; \sE \otimes \halfDenM)$ is disjoint from the conormal bundle $\varSigma^{\perp *}$ of $\varSigma \subset \sM$ then the restriction operator can be extended to a sequentially continuous linear operator 
    \begin{subequations}
    	\begin{eqnarray}
            && \iota_{\chi}^{*} : I^{m} (M, \Lambda; \sE \otimes \halfDenM) \to I^{m + \codim \varSigma / 4} (\varSigma, \Lambda |_{\tangent M_{\varSigma}}; \sE_{\varSigma} \otimes \halfDen \varSigma), ~ u \mapsto \iota_{\chi}^{*} u, 
            \nonumber \\ 
            && \symb{\iota_{\chi}^{*} u} (x', \xi') \asymp (2 \pi)^{- \codim \varSigma / 4} \symb{u} (x', \xi'), 
    	\end{eqnarray}
    \end{subequations}
    where the last equation is in the sense of modulo Keller-Maslov contribution. 
\end{example}
%
%
%
%
%
%
%
%
%
%
\subsubsection{Products of operators with vanishing principal symbol}
By far, we have learned to compute the principal symbol of the product of Lagrangian distributions given those of individual distributions. 
An innate question in this context is that what if any of the principal symbols  vanishes identically. 
We are going to address a special case of this situation: product between a pseudodifferential operator and a Lagrangian distribution whenever the principal symbol of the former vanishes. 
This is motivated by its direct applicability at some later parts of this thesis. 
It turns out that the Lie derivative plays an essential role in this regard, so we spell out a few rudimentary formulae. 

%
%
%
\begin{center}
	\begin{tikzpicture}
		\node (a) at (0, 0) {$\coTanM$};
		\node (b) at (3, 0) {$M$};
		\draw[->] (a) -- (b); 
		\node (c) at (0, 1.5) {$\pi^{*} \sE$}; 
		\node (d) at (3, 1.5) {$\sE$}; 
		\draw[->] (c) -- (d); 
		\draw[->] (c) -- (a); 
		\draw[->] (d) -- (b); 
		\node[above] at (1.5, 0) {$\pi$}; 
		\node[above] at (1.5, 1.5) {$\cong$}; 
        \node (e) at (6, 0) {$C^{\infty} (M; \sE \otimes \halfDen)$};
		\node (f) at (12, 0) {$C^{\infty} (M; \sE \otimes \halfDen)$};
        \node (g) at (6, 1.5) {$C^{\infty} (\coTanM; \pi^{*} \sE)$}; 
		\node (h) at (12, 1.5) {$C^{\infty} (\coTanM; \pi^{*} \sE)$}; 
		\draw[->] (e) -- (f); 
		\draw[->] (e) -- (g); 
		\draw[->] (h) -- (f); 
        \draw[->] (g) -- (h);
        \node[above] at (9, 0) {$\connectionE_{X}$}; 
        \node[above] at (9, 1.5) {$\nabla^{\ms \Hom{\pi^{*} \sE, \pi^{*} \sE}}_{X_{p}}$}; 
        \node[left] at (6, 0.75) {$\pi^{*}$}; 
        \node[right] at (12, 0.75) {$(\pi^{*})^{-1}$}; 
	\end{tikzpicture}
	\captionof{figure}[Pullback of $\sE \to M$ via $\pi: \coTanM \to M$]{
        Pullback of a vector bundle $\sE \to M$ via the canonical projection $\pi: \coTanM \to M$ and the $P$-compatible connection on $\sE$.
        }
	\label{fig: pullback_vector_bundle}
\end{center}
%
%
%

Let $(\sE \to M, \connectionE_{X})$ be a vector bundle over a manifold $M$, equipped with a covariant derivative $\connectionE_{X}$ with respect to \textit{some} vector field $X$ on $M$ and let $\varXi_{s} : M \to M$ be the flow of $X$. 
This covariant derivative induces a parallel transport map 
\begin{equation}
    \hat{\varXi}_{s} (x) : \sE_{x} \to \sE_{\varXi_{s} (x)}
\end{equation}
along the integral curve $c_{x} (s) = \varXi_{s} (x)$ of $X$ passing through $x \in M$. Then the \textbf{Lie derivative} $\pounds_{X}$ along $X$ is defined by 
\begin{equation} \label{eq: def_Lie_derivative}
    \pounds_{X} u 
    := \lim_{s \to 0} \frac{\hat{\varXi}_{-s} \Big( u \big( \varXi_{s} (x) \big) \Big) - u (x)}{s} 
    := \frac{\rd}{\rd s} \Big|_{s=0} \varXi_{s}^{*} u,  
\end{equation}
where the pullback map on $\sE$, induced by $\varXi_{s}$ and the linear isomorphism $\hat{\varXi}_{s}$, is given by 
\begin{equation} \label{eq: def_pullback_VB_parallel_transport}
    \varXi_{s}^{*} : C^{\infty} (M; \sE) \to C^{\infty} (M; \sE), ~ u \mapsto \varXi_{s}^{*} u := \hat{\varXi}_{-s} \circ u \circ \varXi_{s}. 
\end{equation}
For convenience, we set  
\begin{equation} 
    \LieDeri_{X} := - \ri \pounds_{X} : C^{\infty} (M; \sE) \to C^{\infty} (M; \sE)
\end{equation}
and call this first-order differential, the \textit{Lie derivative} as well. 
Locally, $X = X^{i} \nicefrac{\partial}{\partial x^{i}}$ and $C^{\infty} (M; \sE \otimes \halfDen) \ni u = u^{r} \rE_{r} \otimes \sqrt{|\rd x|}$. 
Thus the Leibniz rule entails 
(see e.g.~\cite[p. 130]{Abraham_AMS_1978},~\cite[$(25.2.11)$]{Hoermander_Springer_2009},~\cite[$(2.60)$ $($p. 377$)$]{Treves_Plenum_1980}) 
\begin{equation} \label{eq: Lie_derivative_density_coordinate}
    \pounds_{X} u 
    = \bigg( X^{i} \parDeri{x^{i}}{u^{r}} + \frac{1}{2} \Div (X) \, u^{r} \bigg) \rE_{r} \otimes \sqrt{|\rd x|} + u^{r} \pounds_{X} (\rE_{r}) \otimes \sqrt{|\rd x|}. 
\end{equation}

In what follows, we consider $P \in \PsiDOME$ having a scalar principal symbol $p$ and subprincipal symbol $\subSymb{P}$. 
Let $X_{p}$ be the Hamiltonian vector field (see~\eqref{eq: def_HVF}) generated by $p$ and $\gamma$ its integral curve (cf. Definition~\ref{def: bicharacteristic_strip}) on the cotangent bundle $\pi : \coTanM \to M$. 
These data enables us to define a natural covariant derivative in a geometric fashion as detailed in Section~\ref{sec: connection_form_subprincipal_symbol}. 
The so-called ``$P$-compatible covariant derivative'' (see Definition~\ref{def: P_compatible_connection}) $\nabla^{\ms \Hom{\pi^{*} \sE, \pi^{*} \sE}}_{X_{p}}$ is a covariant derivative on the homomorphism bundle $\Hom{\pi^{*} \sE, \pi^{*} \sE} \to \coTanM$ over $\coTanM$ along $\gamma$. 
The projection $c$ of $\gamma$ on $M$ is the integral curve of the vector field $X := \rd \pi (X_{p})$ on $M$ induced by $X_{p}$, where $\rd \pi : \tangent (\coTanM) \to \tangent M$ is the tangent map of $\pi$.  
Thus we have a covariant derivative $\connectionE_{X}$ on $\sE$ along $c$, as depicted in Figure~\ref{fig: pullback_vector_bundle}. 
\newline 

The following result is closely related to~\cite[Thm. 5.9]{Sandoval_CPDE_1999}.   

%
%
%
\begin{theorem} \label{thm: Hoermander_thm_25_2_4}
    As in the terminologies of Definition~\ref{def: FIO}, let $\fA \in I^{m} \big( M \times N, C'; \mathrm{Hom} (\sF,$ $ \sE) \otimes \halfDenMN \big)$ and $a$ the principal symbol of $\fA$.  
	Suppose that $P \in \PsiDO{m'}{M; \sE \otimes \halfDenM}$ is properly supported with a scalar principal symbol $p$ and a subprincipal symbol $\subSymb{P}$, and that $p$ vanishes identically on the projection of $C$ in $\dotCoTanM$. 
	Then, $P \fA \in I^{m + m' - 1} \big( M \times N, C'; \Hom{\sF, \sE} \otimes \halfDenMN \big)$ and its principal symbol is  
	\begin{equation}
		\symb{P \fA} = \left( \LieDeri_{X_{p}} + ~ \subSymb{P} \right) a, 
	\end{equation}
	where $\LieDeri_{X_{p}}$ is $- \ri$ times the Lie derivative along the Hamiltonian vector field $X_{p}$ of $p$, lifted to a function on $\dotCoTanM \times \dotCoTanN$, so $X_{p}$ is tangential to $C$. 
\end{theorem}
%
%
%
\begin{proof}
	Since the principal symbol is locally defined, our strategy is to make use of the partition of unity to boil down the statements in the level of a matrix of scalar operators as demonstrated in Sections~\ref{sec: PsiDO},~\ref{sec: Lagrangian_distribution} and then employ the scalar-version of this 
	statement~\cite[Thm. 5.3.1]{Duistermaat_ActaMath_1972} 
	(see also, e.g.~\cite[Thm. 25.2.4]{Hoermander_Springer_2009},~\cite[pp. 451-454]{Treves_Plenum_1980}). 
	The differences in the proof between the scalar and the bundle versions are essentially bookkeeping, yet we provide the details for completeness.    
    Locally, the transition functions of $\Maslov$ are constant (see Section~\ref{sec: Keller_Maslov_bundle}) so the Maslov factor can be ignored while computing the Lie derivative. 
    Moreover, $a$ is a matrix-valued density (Definition~\ref{def: symbol_FIO}) on $C'$ and thus the last term in~\eqref{eq: Lie_derivative_density_coordinate} does not contribute. 
    It is always possible to choose a local coordinate $(x^{i}, \xi_{i}; y^{j}, \eta_{j})$ around $\xxiyeta \in \dotCoTanM \times \dotCoTanN$ such that (see~\eqref{eq: def_canonical_nondegenerate_phase_function_canonical_relation}-\eqref{eq: canonical_relation_clean_phase_function}) $x^{i} \xi_{i} + y^{j} \eta_{j} - H$ is a non-degenerate phase function for $C$, where $H$ is a smooth positively homogeneous $\R$-valued function of degree $1$ on an open conic neighbourhood of $(\xi, \eta) \in \dotCoTan_{x} M \times \dotCoTan_{y} N$. 
    Then (cf.~\eqref{eq: Lagrangian_distribution_local} and see e.g.~\cite[Lem. 25.2.5]{Hoermander_Springer_2009}),  
	\begin{equation}
	   	\fA^{l}_{k} (\by, \bz) 
	    = 
        (2 \pi)^{- 3 (\dM + \dN) / 4} \int_{\R^{\dM + \dN}}  
	    \re^{\ri (\bfeta \cdot \by + \bzeta \cdot \bz - H (\bfeta, \bzeta))} \, \fa^{l}_{k} (\bfeta, \bzeta) \, \rd \bfeta \, \rd \bzeta, 
	\end{equation}
    where $\fa^{l}_{k} \in S^{m - (\dM + \dN) / 4} (\R^{\dM + \dN})$ having support in a conic neighbourhood of $(\by = \partial H / \partial \bfeta, \bz = \partial H / \partial \bzeta)$. 
    For any $v \in \comSecNF$, by definition of Fourier transform:  
    \begin{eqnarray}
    	\big( \Fourier (A^{l}_{k} v^{k}) \big) (\bxi) 
        & = & 
        \int_{\R^{\dM}} \rd \by \, \re^{- \ri \by \cdot \bxi} (A^{l}_{k} v^{k}) (\by) 
    	\nonumber\\ 
    	& = & 
        \int_{\R^{\dM}} \rd \by \, \re^{- \ri \by \cdot \bxi} \int_{\R^{\dN}} \fA^{l}_{k}  (\by, \bz) \, v^{k} (\bz) \, \rd \bz 
    	\nonumber\\ 
    	& = & 
        (2 \pi)^{- 3 (\dM + \dN) / 4} \int_{\R^{\dM}} \rd \by \int_{\R^{\dN}} \rd \bz \int_{\R^{\dM + \dN}} \rd \bfeta \, \rd \bzeta 
    	\nonumber\\ 
    	&& \re^{\ri \by \cdot (\bfeta - \bxi) + \ri \bz \cdot \bzeta - \ri H (\bfeta, \bzeta)} 
    	\, \fa^{l}_{k} (\bfeta, \bzeta) \, v^{k} (\bz) 
    	\nonumber\\ 
    	& = & 
        (2 \pi)^{\dM - \frac{3}{4} (\dM + \dN)} \int_{\R^{\dN}} \re^{- \ri H (\bxi, \bzeta)}  
    	\fa^{l}_{k} (\bxi, \bzeta) \, (\Fourier v^{k}) (- \bzeta) \, \rd \bzeta, 
    \end{eqnarray}  
    where we have made use of the Fubini's theorem for oscillatory integrals~\cite[$(1.2.4)$]{Hoermander_ActaMath_1971}. 
    Since $P$ is properly supported, by an application of Theorem~\ref{thm: characterisation_psPsiDO_Euclidean} we obtain  
    \begin{equation}
    	(P^{r}_{l} A^{l}_{k} v^{k}) (\bx)
    	=
    	\frac{1}{(2 \pi)^{\frac{3}{4} (\dM + \dN)}} 
        \int_{\R^{\dM}} \rd \bxi \int_{\R^{\dN}} \rd \bzeta \, \re^{\ri \bx \cdot \bxi - \ri H (\bxi, \bzeta)} 
        (\totSymb{P})^{r}_{l} (\bx, \bxi) \, \fa^{l}_{k} (\bxi, \bzeta) (\Fourier v^{k}) (- \bzeta),
    \end{equation}
    where $(\totSymb{P})^{r}_{l}$ is the total symbol of $P$ in the chosen coordinate-charts and bundle-charts. 
    The preceding equation entails that the Schwartz kernel of $PA$ is 
    \begin{equation}
    	(P^{r}_{l} \fA^{l}_{k}) (\bx, \by) 
       	= 
        (2 \pi)^{- \frac{3}{4} (\dM + \dN)}
        \int_{\R^{\dM}} \rd \bxi \int_{\R^{\dN}} \rd \bfeta \, \re^{\ri (\bx \cdot \bxi + \by \cdot \bfeta - H)} (\totSymb{P})^{r}_{l} (\bx, \bxi) \, \fa^{l}_{k} (\bxi, \bfeta).  
    \end{equation}
    One can decompose $(\totSymb{P})^{r}_{l} = p^{r}_{l} + \tilde{\fp}^{r}_{l}$ as a combination of its principal part $p^{r}_{l}$ and the remainder terms $\tilde{\fp}^{r}_{l}$ (cf.  Definition~\ref{def: polyhomogeneous_symbol_Euclidean}). 
    We have assumed that $p = 0$ on the projection of $C$ on $\dotCoTanM$.
    By means of Taylor's series one can write 
    \begin{equation} \label{eq: FIO_II_5_3_4}
        p^{r}_{l} (\bx, \bxi) 
        = \left( \bx^{i} - \parDeri{\bxi_{i}}{H} \right) f^{r}_{l, i} (\bx, \bxi, \bfeta)    
        = \parDeri{\bxi_{i}}{(\bx \cdot \bxi + \by \cdot \bfeta - H)} f^{r}_{l, i} (\bx, \bxi, \bfeta), 
    \end{equation}
    where $f^{r}_{l, i}$ is homogeneous of degree $m'$ with respect to $(\bxi, \bfeta)$, given by the mean value theorem: 
    \begin{equation}
    f^{r}_{l,i} (\bx, \bxi, \bfeta) 
    = 
    \int_{0}^{1} \frac{\partial p^{r}_{l}}{\partial \bx^{i}} \Big( \lambdaup \, \bx + (1 - \lambdaup) \frac{\partial H}{\partial \bxi}, \bxi \Big) \rd \lambdaup.
    \end{equation}
    In particular, $f^{r}_{l,i} \upharpoonright C = (\partial_{\bx} p^{r}_{l}) (\partial_{\bxi} H, \bxi)$. 
    It is always possible to assume that $\fa$ vanishes in a neighbourhood of $0$ and then an integration by parts gives for $P \fA$:   
    \begin{eqnarray}
       	P^{r}_{l} \fA^{l}_{k}  
        & = &   
        (2 \pi)^{- \frac{3}{4} (\dM + \dN)}  
        \int_{\R^{\dM + \dN}} \re^{\ri (\bx \cdot \bxi + \by \cdot \bfeta - H)} \,  
        (\tilde{\fp}^{r}_{l} + p^{r}_{l}) \fa^{l}_{k} 
       	\, \rd \bxi \, \rd \bfeta
        \nonumber \\ 
       	& = &   
        (2 \pi)^{- \frac{3}{4} (\dM + \dN)}  
        \int_{\R^{\dM + \dN}} \re^{\ri (\bx \cdot \bxi + \by \cdot \bfeta - H)} \,  
       	\Big( \tilde{\fp}^{r}_{l} \fa^{l}_{k} 
       	- 
        \rD_{\bxi_{i}} (f^{r}_{l, i} \fa^{l}_{k}) \Big) \rd \bxi \, \rd \bfeta. 
    \end{eqnarray}
    Since $a^{l}_{k} = \fa^{l}_{k} (\bxi, \bfeta) \sqrt{| \rd \bxi| \, | \rd \bfeta |}$ is the principal symbol of $\fA$, hence 
    \begin{eqnarray}
        (\symb{P \fA})^{r}_{k}  
       	& = &  
       	\Big( \tilde{\fp}^{r}_{l} a^{l}_{k} 
        - \rD_{\bxi_{i}} (f^{r}_{l, i} a^{l}_{k}) \Big)_{\bx = \partial H / \partial \bxi}   
       	\sqrt{| \rd \bxi| \, | \rd \bfeta|} 
        \nonumber \\ 
        & = &  
       	\Big( \tilde{\fp}^{r}_{l} a^{l}_{k} 
        - \rD_{\bxi_{i}} (f^{r}_{l, i}) a^{l}_{k} 
        - f^{r}_{l, i} \rD_{\bxi_{i}} (a^{l}_{k}) \Big)_{\bx = \partial H / \partial \bxi}   
       	\sqrt{| \rd \bxi| \, | \rd \bfeta|} 
    \end{eqnarray}
    whenever $(\bxi, \bfeta)$ are taken as coordinates on $C$. 
    We will suppress the bundle indices from afterwords as there is no scope of confusion. 
    \newline 

    By definition, $X_{p} = \nicefrac{\partial p}{\partial \bxi_{i}} \nicefrac{\partial}{\partial \bx^{i}} - \nicefrac{\partial p}{\partial \bx^{i}} \nicefrac{\partial}{\partial \bxi_{i}}$. 
    So, $\Char{P}$ is a hypersurface near $\xxi$ and $X_{p}$ spans the symplectic-orthogonal (in $\dotCoTanM$) of its tangent space. 
    Now, we pullback $p \xxi$ via the projector $\dotCoTanM \times \dotCoTanN \to \dotCoTanM$ to consider it as a function $p \xxiyeta$ on $\dotCoTanM \times \dotCoTanN$. 
    The restriction of this function $p \xxiyeta$ on $C'$ is only a function of $(\xi, \eta)$ as $x^{i} = \partial H / \partial \xi_{i}$ and $y^{j} = \partial H / \partial \eta_{j}$ on $C'$. 
    Thus, $X_{p}$ on $C$ must be of the form 
    \begin{equation}
        - \parDeri{\bx^{i}}{} \big( p (\bx, \bxi) \big) \parDeri{\bxi_{i}}{} 
        = - f_{i} (\bx, \bxi, \bfeta) \parDeri{\bxi_{i}}{}, 
        \quad 
        \bx =  \parDeri{\bxi}{H}, 
    \end{equation}
    in our chosen parametrisation of $C$ and it is pointwise tangent to $C$.
    We compute 
    \begin{eqnarray}
        \LieDeri_{X_{p}} a 
        & = &  
        - f_{i} \Big(\parDeri{\bxi}{H}, \bxi, \bfeta \Big) \frac{\partial a}{\ri \partial \bxi_{i}} 
        - \frac{1}{2} \frac{\partial}{\ri \partial \bxi_{i}} \bigg( f_{i} \Big(\parDeri{\bxi}{H}, \bxi, \bfeta \Big) \bigg) \, a 
        \nonumber \\ 
        & = &  
        - f_{i} \Big(\parDeri{\bxi}{H}, \bxi, \bfeta \Big) \, \rD_{\bxi_{i}} (a) 
        - (\rD_{\bxi_{i}} f_{i}) \Big(\parDeri{\bxi}{H}, \bxi, \bfeta \Big) \, a 
        \nonumber \\ 
        && 
        + (\rD_{\bxi_{i}} f_{i}) \Big(\parDeri{\bxi}{H}, \bxi, \bfeta \Big) \, a 
        + \frac{\ri}{2} \frac{\partial}{\partial \bxi_{i}} \bigg( f_{i} \Big(\parDeri{\bxi}{H}, \bxi, \bfeta \Big) \bigg) \, a 
        \nonumber \\ 
        & = &  
        - f_{i} \Big(\parDeri{\bxi}{H}, \bxi, \bfeta \Big) \, \rD_{\bxi_{i}} (a) 
        - (\rD_{\bxi_{i}} f_{i}) \Big(\parDeri{\bxi}{H}, \bxi, \bfeta \Big) \, a 
        \nonumber \\ 
        && 
        - \ri \frac{\partial f_{i}}{\partial \bxi_{i}} \Big(\parDeri{\bxi}{H}, \bxi, \bfeta \Big) \, a 
        + \frac{\ri}{2} \frac{\partial f_{i}}{\partial \bxi_{i}} \Big(\parDeri{\bxi}{H}, \bxi, \bfeta \Big) \, a 
        + \frac{\ri}{2} \frac{\partial^{2} H}{\partial \bxi_{i} \partial \bxi_{j}} \parDeri{\bx^{j}}{f_{i}} \Big(\parDeri{\bxi}{H}, \bxi, \bfeta \Big) \, a 
        \nonumber \\ 
        & = &  
        - f_{i} \Big(\parDeri{\bxi}{H}, \bxi, \bfeta \Big) \, \rD_{\bxi_{i}} (a) 
        - (\rD_{\bxi_{i}} f_{i}) \Big(\parDeri{\bxi}{H}, \bxi, \bfeta \Big) \, a 
        \nonumber \\ 
        && 
        - \frac{\ri}{2} \Bigg( \frac{\partial f_{i}}{\partial \bxi_{i}} \Big(\parDeri{\bxi}{H}, \bxi, \bfeta \Big) 
        - \frac{\partial^{2} H}{\partial \bxi_{i} \partial \bxi_{j}} \parDeri{\bx^{j}}{f_{i}} \Big(\parDeri{\bxi}{H}, \bxi, \bfeta \Big) \Bigg) \, a 
        \nonumber \\ 
        & = &  
        - f_{i} \Big(\parDeri{\bxi}{H}, \bxi, \bfeta \Big) \, \rD_{\bxi_{i}} (a) 
        - (\rD_{\bxi_{i}} f_{i}) \Big(\parDeri{\bxi}{H}, \bxi, \bfeta \Big) \, a 
        - \frac{\ri}{2} \frac{\partial^{2} p}{\partial \bx^{j} \partial \bxi_{j}} \Big(\parDeri{\bxi}{H}, \bxi \Big)  
    \end{eqnarray}
    and insert this into the hindmost expression of $\symb{P \fA}$ to reach our goal:
    \begin{equation}
        \symb{P \fA} 
        = \bigg( \LieDeri_{X_{p}} + \fp_{m' - 1} \Big(\parDeri{\bxi}{H}, \bxi \Big) + \frac{\ri}{2} \frac{\partial^{2} p}{\partial \bx^{j} \bxi_{j}} \Big(\parDeri{\bxi}{H}, \bxi \Big) \bigg) a
        + \ldots 
        = (\LieDeri_{X_{p}} + \subSymb{P}) a 
        \mod S^{m + m' - 2} (\ldots).   
    \end{equation}
    Here $\fp_{m' - 1}$ is the pullback of the next-to-leading order term in the asymptotic expression of the total symbol of $P$ to $C'$ and we have used Definition~\ref{def: subprincipal_symbol} together with 
    \begin{equation}
        \frac{\partial^{2} p}{\partial \bx^{j} \bxi_{j}} (\bx, \bxi) 
        = 
        \parDeri{\bxi_{j}}{f_{j}} (\bx, \bxi, \bfeta)
        - \frac{\partial^{2} H}{\partial \bxi_{j} \partial \bxi_{i}} (\bxi, \bfeta) \, \parDeri{\bx^{j}}{f_{i}} (\bx, \bxi, \bfeta) 
        + \Big( \bx^{i} - \parDeri{\bxi_{i}}{H} (\bxi, \bfeta) \Big) \frac{\partial f_{i}}{\partial \bx^{j} \partial \bxi_{j}} (\bx, \bxi, \bfeta)  
    \end{equation}
    when evaluated at $\bx = \partial H / \partial \bxi$. 
\end{proof}
%
%
%
%
%
%
%
%
%
%
\subsection{Lagrangian distributions defined by symplectomorphisms}
\label{sec: FIO_symplecto}
The simplest canonical relation $\varGamma' \subset \dotCoTanM \times \dotCoTanN$ is that defined as the twisted graph $\varGamma'$ of a homogeneous symplectomorphism (see~\eqref{eq: def_symplecto} and~\eqref{eq: group_action_symplecto}) from $\dotCoTanN$ to $\dotCoTanM$, which enforces $d := \dim M = \dim N$ 
(see Example~\ref{exm: graph_symplectomorphism_Lagrangian_submf}). 
In this case, it is always possible to choose coordinates $(\by_{0}^{i})$ on the image of $y_{0} \in N$ under some local chart such that 
(see e.g.~\cite[Prop. 25.3.3]{Hoermander_Springer_2009})
\begin{equation}
    \varGamma_{\varphi} = \{ (\bx_{0}, \rd_{\bx_{0}} \varphi; \grad_{\bfeta^{0}} \varphi, \bfeta^{0}) \}, 
	\qquad 
	\det \frac{\partial^{2} \varphi}{\partial \bx \partial \bfeta} (\bx_{0}, \bfeta^{0}) \neq 0 
\end{equation} 
in an open conic neighbourhood of $(\bx_{0}, \bxi^{0}; \by_{0}, \bfeta^{0}) \in (\Rd \times \dotRd) \times (\Rd \times \dotRd)$ making use of a non-degenerate $\R$-valued phase function $\varphi (\bx, \bfeta)$ on an open conic neighbourhood of $(\bx_{0}, \bfeta^{0})$. 
\newline 

This class of Lagrangian submanifolds results in the simplest class of Fourier integral operators $A$ via the corresponding space of Lagrangian distributions $\fA \in \LagrangianDistGraphHomFE$ whose Euclidean representative is of the form~\cite[pp. 169-173]{Hoermander_ActaMath_1971}  
(see e.g.~\cite[(25.3.2)$^{\prime}$]{Hoermander_Springer_2009},~\cite[$(65.25)$]{Eskin_AMS_2011}) 
\begin{equation} \label{eq: Lagrangian_distribution_symplecto_local}
	\fA_{k}^{r} (\bx, \by) \equiv \int_{\Rd} \re^{\ri (\varphi (\bx, \bfeta) - \by \cdot \bfeta)} \fa_{k}^{r} (\bx, \by; \bfeta) \frac{\rd \bfeta}{(2 \pi)^{d}},     
\end{equation}  
where the total symbol $\fa_{k}^{r}$ is of homogeneous of degree $m$, having support in the interior of a small conic neighbourhood of $(\bx_{0}, \bfeta^{0})$ contained in the domain of definition of $\varphi$. 
\newline 

In order to describe the principal symbol $\symb{\fA}$ of $\fA$, one notes that the natural half-density $\dVol_{\ms \varGamma}$ (see~\eqref{eq: symplectic_form_graph_symplecto}) on $\varGamma$ can be factored out from $\symb{\fA} \xxiyeta \in S^{m + (d+d)/4} \big( \varGamma; \Maslov \otimes \widetilde{\mathrm{Hom}} (\sF, \sE) \otimes \halfDen \varGamma \big)$ so that the order $m + (d + d) / 4$ of the half-density-valued principal symbol is reduced to $m$ as $\halfDen \varGamma$ is of order $d/2$.  
Therefore, the principal symbol map  
\begin{equation} \label{eq: symbol_FIO_symplecto}
	\symb{} : I^{m - [1]} \big( M \times N, \varGamma'; \halfDenMN \otimes \Hom{\sF, \sE} \big) \to S^{m - [1]} \big( \varGamma; \Maslov \otimes \widetilde{\mathrm{Hom}} (\sF, \sE) \big)  
\end{equation}
is locally given by  
(see e.g.~\cite[p. 27]{Hoermander_Springer_2009}) 
\begin{equation} \label{eq: symbol_FIO_symplecto_local}
	\sigma_{A_{k}^{r}} (\bx, \bfeta) := a_{k}^{r} \left( \bx, \parDeri{\bfeta}{\varphi}; \bfeta \right) \, \left| \det \frac{\partial^{2} \varphi}{\partial \bx \partial \bfeta} \right|^{-1/2} \m 
    \mod S^{m-1} (\cdot),     
\end{equation} 
where the Keller-Maslov bundle $\Maslov$ is trivialised by the phase function in~\eqref{eq: Lagrangian_distribution_symplecto_local}. 
\newline 

We remark that the composition 
$\fA \fB \in I^{m + m'} \big( M \times N, (\varGamma \circ \varLambda)'; \Hom{\sF, \sE} \otimes \halfDenMN \big)$ 
between 
$\fA \in I^{m} \big( M \times O, \varGamma', \Hom{\sG, \sE} \otimes \halfDen (M \times O) \big)$ 
and a properly supported 
$\fB \in I^{m'} \big( O \times N, \varLambda'; \Hom{\sF, \sG} \otimes \halfDen (O \times N) \big)$ 
is always well-defined, where $\varGamma \circ \varLambda$ is the graph of the composition of symplectomorphisms $\dotCoTanN \ni \yeta \mapsto (z, \zeta) \mapsto \xxi \in \dotCoTanM$. 
Furthermore, the expression of $\symb{\fA \fB}$ simplifies quite a bit compared to~\eqref{eq: symbol_product_Lagrangian_dist} and it is given by the standard composition of homomorphisms~\cite[p. 180]{Hoermander_ActaMath_1971}  
(see also, e.g.~\cite[Thm. 65.7]{Eskin_AMS_2011},~\cite[$(6.11)$ $($p. 465$)$]{Treves_Plenum_1980}) 
\begin{equation} \label{eq: product_symbol_FIO_symplecto}
    \symb{\fA \fB} \xxiyeta 
    = 
    \sum_{(z, \zeta) \,|\, (x, \xi; z, \zeta) \in \varGamma, (z, \zeta; y, \eta) \in \varLambda} \symb{\fA} (x, \xi; z, \zeta) \, \big( \symb{\fB} (z, \zeta; y, \eta) \big). 
\end{equation} 
If the respective vector bundles are hermitian then the algebra of Lagrangian distributions is a $*$-algebra. 
%
%
%
%
%
%
%
%
%
%
\subsubsection{Parametrices and elliptic operators}
We will now introduce the notion of ellipticity~\cite[p. 186]{Duistermaat_ActaMath_1972} 
for a Fourier integral operator and show that an approximate inverse always exists for such an operator.  

%
%
%
\begin{definition} \label{def: elliptic_FIO}
    As in the terminologies of Definitions~\ref{def: FIO} and~\ref{def: symbol_FIO}, and~\eqref{eq: symbol_FIO_symplecto}, let $\varGamma$ be the graph of a homogeneous symplectomorphism from $\dotCoTanN$ to $\dotCoTanM$. 
    Then, a Lagrangian distribution $\fA \in \LagrangianDistGraphHomFE$ is called \textbf{non-characteristic} at $\xxiyetaNot \in \varGamma$ if its principal symbol has an inverse $\in S^{-m} \big( \varGamma; \Maslov^{-1} \otimes \widetilde{\mathrm{Hom}} (\sE, \sF) \big)$ in a conic neighbourhood of $\xxiyetaNot$. 
	The Lagrangian distribution $\fA$ is called \textbf{elliptic} if it is non-characteristic at every point of $\varGamma$. 
    The complement of non-characteristic points is called the \textbf{characteristic set} $\Char{\fA}$ of $\fA$~\cite[Def. 25.3.4]{Hoermander_Springer_2009}.   
\end{definition}
%
%
%

We note that~\eqref{eq: Lagrangian_distribution_symplecto_local} and~\eqref{eq: symbol_FIO_symplecto_local} 
imply that the non-characteristic points belong to $\WFPrime{\fA}$. 
If $\varGamma^{-1}$ is also a graph and $\fA \in \LagrangianDistGraphHomFE$ is elliptic and properly supported, then $\fA$ has a \textbf{unique parametrix} $\fE \in I^{-m} \big( N \times M, \varGamma^{-1 \prime}; \Hom{\sE, \sF} \otimes \halfDen (N \times M) \big)$. 
In other words, 
\begin{equation} \label{eq: parametrix_elliptic_FIO}
	EA - I_{\sF} \in \PsiDO{- \infty}{N; \sF \otimes \halfDenN}, 
	\qquad 
	AE - I_{\sE} \in \PsiDO{- \infty}{M; \sE \otimes \halfDenM}. 
\end{equation}
A proof of the microlocal version of this claim is going to be presented shortly after one devises a variant of Lagrangian distributions where the closedness assumption on the Lagrangian submanifold has been relaxed. 

%
%
%
\begin{definition} \label{def: Lagrangian_distribution_without_closed}
	As in the terminologies of Definition~\ref{def: Lagrangian_distribution}, let $C \subset \dotCoTanM \times \dotCoTanN$ be a homogeneous canonical relation which is not necessarily closed in $\dotCoTanMN$ and $\cK \subset C$ a conic subset which is closed in $\dotCoTanMN$. 
    By $I^{m} \big( M \times N, \cK'; \Hom{\sF, \sE} \otimes \halfDenMN \big)$, one denotes the set of all matrices with entries~\eqref{eq: Lagrangian_distribution_local} together with the additional condition that, for each $\alpha, \beta$ as in~\eqref{eq: Lagrangian_distribution_mf_local}, the restriction of $\fa_{\alpha, \beta; k}^{\ \; \; \; \; r}$ to some conic neighbourhood in $\R^{\dM} \times \R^{\dN} \times \R^{n_{\alpha, \beta}}$ of the pullback of $C \setminus \cK$ by the Lagrangian fibration (see~\eqref{eq: def_canonical_relation_fibration}) $\jmath: \sC_{\alpha, \beta} := \eqref{eq: fibre_critical_mf_Euclidean} \to C'_{\alpha, \beta} := \eqref{eq: Euclidean_rep_canonical_relation}$, is in the class $S^{-\infty} (\cdot)$. 
    \newline 

    As in the terminologies of Definition~\ref{def: symbol_FIO}, 
    $S^{m} \big( \cK, \Maslov \otimes \halfDenC \otimes \widetilde{\mathrm{Hom}} (\sF, \sE) \big)$ denotes the set of all $a \in S^{m} \big( C, \Maslov \otimes \halfDenC \otimes \widetilde{\mathrm{Hom}} (\sF, \sE) \big)$ such that $a \in S^{-\infty} \big( C \setminus \cK, \Maslov \otimes \halfDenC \otimes \widetilde{\mathrm{Hom}} (\sF, \sE) \big)$~\cite[p. 187]{Duistermaat_ActaMath_1972}.  
\end{definition}
%
%
%

Note, the definition depends on both $\cK$ and $C$ yet suppressed to lighten the notion.
The analogue of the principal symbol isomorphism (Definition~\ref{def: symbol_FIO}) reads~\cite[p. 187]{Duistermaat_ActaMath_1972}    
\begin{equation} \label{eq: def_symbol_FIO_microlocal}
    I^{m - [1]} \big( M \times N, \cK'; \Hom{\sF, \sE} \otimes \halfDenMN \big) 
    \cong 
    S^{m + \frac{\dM + \dN}{4} - [1]} \big( \cK, \Maslov \otimes \halfDenC \otimes  \widetilde{\mathrm{Hom}} (\sF, \sE) \big). 
\end{equation}  

%
%
%
\begin{theorem} \label{thm: existence_parametrix_FIO}
	As in the terminologies of Definitions~\ref{def: FIO} and~\ref{def: Lagrangian_distribution_without_closed}, suppose that $\varGamma$ is the graph of a homogeneous symplectomorphism $\varkappa$ from an open conic subset $\cV \subset \dotCoTanN$ into $\dotCoTanM$ and that $\cK \subset \varGamma$ is a conic subset which is closed in $\dotCoTanMN$.  
    Then, for any conic subset $\cU \subset \cV$ such that $\cU$ (resp. $\varkappa (\cU)$) is closed in $\dotCoTanN$ (resp. $\dotCoTanM$), if $\fA \in I^{m} \big( M \times N, \cK'; \Hom{\sF, \sE} \otimes \halfDenMN \big)$ is non-characteristic on $\{ \big( \varkappa (\cU), \cU \big) \} \subset \varGamma$, then there exists an elliptic $\fE \in I^{-m} \big( N \times M, \cK^{-1 \prime}; \Hom{\sE, \sF} \otimes \halfDen (N \times M) \big)$ such that 
	\begin{equation} \label{eq: microlocal_parametrix_elliptic_FIO}
		\ES (EA - I_{\sF}) \cap \cU = \emptyset, 
		\qquad 
		\ES (AE - I_{\sE}) \cap \varkappa (\cU) = \emptyset, 
	\end{equation}
    where $I_{\sE} \in \PsiDO{0}{M; \sE \otimes \halfDenM}$ and $I_{\sF} \in \PsiDO{0}{N; \sF \otimes \halfDenN}$ are identity operators.  
	The parametrix $\fE$ is unique in the sense that $\big( \cU, \varkappa (\cU) \big) \nsubset \WF' (\fE - \fE')$ for any other parametrix $\fE'$ of $\fA$.
\end{theorem}
%
%
%
\begin{proof}
	The arguments used to prove the 
    scalar-version~\cite[Prop. 5.1.2]{Duistermaat_ActaMath_1972} 
    (see also, e.g.~\cite[Thm. 4.2.5]{Duistermaat_Birkhaeuser_2011}) 
    flows over the bundle case albeit provided here for completeness.   
	Let $a$ be a principal symbol of $\fA$ in the sense of~\eqref{eq: def_symbol_FIO_microlocal}.  
    By hypotheses, there exists a $b (\cdot) \in S^{-m} \big( \cK, \Maslov^{-1} \otimes \halfDenC \otimes \widetilde{\mathrm{Hom}} (\sE, \sF) \big)$ such that $ba = \symb{I_{\ms \sF}}$ in a conic neighbourhood of $\varDelta \cU$. 
    By utilising appropriate 
    microlocal partition of unities\footnote{As 
        in the terminologies of Definition~\ref{def: PsiDO}, a family $\{\varPsi_{\alpha}\}_{\alpha}$ of $\varPsi_{\alpha} \in \PsiDO{0}{M; \sE \otimes \halfDen}$ is called a \textbf{microlocal partition of unity} if the supports of the Schwartz kernels of $\varPsi_{\alpha}$'s  are locally finite, $\WF (\varPsi_{\alpha}) \subset \cU_{\alpha}$ and $\sum \varPsi_{\alpha} = I$. 
        Here, $(\cU_{\alpha})_{\alpha}$ is an open cover of the cosphere bundle $\bbS^{*} M$ such that both $\halfDen$ and $\sE$ admit trivialisations over its projection on $M$.  
        A microlocal partition of unity always exists 
        (see e.g.~\cite[Lem. 6.10, 6.11]{Hintz_2019} for a scalar version). 
    } \label{foot: microlocal_partition_unity}
    $\varPsi$ (resp. $\varPhi$) subordinated to $\cU$ (resp. $\varkappa (\cU)$), one chooses $b = 0$ outside of a sufficiently small conic neighbourhood of $\{ \big( \varkappa (\cU), \cU \big) \}$ where $\cU$ is identified with $\varDelta \cU$ (via projection). 
    Then we have a properly supported $\fB_{0} \in I^{-m} \big( N \times M, \cK^{-1 \prime}; \Hom{\sE, \sF} \otimes \halfDen (N \times M) \big)$ whose principal symbol is $\symb{\varPsi} b \symb{\varPhi} \equiv b$, and by the composition of Lagrangian distributions (\eqref{eq: product_Lagrangian_dist} and~\eqref{eq: symbol_product_Lagrangian_dist}), there exists a properly supported $R \in \PsiDO{-1}{N; \sF \otimes \halfDenN}$ such that $B_{0} A = I_{\sF} - R$.  
	Next, we want to invert $I_{\sF} - R$ by making use of the Neumann series: 
	$(I_{\sF} - R)^{-1} = \sum_{k \in \N_{0}} R^{k}$.  
	Set   
    $\fB_{k} := R^{k} \fB_{0} \in I^{-m - k} \big( N \times M, \cK^{-1 \prime}; \Hom{\sE, \sF} \otimes \halfDen (N \times M) \big)$ 
	and then 
	\begin{equation}
		I_{\sF} 
		= 
		(I_{\sF} - R) \sum_{k=0}^{\infty} R^{k}  
		= 
		\sum_{k=0}^{N-1} R^{k} (I_{\sF} - R) + \sum_{k=N}^{\infty} R^{k} (I_{\sF} - R) 
		= 
		\sum_{k=0}^{N-1} B_{k} A + R^{N}. 
	\end{equation}
	Let $\fE$ be defined by the asymptotic summation: $\fE :\sim \sum_{k \in \N_{0}} \fB_{k}$, then inserting the last equation one obtains   
	\begin{equation}
		EA - I_{\sF} 
		=
		\left(E - \sum_{k=0}^{N-1} B_{k} \right) A - R^{N} 
		\sim  
		\sum_{k = N}^{\infty} B_{k} A - R^{N} 
		\in \PsiDO{- N}{N; \sF \otimes \halfDenN}  
	\end{equation}
	for every $N \in \N$, which in turn proves the first part of the theorem as a right parametrix can be constructed analogously. 
	\newline 

	To prove the uniqueness, we suppose that $\fE'$ is another right parametrix for $\fA$. 
	Then  
	\begin{equation}
        \big\{ \big( \cU, \varkappa (\cU) \big) \big\} \nsubset \WF' \fE = \WF' (\fE \fA \fE') = \WF' \fE' 
	\end{equation}
	and similarly for the left parametrix. 
\end{proof}
%
%
%
%
%
%
%
%
%
%
\subsubsection{Egorov's theorem}
We are now going to present a key theorem for microlocalisation (see Theorem~\ref{thm: microlocalisation_P}) related to the constructions presented in this section.

%
%
%
\begin{theorem}[Egorov's theorem] \label{thm: Egorov}
	As in the terminologies of Definition~\ref{def: FIO}, let $\varGamma$ be the graph of a homogeneous symplectomorphism $\varkappa: \dotCoTanN \to \dotCoTanM$. 
    Suppose that $\fA \in \LagrangianDistGraphHomFE, \fB \in I^{-m} \big( N \times M, \varGamma^{-1 \prime}; \Hom{\sE, \sF} \otimes \halfDen (N \times M) \big)$ and that $P \in \PsiDO{m'}{M; \sE \otimes \halfDenM}$ having a scalar principal symbol $\symb{P}$. 
    Then, $\fB P \fA \in \PsiDO{m'}{N; \sF \otimes \halfDenN}$ is properly supported whose principal symbol is given by 
	\begin{equation}
        \symb{\fB P \fA} = \symb{\fB \fA} (\varkappa^{*} \symb{P}),  
	\end{equation}
    where $\varkappa^{*} : C^{\infty} \big( \dotCoTanM, \Hom{\sE, \sE} \big) \to C^{\infty} \big( \dotCoTanN, \Hom{\sF, \sF} \big)$ is the pullback map via $\varkappa$. 
\end{theorem}
%
%
%

\begin{proof}
	We will use the same strategy employed for proving the scalar version~\cite[Thm. 25.3.5]{Hoermander_Springer_2009}. 
	By repeated applications of the composition of Lagrangian distributions,  we have 
    $P \fA \in I^{m+m'} \big( M \times N, \varGamma'; \Hom{\sF, \sE} \otimes \halfDenMN \big)$ and  
	$\fB P \fA \in \PsiDO{m'}{N; \sF \otimes \halfDenN}$.  
    In order to compute $\symb{P \fA}$, one lifts $\symb{P}$ to $\varGamma$ via the pullbacks of the projector $\dotCoTanM \times \dotCoTanN \to \dotCoTanM$ followed by the inclusion $\varGamma \hookrightarrow \dotCoTanM \times \dotCoTanN$. 
    Then $\symb{P \fA} = \symb{P} \symb{\fA}$. 
    Equivalently, one can lift $\varkappa^{*} \symb{P}$ to $\dotCoTanM \times \dotCoTanN$ via the projector $\dotCoTanM \times \dotCoTanN \to \dotCoTanN$ and then consider it as a homomorphism on $\varGamma$ as before. 
    Thus $P \fA - \fA R \in I^{m + m' - 1} \big( M \times N, \varGamma'; \Hom{\sF, \sE} \otimes \halfDenMN \big)$ if $R \in \PsiDO{m'}{N; \sF \otimes \halfDenN}$ having the principal symbol $\varkappa^{*} \symb{P}$. 
    As $\symb{P}$ is scalar, therefore $\fB P \fA - \fB \fA R \in \PsiDO{m' - 1}{N; \sF \otimes \halfDenN}$ which entails the claim. 
\end{proof}
%
%
%
%
%
%
%
%
%
%
\section{Literature}
Historically, the theory \textit{pseudodifferential operator}s was developed from the theory of singular integral operators by
Joseph J. \textsc{Kohn} and Louis \textsc{Nirenberg}~\cite{Kohn_CPAM_1965} 
and enriched, in particular, by
Lars \textsc{H\"{o}rmander}~\cite{Hoermander_CPAM_1965}. 
\newline 

\textit{Fourier integral operator}s are an enormous generalisation of pseudodifferential operators,  originated from the investigation of the singularities of solutions of hyperbolic differential equations by
Peter \textsc{Lax}~\cite{Lax_DukeMathJ_1957} 
and 
in the study of asymptotic analysis in the context of geometric optics and quantum mechanics by 
Victor \textsc{Maslov}~\cite[Chap. 8]{Maslov_Springer_1981}. 
Their somewhat local formulation was later systematically developed and globalised in two seminal articles by  
Lars \textsc{H\"{o}rmander}~\cite{Hoermander_ActaMath_1971} 
and by 
Johannes J. \textsc{Duistermaat} and Lars \textsc{H\"{o}rmander}~\cite{Duistermaat_ActaMath_1972} 
for scalar operators while the bundle version is available in 
the fourth volume of H\"{o}rmander's treatise~\cite{Hoermander_Springer_2009}. 
\newline 

\textit{Egorov}'s theorem is named after  
Yu. V. \textsc{Egorov}~\cite{Egorov_UMN_1969} 
who proved it at first.   
Since then various generalisation of this result has been reported sporadically, 
e.g.~\cite{Dencker_JFA_1982},~\cite[Thm. 3.2]{Bolte_CMP_2004},~\cite[Thm. 1.7]{Kordyukov_MPAG_2005},~\cite[Thm. 6.1]{Kordyukov_JGP_2007},~\cite[Prop. 3.3]{Jakobson_CMP_2007},~\cite[Prop. A.3]{Gerard_CMP_2015}.  
Theorem~\ref{thm: Egorov} is slightly general than the statements available in literature to the best of our knowledge. 
\newline 

We refer, for instance, the 
monographs~\cite{Hoermander_Springer_2007, Treves_Plenum_1980, Taylor_PUP_1981, Shubin_Springer_2001} 
and the expository 
lecture notes~\cite{Melrose_2007, vandenBan_2017, Hintz_2019} for an in-depth discussion on pseudodifferential operators and microlocal analysis. 
For Fourier integral operators, the 
textbooks~\cite{Duistermaat_Birkhaeuser_2011, Treves_Plenum_1980, Grigis_CUP_1994} 
are referred for details. 

%% file: Feynman_propagator.tex
\chapter{Feynman Propagators}
\label{ch: Feynman_propagator}
\textsf{The subject matter of this chapter is a microlocal construction of Feynman propagators for normally hyperbolic operators and for Dirac-type operators on a globally hyperbolic spacetime. 
These results are shown by constructing Feynman parametrices and by employing the well-posed Cauchy problem. 
A vector bundle generalisation of the microlocalisation of a pseudodifferential operator of real principal type is achieved in order to give the microlocal construction of Feynman parametrices.}
%
%
%
%
%
%
%
%
%
%
\section{Lorentzian manifolds} 
\label{sec: spacetime}
By a \textbf{Lorentzian manifold} $\spacetime$ we mean a Hausdorff, second countable, connected, and smooth $d \geq 2$-dimensional manifold $\sM$ endowed with a smooth Lorentzian metric $\fg$ of signature $(+, - \ldots -)$. 
The global existence of Lorentzian metrics depends on the topology of $\sM$: \textit{they always exist on any non-compact $\sM$, whereas a compact $\sM$ admits a Lorentzian metric if and only if its Euler characteristic vanishes}. 
Hausdorff condition together with the existence of a Lorentzian metric entails that $\sM$ is paracompact. 
\newline 

Since $\fg$ is a (symmetric) non-degenerate bilinear map $\fg : C^{\infty} (\sM; \tansM) \times C^{\infty} (\sM; \tansM)$ $\to C^{\infty} (\sM, \R)$, it induces a unique vector bundle (linear) isomorphism $C^{\infty} (\sM; \tansM) \ni X \mapsto X^{\flat} \in C^{\infty} (\sM; \coTansM)$ defined pointwise via $X_{x}^{\flat} (\cdot) := \fg_{x} (X_{x}, \cdot)$ for any $x \in \sM$. 
Moreover, it induces another unique symmetric, non-degenerate, indefinite bilinear form $\fg^{-1}$ on $\coTansM$, pointwise given by $\fg_{x}^{-1} (\xi_{x}, \eta_{x}) := \fg_{x} (\xi_{x}^{\sharp}, \eta_{x}^{\sharp})$ where $\sharp : C^{\infty} (\sM, \coTansM) \to C^{\infty} (\sM; \tansM)$ is the inverse of $\flat$. 
%
%
%
%
%
%
%
%
%
%
%
\subsection{Causal structure}
\label{sec: causal_structure_Lorentzian_mf}
The causal structure is one of the drastic differences between a Lorentzian and a Riemannian manfold. 
The signature of $\fg$ allows one to make the following classification. 

%
%
%
\begin{definition}
    A cotangent vector $\xi \in \coTan_{x} \sM$ on any point $x$ in a Lorentzian manifold $\spacetime$ is called \textbf{timelike} if $\fg_{x}^{-1} (\xi, \xi) > 0$, \textbf{spacelike} if $\fg_{x}^{-1} (\xi, \xi) < 0 $ or $\xi = 0$, \textbf{lightlike} if $\fg_{x}^{-1} (\xi, \xi) = 0 $ whenever $\xi \neq 0$, and \textbf{causal} if $\xi$ is either timelike or lightlike.
\end{definition}
%
%
%

The attribution of any covector according to this classification is called its causal character.  
The set $\dotCoTan_{0, x} \sM \subset \dotCoTan_{x} \sM$ of all lightlike covectors at $x \in \sM$ is called the \textbf{lightcone} at $x$ and it is an inner product space. 
One notices that, if $\xi, \eta \in \dotCoTan_{x} \sM$ are timelike then $\fg_{x}^{-1} (\xi, \eta) \lessgtr 0$. 
Thus, at each $x \in \sM$, the cotangent vector space $\dotCoTan_{x} \sM$ has two possible open convex cones of timelike covectors whose boundaries consist of the lightlike covectors, and there is no intrinsic way to distinguish one from another albeit it is crucial to do so for some purpose. 
$\dotCoTan_{x} \sM$ is called \textbf{time-oriented} if one makes a choice (termed \textit{future}) for one of these cones; the other cone is then termed as \textit{past}. 
Whether such a designation between future and past is possible in a continuous fashion as $x$ varies over $\sM$ gives the notion of time-orientability of $\spacetime$, which is logically independent of the topological orientation of $\sM$. 
If $\mathfrak{t}$ is a smooth function on $\sM$ that assigns to each $x \in \sM$ either a future or equivalently a past cone, then it is called a \textbf{time-orientation} of $\spacetime$. 
If $\sM$ admits a time-orientation then $\spacetime$ is said to be time-orientable. 

%
%
%
\begin{definition} \label{def: spacetime}
    A time-oriented and orientated Lorentzian manifold is called a \textbf{spacetime} $\spacetime$.
\end{definition}
%
%
%

Note, a Lorentzian manifold is time-orientable if and only if it admits a smooth (highly non-unique) \textit{global timelike vector field} $X$. 
Let us pick a time-orientation once and for all such that $X$ is future-directed at all points and then any causal $\xi \in \dotCoTan_{x} \sM$ is future directed if and only if $\fg_{x}^{-1} (\xi, X^{\flat}) > 0$. 
Therefore, we have a consistent notion of future and past of any point in spacetime. 
To be precise, any differentiable curve $c : \I \to \sM$ is called spacelike (resp. future/past directed causal) if its tangent $\nicefrac{\rd c}{\rd s}$ is spacelike (resp. future/past directed causal), where $\I \subset \R$ is an interval and $s \in \I$. 
\newline 

We now introduce the causal structure of spacetime $\spacetime$. 
If $x, y \in \sM$, then 
\begin{itemize}
    \item 
    $x~ \rhd~y$ means that there is a future directed timelike curve in $\sM$ from $x$ to $y$; 
    \vspace*{-0.25cm}
    \item 
    $x~\unrhd~y$ means that either there is a future directed causal curve in $\sM$ from $x$ to $y$, or $x = y$.
\end{itemize}

The \textbf{causal future} $J^{+} (x) := \{ y \in \sM \,|\, x~ \unrhd~ y \}$ of a point $x \in \sM $ is defined as the set of all points that can be reached by a future directed causal curve in $\sM$ emanating from $x$ and $x$ itself. 
Analogously, $J^{-} (x) := \{ y \in \sM \,|\, y~ \unrhd~ x \}$. 
Then the causal future (past) of a set $K \subset \sM$ is defined by   
\begin{equation}
    J^{\pm} (K) 
    := \bigcup_{x \in K} J^{\pm} (x)
\end{equation}
and one also uses the notation
\begin{equation} \label{eq: causal_future_past}
	J (K) := J^{+} (K) \cup J^{-} (K). 
\end{equation}
If one replaces causal curves with timelike curves then it yields \textbf{chronological} future/past. 
%
%
%
%
%
%
%
%
%
%
\subsection{Globally hyperbolic spacetimes} 
\label{sec: globally_hyperbolic_spacetime}
A generic spacetime is too general to consider as it may contain pathological mathematical and physical properties. 
For instance, the initial-value problem for a wave equation is not well-posed on an arbitrary spacetime. 
Also, there may exist closed causal curves favouring the ``grandfather paradox''. 
Therefore, we look for a subset of spacetimes excluding all such undesirable possibilities, as introduced below. 

%
%
%
\begin{definition} \label{def: globally_hyperbolic_spacetime}
    A spacetime $\spacetime$ is called \textbf{globally hyperbolic} if and only if it is causal and the diamonds $J^{+} (x) \cap J^{-} (y)$ are compact for all $x, y \in \sM$.  
\end{definition}
%
%
%

In what follows, any globally hyperbolic spacetime admits a specific global orthogonal splitting. 
Such a characterisation is of utmost importance for our thesis and we introduce the following terminologies in order to formulate the characterisation.  

%
%
%
\begin{definition} \label{def: Cauchy_hypersurface}
    A \textbf{Cauchy hypersurface} $\varSigma$ in a Lorentzian manifold $\spacetime$ is a subset $\varSigma \subset \sM$ which is met exactly once by every inextensible timelike curve.
\end{definition}
%
%
%

\begin{definition} \label{def: time_function}
	Any $\ft \in C^{\infty} (\sM, \R)$ on a spacetime $\spacetime$ 
	is called a 
	\begin{itemize}
		\item 
		\textbf{time function} if $\ft$ is strictly increasing along all future-directed causal curves; 
		\item 
        \textbf{temporal function} if $\ft$ is smooth having a future-directed timelike gradient; 
	\end{itemize}
    A time (and/or temporal) function is called \textbf{Cauchy} if its level sets $\varSigma_{t} := \ft^{-1} (t)$ are Cauchy hypersurfaces for any $t \in \R$. 
\end{definition}
%
%
%

We now inscribe the global orthogonal splitting of a globally hyperbolic spacetime following the survey article~\cite[Thm. 3.78]{Minguzzi_EMS_2008}. 

%
%
%
\begin{theorem}[Geroch~\cite{Geroch_JMP_1970}-Bernal-S\'{a}nchez~\cite{Bernal_CMP_2003, *Bernal_CMP_2005, *Bernal_LMP_2006, *Bernal_CQG_2007} splitting theorem] \label{thm: globally_hyperbolic_Cauchy}
    Let $\spacetime$ be a spacetime. 
    It is globally hyperbolic if and only if it admits  a smooth spacelike Cauchy hypersurface $\varSigma$. 
    \newline 

    In this case, it allows a Cauchy temporal function $\ft$ and, thus, it is isometrically diffeomorphic to the smooth product manifold 
    \begin{equation}
        \spacetime \cong (\R \times \varSigma, \betaup \, \rd t^{2} - \fh_{t}),   
    \end{equation}
    where $\betaup \in C^{\infty} (\R \times \varSigma, \R_{+})$ is the lapse function, 
    $t \in C^{\infty} (\R \times \varSigma, \R)$ is the natural projection, 
    each level set $\varSigma_{t}$ of $t = \cst$ is a spacelike Cauchy hypersurface, and $\fh_{t}$ is a Riemannian metric on each $\varSigma_{t}$, which varies smoothly with $t$. 
    \newline 

    Moreover, if $\mathsf{\Sigma}$ is a prescribed topological Cauchy hypersurface then there exists a smooth Cauchy function $\ft$ such that $\mathsf{\Sigma}$ is one of its level sets ($\mathsf{\Sigma} = \mathsf{\Sigma}_{0}$). 
    In addition, if $\mathsf{\Sigma}$ is 
    \begin{itemize}
    	\item 
        acausal then $\ft$ becomes a smooth Cauchy time function; 
        \item 
        spacelike then $\ft$ can be modified to obtain a Cauchy temporal function on $\sM$ such that $\mathsf{\Sigma} = \ft^{-1} (0)$. 
    \end{itemize}

    Furthermore, if $S \subset \sM$ is a compact achronal set then it can be extended to a Cauchy hypersurface. 
    Additionally, if $S$ is acausal and a smooth spacelike submanifold with boundary then it can be extended to a spacelike Cauchy hypersurface $\varSigma \supset S$. 
\end{theorem}
%
%
%

\begin{remark} \label{rem: Cauchy_hypersurface_embedded_submf}
    Globally hyperbolic spacetimes $\spacetime$ can always be foliated by a Cauchy hypersurface $\varSigma$, but the foliation is \textit{non-unique}. 
	Necessarily, $\varSigma$ is a closed subset and an \textit{embedded} topological submanifold of codimension $1$. 
\end{remark}
%
%
%

\begin{example} \label{ex: RW_spacetime}
	Let $\I \subset \R$ be open and let $(\varSigma, \fh)$ be a connected Riemannian manifold. 
	Then $(\sM := \I \times \varSigma, \fg := \rd t^{2} - f^{2} \fh)$ is a globally hyperbolic spacetime if and only if $(\varSigma, \fh)$ is \textit{complete}, where $f \in C^{\infty} (\I, \R_{+})$ 
    (see e.g.~\cite[Lem. A.5.14]{Baer_EMS_2007}). 
\end{example}
%
%
%

Physically speaking, our understanding of the physical universe is currently best explained by the ``$\Lambda$CDM model'' where the spacetime geometry is given by the  
Robertson-Walker spacetime (see e.g.~\cite[Sec. 5.3]{Hawking_CUP_1973},~\cite[Def. 7, p. 343]{ONeill_Academic_1983}).  
This so-called standard model of cosmology describes, for instance, the cosmic microwave background, expansion of the universe, and the cosmological Redshift pretty accurately   
(see e.g.~\cite{Green_CQG_2014}). 
However, the standard model has a few limitations and the 
de Sitter spacetime (see e.g.~\cite[Sec. 5.2]{Hawking_CUP_1973},~\cite[p. 229]{ONeill_Academic_1983}) 
offers a solution of some of these open problems.  
Example~\ref{ex: RW_spacetime} encompasses both these spacetimes together with, of course, the Minkowski spacetime. 
Besides, a number of black hole solutions of Einstein equation are globally hyperbolic, e.g., the interior and exterior of Schwarzschild spacetime, parts of regions of the Reissner–Nordstr\"{o}m, and the Kerr spacetimes. 
There are, of course, solutions of Einstein equation which are not globally hyperbolic at all. 
One of such spacetimes received huge attention in the context of superstring theory is the 
anti-de Sitter solution (see e.g.~\cite[Sec. 5.2]{Hawking_CUP_1973},~\cite[p. 229]{ONeill_Academic_1983}).  
Roughly speaking, globally hyperbolic spacetimes are interesting because they offer prediction (or retrodiction) of the entire future (or past) history of the universe from conditions at an instant of time. 

%
%
%
\begin{remark}
	If $\spacetime$ is a globally hyperbolic spacetime then the set $\coLightBun \to \sM$ of lightcones (lightlike geodesics on $\coTansM$) has a smooth bundle structure 
	(see e.g.~\cite{Low_JMP_1989}). 
\end{remark}
%
%
%
%
%
%
%
%
%
%
\section{Normally hyperbolic operators and parametrices}
\subsection{Normally hyperbolic operators (NHOs)}
\label{sec: NHOp}
%
%
%

\begin{definition} \label{def: NHOp}
	Let $\sE \to \sM$ be a vector bundle over a Lorentzian manifold $\spacetime$.  
	A second-order linear differential operator acting on half-density-valued smooth sections of $\sE$ 
	(see e.g.~\cite[Def. 3.1]{Baum_AGAG_1996},~\cite[Sec. 1.5]{Baer_EMS_2007})  
    \begin{equation}
    	\square : \secsME \to \secsME
    \end{equation}
    is called a \textbf{normally hyperbolic operator} if its principal symbol is given by the spacetime metric $\fg^{-1}$ on the cotangent bundle $\coTansM$ 
	\begin{equation}
        \symb{\square} \xxi := \fg_{x}^{-1} (\xi, \xi) ~\one_{(\End \sE)_{x}} 
    \end{equation}
	for all $\xxi \in C^{\infty} (\sM; \coTansM)$. 
\end{definition}
%
%
%

In a local coordinate chart $\big(U, (x^{i}) \big)$ of $\sM$, after trivialising $\sE$, $\square$ is given by  
\begin{equation}
   	\square \upharpoonright U  
   	= 
   	- \fg^{ij} \frac{\partial^{2}}{\partial x^{i} \partial x^{j}} 
   	+ \tA^{i} \frac{\partial}{\ri \partial x^{i}} 
   	+ \tB, 
   	\quad i,j = 1, \ldots, d, 
\end{equation}
where $\tA^{i} \in C^{\infty} (U; \tangent U \otimes \End \sE_{U}), \tB \in C^{\infty} (U, \End \sE_{U})$. 
\newline 

Normally hyperbolic operators arise naturally in the context of geometric analysis and quantum field theories in curved spacetime, and we list below a few. 

%
%
%
\begin{example} \label{exm: covariant_Klein_Gordon_op}
	Let $\sE \to \sM$ be a trivial-$\C$-line bundle, that is, its sections are just $\C$-valued smooth functions $C^{\infty} (\sM)$ on $\spacetime$. 
    The \textbf{covariant Klein-Gordon operator} is defined by 
    (see e.g.~\cite[Sec. 3.1, Exm. 1]{Baum_AGAG_1996},~\cite[Exm. 1.5.1]{Baer_EMS_2007}) 
    \begin{equation} \label{eq: def_covariant_Klein_Gordon_op}
        \square := - \Div \circ \grad - \rrm^{2} - \lambda \ric : C^{\infty} (\sM) \to C^{\infty} (\sM), 
    \end{equation}
	where $\rrm^{2} \in \R_{+}$ is a parameter, physically interpreted (in appropriate situation) as the mass-squared of a linear Klein-Gordon field, $\lambda$ is a coupling, $\ric$ is the Ricci scalar of $\sM$, $\grad f := (\rd f)^{\sharp}$ is the gradient of a function $f \in C (\sM)$, and $\Div (X) := \tr (\connectionLC X)$ is the divergence of a vector field $X \in C (\sM; \tangent \sM)$. 
    Here $\connectionLC$ is the Levi-Civita connection of $\sM$ and the \textbf{endomorphism trace} 
    \begin{equation} \label{eq: def_endo_trace}
        \tr : \End \sE \to \C
    \end{equation}
    is the finite-dimensional trace on the fibres $\sE_{x}$ of any vector bundle $\sE \to M$ over a generic  manifold $M$, defined by the \textit{composition of the canonical isomorphism $\End (\sE_{x}) \cong \sE^{*}_{x} \otimes \sE_{x}$ followed by the contraction mapping} $\sE^{*}_{x} \otimes \sE_{x} \to \C$. 
    \newline 

    The particular case of the preceding expression of $\square$ when $\rrm, \lambda = 0$ is called the \textbf{d'Alembertian} or the \textbf{relativistic wave operator} which is locally given by 
    \begin{equation} \label{eq: d_Alembertian_coordinate} 
		\square = - \fg^{ij} \connectionLC_{i} \connectionLC_{j}, 
    	\quad i,j = 1, \ldots, d. 
    \end{equation}
    In case $\spacetime := (\R \times \varSigma, \rd t^{2} - \fh) $ is a ultrastatic manifold then 
    \begin{equation} \label{eq: def_d_Alembertian_Riemann_mf}
    	\square = - \frac{\partial^{2}}{\partial t^{2}} - \Delta  
    \end{equation}
	is known as the \textbf{wave operator} of a Riemannian manifold $(\varSigma, \fh)$ where $t$ is the natural coordinate on $\R$ and $\Delta$ is the Laplace-Beltrami operator of $(\varSigma, \fh)$ in the geometers' convention: 
    \begin{equation} \label{eq: Laplacian_coordinate}
    	\Delta 
    	:= - \Div_{\fh} \circ \grad 
    	= - \frac{1}{\sqrt{\det \fh}} \frac{\partial}{\partial x^{i}} \Big( \fh^{ij} \sqrt{\det \fh} \frac{\partial}{\partial x^{j}} \Big), 
    	\quad 
    	i,j = 2, \ldots, d. 
    \end{equation}
\end{example}
%
%
%

\begin{example} \label{exm: Yamabe_op}
	Suppose that $\sE \to \sM$ is a trivial-$\C$-line bundle over a $d \geq 3$-dimensional Lorentzian manifold $\spacetime$ and that $\square$ resp. $\ric$ are the d'Alembertian resp. the Ricci scalar of $\spacetime$. 
	The \textbf{Yamabe operator} is defined by 
    (see e.g.~\cite[Sec. 3.1, Exm. 2]{Baum_AGAG_1996},~\cite[Def. 3.5.9]{Baer_EMS_2007})  
    \begin{equation}
        \mathrm{Y} := \square + \frac{d-2}{4 (d-1)} \ric : C^{\infty} (\sM) \to C^{\infty} (\sM).    
    \end{equation}
\end{example}
%
%
%

\begin{example} \label{exm: Hodge_d_Alembert_op}
	Let $\sE := \wedge^{k} \coTansM$ be the bundle of $k$-forms over an oriented Lorentzian manifold $\spacetime$. 
    The \textbf{Hodge-d'Alembert operator} is defined by 
    (see e.g.~\cite[Sec. 3.1, Exm. 3]{Baum_AGAG_1996},~\cite[Exm. 1.5.3]{Baer_EMS_2007})  
    \begin{equation}
        \Box := \rd^{*} \circ \rd + \rd \circ \rd^{*} : C^{\infty} \big( \sM; \wedge^{k} \coTansM \big) \to C^{\infty} \big( \sM; \wedge^{k} \coTansM \big),     
    \end{equation}
    where $\rd$ and $\rd^{*}$ are the exterior and the coexterior derivatives, respectively. 
\end{example}
%
%
%

\begin{example} \label{exm: connection_d_Alembertian}
	Let $(\sE \to \sM, \connectionE)$ be a vector bundle equipped with a connection $\connectionE$ and $(\coTansM \to \sM, \connectionLC)$ is the cotangent bundle endowed with the Levi-CIvita connection $\connectionLC$, over a Lorentzian manifold $\spacetime$. 
	Then, we have an induced connection 
    \begin{equation}
    	\connectionCoTanME := \connectionLC \otimes \one_{\sE} + \one_{\coTansM} \otimes \connectionE
    \end{equation}
    on $\coTansM \otimes \sE$. 
    The \textbf{connection d'Alembertian} is defined by the composition of the following three maps 
    (see e.g.~\cite[Exm. 1.5.2]{Baer_EMS_2007})
    \begin{equation} \label{eq: def_connection_d_Alembertian}
        \secsME 
        \stackrel{\connectionE}{\to} C^{\infty} (\sM; \coTansM \otimes \sE) 
        \xrightarrow{\connectionCoTanME} C^{\infty} (\sM; \coTansM \otimes \coTansM \otimes \sE)  
        \xrightarrow{\tr_{\fg} \otimes I} \secsME, 
    \end{equation}
    where $\tr_{\fg} : C^{\infty} (\sM; \coTansM \otimes \coTansM) \to C^{\infty} (\sM, \R)$ denotes the \textbf{metric trace}: 
    \begin{equation} \label{eq: def_metric_trace}
        \tr_{\fg} \big( \xxi \otimes (x, \eta) \big) := \fg_{x}^{-1} (\xi, \eta) 
    \end{equation}
	for all $\xi, \eta \in \coTan_{x} \sM$. 
	The minus of a connection d'Alembertian, 
    \begin{equation}
    	\square := - \tr_{\fg} \big( \connectionCoTanME \circ \connectionE \big)
    \end{equation}
    is a normally hyperbolic operator. 
    \newline 

    One observes that the scalar d'Alembertian~\eqref{eq: d_Alembertian_coordinate} $- \tr_{\fg} (\connectionLC \circ \rd u)$ is a special case of the connection d'Alembertian when $u$ is a section of any trivial $\C$-line bundle.  
    Note that the first two composition of~\eqref{eq: def_connection_d_Alembertian} 
    \begin{equation} \label{eq: def_Hessian}
    	\Hess := \connectionCoTanME \circ \connectionE : \secsME \to C^{\infty} (\sM; \coTansM \otimes \coTansM \otimes \sE) 
    \end{equation}
    is called the \textbf{Hessian} and the corresponding differential operator $\Hess_{X, Y} = \connectionE_{X} \connectionE_{Y} - \connectionE_{\connectionLC_{X} Y}$ is called the \textbf{second covariant derivative} on $\sE$, where $X, Y$ are any tangent vectors on $\sM$ 
    (see e.g.~\cite[p. 66]{Berline_Springer_2004}).  
\end{example}
%
%
%

\begin{example} \label{exm: Bochner_d_Alembertian}
	Let $\big( \sE \to \sM, (\cdot|\cdot), \connectionE \big)$ be a vector bundle over a Lorentzian manifold $\spacetime$, where $(\cdot|\cdot)$ and $\connectionE$ are a sesquilinear form and a connection on $\sE$, respectively.  
	The \textbf{Bochner-d'Alembertian}  
	(see e.g.~\cite[Def. 2.4]{Berline_Springer_2004}) 
    \begin{equation} \label{eq: def_Bochner_d_Alembertian}
         \square := \nabla^{{\ms \sE} *} \circ \connectionE : \secsME \to \secsME,  
    \end{equation}
    where $\nabla^{{\ms \sE} *}$ is the connection induced by $\connectionE$ on the dual bundle (Section~\ref{sec: adjoint}) $\sE^{*}$ of $\sE$, differs from the corresponding connection d'Alembertian only by a minus sign: $\nabla^{{\ms \sE} *} \circ \connectionE = - \tr_{\fg} \big( \connectionCoTanME \circ \connectionE \big)$. 
    In any orthonormal frame $\{ e_{i} \}$ of $\tangent \! \sM$, it is given by 
    \begin{equation} \label{eq: Bochner_d_Alembertian_coordinate} 
    	\square = - \fg^{ij} \Big( \connectionE_{e_{i}} \connectionE_{e_{j}} - \Gamma_{\ ij}^{\LC \, k} \connectionE_{e_{k}} \Big), 
    	\quad 
    	i,j,k = 1, \ldots, d := \dim \sM, 
    \end{equation}
    where $\Gamma^{\LC}$ is the Levi-Civita connection $1$-form with respect to $e_{i}$: $\connectionLC_{e_{i}} e_{j} = \Gamma^{\LC \, k}_{\ ij} e_{k}$. 
\end{example}
%
%
%

There are some operators which are not normally hyperbolic yet can be related with a normally hyperbolic one with additional constraints.  

%
%
%
\begin{example} \label{exm: Proca_op}
	The \textbf{Proca operator} acting on the space of smooth covectors on a Lorentzian manifold  $\spacetime$ is defined as   
    \begin{equation}
        \mathrm{P} := \rd^{*} \circ \rd + \rrm^{2} : C^{\infty} (\sM; \coTansM) \to C^{\infty} (\sM; \coTansM), 
        \quad \mathrm{m} \in \dot{\R},  
    \end{equation}
	which is \textit{not} normally hyperbolic. 
    Nevertheless, this is equivalent to the normally hyperbolic operator $\rd^{*} \circ \rd + \rd \circ \rd + \rrm^{2}$ whenever the \textbf{Lorenz constraint}  $\rd^{*} \xi = 0$ is imposed.  
\end{example}
%
%
%

It is evident from the preceding examples that $\square$ differs from some d'Alembertian only by a smooth term. 
This is actually true in general. 

%
%
%
\begin{remark} \label{rem: Weitzenboeck_connection}
	By the Weitzenb\"{o}ck formula, given a normally hyperbolic operator $\square$ on any vector bundle $\sE \to \sM$ over a Lorentzian manifold $\spacetime$, there exists a \textit{unique} bundle connection $\connectionE$, called the \textbf{Weitzenb\"{o}ck connection} and a \textit{unique} potential $V \in C^\infty(\sM; \End \sE)$ such that 
    (see e.g.~\cite[Prop. 3.1]{Baum_AGAG_1996},~\cite[Lem. 1.5.5, 1.5.6]{Baer_EMS_2007}): 
    \begin{equation} \label{eq: Weitzenboeck_connection}
        \square = - \tr_{\fg} \big( \connectionCoTanME \circ \connectionE \big) + V.  
    \end{equation}
	Thus, in an orthonormal tangent frame $\{ \partial_{i} := \partial / \partial x^{i} \}$ and a bundle frame $\{ \rE_{r} \}$, a straightforward calculation using~\eqref{eq: Bochner_d_Alembertian_coordinate} and $\connectionE_{i} = \partial_{i} + \Gamma_{i}$,  $\Gamma$ being the Weitzenb\"{o}ck-connection $1$-form with respect to $\partial_{i}$: $\connectionE_{i} \rE_{r} = \Gamma_{ir}^{r'} \rE_{r'}$ entails 
    \begin{equation} \label{eq: NHOp_coordinate} 
    	(\square u)^{r}  
    	= 
    	\fg^{ij} \biggr(- \frac{\partial^{2} u^{r}}{\partial x^{i} \partial x^{j}} 
    	- 2 \Gamma_{ir'}^{r} \frac{\partial u^{r'}}{\partial x^{j}} + \Gamma^{\LC \, k}_{\ ij} \frac{\partial u^{r}}{\partial x^{k}} 
    	- \Big( \frac{\partial \Gamma_{jr'}^{r}}{\partial x^{i}} - \Gamma^{r}_{ir''} \Gamma^{r''}_{jr'} + \Gamma^{\LC, k}_{\ ij} \Gamma^{r}_{kr'} \Big) u^{r'} \bigg) + V u^{r},  
    \end{equation}
    where $i,j, k = 1, \ldots, d := \dim \sM$ and $r,r',r'' = 1, \ldots, \rk \sE$, and $\Gamma^{\LC}$ is the Levi-Civita connection $1$-form with respect to $\partial_{i}$. 
    Therefore, the subprincipal symbol of any normally hyperbolic operator is given by   
    \begin{eqnarray} \label{eq: subprincipal_symbol_NHOp}
    	\subSymb{\square} \xxi 
    	& \stackrel{\eqref{eq: subprincipal_symbol_bundle_Lor_mf}}{=} & 
    	- 2 \ri \fg^{ij} \Gamma_{i} \xi_{j} + \ri \fg^{ij} \Gamma^{\LC \, k}_{\ ij} \xi_{k} + \frac{\mathrm{i}}{2} \frac{\partial^{2} (\fg^{\mu \nu} \xi_{\mu} \xi_{\nu})}{\partial x^{i} \partial \xi_{i}} + \frac{\mathrm{i}}{2} \Gamma^{\LC \, j}_{\ ji} \frac{\partial \fg^{\mu \nu} \xi_{\mu} \xi_{\nu}}{\partial \xi_{i}} 
    	\nonumber \\ 
    	& = & 
    	- 2 \ri \fg^{ij} \Gamma_{i} \xi_{j} + \ri \Big( \fg^{ij} \Gamma^{\LC \, k}_{ij} + \frac{\partial \fg^{ik}}{\partial x^{i}} + \Gamma^{\LC \, j}_{\ ji} \fg^{ik} \Big) \xi_{k} 
    	\nonumber \\ 
    	& = & 
    	- 2 \ri \fg^{ij} \Gamma_{i} \xi_{j}, 
    \end{eqnarray}
    because of the identities  
    (see e.g.~\cite[$(3.4.9)$]{Wald_Chicago_1984}) 
    \begin{equation}
    	\fg^{ij} \Gamma^{\LC \, k}_{\ ij} = - \frac{1}{\sqrt{|\det \fg|}} \frac{\partial (\sqrt{|\det \fg|} \fg^{ik})}{\partial x^{i}}, 
    	\quad 
    	\Gamma^{\LC \, j}_{\ ji} = \frac{1}{2} \frac{1}{|\det \fg|} \frac{\partial |\det \fg|}{\partial x^{i}}.  
    \end{equation}

    If $\sE$ carries a bundle metric $(\cdot|\cdot)$ such that $\square$ is symmetric then $\connectionE$ is compatible with $(\cdot|\cdot)$ as well. 
    We do not, however, assume this unless stated otherwise. 
\end{remark}
%
%
%
%
%
%
%
%
%
%
\subsection{Green's operators}
\label{sec: Green_op}
\begin{definition} \label{def: Green_op}
	Let $\sE \to \sM$ be a vector bundle over a time-oriented Lorentzian manifold $\spacetime$ and $L$ a differential operator on $\sE$.  
	A \textbf{Green's operator} $G^{\times}$ for $L$ is a linear mapping 
	(see e.g.~\cite[Def. 3.4.1]{Baer_EMS_2007}) 
	\begin{equation}
		G^{\times} : \comSecsME \to \secsME ~|~ 
		L \circ G^{\times} = I, ~
		G^{\times} \circ (L \upharpoonright \comSecsME) = I. 
	\end{equation}
	A Green's operator is called 
	\begin{itemize}
		\item 
		\textbf{retarded} $\GreenOpRet$ if $\supp (\GreenOpRet u) \subset J^{+} (\supp u)$ and 
		\item 
		\textbf{advanced} $\GreenOpAdv$ if $\supp (\GreenOpAdv u) \subset J^{-} (\supp u)$, 
	\end{itemize}
	where $J^{\pm} (K) := \eqref{eq: causal_future_past}$ denotes the causal future/past of a set $K \subset \sM$. 
	The \textbf{Jordan-Lichnerowicz-Pauli} or the \textbf{causal} Green's operator is defined as the antisymmetric combination  
	\begin{equation}
		G := \GreenOpRet - \GreenOpAdv. 
	\end{equation}
\end{definition}
%
%
%

In physics, Green's operators are usually called \textit{propagators} and we will use this terminology as well.  
Propagators are closely related to the concept of fundamental solutions. 

%
%
%
\begin{definition} \label{def: fundamental_solution} 
	As in the terminologies of Definition~\ref{def: Green_op}, a \textbf{fundamental solution} of $L$ at $x \in \sM$ is an $\sE_{x}^{*}$-valued distribution on $\sE$ such that  
	(see e.g.~\cite[Def. 2.1.1]{Baer_EMS_2007},~\cite[Def. 3.1]{Guenther_AP_1988}) 
	\begin{equation}
		F_{x}^{\times} : \comSecsMEStar \to \sE_{x}^{*}, ~\phi \mapsto 
		F_{x}^{\times} (L^{*} \phi) = \phi (x),  
	\end{equation}
	where $L^{*}$ is the formal adjoint of $L$. 
	In other words, $F_{x}^{\times} \in \cD' \big( \sM; \sE; \sE_{x}^{*} \big)$ such that 
    \begin{equation}
    	L F_{x}^{\times} = \delta_{x}. 
    \end{equation}
    A fundamental solution is called \textbf{retarded} $\fundaSolRet$ if $\supp \fundaSolRet_{x} \subset J^{+} (x)$ and \textbf{advanced} $\fundaSolAdv$ if $\supp \fundaSolAdv_{x} \subset J^{-} (x)$. 
    The antisymmetric combination of the preceding two defines the \textbf{Jordan-Lichnerowicz-Pauli} or the \textbf{causal} fundamental solution $F := \fundaSolRet - \fundaSolAdv$. 
\end{definition}
%
%
%

To explain these concepts, let us look at the simplest non-trivial differential operator, the partial derivative.  

%
%
%
\begin{example} \label{ex: advanced_retarded_fundamental_sol_partial_derivative}
	Let $x \in \Rd$ and we write it as $x = (x^{1}, x')$ where $x^{1} \in \R$ and $x' = (x^{2}, \ldots, x^{d}) \in \R^{d-1}$. 
	Since the 	distributional derivative of 
	Heaviside step function\footnote{We 
        recall that it is defined as $\Theta_{x_{0}} (x) := 1$ for $x > x_{0}$ and $\varTheta_{x_{0}} (x) := 0$ for $x < x_{0}$.
        } 
	$\Theta_{x_{0}}$ is the Dirac delta distribution $\delta_{x_{0}}$ concentrated at $x_{0}$, the advanced, retarded, and causal fundamental solutions of $\rD_{1} := - \ri \partial / \partial x^{1}$ are 
    \begin{subequations}
    	\begin{eqnarray}
    		\fundaSolAdv_{1} (x) & = & - \ri \Theta (- x^{1}) \otimes \delta (x'), 
            \\ 
            \fundaSolRet_{1} (x) & = & - \ri \Theta (x^{1}) \otimes \delta (x'), 
            \\ 
            F_{1} (x') & = & - \ri \delta (x'). 
    	\end{eqnarray}
    \end{subequations}
	Their Schwartz kernels are given by  
	\begin{subequations}
		\begin{eqnarray}
			\fundaSolAdvKernel_{1}  (x - y)
			& = & 
			- \ri \Theta (y^{1} - x^{1}) \otimes \delta (x' - y'), 
			\\ 
			\fundaSolRetKernel_{1} (x - y)
			& = & 
			- \ri \Theta (x^{1} - y^{1}) \otimes \delta (x' - y'), 
			\\ 
			\fF_{1} (x' - y')
			& = & 
			- \ri \delta (x' - y') 
			= 
			\frac{- \ri}{(2 \pi)^{d-1}} \int_{\R^{d-1}} \re^{\ri (x^{2} - y^{2}) \theta_{2} + \ldots + \ri (x^{d} - y^{d}) \theta_{d}} \, \rd \theta_{2} \ldots \rd \theta_{d}. \quad 
		\end{eqnarray}
	\end{subequations}
\end{example}
%
%
%

One reads off the following facts from the Fourier integral representations of the preceding fundamental solutions.   

%
%
%
\begin{proposition} \label{prop: causal_propagator_partial_derivative_FIO}
    As in the terminologies of Example~\ref{ex: advanced_retarded_fundamental_sol_partial_derivative}, let $\chi$ be a smooth function on $\Rd \times \Rd$ vanishing near the diagonal. 
    Then~\cite[Prop. 6.1.2]{Duistermaat_ActaMath_1972}  
    (see also~\cite[Prop. 26.1.2]{Hoermander_Springer_2009})
    \begin{eqnarray} 
        && 
        \WF'  \mathsf{F}_{1}^{\adv, \ret} 
        = 
        \varDelta \, \dotCoTan \Rd \bigcup C_{1}^{\adv, \ret},  
        \label{eq: Hoermander_Prop_26_1_2_i}
        \\ 
    	&&	
        C_{1}^{\adv, \ret} := \{ \xxiyeta \in C_{1} \, | \, x^{1} \lessgtr y^{1} \},  
        \\ 
        && 
        C_{1} := \{ \xxiyeta \in \coTan \Rd \times \coTan \Rd \, | \, x' = y', \xi' = \eta' \neq 0, \xi_{1} = 0 = \eta_{1} \}, 
        \\ 
        && 
        \fF_{1}, \chi \mathsf{F}_{1}^{\adv, \ret} \in I^{- 1/2} (\Rd \times \Rd, C'_{1}) 
        \label{eq: Hoermander_Prop_26_1_2_ii}, 
    \end{eqnarray}
	where $\varDelta \, \dotCoTan \Rd$ is the diagonal in $\dotCoTan \Rd \times \dotCoTan \Rd$ and $I^{-1/2} (\Rd \times \Rd, C'_{1})$ is the Lagrangian distribution (Definition~\ref{def: Lagrangian_distribution_Euclidean}) on $\Rd \times \Rd$ associated with the canonical relation $C_{1}$. 
\end{proposition}
%
%
%

Suppose that $F_{x}^{\ret, \adv}$ is a family of retarded and advanced fundamental solutions for the adjoint $\square^{*}$ operator of a normally hyperbolic operator $\square$ on a vector bundle $\sE \to \sM$ over a time-oriented Lorentzian manifold $\spacetime$ and that $F_{x}^{\ret, \adv}$ depends smoothly on $x \in \sM$ in the sense that $x \mapsto F_{x}^{\ret, \adv} (\phi)$ is smooth for any $\phi \in \comSecsMEStar$. 
If it satisfies the differential equation $\square (F_{\cdot}^{\ret, \adv} \phi) = \phi$ then 
(see e.g.~\cite[Prop. 3.4.2]{Baer_EMS_2007}) 
\begin{equation}
    (G^{\adv, \ret} \phi) (x) := F_{x}^{\ret, \adv} (\phi) 
\end{equation}
defines advanced resp. retarded Green's operators for $\square$. 

%
%
%
\begin{remark} \label{rem: exist_unique_advanced_retarded_Green_op_NHOp}
	It is a classical result that any normally hyperbolic operator $\square$ on a vector bundle $\sE \to \sM$ over a globally hyperbolic spacetime $\spacetime$ admits \textit{unique retarded} and \textit{advanced} Green's operators  
    (see e.g.~\cite[Cor. 3.4.3]{Baer_EMS_2007},~\cite[Prop. 4.1, Rem. 4.3 b]{Guenther_AP_1988}):   
    \begin{equation}
    	G^{\adv, \ret} : \comSecsME \to C_{\mathrm{sc}}^{\infty} (\sM; \sE), 
    \end{equation}
	where $C_{\mathrm{sc}}^{\infty} (\sM; \sE)$ is the space of spatially compact smooth sections of $\sE$ (Section~\ref{sec: convention}). 
	\newline 

    Furthermore, the Green's operators $(G^{\adv, \ret})^{*}$ for the formal adjoint $\square^{*}$ of $\square$ are related to $G^{\adv, \ret}$ by 
    (see e.g.~\cite[Lem. 3.4.4]{Baer_EMS_2007})
    \begin{equation}
        \big( (G^{\adv, \ret})^{*} \phi \big) (u) = \phi \big( G^{\ret, \adv} u \big) 
    \end{equation}
    for any $\phi \in \comSecsMEStar$ and $u \in \comSecsME$. 
	Hence, the causal propagator
	\begin{equation} \label{eq: def_causal_propagator}
		G := G^{\ret} - G^{\adv} : \comSecsME \to C_{\mathrm{sc}}^{\infty} (\sM; \sE)
	\end{equation}
	satisfies $(G^{*} \phi) (u) = - \phi (G u)$. 
\end{remark}
%
%
%

A natural question in this context is, are there \textit{more} Green's operators for a given $\square$?   
In what follows, the answer leads to a special kind of propagators which are uniquely characterised in terms of their wavefront sets instead of their supports, in contrast to the advanced and retarded propagators. 
%
%
%
%
%
%
%
%
%
%
\subsection{Parametrices} 
\label{sec: parametrix_NHOp}
In order to investigate the preceding question, it is useful to introduce a weaker notion of Green's operator, as defined below. 

%
%
%
\begin{definition} \label{def: parametrix_PsiDO}
	Let $\sE, \halfDen \to M$ be a vector bundle resp. the half-density bundle over a manifold $M$ and $P$ a properly supported pseudodifferential operator (Definitions~\ref{def: PsiDO} and~\ref{def: psPsiDO_Euclidean}) on $\halfDen$-valued sections of $\sE$. 
    A \textbf{right parametrix} of $P$ is a continuous operator 
    \begin{equation}
    	\parametrixRight : \comSecE \to \secE \quad | \quad P \parametrixRight = I + R, 
    \end{equation}
    where $I$ is the identity pseudodifferential operator and $R$ is a smoothing operator. 
    Similarly, a \textbf{left parametrix} is defined by $\parametrixLeft P = I + L$ where $L$ is a smoothing operator. 
	One says $E$ is a \textbf{parametrix} of $P$ when it is both a right and a left parametrix. 
\end{definition}
%
%
%

Now, let $\sE \to \sM$ be a vector bundle over a globally hyperbolic spacetime $\spacetime$ and let $E$ be a parametrix of a normally hyperbolic operator $\square$ on $\sE$: 
\begin{equation}
    \square E \equiv I, 
	\quad 
	E \big( \square \upharpoonright \comSecsME \big) \equiv I,  
\end{equation}
where $\equiv$ means modulo smoothing operators. 
Remark~\ref{rem: exist_unique_advanced_retarded_Green_op_NHOp} entails that $\square$ admits at least two parametrices: the advanced $\parametrixAdv$ and the retarded $\parametrixRet$. 
We would like to know whether there are more parametrices.  
This question has been investigated in a great detail by 
Duistermaat and H\"{o}rmander~\cite[Sec. 6.5, 6.6]{Duistermaat_ActaMath_1972} 
for any scalar pseudodifferential operator of real principal type (Definition~\ref{def: real_principal_type_PsiDO_mf}). 
To translate their analysis in the present context, we revisit Example~\ref{exm: NHOp_real_principal_type} and consider the following relation.  

%
%
%
\begin{definition} \label{def: geodesic_relation}
	Let $\coLightBun \to \sM$ be the bundle of lightlike covectors over a globally hyperbolic spacetime $\spacetime$. 
	The \textbf{geodesic relation} $C$ is defined as the set of all covector pairs in the product manifold $\coLightBun \times \coLightBun$ lying on the same geodesic on $\coLightBun$~\cite[$(6.5.2)$]{Duistermaat_ActaMath_1972} 
    (see also~\cite[$(26.1.8)$]{Hoermander_Springer_2009}):   
	\begin{equation} \label{eq: def_geodesic_relation} 
        C := \big\{ (x, \xi; y, \eta) \in \coLightBun \times \coLightBun \,|\, \exists ! s \in \R : \xxi = \varPhi_{s} (y, \eta) \big\}.  
    \end{equation}
	The \textbf{forward} resp. \textbf{backward} \textbf{geodesic relations} $C^{\pm}$ are defined as the set of all $\xxiyeta$ such that $\xxi$ lies after resp. before $\yeta$ on a geodesic:   
    \begin{equation} \label{eq: def_forward_backward_geodesic_relation} 
        C^{\pm} := \big\{ \xxiyeta \in \coLightBun \times \coLightBun \,|\, \exists ! s \in \R_{\gtrless 0} : \xxi = \varPhi_{s} (y, \eta) \big\}.    
    \end{equation}
    Here, $\varPhi_{s}$ is the geodesic flow on the cotangent bundle and $s \in \R$ is the flow-parameter. 
\end{definition}
%
%
%

\begin{remark} \label{rem: geodesic_relation_canonical_relation}
	The geodesic relation is a conic Lagrangian submanifold (see Definition~\ref{def: Lagrangian_submf}) of $\coLightBun \times \coLightBun$ which is closed in $\dotCoTan_{0} (\sM \times \sM)$ and consequently a homogeneous canonical relation (see Definition~\ref{def: canonical_relation}) from $\coLightBun$ to $\coLightBun$~\cite[Prop. 6.5.2]{Duistermaat_ActaMath_1972} 
    (see also~\cite[Prop. 26.1.13]{Hoermander_Springer_2009}).
\end{remark}
%
%
%

With these terminologies in hand, the wavefront sets of the Schwartz kernels $\parametrixKernelAdv, \parametrixKernelRet$ of $\GreenOpAdv, \GreenOpRet$ are given by 
\begin{eqnarray} \label{eq: def_advanced_retarded_canonical_relation_NHOp}
	\WF' \parametrixKernelAdv
    & \subset &
    \varDelta \dotCoTansM \cup C^{\adv}, 
    \quad 
    C^{\adv} := \sbracket{(x, \xi, y, \eta) \in C \,|\, x \in J^{-} (y)}, 
    \\ 
    \WF' \parametrixKernelRet  
    & \subset & 
    \varDelta \dotCoTansM \cup C^{\ret}, 
    \quad 
    C^{\ret} := \sbracket{(x, \xi, y, \eta) \in C \,|\, x \in J^{+} (y)}. 
\end{eqnarray}
Duistermaat-H\"{o}rmander's~\cite[p. 218]{Duistermaat_ActaMath_1972}  
groundbreaking analysis shows that, if $d \geq 3$ then there are two more parametrices of $\square$ whose wavefront sets are off-diagonally given by the forward and backward geodesic relations.  
These parametrices are called the \textbf{distinguished parametrices} in pure mathematics, whereas they are known as the \textbf{Feynman parametrices} in the theoretical physics literature. 

%
%
%
\begin{definition} \label{def: Feynman_parametrix}
	Let $\sE \to \sM$ be a vector bundle over a globally hyperbolic spacetime $\spacetime$ and $\square$ a normally hyperbolic operator on $\sE$.  
	The \textbf{Feynman} $\parametrixFeyn$ and the \textbf{anti-Feynman} $\parametrixAntiFeyn$ parametrices of $\square$ are those parametrices of $\square$ whose Schwartz kernels satisfy~\cite[p. 229]{Duistermaat_ActaMath_1972} 
	\begin{equation} \label{eq: def_Feynman_parametrix}
		\WFPrime{\parametrixKernelFeyn} \subset \varDelta \, \dotCoTansM \cup C^{+}, 
		\qquad 
		\WFPrime{\parametrixKernelAntiFeyn} \subset \varDelta \, \dotCoTansM \cup C^{-}, 
	\end{equation}
	where $\varDelta \, \dotCoTansM$ is the diagonal to the punctured cotangent manifold $\dotCoTansM \times \dotCoTansM$ and $C^{\pm}$ are the forward (resp. backward) geodesic relations. 
\end{definition}
%
%
%

Therefore, given $\square$ on $d \geq 3$, one has at \textit{most} $4$ independent parametrices $\parametrixAdv, \parametrixRet$, $\parametrixFeyn, \parametrixAntiFeyn$. 
Evidently, the hindmost two are characterised by their wavefront sets rather supports. 

%
%
%
\begin{remark} \label{rem: sixteen_parametrix_dim_two_NHOp}
    There are $16$ parametrices of $\square$ in $d=2$ because $\coLightBun$ has $4$ connected components, whereas $\dotCoTan_{0, \pm} \sM$ are the two connected components of $\coLightBun$ in $d \geq 3$.  
\end{remark}

%
%
%
\begin{center}
    \begin{tikzpicture}
        \fill[
            left color=red!70!white,
            right color=red!70!white,
            middle color=red!20!white,
            shading=axis,
            opacity=1.0
        ]
        (-4, 0) -- (-2, -2) -- (0, 0) arc (360: 180: 2cm and 0.40cm);
        \fill[
            left color=green!70!white,
            right color=green!70!white,
            middle color=green!20!white,
            shading=axis,
            opacity=1.0
        ]
        (-4, -4) -- (-2, -2) -- (0, -4) arc (360: 180: 2cm and 0.40cm);
        \draw[color=red, thick] (0, 0) arc (360: 0: 2cm and 0.40cm);
        \draw[red, thick] (-4, 0) -- (-2, -2) -- (0, 0); 
        \draw[color=green!50!black, thick] (0, -4) arc (360: 180: 2cm and 0.40cm);
        \draw[color=green!50!black, thick, dashed] (-4, -4) arc (180: 0: 2cm and 0.40cm);
        \draw[black!50!green, thick] (-4, -4) -- (-2, -2) -- (0, -4); 
        \draw[-latex, thick] (-1.75, -1.5) -- (-1.25, -0.75);
        \filldraw (-1.75, -1.5) circle (1.5pt); 
        \draw[-latex, thick] (-2.5, -1.25) -- (-3, -0.5);
        \filldraw (-2.5, -1.25) circle (1.5pt);
        \draw[-latex, thick] (-1.25, -3.0) -- (-0.75, -3.75); 
        \filldraw (-1.25, -3.0) circle (1.5pt);
        \draw[-latex, thick] (-2.5, -3.25) -- (-3.0, -4.0); 
        \filldraw (-2.5, -3.25) circle (1.5pt);
        \draw[-latex, thick, blue] (-2, -2) -- (-2.5, -2.75); 
        \filldraw[blue] (-2, -2) circle (1.5pt);
        %
        \node[right] at (-1.25, -0.75) {$\xi$};
        \node[left] at (-1.7, -1.5) {$x$}; 
        \node[left] at (-2, -2) {\textcolor{blue}{$y$}};
        \node[right] at (-2.45, -2.75) {\textcolor{blue}{$\eta$}}; 
     \end{tikzpicture} 
     \hspace*{2cm}
     \begin{tikzpicture}
        \fill[
            left color=red!70!white,
            right color=red!70!white,
            middle color=red!20!white,
            shading=axis,
            opacity=1.0
        ]
        (-4, 0) -- (-2, -2) -- (0, 0) arc (360: 180: 2cm and 0.40cm);
        \draw[color=red, thick] (0, 0) arc (360: 0: 2cm and 0.40cm);
        \draw[red, thick] (-4, 0) -- (-2, -2) -- (0, 0); 
         \filldraw (-1.46, -1.15) circle (1.5pt); 
         \draw[>=latex, thick] (-2.15, -1.75) edge[<->] (-3, -0.5);
         \filldraw (-2.55, -1.15) circle (1.5pt); 
         \draw[>=latex, thick] (-1.9,-1.75) edge[<->] (-1,-0.5);
        \draw[-latex, thick, white] (-1.25, -3.0) -- (-0.75, -3.75); 
        \draw[-latex, thick, white] (-1.25, -3.0) -- (-1.75, -2.25); 
     \end{tikzpicture}
     \newline 
     \begin{tikzpicture}
        \fill[
            left color=green!70!white,
            right color=green!70!white,
            middle color=green!20!white,
            shading=axis,
            opacity=1.0
        ]
        (-4, -4) -- (-2, -2) -- (0, -4) arc (360: 180: 2cm and 0.40cm);
        \draw[color=green!50!black, thick] (0, -4) arc (360: 180: 2cm and 0.40cm);
        \draw[color=green!50!black, thick, dashed] (-4, -4) arc (180: 0: 2cm and 0.40cm);
        \draw[black!50!green, thick] (-4, -4) -- (-2, -2) -- (0, -4); 
        \draw[-latex, thick] (-1.25, -3.0) -- (-0.75, -3.75); 
        \draw[-latex, thick] (-1.25, -3.0) -- (-1.75, -2.25); 
        \filldraw (-1.25, -3.0) circle (1.5pt);
        \draw[-latex, thick] (-2.5, -3.25) -- (-3.0, -4.0); 
        \filldraw (-2.5, -3.25) circle (1.5pt); 
        \draw[-latex, thick] (-2.5, -3.25) -- (-2.0, -2.45); 
        %
     \end{tikzpicture}
     \hspace*{2cm}
     \begin{tikzpicture}
        \fill[
            left color=red!70!white,
            right color=red!70!white,
            middle color=red!20!white,
            shading=axis,
            opacity=1.0
        ]
        (-4, 0) -- (-2, -2) -- (0, 0) arc (360: 180: 2cm and 0.40cm);
        \fill[
            left color=green!70!white,
            right color=green!70!white,
            middle color=green!20!white,
            shading=axis,
            opacity=1.0
        ]
        (-4, -4) -- (-2, -2) -- (0, -4) arc (360: 180: 2cm and 0.40cm);
        \draw[color=red, thick] (0, 0) arc (360: 0: 2cm and 0.40cm);
        \draw[red, thick] (-4, 0) -- (-2, -2) -- (0, 0); 
        \draw[color=green!50!black, thick] (0, -4) arc (360: 180: 2cm and 0.40cm);
        \draw[color=green!50!black, thick, dashed] (-4, -4) arc (180: 0: 2cm and 0.40cm);
        \draw[black!50!green, thick] (-4, -4) -- (-2, -2) -- (0, -4); 
        %
        \draw[-latex, thick] (-1.25, -0.75) -- (-1.75, -1.5);
        \filldraw (-1.25, -0.75) circle (1.5pt); 
        \draw[-latex, thick] (-3, -0.5) -- (-2.5, -1.25);
        \filldraw (-3, -0.5) circle (1.5pt);
        \draw[-latex, thick] (-0.75, -3.75) -- (-1.25, -3.0); 
        \filldraw (-0.75, -3.75) circle (1.5pt);
        \draw[-latex, thick] (-3.0, -4.0) -- (-2.5, -3.25); 
        \filldraw (-3.0, -4.0) circle (1.5pt);
        \node at (0, -3) {$\textcolor{green!50!black}{\dotCoTan_{0, -} \sM}$};
         \node at (0, -1) {$\textcolor{red}{\dotCoTan_{0, +} \sM}$}; 
     \end{tikzpicture}
     
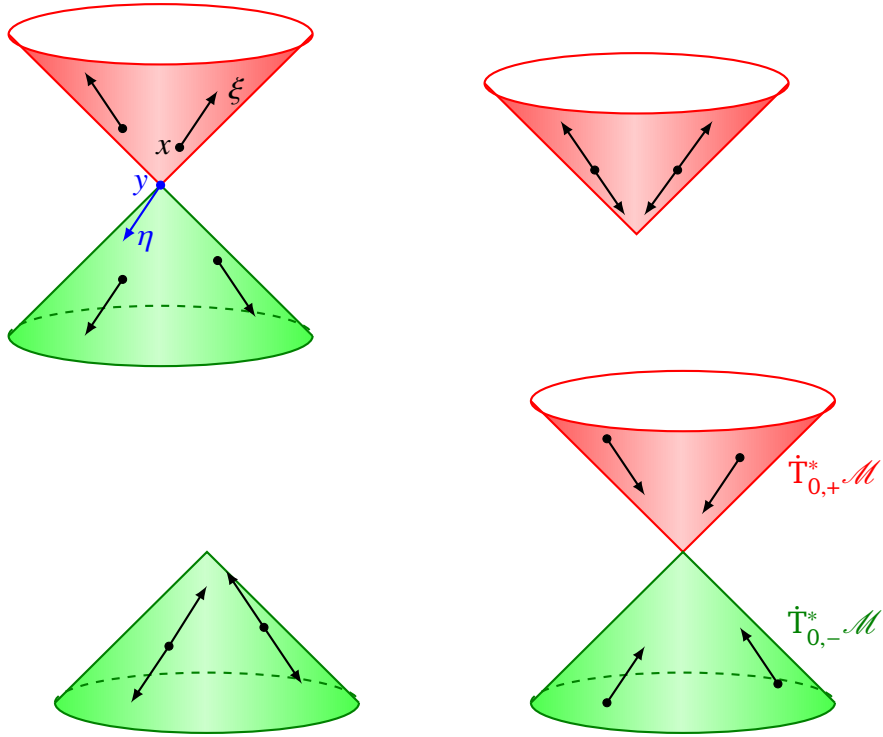
\captionof{figure}[$C^{+}, C^{\ret}, C^{\adv}, C^{-}$ in $d \geq 3$]{Schematic 
         illustration of $C^{+}, C^{\ret}, C^{\adv}, C^{-}$ in spacetime dimensions $d \geq 3$. 
         In these diagrams $y$ is fixed at the vertices of the cones and only the directions of $\xi$ are shown while those of $\eta$ being left undisplayed for simplicity except for the top left one  
        (redrawn from~\cite[Fig. 1]{Radzikowski_CMP_1996}).}
 \end{center}
%
%
%
 
%
%
%
%
%
%
%
%
%
%
\section{Feynman propagators for a NHO} 
\subsection{Primary results}
\label{sec: result_Feynman_propagator_NHOp}
To begin with, we provide a bundle version of the classic result by 
Duistermaat-H\"{o}rmander ~\cite[Thm. 6.5.3]{Duistermaat_ActaMath_1972} 
(see also~\cite[Thm. 26.1.14]{Hoermander_Springer_2009}) 
on the existence and uniqueness of Feynman parametrices of a normally hyperbolic operator. 
Recall that $H_{\mathrm{loc}}^{s} (M; \sE \otimes \halfDen)$ is the space of $\halfDen$-valued sections of $\sE$ that are locally in the Sobolev space $H^{s}$ with respect to any smooth chart and any smooth bundle chart. 
The space of sections in $H_{\mathrm{loc}}^{s} (M; \sE \otimes \halfDen)$ of compact support is denoted by $H_{\mathrm{c}}^{s} (M; \sE \otimes \halfDen)$. 
As usual, the space $H_{\mathrm{loc}}^{s} (M; \sE \otimes \halfDen)$ is equipped with the locally convex topology of convergence locally in $H^{s} (M; \sE \otimes \halfDen)$. 
The space $H_{\mathrm{c}}^{s} (M; \sE \otimes \halfDen)$ is the union $\bigcup_{K \Subset M} H_{\mathrm{c}}^{s} (K; \sE \otimes \halfDen)$ and it is equipped with the inductive limit topology, where the union runs over all compact subset $K$ of $M$. 
For details, we refer, for instance~\cite[App. B1]{Hoermander_Springer_2007}. 

%
%
%
\begin{theorem}[Existence and uniqueness of Feynman parametrices]  
\label{thm: exist_unique_Feynman_parametrix_NHOp}
	Let $\sE \to \sM$ be a vector bundle over a globally hyperbolic spacetime $\spacetime$ and $\square$ a normally hyperbolic operator on $\sE$. 
	Then, there exist 
	unique\footnote{In the sense of parametrices, i.e., modulo smoothing operators.}
	Feynman $\parametrixFeyn$ and anti-Feynman $\parametrixAntiFeyn$ parametrices (given by~\eqref{eq: Feynman_parametrices_NHOp}) of $\square$.  
	Moreover, for every $k \in \R$, $\parametrixFeyn$ and $\parametrixAntiFeyn$ extend to continuous maps from  $H_{\mathrm{c}}^{k} (\sM; \sE)$ to  $H_{\mathrm{loc}}^{k+1} (\sM; \sE)$ and the difference of their Schwartz kernels  
	\begin{equation} \label{eq: diff_Feyn_anti_Feyn_parametrix}
		\parametrixKernelFeyn - \parametrixKernelAntiFeyn \in I^{-3/2} \big( \sM \times \sM, C'; \Hom{\sE, \sE} \big),  
	\end{equation} 
    where $I^{-3/2} (\ldots)$ is the space of Lagrangian distributions (Definition~\ref{def: Lagrangian_distribution}) associated to the geodesic relation $C$ (Definition~\ref{def: geodesic_relation}). 
	Furthermore, $\parametrixKernelFeyn - \parametrixKernelAntiFeyn$ is non-characteristic (Definition~\ref{def: elliptic_FIO}) at every point of $C$.  
\end{theorem}
%
%
%

A special case of this result for the massive Klein-Gordon operator (Example~\ref{exm: covariant_Klein_Gordon_op}) was given by 
Radzikowski~\cite[Prop. 4.2-4.4]{Radzikowski_CMP_1996} 
as a direct consequence of~\cite[Thm. 6.5.3]{Duistermaat_ActaMath_1972}. 
Employing the distinguished global phase function approach of Fourier integrals operators~\cite{Laptev_CPAM_1994}, 
Capoferri \textit{et al}.~\cite[Thm. 5.2]{Capoferri_JMAA_2020} 
have constructed these parametrices for scalar wave operators with time-independent smooth potential in spatially compact globally hyperbolic ultrastatic spacetime.  
Recently, 
Lewandowski~\cite[Prop. 3.5]{Lewandowski_JMP_2022} 
has published a direct construction for a real vector bundle utilising the Hadamard series expansion along with the presentation of~\cite{Baer_EMS_2007}. 
In contrast to these literature, our proof is purely microlocal as in the original 
treatment~\cite[Thm. 6.5.3]{Duistermaat_ActaMath_1972}.
This only requires a bundle version of microlocalisation, as developed in due course (Theorems~\ref{thm: microlocalisation_P} and~\ref{thm: microlocalisation_NHOp})  which is along the lines of 
Dencker's~\cite{Dencker_JFA_1982} proof of propagation of the polarization sets albeit our presentation is more geometric. 
\newline 

The Feynman parametrix $\parametrixFeyn$ can be turned into a Green's operator $\GreenOpFeyn$ --- the Feynman propagator for $\square$. 
As elucidated in Appendix~\ref{ch: Hadamard_state}, a Feynman propagator is intimately related to a quantum state and hence this is equivalent to construct a Hadamard $2$-point distribution provided certain 
positivity\footnote{This is related to the positivity of Wightman distributions in quantum field theory.} 
property (as described below) is satisfied. 
This remarkable property was first observed by 
Duistermaat-H\"{o}rmander~\cite[Thm. 6.6.2]{Duistermaat_ActaMath_1972} 
for scalar operators.  
Namely, that one can modify $\parametrixFeyn$ by a smoothing operator $R$ such that
\begin{equation}
    \hat{W} + R := - \ri (\parametrixFeyn - \GreenOpAdv) + R
\end{equation}
is non-negative in the sense that $(\hat{\fW} + \fR) (\bar{u} \otimes u) \geq 0$ for any $u \in C_{\mathrm{c}}^{\infty} (\sM)$. 
This means, of course, that there exists a Feynman parametrix that has this property with $R = 0$. 
Thus, there exists a Feynman parametrix that satisfies $\hat{W} \geq 0$ in the above sense. 
Unlike functions, there is \textit{not} any notion of positive operators on a vector bundle unless the bundle is hermitian. 
Hence, we equip $\sE$ with a hermitian form and generalise this positivity property in the case of vector bundles. 

%
%
%
\begin{restatable}{proposition}{positivityFeynmanMinusAdvNHOp} 
\label{prop: positivity_Feynman_minus_adv_NHOp}
	Let $\sE \to \sM$ be a vector bundle over a globally hyperbolic spacetime $\spacetime$ and $\square$ a normally hyperbolic operator on $\sE$. 
	We assume that $\sE$ is endowed with a (non-degenerate) sesquilinear form $(\cdot|\cdot) \in C^{\infty} (\bar{\sE}^{*} \otimes \sE)$ such that $\square$ is symmetric. 
	Then, there exists a Feynman parametrix $\parametrixFeyn$ of $\square$ such that $\hat{W} := - \ri (\parametrixFeyn - \GreenOpAdv)$ is symmetric, where $\GreenOpAdv$ is the advanced propagator for $\square$. 
	Additionally, if $(\cdot|\cdot)$ is hermitian (positive-definite) then $\hat{W}$ can be chosen non-negative. 
\end{restatable}
%
%
%

Duistermaat-H\"{o}rmander have proven their version by deploying a refined microlocalisation of scalar pseudodifferential operators~\cite[Lem. 6.6.4]{Duistermaat_ActaMath_1972}. 
We provide such a microlocalisation for normally hyperbolic operators in~\eqref{eq: microlocalisation_NHOp_refined}. 
\newline 

Finally, we turn the Feynman parametrix $\parametrixFeyn$ into a Feynman propagator $\GreenOpFeyn$ utilising the well-posedness of Cauchy problem for $\square$ on a globally hyperbolic spacetime, which in turn gives the existence of Hadamard bisolutions of $\square$. 

%
%
%
\begin{restatable}[Existence of Feynman propagators]{theorem}{existenceFeynmanGreenOpNHOp}
\label{thm: existence_Feynman_propagator_NHOp}
	As in the setup of Proposition~\ref{prop: positivity_Feynman_minus_adv_NHOp}, there exists a Feynman propagator $\GreenOpFeyn$ for $\square$ on $\big( \sE \to \sM, (\cdot|\cdot) \big)$ such that $W := - \ri (\GreenOpFeyn - \GreenOpAdv)$ is symmetric with respect to the sesquilinear form $(\cdot|\cdot)$.   
	In addition, if $(\cdot|\cdot)$ is hermitian then $\GreenOpFeyn$ can be chosen such that $\fW (\bar{u}^{*} \otimes u) \geq 0$ for any $u \in \comSecsME$. 
	Thus, $\fW$ defines a Hadamard state. 
	Here $\secsME \ni u \mapsto \bar{u}^{*} \in C^{\infty} (\sM; \bar{\sE}^{*})$ is the fibrewise linear mapping induced by $(\cdot|\cdot)$.
\end{restatable}
%
%
%

Since 
Duistermaat-H\"{o}rmander~\cite[p. 229]{Duistermaat_ActaMath_1972} 
have considered a much wider class of operators --- scalar pseudodifferential operators of real principal type, the pivotal step of determining the appropriate smoothing operators required to obtain Feynman propagators from respective parametrices was, however, left open.
This indeterminacy can be fixed in various ways on special spacetimes even in the absence of the timelike Killing vector field. 
Such constructions have appeared in the literature on microlocal analysis. 
For example, 
Gell-Redman \textit{et al}.~\cite[Thm. 3.6]{GellRedman_CMP_2016} 
have treated the scalar wave operator in spaces with non-trapping Lorentzian scattering matrices.  
Vasy has constructed Feynman propagators by making assumptions on global dynamics for 
(i) scalar symmetric operators with real principal symbol on closed manifolds~\cite[Thm. 1]{Vasy_AHP_2017}, 
and in spaces of Lorentzian scattering matrices for 
(ii) wave operators in Melrose's b-pseudodifferential algebraic framework~\cite[Thm. 7]{Vasy_AHP_2017} 
and 
(iii) Klein-Gordon operators in Melrose's scattering pseudodifferential algebraic formalism~\cite[Thm. 10, 12]{Vasy_AHP_2017}.
His idea is to identify the appropriate spaces where these operators are invertible and then define the Feynman propagators as the inverse of those operators satisfying the required properties --- a generalisation of Feynman's original (see Example~\ref{ex: Feynman_propagator_Minkowski} and the comment afterward) "$\pm \ri \epsilon$" prescription. 
As a consequence, he has also achieved respective positivity property for Feynman parametrices~\cite[Cor. 5, 9, 11, 13]{Vasy_AHP_2017}. 
A special case of Theorem~\ref{thm: existence_Feynman_propagator_NHOp} for the Klein-Gordon operator minimally coupled to a static electromagnetic potential on a static spacetime has been proven by 
Derezi\'{n}ski and Siemssen~\cite[Thm. 7.7]{Derezinski_RMP_2018}. 
In the spirit of the limiting absorption principle, they have shown that the Feynman propagator can be considered as the boundary value of the resolvent of the Klein-Gordon operator.
Employing Proposition~\ref{prop: positivity_Feynman_minus_adv_NHOp} and ideas from~\cite[Sec. 3.3]{Gerard_CMP_2015}, 
Lewandowski~\cite[Thm. 4.3]{Lewandowski_JMP_2022} 
has recently given a construction of Hadamard states on a Riemannian vector bundle. 
\newline 

In order to proof these assertions, we need some preparation, as presented in the next two subsections. 
Some of the results of the following sections are interesting in a broader context, so we will work on a vector bundle $\sE \to M$ over a generic manifold $M$ for greater generality and come back to the context of a globally hyperbolic spacetime $\spacetime$ afterward. 
%
%
%
%
%
%
%
%
%
%
\subsection{Connection $1$-forms and subprincipal symbols}
\label{sec: connection_form_subprincipal_symbol}
Let $\sE, \halfDen \to M$ be a vector bundle resp. the half-density bundle, over a manifold $M$.  
This has been observed that the subprincipal symbol (Definition~\ref{def: subprincipal_symbol}) $\subSymb{P}$ of $P \in \PsiDOME, m \in \R$ with scalar principal symbol $p$ transforms like a 
partial\footnote{We 
    recall that the \textbf{partial connection} $\connectionE$ on a vector bundle $\sE \to M$ over a manifold $M$ with respect to a foliation $N$ of $M$ is defined by the covariant derivative $\connectionE_{X}$ satisfying the standard properties for any vector field $X$ on $N$ 
    (see e.g.~\cite[p. 24]{Kamber_Springer_1975}. 
    Note, not all partial connections arise as restrictions of some ``full'' connection, for instance, the Bott connection associated with any involutive distribution of a manifold.}-
\textit{connection $1$-form} under change of bundle charts 
(see e.g.~\cite[Prop. 3.1]{Jakobson_CMP_2007}). 
This is perhaps not surprising given that the subprincipal symbol appears as a constant term in transport equations. 
Connections that are naturally defined from transport equations have appeared first in the work of Dencker on propagation of polarisation sets. 
They are often referred to as the 
\textit{Dencker connections}~\cite[p. 365-366]{Dencker_JFA_1982} 
in the mathematical physics literature. 
For the Maxwell system, Dencker found that this connection equals the Levi-Civita connection. 
For the spin-Dirac operator (see Example~\ref{ex: spin_Dirac_op}), it was verified, for example, 
in~\cite{Hollands_adiabatic_CMP_2001} 
that the Dencker connection is indeed the spin-connection.
In this section, we will systematically explain the precise relation between a geometrically defined connection on $\sE$ and the parallel transporter induced by $\subSymb{P}$. 
As already Remarked~\ref{rem: subprincipal_symbol_bundle_chart_change}, $\subSymb{P}$ depends on the choice of trivialisations of $\sE$ and $\halfDen$.
We are going to see below that under a change of bundle frames, this matrix will transform like a connection $1$-form along the Hamiltonian vector field $X_{p}$ in the bundle $\pi^{*} \sE \to \dotCoTanM$ where $\pi : \dotCoTanM \to M$ is the punctured cotangent bundle of $M$. 
\newline 

To begin with, let us summarise the transformation properties of the subprincipal symbol. 
In case $P$ has a scalar principal symbol $p$ and $Q \in \PsiDO{m'}{M; \sE \otimes \halfDen}, m' \in \R$ properly supported with principal symbol $q$, we have the following multiplication formula~\cite[(1.4)]{Duistermaat_InventMath_1975}  
(see also, e.g.~\cite[Prop. 2.1.13]{Safarov_AMS_1997})
\begin{eqnarray} \label{eq: subprincipal_symbol_product}
	\subSymb{PQ} 
	& = & 
	\subSymb{\sum_{l} P^{r}_{l} Q^{l}_{k}} 
    =
	\sum_{l} \subSymb{P^{r}_{l} Q^{l}_{k}} 
	\stackrel{ \textrm{\cite[$(1.4)$]{Duistermaat_InventMath_1975}} }{=}  
	\sum_{l} \subSymb{P^{r}_{l}} q^{l}_{k} + p^{r}_{l} \subSymb{Q^{l}_{k}} + \frac{1}{2 \ri} \big\{ p^{r}_{l}, q^{l}_{k} \big\}   
 	\nonumber \\ 
 	& = & 
	\subSymb{P} \, q + p \, \subSymb{Q} + \frac{1}{2 \ri} \{ p, q \},   
\end{eqnarray}
where we have used the linearity of $\subSymb{}$ in the second line and $\{\cdot, \cdot\}$ is the Poisson bracket (see~\eqref{eq: def_Poisson_bracket}). 
Consequently, we also have~\cite[$(1.3)$]{Duistermaat_InventMath_1975} 
\begin{equation} \label{eq: subprincipal_symbol_power}
	\subSymb{P^{k}} = k \, p^{k-1} \, \subSymb{P}  
\end{equation}
for any $k \in \N$. 
In addition, if $P$ is elliptic (Definition~\ref{def: char_PsiDO}) and $E$ is a parametrix (Definition~\ref{def: parametrix_PsiDO}) of $P$ then   
\begin{equation} \label{eq: subprincipal_symbol_inverse}
	\subSymb{E} = - \, p^{-2} \, \subSymb{P},   
\end{equation}
and in this sense the above formula also holds for any negative integer $k$. 
\newline 

We will now show that the subprincipal symbol indeed has the claimed transformation property under a change of bundle charts. 
In fact, we show a microlocal version of this statement 
(see also~\cite{Hintz_JST_2017}). 

%
%
%
\begin{proposition} \label{prop: P_connection_transformation}
	Let $\halfDen \to M$ be the half-density bundle over a manifold $M$ and let $m, m' \in \R, N \in \N$. 
	Assume $P \in \PsiDO{m}{M, \matN \otimes \halfDen}$ having a scalar principal symbol $p$. 
	Suppose that $Q \in \PsiDO{m'}{M, \matN \otimes \halfDen}$ is non-characteristic at some $\xxiNot \in \dotCoTanM$ and that $E$ is a microlocal parametrix of $Q$ in an open conic neighbourhood $\cU$ of $\xxiNot$. 
	Then the subprincipal symbol of $EPQ$ is 
	\begin{equation}
		\forall \xxi \in \cU : \subSymb{EPQ} \xxi = (q^{-1} \subSymb{P} q) \xxi - \ri \big( q^{-1} X_{p} (q) \big) (x,\xi), 
	\end{equation}
	where $q$ is the principal symbol of $Q$ and $X_{p}$ is the Hamiltonian vector field generated by $p$.  
\end{proposition}
%
%
%
\begin{proof}
    Since $P$ has a scalar principal symbol, the order of the commutator $[P,Q]_{-}$ is $m+m'-1$ and the principal symbol of $EPQ$ is microlocally that of $P$.  
    These imply that the order of $E [P,Q]_{-}$ is $m-1$ and that $\xxiNot \notin \ES (EPQ - P - E [P,Q]_{-})$. 
 	Thus, deploying the order, 
 	\begin{equation}
 		\subSymb{EPQ} = \subSymb{P}  + \symb{E [P,Q]_{-}}
 	\end{equation}
    on $\cU$. 
    It remains to compute the principal symbol of $[P,Q]_{-}$. 
    Since the principal symbol of $P$ is scalar-valued, we obtain the commutator relation 
 	(see e.g.~\cite[$(14)$]{Jakobson_CMP_2007}) 
 	\begin{equation} \label{eq: symbol_commutator_PsiDO}
 		\symb{[P,Q]_{-}} = - \ri \{ p, q \} + [\subSymb{P}, q] _{-}, 
 	\end{equation}
 	as a consequence of the multiplication formula~\eqref{eq: subprincipal_symbol_product}. 
	Therefore, on $\cU$, we get 
	$
	\symb{Q^{-1} [P,Q]_{-}} 
	= q^{-1} \symb{[P, Q]_{-}} 
	= - \ri q^{-1} \{ p, q \} + q^{-1}  \subSymb{P} q - \subSymb{P}
	$ 
	and the result follows.
\end{proof}
%
%
%

As a simple consequence of the above product formula we also have  

%
%
%
\begin{proposition} \label{prop: P_orderonerescale}
	As in the terminologies of Proposition~\ref{prop: P_connection_transformation}, let $Q$ has a scalar principal symbol $q$ and vanishing subprincipal symbol.
	If $P$ has a scalar principal symbol $p$ and a subprincipal symbol $\subSymb{P}$, then 
	\begin{equation}
		\subSymb{QPQ} = q^2 \subSymb{P}.
	\end{equation}
\end{proposition}
%
%
%

With these preparations in hand, we introduce the following notion of a connection $1$-form induced by the subprincipal symbol of a pseudodifferential operator. 

%
%
%
\begin{definition}\label{def: P_compatible_connection}
    Let $(\sE, \sF \to M; \connectionE, \nabla^{\ms \sF})$ be two vector bundles over a manifold $M$, equipped with their connections $\connectionE, \nabla^{\ms \sF}$ and let $\halfDen \to M$ be the half-density bundle. 
    A connection $\nabla^{\ms \sE \otimes \sF}$ on the tensor product bundle $\sE \otimes \sF \to M$ will be called $P \in \PsiDOMEF$-\textbf{compatible} if and only if 
	\begin{equation}
        \forall \xxi \in \Char{P} : 
		\mathsf{\Gamma} \big( (\rd_{x} \pi) X_{p} \big) \xxi = \ri \subSymb{P} \xxi, 
	\end{equation}
	where $\subSymb{P}$ resp. $\Char{P}$ are the subprincipal symbol resp. the characteristic set of $P$, $\mathsf{\Gamma}$ is the connection $1$-form of $\nabla^{\ms \sE \otimes \sF}$ and $\pi: \dotCoTanM \to M$ is the punctured cotangent bundle. 
	In other words, $\nabla^{\ms \sE \otimes \sF}$ indices the covariant derivative 
	\begin{equation}
        \nabla_{X_{p}}^{\ms \pi^{*} \Hom{\sE, \sF}} = X_{p} + \mathsf{\Gamma} \big( (\rd \pi) X_{p} \big) 
	\end{equation}
    on the bundle $\pi^{*} \Hom{\sE, \sF} \to \dotCoTanM$ along the Hamiltonian vector field $X_{p}$ generated by the principal symbol $p$ of $P$. 
\end{definition}
%
%
%

This definition makes sense because both quantities, $\mathsf{\Gamma}$ and $\subSymb{P}$ have the same transformation law under change of bundle frame. 
Hence $X_{p} + \ri \subSymb{P}$ has an invariant meaning.

%
%
%
\begin{proposition} \label{prop: L_P_L_compatible_connection}
	Let $\sE, \halfDen \to M$ be a vector bundle resp. the half-density bundle over a manifold $M$ and $m, m' \in \R, N \in \N$.  
	Suppose that $P \in \PsiDOME$ has a scalar principal symbol $p$ and that $\sE$ admits a $P$-compatible connection $\connectionE$. 
	Assume that $Q \in \PsiDO{m'}{M, \matN \otimes \halfDen}$ having a scalar principal symbol $q$ and vanishing subprincipal symbol. 
	Then $\connectionE$ is $QPQ$-compatible.
\end{proposition}
%
%
%

\begin{proof}
	We fix any bundle frame to check this.
	Then 
    $X_{qpq} 
		= \{qpq, \cdot \} 
		= \{q^{2} p, \cdot \} 
		=  q^{2} \{p , \cdot \} + 2 p q \{q, \cdot \}
    $.
	When restricted to $\Char{P}$ this equals $q^{2} \{p , \cdot \} = q^{2} X_{p}$.
	Hence, on $\Char{P}$ we have 
	\begin{equation}
		\connectionEndPiE_{ X_{qpq} } - X_{qpq} 
		= q^{2}  \Big( \connectionEndPiE_{X_{p}} - X_{p} \Big) 
		= q^{2} \subSymb{P} 
		= \subSymb{QPQ}.
	\end{equation}
\end{proof}
%
%
%

The primary observation is now that the Weitzenb\"{o}ck connection (Remark~\ref{rem: Weitzenboeck_connection}) defined by a normally hyperbolic operator $\square$ is compatible with $\square$ in the above sense. 

%
%
\begin{theorem} \label{thm: NHOp_compatible_connection}
	Let $\sE \to \sM$ be a vector bundle over a Lorentzian manifold $\spacetime$ and $\square$ a normally hyperbolic operator on $\sE$. 
	If $\connectionE$ is the corresponding Weitzenb\"{o}ck connection then it is $\square$-compatible.
\end{theorem}
%
%
%

\begin{proof}
    Since the computation is local, we choose a coordinate basis $\{(x^{i}, \partial_{i}) \}$ of $\tangent_{x} \sM$ to compute the Hamiltonian vector field $X_{\fg}$ generated by the principal symbol $\fg^{-1}$ of $\square$:  
	\begin{equation} \label{eq: HVF_NHOp}
		X_{\fg} \xxi  
		= 
		\frac{\partial (\fg^{\mu \nu} \xi_{\mu} \xi_{\nu})}{\partial \xi_{i}} \frac{\partial}{\partial x^{i}} - \frac{\partial (\fg^{\mu \nu} \xi_{\mu} \xi_{\nu})}{\partial x^{i}} \frac{\partial}{\partial \xi_{i}} 
		=  2 g^{ij} \xi_{j} \frac{\partial}{\partial x^{i}} - \frac{\partial \fg^{\mu \nu}}{\partial x^{i}} \xi_{\mu} \xi_{\nu} \frac{\partial}{\partial \xi_{i}}.
	\end{equation}
    Hence, $(\rd_{x} \pi) X_{\fg} = 2 g^{ij} \xi_{j} \partial_{i}$ and picking any bundle frame $\{ \rE_{r} \}$, we can see that 
	\begin{equation}
        \mathsf{\Gamma} \big( (\rd_{x} \pi) X_{\fg} \big)  \xxi 
        = 2 g^{ij} \mathsf{\Gamma}_{i} \xi_{j}  
		\stackrel{\eqref{eq: subprincipal_symbol_NHOp}}{=} \ri \subSymb{\square} \xxi.  
	\end{equation}
\end{proof}
%
%
%

A straightforward application of this theorem and Proposition~\ref{prop: L_P_L_compatible_connection} yields 

%
%
%
\begin{corollary} \label{cor: L_NHOp_L_compatible_connection}
	As in the terminologies of Theorem~\ref{thm: NHOp_compatible_connection}, let $N \in \N$. 
	Suppose that $Q \in \PsiDO{-1/2}{\sM, \matN}$ having a scalar principal symbol and vanishing subprincipal symbol. 
	Then $\connectionE$ is $Q \square Q$-compatible, i.e. the Weitzenb\"ock connection is compatible with the first-order operator $Q \square Q$.
\end{corollary}
%
%
%
%
%
%
%
%
%
%
\subsection{Microlocal conjugate of a NHO}
\label{sec: microlocalisation_NHOp}
Two pseudodifferential operators are called microlocally conjugate if they can be conjugated to one another by an elliptic Fourier integral operator once they have been appropriately localised in cotangent space. 
The key point is that under some natural assumptions any first-order pseudodifferential operator can be microlocally conjugated to a vector field. 
This is originally due to 
Duistermat-H\"{o}rmander~\cite[Prop. 6.1.4]{Duistermaat_ActaMath_1972} 
for scalar operators of real principal type (Definition~\ref{def: real_principal_type_PsiDO_mf}), which has been extended to vector bundles by 
Dencker~\cite{Dencker_JFA_1982} 
who formulated microlocal conjugate of a system of classical pseudodifferential operators locally of real principal type (Definition~\ref{def: real_principal_type_PsiDO}).  
More precisely, Dencker transformed the system of operators to a scalar pseudodifferential operator $P$ (modulo smoothing operators) with vanishing subprincipal symbol by conjugating with a system of elliptic pseudodifferential operators.  
Then, he depicted the microlocal conjugation of $P$ with $\rD_{1} := - \ri \one_{\matk} \partial_{1}$ by appropriate elliptic Fourier integral operators associated with the graph of a symplectomorphism locally connecting $\Char P$ to $\Char \rD_{1}$.  
\newline

In this section we will explain microlocalisation in an intrinsic geometric language and formulate Dencker's result in a slightly more general form. 
We are going to show that any first-order pseudodifferential operator $P$ on $\sE$ of real principal type having a scalar principal symbol is microlocally conjugate to $\rD_{1}$ in the sense of~\eqref{eq: microlocal_conjugate_P} in the theorem below. 

%
%
%
\begin{theorem}[Microlocalisation] \label{thm: microlocalisation_P}
    Let $\sE, \halfDen \to M$ be a vector bundle of rank $k$ and the bundle of half-densities over a $d$-dimensional manifold $M$, respectively.	
    Suppose that $P \in \PsiDO{1}{M; \sE \otimes \halfDen}$ is a properly supported first-order pseudodifferential operator on $\sE$ with real scalar principal symbol $p$ such that 
	\begin{enumerate}[label=(\alph*)]
		\item \label{con: p_x_xi_zero}
		$p \xxiNot = 0$ for some element $\xxiNot$ in the punctured cotangent bundle $\dotCoTanM$ of $M$; 
		\item \label{con: HVF_radial}
		the Hamiltonian vector field $X_{p}$ of $p$ and the radial direction are linearly independent at $\xxiNot$.
    \end{enumerate}
	Then for any $m \in \R$, there exist 
	\begin{enumerate}[label=(\roman*)] 
		\item 
		a homogeneous symplectomorphism $\varkappa$ from an open conic neighbourhood $\tilde{\cU}$ of $(0, \eta_{1}$ $ \rd y^{1})$ in $\dotCoTanRd$ to an open conic coordinate chart $\big( \cU, (x^{i}, \xi_{i}) \big)$ of $\xxiNot$ in $\dotCoTanM$ such that 
        \begin{equation} \label{eq: choice_symplecto}
            \varkappa^{*} p = \xi_{1} \mathbbm{1}_{\Hom{\sE, \sE}}; 
        \end{equation}
		\item 
		properly supported Lagrangian distributions 
        $\fB \in I^{m} \big( M \times \Rd, \varGamma'; \mathrm{Hom} (\C^{k}, \sE) \otimes \halfDen (M \times \Rd) \big)$   
        and 
        $\tilde{\fB} \in I^{-m} \big( \Rd \times M, \varGamma^{-1 \prime}; \Hom{\sE, \C^{k}} \otimes \halfDen (\Rd \times M) \big)$ 
        so that $B \tilde{B}, \tilde{B} B$ both are zero-order pseudodifferential operators and  
        \begin{subequations}
            \begin{eqnarray}
                \WF' \fB \subset \cU_{(x_{0}, \xi^{0}; 0, \eta_{1} \rd y^{1})}, 
                && 
                \WF' \tilde{\fB} \subset \tilde{\cU}_{(0, \eta_{1} \rd y^{1}; x_{0}, \xi^{0})},
                \label{eq: microlocalisation_1}
                \\ 
                \xxiNot \notin \ES (B \tilde{B} - I_{\sE}),  
                && 
                (0, \eta_{1} \rd y^{1}) \notin \ES (\tilde{B} B - I), 
                \label{eq: microlocalisation_2}
                \\ 
                (x_{0}, \xi^{0}) \notin \ES (B \rD_{1} \tilde{B} - P), 
                &&
                (0, \eta_{1} \rd y^{1}) \notin \ES (\tilde{B} P B - \rD_{1}),      
                \label{eq: microlocal_conjugate_P}
            \end{eqnarray}
        \end{subequations}
        where $\varGamma$ is the graph of $\varkappa$, $\rD_{1} := - \ri \one_{\matk} \partial / \partial y^{1} : C_{\mathrm{c}}^{\infty} (\Rd, \CN) \to C_{\mathrm{c}}^{\infty} (\Rd, \CN)$, and $\cU_{(x_{0}, \xi^{0}; 0, \eta_{1}\rd y^{1})}$ resp. $\tilde{\cU}_{(0, \eta_{1}\rd y^{1}; x_{0}, \xi^{0})}$ are small conic neighbourhoods of $(x_{0}, \xi^{0}; 0$, $\eta_{1} \rd y^{1}) \in \dotCoTanM \times \dotCoTanRd$ resp. $(0, \eta_{1}\rd y^{1}; x_{0}, \xi^{0}) \in \dotCoTanRd \times \dotCoTanM$.
	\end{enumerate}
	In additional, if there are more structures on $\sE$ then the Fourier integral operators $B, \tilde{B}$ can be simplified further as follows. 
	\begin{enumerate}[label=(\Alph*)]
		\item 
		If $\sE \to M$ is endowed with a sesquilinear form $(\cdot|\cdot)$ with respect to which $P$ is symmetric, then $\tilde{B}$ can be chosen as the formal adjoint of $B$ provided that $\CN$ is endowed with a standard sesquilinear form of the same signature as $(\cdot|\cdot)$. 
		\item 
		If $\sE \to M$ is equipped with a $P$-compatible connection $\connectionE$ then the principal symbols of $\fB$ resp. $\tilde{\fB}$ can be chosen $\one$ near $(x_{0}, \xi^{0}; 0, \eta_{1} \rd y^{1})$ resp. $(0, \eta_{1} \rd y^{1}; x_{0}, \xi^{0})$ with respect to a frame that is parallel along $X_{p}$. 
		\item 
		On the vector bundle $\big( \sE \to M, (\cdot|\cdot), \connectionE \big)$ where $\connectionE$ is a $P$-compatible connection and $P$ is symmetric with respect to the sesquilinear form $(\cdot|\cdot)$, we can choose $\fB$ such that its principal symbol equals $\one$ near $(x_{0}, \xi^{0}; 0, \eta_{1} \rd y^{1})$ with respect to a frame that is unitary and parallel along $X_{p}$, and  $\tilde{\fB}= \fB^{*}$.
	\end{enumerate}
\end{theorem}
%
%
%

A schematic of this notion has been portrayed in Figure~\ref{fig: microlocalisation}.  

%
%
%
\begin{center}
	\begin{tikzpicture}
		\node (a) at (0,0) {$\cU$};
		\node (b) at (4,0) {$\cU'$}; 
		\node (c) at (0,2) {$\pi^{*} \Hom{\sE, \sE}_{\cU}$};
		\node (d) at (4,2) {$\matk$}; 
		\node[below] at (2, 0) {$\varkappa$}; 
		\node[left] at (0,1) {$\symb{P}$}; 
		\node[above] at (2.1, 2) {$\hat{\varkappa}$}; 
		\node[right] at (4,1) {$\symb{\rD_{1}}$};  
		\draw[->] (b) -- (a); 
		\draw[->] (d) -- (c); 
		\draw[->] (a) -- (c); 
		\draw[->] (b) -- (d);
	\end{tikzpicture}
	\hfil
	\begin{tikzpicture}
		\node (a) at (0,0) {$C_{\rc}^{\infty} (M; \sE \otimes \halfDen)$};
		\node (b) at (5,0) {$C_{\textrm{c}}^{\infty} (\Rd, \CN)$}; 
		\node (c) at (0,2) {$C_{\rc}^{\infty} (M; \sE \otimes \halfDen)$};
		\node (d) at (5,2) {$C_{\textrm{c}}^{\infty} (\Rn, \CN)$}; 
		\node[below] at (2.75, 0) {$B$}; 
		\node[left] at (0,1) {$P$}; 
		\node[above] at (2.75, 2) {$\tilde{B}$}; 
		\node[right] at (5,1) {$\rD_{1}$}; 
		\draw[->] (b) -- (a); 
		\draw[->] (a) -- (c); 
		\draw[->] (c) -- (d); 
		\draw[->] (b) -- (d);
	\end{tikzpicture}
	\captionof{figure}[A schematic diagram of microlocalisation]{A 
        schematic diagram of microlocalisation. The diagram on the right commutes in a microlocal sense~\eqref{eq: microlocal_conjugate_P} and the map $\hat{\varkappa}$ is defined by  $\hat{\varkappa} (\cdot) = \symb{\fB} (\cdot) \symb{\tilde{\fB}}$ with all other symbols as defined in Theorem~\ref{thm: microlocalisation_P}.}
	\label{fig: microlocalisation}
\end{center}
%
%
%

%
%
%
\begin{proof}
	We will prove the proposition imitating the strategy used for the scalar version~\cite[Prop. 6.1.4, Lem. 6.6.4]{Duistermaat_ActaMath_1972} 
	(see also~\cite[Prop. 26.1.3]{Hoermander_Springer_2009}). 
	The existence of $\varkappa$ satisfying~\eqref{eq: choice_symplecto} is guaranteed by the  
	homogeneous Darboux Theorem~\ref{thm: homogeneous_Darboux}	which prerequisites our hypotheses~\ref{con: p_x_xi_zero} and~\ref{con: HVF_radial}.  
	Suppose that $b \in S^{m} \big( \varGamma; \Maslov \otimes \widetilde{\mathrm{Hom}} (\CN, \sE) \big)$ has an inverse in a conic neighbourhood of $(x_{0}, \xi^{0}; 0, \eta_1 \rd y^{1}) \in \varGamma$. 
	Then we can obtain a properly supported  $\fB_{1} \in I^{m} \big( M \times \Rd, \varGamma'; \Hom{\CN, \sE} \otimes \halfDen (M \times \Rd) \big)$ such that $\WFPrime{\fB_{1}} \subset \cU_{(x_{0}, \xi^{0}; 0, \eta_1 \rd y^{1})}$ and $\fB_{1}$ is non-characteristic (Definition~\ref{def: elliptic_FIO}) at $(x_{0}, \xi^{0}; 0, \eta_1 \rd y^{1})$ from the construction given in 
	Section~\ref{sec: FIO_symplecto}. 
    \newline 

	By Theorem~\ref{thm: existence_parametrix_FIO}, there exists a unique microlocal parametrix $\tilde{\fB}_{1} \in I^{- m} \big( \Rd \times M, \varGamma^{-1 \prime}$; $\Hom{\sE, \CN} \otimes \halfDen (\Rd \times M) \big)$ such that~\eqref{eq: microlocalisation_2} is satisfied . 
	Since $\fB_{1}$ and $\tilde{\fB}_{1}$ have reciprocal principal symbols to each other on $\cU_{(x_{0}, \xi^{0}; 0, \eta_1 \rd y^{1})}$, $\tilde{B}_{1} P B_{1}$ has the principal symbol $\eta_{1} \one_{\matk}$ on $\cU'$ due to~\eqref{eq: choice_symplecto} and the Egorov Theorem~\ref{thm: Egorov}. 
	Furthermore 
	\begin{equation} \label{eq: WF_B1_P_A1_D_Q}
		(0, \eta_1 \rd y^{1}) \notin \ES (\tilde{B}_{1} P B_{1} - \rD_{1} - Q)
	\end{equation}
    for some $Q \in \PsiDON{0}$. 
    \newline

	To find the operators $B$ and $\tilde{B}$ with the claimed properties, we construct properly supported elliptic $B_{2}, \tilde{B}_{2} \in \PsiDON{0}$ such that 
	\begin{eqnarray}
		&& 
		\tilde{B}_{2} B_{2} \equiv I 
		\Leftrightarrow 
		\tilde{B}_{2} B_{2} = I \mod \PsiDON{-\infty}, 
		\label{eq: 26_1_03_Hoermander}
		\\ 
		&& 
		\tilde{B}_{2} (\rD_{1} + Q) B_{2} \equiv \rD_{1} 
		\Leftrightarrow 
		\tilde{B}_{2} (\rD_{1} + Q) B_{2} = \rD_{1} \mod \PsiDON{-\infty},    
		\label{eq: 26_1_3_Hoermander}
	\end{eqnarray}
	and one sets $B := B_{1} B_{2}$ and $\tilde{B} := \tilde{B}_{2} \tilde{B}_{1}$. 
	Then~\eqref{eq: 26_1_03_Hoermander},~\eqref{eq: 26_1_3_Hoermander},  and~\eqref{eq: WF_B1_P_A1_D_Q} imply 
	\begin{eqnarray}
		(0, \eta_1 \rd y^{1}) \notin \ES \big( \tilde{B}_{2} (\tilde{B}_{1} B_{1} - I) B_{2} \big)  
		= 
		\ES (\tilde{B} B - I),  
		\\ 
		(0, \eta_1 \rd y^{1}) \notin \ES \big( \tilde{B}_{2} (\tilde{B}_{1} P B_{1} - \rD_{1} - Q) B_{2} \big)  
		= 
		\ES (\tilde{B} P B - \rD_{1}), 
	\end{eqnarray}
	which proves the second half of~\eqref{eq: microlocalisation_2} and~\eqref{eq: microlocal_conjugate_P}, and the first half follows immediately after we multiply from left and right by $B$ and by $\tilde{B}$.
    \newline 

	It therefore remains to construct $B_2$ and $\tilde{B}_2$ such that~\eqref{eq: 26_1_03_Hoermander} and~\eqref{eq: 26_1_3_Hoermander} hold.
	By the existence of a parametrix 
	(see e.g.~\cite[Thm. 18.1.24]{Hoermander_Springer_2007}), 
	for every elliptic $B_{2} \in \PsiDON{0}$ there exists $\tilde{B}_{2} \in \PsiDON{0}$ such that $\tilde{B}_{2} B_{2} - I$ and $B_{2} \tilde{B}_{2} - I$ are smooth. 
	Multiplying~\eqref{eq: 26_1_3_Hoermander} by $B_{2}$ from the left we arrive at the equivalent condition for~\eqref{eq: 26_1_3_Hoermander} that 
	\begin{equation} \label{eq: 26_1_3_Hoermander_equiv} 
        [\rD_{1}, B_{2}]_{-} + Q B_{2} \equiv 0   
		\tag{\ref{eq: 26_1_3_Hoermander}$^{\prime}$} 
	\end{equation}
	for some elliptic $B_{2}$. 
    \newline 

	We will now construct a solution of~\eqref{eq: 26_1_3_Hoermander_equiv} order by order, starting with the principal symbol. 
	The principal symbol of~\eqref{eq: 26_1_3_Hoermander_equiv} vanishes provided  
	\begin{equation}
		- \ri \left\{ \symb{\rD_{1}}, \symb{B_{2}} \right\} + \symb{Q} \symb{B_{2}} = 0,   
	\end{equation} 
	as the the subprincipal symbol of $\rD_{1}$ vanishes; cf.~\eqref{eq: symbol_commutator_PsiDO}.  
	If $q$ is the principal symbol of $Q$, then the preceding equation yields 
	\begin{equation} \label{eq: microlocalisation_first_transport}
		\parDeri{y^{1}}{b_0} = - \ri q b_0   
	\end{equation}
    for the principal symbol $b_0$ of $B_2$. 
    This is a first-order differential equation, hence a unique solution exists given an initial condition $b_0 (y^{1}=0, \cdot) = \one_{\C^{k \times k}}$ and this solution depends smoothly on $q$.
	By construction  $b_0 \yeta$ is homogeneous of degree zero and $\det (b_0 \yeta)$ is non-vanishing. 
    Therefore $b_0^{-1}$ exists and it is homogeneous of degree zero.  
	Defining a properly supported $B_{2, 0} \in \PsiDON{0}$ with homogeneous principal symbol $b_0$, we now successively construct $B_{2, l} \in \PsiDON{-l}$ 
	so that, for every $n \in \N$:  
	\begin{equation} \label{eq: microlocalisation_k_transport} 
        [\rD_{1}, B_{2, 0} + \ldots + B_{2, n}]_{-} 
		+ 
		Q (B_{2, 0} + \ldots + B_{2, n}) = R_{n+1} \in \PsiDON{- (n+1)}.   
	\end{equation} 
	This is equivalent to the corresponding principal symbols $b_{n}$ of $B_{2, n}$ and  $r_{n}$ of $R_{n}$ of degree $-n$, to satisfy	
	\begin{equation}
		- \ri \parDeri{y^{1}}{b_{n}} + q b_{n} = - r_{n}. 
	\end{equation} 	
	This equation can be solved by the Duhamel principle and the solution reads 
	\begin{equation} \label{eq: microlocalisation_sol_k_transport}
		b_{n} \yeta = - \ri  \, b_0 \yeta \int_{0}^{y^{1}} 
		b_0^{-1} (t, y^{2}, \ldots, y^{d}; \eta) \, r_{n} (t, y^{2}, \ldots, y^{d}; \eta) \, \rd t.   
	\end{equation}
	Then, using asymptotic summation (see Definition~\ref{def: asymptotic_summation}) of the symbols of $B_{2,n}$, one can now construct an operator $B_2$ satisfying \eqref{eq: 26_1_3_Hoermander_equiv}. 
	\newline 

	\noindent 
	\textbf{The case when $P$ is symmetric:} 
	\\ 
	We are now going to prove that the operators $B$ and $\tilde{B}$ can be chosen microlocally unitary in case $P$ is symmetric with respect to $(\cdot|\cdot)$.
	We endow the space $C_{\mathrm{c}}^{\infty} (\Rd, \CN)$ with a standard sesquilinear scalar product of the same signature as $(\cdot|\cdot)$. 
    The operator $\rD_{1}$ is then symmetric.  
	Furthermore, we make the choice $m = 0$. 
	Acting with $\fB \in I^{0} \big( M \times \Rd, \varGamma'; \Hom{\CN, \sE} \otimes \halfDen (M \times \Rd) \big)$ from the left of~\eqref{eq: microlocal_conjugate_P} gives the equivalent microlocal conjugate relation   
	\begin{equation}
		(x_{0}, \xi^{0}; 0, \eta_1 \rd y^{1}) \notin \WFPrime (P \fB - \fB \rD_{1}) 
	\end{equation}
    between $P$ and $\rD_{1}$. 
	Taking its adjoint we have $(0, \eta_1 \rd y^{1}; x_{0}, \xi^{0}) \notin \WFPrime (\fB^{*} P - \rD_{1} \fB^{*})$ where  $\fB^{*} \in I^{0} \big( \Rd \times M, \varGamma^{-1 \prime}; \Hom{\sE^{*}, \CN} \otimes \halfDen (\Rd \times M) \big)$ and consequently   
	\begin{subequations}
		\begin{eqnarray}
			&& (0, \eta_1 \rd y^{1}) \notin \ES (B^{*} P B - \rD_{1} B^{*} B), 
			\\ 
			&& (0, \eta_1 \rd y^{1}) \notin \ES (B^{*} P B - B^{*} B \rD_{1}).
		\end{eqnarray}  
	\end{subequations}
    Thus $(0, \eta_1 \rd y^{1}) \notin \ES{[B^{*} B, \rD_{1}]_{-}}$ and $B^{*} B$ is non-characteristic at $(0, \eta_1 \rd y^{1})$. 
    Since $B^{*} B$ is a pseudodifferential operator and $\ri \rD_{1}$ is differentiation with respect to $y^1$, the total symbol of the commutator $[B^{*} B, \rD_{1}]_{-}$ is $-\ri$ times the $y^1$-derivative of the total symbol of $B^{*} B$.
	Thus, on the the level of symbols, the above implies that the total symbol of $B^{*} B$ is the sum of a term independent of $y^1$ and a term that is rapidly decaying in a conic neighborhood $\cU'$ of $(0, \eta_1 \rd y^{1})$. 
	It is bounded below on $\cU'$ because of the ellipticity of $B^{*} B$ there. 
    By Proposition~\ref{prop: FIO_2_2_2_matrix} (presented immediately after this proof), one can find a properly supported symmetric $\varPsi \in \PsiDONU{0}$ whose principal symbol is the same as that of $B^{*} B$ and such that  
	\begin{subequations}
		\begin{eqnarray}
			&& 
			(0, \eta_1 \rd y^{1}) \notin \Char{\varPsi},
			\label{eq: microlocalisation_varPsi_1}
			\\ 
			&& 
			(0, \eta_1 \rd y^{1}) \notin \ES (\varPsi^{*} \varPsi - B^{*} B), 
			\label{eq: microlocalisation_varPsi_2}
			\\ 
			&& 
            (0, \eta_1 \rd y^{1}) \notin \ES{[\varPsi, \rD_{1}]_{-}}, 
			\label{eq: microlocalisation_varPsi_3} 
		\end{eqnarray}
	\end{subequations}
	where $U'$ is the projection of $\cU'$ on $\Rd$. 
	Note that Proposition~\ref{prop: FIO_2_2_2_matrix} (stated for operators acting on the same bundle) can be used only in a local trivialisation of $\sE$ by an orthonormal frame that identifies the fibre of $\sE$ with $\CN$ in such a way that the sesquilinear forms are identified. It is at this stage that we must require the \textit{sesquilinear forms to have the same signature}. 
    \newline 

	Here the last property follows from the fact that the construction of the full symbol of $\varPsi$ in Proposition~\ref{prop: FIO_2_2_2_matrix} involves only multiplication and asymptotic summation of symbols. 
	The property of a symbol being a sum of two terms, one independent of $y_1$ and another rapidly decaying in a conic neighborhood of $\cU'$, is preserved under these operations. 
	The full symbol of $\varPsi$ is therefore also of this form. 
    \newline 

	The ellipticity~\eqref{eq: microlocalisation_varPsi_1} entails a unique microlocal parametrix $\varPhi$ for $\varPsi$. 
	In other words, there exists a properly supported $\varPhi \in \PsiDONU{0}$ such that 
	\begin{subequations}
		\begin{eqnarray}
			&&  
			(0, \eta_1 \rd y^{1}) \notin \Char{\varPhi},
			\label{eq: microlocalisation_varPhi_1} 
			\\ 
			&& 
			(0, \eta_1 \rd y^{1}) \notin \ES (\varPhi \varPsi - I) 
			\Leftrightarrow (0, \eta_1 \rd y^{1}) \notin \ES (\varPsi \varPhi - I).
			\label{eq: microlocalisation_varPhi_2} 
		\end{eqnarray}
	\end{subequations}
	Note,~\eqref{eq: microlocalisation_varPhi_1},~\eqref{eq: microlocalisation_varPsi_2}, and~\eqref{eq: microlocalisation_varPsi_3} imply that 
	\begin{subequations}
		\begin{eqnarray}
			&& (0, \eta_1 \rd y^{1}) \notin \ES (\varPhi^{*} \varPsi^{*} \varPsi \varPhi - \varPhi^{*} B^{*} B \varPhi), 
			\\ 
            && (0, \eta_1 \rd y^{1}) \notin \ES (\varPhi [\varPsi, \rD_{1}]_{-} \varPhi).   
		\end{eqnarray}
	\end{subequations}
	This, accounting~\eqref{eq: microlocalisation_varPhi_2} entails that 
	\begin{subequations}
		\begin{eqnarray}
			&& (0, \eta_1 \rd y^{1}) \notin \ES \big( I - (B \varPhi)^{*} B \varPhi \big), 
			\\ 
            && (0, \eta_1 \rd y^{1}) \notin \ES [\varPhi, \rD_{1}]_{-},  
		\end{eqnarray}
	\end{subequations} 
    which completes the proof since $(x_{0}, \xi^{0}; 0, \eta_1 \rd y^{1}) \notin \WF' (\fB [\rD_{1}, \varPhi]_{-})$ and therefore with $\fB^{\backprime} := \fB \varPhi$ we have $(x_{0}, \xi^{0}; 0, \eta_1 \rd y^{1}) \notin \WFPrime (P \fB^{\backprime} - \fB^{\backprime} \rD_{1})$. 
	\newline 

	\noindent 
	\textbf{The case of connection $P$-compatibility:} 
	\\
	We assume that $\connectionE$ is $P$-compatible. 
	Since the construction of the symbols of $\fB$ and $\tilde{\fB}$ are local, we can fix a local frame and local coordinates. 
	We will reduce the general situation to the case when the subprincipal symbol of $P$ vanishes near $\xxiNot$. 
    \begin{enumerate}
    	\item[(i)] 
		\textit{The case of vanishing subprincipal symbol}:  
        We assume that $P$ is given in local coordinates with respect to some bundle frame and that in this frame $\subSymb{P} = 0$ near $\xxiNot$. 
        Let $\tilde{P}$ be a scalar operator whose principal symbol is $p$ and subprincipal symbol vanishes.  
        Using Weyl-quantisation 
        (see e.g.~\cite[Sec. 18.5]{Hoermander_Springer_2007}), 
        $B_1$ and $\tilde{B}_1$ can be constructed in such a way that Egorov's theorem holds up to the subprincipal symbol level~\cite[Thm. 1]{Silva_PublMat_2007}. 
		We can choose these operators in such a way that their principal symbols are constant $\one$ and that the subprincipal symbol of $\tilde{B}_1 \tilde P B_1$ (resp. $B_{1} \rD_{1} \tilde{B}_{1}$)  vanishes near $(0, \eta_1 \rd y^{1})$ (resp. $\xxiNot$). 
        We will now use the same construction as before starting with $B_1$ and $\tilde B_1$. Then the principal symbol $q$ of the remainder term $Q$ vanishes. 
        The construction of $B_2$ and $\tilde{B}_2$ then yields operators with total symbols that are scalar and constant principal symbols equal to $\one$. 
        We conclude that in case the subprincipal symbol of $P$ vanishes near $\xxiNot$ and $B, \tilde{B}$ can be chosen as scalar operators. 
        The principal symbols that are constant along the flow lines of $X_p$ and $X_{\symb{\rD_1}}$. 
        \item[(ii)] 
		\textit{The case of non-vanishing subprincipal symbol}: 
        We will microlocally transform $P$ to a scalar pseudodifferential operator $\tilde{P}$. 
        To be specific, for any properly supported $\tilde{P} \in \PsiDO{1}{M; \halfDen}$, we want to have a  $\hat{B} \in \PsiDO{0}{M; \sE \otimes \halfDen}$ such that $\hat{B}$ is non-characteristic at $\xxiNot$ and $\xxiNot \notin \ES (P \hat{B} - \hat{B} \tilde{P} I)$. 
        We construct this operator locally and hence fix a local frame of $\sE$ near the point $x_0$. 
        In the pullbacked bundle $\pi^* \sE \to \dotCoTanM$ where $\pi : \dotCoTanM \to M$, we can also construct a local frame that is parallel along $X_p$ with respect to $\nabla_{X_{p}}^{\ms \pi^{*} \sE}$. 
		This local parallel frame is constructed by choosing a local transverse to $X_{p}$ and use the original frame on this transverse. 
		Parallel transport along the flow lines of $X_{p}$ then gives the desired frame. 
		The change of frame matrix from the original frame to the parallel frame is then a local section $b$ of $\Hom{\pi^{*} \sE, \pi^{*} \sE} \to \dotCoTanM$. 
        By construction, this frame is homogeneous of degree zero. 
        We now choose an elliptic zero-order pseudodifferential operator $\hat B$ whose principal symbol $\symb{\hat B}$ equals $b$ on $\cU$. 
        Let $\check{B}$ be a parametrix of $\hat B$. 
        The subprincipal symbol of $\check{B} P \hat B$ is equal to 
        $
        - \ri X_{p} \symb{\hat{B}} + [\subSymb{P}, \symb{\hat{B}}]_{-} 
        $ 
        on $\cU$, by an application of~\eqref{eq: symbol_commutator_PsiDO}. 
        By $P$-compatibility, this is exactly the formula for the connection $1$-form in the parallel local bundle frame and it therefore vanishes.
        Proposition~\ref{prop: P_connection_transformation} implies that this is precisely the formula for the subprincipal symbol of $\check{B} P \hat{B}$ on $\cU$ and hence, $\tilde P \equiv \check{B} P \hat{B}$  has vanishing subprincipal symbol on $\cU$. 
        \newline 

        This reduces the problem to the case of vanishing subprincipal symbol discussed in the hindmost paragraph and let $B_3$ and $\tilde{B}_3$ are the corresponding scalar Fourier integral operators. 
        This means $B$ and $\tilde{B}$ are of the form $B := \hat{B} B_3$ and $\tilde{B} := \tilde{B}_3 \check{B}$. 
        Since the principal symbol of $\hat{\fB}$ is the transition function to a parallel frame and $B_3$ is a scalar operator, these imply that the principal symbols of $\fB$ and $\tilde{\fB}$ can be chosen $\one$ with respect to a parallel frame along $X_p$.
    \end{enumerate}
	\textbf{The case of symmetric $P$ with connection $P$-compatibility:} 
	\\
	Finally, suppose that $P$ is symmetric with respect to $(\cdot|\cdot)$ and that $\connectionE$ is $P$-compatible. 
	The above construction of $\hat B$ can be repeated using an orthonormal frame of $\sE$ and another orthonormal frame of $\pi^* \sE$ that is parallel along $X_p$. 
    One obtains an operator that has principal symbol $\one$ near $\xxiNot$ with respect to the parallel orthonormal frame. 
	Now one repeats the construction of $\varPsi$ as before and note that the principal symbol can be kept as $\one$ in that way to turn $B$ into a microlocally unitary operator in the sense that $\tilde B = B^*$ with the desired properties. 
\end{proof}
%
%
%

As a supplement of the preceding proof we present a variant of H\"{o}rmander's square root construction~\cite[Prop. 2.2.2]{Hoermander_ActaMath_1971} for vector bundles. 

%
%
%
\begin{proposition} \label{prop: FIO_2_2_2_matrix}
	Let $\big(\sE \to M, (\cdot|\cdot) \big)$ be a vector bundle over a manifold $M$ where $(\cdot|\cdot)$ is a non-degenerate sesquilinear form, and $\halfDen \to M$ the half-density bundle over $M$. 
	Suppose that $P \in \PsiDOME$ is a properly supported pseudodifferential operator of order $m \in \R_{+}$, symmetric with respect to $(\cdot|\cdot)$, elliptic in a conic neighbourhood $\cU$ of $\xxiNot \notin \ES{P}$ in the punctured cotangent bundle $\dotCoTanM$ of $M$, and that its principal symbol is given by  
	\begin{equation}
        p = q^{*} q 
	\end{equation} 
	for some $q \in C^{\infty} (\cU, \End \sE)$. 
    Then, one can find a properly supported symmetric (with respect to $(\cdot|\cdot)$) $Q \in \PsiDO{m/2}{U; \sE \otimes \halfDen}$ such that 
	\begin{equation}
		\xxiNot \notin \ES (P - Q^{*} Q), \qquad \xxiNot \notin \ES{Q},   
	\end{equation}
	where $U$ is the base projection of $\cU$. 
\end{proposition}
%
%
%
\begin{proof}
	The hypothesis on $p$ implies that $q$ is homogeneous of degree $m/2$ which in turns entails that $q \in S^{m/2} (\cU, \End \sE)$; see  Definition~\ref{def: polyhomogeneous_symbol_Euclidean} and Appendix~\ref{sec: symbol}. 
    We define a properly supported $Q_{0} \in \PsiDO{m/2}{U; \sE \otimes \halfDen}$ whose principal symbol is $q$. 
	Without lose of generality $Q_{0}$ can be taken symmetric, otherwise one can just replace $Q_{0}$ by $(Q_{0} + Q_{0}^{*}) / 2$ without changing the principal symbol. 
    Then $P - Q_{0}^{*} Q_{0} \in \PsiDO{m-1}{U; \sE \otimes \halfDen}$. 
    \newline 

    Now we are left with estimations for lower order terms and will show that it is always possible to obtain properly supported and symmetric $Q_{k} \in \PsiDO{m/2-k}{U; \sE \otimes \halfDen}$ for all $k \in \N$, such that  
	\begin{equation}
        R_{k} := P - (Q_{0} + \ldots + Q_{k})^{*} (Q_{0} + \ldots + Q_{k}) \in \PsiDO{m-1-k}{U; \sE \otimes \halfDen}. 
	\end{equation}
	We proceed inductively. 
	Clearly, $k = 1$ has been checked.    
	Observe that, if $Q_{k}$s have been chosen accordingly, then   
	\begin{equation}
		P - |Q_{0} + \ldots + Q_{k}|^{2} 
		= 
        R_{k} - (Q_{k}^{*} Q_{0} + Q_{0}^{*} Q_{k}) - Q_{k}^{*} Q_{k} + \ldots \in \PsiDO{m-1-k}{U; \sE \otimes \halfDen}. 
	\end{equation}
	Since $R_{k}$ is symmetric, the principal symbol of $R_{k} - R_{k}^{*}$ is $2 \ri \Im (\symb{R_{k}})$ modulo $S^{m - 2 - k} (\cU, \mathrm{End}$ $ \sE)$. 
	Then the desired operators $B_{k}$ are achieved if one sets   
	\begin{equation}
		\forall \xxi \in \cU : 
		\symb{R_{k}} \xxi = \symb{Q_{k}}^{*} \circ q \, \xxi + q^{*} \circ \symb{Q_{k}} \, \xxi,   
	\end{equation}
	that is, in the region where $q (x, \hat{\xi})$ is invertible in $\End \sE_{x}$ for each $x \in U$ and $\hat{\xi} := \xi / \| \xi \| \in \bbS^{\rk (\sE) - 1}$. 
	Finally, the result entails by constructing the total symbol of $Q$ in the sense of asymptotic summation; in other words: $Q :\sim \sum_{k=0}^{\infty} Q_{k}$.
\end{proof}
%
%
%

\begin{remark} \label{rem: microlocalisation_real_principal_type}
	The hypotheses of the preceding theorem are satisfied by any pseudodifferential operator $Q$ of real principal type by Remark~\ref{rem: real_principal_type_PsiDO_mf_HVF_radial}, and hence Theorem~\ref{thm: microlocalisation_P} also holds for $Q$ by a redefinition of $B$. 
\end{remark}
%
%
%

We will now show microlocalisation of a normally hyperbolic operator $\square$ as a consequence of the preceding result by replacing the generic manifold $M$ by a Lorentzian (not necessarily globally hyperbolic) spacetime $\spacetime$. 
In particular, the role of the characteristic points $\xxiNot$ of $P$ in Theorem~\ref{thm: microlocalisation_P}~\ref{con: p_x_xi_zero} will be played by any lightlike covector on $\sM$ and thus, the bicharacteristics of $p$ are actually the lightlike geodesics in $\coTansM$. 

%
%
%
\begin{theorem} 
\label{thm: microlocalisation_NHOp}
	Let $\sE \to \sM$ be a vector bundle of rank $k$ over a $d$-dimensional Lorentzian manifold $\spacetime$ and $\Box : \comSecsME \to \comSecsME$ a normally hyperbolic operator. 
    Denote by $\connectionE$, the Weitzenb\"{o}ck connection associated with $\Box$ and by $\varGamma$, the graph of a homogeneous symplectomorphism from an open conic neighbourhood $\cU'$ of $(0, \eta_1 \rd y^{1})$ in $\dotCoTanRd$ to an open conic coordinate chart $\big( \cU, (x^{i}, \xi_{i}) \big)$ centered at any lightlike covector $\xxiNot$ on $\sM$.
	Then for any $m \in \R$, one can find properly supported Lagrangian distributions 
	$\fA \in I^{-1/2 + m} \big( \sM \times \Rd, \varGamma'; \Hom{\CN, \sE} \big)$  
	and 
    $\tilde{\fA} \in I^{-1/2 - m} \big( \Rd \times \sM, \varGamma^{-1 \prime}; \Hom{\sE, \CN} \big)$ so that $A \tilde{A}, \tilde{A} A$ both are pseudodifferential operators of order $-1$ and 
	\begin{subequations}
		\begin{eqnarray}
			\WFPrime{\fA} \subset \cU_{(x_{0}, \xi^{0}; 0, \eta_1 \rd y^{1})}, 
			&&  
			\WFPrime{\tilde{\fA}} \subset \cU'_{(0, \eta_1 \rd y^{1}; x_{0}, \xi^{0})},
			\label{eq: microlocalisation_NHOp_1}
			\\ 
			\xxiNot \notin \ES (A \tilde{A} - I_{\sE}), 
			&& 
			(0, \eta_1 \rd y^{1}) \notin \ES (\tilde{A} A - I), 
			\label{eq: microlocalisation_NHOp_2}
			\\ 
			(x_{0}, \xi^{0}) \notin \ES (A \rD_{1} \tilde{A} - \Box), 
			&& 
			(0, \eta_1 \rd y^{1}) \notin \ES (\tilde{A} \Box A - \rD_{1}).      
			\label{eq: microlocal_conjugate_NHOp}
		\end{eqnarray}
	\end{subequations}
    Here $\rD_{1} = - \ri \one_{\matk} \partial / \partial y^{1} : C_{\rc}^{\infty} (\Rd, \CN) \to C_{\rc}^{\infty} (\Rd, \CN)$ and $\cU_{(x_{0}, \xi^{0}; 0, \eta_1 \rd y^{1})}$ resp. $\cU'_{(0, \eta_1 \rd y^{1}; x_{0}, \xi^{0})}$ are small conic neighbourhoods of $(x_{0}, \xi^{0}; 0, \eta_{1} \rd y^{1}) \in \dotCoTansM \times \dotCoTanRn$ resp. $(0, \eta_1 \rd y^{1}; x_{0}, \xi^{0}) \in \dotCoTanRn \times \dotCoTansM$. 
    \newline 

	In addition, if $\sE \to \sM$ is endowed with a sesquilinear form $(\cdot|\cdot)$ with respect to which $\square$ is symmetric, then $A$ can be chosen as a scalar operator with respect to a unitary bundle frame that is parallel along the geodesic flow and microlocally $A = \tilde{A}^*$ if $\CN$ is endowed with a standard sesquilinear form of the same signature as $(\cdot|\cdot)$. 
\end{theorem}
%
%
%
\begin{proof}
	The strategy is, as usual, to reduce $\square$ to a first-order operator so that  Theorem~\ref{thm: microlocalisation_P} can be applied. 
	To do so, choose a properly supported elliptic symmetric pseudodifferential operator $L$ on $\sE$ of order $-1/2$ having a scalar principal symbol $l$ and vanishing subprincipal symbol. 
	Such an operator always exists, as it can be constructed locally and then patched to a global operator using a suitable partition of unity.  
	Then $L \square L$ is a first-order pseudodifferential operator whose principal symbol  $l^{2} \xxi \, \fg_{x}^{-1} (\xi, \xi) \, \one_{\End{\sE}}$ vanishes on $\xxiNot$. 
	By Theorem~\ref{thm: NHOp_compatible_connection} and Corollary~\ref{cor: L_NHOp_L_compatible_connection}, $\nabla^{E}$ is compatible with both $\square$ and $L\square L$, respectively. 
	Therefore the hypotheses of Theorem~\ref{thm: microlocalisation_P} are satisfied and the conclusion entails by putting $\fA := L \fB, \tilde{\fA} := \tilde{\fB} L$ where $\fB, \tilde{\fB}$ are the Lagrangian distributions as constructed in Theorem~\ref{thm: microlocalisation_P}. 
	In case $\square$ is symmetric, we choose $L$ symmetric as well.  
\end{proof}
%
%
%

\begin{remark}
\label{rem: global_microlocalisation_NHOp} 
	We have microlocally conjugated $\square$ to $\rD_{1}$ around a lightlike covector in the preceding theorem. 
	It is also possible to refine the construction as follows.
	Let $\I$ be a compact interval in $\R$ and $\gamma : \I \to \dotCoTanM$ an integral curve of $X_{p}$:    
	\begin{equation}
		p \circ \gamma = 0. 
	\end{equation} 
	If the composition of $\gamma$ and the projection $\dotCoTanM \to \dotCoTanM / \R_{+}$ is injective, then one can find a conic neighbourhood $\cV'$ of $\I \times \{(0, \eta_{1} \rd y^{1}) \}$ and a smooth homogeneous symplectomorphism $\varrho$ from $\cV'$ to an open conic neighbourhood $\varrho (\cV') \subset \dotCoTanM$ of $\gamma (\I)$ such that~\cite[Prop. 26.1.6]{Hoermander_Springer_2009}  
	\begin{equation}
		\varrho \big( \I \times \{ (0, \eta_{1} \rd y^{1}) \} \big) = \gamma (\I), 
		\qquad \varrho^{*} p = \xi_{1} \one_{\End \sE}. 
	\end{equation}
	Imitating the proof of Theorem~\ref{thm: microlocalisation_P}, one can microlocalise $P$ to $\rD_{1}$ on $\gamma (\I)$ 
	(see~\cite[Prop. 26.1.3$'$]{Hoermander_Springer_2009} for the scalar version).  
	\newline 

	As a consequence, if $\Gamma$ is the graph of a homogeneous symplectomorphism from a conic neighbourhood of $\I \times \{ (0, \eta_{1} \rd y^{1}) \}$ in $\dotCoTanRd$ to a conic neighbourhood of lightlike geodesic $\gamma (\I)$ in $\dotCoTansM$, then the proof of Theorem~\ref{thm: microlocalisation_NHOp} shows that, for any $m \in \R$, there exists Lagrangian distributions 
	$\fA \in I^{-1/2 + m} \big( \sM \times \Rd, \Gamma'; \Hom{\CN, \sE} \big)$  
	and 
    $\tilde{\fA} \in I^{-1/2 - m} \big( \Rn \times \sM, \Gamma^{-1 \prime}; \Hom{\sE, \C^{N}} \big)$ so that $A \tilde{A}, \tilde{A} A$ both are pseudodifferential operators of order $-1$ and   
	\begin{subequations}
		\begin{eqnarray}
			\WF' \fA \subset \cV_{\Gamma}, 
			&& 
			\WF' \tilde{\fA} \subset \cV'_{\Gamma^{-1}},
			\label{eq: microlocalisation_NHOp_1_global}
			\\ 
			\gamma (\I) \cap \ES (A \tilde{A} - I_{\sE}) = \emptyset, 
			&& 
			\I \times \{ (0, \eta_{1} \rd y^{1}) \} \cap \ES (\tilde{A} A - I) = \emptyset, 
			\label{eq: microlocalisation_NHOp_2_global}
			\\ 
			\gamma (\I) \cap \ES (A \rD_{1} \tilde{A} - \Box) = \emptyset, 
			&& 
			\I \times \{ (0, \eta_{1} \rd y^{1}) \} \cap \ES (\tilde{A} \Box A - \rD_{1}) = \emptyset,      
			\label{eq: microlocal_conjugate_NHOp_global}
		\end{eqnarray}
	\end{subequations}
	where $\cV_{\Gamma}, \cV'_{\Gamma^{-1}}$ are small conic neighbourhoods of $\Gamma$ restricted to $\gamma (\I)$ and its inverse, respectively. 
\end{remark}
%
%
%

We close this section by a simple application of Theorem~\ref{thm: microlocalisation_P} to derive a bundle version of 
H\"{o}rmander's propagation of singularity theorem~\cite{Hoermander_Nice_1970},~\cite[Thm. 6.1.1']{Duistermaat_ActaMath_1972}.  
Since 
\begin{equation}
	C^{\infty} (M; \sE \otimes \halfDen) = \bigcap_{s \in \R} H_{\mathrm{loc}}^{s} (M; \sE \otimes \halfDen), 
\end{equation}
such a refinement of the usual notion of (smooth) wavefront set is captured by the Sobolev wavefront set. 

%
%
%
\begin{definition} \label{def: Sobolev_WF}
	Let $\sE, \halfDen \to M$ be a vector bundle resp. the half-density bundle over a manifold $M$ and $s \in \R$. 
	The \textbf{Sobolev wavefront set} $\WF^{s} u$ of a distribution $u \in \dualComSecE$ relative to Sobolev space $H_{\mathrm{loc}}^{s} (M; \sE \otimes \halfDen)$ is defined by~\cite[p. 201]{Duistermaat_ActaMath_1972}  
    \begin{equation} \label{eq: def_Sobolev_WF}
        \WF u
        := 
        \bigcap_{\substack{P \in \PsiDO{0}{M; \sE \otimes \halfDen} \\ P u \in H_{\mathrm{loc}}^{s} (M; \sE \otimes \halfDen)}} \Char P 
        = 
        \bigcap_{\substack{P \in \PsiDO{s}{M; \sE \otimes \halfDen} \\ P u \in L_{\mathrm{loc}}^{2} (M; \sE \otimes \halfDen)}} \Char P,  
    \end{equation}
    where $\Char P$ is the characteristic set of $P$ and the intersection is over all properly supported $P$. 
\end{definition}
%
%
%

Locally this means, for any open set $U \subset \Rd$, $\xxiNot \notin \WF^{s} u$ is not in the Sobolev wavefront set of a distribution $u \in \dualComSecE$ if and only if there exists a compactly supported section $f$ on $U$, non-vanishing at $x_{0} \in U$ such that the weighted Fourier transform $(1 + |\xi|^{2})^{s} \, | \Fourier (fu) (\xi) |^{2}$ is integrable in a conic neighbourhood of $ \xi^{0}$. 
An accumulation of important properties of this finer wavefront set is available, for example, in~\cite[App. B]{Junker_AHP_2002}. 
\newline 

H\"{o}rmander's theorem has been generalised for system of pseudodifferential operators by
Taylor~\cite[Thm. 4.1 $($p. 135$)$, 2.2 $($p. 154$)$]{Taylor_PUP_1981} (for the Sobolev wavefront set)
and by
Dencker~\cite[Thm. 4.2]{Dencker_JFA_1982} (for the polarisation set). 
Our result, as stated below, is similar to that of Taylor but formulated in a geometric fashion and proved in a different way employing microlocalisation on vector bundles developed in Theorem~\ref{thm: microlocalisation_P}.   

%
%
%
\begin{theorem}[Propagation of Sobolev regularity] 
\label{thm: propagation_singularity_Sobolev_WF}
	Let $\sE, \halfDen \to M$ be a vector bundle resp. the half-density bundle over a manifold $M$ and $u \in \dualComSecE$. 
	Suppose that, for some $m \in \R$, $\varPsi \in \PsiDO{m}{M; \sE \otimes \halfDen}$ satisfies the hypotheses~\ref{con: p_x_xi_zero} and~\ref{con: HVF_radial} of Theorem~\ref{thm: microlocalisation_P} and that $\cI$ is an interval on an integral curve of the Hamiltonian vector field generated by the principal symbol of $\varPsi$ such that $\cI \cap \WF^{s} (\varPsi u) = \emptyset$. 
	Then either $\cI \cap \WF^{s + m -1} u = \emptyset$ or $\cI \subset \WF^{s + m -1} u$. 
\end{theorem}
%
%
%

\begin{proof}
    Once Theorem~\ref{thm: microlocalisation_P} is at our disposal, the rest of the proof is the same as its scalar version~\cite[Thm. 6.1.1']{Duistermaat_ActaMath_1972}  
	(see also~\cite[Thm. 26.1.4]{Hoermander_Springer_2009}). 	
	For completeness we give the details here.
	By conjugating $\varPsi$ with appropriate Fourier integral operators, we reduce the statement to the case $m = 1$. 
	In other words, one can just consider $P$ as in Theorem~\ref{thm: microlocalisation_P} instead of $\varPsi$ to conclude the assertion by utilising the Sobolev continuity properties of pseudodifferential 
	(see e.g.~\cite[p. 92]{Hoermander_Springer_2007}) 
	and Fourier integral operators 
	(see e.g.~\cite[Cor. 25.3.2]{Hoermander_Springer_2009}). 
	Via microlocalisation (Theorem~\ref{thm: microlocalisation_P}), the analysis further boils down to $P = \rD_{1}, s = 0$ and $\xxiNot = (0, \eta_1 \rd y^{1})$ by choosing the order of the Fourier integral operator $B$ in Theorem~\ref{thm: microlocalisation_P} equals to $-s$. 
	Then the claim follows from the form of the advanced and retarded fundamental solutions (Example~\ref{ex: advanced_retarded_fundamental_sol_partial_derivative}) of $\rD_{1}$ and as these map from $L_{\mathrm{c}}^{2} (\Rd, \CN)$ to $L_{\mathrm{loc}}^{2} (\Rd, \CN)$ (cf. Proposition~\ref{prop: causal_propagator_partial_derivative_FIO}).  
\end{proof}
%
%
%

By Remark~\ref{rem: microlocalisation_real_principal_type}, the preceding theorem can be restated as 

%
%
%
\begin{theorem}
	Let $\sE, \halfDen \to M$ be a vector bundle resp. the half-density bundle over a manifold $M$ and $u \in \dualComSecE$. 
	Suppose that $P$ is an $m \in \R$-order pseudodifferential operator on $\sE \otimes \halfDen$ of real principal type and that $\cI$ is an interval on an integral curve of the Hamiltonian vector field generated by $q$ as in Definition~\ref{def: real_principal_type_PsiDO} such that $\cI \cap \WF^{s} (P u) = \emptyset$. 
	Then either $\cI \cap \WF^{s + m -1} u = \emptyset$ or $\cI \subset \WF^{s + m -1} u$. 
\end{theorem}
%
%
%
%
%
%
%
%
%
%
\subsection{Proof of Theorem~\ref{thm: exist_unique_Feynman_parametrix_NHOp}}
To begin with, we show that the Feynman parametrix is unique, if it exists. 
\vspace*{-0.25cm}
%
%
%
%
%
%
%
%
%
%
\subsubsection{Uniqueness of the Feynman parametrix} 
\label{sec: uniqueness_Feynman_parametrix}
Suppose that $\parametrixLeft$ (resp. $\parametrixRight$) is a left (resp. right) Feynman parametrix of $\square$, i.e., the off-diagonal contribution of $\WF' (\parametrixLeft, \parametrixRight)$ is given by $C^{+}$; cf. Definition~\ref{def: Feynman_parametrix}. 
To prove that $\parametrixLeft - \parametrixRight$ is a smoothing operator, we would like to argue that $\parametrixLeft \square \parametrixRight$ is congruent both to $\parametrixLeft$ and to $\parametrixRight$ modulo smoothing operators. 
But $\parametrixLeft$ and $\parametrixRight$ are not properly supported which makes this a non-trivial task. 
To circumvent this difficulty, one employs the fact that $\parametrixLeft Q  \parametrixRight$ is defined when $Q$ is a pseudodifferential operator having Schwartz kernel of compact support in $\sM \times \sM$ and then $Q : \dualComSecsME \to \dualSecsME$. 
\newline 

If $(x, \xi; y, \eta) \in \WFPrime (\parametrixKernelLeft Q \parametrixKernelRight)$ but $\xxi, \yeta$ are in the complement of $\ES Q$, then $(x, \xi; z, \zeta), (z, \zeta; y, \eta) \in C^{+}$ for some $(z, \zeta) \in \ES Q$. 
This follows from the behaviour of wavefront sets under composition of kernels, as available in, for instance~\cite[Thm. 8.2.10]{Hoermander_Springer_2003},
and the fact that $\WFPrime \fQ \subset \varDelta \dotCoTansM$. 
This entails that $(x, \xi), (y, \eta)$, and $(z, \zeta)$ are on the same lightlike geodesic (strip) $\gamma (s)$  with $(z, \zeta)$ in between the other two points. 
\newline 

Since $\sM$ is assumed to be globally hyperbolic, $J (K) := \eqref{eq: causal_future_past}$ is compact for any compact $K \subset \sM$~\cite[Corollary of Prop. 6.6.1]{Hawking_CUP_1973}. 
Therefore, if the projections of endpoints of $\gamma$ on $\sM$ are in $K$ then the projection $c (s)$ of $\gamma (s)$ on $\sM$ stays over $J (K)$, i.e., $c (s) \in J (K)$, as shown in Figure~\ref{fig: uniqueness_Feynman_parametrix}.  
Consequently, if $\singsupp Q \cap J (K) = \emptyset$ then $\WFPrime (\parametrixKernelLeft \fQ \parametrixKernelRight)$ cannot have any point over $K \times K$. 
Let $\chi \in \comSecsME$ such that $\chi$ is identically $1$ on $J (K)$ and define $\mathrm{M}_{\chi}$ be the corresponding multiplication operator. 
We observe that $[\mathrm{M}_{\chi}, \square]_{-}$ vanishes identically on $J (K)$. 
Thus $\WFPrime (\parametrixKernelLeft [\mathrm{M}_{\chi}, \square]_{-} \parametrixKernelRight)$ contains no  point over $K \times K$ and so it true for 
\begin{eqnarray}
	(\dotCoTan K \times \dotCoTan K) \cap \WFPrime (\parametrixKernelLeft - \parametrixKernelRight)
	& = & 
	(\dotCoTan K \times \dotCoTan K) \cap \WFPrime (\parametrixKernelLeft \mathrm{M}_{\chi} \square \parametrixKernelRight - \parametrixKernelLeft \square \mathrm{M}_{\chi} \parametrixKernelRight) 
	\nonumber \\ 
	& = &    
    (\dotCoTan K \times \dotCoTan K) \cap \WFPrime (\parametrixKernelLeft [\mathrm{M}_{\chi}, \square]_{-} \parametrixKernelRight) 
	\nonumber \\ 
	& = & 
	\emptyset. 
\end{eqnarray}
Since $K$ is arbitrary, we conclude that $\parametrixLeft - \parametrixRight$ is a smoothing operator and similarly for the anti-Feynman parametrix.  

%
%
%
\begin{center}
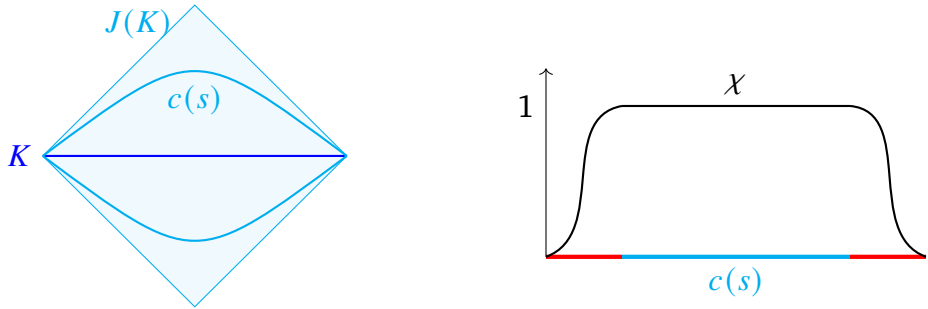

	\begin{tikzpicture}
		\filldraw[color = cyan, fill = cyan!5!white] (0,0) -- (2,2) -- (4,0) -- (2,-2) -- (0,0); 
		\draw[blue, thick] (0,0) -- (4,0); 
		\draw[cyan, thick] (0,0) .. controls (2,1.5) .. (4,0); 
		\draw[cyan, thick] (0,0) .. controls (2,-1.5) .. (4,0);
		\node[left] at (0,0) {$\textcolor{blue}{K}$}; 
		\node[left] at (1.85,1.75) {$\textcolor{cyan}{J (K)}$}; 
		\node at (2,0.75) {\textcolor{cyan}{$c (s)$}};
	\end{tikzpicture}
	\hfil
	\begin{tikzpicture}
		\draw[cyan, ultra thick] (1,0) -- (4,0); 
		\draw[red, ultra thick] (0,0) -- (1,0); 
		\draw[red, ultra thick] (4,0) -- (5,0); 
		\draw[thick] (1,2) -- (4,2);
		\draw[thick] (1,2) to [out=190, in=20] (0,0); 
		\draw[thick] (4,2) to [out=355, in=160] (5,0); 
		\draw[->] (0,0) -- (0, 2.5);
		\node[above] at (2.5, 2) {$\chi$};
		\node[below] at (2.5,0) {$\textcolor{cyan}{c (s)}$}; 
		\node[left] at (0,2) {$1$};
	\end{tikzpicture}
	\captionof{figure}[A consequence of global hyperbolicity of spacetime]{A 
        schematic visualisation of a consequence of the global hyperbolicity of spacetime. 
        Here time flows from left to right.}
    \label{fig: uniqueness_Feynman_parametrix}	
\end{center}
%
%
%
\subsubsection{Construction of the Feynman parametrix} 
Notice that $\square \parametrixRight = I + R$ is equivalent to $(\parametrixRight)^{*} \square^{*} = I + R^{*}$ and, since the adjoint $\square^{*}$ of $\square$ is also normally hyperbolic, the existence of left parametrices with the listed properties in this theorem follows from the existence of right parametrices of $\square$. 
It is, therefore, sufficient to construct a right parametrix of $\square$ with the required regularities. 
%
%
%

%
%
%
\begin{lemma} \label{lem: Lem_26_1_15_Hoermander}
	As in the terminologies of Theorem~\ref{thm: microlocalisation_NHOp}: $\square$ and $\rD_{1}$ are microlocally conjugate to each other, $\fA, \tilde{\fA}$ are the conjugating Lagrangian distributions with $m = 0$, and $\cU$ is the conic neighbourhood of any lightlike covector $\xxiNot$ on $\sM$.  
    Suppose that $\fF_{1}^{\ret, \adv}$ are the Schwartz kernels of the retarded and advanced fundamental solutions of $\rD_{1}$ and that $\chi \in C_{\mathrm{c}}^{\infty} (\Rd \times \Rd, \R)$ is identically $1$ in a neighbourhood of the diagonal and vanishes outside another sufficiently small neighbourhood of the diagonal.  
    If $T \in \PsiDO{0}{\sM; \sE}$ with $\ES T \subset \cU$ and $\fF^{\pm} := \fA (\chi \fF_{1}^{\ret, \adv}) \tilde{\fA} T$, then  
	\begin{subequations}
		\begin{eqnarray}
			&& 
			\WFPrime \fF^{\pm} \subset \varDelta \dotCoTansM \cup C^{\pm}, 
			\label{eq: Hoermander_Lem_26_1_15_i}
			\\ 
			&& 
            \square \fF^{\pm} = T + \fR^{\pm}, 
			\quad 
            \fR^{\pm} \in I^{-3/2} \big( \sM \times \sM, C^{\pm \prime}; \Hom{\sE, \sE} \big), 
			\label{eq: Hoermander_Lem_26_1_15_ii}
			\\ 
			&& 
            \fF^{+} - \fF^{-} \in I^{-3/2} \big( \sM \times \sM, C'; \Hom{\sE, \sE} \big).   
			\label{eq: Hoermander_Lem_26_1_15_iii} 
		\end{eqnarray}
	\end{subequations}
\end{lemma}
%
%
%

\begin{proof}
	The properties~\eqref{eq: Hoermander_Prop_26_1_2_i} and~\eqref{eq: Hoermander_Prop_26_1_2_ii} immediately entail~\eqref{eq: Hoermander_Lem_26_1_15_i} and~\eqref{eq: Hoermander_Lem_26_1_15_iii}. 
	To show~\eqref{eq: Hoermander_Lem_26_1_15_ii}, by definition:  
	\begin{equation} \label{eq: Hoermander_26_1_11}
        \square \fF^{\pm} 
		= (\square \fA - \fA \rD_{1}) (\chi \fF_{1}^{\ret, \adv}) \tilde{\fA} T + \fA \rD_{1} (\chi \fF_{1}^{\ret, \adv}) \tilde{\fA} T. 
	\end{equation} 
	Since $\square$ and $\rD_{1}$ are microlocal conjugates to one another,~\eqref{eq: microlocal_conjugate_NHOp} can be re-expressed as  
	\begin{equation} 
		(x_{0}, \xi^{0}; 0, \eta_{1} \rd y^{1}) \notin \WFPrime (\square \fA - \fA \rD_{1}). 
		\tag{\ref{eq: microlocal_conjugate_NHOp}$^{\prime}$}   
	\end{equation}
	Then, there is a conic neighbourhood $\cV$ of $(0, \eta_{1} \rd y^{1})$ such that $(\square \fA - \fA \rD_{1}) v$ is smooth for any $v \in \cD' (\Rd, \CN)$ when $\WF{v} \subset \cV$. 
	Since $\WF' (\chi \fF_{1}^{\ret, \adv})$ can be made arbitrary close to the diagonal in $\dotCoTanRd \times \dotCoTanRd$ by choosing the support of $\chi$ close to the diagonal in $\Rd \times \Rd$, we can pick an appropriate $\chi$ and a conic neighbourhood $\mathcal{W}$ of $(0, \eta_{1} \rd y^{1})$ such that $\WF \big( (\chi \fF_{1}^{\ret, \adv}) v \big) \subset \cV$ provided $\WF{v} \subset \mathcal{W}$. 
	If $\ES{T} \subset \varkappa (\mathcal{W})$ then the first term in the right-hand side of~\eqref{eq: Hoermander_26_1_11} is smooth and we are left with the last term of that equation. 
	By definition: $\rD_{1} (\chi \fF_{1}^{\ret, \adv}) = I + \tilde{\fR}^{\pm}$ where 
    $\tilde{\fR}^{\pm} := \rD_{x^{1}} \big(\chi (x, y) \big) \, \fF_{1}^{\adv, \ret} 
    \in 
    I^{- 1/2} (\Rd \times \Rd, C'_{1}; \matk)$ (cf. Proposition~\ref{prop: causal_propagator_partial_derivative_FIO}).  
    Since $\fA \tilde{\fA} T - T = (\fA \tilde{\fA} - I) T$ is smooth as long as $\ES{T}$ is sufficiently close to $\xxiNot$, it follows that $\square \fF^{\pm} - T = \fR^{\pm}$ where $\fR^{\pm} - \fA \tilde{\fR}^{\pm} \tilde{\fA} T$ is smooth, which concludes the proof. 
\end{proof} 
%
%
%

To end the proof of Theorem~\ref{thm: exist_unique_Feynman_parametrix_NHOp}, we choose a locally finite covering $\{ \sU_{\alpha} \}$ of $\dotCoTansM$ by open cones $\sU_{\alpha}$ such that Lemma~\ref{lem: Lem_26_1_15_Hoermander} is applicable when $\ES T_{\alpha} \subset \sU_{\alpha}$. 
Denoting the projection of $\sU_{\alpha}$ on $\sM$ by $U_{\alpha}$ and picking  $\sU_{\alpha}$ so that $U_{\alpha}$ are also locally finite, we set 
\begin{equation}
	I = \sum_{\alpha} T_{\alpha}, 
	\quad 
	\ES{T_{\alpha}} \subset \sU_{\alpha}, 
	\quad 
	\supp{\fT_{\alpha}} \subset U_{\alpha} \times U_{\alpha}. 
\end{equation}
Then for every $\alpha$, $\fA (\chi \fF_{1}^{\ret, \adv}) \tilde{\fA} T_{\alpha}$ can be chosen according to the Lemma~\ref{lem: Lem_26_1_15_Hoermander} with $\supp \big( \fA (\chi \fF_{1}^{\ret, \adv}) \tilde{\fA} T_{\alpha} \big) \subset U_{\alpha} \times U_{\alpha}$, and hence the sum 
\begin{equation} 
	\fF^{\pm} := \sum_{\alpha} \fA (\chi \fF_{1}^{\ret, \adv}) \tilde{\fA} T_{\alpha} 
\end{equation}
is well-defined which satisfies the claimed properties~\eqref{eq: def_Feynman_parametrix} and~\eqref{eq: diff_Feyn_anti_Feyn_parametrix} in the statement of Theorem~\ref{thm: exist_unique_Feynman_parametrix_NHOp}.  
Note that, if $\chi$ is taken as a function of $y-z$ then $\chi \fF_{1}^{\ret, \adv}$ is a convolution by a measure of compact support and therefore it continuously maps $H^s_{\mathrm{c}}(\Rd, \CN)$ to $H^s_{\mathrm{c}}(\Rd, \CN)$.
By the mapping properties of pseudodifferential operators and Fourier integral operators~\cite[Thm. 4.3.1]{Hoermander_ActaMath_1971}   
(see also~\cite[Cor. 25.3.2]{Hoermander_Springer_2009}), the factors $A, \tilde{A}, T_{\alpha}$ are 
continuous maps from $H^s_{\mathrm{c}}$ to $H^s_{\mathrm{c}}$. 
Since the sum of kernels is locally finite, this shows that
$E^{\Feyn, \aFeyn}$ continuously maps $H_{\mathrm{c}}^{s} (\sM; \sE)$ into $H_{\mathrm{loc}}^{s} (\sM; \sE)$ for all $s \in \R$.  
\newline 

Until now, we have just shown that $\square \fF^{\pm} = I + \fR^{\pm}$. 
By Lemma~\ref{lem: Hoermander_lem_26_1_16} (stated and proven just after the on going proof), we can choose $\fG^{\pm} \in I^{-3/2} \big( \sM \times \sM, C^{\pm \prime}; \Hom{\sE, \sE} \big)$ so that $\square \fG^{\pm} - \fR^{\pm}$ is smooth. 
Moreover, $G^{\pm}$ extend to continuous mappings from $H_{\mathrm{c}}^{s} (\sM; \sE)$ to $H_{\mathrm{loc}}^{s} (\sM; \sE)$ for every $s \in \R$. 
This follows from mapping properties of Fourier integral operators, for example, Theorem 25.3.8 in~\cite{Hoermander_Springer_2009},  
bearing in mind that the corank of the symplectic form $\sigmaup_{\ms \varGamma}$ (see~\eqref{eq: symplectic_form_graph_symplecto}) on $\varGamma$  is two. 
Since locally a Fourier integral operator on a bundle is a matrix of scalar Fourier integral operators (Section~\ref{sec: Lagrangian_distribution}), this theorem can be applied here directly.
Hence, 
\begin{equation} \label{eq: Feynman_parametrices_NHOp}
    E^{\Feyn, \aFeyn} := F^{\pm} - G^{\pm}
\end{equation}
is a right parametrix which has this continuity property. 
\newline 

Furthermore, the construction shows that $\fF^{+} - \fF^{-}$ and thus $\parametrixKernelFeyn - \parametrixKernelAntiFeyn$ is non-characteri-stic on $\varDelta \, \coLightBun$, as can be seen from the integral representation (Example~\ref{ex: advanced_retarded_fundamental_sol_partial_derivative}) of $F_{1}$ and wavefront set properties of $\fA, \tilde{\fA}$.  
Since $\square (\parametrixKernelFeyn - \parametrixKernelAntiFeyn)$ is smooth, it follows that its principal symbol satisfies a first-order homogeneous differential equation along the lightlike geodesic strip by Theorem~\ref{thm: Hoermander_thm_25_2_4} and Remark~\ref{rem: Weitzenboeck_connection}.   
Therefore $\parametrixKernelFeyn - \parametrixKernelAntiFeyn$ is non-characteristic everywhere: $\WF' (\parametrixKernelFeyn - \parametrixKernelAntiFeyn) = C$ and then it can be concluded that $\WF' \fE^{\Feyn, \aFeyn} \supset C^{\pm}$ because of $\WF' \fE^{\Feyn, \aFeyn} \subset \varDelta \, \dotCoTansM \cup C^{\pm}$. 
Finally, the fact 
\begin{equation}
	\varDelta \, \dotCoTansM = \WFPrime{\fI} = \WFPrime (\square \fE^{\Feyn, \aFeyn}) \subset \WFPrime \fE^{\Feyn, \aFeyn}
\end{equation}
together with Lemma~\ref{lem: Hoermander_lem_26_1_16} complete the proof. 

%
%
%
\begin{lemma} \label{lem: Hoermander_lem_26_1_16}
	Let $\sE \to \sM$ be a vector bundle over a Lorentzian manifold $\spacetime$ and $\square$ a normally hyperbolic operator on $\sE$. 
	Suppose that $C^{\pm}$ are the forward/backward geodesic relations. 
	If $\fB \in I^{m} \big( \sM \times \sM, C^{\pm \prime}; \Hom{\sE, \sE} \big)$ then one can find a properly supported $\fA \in I^{m} \big( \sM \times \sM, C^{\pm \prime}; \Hom{\sE, \sE} \big)$ such that $\square \fA - \fB$ is smoothing.   
\end{lemma}
%
%
%

\begin{proof}
	The assertion has been proven for scalar pseudodifferential operators~\cite[Thm. 5.3.2]{Duistermaat_ActaMath_1972} 
    (see also~\cite[Lem. 26.1.16]{Hoermander_Springer_2009}) 
    which we tailor for normally hyperbolic operators. 
    Let $\fA_{0}$ be a properly supported Lagrangian distribution of order $m$ on $\sM \times \sM$ whose principal symbol is $\fa_{0}$ and that of $\fB$ is $\fb$. 
	Since the principal symbol of $\square$ vanishes on lightlike covectors, Theorem~\ref{thm: Hoermander_thm_25_2_4} implies that $\square \fA_{0} - \fB$ is of order $m-1$ on lightlike covectors on $\sM$, provided    
	\begin{equation}
		\left( - \ri \, \pounds_{X_{\fg}} + \subSymb{\square} \right) \fa_{0} = \fb
	\end{equation}
	is satisfied. 
    Identifying half-densities with functions according to Remark~\ref{exm: density}~\ref{exm: density_Lorentzian_mf} and making use of the fact that the Keller-Maslov bundle $\Maslov$ is a trivial complex vector bundle (see Definition~\ref{def: Keller_Maslov_bundle}), we write (cf.~\eqref{eq: symbol_halfdensity_function}) $\fa_{0} = a_{0} \m$ and $\fb = b \m$ where $\m$ is a non-vanishing section of $\Maslov$ and $a_{0}, b$ are scalar symbols of degree $m$. 
    The preceding equation then reads $(- \ri X_{\fg} + \subSymb{\square}) a_{0} = b$ by~\eqref{eq: Lie_derivative_density_coordinate}. 
	In other words, inserting~\eqref{eq: HVF_NHOp} and~\eqref{eq: subprincipal_symbol_NHOp} yields 
	\begin{equation} 
        \left( - \frac{\partial \fg^{\mu \nu}}{\partial x^{i}} \xi_\mu \xi_\nu \frac{\partial}{\partial \xi_{i}} 
        + 2 \fg^{ij} \xi_{j} \frac{\partial}{\partial x^{i}} 
        + 2 g^{ij} \Gamma_{i}^i \xi_{j}  \right) a_{0} = b.  
	\end{equation}
	This equation has a unique solution for a given initial condition and the essential support of the solution is contained in that of $b$, that is $\WFPrime{\fA_{0}} \subset C^{\pm}$. 
	Now the conclusion follows from applying the argument iteratively to construct Lagrangian distributions $\fA_{k}$ of order $m - k$ such that $\square (\fA_{0} + \ldots + \fA_{k}) - \fB \in I^{m-k-1} \big( \sM \times \sM, C^{\pm \prime}; \Hom{\sE, \sE}\big)$ and setting $\fA :\sim \fA_{0} + \ldots + \fA_{k} + \ldots$. 
\end{proof} 
%
%
%
%
%
%
%
%
%
%
\subsection{Proof of Proposition~\ref{prop: positivity_Feynman_minus_adv_NHOp}} 
Taking the adjoint $\parametrixKernelFeyn {}^{*} \square \equiv I$ of the Feynman parametrix $\parametrixKernelFeyn$ we find that $\WFPrime{\parametrixKernelFeyn {}^{*}} \subset \varDelta \, \dotCoTansM \cup C^{-}$ by Definition~\ref{def: Feynman_parametrix}.  
Thus $\parametrixKernelFeyn {}^{*} - \parametrixKernelAntiFeyn$ is smooth and so is $(\parametrixKernelFeyn - \GreenOpAdvKernel) (\parametrixKernelFeyn - \GreenOpAdvKernel)^{*}$. 
Since $\WF' \big(- \ri (\parametrixKernelFeyn - \GreenOpAdvKernel) \big) \cap \WF' \big( \ri (\parametrixKernelFeyn - \GreenOpAdvKernel)^{*} \big) = \emptyset$, each of these is in $C^{\infty} \big( \sM \times \sM; \Hom{\sE, \sE} \big)$. 
Therefore $- \ri (\parametrixKernelFeyn - \GreenOpAdvKernel)$ can be normalised so that it is symmetric.  
\newline 
%
%
%

To address the positivity, first consider the operator $\rD_{1}$ as in Section~\ref{sec: microlocalisation_NHOp} and observe that $\ri F_{1}$ is positive as evident from Example~\ref{ex: advanced_retarded_fundamental_sol_partial_derivative}:   
\begin{equation} \label{eq: positivity_causal_fundamental_sol_partial_derivative}
    \ri \fF_{1} (\bar{u}^{*} \otimes u) \geq 0. 
\end{equation}

To account $\hat{W}$, one notes that the Lagrangian distributions $\fA$ and $\fA^{*}$ in Theorem~\ref{thm: microlocalisation_NHOp} preserve wave front sets and thus, one obtains $(x_{0}, \xi^{0}; 0, \eta_{1} \rd y^{1}) \notin \WFPrime (\hat{\fW} - \fA \fF_{1} \fA^{*})$ from parametrix differences.  
But it is not sufficient to construct a positive $\hat{\fW}$ on a conic neighbourhood in $\dotCoTansM \times \dotCoTansM$ of the inverse image in $C$ of a neighbourhood of any point in $\coLightBun / C$. 
This compels a refined version microlocalisation and hence respective homogeneous Darboux theorem. 
\newline 

The space of lightlike (unparametrised) geodesic strips is a $2d - 3$-dimensional smooth conic manifold which can be naturally identified with the cotangent unit sphere bundle $\bbS^{*} \varSigma$ of any Cauchy hypersurface $\varSigma$ with respect to the induced metric $\fh$ on $\varSigma$~\cite{Low_NonlinearAnal_2001}. 
We choose $L$, an elliptic symmetric operator of order $-1/2$ with scalar principal symbol and vanishing subprincipal symbol as in Theorem~\ref{thm: microlocalisation_NHOp}. 
One can always choose the principal symbol of $L$ such that $X_{L \square L}$ is a complete vector field~\cite[Thm. 6.4.3, p. 234]{Duistermaat_ActaMath_1972} 
(see also~\cite[Lem. 26.1.11, 26.1.12]{Hoermander_Springer_2009} 
and~\cite[Prop. 2.1]{Strohmaier_AdvMath_2021}).  
Since the bundle of forward/backward lightlike covectors $\dotCoTan_{0,\pm} \sM$ are connected components of the lightcone bundle $\coLightBun$, there exists a homogeneous diffeomorphism $\coLightBun \to (\bbS^{*} \varSigma \oplus \bbS^{*} \varSigma) \times \R \times \R_{+}$ such that $X_{L \square L}$ is mapped to the vector field $\partial / \partial s$ when the variables in $\bbS^{*} \varSigma \times \R \times \R_{+}$ are denoted by $(\cdot, s, \cdot)$. 
In this setting, a refined Darboux theorem yields that, for every lightlike covector $\xxiNot$ on $\sM$, there is an open conic coordinate chart $\big( \cU, (x^{i}, \xi_{i}) \big)$ and a homogeneous symplectomorphism   
\begin{equation} \label{eq: def_bijective_symplecto}
	\kappa : \cU \to \cU', \; \xxi \mapsto 
    \kappa \xxi := \big( y^{1} \xxi, \ldots, y^{d} \xxi; \eta_{1}  \xxi, \ldots, \eta_{d} \xxi \big), 
\end{equation}
bijectively mapping $\cU$ to an open conic neighbourhood $\cU'$ of $\kappa \xxiNot := (0, \eta_{1} \rd y^{1})$ in $\dotCoTanRd$ such that 
\begin{equation}
	\big( \kappa^{-1} \big)^{*} p = \xi_{1} \, \one_{\End{\sE}}, 
\end{equation}
$\Char{\rD_{1}} := \{ \yeta \in \Rd \times \dot{\R}^{d} \,|\, y^{1} = 0 \}$ is symmetric with respect to the plane $y^{1} = 0$ and convex in the direction of $y^{1}$-axis, and $\cU' \cap \Char{\rD_{1}}$ is invariant under the translation along the $y^{1}$-axis. 
Since $\rd \kappa^{-1} (X_{L^{*} \square L}) = \one_{\matk} \partial / \partial y^{1}$, it follows that $\cU \cap \Char{\square}$ is invariant under $X_{\square}$~\cite[Lem. 6.6.3]{Duistermaat_ActaMath_1972}. 
\newline 

Suppose that $\cV_{1}, \cV_{2}$ are closed conic neighbourhoods of $\xxiNot$ such that $\cV_{1} \subset \cU$ and $\cV_{2} \subset \mathrm{int} (\cV_{1})$ while $\cV_{\iota} \cap \coLightBun, \iota = 1, 2$ are invariant under the geodesic flow and that $\cV'_{\iota} := \kappa (\cV_{\iota})$ are convex in the $y^{1}$-direction and symmetric about the plane $y^{1} = 0$. 
Remark~\ref{rem: global_microlocalisation_NHOp} then entails that there exists a properly supported Lagrangian distribution $\fA \in I^{-1/2} \big( \sM \times \Rd, \varGamma'_{1}; \Hom{\CN, \sE} \big)$ so that $A^{*} A$ and $A A^{*}$ both are pseudodifferential operators of order $-1$ and  
\begin{subequations} \label{eq: microlocalisation_NHOp_refined}
	\begin{eqnarray}
		&& 
		\Char{\fA} \cap \varGamma'_{2} = \emptyset,
		\label{eq: microlocalisation_NHOp_1_refined} 
		\\ 
		&& 
		\ES (A^{*} A - I) \cap \cV'_{2} = \emptyset  
		\Leftrightarrow
		\ES (A A^{*} - I_{\sE}) \cap \cV_{2} = \emptyset, 
		\label{eq: microlocalisation_NHOp_2_refined}
		\\ 
		&& 
		\WFPrime (\square \fA - \fA \rD_{1}) \cap \varGamma_{2} = \emptyset, 
		\label{eq: microlocal_conjugate_NHOp_refined}
	\end{eqnarray}
\end{subequations} 
where $\varGamma_{1}$ and $\varGamma_{2}$ are the graphs of the restriction of $\kappa$ to $\cV_{1}$ and to $\cV_{2}$, respectively.  
\newline 

It follows from Theorem~\ref{thm: exist_unique_Feynman_parametrix_NHOp} and~\eqref{eq: microlocalisation_NHOp_refined} that we can microlocally conjugate $\hat{\fW}$ to $\ri \fF_{1}$ 
(see~\cite[p. 237]{Duistermaat_ActaMath_1972} for the scalar version):  
\begin{equation} \label{eq: FIO_6_6_14_bundle}
    \WFPrime (\hat{\fW} - \fA \ri \fF_{1} \fA^{*}) \cap (\cV_{2} \times \cV_{2}) = \emptyset, 
	\quad 
    \cV_{2} \cap \coLightBun \subset \{\dotCoTan_{0, \pm} \sM\}.  
\end{equation} 

%
%
%
\begin{lemma} \label{lem: FIO_lem_6_6_5_bundle}
    For a normally hyperbolic operator $\square$ as in Proposition~\ref{prop: positivity_Feynman_minus_adv_NHOp} and Lagrangian distributions $\fA$ given by~\eqref{eq: microlocalisation_NHOp_refined}, let $\varPsi \in \PsiDO{0}{\sM; \sE}$ such that $\ES{[\varPsi, \square]_{-}} \cap \coLightBun = \emptyset, \ES{\varPsi} \subset \cV_{2}$ and $\coLightBun \cap \cV_{2} \subset \dotCoTan_{0, \pm} \sM$. 
    Then 
	\begin{equation}
        \hat{\fW} \varPsi \varPsi^{*} \equiv \varPsi \fA \ri \fF_{1} \fA^{*} \varPsi^{*} 
		\Leftrightarrow 
		\hat{\fW} \varPsi \varPsi^{*} = \varPsi \fA \ri \fF_{1} \fA^{*} \varPsi^{*} \mod C^{\infty} \big( \sM \times \sM; \Hom{\sE, \sE} \big).
	\end{equation}
\end{lemma}
%
%
%

\begin{proof}
	The idea of the proof is to use microlocalisation to reduce this to the special case $\rD_1$.
	To reach the conclusion our main task is then to compute the relevant commutators which is straightforward and exactly analogous to the scalar counterpart~\cite[Lem. 6.6.5]{Duistermaat_ActaMath_1972}. 
	So, we only sketch the main steps for completeness. 
	The ellipticity~\eqref{eq: microlocalisation_NHOp_1_refined} allows one to relate any properly supported $\varPsi \in \PsiDO{0}{\sM; \sE}$ and $\varPsi' \in \PsiDON{0}$ with $\ES{\varPsi} \subset \cV_{2}$ and $\ES{\varPsi'} \subset \cV'_{2}$ by 
	\begin{equation}
		\varPsi' = A^{*} \varPsi A \mod \PsiDON{-\infty}, 
		\quad 
		\varPsi = A \varPsi' A^{*} \mod \PsiDO{-\infty}{\sM; \sE}. 
	\end{equation} 
    A direct computation shows $\ES{[\varPsi', \rD_{1}]_{-}} \cap \Char{\rD_{1}} = \emptyset$ which implies that the derivative of the symbol of $\varPsi'$ with respect to $y^{1}$ is of order $-\infty$ in a neighbourhood of $\Char{\rD_{1}}$. 
    Denoting the convolution by the Dirac measure at $(h, 0, \ldots, 0)$ by $\tau_{h}$, we can rewrite $\fF_{1}$ as $\fF_{1} = \int \tau_{h} \rd h$ and observe that $\tau_{h} \varPsi'  - \varPsi' \tau_{h}$ is of order $-\infty$ in a conic neighbourhood of $\Char{\rD_{1}}$. 
	Since $\Char{\rD_{1}}$ is invariant under translation in the $y^{1}$-direction, one has  
	\begin{equation}
        \forall v \in \cE' (\Rd, \CN) : \WF ([\tau_{h}, \varPsi']_{-} v) \cap \Char{\rD_{1}} = \emptyset. 
	\end{equation}
	The result follows after the integration with respect to $h$ and a few straightforward algebraic manipulations.
\end{proof}
%
%
%
	
We observe that $\varPsi \fA \ri \fF_{1} \fA^{*} \varPsi^{*}$ is non-negative by~\eqref{eq: positivity_causal_fundamental_sol_partial_derivative} since $A A^{*}$ is non-negative with respect to the hermitian form on $\sE$. 
Thus the only thing left to conclude Proposition~\ref{prop: positivity_Feynman_minus_adv_NHOp} is that the identity can be expressed as a sum of operators of the form $\varPsi \varPsi^{*}$ discussed in Lemma~\ref{lem: FIO_lem_6_6_5_bundle}, provided the prerequisites in this lemma are satisfied.  
Aiming to show this, we consider a closed conic neighbourhood $\cV_{3}$ of any lightlike covector $\xxiNot$ where $\cV_{3} \subset \mathrm{int} (\cV_{2})$ and $\cV_{3} \cap \coLightBun$ is invariant under the geodesic flow.
One can prove that there exists $\varPsi \in \PsiDO{0}{\sM; \sE}$ with $\ES{\varPsi} \subset \cV_{2}$ such that $\varPsi$ has a non-negative principal symbol equals to $\one$ on $\cV_{3}$ and $\ES{[\varPsi, \square]_{-}} \cap \coLightBun = \emptyset$~\cite[Lem. 6.6.6]{Duistermaat_ActaMath_1972}.
In other words, the hypotheses in Lemma~\ref{lem: FIO_lem_6_6_5_bundle} can be satisfied for $\varPsi$.
Taking a suitable cover of $\coLightBun$ one can use these operators to construct a family of operators $\varPsi_{\alpha} \in \PsiDO{0}{\sM; \sE}$ satisfying
\begin{subequations}
	\begin{eqnarray}
		&& \ES{\varPsi_{\alpha}} \bigcap \coLightBun \subset \dotCoTan_{0, \pm} \sM, 
		\label{eq: WF_POU_1}
		\\ 
		&& \mathrm{ES} \bigg( I - \sum_{\alpha} \varPsi_{\alpha} \varPsi_{\alpha}^{*} \bigg) \bigcap  \dotCoTan_{0, \pm} \sM = \emptyset, 
		\label{eq: WF_POU_2}
		\\ 
        && \ES{[\square, \varPsi_{\alpha}]_{-}} \bigcap \coLightBun = \emptyset. 
		\label{eq: WF_commutator_NHOp_POU} 
	\end{eqnarray}
\end{subequations}
The construction of such a family that provides a microlocal partition of unity (see the Footnote~\ref{foot: microlocal_partition_unity} in the proof of Theorem~\ref{thm: existence_parametrix_FIO} for this concept) is carried out in detail in~\cite[p. 238]{Duistermaat_ActaMath_1972}. 
The proof follows the usual strategy of summing the operators and multiplying from left and right by the root of a parametrix. 
One can verify that the hypotheses in Lemma~\ref{lem: FIO_lem_6_6_5_bundle} are preserved under this construction. 
To finalise the proof of Proposition~\ref{prop: positivity_Feynman_minus_adv_NHOp}, we simply write 
\begin{equation} \label{eq: Feynman_parametrix_P_mod_partial_deri}
	\hat{\fW} 
    = \sum_{\alpha} \hat{\fW} \varPsi_{\alpha} \varPsi_{\alpha}^{*} 
    \equiv \sum_{\alpha} \varPsi_{\alpha} \fA_{\alpha} \ri \fF_{1} \fA_{\alpha}^{*} \varPsi_{\alpha}^{*}, 
\end{equation}
which terminates the proof as every term on the right hand side of~\eqref{eq: Feynman_parametrix_P_mod_partial_deri} is non-negative with respect to the hermitian form on $\sE$. 
%
%
%
%
%
%
%
%
%
%
\subsection{Feynman propagators from Feynman parametrices}
We postpone to proof Theorem~\ref{thm: existence_Feynman_propagator_NHOp} in order to make the following observations. 
Let us denote a point in $\sM$ as $(t, x) \in \R \times \varSigma$. 
Then we can define the distributions $f_{1} \otimes \delta_{\varSigma}, f_{0} \otimes \delta'_{\varSigma} \in \cE'(\sM; \sE)$ by
\begin{equation} \label{eq: def_distributional_Cauchy_data_NHOp}
	f_{1} \otimes \delta_\varSigma (\phi)  
	:= 
	\int_\varSigma \phi (x) \, f_{1} (x) \, \dVolh (x), 
	\quad  
	f_{0} \otimes \delta'_\varSigma (\phi)  
	:= 
	\int_\varSigma (\partial_{\rN} \phi) (x) \, f_{0} (x) \, \dVolh (x), 
\end{equation}
where $\phi \in C_{\mathrm{c}}^{\infty} (\sM; \bar{\sE}^{*})$, $\rN$ is a future directed unit normal vector field on $\sM$ along $\varSigma$, and $\dVolh$ is the Riemannian volume form. 
We denote by $\cW_\varSigma \subset \cE' (\sM; \sE)$, the span of the set of distributions of this form. 
By a duality argument both retarded and advanced Green's operators extend to continuous maps $\mathcal{E}' (\sM; \sE) \to \mathcal{D}' (\sM; \sE)$, and if $G$ is the causal propagator and
$f := f_{0} \otimes \delta'_\varSigma + f_{1} \otimes \delta_\varSigma$, then $u = G(f)$ is a smooth solution of $\Box u=0$ with Cauchy data $(f_0, f_1)$ on $\varSigma$ 
(see e.g.~\cite[Lem 3.2.2]{Baer_EMS_2007},~\cite[Thm. 3.20]{Jubin_LMP_2016}). 
\newline 

Assume that $\fw \Box$ and $\bar{\Box}^{*} \fw$ are smooth for any bidistribution $\fw \in \cD' (\sM \times \sM; \sE \boxtimes \bar{\sE}^*)$.  
Let us denote by $\cE'_{\varSigma^{\perp *}} (\sM; \sE)$, the set of compactly supported distributions with wavefront set contained in the conormal bundle (see Example~\ref{exm: conormal_bundle_Lagrangian_submf}) $\varSigma^{\perp *}$ of $\varSigma$. 
This space is endowed with a natural 
locally convex topology~\cite[p. 125]{Hoermander_ActaMath_1971} 
with respect to which $C_{\mathrm{c}}^{\infty} (\sM; \sE)$ is sequentially dense in  $\cE'_{\varSigma^{\perp *}} (\sM; \sE)$ 
(see e.g.~\cite[Thm. 8.2.3]{Hoermander_Springer_2003}  
or the 
exposition~\cite[Sec. 4.3.1]{Strohmaier_Springer_2009} 
for details). 
The bidistribution $\fw$ can be defined as a sequentially continuous bilinear form on $\cE'_{\varSigma^{\perp *}} (\sM; \sE)$. 
Since $\cW_\varSigma \subset \cE'_{\varSigma^{\perp *}} (\sM; \sE)$, the form $\fw$ is then also defined on $\cW_\varSigma$.

%
%
%
\begin{lemma} \label{lem: positive_bisolution_NHOp}
	Let $\big( \sE \to \sM, (\cdot|\cdot) \big)$ be a hermitian vector bundle  over a globally hyperbolic manifold $(M, g)$ whose Cauchy hypersurface is $\varSigma$. 
	Suppose that $\square$ is a normally hyperbolic operator on $\sE$ and that $\fw \in \cD'(\sM \times \sM, \sE \boxtimes \bar{\sE}^*)$ is a bisolution of $\square$. 
	Then, $\fw$ is non-negative with respect to $(\cdot|\cdot)$ if and only if $\fw (\bar{f}^{*} \otimes f) \geq 0$ for any $f \in \cW_\varSigma$, where $\cW_{\varSigma}$ is the span of the distribution Cauchy data~\eqref{eq: def_distributional_Cauchy_data_NHOp} of $\square$.
\end{lemma}
%
%
%

\begin{proof}
    We begin with the fact that $\fw$ extends to a separately sequentially continuous bilinear form on $\cE'_{\varSigma^{\perp *}} (\sM; \bar{\sE}^*) \times \cE'_{\varSigma^{\perp *}} (\sM; \sE)$ as it is a bisolution; let us denote the extension by the same letter. 
    As this form is a bisolution it vanishes on 
    $\bar{\Box}^* C^\infty_{\mathrm{c}} (\sM; \bar{\sE}^*) \times C^\infty_{\mathrm{c}} (\sM; \sE)$ 
 	and on 
 	$C^\infty_{\mathrm{c}} (\sM; \bar{\sE}^*) \times \Box C^\infty_{\mathrm{c}} (\sM; \sE)$. 
 	By sequential continuity and sequential density of $C^\infty_{\mathrm{c}}$ in $\cE'_{\varSigma^{\perp *}}$, it follows that $\fw$ vanishes on 
 	$\bar{\Box}^* \cE'_{\varSigma^{\perp *}} (\sM; \bar{\sE}^{*}) \times \cE'_{\varSigma^{\perp *}} (\sM; \sE)$ 
 	and on
 	$\cE'_{\varSigma^{\perp *}} (\sM; \bar{\sE}^{*}) \times \Box \cE'_{\varSigma^{\perp *}} (\sM; \sE)$. 
 	We will show that any 
 	$u \in \comSecsME$ can be written as $u = f + \Box \phi$ for some $f \in \cW_\varSigma$ and some $\phi \in \cE'_{\varSigma^{\perp *}} (\sM; \sE)$, and the analogous statement holds for the conjugate-dual bundle and the conjugate-dual normally hyperbolic operator. 
 	This obviously implies the statement. 
    \newline 

 	To see that this is indeed the case, observe that $\Box (G u) = 0$ and therefore, if $(f_0, f_1)$ is the Cauchy data of $G(u)$ on $\varSigma$ and $f$ the corresponding element in 
 	$\cW_\varSigma$, then we have $G (u) = G (f)$. 
 	This implies that 
    \begin{equation}
    	G (u - f) = 0 \Leftrightarrow \GreenOpRet (u - f) = \GreenOpAdv (u - f). 
    \end{equation}
 	Thus $\phi := \GreenOpRet (u - f)$ must be a compactly supported distribution because
    $\supp \phi = \supp \big( \GreenOpRet (u-f) \big) \cap \supp \big( \GreenOpAdv (u-f) \big) \subset J^{+} \big( \supp (u-f) \big) \cap J^{-} \big( \supp (u-f) \big)$ is compact due to the global hyperbolicity (Definition~\ref{def: globally_hyperbolic_spacetime}) of $\sM$. 
 	By propagation of singularities (Theorem~\ref{thm: propagation_singularity_Sobolev_WF}) this distribution must be again in $\cE'_{\varSigma^{\perp *}} (\sM; \sE)$ and we have $u - f = \Box \phi$. 
 	The same proof works for the normally hyperbolic operator $\bar{\Box}^*$ on $\bar{\sE}^*$.
\end{proof}
%
%
%

Finally, we are prepared to show the existence of Feynman propagators for $\square$. 
Moreover, if there exists a hermitian form on $\sE$ then Feynman propagators can be chosen non-negative with respect to this form and hence implying existence of $2$-point bidistributions corresponding to Hadamard states. 

%
%
%
\begin{proof}[Proof of Theorem~\ref{thm: existence_Feynman_propagator_NHOp}]
    Let $\hat{W}$ be the parametrix of $\square$ constructed in Proposition~\ref{prop: positivity_Feynman_minus_adv_NHOp}. 
    Since the wavefront set of $\hat{\fW}$ does not intersect with the conormal bundle $\varSigma^{\perp *}$ of any Cauchy hypersurface $\varSigma$ in $\sM$, we fix a $\varSigma$ and restrict $\hat{\fW}$ and its normal derivative to create distributional Cauchy data on $\varSigma \times \varSigma$.
    Suppose that $\sE$ is endowed with a hermitian form $(\cdot|\cdot)$ such that $\hat{\fW}$ is non-negative with respect to $(\cdot|\cdot)$.  
    Then the restriction $\hat{\fW}_{\varSigma}$ of $\hat{\fW}$ on $\varSigma$ is non-negative as well. 
    Next we construct $\fW \in \cD' (\sM \times \sM; \sE \boxtimes \bar{\sE}^*)$ employing the solution operator of $\square$ with the same distributional Cauchy data as $\hat{\fW}$ on $\varSigma \times \varSigma$. 
 	This will be a bisolution of $\square$ with the property 
 	\begin{equation}
        \left( \forall u \in \mathcal{E}'_{\varSigma^{\perp *}} (\sM; \sE) : \fW (\bar{u}^{*} \otimes u) \geq 0 \right)
        \Leftrightarrow W \geq 0, 	
	\end{equation}
 	by the sequential density of $C_{\mathrm{c}}^{\infty} (\varSigma; \sE_{\varSigma})$ in  $\mathcal{E}'_{\varSigma^{\perp *}} (\sM; \sE)$ and Lemma~\ref{lem: positive_bisolution_NHOp}.
    Note that $\fR = \fW - \hat{\fW}$ is smooth since it solves the inhomogeneous problem with zero Cauchy data, which entails that 
    \begin{equation}
    	\GreenOpFeyn := \parametrixFeyn + \ri R
    \end{equation}
	is the \textit{Feynman propagator} with the required \textit{positivity} property.
\end{proof}
%
%
%
%
%
%
%
%
%
%
\subsection{Causal propagator} 
\label{sec: causal_propagator_NHOp_FIO}
In this section we will inscribe the causal propagator $G := \eqref{eq: def_causal_propagator}$for a normally hyperbolic operator $\square$ on a vector bundle $\sE \to \sM$ over a globally hyperbolic spacetime as a Fourier integral operator. 
The result has been proven by 
Duistermaat-H\"{o}rmander~\cite[Thm. 6.6.1]{Duistermaat_ActaMath_1972}  
for any scalar pseudodifferential operator of real principal type, but a systematic derivation for the bundle version is not so easy to find.  

%
%
%
\begin{remark} \label{rem: volume_geodesic_relation}
    There is a natural density $\dVol_{\ms C}$ on the geodesic relation $C$ and hence on the forward/backward geodesic relations $C^{\pm}$, as originally due to 
    Duistermaat and H\"{o}rmander~\cite[p. 230]{Duistermaat_ActaMath_1972} 
    for a generic manifold $M$, which simplifies considerably for a globally hyperbolic spacetime $\spacetime$ as reported by 
    Strohmaier and Zelditch~\cite[$(52)$, Rem. 7.1]{Strohmaier_AdvMath_2021}. 
    By Definition~\ref{def: geodesic_relation}: for each $\xxiyeta \in C$ there is a unique $s \in \R$ such that $\xxi = \varPhi_{s} \yeta$ so that $C$ can be identified with an open subset of $\R \times \coLightBun$, where $\varPhi_{s}$ is the geodesic flow on $\coTansM$ and $s \in \R$ is the flow parameter. 
    $\coLightBun$ is a conic contact manifold where the Hamiltonian  
    \begin{equation} \label{eq: def_Hamiltonian_metric}
        H_{\fg} : C^{\infty}(\sM; \coTansM) \to \R, ~ \xxi \mapsto H_{\fg} \xxi := \frac{1}{2} \fg_{x}^{-1} (\xi, \xi) 
    \end{equation}
	induced by the spacetime metric $\fg$ vanishes identically and its Hamiltonian reduction is the the manifold of \textbf{scaled-lightlike geodesics} $\sN$~\cite[pp. 10-12]{Penrose_1972},~\cite[Thm. 2.1]{Khesin_AdvMath_2009}. 
	Following Strohmaier-Zelditch, we have used the prefactor $1/2$ in the preceding definition so that the \textit{relativistic Hamiltonian flow} is \textit{identical} to the \textit{geodesic flow}. 
    Denoting by $\tilde{s}$, the dilation parameter on $\coLightBun$, the natural half-density on $C$ is given by 
    \begin{equation}
        \sqrt{|\dVol_{\ms C}|} := \sqrt{|\rd s|} \otimes \sqrt{|\rd \tilde{s}|} \otimes \sqrt{| \dVol_{\ms \sN}|}. 
    \end{equation}  
    Note, this density differs from that by Duistermaat-H\"{o}rmander by a factor of $2$ because they used the Hamiltonian flow of $\fg^{-1}$ to parametrise $\sN$, in contrast to the flow of the Hamiltonian vector field $X_{\fg/2}$ generated by $H_{\fg}$. 
    Moreover, 
    \begin{equation}
        \pounds_{X_{\fg}} \dVol_{\ms C} = 0. 
    \end{equation}
    The densities on $C^{\pm}$ follow from the fact $\coLightBun = \dotCoTan_{0, +} \sM \sqcup \dotCoTan_{0, -} \sM$ in $d \geq 3$ and the case $d=2$ is as mentioned in Remark~\ref{rem: sixteen_parametrix_dim_two_NHOp}. 
\end{remark}
%
%
%

Employing 
Theorems~\ref{thm: exist_unique_Feynman_parametrix_NHOp} and~\ref{thm: Hoermander_thm_25_2_4}, 
Definition~\ref{def: P_compatible_connection}, 
and 
Remark~\ref{rem: volume_geodesic_relation}, 
we imitate the proof for the scalar 
version~\cite[Thm. 6.6.1]{Duistermaat_ActaMath_1972},  
to obtain

%
%
%
\begin{theorem} \label{thm: causal_propagator_NHOp_FIO}
    Let $\sE \to \sM$ resp. $\pi: \dotCoTansM \to \sM$ be a vector bundle resp. the punctured cotangent bundle over a globally hyperbolic spacetime $\spacetime$ and $\square$ a normally hyperbolic operator on $\sE$ whose principal symbol is the spacetime metric $\fg^{-1}$ on $\coTansM$.  
	Then, the Schwartz kernel $\fG$ of the causal propagator for $\square$ is the Lagrangian distribution 
	\begin{subequations}
		\begin{eqnarray}
			&& \fG \in I^{-3/2} \big( \sM \times \sM, C'; \Hom{\sE, \sE} \big), 
			\\ 
			&& \symb{\fG} = \frac{\ri}{2} \sqrt{2 \pi |\rd \mathsf{v}_{\ms C}|} \otimes \m w, 
			\\ 
			&& \Char{\fG} \cap C = \emptyset,  
		\end{eqnarray}
	\end{subequations}
    where $\dVol_{\ms C}$ is the natural volume form on the geodesic relation $C$, $\m$ is a section of the Keller-Maslov bundle $\Maslov \to C$ 
    (as constructed in~\cite[pp. 231-232]{Duistermaat_ActaMath_1972}), and $w$ is the unique element of $C^{\infty} \big( C; \pi^{*} \Hom{\sE, \sE}_{C} \big)$ that is diagonally the identity endomorphism and off-diagonally covariantly constant  
	\begin{equation} \label{eq: symbol_causal_propagator_NHOp_covarinatly_constant}
		\nabla_{X_{\fg / 2}}^{\ms \pi^{*} \Hom{\sE, \sE}} w = 0    
	\end{equation}
    with respect to the $\square$-compatible Weitzenb\"{o}ck covariant derivative $\nabla_{X_{\fg / 2}}^{\ms \pi^{*} \Hom{\sE, \sE}}$ along the geodesic vector field $X_{\fg / 2}$. 
\end{theorem}
%
%
%

\begin{remark} \label{rem: HVF_canonical_relation}
    By Definition~\ref{def: P_compatible_connection}, $X_{\fg}$ acts on $C^{\infty} \big( \coTansM; \pi^{*} \Hom{\sE, \sE} \big)$. 
	Thus~\eqref{eq: symbol_causal_propagator_NHOp_covarinatly_constant} 
    must be interpreted in the sense of the induced Hamiltonian vector field on $C$ by the vector field $(X_{\fg}, 0)$ on $(\coTansM \times \coTansM)$~\cite[Rem. 1, p. 216]{Duistermaat_ActaMath_1972} 
	(see also~\cite[Rem. 1, p. 69]{Hoermander_Springer_2009}).   
\end{remark}
%
%
%
%
%
%
%
%
%
%
\section{Feynman propagators for a Dirac operator}
In this section we will extend our analysis for Dirac-type operators. 
These are the first-order differential operators whose squares are normally hyperbolic, and arguably the most fundamental first-order operators in geometric analysis as studied in the context of the celebrated Atiyah-Singer index theorem. 
%
%
%
%
%
%
%
%
%
%
\subsection{Dirac-type operators}
\label{sec: Dirac_type_op}
%
%
%
\begin{definition} \label{def: Dirac_type_op}
	Let $\sE \to \sM$ be a vector bundle over a Lorentzian manifold $\spacetime$. 
	A \textbf{Dirac-type operator} acting on half-density-valued smooth sections of $\sE$ is a first-order differential operator  
    \begin{equation}
    	D : \secsME \to \secsME
    \end{equation}
    whose principal symbol $\symb{D}$ satisfies the Clifford relation 
    (see e.g.~\cite[Def. 2.18]{Baer_Springer_2012}) 
	\begin{equation} 
        \symb{D} \xxi^{2} = \fg_{x}^{-1} (\xi, \xi) \one_{\End \sE_{x}}  
	\end{equation}
	for all $\xxi \in C^{\infty} (\sM; \coTansM)$. 
\end{definition}
%
%
%

In other words, $D^{2}$ is a normally hyperbolic operator and by the polarisation identity 
(see e.g.~\cite[Rem. 2.19]{Baer_Springer_2012}) 
\begin{equation} \label{eq: anticommutation_symbol_Dirac_type}
	[\symb{D} (x, \xi), \symb{D} (x, \eta)]_{+} = 2 \, \fg_{x}^{-1} (\xi, \eta) \, \one_{\End{\sE}_{x}} 
\end{equation}
for any $\xxi, (x, \eta) \in C^{\infty} (\sM; \coTansM)$ and where $[\cdot, \cdot]_{+}$ is the anticommutator bracket. 
The principal symbol $\symb{D}$ defines a Clifford action of $\coTansM$ on $\sE$ by 
\begin{equation} \label{eq: def_symbol_Dirac_type_op}
	\forall f \in C^{\infty} (\sM) : 
	\symb{D} (\rd f) := \ri \, [D, f]_{-}, 
\end{equation}
which is independent of any chosen element $f \in C^{\infty} (\sM)$. 
This turns $\sE \to \sM$ into a bundle of Clifford modules $(\sE \to \sM, \symb{D})$ over $\spacetime$. 
Furthermore, the hindmost equation yields the Leibniz rule for $D$: 
\begin{equation}
	D (fu) = \symb{D} (\rd f) \, (u) + f D (u) 
\end{equation}
for any $f \in C^{\infty} (\sM)$ and any $u \in \secsME$. 
The Clifford mapping~\eqref{eq: anticommutation_symbol_Dirac_type} defines the (pointwise) Clifford multiplication  
\begin{equation} \label{eq: def_Clifford_multiplication}
    \fc: C^{\infty} (\sM; \coTansM \otimes \sE) \to \secsME, ~ \xxi \otimes u \mapsto 
    \fc \big( \xxi \otimes u \big) := \symb{D} \xxi \, (u). 
\end{equation}

%
%
%
\begin{assumption}\label{asp: Direc_op_symmetric}
	We endow the vector bundle $\sE \to \sM$ over a Lorentzian manifold $\spacetime$ with a sesquilinear form $(\cdot|\cdot)$ such that $D$ is symmetric. 
\end{assumption}  
%
%
%

Given any $D$, there exists a unique $(\cdot|\cdot)$-compatible Weitzenb\"{o}ck connection $\connectionE$ on $\sE$ and a unique potential $V \in C^\infty(\sM; \End \sE)$ such that 
\begin{equation} \label{eq: Weitzenboeck_formula_Dirac_type}
	D^{2} = -\tr_{\fg} (\connectionCoTanME \circ \connectionE) + V,   
\end{equation}
by the Weitzenb\"{o}ck formula (Remark~\ref{rem: Weitzenboeck_connection}), where $\tr_{\fg} := \eqref{eq: def_metric_trace}$ is the metric trace. 
\newline 

The most general Dirac-type operator on the Clifford module bundle $\big( \sE \to \sM, (\cdot|\cdot), \symb{D} \big)$ then differs from the composition of the following two maps 
\begin{equation}
	\secsME \stackrel{\connectionE}{\to} C^{\infty} (\sM; \coTansM \otimes \sE) \stackrel{\fc}{\to} \secsME 
\end{equation}
only by a zero-order term, i.e., 
\begin{equation} \label{eq: Dirac_op}
	D = - \ri \, \fc \circ \connectionE + U,    
\end{equation}
where the potential term $U \in C^{\infty} (\sM; \End \sE)$ is defined by the preceding equation. 
Note, Assumption~\ref{asp: Direc_op_symmetric} entails that $U$ and $- \ri \connectionE_{X}$ are symmetric with respect to $(\cdot|\cdot)$ provided that $X$ is divergence free. 
\newline 

The Weitzenb\"{o}ck connection $\connectionE$ induces a connection $\connectionEndE$ on the vector bundle $\Hom{\sE, \sE} \to \sM$ by 
\begin{equation}
	\connectionEndE_{X} \big( \symb{D} \xxi \big) (u) := \connectionE_{X} \big( \symb{D} \xxi \, (u) \big) - \symb{D} \xxi \big( \connectionE_{X} (u) \big) 
\end{equation}
for all $X \in \tangent_{x} \sM$ and all $\xi \in \coTan_{x} \sM$ at any $x \in \sM$. 
In general, $\connectionEndE$ is \textit{not} a Clifford connection because using the Levi-Civita connection we have 
\begin{eqnarray}
    X \big( \fg_{x}^{-1} (\xi, \eta) \big) 
	& = &  
    \fg_{x}^{-1} \big( \connectionLC_{X} \xi, \eta \big) + \fg_{x}^{-1} \big( \xi, \connectionLC_{X} \eta \big) 
	\nonumber \\ 
	& = &  
	\frac{1}{2} \big[ \symb{D} \big( x, \connectionLC_{X} \xi \big), \symb{D} (x, \eta) \big]_{+} + \frac{1}{2} \big[ \symb{D} \xxi, \symb{D} \big( x, \connectionLC_{X} \eta \big) \big]_{+},   
\end{eqnarray}
and then by taking the covariant derivative of~\eqref{eq: anticommutation_symbol_Dirac_type} we obtain 
\begin{equation}
	\big[ \connectionEndE_{X} \big( \symb{D} \xxi \big) - \symb{D} \big( x, \connectionLC_{X} \xi \big), \symb{D} (x, \eta) \big]_{+} 
	= 
	\big[ \connectionEndE_{X} \big( \symb{D} (x, \eta) \big) - \symb{D} \big( x, \connectionLC_{X} \eta \big), \symb{D} \xxi \big]_{+}  
\end{equation}
for all $\eta \in \coTan_{x} \sM$ and others as before. 
In particular, we can always locally choose an orthonormal frame $\{ e_{i} \}$ and its coframe $\{ \varepsilon^{i} \}$ on $\sM$ and express
\begin{equation}
	D = - \ri \slashed{\nabla}^{\ms \sE} + U, 
	\quad 
	\slashed{\nabla}^{\ms \sE} := \symb{D} (\varepsilon^{i}) \connectionE_{e_{i}}, 
\end{equation}
where the Feynman-slash notation $\slashed{\nabla}^{\ms \sE}$ has been used, i.e., $\connectionE$ is composed with Clifford multiplication and then traced over.
A straightforward computation yields  
\begin{equation}
	D^{2} 
	= 
    - \slashed{\nabla}^{\ms \sE} \slashed{\nabla}^{\ms \sE} - \ri \slashed{\nabla}^{\ms \sE} (U) + U^{2} - \Big( \slashed{\nabla}^{\ms \Hom{\sE, \sE}} \big( \symb{D} (\varepsilon^{j}) \big) + \ri U \symb{D} (\varepsilon^{j}) + \ri \symb{D} (\varepsilon^{j}) U \Big) \nabla_{e_{j}}^{\ms \sE}. 
\end{equation}
But 
\begin{eqnarray}
	\slashed{\nabla}^{\ms \sE} \slashed{\nabla}^{\ms \sE} 
	& = & 
	\frac{1}{2} (\slashed{\nabla}^{\ms \sE} \slashed{\nabla}^{\ms \sE} + \slashed{\nabla}^{\ms \sE} \slashed{\nabla}^{\ms \sE}) 
	\nonumber \\ 
	& = & 
	\frac{1}{2} \big( \symb{D} (\varepsilon^{i}) \, \symb{D} (\varepsilon^{j}) (\connectionE_{e_{i}} \connectionE_{e_{j}} - \Gamma_{ij}^{k} \connectionE_{e_{k}}) + \symb{D} (\varepsilon^{j}) \, \symb{D} (\varepsilon^{i}) (\connectionE_{e_{j}} \connectionE_{e_{i}} - \Gamma_{ji}^{k} \connectionE_{e_{k}}) \big) + \slashed{\Gamma}^{k} \connectionE_{e_{k}}
	\nonumber \\ 
	& = & 
	\frac{1}{2} \big( \symb{D} (\varepsilon^{i}) \, \symb{D} (\varepsilon^{j}) \Hess_{e_{i}, e_{j}}^{\ms \sE} + \, \symb{D} (\varepsilon^{j}) \, \symb{D} (\varepsilon^{i}) \Hess_{e_{j}, e_{i}}^{\ms \sE} \big) + \slashed{\Gamma}^{k} \connectionE_{e_{k}}
	\nonumber \\ 
	& = & 
	\frac{1}{2} \big( \symb{D} (\varepsilon^{i}) \, \symb{D} (\varepsilon^{j}) (\Hess_{e_{j}, e_{i}}^{\ms \sE} + \, \sR_{e_{i} e_{j}}) + \symb{D} (\varepsilon^{j}) \, \symb{D} (\varepsilon^{i}) \Hess_{e_{j}, e_{i}}^{\ms \sE} \big) + \slashed{\Gamma}^{k} \connectionE_{e_{k}}
	\nonumber \\ 
	& = & 
	\fg^{ij} \Hess_{e_{j}, e_{i}}^{\ms \sE} + \frac{1}{2} \slashed{\sR} + \slashed{\Gamma}^{k} \connectionE_{e_{k}}, 
\end{eqnarray}
where $\sR$ resp. $\Gamma$ are the curvature resp. connection $1$-form of $\connectionE$. 
Here, we have used~\eqref{eq: Bochner_d_Alembertian_coordinate}, ~\eqref{eq: def_Hessian},~\eqref{eq: anticommutation_symbol_Dirac_type} and $\sR_{X, Y} = \Hess_{X, Y} - \Hess_{Y, X}$. 
Inserting the preceding equation into the last equation for $D^{2}$ yields 
\begin{eqnarray}
	D^{2} 
	& = &  
	- \tr_{\fg} (\connectionCoTanME \circ \connectionE) + \kR + U^{2} - \ri \slashed{\nabla}^{\ms \sE} (U)  
	\nonumber \\ 
	&& 
	- \Big( \slashed{\nabla}^{\ms \Hom{\sE, \sE}} \big( \symb{D} (\varepsilon^{j}) \big) + \slashed{\Gamma}^{j} + \ri U \symb{D} (\varepsilon^{j}) + \ri \symb{D} (\varepsilon^{j}) U \Big) \nabla_{e_{j}}^{\ms \sE}, 
\end{eqnarray}
where 
\begin{equation}
	\kR : \secsME \to \secsME, ~ u \mapsto 
	\kR u := - \frac{1}{2} \big( \symb{D} (X^{\flat}) \, \symb{D} (Y^{\flat}) \, \sR_{X, Y} \big) u 
\end{equation}
is the \textbf{Weitzenb\"{o}ck curvature} of $\connectionE$. 
Since $\connectionE$ is a Weitzenb\"{o}ck connection, comparing with~\eqref{eq: Weitzenboeck_formula_Dirac_type} we equate the coefficients of $\connectionE_{e_{j}}$ in the foregoing equation of $D^{2}$ to zero to obtain the covariant derivative of $\symb{D}$:   
\begin{equation} \label{eq: covariant_derivative_symbol_Dirac} 
    \slashed{\nabla}^{\ms \Hom{\sE, \sE}} \big( \symb{D} (\varepsilon^{j}) \big) = - \slashed{\Gamma} - \big[ U, \symb{D} (\varepsilon^{j}) \big]_{+}, 
	\quad 
	j = 1, \ldots, d. 
\end{equation}

%
%
%
\begin{remark} \label{rem: Clifford_connection}
	A Clifford connection $\tilde{\nabla}^{\ms \Hom{\sE, \sE}}$, characterised by 
	\begin{equation}
		\forall (x, X) \in C^{\infty} (\sM; \tangent \sM), \forall (x, \eta) \in C^{\infty} (\sM; \coTansM) : 
		\big[ \tilde{\nabla}_{X}^{\ms \Hom{\sE, \sE}}, \symb{D} (x, \eta) \big]_{-} = \symb{D} \big(x, \nabla_{X}^{\ms \mathrm{LC}} \eta \big), 
	\end{equation}
    always exists and then the Clifford module bundle $(\sE \to \sM, \symb{D}, \tilde{\nabla}^{\ms \sE})$ is called the \textbf{bundle of Dirac modules}, where the most general Dirac-type operator has the form 
	\begin{equation}
        \tilde{D} := - \ri \, \fc \circ \tilde{\nabla}^{\ms \sE} + \tilde{U}
	\end{equation}
	for all $\tilde{U} \in C^{\infty} \big( \sM; \End{\sE} \big)$ such that (cf.~\eqref{eq: covariant_derivative_symbol_Dirac}) 
	\begin{equation} \label{eq: potential_compatible_Dirac}
		[\tilde{U}, \symb{D} (\cdot)]_{-} = 0. 
	\end{equation}
	The operator $\tilde{D}$ has the same principal symbol as $D$ so that they differ only by a smooth term. 
    Since $\tilde{\nabla}$ is compatible with the Clifford multiplication~\eqref{eq: def_Clifford_multiplication}, $\tilde{D}$ is sometimes called \textbf{compatible Dirac-type operator}. 
\end{remark}
%
%
%

\begin{example} \label{ex: spin_Dirac_op}
	Let $\spacetime$ be a $d$-dimensional spin-spacetime, i.e., a spacetime $\spacetime$ whose  tangent bundle admits a spin-structure. 
	This means, there exists a pair $(\sP, \varTheta)$ where $\PFBSpin$ is a principal $\SpinO$-bundle and $\varTheta : \sP \to \sQ$ is a twofold covering of the bundle $\PFBSO$ of positively oriented and time-oriented tangent frames, such that the following diagram commutes 
    (see e.g.~\cite{Baum_1981, Baer_MathZ_2005})   
    \begin{center} 
        \begin{tikzpicture}  
            \node (a) at (0,0) {$\sP \times \rSpin_{0}$};  
            \node (b) at (4, 0) {$\sP$};
            \node (c) at (0, -2) {$\sQ \times \rSO_{0}$};
            \node (d) at (4, -2) {$\sQ$};
            \node (e) at (6.5, -1) {$\sM$};
            \node at (0.8, -1) {$\varTheta \times \Ad$};
            \node at (4.3, -1) {$\varTheta$}; 
            \draw[->] (a) -- (b);
            \draw[->] (c) -- (d);
            \draw[->] (a) -- (c);
            \draw[->] (b) -- (d);
            \draw[->] (b) -- (e); 
            \draw[->] (d) -- (e);
        \end{tikzpicture}
        \captionof{figure}[Lorentzian spin-structure]{Lorentzian 
            spin-structure on a spacetime. 
            The horizontal lines are the group actions on the principal bundles.}
    \end{center}
	Here, $\SpinO$ resp. $\SOO$ are the connected components to the identity of the Lorentzian spin group $\rSpin (1, d-1)$ resp. the special orthogonal group $\rSO (1, d-1)$ and $\Ad : \SpinO \to \SOO$ is the twofold covering map.  
	\newline 

	A spin structure always exists locally but its \textit{global} existence depends on the topology of $\sM$.  
    A necessary and sufficient condition is that the \textit{second Stieffel-Whitney class of $\sM$ vanishes}~\cite{Haefliger_1956, Borel_AJM_1959} (Riemannian signature); 
    see~\cite{Geroch_JMP_1968, Geroch_spin_JMP_1970, Clarke_GRG_1971,  Baum_1981, Alagia_RUMA_1985} for the Lorentzian setting. 
    Under this topological restriction, the number of inequivalent spin structures on $\sM$ is equal to the number of elements in the cohomology $\rH^{1} (\sM; \Z_{2})$~\cite[Lem.]{Milnor_Enseign_1963} 
    (Riemannian signature; see~\cite{Isham_spin_PRSL_1978} for $4$-dimensional spacetimes). 
    Then, for a fixed spin-structure $(\sP, \varTheta)$, the \textbf{spinor bundle} is defined as the associated vector bundle 
    \begin{equation}
    	\sS := \sP \times_{\tilde{\gammaup}} \cS_{d},  
    \end{equation}
    where $\cS_{d}$ is the complex spinor module. 
	When $d$ is \textit{even}, $\dim_{\C} \cS_{d} = 2^{d/2}$, the spinor module decomposes $\cS_{d} = \cS_{d}^{+} \oplus \cS_{d}^{-}$ into the submodules $\cS_{d}^{\pm}$ of positive(negative) chirality and 
	\begin{equation}
		\tilde{\gammaup} = \tilde{\gammaup}^{+} \oplus \tilde{\gammaup}^{-} : \SpinO \to \Aut \cS_{d}^{+} \times \Aut \cS_{d}^{-} \subset \Aut \cS_{d}
	\end{equation}
	is the \textit{spinor representation}. 
	For an \textit{odd} $d$, $\dim_{\C} \cS_{d} = 2^{(d-1)/2}$ and 
	\begin{equation}
		\tilde{\gammaup} : \SpinO \to \Aut \cS_{d} 
	\end{equation}
	is the \textit{spinor representation}. 
    Sections of a spinor bundle are called the \textbf{spinors}. 
    In other words, for each $x \in \sM$, the fibre $\sS_{x}$ of the spinor bundle $\sS$ over $x$ consists of equivalence classes of pairs $\equiBra{p, s}$ subject to the free left $\SpinO$-group action on $\sP \times \cS_{d}$:   
    \begin{equation}
    	\forall p \in \sP, \forall u \in \cS_{d}, \forall g \in \SpinO : 
		g \bullet \equiBra{p, u} = \equiBra{p g^{-1}, \tilde{\gammaup} (g) \, u}. 
    \end{equation}
    Since $\tansM$ sets inside $\SpinO$, $\tilde{\gammaup}$ induces an $\End \cS_{d}$-valued map $\gammaup$ on $\tansM$: 
    \begin{equation}
        \gammaup (x, X)^{2} := \fg_{x} (X, X) \, \one_{\End \cS_{d}}. 
    \end{equation}
    This enables us to define the Clifford multiplication 
    \begin{equation}
        \fc : C^{\infty} (\sM; \tangent \sM \otimes \sS) \to C^{\infty} (\sM; \sS), ~  
    	\fc \big( \equiBra{\varTheta (p), X} \otimes \equiBra{p, s} \big) := \equiBra{p, \gammaup (X) \, s}, 
    \end{equation}
    where we have expressed $\tansM = \sQ \times_{\rhoup} \Rd$ as a real vector bundle associated to $\sQ$ by the standard representation $\rhoup$ of $\SOO$ on $\Rd$. 
    \newline 

    Employing the Levi-Civita connection $1$-from $\Gamma^{\ms \mathrm{LC}}$ on $\sQ$, we define the spinor connection $1$-form $\Gamma^{\ms \sS} := \Ad_{*}^{-1} \circ \varTheta^{*} (\Gamma^{\ms \mathrm{LC}})$ which induces a connection $\nabla^{\ms \sS}$ on $\sS$. 
	This spin connection is a \textit{metric} and \textit{Clifford}  connection. 
    Moreover, it leaves the positive-negative chiral splitting invariant and satisfies the Leibniz rule. 
	The \textbf{massive spin-Dirac operator} is defined by~\cite[Def. 3.1]{Baum_1981}  
	\begin{equation}
		D := - \ri \fc \circ \nabla^{\ms \sS} + \rrm : C^{\infty} (\sM; \sS) \to C^{\infty} (\sM; \sS), 
	\end{equation} 
	where $\rrm$ is a parameter which can be interpreted as the physical mass of a spinor in appropriate spacetime geometry. 
\end{example}
%
%
%

\begin{remark} \label{rem: spin_structure_4_d_spacetime}
	If $\spacetime$ is a $4$-dimensional non-compact spacetime then the sufficient and necessary condition to admit a spin-structure entails the existence of a global system of orthonormal tetrads~\cite[Thm.]{Geroch_JMP_1968}. 
	In particular, a number of exact solutions of Einstein equation, e.g., the Schwartzschild, the Robertson-Walker, the G\"{o}del, the plane waves, and any \textit{globally hyperbolic spacetime} admit at least a spin-structure~\cite{Geroch_spin_JMP_1970}. 
    Furthermore, the spinor bundle is trivial yet there may exists inequivalent spin-structures $(\sP, \varTheta)$ (given by number of elements in $\rH^{1} (\sM; \Z_{2})$) encoded in $\varTheta$~\cite{Isham_spin_PRSL_1978}. 
\end{remark}
%
%
%
%
%
%
%
%
%
%
\subsection{Green's operators}
\label{sec: Green_op_Dirac}
Since $D^{2}$ is normally hyperbolic, it admits unique advanced $\GreenOpAdv$ and retarded $\GreenOpRet$ Green's operators on any globally hyperbolic spacetime $\spacetime$~\cite[Thm. 1]{Muehlhoff_JMP_2011}. 
Hence  
\begin{equation}
	F^{\adv, \ret} := D G^{\adv, \ret} : \comSecsME \to C_{\mathrm{sc}}^{\infty} (\sM; \sE)   
\end{equation}
are the unique advanced resp. retarded Green's operators for $D$. 
As before, 
\begin{equation} \label{eq: def_causal_propagator_Dirac}
	F := \fundaSolRet - \fundaSolAdv: \comSecsME \to C_{\mathrm{sc}}^{\infty} (\sM; \sE)  
\end{equation}
defines the Pauli-Jordan-Lichnerowicz Green's operator, also known as the \textit{causal propagator} for $D$.  
Assumption~\ref{asp: Direc_op_symmetric} and Remark~\ref{rem: exist_unique_advanced_retarded_Green_op_NHOp} entail that 
\begin{equation} \label{eq: adv_ret_Green_op_Dirac_type_selfadjoint_sesquilinear}
	\forall u, v \in \comSecsME : 
	(F^{\adv, \ret} u | v) = (u | F^{\ret, \adv} v), 
	\quad 
	(F u|v) = - (u|Fv). 
\end{equation}

%
%
%
\begin{lemma} \label{lem: causal_propagator_Dirac}
	Let $\sE \to \sM$ resp. $\pi : \dotCoTansM \to \sM$ be a vector bundle resp. the punctured cotangent bundle over a globally hyperbolic spacetime $\spacetime$ and $D$ a Dirac-type operator on $\sE$ whose principal symbol is denoted by $\symb{D}$. 
	The Schwartz kernel $\fF$ of the causal propagator for $D$ is then a Lagrangian distribution:   
	\begin{subequations} \label{eq: causal_propagator_Dirac}
		\begin{eqnarray}
			&& 
			\fF \in I^{-1/2} \big( \sM \times \sM, C'; \Hom{\sE, \sE} \big), 
			\label{eq: causal_propagator_Dirac_FIO}
			\\  
			&& 
            \symb{\fF} = \ri \sqrt{\frac{\pi}{2}} \symb{D} \circ w \sqrt{| \dVol_{\ms C} |} \otimes \mathbbm{l},   
			\label{eq: symbol_causal_propagator_Dirac}
		\end{eqnarray}
	\end{subequations}
    where $\dVol_{\ms C}$ is the natural volume form (Remark~\ref{rem: volume_geodesic_relation}) on the geodesic relation (Definition~\ref{def: geodesic_relation}) $C$, $\bbl$ is a section of the Keller-Maslov bundle $\bbL_{C} \to C$ 
    (as constructed in~\cite[pp. 231-232]{Duistermaat_ActaMath_1972}), and $w$ is the unique element of $C^{\infty} \big( C; \pi^{*} \Hom{\sE, \sE}_{C} \big)$ that is diagonally the identity endomorphism and off-diagonally covariantly constant  
	\begin{equation} \label{eq: symbol_causal_propagator_Dirac_covarinatly_constant}
        \connectionEndPiE_{X_{\fg/2}} w = 0   
	\end{equation}
	with respect to the $D^{2}$-compatible Weitzenb\"{o}ck covariant derivative (Definition~\ref{def: P_compatible_connection}) $\connectionEndPiE_{X_{\fg/2}}$ along the geodesic vector field $X_{\fg/2}$ (Remark~\ref{rem: HVF_canonical_relation}). 
	\newline 

	If $(\sE \to \sM, \symb{D}, \tilde{\nabla}^{\ms \sE})$ is a bundle of Dirac modules with the corresponding Dirac-type operator $\tilde{D}$ then the preceding holds with the replacement of $\connectionEndPiE$ in~~\eqref{eq: symbol_causal_propagator_Dirac_covarinatly_constant} by the $\tilde{D}^{2}$-compatible connection $\tilde{\nabla}^{\ms \pi^{*} \Hom{\sE, \sE}}$ induced by the Clifford connection $\tilde{\nabla}^{\ms \sE}$.    
\end{lemma}
%
%
%
\begin{proof}
    This follows from Theorem~\ref{thm: causal_propagator_NHOp_FIO} and the facts that $D$ is a first-order differential operator, $\Char D = \coLightBun$, and $\fF (x, \cdot) = D_{x} \fG (x, \cdot)$. 
\end{proof}
%
%
%
%
%
%
%
%
%
%
\subsection{Cauchy problem}
\label{sec: Cauchy_problem_Dirac}
Since $D^{2}$ is a normally hyperbolic operator and $\spacetime$ is a globally hyperbolic spacetime, the Cauchy problem for $D$ is well-posed~\cite[Thm. 2]{Muehlhoff_JMP_2011}. 
In other words, for an arbitrary but fixed Cauchy hypersurface $\varSigma_{t} \subset \sM$, the mapping  
\begin{equation} \label{eq: sol_Dirac_Cauchy_data}
	\cR_{t} : \ker D \to C_{\mathrm{c}}^{\infty} (\varSigma_{t}; \sE_{\varSigma_{t}}), ~ u \mapsto \cR_{t} (u) := u \upharpoonright \varSigma_{t}, 
	\quad \supp{u} \subset J \big( \supp (u \! \upharpoonright \! \varSigma_{t}) \big)  
\end{equation}
is a homeomorphism, where $J$ is defined by~\eqref{eq: causal_future_past}. 
In order to endow $\ker D$ with a hermitian inner product,  we make use of the following 

%
%
%
\begin{assumption} \label{asp: hermitian_form_Dirac_type_op}
	Given an arbitrary but fixed future-directed unit normal covector field $(\cdot, \zeta)$ on a globally hyperbolic spacetime $\spacetime$ along any Cauchy hypersurface $\varSigma$, 
	\begin{equation} \label{eq: hermitian_form_Dirac_type_op}
		\langle \cdot | \cdot \rangle := \, \big( \symb{D} (\cdot, \zeta ) \cdot \big| \cdot \big)   
	\end{equation}
	is a fibrewise hermitian form on the bundle of Clifford modules $\big( \sE \to \sM, \symb{D}, (\cdot|\cdot) \big)$. 
\end{assumption}
%
%
%

The hermitian form $\langle \cdot | \cdot \rangle$ \textit{depends} on the choice of $(\cdot, \zeta)$. 
By global hyperbolicity (Theorem~\ref{thm: globally_hyperbolic_Cauchy}) of $\spacetime$, there exists a global Cauchy temporal function $\ft$ such that each Cauchy hypersurface $\varSigma_{t} := \ft^{-1} (t)$ is a level set of $\ft$, where $t \in \R$. 
Then 
\begin{equation} \label{eq: def_unit_normal_covector_GHST}
	\zeta := \frac{\rd \ft}{\| \rd \ft \|}
\end{equation}
is a unit normal covector field on $\sM$ along $\varSigma_{t}$ and we choose the time-orientation such that $\zeta$ is future-directed. 
\newline 

Employing the preceding topological isomorphism $\cR_{t}$, $\ker D$ can be equipped with a hermitian inner product  
(see e.g.~\cite[Lem. 3.17]{Baer_Springer_2012}):   
\begin{equation} \label{eq: def_hermitian_form_Dirac_type_op}
	\langle u | v \rangle := \int_{\varSigma} \Big( \fc \big( (x, \zeta) \otimes (u \! \upharpoonright \! \varSigma) \big) \big| v \! \upharpoonright \! \varSigma \Big)_{x} \dVolh (x), 
\end{equation}
where $\dVolh$ is the Riemannian volume element on $\varSigma$. 
We remark that $\scalarProdTwo{\cdot}{\cdot}$ is independent of the chosen Cauchy hypersurface $\varSigma$ due to the Green-Stokes formula (see~\eqref{eq: Stokes_formula_Dirac}).  
Therefore, $(\sE \to \sM, \langle \cdot|\cdot \rangle)$ is a \textit{hermitian vector bundle}, of course \textit{depending} on $(\cdot, \zeta)$, but \textit{independent} of chosen $\varSigma$.   
\newline 

We infer from the preceding equation that, given an initial data $k \in C_{\mathrm{c}}^{\infty} (\varSigma; \sE_{\varSigma})$, any smooth solution $u$ of the Dirac equation~\eqref{eq: Dirac_eq} can be expressed as 
(see~\cite[Prop. 2.4 (b)]{Dimock_AMS_1982} for a spin Dirac operator)
\begin{equation} \label{eq: smooth_sol_Dirac_op_Cauchy_data}
    u = - \ri \, F \circ (\iota_{\ms \varSigma}^{*})^{-1} \, \fc (\zeta  \otimes k), 
\end{equation} 
where the restriction operator 
\begin{equation} \label{eq: def_restriction_op_GHST}
    \iota_{\ms \varSigma}^{*} : \comSecsME \to C_{\mathrm{c}}^{\infty} (\varSigma; \sE_{\varSigma})
\end{equation}
has been discussed elaborately in Example~\ref{exm: restriction_op_FIO} and we have taken into account Remark~\ref{rem: Cauchy_hypersurface_embedded_submf}.  
All the maps in the exact complex  
\begin{equation} \label{eq: exact_sequence_causal_propagator_Dirac}
	0 
	\to \comSecsME
	\xrightarrow{~ D ~} \comSecsME 
	\xrightarrow{~ F ~} C_{\mathrm{sc}}^{\infty} (\sM; \sE)  
	\xrightarrow{~ D~} C_{\mathrm{sc}}^{\infty} (\sM; \sE) 
\end{equation}
are sequentially continuous as a consequence of $D$ being a local operator and the following exact complex being sequentially continuous  
(see e.g.~\cite[Prop. 3.4.8]{Baer_EMS_2007})  
\begin{equation}
	0 
	\to \comSecsME
	\xrightarrow{~ D^{2} ~} \comSecsME 
	\xrightarrow{~ G~} C_{\mathrm{sc}}^{\infty} (\sM; \sE)  
	\xrightarrow{~ D^{2} ~} C_{\mathrm{sc}}^{\infty} (\sM; \sE).    
\end{equation}
%
%
%
%
%
%
%
%
%
%
\subsection{Primary result}
\label{sec: result_Feynman_propagator_Dirac}
Since $D^{2}$ is a normally hyperbolic operator and the first hypothesis in Theorem~\ref{thm: existence_Feynman_propagator_NHOp} is satisfied by Assumption~\ref{asp: Direc_op_symmetric}, one obtains Feynman propagator $\fundaSolFeyn$ for $D$ immediately by setting 
\begin{equation}
	 \fundaSolFeyn := D \GreenOpFeyn, 
\end{equation}
where $\GreenOpFeyn$ is the Feynman propagator for $D^{2}$. 
Since $\Char D = \coLightBun$, one has   
\begin{equation} \label{eq: WF_advanced_retarded_Feynman_antiFeynman_Green_op_Dirac_type}
	\WFPrime (\fF^{\adv, \ret, \Feyn, \aFeyn}) 
    = 
    \WFPrime (\fG^{\adv, \ret, \Feyn, \aFeyn}) 
    = \varDelta \dotCoTansM \bigcup C^{\adv, \ret, +, -}, 
\end{equation}
where the canonical relations are defined in~\eqref{eq: def_advanced_retarded_canonical_relation_NHOp} and in Definition~\ref{def: geodesic_relation}. 
Nevertheless, the existence of Hadamard states cannot be deduced due to the positivity issue as $(\cdot|\cdot)$ is \textit{not} definite. 
One may wish to consider symmetricness of $D$ with respect to some hermitian form $\langle \cdot | \cdot \rangle$ so that the hypotheses of Theorem~\ref{thm: existence_Feynman_propagator_NHOp} are satisfied in such a way that the existence of Hadamard states can be concluded. 
The natural inner product on spinors on a spacetime is not positive-definite rather that is \textit{indefinite} and therefore the positivity property of Feynman propagators for a spin-Dirac operator (Example~\ref{ex: spin_Dirac_op}) cannot be inferred directly deploying Theorem~\ref{thm: existence_Feynman_propagator_NHOp}. 
In other words, the requirement of symmetricness with respect to $\langle \cdot | \cdot \rangle$ turns out to be too restrictive to encompass all Dirac-type operators on globally hyperbolic manifolds (see Sections~\ref{sec: twisted_Dirac_op} and~\ref{sec: Rarita_Schwinger_op} for concrete examples) and so we refrain to impose this condition. 

%
%
%
\begin{restatable}{theorem}{HadamardBisolutionDiracTypeOp}
	\label{thm: Hadamard_bisolution_Dirac_type_op}	
	Let $\sE \to \sM$ be a vector bundle over a globally hyperbolic spacetime $\spacetime$, endowed with a non-degenerate sesquilinear form $(\cdot|\cdot)$, and $D$ a Dirac-type operator on $\sE$ whose principal symbol is given by $\symb{D}$. 
	Suppose that $\big( \sE \to \sM, \symb{D}, (\cdot|\cdot) \big)$ is a bundle of Clifford modules over $\spacetime$ and that $D$ is symmetric with respect to $(\cdot|\cdot)$. 
	Then, there exists a Feynman propagator $\fundaSolFeyn$ for $D$ such that $W := - \ri (\fundaSolFeyn - \fundaSolAdv)$ is symmetric, where $\fundaSolAdv$ is the advanced propagator for $D$. 
	In addition, if there exists a hermitian form $\langle \cdot | \cdot \rangle$ on $\sE$ in the sense of Assumption~\ref{asp: hermitian_form_Dirac_type_op}, then $\fundaSolFeyn$ can be chosen such that $W$ is non-negative with respect to $(\cdot|\cdot)$ and hence defines a Hadamard bisolution of $D$.
\end{restatable}
%
%
%

\begin{proof}
	To begin with, recall that $\Char D^{2} = \coLightBun$ has only two connected components in three or higher spacetime   
    dimensions\footnote{In 
        two dimensions, $\coLightBun$ has four connected components; cf. Remark~\ref{rem: sixteen_parametrix_dim_two_NHOp}.}:
    the forward lightcone $\dotCoTan_{0,+} \sM$ and the backward lightcone $\dotCoTan_{0,-} \sM$.
    Therefore, $C$ has four different 
    orientations\footnote{By 
        an orientation of the geodesic relation $C$, it is meant that any splitting of $C \setminus \varDelta \, \coLightBun$ into a disjoint union of open $C^{1}, C^{2} \subset C$ which are inverse relations~\cite[p. 218]{Duistermaat_ActaMath_1972}.
	}
    \cite[p. 218]{Duistermaat_ActaMath_1972} 
    (see also~\cite[p. 540]{Radzikowski_CMP_1996}):  
    $(C^{+}, C^{-}), (C^{-}, C^{+}), (C^{\adv}, C^{\ret}), (C^{\ret}, C^{\adv})$ corresponding to $\{ \coLightBun, \dotCoTan_{0,+} \sM, \dotCoTan_{0,-} \sM, \emptyset \}$.  
    This paves a way to have a (non-unique) microlocal decomposition of the causal propagator $G$ for $D^{2}$ as
    \begin{subequations}
        \begin{eqnarray}
            && 
            \fG \equiv \fE^{+} - \fE^{-},  
            \\ 
            && 
			\fE^{\pm} D^{2} \equiv D^{2} \fE^{\pm} \equiv 0,    
            \\ 
            && 
            \fE^{\pm} \in I^{-3/2} \big( \sM \times \sM, C^{\pm \prime}; \Hom{\sE, \sE} \big),  
        \end{eqnarray}
    \end{subequations}
	where $\equiv$ means modulo smoothing kernels. 
    Clearly, this is a slightly stronger version of Theorem~\ref{thm: exist_unique_Feynman_parametrix_NHOp}, which can be proven exactly in the same fashion; cf Theorem~\ref{thm: causal_propagator_NHOp_FIO}. 
    \newline 

	To achieve the analog of the proceeding results for $D$, we remember $\Char D = \Char D^{2}$ and decompose $\fF$ as $\fS^{+} + \fS^{-}$ where
    \begin{equation}
    	\fS^{\pm} := \sum_{\alpha} Q_{\alpha}^{\pm} \fF Q_{\alpha}^{\pm}
    \end{equation}
    and $(Q_{\alpha}^{\pm})^{*} Q_{\alpha}^{\pm}$ are microlocal partition of unity such that $D Q_{\alpha}^{\pm} \equiv Q_{\alpha}^{\pm} D$ and $\ES{Q_{\alpha}^{\pm}} = \varDelta \, \dotCoTan_{0, \pm} \sM$. 
    Consequently,  
    \begin{subequations}
        \begin{eqnarray}
            && 
            \fF \equiv \fS^{+} + \fS^{-},    
            \label{eq: microlocal_decomposition_causal_propagator_Dirac_type_op}
            \\ 
            && 
            \fS^{\pm} D \equiv D \fS^{\pm} \equiv 0,  
            \label{eq: distinghuished_parametrices_Dirac_type_op} 
            \\ 
            && 
            \fS^{\pm} \in I^{-1/2} \big( \sM \times \sM, C^{\pm \prime}; \Hom{\sE, \sE} \big), 
        \end{eqnarray}
    \end{subequations}
	by an application of~\eqref{eq: product_Lagrangian_dist}. 
    The operators $Q_{\alpha}^{\pm}$ can be constructed by choosing any $\fq_{\alpha}^{\pm} \in S^{0 - [\infty]} \big( \dotCoTansM, \Hom{\sE, \sE} \big)$ as its left total symbol such that $\esssupp \fq_{\alpha}^{\pm} = \dotCoTan_{0, \pm} \sM$. 
    In other words $\fq_{\alpha}^{\pm}$ is of order $-\infty$ near $y \in \sM$ and in a conic neighbourhood of $\eta$ for all $\yeta$ in the complement of $\esssupp \fq^{\pm}$. 
    Furthermore, we pick $\supp \fq_{\alpha}^{\pm}$ slightly away from the projection of $\ES{Q_{\alpha}^{\pm}}$ on $\sM$ so that $[D, Q_{\alpha}^{\pm}]_{-}$ is smooth. 
	
    %
    %
    %
    \begin{lemma}
        As in the terminologies of Theorem~\ref{thm: Hadamard_bisolution_Dirac_type_op}, let $F$ be the causal propagator for $D$. 
        Then $\ri F$ is non-negative with respect to $(\cdot|\cdot)$. 
    \end{lemma}
    %
    %
    %
	
    \begin{proof} 
        Let $u, v \in \secsME$ having compact $\supp{u} \cap \supp{v}$ and let $K \subset \sM$ be a compact set having smooth boundary $\partial K$ with an outward unit normal covector field $\zeta$ and volume element $\dVol$. 
        Then, the Green-Stokes formula yields 
        (see e.g.~\cite[Prop. 9.1, $($p. 178$)$]{Taylor_I_Springer_2011}),~\cite[$(1.7)$]{Dimock_AMS_1982}) 
        \begin{equation} \label{eq: Stokes_formula_Dirac}
            \int_{x \in K} \big( (Du| v)_{x} - (u | Dv)_{x} \big) \dVolg (x) 
            = 
            - \ri \int_{x' \in \partial K} \big( \symb{D} (\cdot, \zeta ) \, u \big| v \big)_{x'} \dVol (x'),  
        \end{equation}
		where $\dVolg$ is the Lorentzian volume (Remark~\ref{exm: density}~\ref{exm: density_Lorentzian_mf}) form on $\sM$. 
        Setting $u := F^{\ret} w$ and $v := Fw$ for any $w \in \comSecsME$,  the above formula entails $(w | \ri F w) \geq 0$ due to the compactness of $\supp (F^{\ret} u) \subset J^{+} (\supp{u})$. 
    \end{proof}
    %
    %
    %
    Therefore, $(u | \ri S^{\pm} u) \geq 0$ by~\eqref{eq: microlocal_decomposition_causal_propagator_Dirac_type_op}. 
    Note, $C^{\pm} = C \cap (\dotCoTan_{0, +} \! \sM \times \dotCoTan_{0, +} \! \sM)$. 
    Inspecting the (twisted) wavefront sets of $\fundaSolFeyn$ and $\fundaSolAdv$ as before in the proof of Proposition~\ref{prop: positivity_Feynman_minus_adv_NHOp}, we have the analogue of Proposition~\ref{prop: positivity_Feynman_minus_adv_NHOp} for Dirac-type operators: 
    \begin{equation}
    	(u| \hat{W} u) \geq 0, 
    	\quad 
    	\hat{W} := \ri (\fundaSolAdv - \fundaSolFeyn). 
    \end{equation}

    As in the proof of Theorem~\ref{thm: existence_Feynman_propagator_NHOp}, in order to turn $\hat{\fW}$ into an exact distributional bisolution $\fW$, we are going to employ the well-posedness of the Cauchy problem for $D$ on a globally hyperbolic spacetime~\cite[Thm. 2]{Muehlhoff_JMP_2011}.  
    For any $k \in \comSecECauchy$, we define distributional Cauchy data $k \otimes \delta_{\varSigma}$ by 
    \begin{equation} \label{eq: def_distributional_Cauchy_data_Dirac_type_op}
        k \otimes \delta_{\varSigma} (\phi) := \int_{x \in \varSigma} \phi (x) \, \fc (\zeta \otimes k)_{x} \, \dVolh (x) 
    \end{equation}
    for any $\phi \in \secsMEStar$. 
	By $\cY_{\varSigma} \subset \cE_{\varSigma^{\perp *}}' (\sM; \sE)$, we denote the set of distributions of the above form where $\cE_{\varSigma^{\perp *}}' (\sM; \sE)$ is the set of compactly supported distributions whose wavefront set is contained in the conormal bundle $\varSigma^{\perp *}$. 
    Any bidistribution $\fw \in \cD' (\sM \times \sM; \sE \boxtimes \bar{\sE}^{*})$ such that $\fw D$ and $D \fw$ are smooth in distributional sense, can be defined as a sequential continuous bilinear form on $\cY_{\varSigma}$.  
	
    %
    %
    %
    \begin{lemma}
        As in the terminologies of Theorem~\ref{thm: Hadamard_bisolution_Dirac_type_op}, suppose that $\fw \in \cD' (\sM \times \sM; \sE \boxtimes \bar{\sE}^{*})$ is a bisolution of $D$. 
		Then $(u | w u) \geq 0$ if and only if $\fw (\bar{u}^{*} \otimes u) \geq 0$ for all $u \in \cY_{\varSigma}$, where $\cY_{\varSigma}$ is the span of the distributional Cauchy data~\eqref{eq: def_distributional_Cauchy_data_Dirac_type_op} of $D$. 
    \end{lemma}
    %
    %
    %
	
    \begin{proof}
		The proof carries over its analogue Lemma~\ref{lem: positive_bisolution_NHOp} after the obvious replacement of $D^{2} (= \square)$, its Cauchy data~\eqref{eq: def_distributional_Cauchy_data_NHOp}, and causal propagator $G$, by $D$, its Cauchy data~\eqref{eq: def_distributional_Cauchy_data_Dirac_type_op}, and causal propagator $F$, respectively. 
    \end{proof} 
\end{proof}
%
%
%
%
%
%
%
%
%
%
\section{Examples}
We list below several operators arising in the context of quantum field theories in a curved spacetime and discuss how they fit into our framework. 
%
%
%
%
%
%
%
%
%
%
\subsection{Covariant Klein-Gordon operator} 
As described in Example~\ref{exm: covariant_Klein_Gordon_op}, the natural inner product on $C^{\infty} (\sM)$ is positive-definite with respect to which $\square$ is symmetric. 
Therefore, Theorem~\ref{thm: existence_Feynman_propagator_NHOp} applies and positive Feynman propagators can be constructed using this method.   
%
%
%
%
%
%
%
%
%
%
\subsection{Connection d'Alembert operator} 
A Feynman propagate exists for $\square$ (Example~\ref{exm: connection_d_Alembertian}) by  Theorem~\ref{thm: existence_Feynman_propagator_NHOp}. 
The positivity issue depends on the positivity of the inner product with respect to which $\square$ is symmetric. 
%
%
%
%
%
%
%
%
%
%
\subsection{Hodge-d'Alembert operator} 
\label{sec: Hodge_d_Alembert_op}
A Feynman propagator exists for this operator defined in Example~\ref{exm: Hodge_d_Alembert_op} by Theorem~\ref{thm: existence_Feynman_propagator_NHOp}. 
Since the natural inner product with respect to which this operator is symmetric, is not positive-definite unless $k=0$ or $k=d$, we cannot conclude non-negativity in a straightforward way.  
For example, for $k=1$ non-negativity is only expected on a subset as is usual for gauge theories. 
%
%
%
%
%
%
%
%
%
%
\subsection{Proca operator} 
\label{sec: Proca_op}
As in the terminologies of Example~\ref{exm: Proca_op}, a Feynman propagator for $\rd^{*} \rd + \rd \rd^{*} + \mathrm{m}^{2}$ can be constructed by this method and the positivity issue is same as in Section~\ref{sec: Hodge_d_Alembert_op}. 
%
%
%
%
%
%
%
%
%
%
\subsection{Twisted spin-Dirac operator} 
\label{sec: twisted_Dirac_op}
Let $(\sS \to \sM, \nabla^{\ms \sS})$ and $(\sE \to \sM, \connectionE)$ be a spinor bundle (Example~\ref{ex: spin_Dirac_op} and Remark~\ref{rem: spin_structure_4_d_spacetime}) and a vector bundle endowed with the Levi-Civita spin connection $\nabla^{\ms \sS}$ and a bundle connection $\connectionE$, over a globally hyperbolic spin-spacetime $\spacetime$.    
These two connections induce another connection $\nabla := \nabla^{\ms \sS} \otimes \one_{\sE} + \one_{\sS} \otimes \nabla^{\ms \sE}$ on the twisted spinor bundle $\sS \otimes \sE \to \sM$ and the \textbf{twisted Dirac operator} is defined by 
(see e.g.~\cite{Baum_1981}) 
\begin{equation} \label{eq: def_twisted_Dirac_op}
	\slashed{D} := - \ri \upgamma \circ \fg^{-1} \circ \nabla : C^{\infty} (\sM; \sS \otimes \sE) \to C^{\infty} (\sM, \sS \otimes \sE), 
\end{equation} 
where $\upgamma : C^{\infty} (\sM; \tansM) \to \End{\sS}$ is the Clifford mapping such that 
\begin{equation}
    \upgamma (X) \, \upgamma (Y) + \upgamma (Y) \, \upgamma (X) = 2 \fg_{x} (X, Y) \one_{\End{\sS}_{x}} 
\end{equation}
for any $X, Y \in \tangent_{x} \sM$ and $x \in \sM$.  
The Clifford multiplication is given by 
\begin{equation}
    C^{\infty} (\sM; \tansM \otimes \sS \otimes \sE) \ni (x, X) \otimes \psi \otimes u \mapsto (\upgamma (X) \psi) \otimes u \in C^{\infty} (\sM; \sS \otimes \sE).
\end{equation}
The Schr\"{o}dinger~\cite{Schroedinger_GRG_2020}-Lichnwerowicz~\cite{Lichnerowicz_1963} 
formula entails  
(see e.g.~\cite[Exm. 1]{Muehlhoff_JMP_2011}) 
\begin{equation}
	\slashed{D}^{2} = \nabla^{*} \nabla + \frac{\ric}{4} + \kF, 
\end{equation}
where $\ric$ is the scalar curvature of $\sM$ and $\kF$ is the Clifford multiplied curvature of $\nabla^{\ms \sE}$. 
\newline 

Clearly, $\slashed{D}^{2}$ is a normally hyperbolic operator and thus $\slashed{D}$ is of  Dirac-type. 
This formula also shows that the \textit{Weitzenb\"ock connection} is the \textit{twisted spin-connection} $\nabla$ on the twisted spinor bundle $(\sS \otimes \sE \to \sM, \nabla)$ and therefore the induced connection on $\Hom{\sS \otimes \sE, \sS \otimes \sE} \to \dotCoTansM$ is to be used for microlocalisation (Theorem~\ref{thm: microlocalisation_P}). 
Since $\slashed{D}$ is symmetric with respect to the natural inner product $(\cdot | \cdot)$ on $\sS \otimes \sE$, a Feynman propagator for $\slashed{D}$ exists by Theorem~\ref{thm: existence_Feynman_propagator_NHOp} yet the existence of Hadamard states cannot be concluded due to the fact that $(\cdot | \cdot)$ is not positive-definite unless $\spacetime$ is a Riemannian manifold. 
However, $\spacetime$ is globally hyperbolic and thus (Theorem~\ref{thm: globally_hyperbolic_Cauchy} and~\eqref{eq: def_unit_normal_covector_GHST}), there exist a smooth spacelike Cauchy hypersurfaces $\varSigma$ and a future-directed unit vector field $\rN$ on $\sM$ normal to $\varSigma$ so that $\upgamma (\rN) \cdot$ defines the Clifford multiplication to achieve the hermitian form $\langle \cdot | \cdot \rangle$ (cf. Assumption~\ref{asp: hermitian_form_Dirac_type_op}). 
Then, Theorem~\ref{thm: Hadamard_bisolution_Dirac_type_op} is applicable and we have Feynman propagators with the desired positivity with respect to $(\cdot|\cdot)$.  
%
%
%
%
%
%
%
%
%
%
\subsection{Rarita-Schwinger operator} 
\label{sec: Rarita_Schwinger_op}
As in the terminologies of Example~\ref{ex: spin_Dirac_op}, let $\sS \to \sM$ be a spinor bundle over a spin-spacetime $\spacetime$ and the Clifford multiplication is given by 
$\coTansM \otimes \sS \ni \xxi \otimes u \mapsto \upgamma (\xi^{\sharp}) u \in \sS$. 
Then, we have the representation theoretic splitting    
\begin{equation}
	\coTansM \otimes \sS = \iota (\sS) \oplus \sS^{3/2}, 
\end{equation}
where $\sS^{3/2}$ is defined as the kernel of the Clifford multiplication and the embedding $\iota$ of $\sS$ into $\coTansM \otimes \sS$ is locally defined by 
\begin{equation}
	\iota (u) := - \frac{1}{d} e_{i} \otimes \upgamma (e^{i}) u, 
\end{equation}
where $\{e_{i}\}$ is an orthonormal basis of $\tansM$. 
Suppose that 
\begin{equation}
	\slashed{D} := \ri (\one \otimes \upgamma) \circ \nabla
\end{equation}
is the twisted Dirac operator on $\coTansM \otimes \sS$, then the \textbf{Rarita-Schwinger operator} is defined as~\cite[Def. 2.25]{Baer_Springer_2012} 
\begin{equation} \label{eq: def_Rarita_Schwinger_op}
	\rR := (I - \iota \circ \upgamma) \circ \slashed{D} : C^{\infty} (\sM; \sS^{3/2}) \to C^{\infty} (\sM; \sS^{3/2}). 
\end{equation}
The characteristic set $\Char{\rR}$ of $\rR$ coincides with the set $\coLightBun$ of lightlike covectors in dimensions $d \geq 3$  
(see e.g.~\cite[Lem. 2.26]{Baer_Springer_2012}) 
and $\rR$ is a symmetric differential operator whose Cauchy problem is \textit{well-posed} when $\spacetime$ is globally hyperbolic albeit $\rR^{2}$ is \textit{not} a normally hyperbolic operator
(see e.g~\cite[Rem. 2.27]{Baer_Springer_2012}). 
Moreover, a hermitian form~\eqref{eq: hermitian_form_Dirac_type_op} does \textit{not} exist in $d \geq 3$~\cite[Example 3.16]{Baer_Springer_2012}. 
\newline 

Originally, Rarita and Schwinger (in Minkowski spacetime) considered the twisted Dirac operator $\slashed{D}$ restricted to $\sS^{3/2}$ but not projected back to $\sS^{3/2}$~\cite[$(1)$]{Rarita_PR_1941} 
(see also~\cite[Sec. 2]{Homma_CMP_2019},~\cite[Prop. 2.7]{Wang_IndianaUMJ_1991} 
for Riemannian 
and~\cite[Rem. 2.28]{Baer_Springer_2012} 
for Lorentzian spin manifolds), 
that is,   
\begin{equation} \label{eq: constraint_Rarita_Schwinger}
	\slashed{D} \upharpoonright C^{\infty} (\sM; \sS^{3/2}) : C^{\infty} (\sM; \sS^{3/2}) \to C^{\infty} (\sM; \coTansM \otimes \sS), 
\end{equation}
in order to ensure the correct number of propagating degrees of freedom for spin-$3/2$ fields     
(see, for instance, the reviews~\cite{Sorokin_AIP_2005, Rahman_2015} for physical motivation and  different approaches used in Physics literature). 
The corresponding Rarita-Schwinger operator 
\begin{equation}
	\rQ := (I - \iota \circ \upgamma) \circ \big( \slashed{D} \upharpoonright C^{\infty} (\sM; \sS^{3/2}) \big)
\end{equation}
is then an overdetermined system and this constrained system limits possible curvatures of the spacetime~\cite{Gibbons_JPA_1976}   
(see also, e.g.~\cite[p. 856]{Homma_CMP_2019},~\cite{Hack_PLB_2013}): 
\begin{equation}
    \Big( \Ric - \frac{\ric}{d} \fg \Big)^{*} u = 0, 
\end{equation} 
where $\Ric$ resp. $\ric$ are the Ricci tensor resp. Ricci scalar curvature of $\sM$.  
In other words, solutions (Rarita-Schwinger field) of $\rQ$ exists only in \textit{Einstein spin manifolds}. 
However, $\rQ$ does \textit{not} admit any Green's operator~\cite[Rem. 2.28]{Baer_Springer_2012}. 
\newline 

Since the Cauchy problem for $\rR$ is well-posed on a globally hyperbolic spacetime $\spacetime$, the Rarita-Schwinger operator~\eqref{eq: def_Rarita_Schwinger_op} admits unique advanced $S^{\adv}$ and retarded $S^{\ret}$ Green's operators. 
Moreover, it also admits Feynman $S^{\Feyn}$ and anti-Feynman $S^{\aFeyn}$ Green's operators as $\Char \rR = \coLightBun$ (entailing four possible orientations of the geodesic relation). 
These can be constructed by assuming the construction of the causal Green's operator $S := S^{\ret} - S^{\adv}$ and following the same strategy (cf. the proof of Theorem~\ref{thm: Hadamard_bisolution_Dirac_type_op}) as in the case of a Dirac-type operator. 
However, the positivity issue is non-conclusive.
%
%
%
%
%
%
%
%
%
%
\subsection{Higher spin operators} 
The straightforward attempts to generalise  
Dirac operator on Minkowski spacetime for arbitrary spin~\cite{Dirac_PRSA_1936} 
in 
curved spacetimes\footnote{Not 
	necessarily be globally hyperbolic.}
leads to difficulties 
(see e.g~\cite[p. 324]{Illge_AnnPhys_1999} for a panoramic view 
and the reviews~\cite{Sorokin_AIP_2005, Rahman_2015}). 
A crucial advancement came through 
Buchdahl operator (in Riemannian manifold)~\cite{Buchdahl_JPA_1982} 
(for a Lorentzian formulation, see ~\cite[Exam. 2.24]{Baer_Springer_2012}) 
whose square turns out to be a normally hyperbolic operator~\cite[p. 8]{Baer_Springer_2012}, 
yet the minimum coupling principle seems to be violated and a "by hand" proposal is required in the original idea of Buchdahl. 
These minor imperfections were cured by 
W\"{u}nsch~\cite{Wuensch_GRG_1985},  
by 
Illege~\cite{Illge_ZAA_1992, Illge_CMP_1993} 
and by 
Illege and Schimming~\cite{Illge_AnnPhys_1999} 
for the massive case deploying the $2$-spinor formalism in $4$-dimensional curved spacetimes. 
In particular, the square of Buchdahl operator (as modified by W\"{u}nsch and by Illge) is a normally hyperbolic operator on certain twisted bundle. 
In contrast, there are a few open questions for the massless case~\cite{Frauendiener_JGP_1999} 
(see also the reviews~\cite{Sorokin_AIP_2005, Rahman_2015} for the contemporary status and other formulations used in Physics literature). 
Hence, the existence of a Feynman propagator is evident either by Theorem~\ref{thm: existence_Feynman_propagator_NHOp} or Theorem~\ref{thm: Hadamard_bisolution_Dirac_type_op} but the issue of positivity is non-conclusive at this stage. 
%
%
%
%
%
%
%
%
%
%
\section{Literature} 
The \textit{Lorentzian geometry of spacetime} 
is quite well-studied and presented in great detail, for instance, in the 
monographs~\cite{Hawking_CUP_1973, ONeill_Academic_1983, Wald_Chicago_1984, Beem_CRC_1996, Choquet_Bruhat_OUP_2008}.  
The notion of global hyperbolicity was introduced by 
Jean \textsc{Leray}~\cite{Leray_Princeton_1953} 
to ensure the well-posedness of the Cauchy problem.    
Later, this class of spacetimes was emphasized in the context of Einstein's general relativity notably by  
Roger \textsc{Penrose}~\cite{Penrose_RivNuovoCim_1969} 
by means of his cosmic censorship conjecture
(see e.g.~\cite[pp. 200 - 209, 299 - 308]{Wald_Chicago_1984}). 
The characterisation of a globally hyperbolic manifold in terms of a Cauchy hypersurface is originally due to 
Robert \textsc{Geroch}~\cite{Geroch_JMP_1970} 
who proved ``a topological version'' of Theorem~\ref{thm: globally_hyperbolic_Cauchy}. 
His topological splitting theorem has been sharpened in a series of four articles by 
Antonio \textsc{Bernal} and Miguel \textsc{S\'{a}nchez}~\cite{Bernal_CMP_2003, *Bernal_CMP_2005, *Bernal_LMP_2006, *Bernal_CQG_2007} 
who, in particular, have achieved an orthogonal smooth splitting in the level of spacetime metric. 
We refer the survey articles~\cite{Minguzzi_EMS_2008, Sanchez_AMS_2011} for details. 
\newline 

\textit{Fundamental solutions} 
of a normally hyperbolic operator have been constructed in great generality by 
Jacques \textsc{Hadamard}~\cite{Hadamard_ActaMath_1908, *Hadamard_Dover_2003} 
and by   
Marcel \textsc{Riesz}~\cite{Riesz_ActaMath_1949, *Riesz_CPAM_1960}. 
Contemporary expositions include, for example~\cite{Friedlander_CUP_1975, Guenther_AP_1988, Baer_EMS_2007}.
\newline 

The \textit{Feynman propagator} 
was introduced by 
Richard \textsc{Feynman}~\cite{Feynman_PR_1949} 
propounding the idea of 
Ernst \textsc{St\"{u}ckelberg}~\cite{Stueckelberg_HPA_1941}
that ``particles propagate forward in time whereas antiparticles propagate backward in time'' 
in the context of quantum field theory by means of the so-called time-ordered product (see Appendix~\ref{sec: timeordered_product}). 
His original definition was not rigorous and can only be understood well in Minkowski spacetime as recalled in Example~\ref{ex: Feynman_propagator_Minkowski}. 
The mathematical underpinning was given by
Johannes J. \textsc{Duistermaat} and Lars \textsc{H\"{o}rmander}~\cite[Sec. 6.6]{Duistermaat_ActaMath_1972}.  
In fact, the notion of distinguished parametrices for scalar pseudodifferential operators of real principal type was introduced in their seminal article where they have identified a geometric notion of pseudoconvexity which allows to prove the uniqueness of such parametrices. 
Although they were motivated and certainly aware of the developments in physics, it was realised only much later by 
Maciej \textsc{Radzikowski}~\cite{Radzikowski_CMP_1996} 
that the expectation value of the time-ordered product with respect to a state satisfying the Hadamard condition is the Feynman propagator. 
\newline 

The \textit{spin-Dirac operator} 
is named after 
Paul \textsc{Dirac}~\cite{Dirac_PRSA_1928, *Dirac_PRSA_1928_P2} 
who introduced this operator in the context of quantum field theory to describe the dynamics of an electron in Minkowski spacetime. 
The bundle language was not developed back then and so a geometric formulation of  Dirac operator was not available until the precise formulation of spin-structure. 
Subsequently, Dirac-type operators appear in the work of 
Michael \textsc{Atiyah}, Vijay \textsc{Patodi} and Isadore \textsc{Singer}~\cite{Atiyah_MathProcCam_1975} 
in their celebrated index theorem. 
We refer~\cite{Baum_1981, Baer_MathZ_2005} 
for a discourse on Lorentzian spin-geometry 
and~\cite{Branson_JFA_1992, Berline_Springer_2004} 
for a detailed account of Riemannian Dirac-type operators.

%% file: Gutzwiller_trace.tex
\chapter{Gutzwiller trace formula} 
\label{ch: Gutzwiller_trace}
\textsf{The Duistermaat-Guillemin trace formula of the time-flow on the kernel of a Dirac operator on a spatially compact globally hyperbolic stationary spacetime has been derived in this chapter.}
%
%
%
%
%
%
%
%
%
%
\section{Trace of a smoothing operator} 
We recall that a \textbf{trace} on an associative algebra $(\sA, \cdot)$ over $\C$ is a linear functional   
\begin{equation}
	\Tr : \sA \to \C \quad \textrm{such that} \quad \forall A, B \in \sA : \Tr ([A, B]_{-}) = 0,  
\end{equation}
where $[A, B]_{-} := A \cdot B - B \cdot A$ is the commutator defined by the product $\cdot$ structure on $\sA$. 
In particular, let $\sE_{\varSigma} \to \varSigma$ be a vector bundle over a compact manifold $\varSigma$ without boundary. 
Then $\PsiDO{-\infty}{\varSigma; \sE_{\varSigma} \otimes \halfDen \varSigma}$ is a Fr\'{e}chet algebra of smoothing operators and 
\begin{equation} \label{eq: def_trace_smooth_op}
	\Tr : \PsiDO{- \infty}{\varSigma; \sE_{\varSigma} \otimes \halfDen \varSigma} \to \C, ~ P \mapsto 
	\Tr P := \int_{M} \tr \big( \fP (x, x) \big),  
\end{equation}
where $\fP (x, x)$ is the restriction of the Schwartz kernel $\fP$ of $P$ to the diagonal in $\varSigma \times \varSigma$ and $\tr$ is the endomorphism trace defined by~\eqref{eq: def_endo_trace}. 
The traciality of this trace is the Fubini's theorem and this trace is the unique trace on $\PsiDO{- \infty}{\varSigma; \sE_{\varSigma} \otimes \halfDen \varSigma}$
(see e.g.~\cite[Sec. 1.1.7, 4.3.2]{Scott_OUP_2010}). 
%
%
%
%
%
%
%
%
%
%
\section{Stationary spacetimes} 
\label{sec: stationary_spacetime}
%
%
%
\begin{definition} \label{def: stationary_spacetime}
	A spacetime (Definition~\ref{def: spacetime}) $\spacetime$ is called  \textbf{stationary} $\SST$ if it admits a global timelike Killing flow $\varXi$. 
\end{definition}
%
%
%

A Killing flow $\varXi_{s} : \sM \to \sM$, per se, keeps the spacetime metric $\fg$ invariant: 
\begin{equation} \label{eq: def_Killing_flow}
	\varXi_{s}^{\ms \coTansM \times \coTansM} \fg = 0 \Leftrightarrow \pounds_{Z} \fg = 0 
\end{equation}
for any $s \in \R$, where $Z$ is the infinitesimal generator of $\varXi_{s}$, called the \textit{Killing vector field}. 
In other words, the Killing flow $\varXi_{s}$ is a \textit{spacetime isometry}. 
In local coordinates $(x^{i})$ on $\sM$, $Z = Z^{i} \partial_{i}$, and the preceding equation entails the \textbf{Killing's equation}: 
\begin{equation}
	\connectionLC_{i} Z_{j} + \connectionLC_{j} Z_{i} = 0, \quad i,j=1, \ldots, d, 
\end{equation}
where $\connectionLC_{X}$ is the Levi-Civita covariant derivative of $\spacetime$. 

%
%
%
\begin{example} \label{exm: KVF_Minkowski}
	Let $\spacetime = (\R^{4}, \etaup)$ be the $4$-dimensional Minkowski spacetime. 
	The set of all isometries of $\R^{4}$ consists of translations and Lorentz-transformations, and it has a group structure --- the \textit{Poincar\'{e} group}. 
	In inertial coordinates $(x^{i}) = (t, x, y, z)$: 
	\begin{itemize}
		\item 
		$(Z^{i}) = (0, 0, 0, 1)$ is the Killing vector field generating the translations in the $z$-direction. 
		\item 
		$(Z^{i}) = (0, 0, -z, y)$ is the Killing vector field generating the rotations  around the $x$-axis.
		\item 
		$(Z^{i}) = (x, t, 0, 0)$ is the Killing vector field generating the boosts along the $x$-axis.
	\end{itemize}
\end{example}
%
%
%

The Killing flow is physically interpreted as the \textit{flow of time} which allows a canonical $1 + (d - 1)$ dimensional non-unique topological (resp. local geometric) decomposition of $\sM$ (resp. $\fg$). 
This can be achieved either by the 
projection formalism~\cite{Geroch_JMP_1971} 
(see also, e.g.~\cite{Harris_CQG_1992, Javaloyes_CQG_2008})  
or by the Kaluza-Klein reduction 
(see e.g.~\cite[Sec. XIV.2]{Choquet_Bruhat_OUP_2008}). 
In what follows, we can (and at some point, will) consider a slightly specialised version of stationary spacetimes. 

%
%
%
\begin{definition} \label{def: std_stationary_spacetime}
	Let $\Cauchy$ be a Riemannian manifold. 
	A \textbf{standard stationary spacetime} is defined as the triplet $(\sM, \fg, Z)$ where  
	(see e.g.~\cite[$(7.1)$]{Sanchez_AMS_2011},~\cite[Def. 3.3]{Sanders_IJMPA_2013})  
    \begin{equation}
    	\sM := \R \times \varSigma
    \end{equation}
	is the spacetime manifold admitting a complete and timelike Killing vector field $Z$ and a local coordinate chart $\big( U, (t, x^{i}) \big)$ such that $Z \upharpoonright U = \partial_{t}$ and the spacetime metric 
	\begin{equation}\label{eq: def_spacetime_metric_SSST}
        \fg \upharpoonright U := \upbeta^{2} \, \rd t^{2} - \fh_{ij} (\rd x^{i} + \upalpha^{i} \, \rd t) (\rd x^{j} + \upalpha^{j} \, \rd t), 
		\quad i,j = 2, \ldots, d,     
    \end{equation}
    where $\upalpha$ is the shift vector field and $\upbeta$ is the lapse function, both independent of $t$.   
\end{definition}
%
%
%

Not all standard stationary spacetimes are globally hyperbolic (Definition~\ref{def: globally_hyperbolic_spacetime}) and as argued in Section~\ref{sec: globally_hyperbolic_spacetime}, we are only interested in those which are. 
Our desired class is achieved pertaining to the Fermat metric 
(see e.g.~\cite[Sec. 7]{Sanchez_AMS_2011} and references therein) 
\begin{equation} \label{eq: def_Fermat_metric_SSST}
	\ff := \frac{1}{\betaup^{2} - \| \alphaup \|_{\fh}^{2}} \fh_{ij} \alpha^{i} \rd x^{j} + \sqrt{\frac{1}{\betaup^{2} - \| \alphaup \|_{\fh}^{2}} \fh_{ij} \rd x^{i} \rd x^{j} + \bigg( \frac{\fh_{ij} \alphaup^{i} \rd x^{j}}{\betaup^{2} - \| \alphaup \|_{\fh}^{2}} \bigg)^{2}   }
\end{equation}
on $\varSigma$ for a given splitting~\eqref{eq: def_spacetime_metric_SSST}. 
This metric induces a (not necessarily symmetric) distance $\fd_{\ff}$ on $\varSigma$ which in turn induces forward $\bar{\B}^{+} (x, \rr) := \{ y \in \sM \,|\, \fd_{\ff} (x, y) \leq \rr \}$ and backward $\bar{\B}^{-} (x, \rr) := \{ y \in \sM \,|\, \fd_{\ff} (y, x) \leq \rr \}$ closed balls. 
Analogously, one has forward (backward) Cauchy sequences and forward (backward) completeness. 
Moreover, one defines the closed symmetrised balls $\bar{\B}_{\mathrm{s}} (x, \rr)$ corresponding to the symmetrised distance $\big( \fd_{\ff} (x, y) + \fd_{\ff} (y, x) \big) / 2$. 
In terms of these terminologies, the analysis of~\cite[Thm. 4.3b, 4.4, Cor. 5]{Caponio_RMI_2011} 
can be summarised as

%
%
%
\begin{theorem}
	The following properties are equivalent for a standard stationary spacetime $(\sM, \fg, Z)$ 
	(see e.g., the reviews~\cite[Thm. 7.2]{Sanchez_AMS_2011},~\cite[Thm. 3.1]{Sanders_IJMPA_2013}).   
	\begin{itemize}
		\item 
		$(\sM, \fg, Z)$ is globally hyperbolic. 
		\item 
		The closed symmetrised balls $\bar{\B}_{\mathrm{s}}$ of the Fermat metric~\eqref{eq: def_Fermat_metric_SSST} associated to one (and then to any) standard stationary splitting of $\spacetime$ are compact. 
	\end{itemize}
	Furthermore, the slices associated to a standard stationary splitting are Cauchy hypersurfaces if and only if the Fermat metric associated to that splitting is both forward and backward complete. 
\end{theorem}
%
%
%

A sufficient condition is given by 

%
%
%
\begin{corollary}
	Let $(\sM, \fg, Z)$ be a standard stationary spacetime. 
	If $\fh$ is complete and $\fd_{0}$ is the distance function on $\Cauchy$ to some (and then any) fixed point $x_{0} \in \varSigma$ such that  
	\begin{equation}
		\sup \left\{ \frac{\norm{\alphaup}{\fh}}{\fd_{0}} (x), \frac{ \sqrt{\betaup^{2} - \norm{\alphaup}{\fh}^{2}} }{\fd_{0}} (x) ~|~ x \in \varSigma, \fd_{0} (x) > 1 \right\} < \infty, 
	\end{equation}
	then $(\sM, \fg, Z)$ is globally hyperbolic and all the slices $\varSigma_{t}$ are Cauchy hypersurfaces.  
	\newline 

	In particular, the preceding finite bound holds if $\alphaup$ is sublinear: $\norm{\alphaup}{\fh} < \cst_{1} \fd_{0} (x) + \cst_{2}$ and $\betaup$ is subquadratic: $\betaup (x) < \cst_{3} \rd_{0} (x) + \cst_{4}$ for some constants 
	(see e.g.~\cite[Cor. 7.1]{Sanchez_AMS_2011}). 
\end{corollary}
%
%
%

\begin{remark} \label{rem: SST}
	As in the terminologies of Definitions~\ref{def: stationary_spacetime} and~\ref{def: std_stationary_spacetime}:  
	\begin{enumerate} [label=(\alph*)]
		\item \label{rem: complete_Killing_flow}
		The completeness means that the Killing flow is globally defined
        \begin{equation} \label{eq: Killing_flow}
        	\varXi : \R \times \sM \to \sM, ~ (s, x) \mapsto \varXi (s, x) =: \varXi_{s} (x), 
            \quad 
            \rd \varXi (s, x; \partial_{s}, 0) = Z \big( \varXi (s, x) \big), 
        \end{equation}
        where $Z$ is the global timelike Killing vector field, locally given by $\partial_{s}$. 
		We remark that some authors define stationary spacetimes without imposing the completeness assumption but this is crucial for the thesis.
		\item \label{rem: SSST}
		If one of the timelike Killing vector fields of a globally hyperbolic stationary spacetime is complete then the spacetime is a standard stationary one~\cite[Thm. 2.3]{Candela_AdvMath_2008}. 
        \item 
		The splitting~\eqref{eq: def_spacetime_metric_SSST} is \textit{not} unique. 
		\item 
        A globally hyperbolic stationary spacetime is called \textbf{spatially compact} if and only if its Cauchy hypersurface is compact. 
        Such spacetimes are \textit{geodesically complete}~\cite[Lem. 1.1]{Anderson_AHP_2000}.  
		\item \label{rem: hypersurface_orthogonality_SST}
		On a generic $d > 2$-dimensional (standard) stationary spacetime the Killing vector field $(\partial_{t} =) Z$ is \textit{not} orthogonal to $\varSigma_{t}$ because the corresponding $1$-form $Z^{\flat}$ does not satisfy the \textit{hypersurface-orthogonality} condition: 
		\begin{equation}
			Z^{\flat} \wedge \rd Z^{\flat} \neq 0.  
		\end{equation}
		In other words, the orthogonal geometric distribution of $Z$ is non-involutive.  
		Physically this means that the neighbouring orbits of $Z$ can twist around each other. 
		In $d=2$ this cannot happen, i.e., every Killing vector field is at least locally hypersurface-orthogonal because there is no freedom to rotate. 
		If one imposes the condition that $\varSigma_{t}$ is orthogonal to the orbits of the spacetime isometry, then $\SST$ is called a \textbf{static spacetime} and one has a \textit{canonical non-unique global time-coordinate} $t$. 
		In this case, $\upalpha$ vanishes identically so that there is no $\rd t \, \rd x^{i}$-type cross terms in~\eqref{eq: def_spacetime_metric_SSST}. 
		Additionally, if we demand that $Z$ has a constant norm, then a static spacetime is called an \textbf{ultrastatic spacetime}. 
		This enforces $\upbeta$ to be the identity function 
		(see e.g.~\cite[pp. 120-122]{Beem_CRC_1996}~\cite[(B.3.6) and App. C.3]{Wald_Chicago_1984}). 
		\item 
		On a $d$-dimensional Lorentzian manifold, there can be at most $d (d+1)/2$ linearly independent Killing vector fields. 
		In our physical universe, currently $d=4$. 
		Hence, it can admit maximum $10$ Killing vector fields and this maximally symmetric spacetime is known as the \textit{de Sitter spacetime}. 
	\end{enumerate}
\end{remark}
%
%
%
%
%
%
%
%
%
%
\section{Gutzwiller trace formula for Dirac operators}
In order to get an intuition of the Gutzwiller trace formula for Dirac-type operators on a stationary spacetime we begin with the simplest case.  
%
%
%
%
%
%
%
%
%
%
\subsection{Poisson summation formula on $\bbS$}
Since $\rH^{1} (\bbS; \Z_{2}) = \Z_{2}$ (cf. Example~\ref{ex: spin_Dirac_op}), a circle $\bbS$ admits \textit{two inequivalent} spin-structures, as described below 
(see e.g.~\cite[Sec. 3]{Baer_Seminar_2000}). 
\begin{itemize}
	\item 
	\textbf{Trivial spin-structure}:  
	One uses $\rSO (1) = \{1\}$ and $\rSpin (1) = \Z_{2} = \{ + 1, -1 \}$. 
	Then the frame bundle $(\sQ, \rSO (1), \bbS)$ is trivial and so is the spin-structure   
	\begin{equation}
		\sP_{\triv} = \bbS \times \rSpin (1), 
		\quad 
		\sS_{\triv} := \sP_{\triv} \times \C, 
	\end{equation}
	where the double covering $\varTheta$ is the projection of the first factor. 
	Hence, the spinors on $\bbS$ are simply ($\C$-valued) functions on $\bbS$ which are identified with $T$-periodic functions on $\R$:  
    \begin{equation}
        C^{\infty} (\bbS; \sS_{\triv}) 
		= C^{\infty} (\bbS) 
		= \{ u \in C^{\infty} (\R) \,|\, \forall T > 0 : u (x + T) = u (x) \}. 
    \end{equation}
    The massless spin-Dirac operator on $\bbS$ is nothing but the derivative operator $\rD_{x}$:  
    \begin{equation} \label{eq: massless_spin_Dirac_op_S_1}
        \slashed{D}_{\bbS} := - \ri \frac{\rd}{\rd x} = \rD_{x}.  
    \end{equation}
    Therefore $\spec \slashed{D}_{\bbS} = \spec \rD_{x}$ and the Poisson summation formula~\eqref{eq: Poisson_summation_formula_derivative_circle} is essentially the Poisson summation formula for $\slashed{D}_{\bbS}$ on $\bbS$ for the trivial spin-structure. 
	\item 
	\textbf{Non-trivial spin-structure}: 
	In this case, for $T > 0$,  
	\begin{equation}
		\sP_{\nonTriv} = \big( [0, T] \times \rSpin (1) \big) / \sim,  
	\end{equation}
	where $\sim$ identifies $0$ with $T$ while it interchanges two elements of $\rSpin (1) = \{ +1, -1 \}$. 
    Thus, sections of $\sS_{\nonTriv}$ are no longer $C^{\infty} (\bbS)$ rather they are anti-periodic:  
    \begin{equation}
    	C^{\infty} (\bbS; \sS_{\nonTriv}) = \{ u \in C^{\infty} (\R) \,|\, \forall T > 0 : u (x + T) = - u (x)  \}.
    \end{equation}
    In contrast to~\eqref{eq: spec_spin_Dirac_trivial_circle}, the spectrum reads 
    \begin{equation} \label{eq: spec_spin_Dirac_nontrivial_circle}
        \spec \slashed{D}_{\bbS} = \{ \lambda_{n} | n \in \Z \}, 
        \quad 
        \lambda_{n} := \left( n + \frac{1}{2} \right) \omega, 
        \quad 
        \omega := \frac{2 \pi}{T}
    \end{equation}
    corresponding to the orthonormal eigenspinors $\phiup_{n} := \re^{\ri \lambda_{n} x} / \sqrt{T}$. 
    One can obtain a Poisson summation formula for non-trivial spin-structure analogously and it is structurally similar to that of trivial spin-structure albeit not exactly the same. 
\end{itemize}

In relativistic language, the Poisson summation formulae for trivial and non-trivial spin-structures on $\bbS$ describe an exact trace formula $\Tr U_{t}$ for $U_{t} := \re^{- \ri t \slashed{D}_{\bbS}}$ on the product spin-manifold $(\R \times \bbS, \rd t^{2} - \rd x^{2})$. 
The thesis is devoted to an enormous generalisation of this result in the sense that $\bbS$ in the preceding examples is replaced by a compact Cauchy hypersurface $\varSigma$ and the role the product metric $\rd t^{2} - \rd x^{2}$ is played by a stationary spacetime metric (cf. Definitions~\ref{def: stationary_spacetime},~\ref{def: std_stationary_spacetime} and Remark~\ref{rem: SST}~\ref{rem: SSST}). 
%
%
%
%
%
%
%
%
%
%
\subsection{Primary results}
\label{sec: result_Gutzwiller_trace}
Let $\sE \to \sM$ be a vector bundle over a $d \geq 2$ dimensional spatially compact globally hyperbolic stationary spacetime (Section~\ref{sec: stationary_spacetime}) $\SST$. 
That is, the spacetime manifold $\sM$ is homeomorphic to the product manifold $\R \times \varSigma$ where $\varSigma \subset \sM$ is a spacelike compact Cauchy hypersurface without boundary and $\spacetime$ admits a \textit{complete} timelike Killing vector field $Z$ whose flow is $\varXi$ (Remark~\ref{rem: SST}~\ref{rem: complete_Killing_flow}).  

%
%
%
\begin{assumption} \label{asp: trace_formula}
	We consider Dirac-type operators (Definition~\ref{def: Dirac_type_op}) $D$ on a vector bundle $\sE \to \sM$ over a globally hyperbolic stationary spacetime $\SST$, subjected to the following assumptions:  
    \begin{enumerate}[label=(\alph*)] 
        \item 
        $\sE$ is endowed with a sesquilinear form $(\cdot|\cdot)$ invariant under the Killing flow $\varXi_{s}^{*}$ such that $D$ is symmetric;  
        \label{asp: Direc_op_symmetric_SST}
        \item 
        Given an arbitrary but fixed future-directed unit normal covector field $(\cdot, \zeta)$ on $\sM$ along any Cauchy hypersurface $\varSigma$, 
        \begin{equation} \tag{\ref{eq: hermitian_form_Dirac_type_op}}
            \langle \cdot | \cdot \rangle := \, \big( \symb{D} (\cdot, \zeta ) \cdot \big| \cdot \big)   
        \end{equation}
		is a fibrewise hermitian form on the bundle of Clifford modules (Section~\ref{sec: Dirac_type_op}) $\big( \sE \to \sM, \symb{D}, (\cdot|\cdot) \big)$, where $\symb{D}$ is the principal symbol of $D$;  
        \label{asp: hermitian_form_Dirac_op_SST}
        \item 
        $D$ commutes with the induced Killing flow $\varXi_{s}^{*}$ on $\sE$ for all $s \in \R$.
        \label{asp: commutator_Dirac_op_Killing_flow}
    \end{enumerate} 
\end{assumption}
%
%
%

Regarding the second assumption, we recall the comments after Assumption~\ref{asp: hermitian_form_Dirac_type_op}. 
In the present setting, we use $Z$ to fix the time-orientation and note that  on a static spacetime, one can choose the Killing covector field $(\partial_{t})^{\flat}$ (up to normalisation) for $\zeta$ 
(see Remark~\ref{rem: SST}~\ref{rem: hypersurface_orthogonality_SST} for details).   
\newline 

As noted before, the Clifford module bundle $\big( \sE \to \sM, \symb{D}, (\cdot | \cdot) \big)$ is naturally furnished with the unique Weitzenb\"{o}ck connection $\connectionE$ (\eqref{eq: Weitzenboeck_formula_Dirac_type}) induced by $D$ and $\fg$. 
This connection induces a parallel transport map $\hat{\varXi}_{s} : \sE_{\cdot} \to \sE_{\varXi_{s} (\cdot)}$ along the spacetime isometry $\varXi_{s} (\cdot)$. 
Subsequently, the following diagram commutes:   

%
%
%
\begin{center}
	\begin{tikzpicture}
		\node (a) at (0, 0) {$\sM$};
		\node (b) at (3, 0) {$\sM$};
		\node (c) at (0, 1.5) {$\varXi_{s}^{*} \sE$}; 
		\node (d) at (3, 1.5) {$\sE$};
		\node (e) at (-3, 0) {$\sM$};
		\node (f) at (-3, 1.5) {$\sE$}; 
		\draw[->] (a) -- (b); 
		\draw[->] (c) -- (d); 
		\draw[->] (c) -- (a); 
		\draw[->] (d) -- (b); 
		\draw[->] (f) -- (e); 
		\draw[->] (e) -- (a); 
		\draw[->] (f) -- (c); 
		\node[above] at (1.5, 0) {$\varXi_{s}$}; 
		\node[above] at (-1.5, 0) {id}; 
		\node[above] at (1.5, 1.5) {$\cong$}; 
		\node[above] at (-1.5, 1.5) {$\hat{\varXi}_{s}$}; 
	\end{tikzpicture}
	\captionof{figure}[Pullback of a Clifford module bundle via a spacetime isometry]{Pullback $\varXi_{s}^{*} \sE$ of a bundle of Clifford modules $\big( \sE \to \sM, \symb{D}, (\cdot \vert \cdot) \big)$ over a stationary spacetime $\SST$ via a spacetime isometry $\varXi_{s} : \sM \to \sM$.}
	\label{fig: bundle_lift_spacetime_isometry}
\end{center}
%
%
%

Let $\pounds_{\! Z}$ be the Lie derivative (given by~\eqref{eq: def_Lie_derivative}) on $\sE$ with respect to the Killing vector field $Z$ and (cf.~\eqref{eq: def_Lie_derivative})
\begin{equation} \label{eq: def_L}
	\LieDeri := - \ri \pounds_{\! Z} : \secsME \to \secsME
\end{equation}
so that the induced Killing flow on $\secsME$ is expressed as $\varXi_{s}^{*} = \re^{\ri s \LieDeri}$.  
By hypotheses, 
\begin{equation}
    [\LieDeri, D]_{-} = 0.  
\end{equation}
Thus, on $\ker D$, the eigensections $\uppsi_{n}$ of $\LieDeri$ are the joint eigensections: 
\begin{equation}
	D \uppsi_{n} = 0, \qquad \LieDeri \uppsi_{n} = \lambda_{n} \uppsi_{n}. 
\end{equation}
We equip $\ker D$ with the hermitian inner product $\scalarProdTwo{\cdot}{\cdot}$ to have the Hilbert space $(\ker D, \scalarProdTwo{\cdot}{\cdot})$. 

%
%
%
\begin{theorem} \label{thm: spectrum_L}
	Under Assumption~\ref{asp: trace_formula}, the spectrum of $\LieDeri :=$~\eqref{eq: def_L} on the Hilbert space $\sH := (\ker D, \langle \cdot | \cdot \rangle := \eqref{eq: def_hermitian_form_Dirac_type_op})$ is purely discrete and comprises infinitely many real eigenvalues that grow polynomially and accumulate at $\pm \infty$.  
\end{theorem}
%
%
%

Thus, we can restrict our attention entirely to smooth sections $C^{\infty} (\sM; \sE)$ of $\sE$ owing to the elliptic regularity and $\tr_{\sH} \re^{\ri t \LieDeri}$ makes sense as a distribution with  (cf.~\eqref{eq: trace_time_evolution_op_tr_Killing_flow}) the identification 
\begin{equation} \label{eq: Tr_Cauchy_evolution_op_Dirac_Tr_Killing_flow}
	\Tr U_{t} = \tr_{\sH} \varXi_{t}^{*},   
\end{equation}
where $U_{t}$ is the time evolution operator of the Dirac equation~\eqref{eq: Dirac_eq}. 
\newline 

Since the classical dynamics is governed by the principal symbol $\symb{D}$ of $D$, its characteristic set $\Char D$ can be viewed as a classical limit of $\ker D$ as per se geometric quantisation. 
Thus, our classical phase space is the lightcone bundle $\coLightBun \to \sM$ over $\SST$ where the metric-Hamiltonian $H_{\fg} :=$~\eqref{eq: def_Hamiltonian_metric} vanishes and the reduced phase space is the conic symplectic manifold $\sN$ of scaled-lightlike geodesics in the cotangent bundle~\cite{Penrose_1972} 
(see also~\cite{Khesin_AdvMath_2009}). 
Being a globally hyperbolic spacetime, $\sM$ does \textit{not} admit any \textit{closed timelike geodesic}. 
Therefore, the notion of Lorentzian analogue of periodic trajectories is defined by means of the spacetime isometry $\varXi_{t}$ induced (reduced) symplectic flow $\varXi_{t}^{\sN}$ on $\sN$, whose Hamiltonian is given by~\cite[Lem. 1.1]{Strohmaier_AdvMath_2021} 
\begin{equation} \label{eq: def_Hamiltonian_scaled_lightlike_geodesic}
	H : \sN \to \R, ~\gamma \mapsto H (\gamma) := \fg \Big( \frac{\rd c}{\rd s}, Z \Big), 
\end{equation} 
where $\R \ni s \mapsto c (s) \in \sM$ is any lightlike geodesic on $\sM$ and the value $\fg (\nicefrac{\rd c}{\rd s}, Z)$ is independent of the cotangent lift $\gamma$ of $c$. 
This Hamiltonian is positive as $Z$ is timelike. 
For any $\tE \in \R_{+}$, we denote the constant $\tE$-energy surface by  
\begin{equation} \label{eq: def_cst_energy_surface_Hamiltonian_scaled_lightlike_geodesic}
	\sN_{\tE} := \{ \gamma \in \sN \,|\, H (\gamma) = \tE \}.   
\end{equation}
The set of \textbf{periods} resp. \textbf{periodic lightlike geodesics} of $\varXi_{s}^{\sN}$ are then given 
by~\cite[$(6)$]{Strohmaier_AdvMath_2021}  
\vspace*{-0.4cm}
\begin{subequations}
	\begin{eqnarray}
		\cP 
		& := & 
		\{ T \in \R_{+} \,|\, \exists \gamma \in \sN : \varXi_{T}^{\sN} (\gamma) = \gamma \}, 
		\label{eq: def_period_Killing_flow}
		\\  
		\cP_{T} 
		& := & 
		\{ \gamma \in \sN \,|\, \varXi_{T}^{\sN} (\gamma) = \gamma \}. 
		\label{eq: def_periodic_pt_Killing_flow}
	\end{eqnarray}
\end{subequations}
Recall, the set of all lengths of periodic geodesics on a manifold counted with multiplicities is called the \textit{length spectrum} of the manifold. 
Our first finding is that $\Tr U_{t}$ determines the Lorentzian length spectrum of $\sN$ in the sense described below. 

%
%
%
\begin{proposition} \label{prop: length_spectrum}
	Let $\sE \to \sM$ be a vector bundle over a spatially compact globally hyperbolic stationary spacetime $\SST$, endowed with a sesquilinear form $(\cdot|\cdot)$, and $D$ a Dirac-type operator on $\sE$ whose principal symbol is $\symb{D}$ so that $\big( \sE \to \sM, \symb{D}, (\cdot|\cdot) \big)$ is a bundle of Clifford modules over $\SST$.   
	If $U_{t}$ is the time evolution operator of $D$ then under Assumption~\ref{asp: trace_formula}, $\Tr U_{t}$ is a distribution on $\R$ and its singular support 
	\begin{equation}
		\singsupp{\Tr U_{t}} \subset \{0\} \cup \cP,  
	\end{equation}
	where $\cP$ is the set of periods of the induced Killing flow $\varXi_{T}^{\sN}$ on the manifold of scaled-lightlike geodesics $\sN$, given by~\eqref{eq: def_period_Killing_flow}.   
\end{proposition}
%
%
%

Lagrangian distributions (Section~\ref{sec: Lagrangian_distribution}) offer the elegant characterisation of $\Tr U_{t}$, as stated below.  

%
%
%
\begin{theorem} \label{thm: trace_formula_t_zero}
    As in the set-up of Proposition~\ref{prop: length_spectrum}, $\Tr U_{T}$ is the Lagrangian distribution $I^{d - 7/4} (\R, \Lambda_{T})$ where $d := \dim \sM$ and 
	\begin{equation} \label{eq: def_canonical_relation_Tr_U_t}
        \Lambda_{T} := \{ (T, \tau) \in \R \times \R_{-} \,|\, \exists \gamma \in \sN : \varXi_{T}^{\sN} (\gamma) = \gamma, \tau = - H (\gamma) \},   
	\end{equation}
	is the Lagrangian manifold $\Lambda_{T}$, where the Hamiltonian $H$ on $\sN$ is given by~\eqref{eq: def_Hamiltonian_scaled_lightlike_geodesic}. 
	Furthermore,    
	\begin{equation}
		(\Tr U_{0}) (t) = u_{0} (t) + v_{0} (t), 
	\end{equation}
	where $v_{0}$ is a distribution on $\R$ which is smooth in the vicinity of $t=0$ and $u_{0}$ is a Lagrangian distribution admitting the following singularity expansion around $t=0$:   
	\begin{eqnarray}
		u_{0} (t) 
		& \sim & 
		r (d-1) \frac{\vol (\sN_{H \leq 1})}{(2 \pi)^{d-1}} \mu_{d-1} (t) + \mathrm{c}_{2} \, \mu_{d-2} (t) + \ldots, 
		\nonumber \\ 
		r 
		& := & 
		\rk{\sE} \int_{\varSigma} \fg_{x}^{-1} (\eta, \zeta) \, \dVolh (x), 
		\nonumber \\ 
		\mu_{d-k} (t) 
		& := & 
		\int_{\R_{\geq 0}} \re^{- \ri t \tau} \tau^{d-1-k} \rd \tau, \quad k = 1, 2, \ldots.   
	\end{eqnarray}
	Here, $\mu_{d-k}$ is a distribution on $\R$ given by the preceding oscillatory integral, $\eta$ is any lightlike covector on $\sM$ and $\zeta$ as in Assumption~\ref{asp: trace_formula}, both restricted to a Cauchy hypersurface $\varSigma$ of $\sM$, $\mathrm{c}_{2}$ is some constant (Dirac-wave trace invariant), and $\vol \sN_{H \leq 1}$ is the volume of $\sN_{H \leq 1}$.
\end{theorem}
%
%
%

As remarked in~\ref{rem: SST}~\ref{rem: SSST}, $\SST$ can be considered as a (globally hyperbolic) \textit{standard} stationary spacetime since $Z$ is complete. 
Then, we can use the spacetime metric~\eqref{eq: def_spacetime_metric_SSST} without loss of generalities and the volume of $\sN_{H \leq 1}$ is given by~\cite[$(15)$]{Strohmaier_AdvMath_2021}
\begin{equation} \label{eq: vol_energy_surface}
	\vol \sN_{H \leq 1} = \vol (\mathbb{B}^{d-1}) \int_{\varSigma} \upbeta (x) \big( \upbeta^{2} (x) - \fh_{x} (\upalpha, \upalpha) \big)^{-d/2} \dVolh (x),  
\end{equation}
where $\vol (\mathbb{B}^{d-1})$ is the volume of the unit ball $\mathbb{B}^{d-1} \subset \R^{d-1}$, and $\upbeta$ and $\upalpha$ are the lapse function and shift vector field, respectively, appearing in spacetime metric~\eqref{eq: def_spacetime_metric_SSST}.
\newline

Since $\LieDeri$ has discrete eigenvalues, let us introduce its eigenvalue counting  
function\footnote{One   
	can define it with negative eigenvalues as well.
} 
\begin{equation}
	\mathsf{N} (\lambda) := \# \{ n \in \N \,|\, 0 \leq \lambda_{n} \leq \lambda \}. 
\end{equation}

%
%
%
\begin{corollary}[Weyl law] \label{cor: Weyl_law}
	As in the terminologies of Theorems~\ref{thm: spectrum_L} and~\ref{thm: trace_formula_t_zero}, the Weyl eigenvalue counting function of $\LieDeri$ has the asymptotics 
	\begin{equation}
		\fN (\lambda) = \Big( \frac{\lambda}{2 \pi} \Big)^{d-1} r \, \vol (\sN_{H \leq 1}) + O (\lambda^{d-2}), 
		\quad \textrm{as} \quad \lambda \to \infty. 
	\end{equation}
\end{corollary}
%
%
%

In order to describe $\Tr U_{T}$ for the non-trivial periods $T \neq 0$, we let $\connectionEndPiE_{X_{\fg/2}}$ be the $D^{2}$-compatible (Weitzenb\"{o}ck) covariant derivative (Definition~\ref{def: P_compatible_connection}) with respect to the geodesic vector field $X_{\fg/2}$ on the bundle $\pi^{*} \Hom{\sE, \sE} \to \dotCoTansM$ where $\pi : \coTansM \to \sM$ is the cotangent bundle. 
By $\Hol := C_{\mathrm{Hol}}^{\infty} \big( \dotCoTansM; \pi^{*} \Hom{\sE, \sE} \big)$ we will denote the \textit{set of all sections of $\pi^{*} \Hom{\sE, \sE}$ invariant under the holonomy group of the  parallel transporter $\cT$ with respect to} $\connectionEndPiE_{X_{\fg/2}}$.  
\newline 

We suppose $\rP_{\gamma}$ is the (linearised) Poincar\'{e} map of a periodic geodesic $\gamma$ 
(see e.g.~\cite[Sec. 7.1]{Abraham_AMS_1978},~\cite{Strohmaier_AdvMath_2021}). 
Then $\gamma$ is called \textbf{non-degenerate} if $1$ is not an eigenvalue of $\rP_{\gamma}$. 
\newline 

With all these devises, the precise characterisation of $\Tr U_{T \neq 0}$ can be stated as follows. 

%
%
%
\begin{theorem} \label{thm: trace_formula_t_T}
	As in the set-ups of Proposition~\ref{prop: length_spectrum} and Theorem~\ref{thm: trace_formula_t_zero}: if the periods $T$ are discrete and the set $\cP_{T} := \eqref{eq: def_periodic_pt_Killing_flow}$ of periodic lightlike geodesics of $\varXi_{T}^{\sN}$ is a finite union of non-degenerate periodic orbits $\gamma$, then  
	\begin{equation}
		(\Tr U_{T}) (t) = \sum_{\gamma \in \cP_{T}} u_{\gamma} (t) + v_{\ms T} (t), 
	\end{equation}
    where $v_{\ms T}$ is a distribution on $\R$ that is smooth in the vicinity of $t=T$ and $u_{\gamma} (t)$'s are Lagrangian distributions having singularities at $t=T_{\gamma}$ with the asymptotic expansion  
	\begin{eqnarray} \label{eq: asymptotic_expansion_Lagrangian_dist_T}
        \lim_{t \to T} (t - T_{\gamma}) u_{\gamma} (t) 
        & \sim & 
        \frac{1}{2 \pi}  
        \int_{ \kF_{T_{\gamma}} } \tr \Big( \symb{D} (\gamma) \, \cT_{\varXi_{\ms T_{\gamma}}^{\ms \sN} (\gamma)} \symb{D} (x, \zeta) \Big) \frac{\re^{- \ri \pi \km (\gamma) / 2} |\rd T_{\gamma}|}{\sqrt{| \det (I - \rP_{\gamma}) |}} \nu
        + \ldots, 
        \nonumber \\ 
        \nu 
        & := & 
        \int_{\R_{\geq 0}} \re^{- \ri (t - T) \tau} \Big( \frac{\tau}{2 \pi} \Big)^{d-2} \rd \tau, 
	\end{eqnarray}
    where $\km (\gamma)$ resp. $\rP_{\gamma}$ are the Maslov index resp. the (linearised) Poincar\'{e} map of $\gamma$, $\mathcal{T}_{\gamma}$ is an element of the holonomy group $\Hol_{\gamma}$ with respect to the $D^{2}$-compatible Weitzenb\"{o}ck connection at base point $\gamma \in \sN$ whose projection on $\sM$ is $x$, and 
	\begin{equation} \label{eq: def_F_T}
        \kF_{T_{\gamma}} := \left\{ \gamma \in \sN \,|\, \tau = - H (\gamma), \varXi_{T_{\gamma}}^{\sN} (\gamma) = \gamma \right\} 
	\end{equation}
    is the fixed point set of $\varXi_{T_{\gamma}}^{\sN}$.  
\end{theorem}
%
%
%

At this point we would like to comment on the issue of the parallel transporter raised in the introduction (Section~\ref{sec: Sandoval_trace_formula}) of this thesis and reproduce  
Sandoval's result~\cite[Thm. 2.8]{Sandoval_CPDE_1999} 
as a special case of this theorem. 
We set $\upalpha = 0$ and $\upbeta = 1$ in~\eqref{eq: def_spacetime_metric_SSST} so that the Dirac equation~\eqref{eq: Dirac_eq} on the ultrastatic spacetime becomes 
\begin{equation} \label{eq: Dirac_Hamiltonian_ultrastatic}
	- \ri \partial_{t} u = H_{D} u, 
    \quad 
    H_{D} := \fc (\rd t) \big( \hat{D} - U \big), 
\end{equation}
where $H_{D}$ is the Dirac Hamiltonian, $\fc := \eqref{eq: def_Clifford_multiplication}$ is the Clifford multiplication, and $U$ is the potential in~\eqref{eq: Dirac_op} allowed by the most general Dirac-type operator. 
Sandoval plumped for $U = 0$ and used the parallel transporter $\check{\cT}$ corresponding to $\check{D}$ instead of $\hat{D}$ and hence $\check{T}$ together with a ``suitable average'' of $W := \hat{D} - \check{D}$ showed up under the $\tr$ in~\eqref{eq: asymptotic_expansion_Lagrangian_dist_T} in her work. 
However, no such ad-hoc choices are necessary in our formulation at all, rather $\cT$ naturally induces the parallel transport $\hat{\cT}$ for the preceding choices.      
\newline 

Note, the Maslov index (see Appendix~\ref{sec: Maslov_index}) $\km (\gamma)$ is essentially the Conley-Zehnder index of $\gamma$~\cite{Meinrenken_JGP_1994} 
and the factor $\re^{- \ri \pi \km (\gamma) / 2}$ is often known as the Maslov factor.  
\newline

In the future, several generalisations have been planned. 
For instance, we intend to address the spectral asymptotics on stationary black holes and explore the semi-classical regime.   
We also wish to extend the study for Hodge-d'Alambertians and connect with interesting applications on relativistic quantum chaos on curved spacetimes. 
%
%
%
%
%
%
%
%
%
%
\subsection{Proof strategy and novelty} 
\label{sec: proof_strategy}
We divide the description into several steps to give a panorama view. 
\vspace*{-0.4cm}
\subsubsection{$U_{t}$ as a Fourier integral operator} 
This pivotal idea was originally due to 
Duistermaat and Guillemin~\cite{Duistermaat_InventMath_1975} 
who worked it out (modulo smoothing operators) in the context of scalar half-wave operators on an ultrastatic (Remark~\ref{rem: SST}~\ref{rem: hypersurface_orthogonality_SST}) spacetime. 
We, however, have not followed their approach directly, instead have expressed $U_{t}$ in terms of the Killing flow $\varXi_{t}^{*}$ and the causal propagator $F := \eqref{eq: def_causal_propagator_Dirac}$ of $D$ by propounding 
Strohmaier and Zelditch's work on d'Alembertian~\cite{Strohmaier_AdvMath_2021}. 
More precisely, there exists unique advanced and retarded Green's operators (Section~\ref{sec: Green_op_Dirac}) for $D$ owing to the global hyperbolicity of $\spacetime$ and hence $F$ together with the restriction operator $\iota_{\ms \varSigma}^{*}$ (Example~\ref{exm: restriction_op_FIO}) pave the way to construct the Cauchy restriction operator $\cR$ as described in Section~\ref{sec: Cauchy_problem_Dirac}. 
By Assumption~\ref{asp: trace_formula}~\ref{asp: commutator_Dirac_op_Killing_flow}, the time flow determines the time evolution of Cauchy data. 
Then the combination of these facts allows us to express $U_{t}$ as a Lagrangian distribution as inscribed in Lemma~\ref{lem: symbol_Tr_U_0}. 

\vspace*{-0.3cm}
\begin{center}
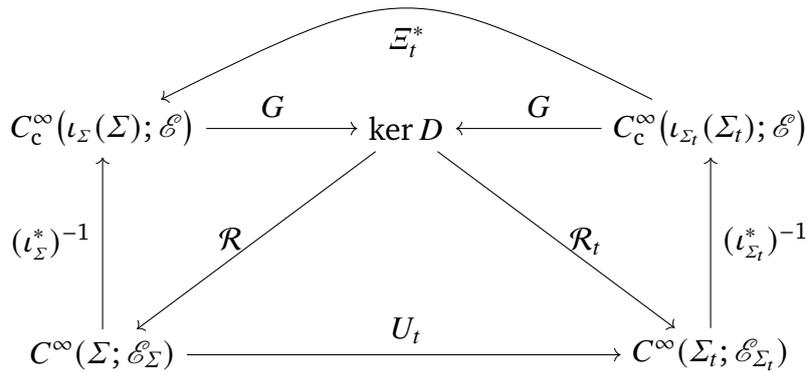

	\begin{tikzpicture}
		\node (a) at (0, 0) {$C^{\infty} (\varSigma; \sE_{\varSigma})$}; 
		\node (b) at (8, 0) {$C^{\infty} (\varSigma_{t}; \sE_{\varSigma_{t}})$}; 
		\node (c) at (0, 3) {$C_{\mathrm{c}}^{\infty} \big( \iota_{\ms \varSigma} (\varSigma); \sE \big)$}; 
        \node (d) at (8, 3) {$C_{\mathrm{c}}^{\infty} \big( \iota_{\ms \varSigma_{t}} (\varSigma_{t}); \sE \big)$}; 
		\node (e) at (4,3) {$\ker D$}; 
		\draw[->] (a) -- (b); 
		\draw[->] (a) -- (c); 
		\draw[->] (b) -- (d); 
		\draw[->] (e) -- (a); 
		\draw[->] (e) -- (b); 
		\draw[->] (c) -- (e); 
		\draw[->] (d) -- (e); 
		\draw[->] (d) .. controls (4,5) .. (c);
		\node[above] at (4,0) {$U_{t}$}; 
		\node[left] at (2, 1.5) {$\cR_{}$}; 
		\node[right] at (6, 1.5) {$\cR_{t}$}; 
        \node[left] at (0, 1.5) {$(\iota_{\ms \varSigma}^{*})^{-1}$}; 
        \node[right] at (8, 1.5) {$(\iota_{\ms \varSigma_{t}}^{*})^{-1}$}; 
		\node[above] at (5.75, 3) {$G$}; 
		\node[above] at (2.25, 3) {$G$}; 
		\node[below] at (4, 4.5) {$\varXi_{t}^{*}$};
	\end{tikzpicture}
	\captionof{figure}[Time evolution map in terms of the causal propagator and the Killing flow]{The 
		time evolution map $U_{t}$ in terms of the causal propagator $F$ for $D$ and the induced Killing flow $\varXi_{t}^{*}$. 
		Here, $\cR_{t}$ denotes the Cauchy restriction map and $\iota_{\ms \varSigma_{t}}^{*}$ the restriction map. 
	}
	\label{fig: time_evolution_causal_propagator_restriction_op_Killing_flow}
\end{center} 

In order to describe $U_{t}$ as a Fourier integral operator one requires to obtain such a description of $F$ and $\iota_{\ms \varSigma}^{*}$. 
Both are well-known for the scalar case and the bundle generalisation is straightforward for the latter. 
But the former demands an intricate treatment due to possible bundle curvature, which has been spelled out in Lemma~\ref{lem: causal_propagator_Dirac}. 
In fact, this is precisely from where the expedient choice (as discussed in Section~\ref{sec: Sandoval_trace_formula}) in Sandoval's work stems. 
We have computed $\symb{\fF}$ in Lemma~\ref{lem: causal_propagator_Dirac} deploying the $D^{2}$-compatible Weitzenb\"{o}ck connection in an intrinsically geometric fashion. 
As a consequence, it solves the relevant leading-order transport equation concisely from where the holonomy group shows up elegantly. 
Furthermore, it closes the ad-hoc consideration in Sandoval's analysis as explained after Theorem~\ref{thm: trace_formula_t_T}.  
\newline 

In contemporary of this work (after the submission of an arXiv-preprint based on this work), 
Capoferri and Murro~\cite[Thm. 1.1]{Capoferri} 
have obtained an oscillatory integral representation of $U_{t}$ (modulo smoothing operators) for the reduced massless Dirac equation (in the sense of Definition 3.5 of their paper) on a $4$-dimensional spatially compact globally hyperbolic spin-spacetime using the global phase function approach~\cite{Laptev_CPAM_1994} of Fourier integral operators. 
An antecedent of this analysis can be traced back to 
Capoferri and Vassiliev~\cite{Capoferri_2020} 
who have constructed $U_{t} (\check{D}) =~\eqref{eq: def_U_t_Riem_Dirac}$ (modulo smoothing operators) as a summation of two invariantly defined oscillatory integrals,
global in space and in time, with distinguished complex-valued phase functions, when $\check{D}$ is a massless spin-Dirac operator on a $3$-dimensional closed Riemannian manifold. 
\newline 

It is worthwhile to mention that one cannot deploy the algorithm used in~\cite{Capoferri, Capoferri_2020} in a straightforward way to derive the Fourier integral description (Lemma~\ref{lem: symbol_Tr_U_0}) of $U_{t}$ in our setting. 
This is primarily because the Dirac Hamiltonian
\begin{equation} \label{eq: def_Dirac_Hamiltonian_SSST}
	H_{D} := - \ri \frac{\fc (\rd t) \big( \fh_{ij} \fc (\rd x^{i}) \nabla^{{\ms \sE} j} + U \big) + \fh_{ij} \upalpha^{i} \nabla^{{\ms \sE} j}}{\upbeta^{2} - \| \upalpha \|_{\fh}^{2} - \fh_{ij} \upalpha^{i} \fc (\rd t) \fc (\rd x^{j})}
\end{equation}
on a stationary spacetime is \textit{not} of Dirac-type (albeit it can be written as a first-order elliptic operator) which is one of the key assumptions of the hindmost literatures. 
\vspace*{-0.2cm}
\subsubsection{Principal symbol of $\Tr U_{t}$} 
If $\sU_{t} (x, y) := \fU_{t} (x, y) \, \rd t \otimes \sqrt{|\dVolh (x)|} \otimes \sqrt{|\dVolh (y)|}$ denotes the Schwartz kernel of $U_{t}$ then one considers the smoothed-out operator $\fU_{\rho}$ (cf.~\eqref{eq: def_Tr_U_t_distribution}) 
\begin{equation} \label{eq: def_smoothed_time_evolution_op}
	\fU_{\rho} := \int_{\R} \fU_{t} \, \cF^{-1} (\rho) \, \rd t, 
	\quad 
	(\cF^{-1} \rho) (t) = \frac{1}{2 \pi} \int_{\spec \LieDeri} \re^{\ri t \lambda} \rho (\lambda) \, \rd \lambda
\end{equation}
for any Schwartz function $\rho \in \cS (\R)$ on $\R$ such that $\supp (\cF^{-1} \rho)$ is compact. 
Let $\fU_{\rho} (x, x)$ be the diagonal embedding of $\fU_{\rho} (x, y)$. 
The distributional trace $\Tr U_{t}$ is then obtained by 
\begin{equation} \label{eq: def_tr_smoothed_Cauchy_evolution_op}
	\Tr U_{\rho} := \int_{\varSigma} \tr \big( \fU_{\rho} (x, x) \big) \, \dVolh (x), 
\end{equation}
where $\tr$ is the endomorphism trace~\eqref{eq: def_endo_trace}. 
\newpage 

In order to compute $\symb{\Tr U_{t}}$ we employ the bundle generalisation of 
Duistermaat and Guillemin's~\cite{Duistermaat_InventMath_1975} 
idea, due to 
Sandoval~\cite{Sandoval_CPDE_1999}. 
One notices that the mapping $\sU_{t} (x, y) \mapsto \fU_{t} (x, x) \, \dVolh (x)$ can be viewed as the pullback 
\begin{equation}
    \varDelta^{*} : 
    C_{\rc}^{\infty} \Big( \R \times \varSigma \times \varSigma; \varOmega \R \boxtimes \big( \Hom{\sE_{\varSigma}, \sE_{\varSigma}} \otimes \halfDen (\varSigma \times \varSigma) \big) \Big)  
    \to 
    C_{\rc}^{\infty} \big( \R \times \varSigma; \varOmega \R \boxtimes (\End \sE_{\varSigma} \otimes \varOmega \varSigma) \big) 
\end{equation}
of $\sU_{t}$ via the diagonal embedding 
\begin{equation} \label{eq: def_diagonal_embedding_R_Cauchy}
	\varDelta : \R \times \varSigma \to \R \times \varSigma \times \varSigma  
\end{equation}
and the fibrewise isomorphism $\hat{\varDelta} \in C^{\infty} \big( \varSigma; \Iso (\sE_{\varSigma}, \sE_{\varSigma}) \big)$ induced by the Weitzenb\"{o}ck connection $\connectionE$ restricted to $\varSigma$. 
It follows (cf. Example~\ref{exm: pullback_Euclidean_FIO}) that $\varDelta^{*}$ is a Fourier integral operator of order $(d-1) / 4$ associated to the canonical relation 
$
    C_{\varDelta^{*}} 
    := 
    \{ (t, \tau; x, \xi + \eta; t, \tau; x, \xi; x, \eta) \in \dotCoTan \R \times \coTansM_{\varSigma} \times \dotCoTan \R \times \dotCoTansM_{\varSigma} \times \dotCoTansM_{\varSigma} \}
$~\cite[$(1.20)$]{Duistermaat_InventMath_1975}. 
Hence, for a fixed $t$, $\tr (\varDelta^{*} \sU_{t})$ is density on $\varSigma$  which can be integrated. 
The integration over $\varSigma$ is the pushforward 
(see e.g.~\cite[pp. 103-104]{Duistermaat_Birkhaeuser_2010})
\begin{equation}
    \uppi_{*} : 
    C_{\rc}^{\infty} (\R \times \varSigma; \varOmega \R \boxtimes \varOmega \varSigma) 
    \to 
    C^{\infty} (\R; \varOmega \R), 
    ~ \tr (\varDelta^{*} \sU_{t}) \mapsto \uppi_{*} \big( \tr (\varDelta^{*} \sU_{t}) \big) 
\end{equation}
of the Cartesian projection 
\begin{equation}
	\uppi : \R \times \varSigma \to \R,  
\end{equation}
which is also a Fourier integral operator of order $1/2 - (d-1)/4$ associated to the canonical relation 
$
C_{\uppi_{*}} 
    := 
    \{ (t, \tau; t, \tau; x, 0) \in \dotCoTan \R \times \dotCoTan \R \times \coTansM_{\varSigma} \}
$~\cite[$(1.22)$]{Duistermaat_InventMath_1975}. 
Therefore 
\begin{equation} \label{eq: Tr_U_t_pushforward_tr_pullback}
	\Tr U_{t} = \uppi_{*} \circ \tr (\varDelta^{*} \sU_{t}).   
\end{equation}
Alternatively, one can also integrate $\varDelta^{*} \sU_{t}$ over $\varSigma$ so that $\uppi_{*} (\varDelta^{*} \sU_{t})$ is an $\End \sE_{\varSigma}$-valued density on $\R$ by reckoning 
$
    \uppi_{*} : 
    C_{\rc}^{\infty} \big( \R \times \varSigma; \varOmega \R \boxtimes (\End \sE_{\varSigma} \otimes \varOmega \varSigma) \big)  
    \to 
    C^{\infty} (\R; \varOmega_{\R} \boxtimes \End \sE_{\varSigma})
$.   
Then taking the endomorphism trace one arrives at 
\begin{equation} \label{eq: Tr_U_t_tr_pushforward_pullback}
    \Tr U_{t} = \tr (\uppi_{*} \circ \varDelta^{*} \sU_{t}). 
\end{equation}
Hence, in this sense 
\begin{equation}
	\uppi_{*} \circ \tr \circ \varDelta^{*} = \tr \circ \uppi_{*} \circ \varDelta^{*}
\end{equation}
and one utilises the clean intersection (see Definition~\ref{def: clean_composition}) between $\uppi_{*} \circ \varDelta^{*}$ and $\sU_{t}$ to compute the principal symbol of $\Tr U_{t}$.  
%
%
%
%
%
%
%
%
%
%
\subsubsection{Spectral theory} 
Albeit finite energy solutions of Dirac equation do not live in $L^{2}$ sections of $\sE$, $\ker D$ can be naturally given a hermitian structure by equipping with a Killing flow invariant hermitian form $\scalarProdTwo{\cdot}{\cdot}$ owing to Assumption~\ref{asp: trace_formula}~\ref{asp: hermitian_form_Dirac_op_SST}.  
It is noteworthy that the non-definiteness of $(\cdot|\cdot)$ is a characteristic feature of any Lorentzian spin-manifold: a positive-definite $(\cdot|\cdot)$ invariant under the spin group only exists for the Riemannian case. 
%
%
%
%
%
%
%
%
%
%
\subsubsection{Weyl law} 
There are a number of approaches to deriving the Weyl law 
(see e.g. the review~\cite{Ivrii_BullMathSci_2016}).  
Amongst those, we will use the 
Fourier-Tauberian argument~\cite{Hoermander_ActaMath_1968}  
(see also, e.g.~\cite[App. B]{Safarov_AMS_1997},~\cite{Safarov_JFA_2001} and references therein). 
The key idea is to relate the Weyl counting function with $\Tr U_{t}$ via the distributional Fourier transform~\eqref{eq: Weyl_counting_function_Tr_U_t} $\rd \fN / \rd \lambda = \cF_{t \mapsto \lambda}^{-1} (\Tr U_{t})$ and compute the right-hand side using the express for $\Tr U_{0}$. 
\section{Dirac operator on a stationary spacetime}
\subsection{Lie derivative}
\label{sec: Lie_derivative}
On a stationary spacetime $\SST$, the cotangent lift $\coTan \varXi_{s}$ of the spacetime isometry $\varXi_{s} : \sM \to \sM$ naturally induces a $\R$ action on the cotangent bundle $\coTansM$: $\coTan \! \varXi (s; y, \eta) := (\coTan_{x} \varXi_{s}) \eta$. 
Furthermore, the pullback (restoring densities for clarity) 
\begin{equation}
	\varXi^{*} : \CComInfinity{\sM; \sE \otimes \halfDen \sM} \to C^{\infty} \big( \R \times \sM; \varOmega \R \boxtimes (\sE \otimes \halfDen \sM) \big)
\end{equation}
via the Killing flow $\varXi : \R \times \sM \to \sM$ and the fibrewise isomorphism $\hat{\varXi}_{s} \in C^{\infty} \big( \sM; \mathrm{Iso} (\sE, \sE) \big)$ as depicted in Figure~\ref{fig: bundle_lift_spacetime_isometry}, is a Fourier integral operator whose Schwartz kernel is given by (cf. Example~\ref{exm: pullback_Euclidean_FIO}) 
\begin{subequations} \label{eq: pullback_kernel}
	\begin{eqnarray}
		&& 
		\mathsf{\Xi} \in I^{-1/4} \Big( \R \times \sM \times \sM, \varGamma'; \varOmega \R \boxtimes \big( \Hom{\sE, \sE} \otimes \halfDen (\sM \times \sM) \big) \Big), 
		\label{eq: pullback_kernel_Lagrangian_dist}
		\\ 
		&& 
		\varGamma := \big\{ \big( s, \tau; x, \xi; y, - \eta \big) \in \dotCoTan \R \times \dotCoTansM \times \dotCoTansM \,|\, \tau = - \xi (Z), \xxi = (\coTan_{x} \varXi_{s}) \yeta \big\}, \qquad \quad 
		\label{eq: def_canonical_relation_pullback}
		\\ 
		&& 
        \symb{\mathsf{\Xi}} := (2 \pi)^{1/4} \one_{\Hom{\sE, \sE}} \sqrt{|\dVol_{\ms \varGamma}|} \otimes \bbl,  
		\label{eq: symbol_pullback}
	\end{eqnarray}
\end{subequations}
where $\dVol_{\ms \varGamma}$ is the volume form on the homogeneous canonical relation $\varGamma \subset \dotCoTan (\R \times \sM) \times \dotCoTansM$ and $\bbl$ is a section of the Keller-Maslov bundle $\bbL_{\varGamma} \to \varGamma$ over $\varGamma$, constructed as below. 
One observes that $\varGamma$ is the graph of $\coTan \varXi_{s}$ and at $s=0$, $\varGamma$ is essentially the conormal bundle $\varGamma_{0} := \{ (0, \tau) \} \times (\varDelta \, \dotCoTansM)'$. 
Since $\dotCoTansM$ is a symplectic manifold (see Example~\ref{exm: cotangent_bundle_symplectic_mf}), it admits volume form induced by the canonical symplectic form on $\dotCoTansM$. 
Then $\dVol_{\ms \varGamma_{0}}$ is obtained via the pullback of the projector $\Pr : \varGamma_{0} \to \R \times \dotCoTansM$, which is invariant under the flow $\varPsi_{s}$ of the Hamiltonian vector field generated by the extended Hamiltonian $\tau + \xi (Z)$ and given by $\dVol_{\ms \varGamma} = \rd s \otimes \rd x \wedge \rd \xi$ in the parametrisation~\eqref{eq: def_canonical_relation_pullback}. 
To construct $\bbL_{\varGamma}$, we re-use the fact that $\varPsi_{s}$ sweeps out $\{ (0, \tau) \} \times (\varDelta \, \dotCoTanM)'$ to $\varGamma$, and hence $\bbL_{\varGamma}$ is constructed by parallelly transporting the sections of $\bbL_{0}$ along the orbits of $\varPsi_{s}$, where $\bbL_{0}$ the Keller-Maslov bundle (see Example~\ref{exm: Maslov_bundle_conormal_bundle}) over $\varGamma_{0}$ consists of a global constant section.  
\newline 

More generally, $\varXi^{*}$ can be extended to a sequentially continuous linear map on $I^{m} (\sM, \sL; \sE)$ 
by
\begin{equation} \label{eq: pullback_Killing_flow_Lagrangian_dist}
	\varXi^{*} : I^{m} (\sM, \sL; \sE) \to I^{m - 1/4} (\R \times \sM, \varGamma' \circ \sL; \sE), 
	\qquad 
	\symb{\varXi^{*} u} \asymp \symb{\mathsf{\Xi}} \diamond \symb{u} 
\end{equation}
for any $u \in I^{m} (\sM, \sL; \sE)$. 
The composition $\diamond$ of principal symbols is presented in detail in 
Appendix~\ref{sec: composition_halfdensity_Maslov_vector_bundle_canonical_relation} and the equation of $\symb{\varXi^{*} u}$ is in the sense of modulo Keller-Maslov part.  
\newline

The induced Killing flow $\varXi_{s}^{*}$ paves the way to define the Lie derivative $\pounds_{\! Z}$~\eqref{eq: def_Lie_derivative} on $\sE$ with respect to the Killing vector field $Z$. 
Note, $\pounds_{\! Z}$ is essentially a generalisation of the 
Lichnerowicz spinor Lie derivative~\cite[Sec. 6]{Lichnerowicz_1963} 
for stationary spacetimes when one does not necessarily have a spin-structure. 
%
%
%
%
%
%
%
%
%
%
\subsection{Classical dynamics}
\label{sec: classical_dyanmics}
The primary tenet of the semiclassical analysis is to connect the relativistic trace formula with its classical dynamics. 
In non-relativistic mechanics, the cotangent bundle $\coTanCauchy$ models the classical phase space, whereas the Hilbert-space quantum dynamics takes place in $L^{2} (\varSigma)$. 
The naive expectation of using this pair or the pair $\big( \coTansM, L^{2} (\sM; \sE) \big)$ does \textit{not} work because the former depends on the choice of Cauchy hypersurface $\varSigma \subset \sM$ and for the latter pair, the (Killing flow invariant) sesquilinear form $(\cdot\vert\cdot)$ on $\sE$ (Assumption~\ref{asp: trace_formula}~\ref{asp: Direc_op_symmetric_SST}) does not induce any $L^{2}$-norm. 
One can, of course, choose an arbitrary hermitian form in order to have a $L^{2}$ space on $\sE$, but then this $L^{2}$-space will depend on the particular choice of the hermitian form as $\sM$ is non-compact. 
In pursuance of defining the correct classical dynamics, one notes $\Char D$ is the lightcone bundle $\coLightBun$. 
Since $\Char{D}$ and $\Char{\square}$ are identical, the classical dynamics in this case coincides with  that in the Strohmaier-Zelditch trace formula (Section~\ref{sec: Strohmaier_Zelditch_trace_formula}) and hence we adopt their formulae.   
The metric-Hamiltonian $H_{\fg}:= \eqref{eq: def_Hamiltonian_metric}$ is a homogeneous function of degree $2$ in the cotangent fibres. 
Referring to this as the dilation, let $\E$ be the Euler vector field which is the generator of this action. 
On $\coLightBun$, clearly $H_{\fg}$ vanishes and we have $[\E, X_{\fg/2}]_{-} = X_{\fg/2}$ where $X_{\fg/2}$ is the Hamiltonian vector field of $H_{\fg}$. 
Then the Hamiltonian reduction of $\coLightBun$ is the manifold of scaled-lightlike geodesics $\sN$. 
That is, if $\spacetime$ is geodesically complete then $\sN$ is the quotient of $\coLightBun$ by the $\R$-group action generated by $X_{\fg/2}$. 
Similarly, by taking the quotient of $\sN$ by the $\R_{+}$-group action generated by $\E$ we obtain the manifold of unparametrised-lightlike geodesics $\tilde{\sN}$. 
If $\spacetime$ is a spatially compact globally hyperbolic spacetime then $\tilde{\sN}$ is a \textit{conic compact contact manifold} whose symplectisation is the conic symplectic manifold $\sN$ induced from the conic contact manifold $\coLightBun$~\cite[pp. 10-12]{Penrose_1972},~\cite[Thm. 2.1]{Khesin_AdvMath_2009}. 
\newline 

We remark that $\sN$ is defined \textit{invariantly} and~\cite[Prop. 2.1]{Strohmaier_AdvMath_2021}
(see also~\cite{Low_JMP_1989, Low_NonlinearAnal_2001})
\begin{equation} \label{eq: symplectic_diffeo_lightlike_geodesic_coTanCauchy}
	\varsigma: \sN \to \dotCoTanCauchy 
\end{equation}
is a homogeneous symplectomorphism. 
Furthermore, the geodesic flow on any spatially compact globally hyperbolic spacetime $\spacetime$ does not necessarily have to be complete. 
For instance, $\sM := (-1, 1) \times \varSigma$ with the spacetime metric $\fg := \rd t^{2} - \fh$ where $(\varSigma, \fh)$ is a compact Riemannian manifold~\cite[Rem. 2.2]{Strohmaier_AdvMath_2021}.   
\newline 

As announced in Theorems~\ref{thm: trace_formula_t_zero} and~\ref{thm: trace_formula_t_T}, periodic geodesics $\gamma$ on $\sN$ play a vital role in our investigation. 
Employing the symplectomorphism~\eqref{eq: symplectic_diffeo_lightlike_geodesic_coTanCauchy}, these can be considered periodic geodesics on $\dotCoTanCauchy$ as well, which matches better with our intuitive expectation.  
However, $\dotCoTanCauchy$ \textit{depends} heavily on the \textit{choice of Cauchy hypersurface} $\varSigma$ and as discussed in Section~\ref{sec: stationary_spacetime}, there is no canonical choice for $\varSigma$ on a (globally hyperbolic) stationary spacetime $\SST$. 
Hence, $\sN$ offers the \textit{natural} and \textit{invariant} manifold to describe the classical dynamics in contrast to $\dotCoTanCauchy$.  
Recall that a geodesic $c : \R \to \R \times \varSigma$ in $\sM$ can be expressed as $c (s) = \big( t (s), c_{\ms \varSigma} (s) \big)$ where $c_{\ms \varSigma}$ is the base projection of $\varsigma \big( \gamma (s) \big)$. 
A periodic geodesic $c$ in $\sM$ with period $T$ then means  
\begin{equation}
	\forall s \in \R : c_{\ms \varSigma} (s + T) = c_{\ms \varSigma} (s).  
\end{equation}
We note that the set of $T$-periodic curves on $\sN$ is precisely the length spectrum of $(\varSigma, \fh)$ for an ultrastatic spacetime. 
For details, see, for 
instance~\cite{Masiello_NonlinearAnal_1992, Sanchez_NonlinearAnal_1999, Sanchez_AMS_1999, Bartolo_NonlinearAnal_2001}. 

%
%
%
\begin{center}
	\begin{tikzpicture}[line cap=round,line join=round,isometric view]
        \def\r{1}
        \def\h{3}
        \draw[-latex] (0,0,0) -- (-\r-1,0,0) node [left]  {$x$};
        \draw[-latex] (0,0,0) -- (0,-\r-1,0) node [right] {$y$};
        \draw[red,dashed]     (135:\r) arc (135:-45:\r);
        \draw[red] (0,0,\h) + (135:\r) arc (135:-45:\r);
        \helix{135}{315}{dashed}
        \helix{495}{540}{dashed}
        \draw[dashed] (0,0,0) -- (0,0,\h);
        \fill[gray!30, opacity=0.4] (0, 0, \h) circle (\r);
        \draw[gray,left color=gray!30,fill opacity=0.5]
        (135:\r) arc (135:315:\r) --++ (0,0,\h) arc (315:135:\r) -- cycle;
        \helix{0}  {135}{}
        \helix{315}{495}{}
        \draw[-latex] (0,0,\h) --++ (0,0,1) node [above] {$z$};
    \end{tikzpicture}
    \captionof{figure}{A periodic geodesic on ($\R \times \bbS, \rd t^{2} - \rd x^{2}$).}
\end{center}
%
%
%

%
%
%
%
%
%
%
%
%
%
\subsection{Time evolution operator} 
\label{sec: Cauchy_evolution_map}
On a spatially compact globally hyperbolic spacetime this mapping is defined by 
\begin{equation} \label{eq: def_Cauchy_evolution_op}
	U_{t', t} := \cR_{t} \circ (\cR_{t'})^{-1} : C^{\infty} \big( \varSigma_{t'}; \sE_{\varSigma_{t'}} \big) \to C^{\infty} \big( \varSigma_{t}; \sE_{\varSigma_{t}} \big),    
\end{equation} 
which is a homeomorphism and extends to a unitary operator (denoted by the same symbol) on the space of square integrable sections on Cauchy hypersurfaces, i.e.,  
\begin{equation}
	U_{t', t} : L^{2} \big( \varSigma_{t'}; \sE_{\varSigma_{t'}} \big) \to L^{2} (\varSigma_{t}; \sE_{\varSigma_{t}})
\end{equation}
is an isometry.  
%
%
%
%
%
%
%
%
%
%
\section{Trace formula} 
\label{sec: trace_Cauchy_evolution_op}
In this section, we will work on the set-up in Section~\ref{sec: result_Gutzwiller_trace}. 
That is, $\varSigma \subset \sM$ is an embedded submanifold and $\iota_{\ms \varSigma} : \varSigma \hookrightarrow \sM$ is proper (Remark~\ref{rem: Cauchy_hypersurface_embedded_submf}). 
As a consequence, $\iota_{\ms \varSigma}^{*} : \comSecsME \to C_{\mathrm{c}}^{\infty} (\varSigma; \sE_{\varSigma}) = \secECauchy$. 

%
%
%
\begin{theorem} \label{thm: trace_Cauchy_evolution_op}
	As in the terminologies of Theorem~\ref{thm: trace_formula_t_zero}, let $\fF$ be the Schwartz kernel of the causal propagator for $D$. 
	The smoothed-out time evolution operator $U_{\rho} :=$~\eqref{eq: def_smoothed_time_evolution_op} is a trace-class operator on the Hilbert space $\sH := (\ker D,~\eqref{eq: def_hermitian_form_Dirac_type_op})$ and its trace is given by 
	\begin{equation} \label{eq: trace_Cauchy_evolution_op}
		\Tr U_{\rho} 
		= 
        - \ri \int_{\varSigma} \tr \int_{\R} \big( \varXi_{-t}^{*} \circ \fF \circ (\iota_{\ms \varSigma}^{*})^{-1} \, \symb{D} (\cdot, \zeta ) \big) (x, y) \, (\cF^{-1} \rho) (t) \, \rd t |_{x = y} \, \dVolh (x),  
	\end{equation}
    where $\varXi_{t}^{*} : \comSecsME \to \comSecsME$ is the induced Killing flow (Figure~\ref{fig: bundle_lift_spacetime_isometry}) and $\iota_{\ms \varSigma}^{*} : \comSecsME \to \secECauchy$ is the restriction map. 
	In this thesis, $\varXi_{-t}^{*} \fF$ is meant to be the pullback of $\fF (\cdot, )$ in the first argument via $\varXi_{-t}$.   
\end{theorem}
%
%
%

\begin{proof}
	By Assumption~\ref{asp: trace_formula}, the hermitian form~\eqref{eq: def_hermitian_form_Dirac_type_op}, the retarded $\fundaSolRet$ and advanced $\fundaSolAdv$ propagators, all are preserved under the action of $\varXi_{t}^{*}$. 
	Hence $F$ is preserved as well.  
	In other words, if $u \in \ker D$ having initial data on some $\varSigma$ then $\varXi_{t}^{*} u \in \ker D$ having Cauchy data on some $\varSigma_{t}$, which means that the time flow $\varXi_{t}$ induces time evolution of Cauchy data.   
	Let us now choose an arbitrary but fixed $\varSigma$. 
	This picks a global time coordinate $t$ on $\sM$ and then the generator of the Killing flow $\varXi_{t}$ is given by $\partial_{t}$. 
	Thus, $U_{t}$ is identified with the evolution of Cauchy data via the induced Killing flow (see Figure~\ref{fig: time_evolution_causal_propagator_restriction_op_Killing_flow} for a schematic illustration): 
	\begin{equation} \label{eq: U_t_Cauchy_restriction_op_Killing_flow}
		U_{t} = \cR \circ \varXi_{-t}^{*} \circ \cR^{-1}.
	\end{equation} 
    We read off $\cR^{-1} = - \ri \, F \circ (\iota_{\ms \varSigma}^{*})^{-1} \symb{D} (\cdot, \zeta )$ from~\eqref{eq: smooth_sol_Dirac_op_Cauchy_data} and observe that the twisted wavefront set (Lemma~\ref{lem: causal_propagator_Dirac}) of $\fF$ contains only lightlike covectors. 
	Therefore, integration over $t$ results in a smooth Schwartz kernel $\fU_{\rho} (x, y)$ and the expression of $\Tr U_{\rho}$ entails from~\eqref{eq: def_tr_smoothed_Cauchy_evolution_op}. 
\end{proof}
%
%
%
%
%
%
%
%
%
%
\section{Spectral theory of $\rL$ on $\ker D$}
\label{sec: spectral_theory_L}
Recall that, a \textbf{strongly continuous one-parameter unitary group} is a family $\{ U_{t} | t \in \R \}$ of unitaries $U_{t}$ on a Hilbert space $\cH$ such that 
\begin{subequations}
	\begin{eqnarray}
		&& 
		\forall t, s \in \R : U_{t} U_{s} = U_{t+s}, 
		\\ 
		&& 
		\forall u \in \cH, \forall t \in \R : \lim_{h \to 0} U_{t+h} u = U_{t} u. 
	\end{eqnarray}
\end{subequations}
In particular, this implies that
\begin{equation}
	U_{-t} = U_{t}^{-1} = U_{t}^{*}.  
\end{equation}
By Stone's theorem 
(see e.g.~\cite[Thm. VIII.8]{Reed_I}), 
every $U_{t}$ has a unique generator $A$, i.e., $U_{t} = \re^{- \ri t A}$. 
If 
\begin{equation}
	U_{\rho} := \int_{\R} \re^{- \ri t A} (\cF^{-1} \rho) (t) \, \rd t 
\end{equation}
is a compact operator for any $\rho \in \cS (\R)$ such that $\supp (\cF^{-1} \rho)$ is compact, then the spectrum of $A$ is discrete and consists of eigenvalues $\lambda_{n}$ of finite algebraic multiplicities $m_{n}$. 
Moreover, whenever $U_{\rho}$ is trace-class then 
\begin{equation}
	\Tr U_{\rho} = \sum_{n} m_{n} \, \rho (\lambda_{n}).   
\end{equation}

%
%
%
\begin{proof}[Proof of Theorem~\ref{thm: spectrum_L}]
    Since $D$ commutes with induced Killing flow $\varXi_{t}^{*}$ for all $t \in \R$ and $U_{t}$ is the time evolution operator of $D$, the discreteness of $\spec \LieDeri$ and the polynomial growth of eigenvalues $\lambda_{n}$ follow from the fact that $\fU_{\rho} =$~\eqref{eq: def_smoothed_time_evolution_op} is a trace-class operator for $\rho \in \cS (\R)$ such that $\supp (\cF^{-1} \rho)$ is compact. 
	All $\lambda_{n}$ are real as a consequence of the selfadjointness of $\rL$ on $\ker D$.  
\end{proof}
%
%
%
%
%
%
%
%
%
%
\section{Proof of the main theorems}
\label{sec: proof_trace_formula_nondegenerate_fixed_pt}
In order to implement the strategy outlined in Section~\ref{sec: proof_strategy}, we recall that that
\begin{equation}
    \uppi_{*} \circ \varDelta^{*} : 
    C_{\mathrm{c}}^{\infty} \big( \R \times \varSigma \times \varSigma; \Hom{\sE_{\varSigma}, \sE_{\varSigma}} \big) 
    \to 
    C^{\infty} (\R; \End \sE_{\varSigma})
\end{equation}
is a zero-order Fourier integral operator whose Schwartz kernel $\fK$ is that of an identity map~\cite[Lem. 5.2]{Sandoval_CPDE_1999} 
(see also~\cite[Lem. 6.2, 6.3]{Duistermaat_InventMath_1975}):   
\begin{subequations}
	\begin{eqnarray}
		&& 
        \fK \in I^{0} \Big( \R \times \R \times \varSigma \times \varSigma, C_{\uppi_{*} \circ \varDelta^{*}}; \mathrm{Hom} \big( \Hom{\sE_{\varSigma}, \sE_{\varSigma}}, \End \sE_{\varSigma} \big) \Big), 
		\\ 
		&& 
		C_{\uppi_{*} \circ \varDelta^{*}} = \big( \varDelta (\R \times \varSigma)  \big)^{\perp *}, 
		\\ 
		&& 
		\symb{\fK} = \Pi^{*} \big( |\rd t \wedge \rd \tau \wedge \rd x \wedge \rd \xi|^{1/2} \big) \one. 
    \end{eqnarray}
\end{subequations}
Here the canonical relation $C_{\fK} = \{ (t, \tau; t, \tau; x, \xi; x, - \xi) \in \dotCoTan \R \times \dotCoTan \R \times \dotCoTanCauchy \times \dotCoTanCauchy \}$ of $\uppi_{*} \circ \varDelta^{*}$ has been identified with the cornormal bundle $\big( \varDelta (\R \times \varSigma)  \big)^{\perp *}$ to the diagonal in $\R \times \varSigma \times \R \times \varSigma$, $\Pi$ is the projector 
$
	C_{\pi_{*} \varDelta^{*}} \ni (t, \tau; x, \xi; x, - \xi; t, - \tau) \mapsto (t, \tau; x, \xi) \in \dotCoTan \R \times \dotCoTanCauchy 
$, 
and the Keller-Maslov bundle $\T \to \big( \varDelta (\R \times \varSigma)  \big)^{\perp *}$ is trivial (see Example~\ref{exm: Maslov_bundle_conormal_bundle}). 
\newline 

Theorem~\ref{thm: trace_Cauchy_evolution_op} implies that $\Tr U_{t}$ exists as a distribution in $\cD' (\R)$ and it can be re-expressed as  
\begin{equation} \label{eq: tr_U_t}
	\Tr U_{t} = \int_{\varSigma} \tr \big( (\varXi_{-t}^{*} \fF) (x, y) \, \symb{D} (y, \zeta ) \big) \big|_{x = y} \, \dVolh (x). 
\end{equation}
Since the geodesic relation (Definition~\ref{def: geodesic_relation}) is disjoint with the conormal bundle $\varSigma^{\perp *}$, restriction of $\varXi_{-t}^{*} \fF$ is well-defined. 
Then applying~\eqref{eq: Tr_U_t_tr_pushforward_pullback}, the preceding equation can be written as 
\begin{equation} \label{eq: Tr_U_t_tr_pushforward_pullback_restriction_causal_propagator}
	\Tr U_{t} 
	=  
	\tr \Big( \uppi_{*} \circ \varDelta^{*} \circ (\iota_{x}^{*} \boxtimes \iota_{y}^{*}) \big( (\varXi_{-t}^{*} \fF) (x, y) \, \symb{D} (y, \zeta) \big) \Big) 
\end{equation}
for a fixed but arbitrary $t \in \R$ and any $x, y \in \sM$. 
Hence, our task is to compute the principal symbol of 
\begin{equation} \label{eq: def_G_t}
	\fF_{t} := (\iota_{x}^{*} \boxtimes \iota_{y}^{*}) \big( (\varXi_{-t}^{*} \fF) (x, y) \big). 
\end{equation} 
We have worked out $\iota_{\ms \varSigma}^{*}$ in Example~\ref{exm: restriction_op_FIO}, so let us begin by describing $\varXi_{-t}^{*} \fF$ as a Lagrangian distribution.  
To begin with, one notes that $w$ appearing in the expression of $\symb{\fF}$ satisfies~\cite[Thm. 6.6.1]{Duistermaat_ActaMath_1972} 
\begin{equation}
    (X_{\fg^{\pm}/2} \pm \ri \subSymb{D^{2}, \pm}) w = 0, 
\end{equation}
by Lemma~\ref{lem: causal_propagator_Dirac} and Definition~\ref{def: P_compatible_connection}, where $\fg^{\pm}$ resp. $\subSymb{\square, \pm}$ are the lifts of $\fg$ resp. $\subSymb{\square}$ to $\dotCoTansM \times \dotCoTansM$ via the projections on the first resp. second copies of $\dotCoTansM$. 
This means that $\ri (2 \pi)^{3/4} \symb{D} \circ w \sqrt{|\rd t|} \otimes \sqrt{|\dVol_{\ms C}|} / 2$ is the principal symbol of $\varXi^{*} \fF$ on each $C^{+}$ and $C^{-}$. 
Employing~\eqref{eq: pullback_kernel},~\eqref{eq: pullback_Killing_flow_Lagrangian_dist}, Lemma~\ref{lem: causal_propagator_Dirac},   
and Appendix~\ref{sec: composition_halfdensity_Maslov_vector_bundle_canonical_relation}, it is then  straightforward to obtain 

%
%
%
\begin{lemma} \label{lem: pullback_Killing_flow_causal_propagator}
	Suppose that $(\sE \to \sM, \symb{D})$ is a bundle of Clifford modules over a globally hyperbolic stationary spacetime $\SST$ and that $\fF$ is the Schwartz kernel of the causal propagator of any Dirac-type operator $D$ on $\sE$.   
	Then, the pullback (with respect to the first argument) of $\fF$ by the Killing flow is a Lagrangian distribution      
	\begin{subequations} \label{eq: pullback_Killing_flow_causal_propagator_kernel}
		\begin{eqnarray}
			&& 
			\varXi^{*} \fF   
			\in 
			I^{-3/4} \big( \R \times \sM \times \sM, \varGamma' \circ C; \Hom{\sE, \sE} \big), 
			\label{eq: pullback_Killing_flow_causal_propagator_Lagrangian_dist}
			\\ 
			&& 
			\varGamma' \circ C  
			= 
			\big\{ \big( t, - \xi (Z); x, \xi; y, \eta \big) \in \dotCoTan \R \times \coLightBun \times \coLightBun \,|\, \xxi = \coTan \! \varXi_{t} \circ \varPhi_{s} \yeta \big\}, \hspace*{1.25cm} 
			\label{eq: def_canonical_relation_pullback_Killing_flow_causal_propagator}
			\\ 
			&& 
			\symb{\varXi^{*} \fF}  
			\asymp 
            \ri (2 \pi)^{\frac{3}{4}} \symb{D} \circ w \big( \coTan \! \varXi_{t} \circ \varPhi_{s} \yeta, \yeta \big) \, |\rd t|^{\frac{1}{2}} \otimes \big| \dVol_{\ms C} \big( \coTan \! \varXi_{t} \circ \varPhi_{s} \yeta, \yeta \big) \big|^{\frac{1}{2}}, \qquad \qquad  
			\label{eq: symbol_pullback_Killing_flow_causal_propagator}
		\end{eqnarray}
	\end{subequations}
	where the principal symbol is modulo Keller-Maslov part, $Z$ is the infinitesimal generator of $\varXi$, and all other symbols are as defined in Lemma~\ref{lem: causal_propagator_Dirac} and~\eqref{eq: def_canonical_relation_pullback}. 
\end{lemma}
%
%
%

Next, we compute the restriction of $\varXi^{*} \fF$ on $\varSigma_{x} \times \varSigma_{y}$ by an application of Example~\ref{exm: restriction_op_FIO}.  

%
%
%
\begin{lemma} \label{lem: restriction_pullback_causal_propagator}
	As in the terminologies of Lemma~\ref{lem: pullback_Killing_flow_causal_propagator}, let $\varSigma$ be a Cauchy hypersurface of $\sM$. 
	Then, the distribution~\eqref{eq: def_G_t} is a Lagrangian distribution 
	\begin{subequations}
		\begin{eqnarray}
			\fF_{t} 
            & \in &   
			I^{-1/4} \big( \R \times \varSigma \times \varSigma, \cC_{t}; \Hom{ \sE_{\varSigma}, \sE_{\varSigma}} \big), 
			\\ 
			\cC_{t} 
            & := &  
            \big\{ \big( (t, - \xi (Z)), \xxi |_{\tangent \varSigma}, \yeta |_{\tangent \varSigma} \big) \in \dotCoTan \R \times \coTansM_{\varSigma} \times \coTansM_{\varSigma} \,|\, 
            \nonumber \\ 
            && 
            \xxi = \coTan \! \varXi_{t} \circ \varPhi_{s} \yeta \big\}, 
			\label{eq: def_canonical_relation_Cauchy_evolution_op}
			\\ 
			\symb{\fF_{t}} 
            & \asymp &   
			\ri (2 \pi)^{\frac{1}{4}}     
            \symb{D} \circ w \big( \coTan \! \varXi_{t} \circ \varPhi_{s} \yeta |_{\tangent \varSigma}, \yeta |_{\tangent \varSigma} \big) 
			\nonumber \\ 
			&& 
			|\rd t|^{\frac{1}{2}} \otimes \big| \dVol_{\ms C} \big( \coTan \! \varXi_{t} \circ \varPhi_{s} \yeta |_{\tangent \varSigma}, \yeta |_{\tangent \varSigma} \big) \big|^{\nicefrac{1}{2}}, 
			\label{eq: symbol_restriction_pullback_causal_propagator}
		\end{eqnarray}
	\end{subequations}
	where the principal symbol is modulo the Keller-Maslov part.    
\end{lemma}
%
%
%

\begin{proof}
	The first two assertions are immediate 
	from~\eqref{eq: restriction_kernel_Lagrangian_dist},~\eqref{eq: def_canonical_relation_restriction},~\eqref{eq: pullback_Killing_flow_causal_propagator_Lagrangian_dist} 
	and~\eqref{eq: def_canonical_relation_pullback_Killing_flow_causal_propagator}. 
	To compute the principal symbol one notes that $\dVol_{\ms \coLightBun} \yeta$ induces a volume element $\dVol_{\ms \coLightBun_{\varSigma}} \yeta := \dVol_{\ms \coLightBun} / \rd y^{1}$ on $\coLightBun_{\varSigma}$ when $\varSigma$ is parametrised by $y^{1} = \cst$ and then $\dVol_{\ms \coLightBun_{\varSigma}} \yeta = \rd \eta_{1} \wedge \rd y' \wedge \rd \eta'$ in the adapted coordinates $(y^{1} = \cst, y', \eta_{1}, \eta')$ on $\coTansM$. 
	Finally, the claim follows from~\eqref{eq: symbol_restriction},~\eqref{eq: symbol_pullback_Killing_flow_causal_propagator},~\eqref{eq: symbol_restriction}, and transporting $\yeta$ to $\xxi$ by the geodesic flow and the Killing flow. 
\end{proof}
%
%
%

If the composition $C'_{\uppi_{*} \varDelta^{*}} \circ \cC_{t}$ is clean then we can compute $\Tr U_{t}$ by the standard composition of Fourier integral operators. 
But, in general, 
\begin{equation} 
	C'_{\uppi_{*} \varDelta^{*}} \circ \cC_{t} 
	= 
	\big\{ (t, \tau) \in \R \times \dot{\R} \,|\, \tau = - \xi (Z), \xxi |_{\tangent \varSigma} = \big( \coTan \varXi_{-t} \circ \varPhi_{t} \yeta \big) |_{\tangent \varSigma} \big\}  
\end{equation}
may \textit{not} be clean 
(see Appendix~\ref{sec: clean_composition}:~\eqref{eq: def_C_star_varLambda},~\eqref{eq: def_fibre_projection_clean_composition},  and~\eqref{eq: varPi} for the symbol ${\scriptsize \text{\FiveStarOpen}}$ and for further details) 
because the fibres of $C'_{\uppi_{*} \varDelta^{*}} {\scriptsize \text{\FiveStarOpen}} \cC_{t} \to C'_{\pi_{*} \varDelta^{*}} \circ \cC_{t}$ can be identified with the fibres over $\tau \in \dotCoTan_{t} \R$, i.e., the set $\{\kF_{t}\} :=$~\eqref{eq: def_F_T} of periodic geodesics, where we have used~\eqref{eq: def_Hamiltonian_scaled_lightlike_geodesic} and~\eqref{eq: symplectic_diffeo_lightlike_geodesic_coTanCauchy}.   
Then, even if $\{ \kF_{t} \}$ happens to be manifolds, the chances of $\dim \kF_{t}$ is a constant for all $t$ is very low; for instance, if all the orbits of $\varXi_{t}^{\sN}$ are periodic with the same period $t = T$, then the composition is clean  
(see e.g.~\cite[p. 289]{Meinrenken_ReptMathPhys_1992}). 
We remark that this, however, is not an issue for the trivial period ($T=0$) and the assumptions made on classical dynamics in Theorem~\ref{thm: trace_formula_t_T} ensures a clean intersection in the case of non-trivial periods ($T \neq 0$). 
\newline 

Let us record for the future computations that  
\begin{eqnarray}
	&& 
  	\dim \sM = d, \dim (\coTansM) = 2d, \dim (\coLightBun) = 2d-1, 
	\nonumber \\ 
	&& 
	\dim \sN = 2d-2, \dim \tilde{\sN} = 2d-3, 
  	\nonumber \\ 
  	&& 
  	\dim \varSigma = d-1, \dim (\coTanCauchy) = 2d-2 = \dim (\coLightBun_{\varSigma}). 
\end{eqnarray} 
%
%
%
%
%
%
%
%
%
%
\subsection{Principal symbol of $\Tr U_{t}$ at $t = 0$}
\label{sec: symbol_t_zero}
We begin with the trivial periodic orbits where a big singularity is expected as $U_{t}$ reduces to an identity operator in this situation. 
To describe $\symb{\Tr U_{0}}$, it is useful to have the notion of the 
symplectic residue~\cite[Def. 6.1]{Guillemin_AdvMath_1985},  
introduced by Guillemin to derive Weyl's law in the context of Weyl algebra quantising a conic symplectic manifold. 
By construction (cf.~\eqref{eq: def_Hamiltonian_scaled_lightlike_geodesic}), $H^{1-d}$ is a homogeneous function of degree $1-d$. 
Then its symplectic residue is defined by 
\begin{equation}
	\res H^{1-d} := \int_{\tilde{\sN}} \vartheta_{H^{1-d}}, 
	\quad 
	\vartheta := \mathbb{E} \lrcorner \dVol_{\ms \sN}, 
\end{equation}
where $\lrcorner$ denotes the interior multiplication by the Euler vector field (Section~\ref{sec: classical_dyanmics}) $\mathbb{E}$ on $\sN$ and $\vartheta_{H^{1-d}}$ is the pull-back of $H^{1-d} \vartheta$ on $\sN$. 
Homogeneity of $H$ entails 
(cf.~\cite[Proof of Lemma 6.3]{Guillemin_AdvMath_1985}) 
\begin{equation} \label{eq: symplectic_residue_volume_lightlike_geodesic}
	\res H^{1-d} = (d-1) \, \vol \sN_{H \leq 1}. 
\end{equation}
On a globally hyperbolic standard stationary spacetime, $\vol \sN_{H \leq 1}$ has been computed by  
Strohmaier-Zelditch~\cite[$(15)$]{Strohmaier_AdvMath_2021} 
and it is given by~\eqref{eq: vol_energy_surface}.  

%
%
%
\begin{lemma} \label{lem: symbol_Tr_U_0}
	As in the terminologies of Theorem~\ref{thm: trace_Cauchy_evolution_op} (and hence Theorem~\ref{thm: trace_formula_t_zero} as well), 
	\begin{equation}
		\fU_{t} \in I^{-1/4} \big( \R \times \varSigma \times \varSigma, \cC_{t}; \Hom{\sE_{\varSigma}, \sE_{\varSigma}} \big) 
	\end{equation}
	is a Lagrangian distribution associated with the canonical relation $\cC_{t} :=$~\eqref{eq: def_canonical_relation_Cauchy_evolution_op}, whose principal symbol is   
    \begin{equation}
    	\symb{\fU_{t}} \xxiyeta \asymp \symb{\fF_{t}} \xxiyeta \, \symb{D} \yeta 
    \end{equation}
    modulo the Keller-Maslov part, where $\symb{\fF_{t}}$ is given by~\eqref{eq: symbol_restriction_pullback_causal_propagator}. 
    \newline 

    Furthermore, $\Tr U_{t}$ is a Lagrangian distribution on $\R$ of order $d - 7/4$ associated with the Lagrangian submanifold $\Lambda_{T}$ and its principal symbol at $T = 0$ is given by  
	\begin{equation}
		\symb{\Tr U_{0}} (\tau) \asymp r \frac{d-1}{(2 \pi)^{d-1}} \vol (\sN_{H \leq 1}) \, |\tau|^{d-2} \sqrt{|\rd \tau|}. 
	\end{equation}
\end{lemma}
%
%
%
\begin{proof}
	The first assertion simply follows from $\fU_{t} (x, y) = \fF_{t} (x, y) \, \symb{D} (y, \zeta)$.
    We then use identification of $\sN$ with $\coLightBun_{\varSigma}$ so that $|\dVol_{\ms \sN}|$ is identified with $|\dVol_{\ms \coLightBun_{\varSigma}}|$. 
	If $\tE$ is a regular value of $H$ then $\sN_{\tE}$ is a codimension one $X_{H} |_{\sN_{\tE}}$ invariant embedded submanifold of $\sN$, which inherits a natural volume form $\dVol_{\ms \sN_{\tE}}$ from $\dVol_{\ms \sN}$, invariant under the action of $X_{H} |_{\sN_{\tE}}$~\cite[Thm. 3.4.12]{Abraham_AMS_1978}: 
	\begin{equation}
		\dVol_{\! \ms \sN} = \dVol_{\! \ms \sN_{\tE}} \wedge \frac{\rd H}{\| \grad H \|} 
		\Leftrightarrow 
		\dVol_{\! \ms \sN_{\tE}} (\ldots) = \dVol_{\! \ms \sN} \bigg( \frac{\grad H}{\| \grad H \|}, \ldots \bigg). 
	\end{equation}

    One observes that at $T = 0$, 
    $\coTan \varXi_{0} \big( \varPhi_{0} \yeta \big) = \yeta$ 
    and 
    $C = \varDelta \, \coLightBun$ 
    so that 
    $\cC'_{0} = \big\{ \big( 0, - \xi (Z) \big) \big\} \times (\varDelta \, \coLightBun) |_{\tangent \varSigma \times \tangent \varSigma}$ 
    and 
    $\big| \dVol_{\ms C} (y, \eta; y, \eta) |_{\tangent \varSigma \times \tangent \varSigma} \big|^{1/2} = |\dVol_{\ms \sN} \yeta|$. 
    The excess of $C'_{\uppi_{*} \varDelta^{*}} \circ \cC_{0}$ is (see Appendix~\eqref{eq: excess_fibre_projection_clean_composition}) 
	\begin{equation}
		e = \dim \kF_{0} = 2d - 3. 
	\end{equation}
    Hence, by the composition of Lagrangian distributions, 
    $\Tr U_{0} \in I^{d - 7/4} (\R, \Lambda_{T})$. 
    One observes~\cite[$(6.6)$]{Duistermaat_InventMath_1975} 
    that $\varphi (t, \tau)$ defined by $(t - T) \tau$ (resp. $0$) for $\tau < 0$ (resp. $\tau \geq 0$), is a phase function for $\Lambda_{T}$. 
    Obviously, $\varphi$ is not smooth around $\tau = 0$ but the required modification will only reflect by some smooth terms in the oscillatory integral described below.
	For brevity, these have been suppressed notationally. 
    Any element of this Lagrangian distribution can be written as scalar multiples $\rc_{1}, \rc_{2}, \ldots$ of oscillatory integrals of the form (cf.~\eqref{eq: def_naive_kernel_FIO_Euclidean})
    \begin{equation}
        (\Tr U_{0}) (t) = (2 \pi)^{-3/4} (2 \pi)^{d-1} \int_{\R_{\geq 0}} \re^{-\ri t \tau} ( \rc_{1} \tau^{d-2} + \rc_{2} \tau^{d-3} + \ldots) \rd \tau, 
	\end{equation}
    because $\symb{\Tr U_{0}}$ must be of order $d-2$ (cf. Definition~\ref{def: symbol_FIO_Euclidean}) as entailed by the order ($d-7/4$) of $\Tr U_{0}$. 
    \newline 

    Computation of $\symb{\Tr U_{0}}$ involves the following four steps as detailed in Section~\ref{sec: composition_halfdensity_Maslov_vector_bundle_canonical_relation}.  
    Briefly speaking, first one obtains the tensor product $\symb{\fK} \otimes \symb{\fU_{0}}$ followed by intersecting with the diagonal. 
    Then we integrate over $\kF_{0}$ and take $\tr$ of the hindmost quantity. 
    Putting all these together and comparing with the last expression of $\Tr U_{0}$ we obtain  
    \begin{equation}
        \symb{\Tr U_{0}} (0, \tau) = r \frac{\res (H^{1-d})}{(2 \pi)^{d-1}} |\tau|^{d-2} \sqrt{|\rd \tau|}, 
    \end{equation}
    where we have used Lemmas~\ref{lem: restriction_pullback_causal_propagator} and~\ref{lem: causal_propagator_Dirac} together with the identity  
	\begin{equation}
		\tr \big( \symb{D} \yeta \, \symb{D} (y, \zeta) \big) 
		= 
		\tr \frac{\symb{D} \yeta \, \symb{D} (y, \zeta) + \symb{D} (y, \zeta) \, \symb{D} \yeta}{2} 
		= 
        \fg_{y}^{-1} (\eta, \zeta) \, \rk (\sE)
	\end{equation} 
	using the cyclicity of trace and~\eqref{eq: anticommutation_symbol_Dirac_type}. 
\end{proof}
%
%
%
%
%
%
%
%
%
%
\subsection{Principal symbol of $\Tr U_{t}$ at $t = T \in \cP$} 
\label{sec: symbol_t_T}
This is the scenario corresponding to the non-trivial periodic orbits when $C_{\uppi_{*} \varDelta^{*}} \circ \cC_{\ms T}$ is clean or equivalently the set of fixed points $\kF_{T}$ is clean.  
Recall, the fixed point sets of $\varXi_{T}^{\sN}$ are the union of periodic orbits $\gamma$ of $\varXi_{T}^{\sN}$ which depends substantially on manifold $\spacetime$ and we are only concerned about non-degenerate orbits as briefly defined below 
(details are available  
in~\cite[Sec. 4]{Duistermaat_InventMath_1975}).   
Let $\upsigma$ be the symplectic form on $\coTansM$ (see Example~\ref{exm: cotangent_bundle_symplectic_mf}). 
The induced Killing flow $\varXi_{s}^{\sN}$ preserves the level sets $\{ H (c) = \tau \}$ since $\varXi_{s}^{\sN}$ is the Hamiltonian $H := \eqref{eq: def_Hamiltonian_scaled_lightlike_geodesic}$ flow on $\sN$. 
Then, its restriction to $\cV := \tangent_{c} \{ H = \tE \}$ for any $\tE$ satisfies $\upsigma |_{\cV} (X_{H}, \cdot) = 0$ and by the 
\textit{non-degeneracy}~\cite[p. 44]{Strohmaier_AdvMath_2021} 
of $\gamma$, it is meant that the \textit{nullspace of $\upsigma |_{\cV}$ is spanned by} $X_{H}$. 
This condition implies that 
$\ker (I - \rd_{\gamma} \varXi_{s}^{\cV}) = \R X_{H}$ 
and 
$\mathrm{img} (I - \rd_{\gamma} \varXi_{s}^{\sV}) = \cV / \R X_{H}$ 
where $\rd_{\gamma} \varXi_{s}^{\sN} : \cV \to \cV$ and $\ker (I - \rd_{\gamma} \varXi_{s}^{\sN})$ is to be understood as the kernel of $I - \rd_{\gamma} \varXi_{s}^{\sN}$ on $\cV$. 
In this setting, the linearised Poincar\'{e} map $\rP_{\gamma}$ is defined as   
\begin{equation}
    I - \rP_{\gamma} : \cV / \ker (I - \rd_{\gamma} \varXi_{s}^{\sN}) \to \cV / \ker (I - \rd_{\gamma} \varXi_{s}^{\sN}) 
\end{equation} 
the linear symplectic quotient map induced by $I - \rd_{\gamma} \varXi_{s}^{\sN}$. 
%
%
%
\begin{lemma} \label{lem: symbol_Tr_U_T}
	As in the terminologies of Theorem~\ref{thm: trace_formula_t_T}, the principal symbol of $\Tr U_{t}$ at $t = T \in \cP$ is   
	\begin{equation}
		\symb{\Tr U_{T}} (\tau) =   
        \frac{\tau^{d-2}}{(2 \pi)^{d-1}} \int_{\kF_{T}} 
        \tr \Big( \symb{D} (\gamma) \, \cT_{\varXi_{\ms T}^{\ms \sN} (\gamma)} \symb{D} (x, \zeta) \Big) 
        \frac{\re^{- \ri \pi \km (\gamma) / 2} |\rd T| \otimes \sqrt{|\rd \tau|}}{\sqrt{|\det (I - \rP_{\gamma})|}}. 
	\end{equation}
\end{lemma}
%
%
%
\begin{proof}
	We begin with the fact that non-degenerate $\gamma$ belongs to a $2$-dimensional cylinder, transversally intersecting the energy hypersurfaces~\cite[p. 576]{Abraham_AMS_1978} which implies that one can use $\tau = - H (\gamma)$ as a coordinate for fibre over $(T, \tau)$. 
	Then the half-density valued density at $(T, \tau)$ is given by~\cite[pp. 60 - 61]{Duistermaat_InventMath_1975}  
	(see also, e.g.~\cite[Thm. 6.1.1]{Guillemin_InternationalP_2013})
	\begin{equation}
		|\rd T| \otimes \sqrt{\frac{|\rd \tau|}{|\det (I - \rP_{\gamma})|}}, 
	\end{equation}
	where $|\rd T|$ is the density on $\gamma$ induced by the flow.  
	\newline 

    To compute $\symb{\Tr \fU_{T}}$, we use identifications $\sN = \coLightBun_{\varSigma}, |\dVol_{\ms \sN}| = |\dVol_{\ms \coLightBun_{\varSigma}}|$ once again. 
    For $T \neq 0$, $C$ is no more given by $\varDelta \, \coLightBun$ and so we read off the off-diagonal expression of $\symb{\fF}$ from Lemma~\ref{lem: causal_propagator_Dirac} to obtain   
	\begin{equation}
		\symb{\fU_{T}} \xxiyeta 
		\asymp 
		(2 \pi)^{\frac{1}{4}} \ri     
        \symb{D} \xxi \, w \xxiyeta  
		|\rd t|^{\frac{1}{2}} \otimes \big| \dVol_{\ms C} \xxiyeta \big|^{\frac{1}{2}} 
	\end{equation}
    modulo the Keller-Maslov contribution, where $\xxi = \coTan \! \varXi_{t} \circ \varPhi_{s} (y^{\backprime}, \eta^{\backprime}) |_{\tangent \varSigma}, \yeta = (y^{\backprime}, \eta^{\backprime}) |_{\tangent \varSigma}$ for $(y^{\backprime}, \eta^{\backprime}) \in \coLightBun$. 
    As before, $\symb{\Tr U_{T}}$ is computed employing~\eqref{eq: Tr_U_t_tr_pushforward_pullback} by performing the four steps used in Lemma~\ref{lem: symbol_Tr_U_0}. 
	These yield the claimed expression 
    \begin{equation}
        \symb{\Tr U_{T}} (T, \tau) 
        \asymp 
        \frac{\tau^{d-2}}{(2 \pi)^{d-1}} \int_{\kF_{T}} \tr \Big( \symb{D} (\gamma) \, \cT_{\varXi_{\ms T}^{\ms \sN} (\gamma)} \symb{D} (x, \zeta) \Big) \frac{|\rd T| \otimes \sqrt{|\rd \tau|}}{\sqrt{| \det (I - \rP_{\gamma}) |}} 
    \end{equation}
	modulo the contribution coming from the Keller-Maslov line bundle $\Maslov \to \Lambda_{T}$. 
	\newpage  

	Thus, we are left with the computation of the $\Maslov$-part which has been calculated in~\cite[$(6.16)$]{Duistermaat_InventMath_1975} 
    (see also~\cite[Sec. 3]{Meinrenken_JGP_1994}). 
	For completeness, we briefly outline the main steps.  
	First, one constructs $\bbL \to \cC_{T}$ from $\bbL_{C}, \bbL_{\varGamma}$, and $\bbL_{\varLambda}$ following the procedure detailed in Section~\ref{sec: Keller_Maslov_bundle}. 
    Let $\varphi$ resp. $\phi$ be generating functions of $C_{\uppi_{*} \varDelta^{*}}$ resp. $\cC_{T}$. 
    Then $\psi := \varphi + \phi$ locally generate $\Lambda_{T}$. 
    Suppose that we partition $[0, 1]$ as $0 = s_{0} < s_{1} < \ldots < s_{\ms N} = 1$ such that $\{ \cL_{\alpha} \}_{k=0, \ldots, N}$ is a conic covering of $\Lambda_{T}$ around $\gamma (s_{\alpha})$ and that $\{ \psi_{\alpha} \}$ are corresponding generating functions. 
    Then, following the steps explained in Section~\ref{sec: Maslov_index} one obtains that $\re^{-\ri \pi \km (\gamma) / 2}$ is the Keller-Maslov contribution. 
\end{proof}
%
%
%
%
%
%
%
%
%
%
\section{Proof of Weyl law}  
As shown in Theorems~\ref{thm: trace_formula_t_zero} and~\ref{thm: trace_formula_t_T} that $\Tr U_{t}$ has singularities at $t=0, T$. 
Our aim is to only cut the $t = 0$ singularity. 
In order to do so, one introduces a Schwartz function $\chi (\lambda)$ on $\R$ such that 
(see e.g.~\cite[p. 133]{Grigis_CUP_1994},~\cite[Sec. 21.1-3]{Shubin_Springer_2001},~\cite[Thm. 7.5]{Wunsch_Zuerich_2008}) 
\begin{enumerate}[label=(\roman*)] 
	\item 
    $\chi (\lambda) > 0$ for all $\lambda$; 
    \item 
    $(\Fourier \chi) (0) = 1$; 
    \item 
    $(\Fourier \chi) (t) = (\Fourier \chi) (-t)$; 
    \item 
    $\supp (\Fourier_{\lambda \mapsto t} \chi) \subset (- \varepsilon, \varepsilon)$, where $\varepsilon > 0$ is sufficiently small.  
\end{enumerate}
Then employing expression of $\Tr U_{0}$ from Theorem~\ref{thm: trace_formula_t_zero} and the aforementioned properties of $\chi$, we obtain 
\begin{eqnarray}
    \Fourier_{t \mapsto \lambda}^{-1} \big( \Fourier_{\cdot \mapsto t} (\chi) \Tr U_{t} \big) 
	& \approx & 
	\int \rd \tau \, \delta (\tau - \lambda) r \frac{\res (H^{1-d})}{(2 \pi)^{d-1}} \tau^{d-2} + \ldots 
	\nonumber \\ 
	& = & 
	r (d-1) \frac{\vol (\sN_{H \leq 1})}{(2 \pi)^{d-1}} \lambda^{d-2} + \ldots, 
\end{eqnarray}
where $f \approx g$ means that $\nicefrac{f}{g} \to 1$ as $t \to 0$  for any functions $f (t), g(t)$. 
That means that 
\begin{equation}
    (\chi * \rd \fN) (\lambda) = r (d-1) \frac{\vol (\sN_{H \leq 1})}{(2 \pi)^{d-1}} \lambda^{d-2} + O (\lambda^{d-3}), 
\end{equation}
where $*$ denotes the convolution. 
Employing 
\begin{equation}
	(\chi * \fN) (\lambda) 
    =
    \int_{- \infty}^{\lambda} (\chi * \rd \fN) (\mu) \, \rd \mu 
    =  
    r \frac{\vol (\sN_{H \leq 1})}{(2 \pi)^{d-1}} \lambda^{d-1} + O (\lambda^{d-2}). 
\end{equation}
The Weyl law entails from 
\begin{equation}
	\fN (\lambda) = (\fN * \chi) (\lambda) - \int \big( \fN (\lambda - \mu) - \fN (\lambda) \big) \, \chi (\mu) \, \rd \mu
\end{equation}
together with the facts $\fN (\lambda - \mu) - \fN (\lambda) \lessapprox \langle \mu \rangle^{d-1} \langle \lambda \rangle^{d-2}$ and $\int \chi (\mu) \, \rd \mu = 1$. 
%
%
%
%
%
%
%
%
%
%
\section{Literature} 
\label{sec: literature}
The idea of the asymptotic trace formula was originally due to   
Martin \textsc{Gutzwiller}~\cite{Gutzwiller_JMP_1971} 
in the context of quantum mechanics, 
which was given a rigour mathematical underpinning by 
Johannes J. \textsc{Duistermaat} and Victor W. \textsc{Guillemin}~\cite{Duistermaat_InventMath_1975}.   
We refer to the references mentioned before Section~\ref{sec: Duistermaat_Guillemin_trace_formula} for important chronological steps in between and, for instance, the survey article~\cite{Verdiere_AIF_2007} 
for a historical detour. 
\newline 

Jens \textsc{Bolte} and Stefan \textsc{Keppeler}~\cite{Bolte_PRL_1998, Bolte_AnnPhys_1999} 
have generalised Gutzwiller's work for spin Dirac operators on Minkowski spacetime 
(see also the elucidating articles~\cite{Bolte_FoundPhys_2001, Bolte_JPA_2004} 
and references therein for chronological developments). 
In contemporary to their reports, mathematical rigour analysis has been reported by 
Mary R. \textsc{Sandoval}~\cite{Sandoval_CPDE_1999} 
promoting the Duistermaat-Guillemin framework for Dirac-type operators on a closed Riemannian manifold. 
Her result on Dirac-wave-trace invariants for trivial periods is closely related to the evaluation of the residues of the eta-invariant (by 
Branson and Gilkey~\cite{Branson_JFA_1992})  
and the behaviour of eigenfunctions in the high energy limit (by 
Jakobson and Strohmaier~\cite{Jakobson_CMP_2007}).  
In the particular case of a massless spin-Dirac operator on a closed $3$-Riemannian manifold,  
Capoferri and Vassiliev~\cite{Capoferri_2020}
have computed the third local Weyl coefficients by advancing the framework by  
Chervova \textit{et al}.~\cite{Chervova_JST_2013} 
(see also the review by 
Avetisyan \textit{et al}.~\cite{Avetisyan_JST_2016})
for asymptotic spectral analysis of a general elliptic first-order system;   
see~\cite[Sec. 11]{Chervova_JST_2013} and~\cite[Chap. 3]{Fang_PhD} 
for a bibliographic overview.
An in-depth investigation of different spectral coefficients of a selfadjoint Laplace-type operator on a hermitian vector bundle over a closed Riemannian manifold has been performed by 
Li and Strohmaier~\cite{Li_JGP_2016}.   
In particular, they have obtained the relevant coefficients for Dirac-type operators with a more general bundle endomorphism than in  
Sandoval~\cite{Sandoval_CPDE_1999} (identity endomorphism)
and in 
Branson-Gilkey~\cite{Branson_JFA_1992} (scalar endomorphism),  
and bridged the results by 
Chervova \textit{et al}.~\cite{Chervova_JST_2013} 
(and their follow up works as mentioned in~\cite{Li_JGP_2016}) 
with the known heat-trace invariants. 
\newline 

The primary and common ingredient of~\cite{Bolte_PRL_1998, Bolte_AnnPhys_1999, Sandoval_CPDE_1999, Capoferri_2020} 
is to determine the time evolution operator (modulo smoothing operators) by solving the transport equations order by order.   
Bolte-Keppeler have also identified terms responsible for spin-magnetic and spin-orbit interactions in the semiclassical expression. 
However, they have finally considered regularised truncated time evolution operator by introducing an energy localisation to deal with the continuous spectrum of Dirac Hamiltonian arising due to the non-compact Minkowski spacetime. 
Such restrictions were absent in the hindmost references as they considered closed Riemannian manifolds and utilised the full power of Fourier integral operator theory in contrast to Bolte-Keppeler who have worked with oscillatory integrals. 
A novel feature in the lines of research by 
Vassiliev and his collaborators~\cite{Chervova_JST_2013, Avetisyan_JST_2016} 
is the second term of the Weyl law (see the references cited in these papers for earlier works). 
Analogous results have been achieved by  
Li-Strohmaier~\cite{Li_JGP_2016} 
employing a different spectral analysis.  
Specifically, the global phase function approach of Fourier integral operators has been deployed in~\cite{Capoferri_2020} 
to construct (modulo smoothing operators) the solution operator of a massless spin-Dirac operator on a $3$-dimensional closed Riemannian manifold. 
They have also provided a closed formula for the principal symbol and an algorithm for computing the subprincipal symbol of this operator.
\newline 

A general relativistic generalisation of Duistermaat-Guillemin-Gutzwiller trace formula has been initiated by 
Strohmaier and Zelditch~\cite{Strohmaier_AdvMath_2021, Strohmaier_IndagMath_2021}  
(see also the review~\cite{Strohmaier_RMP_2021}) 
who have studied the d'Alembertian on a globally hyperbolic spatially compact stationary spacetime as briefed in Section~\ref{sec: Strohmaier_Zelditch_trace_formula}. 
Their crucial step was to set up a relativistic description of the classical and the quantum dynamics and advance the celebrated 
Duistermaat-Guillemin~\cite{Duistermaat_InventMath_1975} 
framework accordingly. 
In particular, they have expressed the time evolution operator by means of the causal propagator of d'Alembertian and Killing flow, and computed its principal symbol by utilising the symbolic calculus of Fourier integral operators based on 
Duistermaat and H\"{o}rmander's~\cite{Duistermaat_ActaMath_1972} 
classic work of distinguished parametrices.
Subsequently, they employed 
Guillemin's symplectic residue~\cite{Guillemin_AdvMath_1985} 
approach at trivial period and tailored the Duistermaat-Guillemin computation for the non-trivial periods. 
Consequently, the Weyl law in the space of lightlike geodesics has been reported by them using the standard Fourier-Tauberian argument.      
Apart from this investigation, their work has been generalised in a bundle setting for a d'Alembertian on a globally hyperbolic stationary Kaluza-Klein spacetime by 
McCormick~\cite{McCormick} 
who utilised some technical results of this thesis. 
%
%
%

%% file: composition_density_canonical_relation.tex
\chapter{Canonical Relations \& Compositions}
\label{ch: canonical_relation_composition}
\textsf{In this chapter we will first briefly recall some standard notions of conic symplectic geometry and then review canonical relations, their clean intersection, and the composition of bundle-valued densities on the them}.
%
%
%
%
%
%
%
%
%
%
\section{Canonical relations}
\subsection{Preliminaries on conic symplectic geometry}
\label{sec: preliminary_conic_symplectic_geometry}
We recall that a \textbf{symplectic manifold} is a pair $(\fM, \sigmaup)$ where $\fM$ is a smooth manifold and $\sigmaup$ is a smooth, non-degenerate, closed $2$-form. 
Symplectic manifolds must be even dimensional and are always \textit{orientable} because they admit the \textbf{Liouville form} $\dVol_{\ms \fM}$ defined in~\eqref{eq: def_Liouville_form}.  
If $(\fM, \sigmaup_{\ms \dM}), (\fN, \sigmaup_{\ms \fN})$ are two symplectic manifolds and $\varkappa : \fN \to \fM$ is a smooth map (resp. diffeomorphism --- in this case $\dim \fM = \dim \fN$) with 
\begin{equation} \label{eq: def_symplecto}
	\varkappa^{*} \sigmaup_{\ms \fM} = \sigmaup_{\ms \fN}
\end{equation}
then $\varkappa$ is called a \textbf{symplectic map} (resp. \textbf{symplectomorphism}). 
The latter is known as a \textit{canonical transformation} in Physics literature. 
If $\dim \fM = \dim \fN$ then any symplectic map $\varkappa$ is a volume preserving $\varkappa^{*} \fv_{\ms \fM} = \fv_{\ms \fN}$ local diffeomorphism. 
Stated differently, the total volume of a symplectic manifold is a global symplectic invariant. 
By the Darboux theorem, \textit{every $2d$-dimensional symplectic manifold $(\fM, \sigmaup)$ is locally symplectomorphic to $(\R^{2d}, \sigmaup_{0})$} where $\sigmaup_{0}$ is the standard symplectic form on $\R^{2d}$. 
This implies that 
\begin{equation}
	\sigmaup \! \upharpoonright \! U = \rd x^{i} \wedge \rd \xi_{i}, i = 1, \ldots, d, 
    \quad 
    \fv_{\ms M} \! \upharpoonright \! U = \rd x^{1} \wedge \ldots \wedge \rd x^{d} \wedge \rd \xi_{1} \wedge \ldots \wedge \rd \xi_{d} = \rd x \rd \xi
\end{equation}
on any symplectic chart $\big( U, (x^{i}, \xi_{i}) \big)$ for $(\fM, \sigmaup)$. 
\newline 

We recall $\sigmaup$ induces a vertical vector bundle isomorphism $\sigmaup^{\flat} : \tangent \fM \to \coTan \fM$ whose inverse is denoted by $\sigmaup^{\sharp} : \coTan \fM \to \tangent \fM$. 
Let $f \in C^{\infty} (\fM, \R)$. 
Then the vector field 
\begin{equation} \label{eq: def_HVF}
	X_{f} := (\rd f)^{\sharp}
\end{equation}
is called the \textbf{Hamiltonian vector field} generated by $f$, which is locally given by 
\begin{equation}
	X_{f} \xxi = \frac{\partial f}{\partial \xi_{i}} \xxi \frac{\partial}{\partial x^{i}} \xxi - \frac{\partial f}{\partial x^{i}} \xxi \frac{\partial}{\partial \xi_{i}} \xxi. 
\end{equation}
The \textbf{Poisson bracket} is defined by
\begin{equation} \label{eq: def_Poisson_bracket}
	\{ f, g \} := \sigmaup (X_{f}, X_{g}) = X_{f} (g)  
\end{equation}
for any $f, g \in C^{\infty} (\fM, \R)$. 
Locally it reads  
\begin{equation}
    \{ f, g \} \xxi = \parDeri{\xi_{i}}{f} \xxi \parDeri{x^{i}}{g} \xxi - \parDeri{x^{i}}{f} \xxi \parDeri{\xi_{i}}{g} \xxi. 
\end{equation}

%
%
%
\begin{definition} \label{def: conic_mf} 
    A $d$-dimensional manifold $(\cM, \fm_{\lambdaup})$ is called \textbf{conic} if there is a given smooth map $\R_{+} \times \cM \ni (\lambdaup, x) \mapsto \fm_{\lambdaup} (x) := \fm (\lambdaup, x) \in \cM$ 
    and for every $x \in \cM$, there is an open neighbourhood $U$ with $\fm_{\ms \R_{+}} (U) = U$ and a diffeomorphism $\kappa : U \to \U$ such that $\lambdaup \kappa = \kappa \, \fm_{\lambdaup}$, where $\U$ is an open cone in $\dotRd$. 
    \newline 

    If a conic manifold $(\fM, \fm_{\lambdaup})$ is also symplectic $(\fM, \sigmaup)$ with the property 
    \begin{equation} \label{eq: group_action_symplecto}
    	\fm_{\lambdaup}^{*} \sigmaup = \lambdaup \sigmaup 
    \end{equation}
    then $(\fM, \fm, \sigmaup)$ is called a \textbf{conic symplectic manifold} 
    (see e.g.~\cite[Definition 21.1.8]{Hoermander_Springer_2007}).   
\end{definition}
%
%
%

Since $\fm_{\lambdaup} \circ \fm_{\lambdaup'} = \fm_{\lambdaup \lambdaup'}$ holds true in $\dotRd$, succinctly speaking, a conic manifold $(\cM, \fm)$ is a manifold endowed with a \textit{free, proper}, and \textit{smooth} $\R_{+}$-\textit{multiplicative group action} $\fm$ on $\cM$ 
(see e.g.~\cite[Def. 2.1.1]{Duistermaat_Birkhaeuser_2011}).  
If $f \in C^{\infty} (\cM)$ then 
\begin{equation} \label{eq: def_radial_VF}
    \sX (f) := \frac{\rd (\fm_{\lambdaup}^{*} f)}{\rd \lambdaup} (1) 
\end{equation}
is called the \textbf{radial vector field}. 
Locally this means $\sX = \xi_{i} \nicefrac{\partial}{\partial \xi_{i}}$. 

%
%
%
\begin{example} \label{exm: cotangent_bundle_symplectic_mf}
	If $M$ is a $d$-dimensional manifold then its cotangent bundle $\coTanM$ can naturally be seen as a symplectic manifold because it carries a \textit{canonical $1$-form}, also known as the \textbf{Liouville $1$-form} $\thetaup$ whose negative differential is \textbf{canonical symplectic form} $\sigmaup := - \rd \thetaup$ on $\coTanM$. 
    In symplectic coordinates $(x^{i}, \xi_{i})$ on $\coTanM$, 
    \begin{equation}
    	\thetaup = \xi_{i} \, \rd x^{i}, \quad \sigmaup = \rd x^{i} \wedge \rd \xi_{i}. 
    \end{equation}
    Whilst $(\coTanM, \sigmaup)$ is a symplectic manifold, the punctured cotangent bundle $(\dotCoTanM, \sigmaup)$ is naturally a conic symplectic manifold where 
    \begin{equation} \label{eq: punctured_coTanM_conic_mf}
        \fm_{\lambdaup} \xxi := (x, \lambdaup \, \xi), 
        \quad 
        \kappa \xxi := \big( (|\xi_{1}| + \ldots + |\xi_{d}|) (x^{1}, \ldots, x^{d}), (\xi_{1}, \ldots, \xi_{d}) \big). 
    \end{equation}
\end{example}
%
%
%

\begin{example} \label{exm: punctured_real_vector_bundle_conic_mf}
	Let $\sV \to M$ be a $\R$-vector bundle of rank $n$ over a manifold $M$. 
    Then $(\dot{\sV}, \fm_{\lambdaup})$ is a conic manifold with respect to the multiplication $\fm_{\lambdaup}$ defined as all dilations in its fibres: 
    \begin{equation} \label{eq: def_dilation_fibre}
        \fm_{\lambdaup} : \dot{\sV}_{x} \to \dot{\sV}_{x}, ~ v (x) \mapsto \fm_{\lambdaup} \big( v (x) \big) := (\lambdaup v) (x) 
    \end{equation}
    for any $\lambdaup \in \R_{+}$. 
    Locally $\sV \cong U \times \Rn$ and $\dot{\sV} \cong U \times \dotRn$ for some open set $U \subset M$ 
    (see e.g.~\cite[p. 27]{Duistermaat_Birkhaeuser_2011}). 
    The bundle $\dot{\sV}$ being a conic then locally translates  
    \begin{equation} \label{eq: dilation_fibre_local}
        \xtheta \in U \times \dotRn \Rightarrow \fm_{\lambdaup} \xtheta = (x; \lambdaup \theta) \in U \times \dotRn.
    \end{equation}
    A set $\cU := \{ \xtheta \in U \times \Rn \}$ is called a \textbf{conic subset} of $U \times \Rn$ if it is stable under all dilation~\eqref{eq: dilation_fibre_local}: $\xtheta \mapsto (x; \lambdaup \theta)$.
    Thereby, a positively homogeneous function $a$ of degree $s$ on $\cU$ means (cf. Definition~\ref{def: conic_mf})  
    \begin{equation} \label{eq: def_positive_homogeneous_function}
        a (x; \lambdaup \theta) = \lambdaup^{s} a \xtheta.  
    \end{equation}
\end{example}

%
%
%
\begin{theorem}[Homogeneous Darboux Theorem] \label{thm: homogeneous_Darboux}
    Let $(\fM, \fm_{\lambdaup}, \sigmaup)$ be a $2d$-dimensional conic symplectic manifold. 
    Suppose that $I$ and $J$ are subsets of $\{ 1, \ldots, d \}$ and that $(x^{i})_{i \in I}$ and $(\xi_{j})_{j \in J}$ are smooth functions on a conic neighbourhood of any $\xxi \in \fM$ such that~\cite[Prop. 6.1.3]{Duistermaat_ActaMath_1972}  
    (see also~\cite[Thm. 21.1.9]{Hoermander_Springer_2009}) 
    \begin{enumerate}[label=(\alph*)]
    	\item \label{con: homogeneous_Darboux_I}
        $x^{i}$ and $\xi_{j}$ are homogeneous of degree $0$ and $1$,  respectively, 
        \item \label{con: homogeneous_Darboux_II}
        $x^{i}$ Poisson commutes with $x^{j}$ for all $i, j \in I$ and so does $\xi_{i}$ with $\xi_{j}$ for all $i, j \in J$ whilst $\{ \xi_{j}, x^{i} \} = \delta_{j}^{i}$ for all $i \in I$ and for all $j \in J$, 
        \item \label{con: homogeneous_Darboux_III} 
        the Hamiltonian vector fields $X_{x^{i}}, X_{\xi_{j}}$ and the radial vector field $\sX$ are linearly independent at $\xxi$, 
        \item \label{con: homogeneous_Darboux_IV} 
        there are arbitrary real numbers $\bx^{i}$ and $\bxi_{j}$ so that $x^{i} \xxi = \bx^{i}$ and $\xi_{j} \xxi = \bxi_{j}$ with the assumption that $\bxi_{k} \neq 0$ for some $k \notin I$.  
    \end{enumerate}
    Then there exists smooth functions $(x^{i'})_{i'  \notin I}$ and $(\xi_{j'})_{j' \notin J}$ in a conic neighbourhood of $\xxi$ such that~\ref{con: homogeneous_Darboux_I},~\ref{con: homogeneous_Darboux_II},~\ref{con: homogeneous_Darboux_IV} hold true for indices running from $1$ to $d$ and~\ref{con: homogeneous_Darboux_III} remains valid if $i' \neq k$. 
    In particular, $(x^{i}; \xi_{i})$ defines a 
    homogeneous\footnote{Cf.~\eqref{eq: group_action_symplecto}. 
    } 
    symplectomorphism 
    \begin{equation}
        \varkappa : \cU \to \varkappa (\cU), \qquad \fm_{\lambdaup}^{*} \varkappa = \lambdaup \, \varkappa
    \end{equation}
    of degree $1$ from an open conic neighbourhood $\cU$ of $\xxi$ in $\fM$ to an open conic neighbourhood $\varkappa (\cU)$ of $(\bx, \bxi) \in \dotCoTan \R^{d}$.  
\end{theorem}
%
%
%
%
%
%
%
%
%
%
\subsection{Lagrangian submanifolds}
\label{sec: Lagrangian_submf}
%
%
%
\begin{definition} \label{def: Lagrangian_submf}
	Let $(\fM, \sigmaup)$ be a symplectic manifold.
    An immersion $\iota : \varLambda \hookrightarrow \fM$ is called \textbf{Lagrangian} if and only if $\iota^{*} \sigmaup = 0$ and $\dim \varLambda = \dim \fM / 2$.  
\end{definition}
%
%
%

\begin{example} \label{exm: Lagrangian_submf_img_exact_form}
	Let $M$ be a manifold. 
    Then the zero covector and the image $\xiup (M)$ of an exact $1$-form $\xiup$ on $M$ are Lagrangian submanifolds of $\coTanM$. 
\end{example}

%
%
%
\begin{example} \label{exm: conormal_bundle_Lagrangian_submf}
	Let $(\varSigma, \iota)$ be a submanifold of a manifold $M$. 
    The quotient vector bundle $\varSigma^{\perp} := \nicefrac{\tangent M_{\varSigma}}{\tangent \varSigma}$ is called the \textbf{normal bundle} whose algebraic dual 
    \begin{equation}
        \varSigma^{\perp *} := \nicefrac{\coTanM_{\varSigma}}{\coTanCauchy}, \quad \coTanM_{\varSigma} := \coTanM  \upharpoonright \varSigma
    \end{equation}
    is called the \textbf{conormal bundle} of $(\varSigma, \iota)$.  
    It is a $\R$-vector bundle over $\varSigma$ of $\rk \varSigma^{\perp *} = \codim \varSigma$ so $\dim \varSigma^{\perp *} = \dim M$. 
    Moreover, it is naturally isomorphic to the annihilator of $\rd \iota$ in $\tangent M_{\varSigma}$,  which is a Lagrangian submanifold of $\coTanM$  
    (see e.g.~\cite[p. 374]{Rudolph_Springer_2013}). 
    \newline 

    If we assume that $(\varSigma, \iota)$ is embedded then its conormal bundle is characterised as 
    \begin{equation}
        \varSigma^{\perp *} := \big\{ \yeta \in \coTanM \,|\, y \in \varSigma, \eta \upharpoonright \tangent_{y} M_{\varSigma} = 0 \big\}. 
	\end{equation}
    Note, $\dot{\varSigma}^{\perp *} \subset \dotCoTanM$ is a conic Lagrangian manifold 
    (see e.g.~\cite[Prop. 3.7.2]{Duistermaat_Birkhaeuser_2011}). 
\end{example}
%
%
%

For any exact $1$-form $\xiup$ on a manifold $M$, we can write $\xiup = \rd \varphi$ for some $\varphi \in C^{\infty} (M, \R)$.  
The function $\varphi$ is an example of the so-called generating function for the Lagrangian submanifold $\xiup (M) = \rd \varphi \, (M)$. 
This is an archetypal result for any Lagrangian immersion of a generic symplectic manifold because 
$\sigmaup$ vanishes identically on $\varLambda$ and so the Poincar\'{e} lemma implies that locally $\sigmaup = \rd \vartheta$ for some potential $\vartheta$ on $M$. 
Applying this lemma once again yields an open subset $\cL \subset \varLambda$ and a $\varphi \in C^{\infty} (\cL, \R)$ depending on $\vartheta$, such that $(\iota^{*} \vartheta) \upharpoonright \cL = \rd \varphi$. 
If we consider the particular class of Lagrangian immersions $\varLambda$ of the cotangent bundle $\pi : \coTanM \to M$ then the set $\Sigma (\varLambda)$ of all critical points of
\begin{equation} \label{eq: def_projection_Lagrangian_M}
	\Pi : \pi \circ \iota : \varLambda \to M, ~\xxi \mapsto \Pi \xxi := x
\end{equation}
is called the \textbf{singular set} and its projection on $M$ is called the \textbf{caustic} of $\varLambda$. 
This set comprises the points where $\varLambda$ is not transversal to cotangent fibration and in its complement $\varLambda$ can be always locally represented as the graph of an exact $1$-form $\rd \varphi$ defined on some open subset of $M$.  
We will see below that this construction can be carried over all $\varLambda$ by introducing some auxiliary parameter $\theta \in \Rn$ 
(see e.g.~\cite{Meinrenken_ReptMathPhys_1992},~\cite[Rem. 8.8.2]{Rudolph_Springer_2013}).  
\newline 

The primary objective of introducing this specific class of submanifolds is to complement the theory of Fourier integral operators (Sections~\ref{sec: FIO_Euclidean},~\ref{sec: Lagrangian_distribution_mf}, and~\ref{sec: Lagrangian_distribution}) where we only require the \textit{conic} Lagrangian manifolds $\varLambda$. 
We remark that such Lagrangian manifolds are precisely those where the canonical $1$-form $\thetaup$ vanishes identically and vice-versa~\cite[$(3.1.3)$]{Hoermander_ActaMath_1971} 
(see also, e.g.~\cite[Prop. 3.7.1]{Duistermaat_Birkhaeuser_2011}). 
\newline 

We are now going to systematically describe the parametrisation of a conic Lagrangian manifold $\varLambda \subset \dotCoTanM$ by means of a generating function, also known as a phase function. 

%
%
%
\begin{definition} \label{def: clean_phase_function}
	Let $M$ be a $d$-dimensional manifold and $\cU \subset M \times \dotRn$ an open conic set where $n$ is not necessarily equal to $d$.  
    Any $\R$-valued smooth function $\varphi \in C^{\infty} (\cU, \R)$ on $\cU$ is called a \textbf{clean phase function} with \textbf{excess} $0 \leq e \leq n$ if~\cite[Sec. 7]{Duistermaat_InventMath_1975}  
    (see also, e.g.~\cite[Def. 21.2.15]{Hoermander_Springer_2007},~\cite[Def. 1]{Meinrenken_ReptMathPhys_1992}) 
    \begin{enumerate} [label=(\alph*)]
    	\item 
        it is positive-homogeneous of degree $1$ in the $\dotRn$-variable $\theta$, 
        \item 
        its differential $\rd \varphi \neq 0$ is non-vanishing, 
        \item 
        its \textbf{fibre-critical set} 
        \begin{equation} \label{eq: def_fibre_critical_set}
            \sC := (\grad_{\theta} \varphi)^{-1} (0) := \Big\{ \xthetaNot \in \cU \,\Big|\, \Big( \parDeri{\theta_{k}}{\varphi} \Big)_{k=1, \ldots, n} \xthetaNot = 0 \Big\}
        \end{equation}
        is a $d+e$-dimensional smooth conic submanifold of $M \times \dotRn$, whose tangent space is 
        \begin{equation}
            \tangent_{\xtheta} \sC := \ker \bigg( \rd \Big( \parDeri{\theta_{k}}{\varphi} \Big)_{k=1, \ldots, n} \xtheta \bigg).  
        \end{equation}
    \end{enumerate}
    If $e = 0$ then $\varphi$ is called a \textbf{non-degenerate phase function}. 
\end{definition}
%
%
%

Given a clean $\varphi$, the differentials $\rd (\partial \varphi / \partial \theta_{k})$ at any point of $\sC$ span a linear space of dimension $n - e$ and the map~\cite[Lem. 7.1]{Duistermaat_InventMath_1975} 
\begin{equation} \label{eq: def_Lagrangian_submf_fibration}
    \jmath_{\varphi} : \sC \to \varLambda_{\varphi} \subset \dotCoTanM, ~\xtheta \mapsto 
    \jmath \xtheta := (x, \rd_{x} \varphi) 
\end{equation}
is locally a \textit{homogeneous fibration} of fibre dimension $e$ and its image is an \textit{immersed conic Lagrangian submanifold} $\varLambda_{\varphi}$. 
If $e = 0$ then $\big( \rd (\partial \varphi / \partial \theta_{1}), \ldots, \rd (\partial \varphi / \partial \theta_{n}) \big)$ are linearly independent at $\xtheta \in \sC$ and $\jmath_{\varphi} :=$~\eqref{eq: def_Lagrangian_submf_fibration} is a \textit{homogeneous immersion}~\cite[p. 134]{Hoermander_ActaMath_1971} 
(see also, e.g.~\cite{Meinrenken_ReptMathPhys_1992}). 
\newline 

We remark that $\varLambda \subset \dotCoTanM$ is a conic Lagrangian manifold if and only if, there exists a conic neighbourhood $\cL$ for every $\xxi \in \varLambda$ such that 
\begin{equation}
    \varLambda \cap \cL = \varLambda_{\varphi} \Leftrightarrow \{ \xxi \} = \{ (x, \rd_{x} \varphi) \} 
\end{equation}
for some non-degenerate phase function $\varphi$ on $\cL$
(see e.g.~\cite[Prop. 3.7.3]{Duistermaat_Birkhaeuser_2011},~\cite[Thm., p. 26]{Weinstein_AMS_1977}).  
In this situation, we say that $\varphi$ \textbf{locally generates} $\varLambda$ 
(see e.g.~\cite[p. 668]{Rudolph_Springer_2013}). 
Such a local generating function always exists and has a canonical form, i.e., 
\textit{every (conic) Lagrangian submanifold is locally generated by a non-degenerate phase function}. 
For any $\xxi \in \varLambda$, local coordinates $(x^{i})$ at $x \in M$ can be chosen such that with the corresponding homogeneous symplectic coordinates $(x^{i}, \xi_{i})$ on $\dotCoTanM$, the Lagrangian plane $\xi_{i} = \cst$ through $\xxi$ is transversal to $\varLambda$. 
Furthermore, there exists a non-degenerate phase function $\varphi$ in an open conic subset of $M \times \dotRd$ having the unique form~\cite[Thm. 3.1.3]{Hoermander_ActaMath_1971} 
(see also, e.g.~\cite[Thm. 21.2.16, 21.2.18]{Hoermander_Springer_2007},~\cite[Prop. 1.4, Cor. 1.1, Rem 1.1 $($pp. 418-419$)$]{Treves_Plenum_1980})
\begin{equation} \label{eq: def_canonical_nondegenerate_phase_function_Lagrangian_submf}
    \varphi = x^{i} \xi_{i} - H 
    \quad \textrm{so that} \quad 
    \varLambda_{\varphi} = \varLambda_{H} = \Big\{ \parDeri{\xi_{i}}{H}, \xi_{i} \Big\}, 
\end{equation}
where $H$ is a smooth positively homogeneous $\R$-valued function of degree $1$ on an open conic neighbourhood of $\xi \in \dotCoTan_{x} M$.  
We remark that the homogeneous immersion~\eqref{eq: def_Lagrangian_submf_fibration} is the identity map in this case. 
\newline 

The fascinating existence result immediately provokes 
\begin{itemize}
	\item 
    Are local phase functions unique?
    \item 
    Is it possible to achieve a global parametrisation of a Lagrange manifold utilising a phase function?
    \item 
    How the number $n$ of fibre variables are chosen?
\end{itemize}
The first question rises the notion of an equivalence class of phase functions. 
Technically speaking, let $\tilde{\cU} \subset M \times \dot{\R}^{\tilde{n}}$ be an open conic set such that there exists a diffeomorphism
\begin{equation} \label{eq: def_fibre_preserving_diffeo}
	\upsilon : \cU \to \tilde{\cU}, ~\xtheta \mapsto \upsilon \xtheta := (x; \tilde{\theta}),  
\end{equation}
which is homogeneous with respect to $\theta$ of degree $1$ and fibre-preserving $\widetilde{\Pr} \circ \upsilon = \Pr$ where $\Pr, \widetilde{\Pr} : \cU, \tilde{\cU} \to M$ are projection operators and $\tilde{\theta}$ is a smooth function of $\xtheta$.  
Then $\varphi$ is said to be (locally) \textbf{equivalent} to a phase function $\tilde{\varphi}$ on $\tilde{\cU}$ if~\cite[p. 134]{Hoermander_ActaMath_1971} 
(see also, e.g.~\cite[p. 27]{Weinstein_AMS_1977}) 
\begin{equation} \label{eq: def_equivalent_phase_function}
	\upsilon^{*} \tilde{\varphi} = \varphi. 
\end{equation}
Given any $\xxi \in \varLambda$ we would like to know the \textit{necessary} and \textit{sufficient} conditions for the existence of equivalent phase functions.   
Let $\varphi$ and $\tilde{\varphi}$ be phase functions for $\varLambda$ with the same excess $e$ on some open conic subsets $U \times \dotRn$ and $\tilde{U} \times \dot{\R}^{\tilde{n}}$, respectively, such that $\Pi \xxi \in U \subset M$. 
Then $\varphi$ and $\tilde{\varphi}$ are \textit{equivalent} on any open conic neighbourhoods $\cU$ and $\tilde{\cU}$ of $\Pi \xxi$ in $U \times \dotRn$ and $\tilde{U} \times \dot{\R}^{\tilde{n}}$, respectively, if and only if~\cite[Prop. 1.5 (p. 421)]{Treves_Plenum_1980} 
(see also~\cite[Prop. 1]{Meinrenken_ReptMathPhys_1992}) 
\begin{subequations} \label{eq: necessarily_sufficient_equivalent_clean_phase_function}
	\begin{eqnarray} 
        \sC_{\varphi} \ni \xtheta \mapsto (x, \rd_{x} \varphi = \xi) 
		& = & 
        (x, \rd_{x} \tilde{\varphi} = \xi)  \mapsfrom (x; \tilde{\theta}) \in \sC_{\tilde{\varphi}}, 
        \\ 
		n 
		& = & 
		\tilde{n}, 
		\\ 
        \sgn \big( \Hess_{\theta} \varphi \, \xtheta \big) 
		& = & 
        \sgn \big( \Hess_{\tilde{\theta}} \varphi \, (x; \tilde{\theta}) \big),   
    \end{eqnarray}
\end{subequations}
where $\sC_{\varphi}$ resp. $\sC_{\tilde{\varphi}}$ are fibre-critical manifolds of $\varphi$ resp. $\tilde{\varphi}$ and $\sgn$ denotes the signature of the Hessian (\eqref{eq: def_Hessian}) matrix $\Hess_{\theta} \varphi$ of $\varphi$ with respect to its fibre-variable $\theta \in \dotRn$. 
These constraints can be alternatively given an operational viewpoint. 

%
%
%
\begin{remark} \label{rem: stable_equivalence}
    Let $\kM (\varLambda)$ be the class of equivalent phase functions for $\varLambda$. 
    Given $\varphi \in \kM (\varLambda)$, the following operations produce another element of $\kM (\varLambda)$~\cite[p. 48]{Bates_AMS_1997} 
    (see also, e.g.~\cite[p. 671]{Rudolph_Springer_2013}): 
    \begin{enumerate}
        \item[(i)] 
        \textit{Addition}: 
        $\tilde{\varphi} := \varphi + \lambdaup \in \kM (\varLambda)$ for any $\lambdaup \in \R$; 
        \item[(ii)] 
        \textit{Composition}: 
        $\tilde{\varphi} = \eqref{eq: def_equivalent_phase_function} \in \kM (\varLambda)$ for any fibre-preserving diffeomorphism~\eqref{eq: def_fibre_preserving_diffeo}; 
        \item[(iii)] 
        \textit{Suspension}: 
        For any non-degenerate bilinear form $\mathsf{q}$ on $\dot{\R}^{\tilde{n}}$, $\tilde{\varphi} := \Pr_{\cU}^{*} \varphi + \Pr_{ \dot{\R}^{\tilde{n}} }^{*} \mathsf{q}$ on $\tilde{\cU} := U \times \dot{\R}^{\tilde{n}}$ is an element of $\kM (\varLambda)$. 
        \item[(iv)] 
        \textit{Restriction}: 
        If $\varphi$ is a phase function on $M \times \dotRn$ then $\tilde{\varphi} := \varphi \upharpoonright \tilde{\cU} \in \kM (\varLambda)$ where $\tilde{\cU} \subset M \times \dotRn$ is an open conic subset containing the preimage of $\jmath_{\varphi} \xxi$. 
    \end{enumerate}
    These four operations generate an equivalence relation among phase functions called \textbf{stable equivalence}. 
    So, we conclude that \textit{a phase function is not unique yet it is so up to stable equivalence}. 
\end{remark} 
%
%
%

The parametrisation of a Lagrangian submanifold by means of an equivalent class of phase functions \textit{cannot} be performed \textit{globally} due to \textit{topological restrictions}.  
When $M$ is compact, these obstructions require cohomological and $k$-theoretic language to formulate, which are somehow tangential to the subject matter so we refrain those precise formulae and refer~\cite{Latour_ASENS_1991, Giroux_Springer_1990} together with the earlier references cited therein. 
The non-compact case is still an open issue to the best of our knowledge. 
Notwithstanding, a global parametrisation can be achieved, for instance, by an equivalence class of $\C$-valued non-degenerate phase functions in the particular case whenever conic Lagrangian submanifolds are given by the graphs of homogeneous symplectomorphisms~\cite[Lem. 1.2, 1.7]{Laptev_CPAM_1994}.
\newline 

Finally, the fibre-dimension is bounded from below by~\cite[Prop. 1.3 (p. 417)]{Treves_Plenum_1980} 
\begin{equation} \label{eq: lower_bound_fibre_dimension_phase_function}
    n \geq d - \rk (\rd_{\xxi} \Pi) + e
\end{equation}
for any $\xxi \in \varLambda$, where $\rd \Pi$ is the differential of~\eqref{eq: def_projection_Lagrangian_M}.  
\newline 

If $(\fM, \sigmaup_{\ms \fM})$ and $(\fN, \sigmaup_{\ms \fN})$ are two symplectic manifolds then their  product manifold $\fM \times \fN$ can be naturally equipped with the symplectic forms 
\begin{equation} \label{eq: def_product_symplectic_form}
    \sigmaup_{\ms \fM} \boxplus \sigmaup_{\ms \fN} := \pr_{\ms \fM}^{*} \sigmaup_{\ms \fM} + \pr_{\ms \fN}^{*} \sigmaup_{\ms \fN}, 
    \quad 
    \sigmaup_{\ms \fM} \boxminus \sigmaup_{\ms \fN} := \pr_{\ms \fM}^{*} \sigmaup_{\ms \fM} - \pr_{\ms \fN}^{*} \sigmaup_{\ms \fN},  
\end{equation}
where $\pr_{\ms \fM}, \pr_{\ms \fN} : \fM \times \fN \to \fM, \fN$ are the Cartesian projectors. 
In symplectic coordinates $(x^{i}, \xi_{i}; y^{j}, \eta_{j})$ for $\fM \times \fN$, these read $\rd x^{i} \wedge \rd \xi_{i} \boxplus \rd y^{j} \wedge \rd \eta_{j}$ resp. $\rd x^{i} \wedge \rd \xi_{i} \boxminus \rd y^{j} \wedge \rd \eta_{j}$. 

%
%
%
\begin{example} \label{exm: graph_symplectomorphism_Lagrangian_submf}
    Let $(\fM, \sigmaup_{\ms \fM}), (\fN, \sigmaup_{\ms \fN})$ be $2d$-dimensional symplectic manifolds and $\varkappa : \fN \to \fM$ a symplectomorphism. 
    Then its graph 
    \begin{equation}
    	\varGamma = \big\{ \big( \varkappa \yeta; \yeta \big) | \yeta \in \fN \big\} \subset (\fM \times \fN, \sigmaup_{\ms \fM} \boxminus \sigmaup_{\ms \fN})
    \end{equation}
    is a Lagrangian submanifold whereas the twisted graph $\varGamma' = \big\{ \big( \varkappa \yeta; (y, - \eta) \big) \big\}$ is a Lagrangian submanifold with respect to the product symplectic form $\sigmaup_{\ms \fM} \boxplus \sigmaup_{\ms \fN}$.  
    Note, \textit{any diffeomorphism $\varkappa$ is a symplectomorphism if and only if its graph $\varGamma \subset (\fM \times \fN, \sigmaup_{\ms \fM} \boxminus \sigmaup_{\ms \fN})$ is a Lagrangian submanifold}. 
\end{example}
%
%
%

A special case of the preceding example for the identity symplectomorphism has some interesting consequences as presented below.  

%
%
%
\begin{example} \label{exm: diagonal_embedding_Lagrangian_submf}
    Let $M$ be a manifold whose punctured cotangent bundle $(\dotCoTanM, \sigmaup_{\ms M})$ is considered as a conic symplectic manifold as in Example~\ref{exm: cotangent_bundle_symplectic_mf} and let $\id : M \to M$ be the identity diffeomorphism. 
    Then its cotangent lift 
    \begin{equation}
        \coTan \id : \dotCoTanM \to \dotCoTanM, ~\xxi \mapsto 
        (\coTan_{x} \! \id) (\xi) \, (X) := \xi \big( \rd_{x} \! \id \, (X) \big) = \xi (X) 
    \end{equation}
    for any $X \in \tangent_{x} M$, is the identity homogeneous symplectomorphism preserving the Liouville form. 
    Its twisted graph     
    \begin{equation}
        \varGamma' 
        = 
        \{ (x, \xi; x, - \xi) \in \dotCoTanM \times \dotCoTanM \} 
        = 
        (\varDelta_{M} M)^{\perp *} 
        =  
        (\varDelta \dotCoTanM)' 
        \subset (\dotCoTanM \times \dotCoTanM, \sigmaup_{\ms M} \boxplus \sigmaup_{\ms M}) 
    \end{equation}
    is a Lagrangian submanifold and it is nothing but the conormal bundle (Example~\ref{exm: conormal_bundle_Lagrangian_submf}) $(\varDelta_{M} M)^{\perp *}$ of the diagonal embedding $\varDelta_{M} : M \to M \times M$. 
    Equivalently, this is the twisted diagonal embedding $(\varDelta \, \dotCoTanM)' \hookrightarrow \dotCoTanM \times \dotCoTanM \cong \dotCoTanMM$ of the product punctured cotangent bundle.  
    The manifold $\varGamma'$ is parametrised by the non-degenerate phase function 
    \begin{equation}
    	\varphi := (x^{i} - y^{i}) \xi_{i}
    \end{equation}
    for any coordinates $(x^{i})$ (resp. $(y^{i})$) around $x$ (resp. $y$) in $M$ and $\xi \in \dotRd$ 
    (see e.g.~\cite[Exm. 5.2, 5.3 (pp. 455-456)]{Treves_Plenum_1980}). 
\end{example}
%
%
%

Example~\ref{exm: graph_symplectomorphism_Lagrangian_submf} exhibits that the standard recipe to construct a Lagrangian submanifold $\varLambda \subset (\fM \times \fN, \sigmaup_{\ms \fM} \boxminus \sigmaup_{\ms \fN})$ as the graph of a symplectomorphism only works provided that $\dim \fM = \dim \fN$ which cannot be guaranteed a priori and relaxing this condition leads to the notion of a canonical relation.  

%
%
%
\begin{definition} \label{def: canonical_relation} 
    If $(\fM, \sigmaup_{\ms \fM})$ and $(\fN, \sigmaup_{\ms \fN})$ are two symplectic manifolds then an immersed Lagrangian submanifold $(C, \iota)$ of $(\fM \times \fN, \sigmaup_{\ms \fM} \boxminus \sigmaup_{\ms \fN}:= \eqref{eq: def_product_symplectic_form})$ is called a \textbf{canonical relation} from $\fN$ to $\fM$.  
    If $\fM, \fN$, and $C$ all are conic, then we call the relation \textbf{homogeneous}~\cite[Def. 4.1.2]{Hoermander_ActaMath_1971}  
    (see also~\cite[Def. 21.2.12]{Hoermander_Springer_2007}). 
\end{definition}
%
%
%

\begin{example}
    Let $(\fM, \sigmaup_{\ms \fM})$ be a symplectic manifold and let $(\coTan \R \times \fM, \sigmaup_{\ms \R} \boxplus \sigmaup_{\ms \fM}, H)$ be a time-dependent Hamiltonian system. 
    The twisted graph  
    (see e.g.~\cite[pp. 441 - 442]{Rudolph_Springer_2013},~\cite[Sec. 4.2]{Meinrenken_ReptMathPhys_1992} for details) 
    \begin{equation}
        \varGamma' = \{ (t, \tau; x, \xi; y, - \eta) \in \coTan \R \times \fM \times \fM \,|\, \xxi = \varPhi_{t} \yeta, \tau + H \xxi = 0 \}  
    \end{equation}
    of $H$-flow $\varPhi_{t}$ which is assumed to be complete, is a canonical relation of $\coTan \R \times \fM \times \fM$.  
\end{example}
%
%
%

\begin{remark} \label{rem: twisted_canonical_relation_product_symplectic_form}
    By Definition~\ref{def: canonical_relation}, $C \subset \fM \times \fN$ is a Lagrangian submanifold with respect to $\sigmaup_{\ms \fM} \boxminus \sigmaup_{\ms \fN}$, then it follows that the twisted canonical relation $C' \subset \fM \times \fN$ is a Lagrangian submanifold with respect to the product symplectic form $\sigmaup_{\ms \fM} \boxplus \sigmaup_{\ms \fN}$. 
\end{remark}
%
%
%

Let $C$ be a homogeneous canonical relation from $\pi_{\ms N} : \dotCoTanN \to N$ to $\pi_{\ms M} : \dotCoTanM \to M$. 
Being a Lagrangian submanifold, it can also be uniquely parametrised by a  stable equivalence class $[\varphi]$ of phase functions $\varphi$ defined on an open conic subset of $M \times N$ as described below. 
For each $\xxiyetaNot \in C$, there exists a homogeneous symplectic coordinate neighbourhood $\cU \times \cV$ on $\dotCoTanM \times \dotCoTanN$ and a clean phase function $\varphi$ on $\pi_{\ms M} (\cU) \times \pi_{\ms N} (\cV) \times \dotRn$ such that 
\begin{equation} \label{eq: fibre_critical_mf_M_N_Rn}
    \sC_{\varphi} := \{ \xythetaNot \in \pi_{\ms M} (\cU) \times \pi_{\ms N} (\cV) \times \dotRn \,|\, \grad_{\theta} \varphi \, \xythetaNot = 0 \} \subset M \times N \times \dotRn 
\end{equation}
is a $\dM + \dN + e$-dimensional submanifold for some $\theta^{0} \in \dotRn$. 
Then, $C$ is locally obtained by the homogeneous fibration   
\begin{equation} \label{eq: def_canonical_relation_fibration}
    \jmath_{\varphi} : \sC_{\varphi} \to C_{\varphi}, ~\xytheta \mapsto \jmath \xytheta := (x, \rd_{x} \varphi; y, \rd_{y} \varphi) 
\end{equation}
of fibre-dimension $e$. 
Obviously, one can opt $e = 0$ so that the preceding mapping is a homogeneous immersion. 
Alternatively, on any homogeneous symplectic chart $\big(\cU \times \cV, (x^{i}, \xi_{i}; y^{j}, \eta_{j}) \big)$ around $\Pi \xxiyeta$, a non-degenerate phase function of the form 
\begin{equation} \label{eq: def_canonical_nondegenerate_phase_function_canonical_relation}
    \psi := x^{i} \xi_{i} + y^{j} \eta_{j} - H  
\end{equation}
can be chosen to parametrise $C$, where $H$ is a smooth positively homogeneous $\R$-valued function of degree $1$ on an open conic neighbourhood of $(\xi, \eta) \in \dotCoTan_{x} M \times \dotCoTan_{y} N$. 
The corresponding homogeneous immersion 
\begin{equation} \label{eq: def_canonical_relation_identity_immersion}
    \jmath_{\ms H} : \sC_{H} \to C_{H}, ~ \xxiyeta \mapsto \left( \parDeri{\xi_{i}}{H}, \xi_{i}; \parDeri{\eta_{j}}{H}, \eta_{j} \right) 
\end{equation}
is then actually the identity map. 
Complying with Remark~\ref{rem: twisted_canonical_relation_product_symplectic_form}, this means 
\begin{equation} \label{eq: canonical_relation_clean_phase_function}
    C'_{\varphi} = \{ (x, \rd_{x} \varphi; y, - \rd_{y} \varphi) \} 
    \quad \textrm{and} \quad  
    C'_{H} = \Big\{ \Big( \parDeri{\xi_{i}}{H}, \xi_{i}; \parDeri{\eta_{j}}{H}, - \eta_{j} \Big) \Big\}. 
\end{equation}
are canonical relations with respect to $\sigmaup_{\ms M} \boxplus \sigmaup_{\ms N}$ 
(see e.g.~\cite[$(21.2.9)$]{Hoermander_Springer_2007}). 
%
%
%
%
%
%
%
%
%
%
\subsection{Clean composition}
\label{sec: clean_composition}
A relation maps a set to another set, thus, if $S$ is a subset of a symplectic manifold $(\fO, \sigmaup_{\ms \fO})$ and $C \subset (\fM \times \fO, \sigmaup_{\ms \fM} \boxminus \sigmaup_{\ms \fO})$ is a canonical relation then $C (S) = \{ \xi \in \fM \,|\, (\xi, \zeta) \in C, \textrm{for some} \, \zeta \in S \}$ can be expressed as $C (S) = \pr_{\ms \fM} \big( C \cap \pr_{\ms \fO}^{-1} S \big)$ where $\pr_{\ms \fM}, \pr_{\ms \fO} : \fM \times \fO \to \fM, \fO$ are the Cartesian projectors. 
The composition of $C$ with another canonical relation $\varLambda \subset (\fO \times \fN, \sigmaup_{\ms \fO} \boxminus \sigmaup_{\ms \fN})$ then reads  
(see e.g.~\cite[p. 290]{Hoermander_Springer_2007}) 
\begin{equation}
    C \circ \varLambda 
    :=
    \{ (\xi, \eta) \in \fM \times \fN \,|\, (\xi, \zeta) \in C, (\zeta, \eta) \in \varLambda, \textrm{for some} \, \zeta \in \fO \} 
    = 
    \varPi (C {\scriptsize \text{\FiveStarOpen}} \varLambda), 
\end{equation}
where we have used 
\begin{equation} \label{eq: def_C_star_varLambda}
	C {\scriptsize \text{\FiveStarOpen}} \varLambda := (C \times \varLambda) \bigcap (\fM \times \varDelta \fO \times \fN) 
\end{equation}
and the Cartesian projector $\varPi : \fM \times \varDelta \fO \times \fN \to \fM \times \fN$. 
\newline 

But, in general, $C \circ \varLambda$ is \textit{not} a manifold and one requires to know ``how'' $C \times \varLambda$ intersects with $\fM \times \varDelta \fO \times \fN$ for which we introduce the terminology of clean intersection 
(see e.g.~\cite[Sec. 5]{Duistermaat_InventMath_1975}).  

%
%
%
\begin{definition} \label{def: clean_intersection}
	Let $S$ and $\varSigma$ be two submanifolds of a manifold $M$. 
    If their intersection $S \cap \varSigma$ is a submanifold satisfying 
    (see e.g.~\cite[Def. C.3.2]{Hoermander_Springer_2007},~\cite[Def. 5.9]{Bates_AMS_1997}) 
    \begin{equation} \label{eq: def_clean_intersection} 
        \forall x \in S \cap \varSigma : 
        \tangent_{x} (S \cap \varSigma) = \tangent_{x} S \cap \tangent_{x} \varSigma  
    \end{equation}
    then the intersection is said to be \textbf{clean}, denoted by $S \pitchfork_{e} \varSigma$, with \textbf{excess} $e \in \NO$ given by 
    \begin{equation} \label{eq: def_excess} 
    	e := \codim S + \codim \varSigma - \codim (S \cap \varSigma).   
    \end{equation} 
    If $e = 0$ then $S$ and $\varSigma$ are said to intersect \textbf{transversally}, symbolically $S \pitchfork \varSigma$. 
\end{definition}
%
%
%

The criterion that $S \pitchfork \varSigma$ at any $x \in S \cap \varSigma$ can be restarted as 
\begin{equation}
    \tangent_{x} S + \tangent_{x} \varSigma = \tangent_{x} M \Leftrightarrow S_{x}^{\perp} \cap \varSigma_{x}^{\perp} = \{0\}. 
\end{equation}

%
%
%
\begin{example}
    Let $(x, y, z) \in \R^{3}$. 
    Then $\{ z = 0 \} \pitchfork \{ z = x \}$ but $\{ z = 0 \} \not\pitchfork \{ z = xy \}$. 
\end{example}
%
%
%

Intuitively speaking, $S \pitchfork \varSigma$ means they intersect ``as little as possible''. 
Thus, $S \cap \varSigma = \emptyset$ is trivially $S \pitchfork \varSigma$. 
%
%
%
\begin{example}
    The $x$-axis intersects with the $y$-axis in $\R^{3}$ cleanly with $e = 2 +2 - 3 = 1$. 
    In $\R^{2}$, the $x$-axis and $\{ y = x^{2} \}$ do not intersect cleanly. 
\end{example}
%
%
%

If $(C \times \varLambda) \pitchfork_{e} (\fM \times \varDelta \fO \times \fN)$ then 
(see e.g.~\cite[Thm. 21.2.14]{Hoermander_Springer_2007},~\cite[Prop. 5.28]{Bates_AMS_1997}) 
\begin{equation}
    \rk \varPi = \frac{\dim \fM + \dim \fN}{2}
\end{equation}
and the range $C \circ \varLambda$ of $\varPi$ is locally an \textit{immersed Lagrangian submanifold} of $(\fM \times \fN, \sigmaup_{\ms \fM} \boxminus \sigmaup_{\ms \fN})$. 
Furthermore, if $U$ is an open set in $\fM \times \fN$ and $\varPi \upharpoonright \varPi^{-1} (U)$ is \textit{proper} then $\varPi (\varPi^{-1} U)$ is an immersed \textit{closed} Lagrangian submanifold of $U$ and $\varPi^{-1} U$ is a fibre space with \textit{compact fibres} over it 
(see e.g.~\cite[Rem. p. 291]{Hoermander_Springer_2007}). 

%
%
%
\begin{definition} \label{def: clean_composition} 
    Let $(\fM, \sigmaup_{\ms \fM}), (\fO, \sigmaup_{\ms \fO})$, and $(\fN, \sigmaup_{\ms \fN})$ be symplectic manifolds and let $C \subset (\fM \times \fO, \sigmaup_{\ms \fM} \boxminus \sigmaup_{\ms \fO})$ and $\varLambda \subset (\fO \times \fN, \sigmaup_{\ms \fO} \boxminus \sigmaup_{\ms \fN})$ be closed canonical relations. 
    Then $C$ and $\varLambda$ is said to have a \textbf{clean composition $C \circ_{e} \varLambda$ with excess $e \in \NO$} if  $C \times \varLambda$ and $\fM \times \varDelta \fO \times \fN$ intersects cleanly with excess $e$. 
    The composition is called \textbf{proper} if the Cartesian projector $\varPi : \fM \times \varDelta \fO \times \fN \to \fM \times \fN$ is proper 
    (see e.g.~\cite[Thm. 21.2.14]{Hoermander_Springer_2007},~\cite[Def. 5.2 $($p. 458$)$]{Treves_Plenum_1980}). 
\end{definition}
%
%
%

We remark that 
(see e.g.~\cite[Thm. 4.2.2]{Guillemin_InternationalP_2013}) 
\begin{equation} \label{eq: excess_fibre_projection_clean_composition}
    e = \dim \kF,   
\end{equation}
where 
\begin{eqnarray} \label{eq: def_fibre_projection_clean_composition}
	\kF & := & \{ \kF_{(\xi, \eta)} \,|\, (\xi, \eta) \in C \circ_{e} \varLambda \}, 
    \nonumber \\ 
    \kF_{(\xi, \eta)} & := & \{ (\xi, \zeta, \zeta, \eta) \in C {\scriptsize \text{\FiveStarOpen}} \varLambda \,|\, \varPi (\xi, \zeta, \zeta, \eta) = (\xi, \eta) \in C \circ_{e} \varLambda \}
\end{eqnarray}
is the set of compact connected fibres $\kF_{\xi, \eta}$ of the smooth fibration 
\begin{equation} \label{eq: varPi}
	\varPi : C {\scriptsize \text{\FiveStarOpen}} \varLambda \to C \circ_{e} \varLambda. 
\end{equation}

In terms of generating functions, if $\varphi (x, z; \theta)$ and $\phi (z, y; \vartheta)$ are non-degenerate phase functions for homogeneous canonical relations $C \subset \dotCoTanM \times \dotCoTan O$ and $\varLambda \subset \dotCoTan O \times \dotCoTanN$ in open conic neighbourhoods of $(x_{0}, z_{0}; \theta^{0}) \in M \times O \times \dot{\R}^{\nM}$ and $(z_{0}, y_{0}; \vartheta^{0}) \in O \times N \times \dot{\R}^{\nN}$ such that $(x_{0}, z_{0}; \theta^{0}) \in \sC_{\varphi}$ and $(z_{0}, y_{0}; \vartheta^{0}) \in \sC_{\phi}$, respectively, with $\rd_{z} \varphi \, (x_{0}, z_{0}; \theta^{0}) + \rd_{z} \phi \, (z_{0}, y_{0}; \vartheta^{0}) = 0$, then 
(see e.g.~\cite[Prop. 21.2.19]{Hoermander_Springer_2007})
\begin{equation} \label{eq: def_phase_function_clean_composition_canonical_relatation}
	\psi (x, z, y; \theta, \vartheta) := \varphi (x, z; \theta) + \phi (z, y; \vartheta) 
\end{equation}
is a clean phase function with excess $e$ parametrising $C \circ_{e} \varLambda$. 
In addition, if $\psi$ is non-degenerate in $(z'; \theta', \vartheta')$ then the $e$ variables $(z''; \theta'', \vartheta'')$ parametrise $\kF$; the meaning of the splitting $z = (z', z'')$ has been explained in  Definition~\ref{def: symbol_FIO_Euclidean}.
%
%
%
%
%
%
%
%
%
%
\subsection{Keller-Maslov bundle} 
\label{sec: Keller_Maslov_bundle}
Given a Lagrangian submanifold $\varLambda$, one can associate a number of bundle structures over $\varLambda$ such that these bundles encode information about how the local pieces of $\varLambda$  generated by phase functions are glued together. 

%
%
%
\begin{definition} \label{def: Keller_Maslov_bundle}
	Let $\dotCoTanM \to M$ be the punctured cotangent bundle over a manifold $M$ and $\varLambda \subset \dotCoTanM$ a conic Lagrangian submanifold. 
    The \textbf{Keller-Maslov bundle} $\M \to \varLambda$ over $\varLambda$ is a complex line bundle obtained from some (as detailed below) principal bundle with structure group $\nicefrac{\Z}{4 \Z}$, which is trivial as a vector bundle~\cite[p. 148]{Hoermander_ActaMath_1971}.  
\end{definition}
%
%
%

\begin{remark} \label{rem: flat_connection_Maslov_bundle}
    Having a discrete structure group, $\Maslov$ carries a natural flat connection with holonomy in $\nicefrac{\Z}{4 \Z}$ 
    (see e.g.~\cite[p. 56]{Bates_AMS_1997}).   
\end{remark}
%
%
%

We are now going to present a construction of $\M$ which is an integral part of the principal symbol (Definition~\ref{def: symbol_FIO}) of a Lagrangian distribution and consequently has played a pivotal part  in the non-zero periodic contribution (Theorem~\ref{thm: trace_formula_t_T}) of the Duis- termaat-Guillemin-Gutzwiller trace formula. 
Let $\varphi$ be a \textit{non-degenerate} phase function for $\varLambda$. 
Then employing the Hessian of $\varphi$ one has the integer-valued map: 
$\sC_{\varphi} \ni \xtheta \mapsto \sgn \big( \Hess_{\theta} \varphi \, \xtheta \big) \in \Z$. 
Note, this function can be somewhat discontinuous as the Hessian can be singular at some elements of $\sC_{\varphi}$. 
We now consider an open conic 
(Leray\footnote{This means that finite intersection of $\cL_{\alpha}$'s are either empty or diffeomorphic to the open ball.}-) 
covering $\{ \cL_{\alpha} \}$ of $\varLambda$, indexed by a countable set with the corresponding non-degenerate phase functions $\varphi_{\alpha}$ and the fibre-critical manifolds $\sC_{\alpha}$. 
It follows from~\eqref{eq: necessarily_sufficient_equivalent_clean_phase_function} that the  mapping~\cite[$(3.2.15)$]{Hoermander_ActaMath_1971}
(see also, e.g.~\cite{Meinrenken_ReptMathPhys_1992},~\cite[Sec. 5. 13]{Guillemin_InternationalP_2013},~\cite[pp. 686-687]{Rudolph_Springer_2013}) 
\begin{eqnarray}
    && 
    \ks_{\alpha \beta} : \sC_{\alpha} \cap \sC_{\beta} \to \Z, ~ (x; \theta^{\alpha}, \theta^{\beta}) \mapsto 
    \nonumber \\ 
    && 
    \ks_{\alpha \beta}(x; \theta^{\alpha}, \theta^{\beta}) 
    := \frac{1}{2} 
    \Big( 
    \sgn \big( \Hess_{\theta} \varphi_{\beta} \, (x; \theta^{\beta}) \big) - n_{\beta} 
    - 
    \sgn \big( \Hess_{\theta} \varphi_{\alpha} \, (x; \theta^{\alpha}) \big) + n_{\alpha} 
    \Big) \qquad 
\end{eqnarray}
is constant for all connected intersections $\sC_{\alpha} \cap \sC_{\beta}$. 
Thereby, 
\begin{equation} \label{eq: def_transition_function_Maslov_bundle}
    g_{\alpha \beta} := \re^{\nicefrac{\ri}{2} \pi \ks_{\alpha \beta} \circ \jmath_{\alpha \beta}^{-1}} : \cL_{\alpha} \cap \cL_{\beta} \to \dot{\C}
\end{equation}
clearly satisfies the cocycle property: $g_{\alpha \beta'} g_{\beta' \beta} = g_{\alpha \beta}$ whenever $\cL_{\alpha} \cap \cL_{\beta'} \cap \cL_{\beta} \neq \emptyset$ with $|g_{\alpha \beta}| = 1$. 
Hence, the collection of non-zero complex numbers $\{ g_{\alpha \beta} \}$ for all $\alpha$ and $\beta$ such that $\varLambda_{\alpha}$ and $\varLambda_{\beta}$ are not disjoint, defines our transition function and the global construction of $\M$ is achieved by taking the disjoint union 
\begin{equation} \label{eq: def_construct_Maslov_bundle}
    \M := \bigsqcup_{\alpha} (\cL_{\alpha} \times \C) / \sim 
\end{equation}
modulo the equivalence relation 
\begin{equation} \label{eq: def_equivalence_relation_Maslov_bundle}
    \big( (x, \xi^{\alpha}), c_{\alpha} \big) \sim \big( (x, \xi^{\beta}), c_{\beta} \big) 
    \Leftrightarrow 
    (x, \xi^{\alpha}) = (x, \xi^{\beta}) \in \cL_{\alpha} \cap \cL_{\beta}, c_{\alpha} = g_{\alpha \beta} c_{\beta}. 
\end{equation}
The constant functions from $\cL_{\alpha}$ to $1_{\C}$ form flat local sections of $\Maslov$: 
\begin{equation}
    \mathtt{e}_{\alpha} : \cL_{\alpha} \to \Maslov, ~ \xxi \mapsto \big[ \big( \xxi, 1_{\C} \big) \big]
\end{equation}
and hence manufacture $\Maslov$ into a line bundle with flat connection over $\varLambda$. 
Any section $\m$ of $\Maslov \to \cL_{\alpha}$ then can be written as $\m = \sum \m^{\alpha} \mathtt{e}_{\alpha}$ for unique functions $\m^{\alpha}$. 
\newline 

Alternatively, one can define a $\sqrt[8]{1_{\C}}$-valued mapping 
\begin{equation} \label{eq: def_section_Maslov_bundle}
    \mathtt{e}_{\varphi} := \re^{\nicefrac{\ri}{4} \pi \sgn (\Hess_{\theta} \varphi) \circ \jmath^{-1}} : \cL_{\alpha} \to \dot{\C} 
\end{equation}
and declare this as the flat section of $\M \upharpoonright \cL_{\alpha}$. 
Then $\M$ is constructed employing transition functions $g_{\alpha \beta}$ valued in the structure group $\nicefrac{\Z}{4 \Z} = \{\pm 1, \pm \ri \}$, which are locally constant on each disjoint covering of $\varLambda$ and absolute value $1$ so that $\mathtt{e}_{\varphi_{\alpha}} = g_{\alpha \beta} \mathtt{e}_{\varphi_{\beta}}$.  
\newline 

Since $\varphi$ is unique only up to the stable equivalence, the above construction must be examined under the four constituent operations mentioned in Remark~\ref{rem: stable_equivalence}. 
Clearly $\sgn (\Hess_{\theta} \varphi_{\alpha})$ is invariant under all those except suspension for which $\sgn (\Hess_{\theta} \varphi_{\alpha}) = \sgn (\Hess_{\tilde{\theta}} \tilde{\varphi}_{\alpha}) + \sgn \mathsf{q}$ --- ensuring the construction for the equivalence class $[\varphi_{\alpha}]$ rather than just a representative $\varphi_{\alpha}$. 
\newline 

So far, the construction has been given for a non-degenerate phase function but it can be well adopted for a clean phase function $\psi \xtheta$ on $U \times \dotRn$ as described below.  
Let $\rk \big( \Hess_{\theta} \psi \, \xthetaNot \big) = e$ for some $0 < e \leq n$ on the fibre-critical manifold $\sC_{\psi} \ni \xthetaNot$. 
We can split (cf. Definition~\ref{def: symbol_FIO_Euclidean}) $\dotRn \ni \theta = (\theta', \theta'') \in \dot{\R}^{n-e} \times \dot{\R}^{e}$ the fibre-variable $\theta$ into two groups and perform a linear change of $\theta$-coordinates such that the matrix $\Hess_{\theta'} \psi \, (x_{0}; \theta^{\prime 0}, \theta^{\prime \prime 0})$ is non-degenerate for a fixed $\theta''^{0} \in \dot{\R}^{e}$ and that the fibration~\cite{Duistermaat_InventMath_1975}
(see also, e.g.~\cite[Thm. 5.8.2, p. 119, 128]{Guillemin_InternationalP_2013}) 
\begin{equation}
    \varsigma : U \times \dotRn \to U \times \dot{\R}^{e}, ~ (x; \theta', \theta'') \mapsto (x; \theta'')
\end{equation}
produces a non-degenerate phase function $\phi (x; \theta', \theta''^{0})$ on $U \times \dot{\R}^{n-e}$ for $\varLambda$: 
\begin{equation} \label{eq: clean_nondegenerate_phase_function}
    \psi = \varsigma^{*} \phi. 
\end{equation}
Consequently,  the $\theta$-variable can be reduced so that $\Hess_{\theta''} \psi \, (x_{0}; \theta'^{0}, \theta''^{0}) = 0$. 
Evidently, $\sC_{\psi} = \varsigma^{*} \sC_{\phi}$ and the signature map $\sgn \big( \Hess_{\theta} \psi \, \xthetaNot \big)$ is simply the pullback of the signature map $\sC_{\varphi} \to \Z$ associated with $\phi$. 
In addition, the immersion $\jmath_{\phi} : \sC_{\varphi} \to \varLambda$ lifts to a fibre-preserving fibration $\jmath_{\psi} : \sC_{\psi} \to \varLambda$ and hence the construction of $\Maslov \to \varLambda$ by means of $\psi$ proceeds analogously.   

%
%
%
\begin{example}[Keller-Maslov bundle $\bbL_{0} \to (\varDelta M)^{\perp *}$ over the conormal bundle $(\varDelta M)^{\perp *}$] 
\label{exm: Maslov_bundle_conormal_bundle}
    As in the terminologies of Example~\ref{exm: diagonal_embedding_Lagrangian_submf}, the non-degenerate phase function $\varphi := (x^{i} - y^{i}) \xi_{i}$ parametrises $(\varDelta M)^{\perp *}$ globally. 
    Then $\Hess_{\xi} \varphi = 0$ and thus we have a global transition function $g_{\alpha \beta} = 1$ resp. a global constant section $\mathtt{e}_{\alpha} = 1$ by~\eqref{eq: def_transition_function_Maslov_bundle} resp.~\eqref{eq: def_section_Maslov_bundle}.  
\end{example}
%
%
%
%
%
%
%
%
%
%
\subsection{Maslov index} 
\label{sec: Maslov_index}  
The exponent in~\eqref{eq: def_section_Maslov_bundle} has a beautiful topological interpretation for which let us recall some notions from symplectic topology. 
Let $(\sV, \sigmaup)$ be a $2d$-dimensional real symplectic vector space and we choose a $\sigmaup$-compatible complex structure $\fJ$ on $\sV$.  
A choice of $\fJ$ defines a positive-definite hermitian form $\scalarProdTwo{\cdot}{\cdot} :=  \sigmaup (\cdot, \fJ \cdot) + \ri \sigmaup (\cdot, \cdot)$ on the complex vector space $\sV$. 
Then the set $\sL (\sV)$ of all Lagrangian subspaces of $(\sV, \sigmaup)$ is manifold of dimension $d (d+1) / 2$, called the \textbf{Lagrange-Gra\ss{}mann manifold} of $(\sV, \sigmaup)$.  
Let us denote by $\rO (\varLambda)$, the isometry group of a given $\varLambda \in \sL (\sV)$ with respect to the scalar product on $\varLambda$ induced by $\scalarProdTwo{\cdot}{\cdot}$. 
It follows that $\rO (\varLambda) \subset \rU (\sV) = \rU (\sV, \scalarProdTwo{\cdot}{\cdot})$ as a Lie group 
(see e.g.~\cite[Prop. 7.5.3/3]{Rudolph_Springer_2013}) when one makes the identification $\rO (\varLambda) = \{ g \in \rU (\sV) | g \varLambda = \varLambda \}$. 
Then, the mapping $\rU (\sV) \ni g \mapsto \det g \in \bbS \subset \C$ values $\pm 1$ on $\rO (\varLambda)$ and every element of a Lagrangian-Gra\ss{}mannian defines the isomorphism 
$\rU (\sV) /\rO (\varLambda) \ni [g] \mapsto g \varLambda \in \sL (\sV)$ 
(see e.g.~\cite[Prop. 7.6.1]{Rudolph_Springer_2013}).  
Thus, we have an induced mapping 
\begin{equation} \label{eq: map_Lagrangian_Grassmannian_circle}
     \mathrm{det}_{\varLambda}^{2} : \sL (\sV) \to \bbS \subset \C. 
\end{equation}

%
%
%
\begin{definition} \label{def: Maslov_index_closed_curve_Lagrangian_Grassmannian}
    Let $\sL (\sV)$ be a Lagrangian-Gra\ss{}mannian of a symplectic vector space $(\sV, \sigmaup)$. 
    The \textbf{Maslov index} $\kn (\gamma)$ of a closed curve $\gamma : \bbS \to \sL (\sV)$ is defined to be the degree of the mapping $\det_{\varLambda}^{2} \circ \gamma : \bbS \to \bbS$ where $\det_{\varLambda}^{2}$ is defined by~\eqref{eq: map_Lagrangian_Grassmannian_circle}~\cite[Sec. 1.3-1.4]{Arnold_FAA_1967}  
    (see also, e.g.~\cite[Def. 7.7.1]{Rudolph_Springer_2013}). 
\end{definition}
%
%
%

\begin{proposition} \label{prop: Maslov_index_homotopy_invariant}
	The Maslov index is an integer-valued homotopy invariant quantity and it is independent of the choice of Lagrangian subspace or complex structure 
    (see e.g.~\cite[Prop. 7.7.2]{Rudolph_Springer_2013}).  
\end{proposition}
%
%
%

In order to define this index for closed curves on a Lagrangian immersion, let 
\begin{equation}
    \chiup: \iota^{*} (\tangent \coTanM) \to \varLambda \times \R^{2d} 
\end{equation}
be a trivialisation of the vertical geometric distribution to $\{0\} \times \Rd$, i.e., $\chiup ( \tangent_{\iota \xxi} \coTan_{x} M) = \{0\} \times \Rd$ for any $\xxi \in \varLambda$. 
Then, the image of $\rd_{\xxi} \iota : \tangent_{\xxi} \varLambda \to \tangent_{\iota \xxi} \coTanM$ is identified with a Lagrangian subspace in $\R^{2d}$ via $\chiup$ and one obtains the following mapping 
(see e.g.~\cite[$(12.6.3)$]{Rudolph_Springer_2013}) 
\begin{equation} \label{eq: def_map_Lagrangian_submf_Lagrangian_Grassmannian_R_2d}
    \kl : \varLambda \to \sL (\R^{2d}), ~ \xxi \mapsto \kl \xxi := \chiup \big( \rd_{\xxi} \iota \, (\tangent_{\xxi} \varLambda) \big).  
\end{equation}

%
%
%
\begin{definition} \label{def: Maslov_index_closed_curve}
    Let $\varLambda$ be a Lagrangian submanifold of a symplectic manifold of dimension $2d$ and $\sL (\R^{2d})$ the Lagrangian-Gra\ss{}mannian  of the symplectic space $\R^{2d}$. 
    Then, the \textbf{Maslov index} of a closed curve $\gamma$ in $\varLambda$ is defined by 
    (see e.g.~\cite[p. 678]{Rudolph_Springer_2013}) 
    \begin{equation}
        \km (\gamma) := \kn (\kl \circ \gamma),  
    \end{equation}
    where $\kn$ is the Maslov index of the closed curve $\kl \circ \gamma$ in $\sL (\R^{2d})$ and $\kl$ is defined by~\eqref{eq: def_map_Lagrangian_submf_Lagrangian_Grassmannian_R_2d}.  
\end{definition}
%
%
%

It turns out that $\km$ is equal to the so-called ``Maslov intersection index''~\cite[Thm. 1.5]{Arnold_FAA_1967}  
(see also, e.g.~\cite[Prop. 12.6.6]{Rudolph_Springer_2013}). 
The formal definition 
(see e.g.~\cite[Def. 12.6.5]{Rudolph_Springer_2013})
and details of the intersection index are tangential to the ongoing discussion so we abstain from that detour and instead remark that \textit{this index is an integer-valued homotopy invariant quantity with fixed end points}~\cite[Thm. 2.2]{Arnold_FAA_1967} 
and it can be computed by means of phase functions as explained below. 
For our purpose, it is sufficient to consider the situation when $\gamma : [0, 1] \to \varLambda$ is a closed curve with end points in the set $\Sigma_{0} (\varLambda) \subset \varLambda$ defined by the points where the differential of~\eqref{eq: def_projection_Lagrangian_M} has rank $d$. 
Suppose that we have picked a partition of $[0, 1]$ by the numbers $0 = s_{0} < s_{1} < \ldots < \ldots < s_{N} = 1$ such that there exist open conic subsets $\cL_{1}, \ldots, \cL_{N}$ of $\varLambda$ so that $\gamma ([s_{\alpha-1}, s_{\alpha}]) \subset \cL_{\alpha}$ and $\cL_{\alpha} \cap \cL_{\beta}$ are connected, and that the immersions $\cL_{\alpha} \hookrightarrow \dotCoTanM$ and~\eqref{eq: def_Lagrangian_submf_fibration} are equivalent.  
Then the Maslov index is given by 
(see e.g.~\cite[$(5)$]{Meinrenken_ReptMathPhys_1992},~\cite[Cor. 12.6.12]{Rudolph_Springer_2013})
\begin{equation}
    \km (\gamma) = \sum_{\alpha = 1}^{N - 1} \ks_{\alpha \alpha + 1}. 
\end{equation}
This reveals that $\M \to \varLambda$ is essentially defined by the representation $\re^{\ri \pi \km (\gamma) / 2} = \ri^{\km (\gamma)}$ of $\piup_{1} (\varLambda)$. 
In other words, sections $\m$ of $\M$ can be identified with functions on the universal cover of $\varLambda$ satisfying the equivalence relation 
(see e.g.~\cite[$(1.3.3)$]{Anne_Cubo_2006}) 
\begin{equation}
    \m \big(\gamma \cdot \xxi \big) = \ri^{- \km (\gamma)} \m \xxi. 
\end{equation}

We close this section by two simple examples to give an intuitive understanding of the abstract formulae. 

%
%
%
\begin{example}
    Let $\R^{2}$ be equipped with the canonical symplectic structure and $\xxi$ its symplectic coordinates. 
    Suppose that $\varLambda := \{ x^{2} + \xi^{2} = 1 \} \hookrightarrow \R^{2}$ and that its orientation is assumed clockwise. 
    For dimensional reasons $\varLambda$ is an embedded Lagrangian submanifold and we consider a closed curve $\R \ni s \mapsto \gamma (s)$ running clockwise through $\varLambda$ exactly once. 
    An element $\check{s}$ of the flow parameter is called a crossing of $\gamma$ if $\gamma (\check{s}) \in \Sigma (\varLambda)$ and $\parDeri{\xi}{x} \big( \gamma (\check{s}) \big) = 0$ whilst $- \sgn \parDeri{\xi}{x}$ equals to the signature of the Hessian of the phase function with respect to $\xi$ outside of $\gamma (\check{s})$. 
    Hence, the coorientation points from the lower half-plane to the upper half-plane at the critical point $(-1, 0)$, whereas at $(1, 0)$, the otherwise; yielding $\km (\gamma) = 2$. 
    \newline 

    The figure eight immersion in $\R^{2}$ can be worked out analogously where it turns out that the Maslov index vanishes 
    (see e.g.~\cite[Exm. 12.6.13]{Rudolph_Springer_2013} for details). 
\end{example}
%
%
%
%
%
%
%
%
%
%
\subsection{Natural volume form} 
\label{sec: volume_canonical_relation}
Let $\varphi$ be a non-degenerate phase function for $C$ whose fibre-critical manifold is $\sC$. 
Then the map 
\begin{equation} \label{eq: surjective_map_M_N_Rn_Rn}
    M \times N \times \dotRn \ni \xytheta \mapsto \bigg( \parDeri{\theta_{1}}{\varphi} \xytheta, \ldots, \parDeri{\theta_{n}}{\varphi} \xytheta \bigg) \in \dotRn
\end{equation} 
is surjective because 
\begin{equation}
    \forall \xytheta \in \sC : 
    \rk \bigg( \rd \Big( \parDeri{\theta_{k}}{\varphi} \Big) \xytheta \bigg) = n, \quad k = 1, \ldots, n,  
\end{equation}
by Definition~\ref{def: clean_phase_function}, forbye, $\tangent_{\xytheta} \sC = \ker \big( (\rd \frac{\partial \varphi}{\partial \theta_{k}}) \xytheta \big)$. 
Thus, $\sC$ can be endowed with the volume form 
(see e.g.~\cite[$(4.1.1)$]{Duistermaat_Birkhaeuser_2011}): 
\begin{equation} \label{eq: def_volume_fibre_critical_mf}
    \dVol_{\ms \sC} (X, Y) \, \Big( \rd \parDeri{\theta_{1}}{\varphi} \wedge \ldots \wedge \rd \parDeri{\theta_{n}}{\varphi} \Big) (Z^{1}, \ldots, Z^{n}) 
    = 
    (\rd x \otimes \rd y \otimes \rd \theta) (X, Y; Z)  
\end{equation} 
for any vector fields $X, Y, Z = (Z^{1}, \ldots, Z^{n})$ on $M, N, \dotRn$, respectively. 
In other words, 
\begin{equation}
	\dVol_{\ms \sC} := \delta (\grad_{\theta} \varphi) 
\end{equation}
is the quotient of the measure $\rd x \otimes \rd y \otimes \rd \theta$ on $M \times N \times \dotRn$ by the pullback of the delta-distribution $\delta$ in $\dotRn$ under the mapping~\eqref{eq: surjective_map_M_N_Rn_Rn}. 
If $(\lambda_{l}) = (\lambda_{1}, \ldots, \lambda_{\dM + \dN})$ denotes arbitrary local coordinates of $\lambda \in \sC$ and we extend them to a neighbourhood in $M \times N \times \dotRn$ then $\dVol_{\ms \sC} (\lambda)$ 
locally\footnote{For instance, see~\cite[$(6.1.1)$]{Hoermander_Springer_2003} 
    for the local definition of pullback of a distribution.} 
means~\cite[p. 143, Prop. 4.1.3]{Hoermander_ActaMath_1971}  
(see e.g.~\cite[p. 123]{Grigis_CUP_1994}) 
\begin{equation} \label{eq: volume_fibre_critical_mf_local}
    \dVol_{\ms \sC} (\lambda) = \left| \det \left( \begin{array}{ccc}
                                                        \dfrac{\partial \lambda_{l}}{\partial x^{i}} 
                                                        & 
                                                        \dfrac{\partial \lambda_{l}}{\partial y^{j}} 
                                                        & 
                                                        \dfrac{\partial \lambda_{l}}{\partial \theta_{k}} 
                                                        \\ 
                                                        \dfrac{\partial^{2} \varphi}{\partial x^{i} \partial \theta_{k}} 
                                                        & 
                                                        \dfrac{\partial^{2} \varphi}{\partial y^{j} \partial \theta_{k}} 
                                                        & 
                                                        \dfrac{\partial^{2} \varphi}{\partial \theta_{k} \partial \theta_{k'}}
                                  \end{array}
                                  \right) \right|^{-1} |\rd \lambda|, 
\end{equation}
where $|\rd \lambda|$ is the Lebesgue density on $\sC$. 
This expression is independent of $(\lambda_{l})$ but depends on the choice of $(x^{i}, y^{j})$ which we keep aside for the moment. 
Thereby 
\begin{equation}
    \dVol_{\ms C} := \big( (\rd \jmath)^{-1} \big)^{*} \dVol_{\ms \sC}
\end{equation}
is a \textit{natural} measure on $C$ defined by transferring $\dVol_{\ms \sC}$ via the immersion $\jmath := \eqref{eq: def_canonical_relation_fibration}$ and hence~\eqref{eq: volume_fibre_critical_mf_local} can be identified as the local expression of $\dVol_{\ms C} (\lambda)$ when $\lambda$ is considered as an element in $C$. 
Since $\varphi$ is non-degenerate, one can pick its canonical form $\psi := \eqref{eq: def_canonical_nondegenerate_phase_function_canonical_relation}$ so that $C$ is parametrised by~\eqref{eq: def_canonical_relation_identity_immersion} and $(\lambda_{l}) = (\xi_{i} = \partial \psi / \partial x^{i}, \eta_{j} = \partial \psi / \partial y^{j})$ can be chosen as local coordinates on $C$, yielding 
\begin{equation}
    \dVol_{C} \xxiyeta = \left| \det \left( \begin{array}{ccc}
                                                                \dfrac{\partial^{2} \psi}{\partial x^{i} \partial \xi_{i}} 
                                                                & 
                                                                \dfrac{\partial^{2} \psi}{\partial y^{j} \partial \eta_{j}} 
                                                                & 
                                                                \dfrac{\partial^{2} \psi}{\partial \xi_{i} \partial \eta_{j}}
                                                            \end{array}
                                                    \right) \right|^{-1} |\rd \xi| |\rd \eta|. 
\end{equation}
Had we insisted a generic $\varphi$ instead of the canonical form $\psi$, $C$ would be parametrised by~\eqref{eq: def_canonical_relation_fibration}. 
Then, one can go for the natural choice of coordinates $(\lambda_{l}) = (\xi_{i} = \partial \varphi / \partial x^{i}, \eta_{j} = \partial \varphi / \partial y^{j})$ on $C$ and obtains 
\begin{equation} \label{eq: volume_canonical_relation_nondegenerate_phase_function}
    \dVol_{\ms C} \xxiyeta = |\det (\Hess \varphi)|^{-1} |\rd \xi| |\rd \eta|. 
\end{equation}

So far, an arbitrary but fixed non-degenerate phase function $\varphi$ has been used to construct $\dVol_{\ms C}$. 
However, as noted before $\varphi$ is unique only up to the stable equivalence. 
So we consider another element $\tilde{\varphi} = \eqref{eq: def_equivalent_phase_function}$ of this equivalence class and then the corresponding Jacobian $\rJ (\tilde{\lambda}, \grad_{\tilde{\theta}} \tilde{\varphi}; x, \tilde{\theta})$ is related to that in~\eqref{eq: volume_fibre_critical_mf_local} by 
(see e.g.~\cite[p. 440]{Treves_Plenum_1980}) 
\begin{equation}
    \rJ (\lambda, \grad_{\theta} \varphi; x, \theta) = \big( \rJ (\tilde{\theta}; \theta) \big)^{2} \rJ (\tilde{\lambda}, \grad_{\tilde{\theta}} \tilde{\varphi}; x, \tilde{\theta}), 
\end{equation} 
where $\tilde{\lambda} \circ \upsilon = \lambda \in \sC_{\varphi}$ for a fibre-preserving diffeomorphism $\upsilon = \eqref{eq: def_fibre_preserving_diffeo}$. 
\newline 

The situation becomes much simpler when the canonical relation $\varGamma$ is the graph of a symplectomorphism from $\dotCoTanN$ to $\dotCoTanM$. 
In this case $\varGamma$ is a symplectic manifold with respect to the symplectic form 
(see e.g.~\cite[p. 25]{Hoermander_Springer_2009}) 
\begin{equation} \label{eq: symplectic_form_graph_symplecto}
	\sigmaup_{\ms \varGamma} := \pr_{\ms M}^{*} \sigmaup_{\ms M} = \pr_{\ms N}^{*} \sigmaup_{\ms N}
\end{equation}
and we have the corresponding Liouville volume form $\fv_{\ms \varGamma} := \eqref{eq: def_Liouville_form}$. 
\newline 

Therefore, one has the bundle $\halfDen C \to C$ of half-densities over a homogeneous canonical relation $C$, whose local sections are $\sqrt{|\dVol_{\ms C}|}$.  
%
%
%
%
%
%
%
%
%
%
\section{Composition of densities on canonical relations}
\subsection{Half-densities} 
\label{sec: composition_halfdensity_canonical_relation}
To define a half-density on the canonical relation $C \circ_{e} \varLambda \subset \dotCoTanM \times \dotCoTanN$ (Definition~\ref{def: clean_composition}),  given half-densities on closed conic canonical relations $C \subset \dotCoTanM \times \dotCoTan O$ and $\varLambda \subset \dotCoTan O \times \dotCoTanN$,  we note that~\cite[Lem. 5.2]{Duistermaat_InventMath_1975} 
(see also, e.g.~\cite[Thm. 7.1.1, $(7.4)$]{Guillemin_InternationalP_2013},~\cite[Thm. 21.6.7]{Hoermander_Springer_2007})  
\begin{equation} \label{eq: isom_tensor_halfdensity_clean_composition}
    \halfDen_{(x, \xi; z, \zeta)} C \otimes \halfDen_{(z, \zeta; y, \eta)} \varLambda \cong \halfDen_{(x, \xi; z, \zeta; y, \eta)} \kF \otimes \halfDen_{\xxiyeta} (C \circ_{e} \varLambda),  
\end{equation} 
where $\kF$ is the set of fibres $\kF_{\xxiyeta} := \eqref{eq: def_fibre_projection_clean_composition}$ of the fibration $\varPi := \eqref{eq: varPi}$. 
This entails that by taking tensor product between half-densities on $C$ and $\varLambda$ followed by intersecting with the diagonal gives a half-density on $\kF$ times a half-density on $C \circ_{e} \varLambda$, and therefore, integrating over the fibre $\kF_{\xxiyeta}$ we achieve the desired  composition of half-densities on canonical relations~\cite[p. 64]{Duistermaat_InventMath_1975}:    
\begin{eqnarray} \label{eq: def_composition_halfdensity_canonical_relation}
	\cdot \diamond_{e} \cdot &:& 
    C^{\infty} \big( C; \halfDen C \big) \times C^{\infty} \big( \varLambda; \halfDen \varLambda \big) \to C^{\infty} \big( C \circ_{e} \varLambda; \halfDen (C \circ_{e} \varLambda) \big), 
	\nonumber\\ 
	&& 
    (\mu, \nu) \mapsto \mu \diamond_{e} \nu \, \xxiyeta := \int_{\kF_{\xxiyeta}} (\mu \otimes \nu) \big( \varPi^{-1} \xxiyeta \big). 
\end{eqnarray}
If $e = 0$ then the isomorphism~\eqref{eq: isom_tensor_halfdensity_clean_composition} reduces to~\cite[p. 179]{Hoermander_ActaMath_1971}  
(see also, e.g.~\cite[$(7.3)$]{Guillemin_InternationalP_2013}) 
\begin{equation}
	\halfDen C \otimes \halfDen \varLambda \cong \halfDen (C \circ \varLambda)
\end{equation}
and hence~\eqref{eq: def_composition_halfdensity_canonical_relation} becomes simply the pointwise evaluation 
(see e.g.~\cite[Thm. 4.2.2]{Duistermaat_Birkhaeuser_2011}) 
\begin{equation}
	\mu \diamond \nu \, \xxiyeta = \sum_{(z, \zeta) \,|\, (x, \xi; z, \zeta) \in C, (z, \zeta; y, \eta) \in \varLambda} \mu (x, \xi; z, \zeta) \, \nu (z, \zeta; y, \eta).  
\end{equation} 
We call the pairs $\xxi, \yeta$ (resp. $(z, \zeta)$) in the preceding equation as the output (rep. input) variables of the composition following the terminology from~\cite[Lem. 8.4]{Strohmaier_AdvMath_2021}. 
%
%
%
%
%
%
%
%
%
%
\subsection{Keller-Maslov bundle-valued half-densities} 
\label{sec: composition_halfdensity_Maslov_canonical_relation}
The principal symbol (Definition~\ref{def: symbol_FIO_mf}) of a scalar Lagrangian distribution associated with a canonical relation $C$ is a Keller-Maslov bundle $\bbL \to C$ valued density on $C$.  
Hence, we are compelled to extend the composition~\eqref{eq: def_composition_halfdensity_canonical_relation} on sections of $\bbL \otimes \halfDen C$. 
The first step, thus, is to construct the Keller-Maslov bundle $\Maslov \to C \circ_{e} \varLambda$ given $\bbL \to C$ and $\T \to \varLambda$. 
This is achieved by utilising the phase function $\psi := \eqref{eq: def_phase_function_clean_composition_canonical_relatation}$ into~\eqref{eq: def_section_Maslov_bundle} and~\eqref{eq: def_transition_function_Maslov_bundle}, 
as it is locally constant on $\kF$.   
Next, one employs the isomorphism~\cite[$(5.7)$]{Duistermaat_InventMath_1975} 
(see also, e.g.~\cite[$(5.29)$]{Guillemin_InternationalP_2013})   
\begin{equation} \label{eq: isom_tensor_Maslov_bundle}
    \bbL \boxtimes \T := \varPi_{\ms C}^{*} \bbL \otimes \varPi_{\ms \varLambda}^{*} \T \cong \varPi^{*} \Maslov,  
\end{equation}
where $\varPi_{\ms C}, \varPi_{\ms \varLambda} : C {\scriptsize \text{\FiveStarOpen}} \varLambda \to C, \varLambda$ are the projections, in order to make the identification  
\begin{equation}
    (\varPi^{*} \m) (x, \xi; z, \zeta; z, \zeta; y, \eta) = (\varPi_{\ms C}^{*} \bbl \otimes \varPi_{\ms \varLambda}^{*} \bbt) (x, \xi; z, \zeta; z, \zeta; y, \eta) 
\end{equation}
for any given sections $\bbl$ resp. $\bbt$ of $\bbL$ resp. $\T$ with a section $\m$ of $\M$. 
The desired extension of~\eqref{eq: def_composition_halfdensity_canonical_relation} is then  
\begin{equation} \label{eq: def_composition_halfdensity_Maslov_canonical_relation}
    C^{\infty} \big( C; \bbL \otimes \varOmega^{\frac{1}{2}} C \big) 
	\times 
	C^{\infty} \big( \varLambda; \T \otimes \varOmega^{\frac{1}{2}} \varLambda \big) 
	\to 
    C^{\infty} \big( C \circ_{e} \varLambda; \Maslov \otimes \varOmega^{\frac{1}{2}} (C \circ_{e} \varLambda) \big), 
    (\bbl \otimes \mu, \bbt \otimes \nu) \mapsto \m \otimes (\mu \diamond_{e} \nu).  
\end{equation}
%
%
%
%
%
%
%
%
%
%
\subsection{Homomorphism bundle-valued half-densities} 
\label{sec: composition_halfdensity_Maslov_vector_bundle_canonical_relation}
Finally, as in the set-up of Section~\ref{sec: algebra_FIO}, we consider any vector bundles $\sE \to M, \sG \to O, \sF \to N$ and so $\Hom{\sG, \sE} \to M \times O$ is a homomorphism bundle, forbye, $\widetilde{\mathrm{Hom}} (\sG, \sE) \to C$ represents the pullback of $\Hom{\sG, \sE}$ to $\coTan (M \times O)$ followed by restriction to $C$. 
One then extends~\eqref{eq: def_composition_halfdensity_Maslov_canonical_relation} by tensoring with $\widetilde{\mathrm{Hom}} (\sG, \sE), \widetilde{\mathrm{Hom}} (\sF, \sG), \widetilde{\mathrm{Hom}} (\sF, \sE)$, respectively to achieve our final composition rule~\cite[$(25.2.10)$]{Hoermander_Springer_2009}  
\begin{eqnarray} \label{eq: def_composition_halfdensity_Maslov_vector_bundle_canonical_relation}
	&& 
	C^{\infty} \big( C; \bbL \otimes \halfDen C \otimes \widetilde{\mathrm{Hom}} (\sG, \sE) \big)
	\times 
	C^{\infty} \big( \varLambda; \T \otimes \halfDen \varLambda \otimes \widetilde{\mathrm{Hom}} (\sF, \sG) \big) 
	\nonumber\\ 
    && \hspace*{6.0cm} \to C^{\infty} \big( C \circ_{e} \varLambda; \Maslov \otimes \halfDen (C \circ_{e} \varLambda) \otimes \widetilde{\mathrm{Hom}} (\sF, \sE) \big), 
	\nonumber\\ 
	&& 
	(\bbl \otimes a, \bbt \otimes b) \mapsto \m \otimes (a \diamond_{e} b), 
    \nonumber \\ 
    && 
    (a \diamond_{e} b) \xxiyeta := \int_{\kF_{\xxiyeta}} (a \boxtimes b) \big( \varPi^{-1} \xxiyeta \big).  
\end{eqnarray}
As before, the hindmost expression reduces to standard composition of homomorphisms when $e = 0$: 
\begin{equation}
	(a \diamond b) \xxiyeta = \sum_{(z, \zeta) \,|\, (x, \xi; z, \zeta) \in C, (z, \zeta; y, \eta) \in \varLambda} a (x, \xi; z, \zeta) \big( b (z, \zeta; y, \eta) \big). 
\end{equation}
%
%
%
%
%
%
%
%
%
%
\section{Literature}
The notions of 
\textit{Maslov index}~\cite[\S 7]{Maslov_Springer_1981}  
and 
\textit{Lagrangian submanifold}~\cite[Def. 4.25]{Maslov_Springer_1981} 
are originally due to Victor \textsc{Maslov} who also studied the local 
generating functions~\cite[Thm. 4.20]{Maslov_Springer_1981} 
for a Lagrangian manifold in some special cases. 
Vladimir \textsc{Arnol'd } put forwarded the former in the modern mathematical footing~\cite{Arnold_FAA_1967} 
while 
Lars \textsc{H\"{o}rmander} generalised the latter ideas by introducing the 
\textit{canonical relation}~\cite[Def. 4.1.2]{Hoermander_ActaMath_1971} 
and the  
\textit{non-degenerate phase function}~\cite[p. 91]{Hoermander_ActaMath_1971}. 
He gave the systematic route to locally generate a generic Lagrangian manifold by means of a phase function and obtained the 
necessarily and sufficient conditions~\cite[Thm. 3.1.16]{Hoermander_ActaMath_1971}  
for two phase functions to parametrise the same Lagrangian manifold. 
Moreover, he advanced the works by
Keller~\cite{Keller_AnnPhys_1958} 
and by 
Maslov~\cite[\S 10]{Maslov_Springer_1981} 
to give complete analytic description of the 
\textit{Keller-Maslov bundle}~\cite[Sec. 3.3]{Hoermander_ActaMath_1971}.  
Furthermore, the composition of Keller-Maslov bundle-valued half-densities on \textit{transversally intersecting} canonical relations was achieved by him~\cite[Thm. 4.2.2]{Hoermander_ActaMath_1971}. 
Subsequently, 
Johannes \textsc{Duistermaat} and Victor \textsc{Guillemin} envisaged the 
\textit{clean phase function}~\cite[p. 71]{Duistermaat_InventMath_1975} 
and propounded H\"{o}rmader's composition formula for \textit{cleanly intersecting} canonical relations.  
In contemporary of these developments a more geometric formulation was offered by 
Alan \textsc{Weinstein} in a series of 
papers~\cite{Weinstein_AdvMath_1971, Weinstein_Nice_1975, Weinstein_AMS_1977} 
(see also, e.g.~\cite{Meinrenken_ReptMathPhys_1992}) 
where the terms \textit{Lagrangian correspondence}, \textit{(Bott-)Morse family}, and \textit{stable equivalence} were coined for canonical relation, (clean) non-degenerate phase function, and equivalent class of phase functions for a given Lagrangian manifold, respectively.      
An in-depth discussion on these aspects of symplectic geometry is available, for instance, in the textbooks~\cite{Abraham_AMS_1978, Rudolph_Springer_2013}   
and we refer 
H\"{o}rmader's treatises~\cite{Hoermander_Springer_2007, Hoermander_Springer_2009}  
together with the expositions~\cite{Duistermaat_Birkhaeuser_2011, Treves_Plenum_1980, Guillemin_InternationalP_2013}
for the conic symplectic geometry.   
%
%
%

%% file: symbol.tex
\chapter{Symbols} 
\label{ch: symbol}
\textsf{We accumulate the background materials on a vector bundle-valued classical symbols on a conic manifold as required for this thesis}. 
\newline 

Throughout this chapter, we have assumed that $d \in \NO := \N \cup \{0\}, n \in \N, m \in \R$ and used  multi-indices $\alpha = (\alpha_{1}, \ldots, \alpha_{n}) \in \NO^{n}, \beta = (\beta_{1}, \ldots, \beta_{d}) \in \NO^{d}$ together with the multi-index notion for partial derivatives:    
\begin{equation}
	\partial_{\theta}^{\alpha} := \frac{\partial^{|\alpha|}}{ \partial \theta_{1}^{\alpha_{1}} \ldots \partial \theta_{n}^{\alpha_{n}} }, 
	\qquad 
    \rD_{x}^{\beta} := (- \ri)^{|\beta|} \frac{\partial^{|\beta|}}{ \partial (x^{1})^{\beta_{1}} \ldots \partial (x^{d})^{\beta_{d}} }, 
\end{equation} 
where $|\alpha| := |\alpha_{1}| + \ldots + |\alpha_{n}|$ and so on for $|\beta|$. 
Recall, the \textbf{Japanese bracket} is defined by  
\begin{equation} \label{eq: def_Japanese_bracket}
	\jbracket{\theta} :=  \sqrt{ 1 + \theta \cdot \theta}, 
\end{equation}
where $\cdot$ is the Euclidean inner product. 
%
%
%
%
%
%
%
%
%
%
\section{Symbols on a Euclidean space} 
\label{sec: symbol_Euclidean}

%
%
%
\begin{definition} \label{def: Kohn_Nirenberg_symbol_Euclidean}
	Let $U$ be an open subset of the Euclidean space $\Rd$.  
	The set $S_{1, 0}^{m} (U \times \Rn)$ of \textbf{Kohn-Nirenberg symbols} of \textbf{order} (at most) $m$ on $U \times \Rn$ is defined as the set of all $\C$-valued smooth functions $a \in C^{\infty} (U \times \Rn)$ on $U \times \Rn$ such that, for every compact set $K \subset U$, the estimation  
	(see e.g.~\cite[$(18.1.1)'$]{Hoermander_Springer_2007}) 
	\begin{equation} \label{eq: def_Kohn_Nirenberg_symbol_Euclidean}
		\forall x \in K, \forall \theta \in \Rn, \forall \alpha \in \NO^{n}, \forall \beta \in \NO^{d} : 
		\big| (\rD_{x}^{\beta} \partial_{\theta}^{\alpha} a) \xtheta \big| 
		\leq 
		\cst_{\alpha, \beta; K} \jbracket{\theta}^{m - |\alpha|} 
	\end{equation}
	holds true for some constant $\cst_{\alpha, \beta; K}$.   
	\newline 

	One sets 
	\begin{equation}
		S_{\bullet}^{\infty} 
		:= 
		\bigcup_{m \in \R} S_{\bullet}^{m}, 
		\quad 
		S_{\bullet}^{- \infty} 
		:= 
		\bigcap_{m \in \R} S_{\bullet}^{m}, 
		\quad 
		S_{\bullet}^{m - [k]} := S_{\bullet}^{m} / S_{\bullet}^{m - k}. 
	\end{equation}
\end{definition}
%
%
%

\begin{example}
    The set of all positively homogeneous smooth functions (\eqref{eq: def_positive_homogeneous_function}) on $U \times \Rn$ of order $m$ is automatically $S_{1, 0}^{m} (U \times \Rn)$.
\end{example}
%
%
%

We equip the set $S_{1, 0}^{m} (U \times \Rn)$ with the seminorm separating points defined as the best constant in~\eqref{eq: def_Kohn_Nirenberg_symbol_Euclidean}: 
\begin{equation} \label{eq: def_seminorm_symbol_U_Euclidean}
	\seminorm{\cdot}{m, N} : S_{1, 0}^{m} (U \times \Rn) \to \R_{\geq 0}, ~ a \mapsto 
	\seminorm{a}{m, N} 
	:= 
	\sup_{\xtheta \in K \times \Rn} \max_{|\alpha| + |\beta| \leq N}
	\frac{ \big| (\Dxbeta \partial_{\theta}^{\alpha} a) \xtheta \big| }{\jbracket{\theta}^{m - |\alpha|}},   
\end{equation}
which induces the metric 
\begin{equation} \label{eq: def_metric_symbol_U_Euclidean}
	\fd_{m} : S_{1, 0}^{m} (U \times \Rn) \times S_{1, 0}^{m} (U \times \Rn) \to \R_{\geq 0}, ~ (a, b) \mapsto \fd_{m} (a, b) 
	:= \sum_{N \in \NO} \frac{1}{2^{N}} \frac{ \seminorm{a - b}{m, N} }{ 1 + \seminorm{a - b}{ m, N} },   
\end{equation}
yielding  

%
%
%
\begin{proposition} \label{prop: Kohn_Nirenberg_symbol_Euclidean_Frechet_space}
	As in the terminologies of Definition~\ref{def: Kohn_Nirenberg_symbol_Euclidean}, $\big( S_{1, 0}^{m} (U \times \Rn), \fd_{m} :=~\eqref{eq: def_metric_symbol_U_Euclidean} \big)$ and $\big( S_{1, 0}^{\pm \infty} (U \times \Rn), \fd_{\pm \infty} \big)$ are Fr\'{e}chet spaces.  
	Furthermore, $S_{1, 0}^{m} (U \times \Rn) \hookrightarrow S_{1, 0}^{m'} (U \times \Rn)$ is a continuous inclusion for all $m \leq m'$, which is never dense if $m < m'$~\cite[Prop. 1.1.11]{Hoermander_ActaMath_1971}.  
\end{proposition}
%
%
%

The above symbol class is too general for the thesis and it turns out that the ``classical symbol'' class, a subclass of $S_{1, 0}^{m}$, is sufficient for our purpose. 
To introduce this subclass, we recall a notion related to the following completeness property of the space of symbols~\cite[Thm. 2.7]{Hoermander_AMS_1967}

%
%
%
\begin{definition} \label{def: asymptotic_summation}
	As in the terminologies of Definition~\ref{def: Kohn_Nirenberg_symbol_Euclidean}, suppose that $(a_{k})_{k \in \NO} \in S_{1, 0}^{m_{k}} (U \times \Rn)$ with $m_{0} > m_{1} > \ldots$ and $m_{k} \to - \infty$ as $k \to \infty$, and that $a \in S_{1, 0}^{\mu_{0}} (U \times \Rn)$ where $\mu_{\ms N} := \max_{k \geq N} m_{k}$. 
	Then $a$ is called the \textbf{asymptotic summation} of $a_{k}$:  
	\begin{equation} \label{eq: def_asymptotic_summation}
		a \sim \sum_{k \in \NO} a_{k} 
		\Leftrightarrow 
		\supp a \subset \bigcup_{k \in \NO} \supp (a_{k}), 
		\;
		\forall N \in \N : \left( a - \sum_{k=0}^{N-1} a_{k} \right) \in S_{1, 0}^{\mu_{\ms N}} (U \times \Rn).
	\end{equation} 
\end{definition}
%
%
%

\begin{remark}
	Asymptotic summation always exists. 
	Furthermore, it is unique modulo $S_{1, 0}^{-\infty} (U \times \Rn)$ and independent of any rearrangement of the series $\sum a_{k}$ 
	(see e.g.~\cite[Prop. 18.1.3]{Hoermander_Springer_2007}). 
\end{remark}
%
%
%

Now, our desired symbol class is formulated as follows. 

%
%
%
\begin{definition} \label{def: polyhomogeneous_symbol_Euclidean}
	Let $U$ be an open subset of the Euclidean space $\Rd$.  
	The set $S^{m} (U \times \Rn)$ of \textbf{polyhomogeneous symbols} of \textbf{order} (at most) $m$ on $U \times \Rn$ is defined as the set of all $a \in S_{1, 0}^{m} (U \times \Rn)$ which can be expressed as 
	(see e.g.~\cite[Def. 18.1.5]{Hoermander_Springer_2007}) 
	\begin{equation} \label{eq: def_polyhomogeneous_symbol}
		a \sim \sum_{k \in \N_{0}} a_{k} \chi_{k},    
	\end{equation}
    where, for each $k$, $a_{k}$ is positively homogeneous of degree $m - k$ when $|\theta| > 1$, $\sim$ denotes the asymptotic summation, and $\chi_{k} \in C^{\infty} (\Rn)$ vanishing identically near $\theta = 0$ with $\chi_{k} = 1$ if $|\theta| \geq 1$.
\end{definition}
%
%
%

As a subspace $S^{m} (U \times \Rn) \subset S_{1, 0}^{m} (U \times \Rn)$ is dense but not closed whereas the residual symbol space $S_{1, 0}^{- \infty} (U \times \Rn) = S^{- \infty} (U \times \Rn) \subset S^{m} (U \times \Rn)$ is closed for any $m$. 
Apart from the subspace topology, $S^{m} (U \times \Rn)$ can be endowed with a seminorm for each $a_{k}$ and a seminorm (cf.~\eqref{eq: def_seminorm_symbol_U_Euclidean}) ensuring asymptotic summability 
\begin{equation} \tag{\ref{eq: def_asymptotic_summation}'}
	\Big| \big(\rD_{x}^{\beta} \partial_{\theta}^{\alpha} (a - \sum_{k=0}^{N-1} a_{k}) \big) \xtheta \Big| \leq \cst_{\alpha, \beta; N} \jbracket{\theta}^{\mu_{\ms N} - |\alpha|}. 
\end{equation}
Thus, $S^{m} (U \times \Rn)$ is equipped with a countable number of seminorms and it becomes a \textit{Fr\'{e}chet space} with respect to the metric induced by these seminorms (cf.~\eqref{eq: def_metric_symbol_U_Euclidean}) 
(see e.g.~\cite[Sec. 2.16]{Melrose_2007}). 
\newline 

We observe, $\rD_{x}^{\beta} \partial_{\theta}^{\alpha} a_{k}$ is of degree $m - k - |\alpha|$ entailing $a_{k} \in S^{m - k} (U \times \Rn)$ provided that $a \in C^{\infty} (U \times \Rn)$ and it vanishes for large $x$. 
Therefore, away from $\theta = 0$, attributed by large $\theta$, \textit{polyhomogeneous symbols are essentially smooth functions with boundary}. 
\newline 

In order to formulate the principal symbol (Definition~\ref{def: symbol_FIO_Euclidean}) of a Fourier integral operator, we are compelled to extend the notion of symbols on a conic set. 

%
%
%
\begin{definition} \label{def: Kohn_Nirenberg_symbol_conic_Euclidean}
	Let $U \subset \Rd$ be an open set.  
	The set $S_{1, 0}^{m} (\cU)$ of \textbf{Kohn-Nirenberg symbols} of \textbf{order} (at most) $m$ on a conic set $\cU \subset U \times \Rn$ is defined as the set of all $a \in C^{\infty} (\cU)$ such that for every compact $\cK \subset \cU$ the estimation~\eqref{eq: def_Kohn_Nirenberg_symbol_Euclidean} holds true in $\{ (x; \lambdaup \theta) \,|\, \xtheta \in \cK \}$ where $\lambdaup \geq 1$  
	(see e.g.~\cite[p. 83]{Hoermander_Springer_2007}). 
\end{definition}
%
%
%

\begin{remark} \label{rem: symbol_conic_generalisation}
	Definition~\ref{def: Kohn_Nirenberg_symbol_conic_Euclidean} reproduces Definition~\ref{def: Kohn_Nirenberg_symbol_Euclidean} by making the choice $\cK := K \times \{ \theta \,|\, |\theta| \leq 1 \}$ where $K$ is as in the later definition 
    (see e.g.~\cite[p. 84]{Hoermander_Springer_2007}).  
	The asymptotic summability carries over $S_{1, 0}^{m} (\cU)$ and so we can define and topolise the corresponding polyhomogeneous class $S^{m} (\cU)$ analogously.  
    Also, by construction:  
    \begin{equation} \label{eq: deri_polyhomogeneous_symbol_Euclidean}
        \rD_{x}^{\beta} : S^{m} (U \times \Rn) \to S^{m} (U \times \Rn), 
        \quad 
        \partial_{\theta}^{\alpha} : S^{m} (U \times \Rn) \to S^{m - |\alpha|} (U \times \Rn). 
    \end{equation}
\end{remark}
%
%
%
%
%
%
%
%
%
%
\section{Symbols on manifolds} 
\label{sec: symbol_mf}
The principal symbol (Definition~\ref{def: symbol_PsiDO_mf}) and the subprincipal symbol  (Definition~\ref{def: subprincipal_symbol_mf}) of a pseudodifferential operator acting on a half-density bundle  (Section~\ref{sec: distribution_density_mf}) $\halfDen \to M$ over a manifold $M$, are invariantly defined homogeneous functions on the punctured cotangent bundle $\dotCoTanM$ as Remarked in~\ref{rem: symbol_mf_coordinate_change} and~\ref{rem: subprincipal_symbol_PsiDO_mf_coordinate_change}. 
Moreover, the principal symbol (Definition~\ref{def: symbol_FIO_mf}) of a Lagrangian distribution is a half-density on a homogeneous canonical relation (Definition~\ref{def: canonical_relation}) $C \subset \dotCoTanM \times \dotCoTanN$ where $N$ is a manifold whose dimension $\dN$ is not necessarily equal to that $\dM := d$ of $M$.  
As evident from the local case, it is sufficient to consider only the smooth functions on $\dotCoTanM$ for large covectors. 
Therefore, guided by the fact that $\dotCoTanM$ (Example~\ref{exm: cotangent_bundle_symplectic_mf}) and $C$ both are examples of conic manifolds (Definition~\ref{def: conic_mf}), together with Remark~\ref{rem: symbol_conic_generalisation}, we are going to describe symbols on the following class of conic manifolds. 

%
%
%
\begin{definition} \label{def: Kohn_Nirenberg_symbol_conic_mf}
    Let $(\dot{\sV}, \fm_{\lambdaup})$ be the conic manifold as described in Example~\ref{exm: punctured_real_vector_bundle_conic_mf}. 
	The set $S_{1, 0}^{m} (\dot{\sV})$ of \textbf{Kohn-Nirenberg symbols} of \textbf{order} (at most) $m$ on this conic manifold is defined as the set of all $a \in C^{\infty} (\dot{\sV})$ such that the functions $\lambdaup^{-m} \fm_{\lambdaup}^{*} a$ are uniformly bounded in $C^{\infty} (\dot{\sV})$ when $\lambdaup \geq 1$ 
	(see e.g.~\cite[p. 13]{Hoermander_Springer_2009},~\cite[Def. 2.1.2]{Duistermaat_Birkhaeuser_2011}). 
\end{definition}
%
%
%

One of the advantages of this definition is that it is well applicable for densities. 
In particular, if $\mu$ is a fixed positive half-density on $\dot{\sV}$ having degree of homogeneity $m'$ (Definition~\ref{def: conic_mf}) then every element $\fa \in S^{m} (\dot{\sV}; \halfDen \dot{\sV})$ can be written in the form 
\begin{equation} \label{eq: symbol_halfdensity_function}
	\fa = a \mu 
\end{equation}
for some scalar $a \in S^{m - m'} (\dot{\sV})$. 
See Example~\ref{exm: order_homogeneous_halfdensity} for an illustration of this concept. 
\newline 

Equivalently, a bottom-up approach can be adapted for which let us consider an atlas $\{ (U_{\alpha}, \kappa_{\alpha}) \}_{\alpha}$ for $M$ and a homogeneous symplectic atlas $\cA := \{ (\cU_{\alpha}, \varkappa_{\alpha}) \}_{\alpha}$ for $\dotCoTanM$ such that $\varkappa_{\alpha}: \cU_{\alpha} \cong \kappa_{\alpha} (U_{\alpha}) \times \dotRd$, as considered in Remark~\ref{rem: symbol_mf_coordinate_change}.  
Then   
\begin{equation} 
	S_{1, 0}^{m} (\cU_{\alpha}) = \big\{ a \in C^{\infty} (\cU_{\alpha}) \,|\, \forall (\cU_{\alpha}, \varkappa_{\alpha}) \in \cA : (\varkappa_{\alpha}^{*})^{-1} a \in S_{1, 0}^{m} \big( \kappa_{\alpha} (U_{\alpha}) \times \dotRd \big) \big\} 
\end{equation}
and consequently, $S_{1, 0}^{m} (\dotCoTanM)$ is constructed utilising the partition of unity. 
We remark that it is sufficient to satisfy these requirements for an atlas due to 

%
%
%
\begin{lemma} \label{lem: symbol_conic_change_variable}
	Let $U \subset \R^{d_{\ms U}}, V \subset \R^{d_{\ms V}}$ be open sets. 
	Suppose that $\cU \subset U \times \dot{\R}^{n_{\ms U}}, \cV \subset V \times \dot{\R}^{n_{\ms V}}$ are open conic sets and that $\kappaup : \cU \to \cV$ is a smooth proper map commuting with multiplication by positive scalars in the second variable. 
	If $a \in S^{m} (V \times \R^{n_{\ms V}})$ has support in the interior of a compactly based conic subset of $\cV$ then $a \circ \kappaup \in S^{m} (U \times \R^{n_{\ms U}})$ provided the composition is defined as $0$ outside $\cU$  
	(see e.g.~\cite[Lem. 25.1.6]{Hoermander_Springer_2007},~\cite[Def. 2.1.2]{Duistermaat_Birkhaeuser_2011}).  
\end{lemma}
%
%
%

The asymptotic summability carries over $S_{1, 0}^{m} (\dotCoTanM)$ 
(see e.g.~\cite[Prop. 2.1.2]{Duistermaat_Birkhaeuser_2011})  
and so we define $S^{m} (\dotCoTanM)$ analogously. 
It is now straightforward to extend the construction for the product manifold $\dotCoTanM \times \dotCoTanN \cong \dotCoTanMN$. 
Given a homogeneous canonical relation $(C, \iota)$ of $\dotCoTanM \times \dotCoTanN$, one achieves    
\begin{equation} \label{eq: symbol_canonical_relation}
	S^{m} (C) := \iota^{*} \, S^{m} (\dotCoTanM \times \dotCoTanN). 
\end{equation}
Since any open conic cover $\{ C_{\alpha, \beta} \}$ of $C$ is locally generated as elucidated in Section~\ref{sec: Lagrangian_submf}  
\begin{equation}
    \cU_{\alpha} \times \cV_{\beta} \xrightarrow{\ker \grad_{\xi, \eta} \psi} \sC_{\alpha, \beta} \xrightarrow{\eqref{eq: def_canonical_relation_identity_immersion}} C_{\alpha, \beta}
\end{equation}
by a non-degenerate phase function $\psi :=$~\eqref{eq: def_canonical_nondegenerate_phase_function_canonical_relation} on a homogeneous symplectic neighbourhood $\cU_{\alpha} \times \cV_{\beta}$ on $\dotCoTanM \times \dotCoTanN$, locally~\eqref{eq: symbol_canonical_relation} means to obtain $S^{m} (C_{\alpha, \beta})$ from $S^{m} (\cU_{\alpha} \times \cV_{\beta})$ via appropriate pullbacks of the preceding mapping.  
Finally, we incorporate $\halfDenC \otimes \Maslov$-valued symbols $\fa \in S^{m + (\dM + \dN) / 4} (C; \halfDenC \otimes \Maslov)$ simply by using~\eqref{eq: symbol_halfdensity_function} 
(see e.g.~\cite[Def. 18.2.10]{Hoermander_Springer_2007}):  
\begin{equation} \label{eq: symbol_canonical_relation_Maslov_density_function}
    \fa = a \, |\dVol_{\ms C}|^{1/2} \otimes \m, \quad a \in S^{m + (\dM + \dN - 2n) / 4} (C), 
\end{equation}
where $\dVol_{\ms C}$ is the natural volume form on $C$, as constructed in Section~\ref{sec: volume_canonical_relation}, and $\m$ is a section of the Keller-Maslov bundle $\Maslov \to C$ as detailed in Section~\ref{sec: Keller_Maslov_bundle}. 
The former resp. latter are of homogeneous of degree $n$ resp. $0$, as explained in Example~\ref{exm: order_homogeneous_halfdensity}. 
%
%
%
%
%
%
%
%
%
%
\section{Symbols on vector bundles} 
\label{sec: symbol}
One of the primary tools to address the subject matter of this thesis is the symbol calculus on vector bundles.   
Loosely speaking, this is achieved by replacing the $\C$-valued smooth maps by a  bundle-valued smooth maps in the discussion of Section~\ref{sec: symbol_mf} yet we spell out the details to supplement the theory of Fourier integral operators on vector bundles. 

%
%
%
\begin{definition} \label{def: Kohn_Nirenberg_symbol_conic_mf_matrix}
    Let $(\dot{\sV}, \fm_{\lambdaup})$ be a conic manifold as described in Example~\ref{exm: punctured_real_vector_bundle_conic_mf} and let $\C^{r \times k}$ denotes the set of all $r$ by $k$ matrices over $\C$.  
    The set $S_{1, 0}^{m} (\dot{\sV}, \C^{r \times k})$ of \textbf{$\C^{r \times k}$-valued Kohn-Nirenberg symbols} on $\dot{\sV}$ of \textbf{order} (at most) $m$ is defined as the set of all $a \in C^{\infty} (\dot{\sV}, \C^{r \times k})$ such that the functions $\lambdaup^{-m} \fm_{\lambdaup}^{*} a$ are uniformly bounded in $C^{\infty} (\dot{\sV}, \C^{r \times k})$ whenever $\lambdaup \geq 1$.   
\end{definition}
%
%
%

This means, for any $a \in C^{\infty} (\cU, \C^{r \times k})$ and for every compact set $\cK \subset \cU$, the estimation  
\begin{equation} \label{eq: def_Kohn_Nirenberg_symbol_Euclidean_matrix}
	\big\| (\rD_{x}^{\beta} \partial_{\theta}^{\alpha} a) \xtheta \big\|  
	\leq 
	\cst_{\alpha, \beta; \cK; \| \cdot \|} \jbracket{\theta}^{m - |\alpha|}  
\end{equation}
is valid for all $\xtheta \in \{ (x; \lambdaup \theta) | \xtheta \in \cK \}$ when $\lambdaup \geq 1$, for any $\alpha \in \NO^{n}, \beta \in \NO^{d}$ and for some constant $\cst_{\alpha, \beta; \cK; \| \cdot \|}$. 
Here $\| \cdot \|$ is some choice of norm on $\C^{r \times k}$ and this particular choice does \textit{not} matter 
(see e.g.~\cite[$(1.5.1.3)$]{Scott_OUP_2010}).   
\newline 

Since $C^{\infty} (\cU, \C^{r \times k})$ is isomorphic to the $r \times k$ matrix of elements  $C^{\infty} (\cU)$, we have $S_{1, 0}^{m} (\cU, \C^{r \times k}) \cong \big( S_{1, 0}^{m} (\cU) \big)_{r \times k}$ for any $m$ and any $r, k \in \N$. 
As in the scalar case, one can equip $S_{1, 0}^{m} (\cU, \C^{r \times k})$ with the Fr\'{e}chet topology and the asymptotic completeness can be shown.  
Consequently, we obtain the polyhomogeneous symbol class $S^{m} (\cU, \C^{r \times k})$ and 
(see e.g.~\cite[Sec. 1.5.2]{Scott_OUP_2010}) 
\begin{equation}
    \forall m \in \R, \forall r, k \in \N : S^{m} (\cU, \C^{r \times k}) \cong \big( S^{m} (\cU) \big)_{r \times k}. 
\end{equation}%

Our quintessential example is the homomorphism bundle $\Hom{\sF, \sE} \to M \times N$ (Section~\ref{sec: Lagrangian_distribution}), given a vector bundle $\sE \to M$ (resp. $\sF \to N$) over a manifold $M$ (resp. $N$), where neither $r := \rk \sE$ and $k := \rk \sF$ nor $\dM$ and $\dN$ are necessarily equal.  
Employing the canonical projections $\pi_{\ms M} : \coTanM \to M, \pi_{\ms N} : \coTanN \to N$, one has the induced homomorphism bundle $(\pi_{\ms M} \boxtimes \pi_{\ms N})^{*} \, \Hom{\sF, \sE} \to \coTanM \times \coTanN$ and the natural isomorphism 
\begin{equation}
	C^{\infty} \big( \coTanM \times \coTanN, \Hom{\sF, \sE} \big) \cong C^{\infty} \big( \coTanM \times \coTanN; (\pi_{\ms M} \boxtimes \pi_{\ms N})^{*} \Hom{\sF, \sE} \big), 
\end{equation}
which will be often exploited without further comments. 
\newline 

Let $\cA := \{ (\cU_{\alpha}, \varkappa_{\alpha}) \}_{\alpha}$ (resp. $\cB := \{ (\cV_{\beta}, \varrho_{\beta}) \}_{\beta}$) be a homogeneous symplectic atlas for $\dotCoTanM$ (resp. $\dotCoTanN$) such that the induced homomorphism bundle $(\pi_{\ms M} \boxtimes \pi_{\ms N})^{*} \, \Hom{\sF, \sE}$ over $\coTanM \times \coTanN$ admits an atlas $\kA := \{ (\cU_{\alpha} \times \cV_{\beta}, \vartheta_{\alpha, \beta}) \}_{\alpha, \beta}$, i.e., 
\begin{equation}
    \vartheta_{\alpha, \beta} : \big( (\pi_{\ms M} \boxtimes \pi_{\ms N})^{*} \, \Hom{\sF, \sE} \big)_{\alpha, \beta} \to \cU_{\alpha} \times \cV_{\beta} \otimes \C^{r \times k}
\end{equation}
is the induced trivialisation of $(\pi_{\ms M} \boxtimes \pi_{\ms N})^{*} \, \Hom{\sF, \sE}$ over $\cU_{\alpha} \times \cV_{\beta}$. 
Thereby
\begin{eqnarray}
    && 
	S^{m} \Big( \cU_{\alpha} \times \cV_{\beta}; \big( (\pi_{\ms M} \boxtimes \pi_{\ms N})^{*} \, \Hom{\sF, \sE} \big)_{\alpha, \beta} \Big) 
	\nonumber \\ 
	&& =   
	\Big\{ a \in C^{\infty} \Big( \cU_{\alpha} \times \cV_{\beta}; \big( (\pi_{\ms M} \boxtimes \pi_{\ms N})^{*} \, \Hom{\sF, \sE} \big)_{\alpha, \beta} \Big) 
	\,\Big|\, 
    \forall (\cU_{\alpha} \times \cV_{\beta}, \vartheta_{\alpha, \beta}) \in \kA : 
	\nonumber \\ 
	&& \qquad 
	(\vartheta_{\alpha, \beta}^{*})^{-1} a \in S^{m} \big( \cU_{\alpha} \times \cV_{\beta}, \C^{r \times k} \big) 
	\Big\} 
	\nonumber \\ 
	&& \cong   
    \big( S^{m} (\cU_{\alpha} \times \cV_{\beta}) \big)_{r \times k}. 
\end{eqnarray}
Assembling these local pieces together deploying the partition of unity, one then achieves $S^{m} \big( \dotCoTanM \times \dotCoTanN, \Hom{\sF, \sE} \big)$.  
\newline

For a canonical relation $\iota: C \to \dotCoTanM \times \dotCoTanN$, one sets 
\begin{eqnarray} \label{eq: symbol_HomFE_canonical_relation}
	S^{m} \big( C; \widetilde{\mathrm{Hom}} (\sF, \sE) \big) 
    & := & \iota^{*} S^{m} \big( \coTanM \times \coTanN, \Hom{\sF, \sE} \big), 
	\nonumber \\ 
	\widetilde{\mathrm{Hom}} (\sF, \sE) 
    & := & \iota^{*} \big( (\pi_{\ms M} \boxtimes \pi_{\ms N})^{*} \, \Hom{\sF, \sE} \big). 
\end{eqnarray}
To end, $S^{m} \big( C; \halfDenC \otimes \Maslov \otimes \widetilde{\mathrm{Hom}} (\sF, \sE) \big)$ is constructed as in~\eqref{eq: symbol_canonical_relation_Maslov_density_function}. 
%
%
%
%
%
%
%
%
%
%
\section{Literature}
The symbol class $S_{1, 0}^{m}$ was named after Joseph J. \textsc{Kohn} and Louis \textsc{Nirenberg} who introduced the class in their study of pseudodifferential operators~\cite{Kohn_CPAM_1965}. 
At any rate, this is not the most general possibility, for instance, the H\"{o}rmander symbol class~\cite{Hoermander_AMS_1967} is a vast generalisation of it, originated in the investigation of  fundamental solutions of 
hypoelliptic\footnote{A differential operator $P$ with smooth coefficients is called \textbf{hypoelliptic} if the equation $P u = f$ only has smooth solutions $u$ when $f$ is smooth.} 
operators of constant strength. 
We refer the monographs~\cite{Shubin_Springer_2001, Hoermander_Springer_2007, Scott_OUP_2010} 
and the expository lecture notes~\cite{Melrose_2007, vandenBan_2017}
for details.

%% file: QFT.tex
\chapter{Hadamard States} 
\label{ch: Hadamard_state}
\textsf{A bare minimum introduction of Feynman propagators by means of the time-ordered product is presented to bridge the quantum field theoretic approach and the microlocal formalism}. 
\newline 

In this chapter we will consider quantum field theory in a classical curved spacetime (Definition~\ref{def: spacetime}), also known as semi-classical gravity in the physics community.  
According to Einstein's general theory of relativity --- the best experimentally tested theory of gravity by far,  our universe is some $4$-dimensional Lorentzian manifold. 
On the other hand, quantum field theory in Minkowski spacetime has proven extremely successful by most experimental outcomes. 
Therefore, semi-classical gravity stems from the urge to extend the standard quantum field theoretic framework to a Lorentzian manifold, in physics terms, to include gravitational effects on quantum fields. 
It should be possible to derive semi-classical gravity by taking a suitable limit of a more fundamental theory wherein the spacetime metric is treated in accord with the principle of quantum physics. 
However, this has not been done --- except in formal and heuristic ways - simply because no present quantum theory of gravity has been developed to the point where such a well-defined limit can be taken. 
\newline 

Quantum field theory in curved spacetime is expected to provide an \textit{accurate description of quantum phenomena in a regime where the effects of spacetime curvature is significant yet the effects of quantum gravity may be neglected}. 
In particular, it is anticipated that this formalism should be applicable to the quantum phenomenology in the early universe and near (and inside of) black holes --- provided that one does not attempt to describe phenomena occurring so near to the singularities that the curvature reaches some extreme scale and the quantum nature of spacetime metric must be taken into account 
(see e.g. the review~\cite{Hollands_PR_2015} and references therein). 
\newline 

There are numerous ways to define a quantum theory and each approach comes with its merits and demerits; we refer the  
reviews~\cite{Ali_RMP_2005, Hollands_PR_2015, Fredenhagen_JMP_2016} 
for compendia of quantisation schemes and related aspects.    
In this thesis, we will use the locally covariant framework~\cite{Brunetti_CMP_2003} 
as it offers a mathematically precise and conceptually sound formulation of quantum field theory in curved spacetime.  
Loosely speaking, the essential idea is to take quantum observables as the fundamental object to define a quantum theory through the association of some appropriate algebra to every extended bounded open region of the spacetime, satisfying \textit{general covariance, Einstein causality}, and \textit{time-slice axiom}.   
These postulates are motivated on physical ground and can be interpreted as follows. 
The first requirement is a ramification of general relativity which tells that quantum fields should transform ``nicely'' under any coordinate change of spacetime geometry. 
Einstein causality ensures finite propagation speed to any physical information and the last demand provides the dynamics of the theory. 
To formulate all these notions precisely, it requires to introduce quite a bit of terminologies from operator algebra and Lorentzian geometry, which are beyond the scope of this thesis, and are not essential to comprehending the idea of Hadamard states and time-ordered products. 
Thus, we refer the 
expositions~\cite{Fewster_Springer_2015, Hollands_PR_2015, Fredenhagen_JMP_2016} 
for interested readers. 
For us, it is sufficient to keep in mind some unital $*$-algebra $\sA$ which satisfies these ``postulates''. 
Concrete characterisation of such an algebra actually defines a quantum field theory, as briefly discussed below. 
%
%
%
%
%
%
%
%
%
%
\section{Algebra of quantum fields}
\label{sec: field_algebra}
Let $\sE \to \sM$ be a vector bundle over a globally hyperbolic spacetime (Definition~\ref{def: globally_hyperbolic_spacetime}) $\spacetime$ and let $\bar{\sE}^{*}$ be the conjugate algebraic-dual bundle (Section~\ref{sec: adjoint}) of $\sE$. 
A \textit{linear quantum field} is some algebra (called the \textbf{field algebra}) $\sF (\sM; \sE)$-valued distribution 
\begin{equation}
	\Phi : \comSecsME \to \sF (\sM; \sE), ~u \mapsto \Phi (u), 
\end{equation}
where $\sF (\sM; \sE)$ can be taken as the algebra $\CCR (\sM; \sE)$ (resp. $\CAR (\sM; \sE)$) of canonical commutation (resp. anticommutation) relations, as defined below.  
The algebra $\CCR (\sM; \sE)$ (resp. $\CAR (\sM; \sE)$) is a unital $*$ (resp. $C^{*}$)-algebra generated over $\C$ by the symbols $\one_{-} := \one_{\CCR}$ (resp. $\one_{+} := \one_{\CAR}$), $\Phi (u)$ and its conjugate algebraic-dual $\bar{\Phi}^{*} (\phi)$, modulo the relations  
\begin{subequations}
	\begin{eqnarray}
		&& 
		\Phi (c u) = c \, \Phi (u), 
		\quad 
		\bar{\Phi}^{*} (c \phi) = \bar{c} \, \Phi^{*} (\phi), 
		\\ 
		&& 
		\Phi (L u) = 0 = \bar{\Phi}^{*} (\bar{L}^{*} \phi), 
		\label{eq: EOM}
		\\ 
		&& 
		[\Phi ([u]), \bar{\Phi}^{*} ([\phi])]_{\mp} = \ri \hbar \, \fG ([\phi] \boxtimes [u]) \one_{\mp}
		\label{eq: CCR_CAR}
	\end{eqnarray}
\end{subequations}
for any $[u] \in \comSecsME / \img{L}$ and $[\phi] \in \comSecsMEConjStar / \img{\bar{L}^{*}}$. 
Here $\hbar$ is the reduced Planck constant and it is $\hbar = 1$ in the \textit{natural units} which we use throughout the thesis. 
The brackets $[\cdot, \cdot]_{\mp}$ are the \textbf{commutator} and the \textbf{anticommutator}, respectively:  
\begin{equation}
	[\Phi (u), \bar{\Phi}^{*} (\phi) ]_{\mp} := \Phi (u) \, \bar{\Phi}^{*} (\phi) \mp \bar{\Phi}^{*} (\phi) \, \Phi (u)   
\end{equation}
and $\fG \in \cD' (\sM \times \sM; \sE \boxtimes \bar{\sE}^{*})$ is the Schwartz kernel of the causal propagator (Definition~\ref{def: Green_op}) for the differential operator $L$ depending on the specific model of $\Phi$, as listed in Table~\ref{tab: EOM_QFT}. 
Quantum fields are called \textit{Bosonic} (resp. \textit{Fermionic}) if $\sF (\sM; \sE)$ is $\CCR (\sM; \sE)$ (resp. $\CAR (\sM; \sE)$). 
%
%
%

\begin{table}
    \centering
	\begin{tabular}{ccc}
		\toprule
		Spin of a quantum field & Equation of motion & Quantisation 
		\\ \hline 
		$0$ & Klein-Gordon (Example~\ref{exm: covariant_Klein_Gordon_op}) & \cite{Dimock_CMP_1980, Kay_PR_1991}
		\\ 
		$1/2$ & Dirac (Example~\ref{ex: spin_Dirac_op}) & \cite{Dimock_AMS_1982, Hollands_CMP_2001, Dappiaggi_RMP_2009, Sanders_RMP_2010, Baer_Springer_2012} 
		\\ 
        $1$ & & 
		\\ 
		Abelian massive & Proca (Example~\ref{exm: Proca_op}) & \cite{Schambach_ReptMathPhys_2018}
		\\ 
		Abelian massless & Maxwell & \cite{Pfenning_CQG_2009, Sanders_CMP_2014, Finster_AHP_2015}
		\\ 
		Non-abelian massless & Yang-Mills & \cite{Hollands_RMP_2008, Fredenhagen_CMP_2013}
		\\ 
		$2$ & Einstein & 
		\\ 
		$3/2$ & Rarita-Schwinger (Section~\ref{sec: Rarita_Schwinger_op}) & \cite{Hack_PLB_2013}
		\\ \bottomrule
	\end{tabular}
	\caption[Models of QFT]{Models 
        of generally covariant quantum field theories. 
        Although we have not discussed Maxwell equation and Yang-Mills equation in this thesis but they can be related with the Hodge-d'Alembertian (Example~\ref{exm: Hodge_d_Alembert_op}) by a suitable ``gauge choice''.   
    } 
	\label{tab: EOM_QFT}
\end{table}
%
%
%

The field algebra is not necessarily the observable algebra because not all quantum fields are necessarily observables, for instance, the Fermi fields, as they do not commute for spacelike separated regions in $\sM$. 
Nevertheless, their appropriate bilinear combinations are observables. 
We will denote the algebra of observables by $\sA (\sM; \sE)$. 
For instance, $\CAR (\sM; \sE) \nsubset \sA (\sM; \sE)$ but $\CCR (\sM; \sE) \subset \sA (\sM; \sE)$. 
%
%
%
%
%
%
%
%
%
%
\section{Quantum states and Hadamard condition}
\label{sec: quantum_state}

%
%
%
\begin{definition} \label{def: algebraic_state}
	Let $\sE \to \sM$ be a vector bundle over a globally hyperbolic spacetime and $\sA (\sM; \sE)$ an algebra of quantum observables. 
    An \textbf{algebraic state} $\omega$ of $\sA (\sM; \sE)$ is a continuous (in a suitable topology) linear functional $\omega : \sA (\sM; \sE) \to \C$ which is positive: $\omega (A^{*} A) \geq 0$ for all $A \in \sA (\sM; \sE)$ and normalised: $\omega (\one) = 1$.
\end{definition}
%
%
%

The standard Hilbert space states used in physics literature are related to algebraic states via the  Gelfand-Najmark-Segal (GNS) theorem --- stating that 
(see e.g.~\cite[Thm. 2]{Drago_LMP_2020}), 
given an algebraic state $\omega$ of any unital $*$-algebra $\sA (\sM; \sE)$ there exists a quadruple $(\sH_{\omega}, \sD_{\omega}, \pi_{\omega}, \varOmega_{\omega})$ where $\sD_{\omega}$ is  a dense subspace of the GNS-Hilbert space $(\sH_{\omega}, \scalarProdTwo{\cdot}{\cdot})$, $\pi_{\omega}$ is an unital algebra representation of $\sA (\sM; \sE)$ on $\sD_{\omega}$ with the property $\pi_{\omega} (A^{*}) \subset \big( \pi_{\omega} (A) \big)^{\dag}$, $\dag$ being the adjoint with respect to $\scalarProdTwo{\cdot}{\cdot}$, and $\varOmega_{\omega}$ is the GNS-vector such that $\sD_{\omega} = \pi_{\omega} \big( \sA (\sM; \sE) \big) \varOmega_{\omega}$ and 
\begin{equation}
	\forall A \in \sA (\sM; \sE) : \omega (A) = \scalarProdTwo{\varOmega_{\omega}}{\pi_{\omega} (A) \varOmega_{\omega}}. 
\end{equation}
By definition, a state is specified by a collection of \textbf{Wightman $n$-point distributions} 
\begin{eqnarray}
	&& 
	\fW_{n} :\comSecsME^{\otimes n} \otimes \comSecsMEConjStar^{\otimes n} \to \C, ~u_{1} \otimes \ldots \otimes u_{n} \otimes \phi_{1} \otimes \ldots \phi_{n} \mapsto 
	\nonumber \\ 
	&& 
	\fW_{n} (u_{1}, \ldots, u_{n}, \phi_{1}, \ldots, \phi_{n}) := \omega \big( \Phi (u_{1}) \ldots \Phi (u_{n}) \, \bar{\Phi}^{*} (\phi_{1}) \ldots \bar{\Phi}^{*} (\phi_{n}) \big). 
\end{eqnarray}
By construction, the field equation~\eqref{eq: EOM}, the canonical commutation (resp. anticommutation) relation~\eqref{eq: CCR_CAR}, and positivity requirement of a state entail 
\begin{subequations}
	\begin{eqnarray}
		&& 
		\fW \big( L (u), \phi \big) = 0 = \fW \big( u, \bar{L}^{*} (\phi) \big) , 
		\\ 
		&& 
		\fW (u, \phi) - \fW (\phi, u) = \ri \fG (\phi \boxtimes u), 
		\\ 
		&& 
		\fW (\bar{u}^{*}, u) \geq 0  
	\end{eqnarray}
\end{subequations}
for the $2$-point Wightman distribution. 
\newline 

If $\spacetime$ admits a non-trivial isometry group $\Iso \sM$ preserving the orientation and the time-orientation, then $\sA (\sM; \sE)$ induces a representation of $\Iso \sM$ by elements $\alpha_{\ms \varXi} \in \Aut \big( \sA (\sM; \sE) \big)$. 
As a consequence, one has an induced automorphism $\alpha_{\varXi}^{*}$ on the states of $\sA (\sM; \sE)$. 
A state $\omega$ is called $\Iso \sM$-\textbf{invariant} if 
\begin{equation} \label{eq: def_invariant_state}
    \forall \alpha_{\varXi} \in \Aut \big( \sA (\sM; \sE) \big) : \alpha_{\varXi}^{*} \omega = \omega.
\end{equation}

%
%
%
\begin{example} \label{ex: Minkowski_vacuum}
    If $(\R^{4}, \upeta)$ is the $4$-dimensional Minkowski spacetime then the \textit{proper orthochronous Poincar\'{e} group} $\bbP_{0}^{\uparrow}$ is the spacetime isometry group (Example~\ref{exm: KVF_Minkowski}). 
    Let $\sE \to \R^{4}$ be a trivial $\R$-line bundle and set $L$ in~\eqref{eq: EOM} to be the \textit{Klein-Gordon operator} 
    \begin{equation}
        \square := - \frac{\partial^{2}}{\partial x^{2}} - \rrm^{2}, \quad \rrm^{2} \in \R_{+}
    \end{equation}
    with the mass-squared term $\rrm^{2}$ (cf. Example~\ref{exm: covariant_Klein_Gordon_op}).
	In physics language, this is called a \textbf{massive real scalar quantum field theory} on \textbf{Minkowski spacetime} whose dynamics is governed by the Klein-Gordon equation.  
    Since $(\R^{4}, \upeta)$ is globally hyperbolic, there exist unique retarded $\GreenKernelRet$ and advanced $\GreenKernelAdv$ Green's Schwartz kernels for $\square$ and hence the causal propagator is given by  
	(see e.g.~\cite[Thm. 6.2.3]{Hoermander_Springer_2003},~\cite[Sec. 4.3.4]{Strohmaier_Springer_2009}): 
	\begin{subequations}
		\begin{eqnarray}
			\fG_{x-y} & = & - \ri (\Fourier \mu_{+} - \Fourier \mu_{-}) (x-y), 
			\label{eq: Green_kernel_causal_Minkowski}
            \\ 
            \mu_{\pm} (\xi) & := & \thetaup (\pm \xi_{1}) \, \delta \big( \upeta (\xi, \xi) - \rrm^{2} \big), 
            \\ 
            \fG^{\ret, \adv}_{x-y} & = & \pm \uptheta \big( \pm (x^{1} - y^{1}) \big) \, \fG_{x-y}, 
			\label{eq: Green_kernel_ret_adv_Minkowski} 
			\\ 
			\GreenKernelAdv_{x-y} & = & \GreenKernelRet_{- (x-y)}, 
		\end{eqnarray}
	\end{subequations}
    where $\mu_{\pm} (\xi)$ are the unique (up to constant multiples) measures on the mass-hyperboloids $\{ \xi \in \R^{4} \,|\, \upeta (\xi, \xi) = \rrm^{2}, \xi_{1} \gtrless 0 \}$ invariant under the proper orthochronous Lorentz group  
	(see e.g.~\cite[Thm. IX. 33, IX. 37]{Reed_II}) 
    and the multiplication by the Heaviside step function  
    $\uptheta$ is well-defined due to the wavefront set (see~\eqref{eq: WF_Green_kernel_causal_Minkowski}) of $\fG$. 
    In standard physics literature, these measures are expressed as $\int_{\R^{4}} \delta \big( \upeta (\xi, \xi) - \rrm^{2} \big) \, \uptheta (\xi_{1} \gtrless 0) \, \rd \xi / (2 \pi)^{4}$ and preceding equations are presented as the 
	formal\footnote{These 
        integrals do not converge, so they must be understood as Fourier transforms of Schwartz distributions.} \label{foot: a}
	expressions 
	(see e.g.~\cite[Sec. 1-3-1]{Itzykson_Dover_2005}) 
	\newpage 
	\begin{subequations}
		\begin{eqnarray}
			\GreenKernelRet_{x-y} 
			& = & 
			- \ri \thetaup (x^{1} - y^{1}) \int_{\R^{3}} \left( \re^{- \ri \tE (x^{1} - y^{1})} - \re^{\ri \tE (x^{1} - y^{1})} \right) \re^{\ri \xi' \cdot (x'-y')} \frac{\rd \xi'}{(2 \pi)^{3} 2 \tE}, 
			\label{eq: Green_kernel_retarded_Minkowski_Fourier}
			\\ 
			\GreenKernelAdv_{x-y} 
			& = & 
			- \ri \thetaup (y^{1} - x^{1}) \int_{\R^{3}} \left( \re^{- \ri \tE (y^{1} - x^{1})} - \re^{\ri \tE (y^{1} - x^{1})} \right) \re^{\ri \xi' \cdot (y'-x')} \frac{\rd \xi'}{(2 \pi)^{3} 2 \tE}, 
			\label{eq: Green_kernel_advanced_Minkowski_Fourier}
			\\ 
			\fG_{x-y} 
			& = & 
			- \ri \int_{\R^{3}} \left( \re^{- \ri \tE (x^{1} - y^{1})} - \re^{\ri \tE (x^{1} - y^{1})} \right) \re^{\ri \xi' \cdot (x'-y')} \frac{\rd \xi'}{(2 \pi)^{3} 2 \tE}, 
			\label{eq: Green_kernel_causal_Minkowski_Fourier}
		\end{eqnarray}
	\end{subequations}
    where we have used the notation $x = (x^{1}, x'), \xi = (\xi_{1}, \xi') \in \R \times \R^{3}$ and $\tE := \sqrt{\xi' \cdot \xi' + \rrm^{2}}$ is the energy where $\cdot$ denotes the Euclidean inner product on $\R^{3}$.  
    \newline 

    In the Minkowski spacetime, there exists a unique $\bbP_{0}^{\uparrow}$-invariant state, the \textbf{vacuum} $\omega_{0}$ whose $2$-point Wightman distribution $\fW_{0}$ is given by 
    (see e.g.~\cite[Thm. IX. 34]{Reed_II})
	\begin{equation} \label{eq: Minkowski_vacuum_Fourier}
		\fW_{0} (x-y) 
		= - \ri (\Fourier \mu) (x-y) 
        = \int_{\R^{3}} \re^{- \ri \tE (x^{1} - y^{1})} \re^{\ri \xi' \cdot (x' - y')} \frac{\rd \xi'}{(2 \pi)^{3} 2 \tE}.  
	\end{equation}
    Evidently, the singularities of the Green's operators and $\fW_{0}$ are  
	(see e.g~\cite[p. 118]{Strohmaier_Springer_2009})
	\begin{subequations}
		\begin{eqnarray}
            \WFPrime \fG^{\substack{\ret \\ \adv}} 
			\! & = & \! 
            \big\{ (x, \xi; x, \xi) \in \dotCoTan \R^{4} \times \dotCoTan \R^{4} \big\} \bigcup \big\{ (x, \xi; y, \xi) \in \dot{\R}^{4} \times \dot{\R}^{4} \times \dot{\R}^{4} \times \dot{\R}^{4} \,|
			\nonumber \\ 
			&& 
			\upeta (x-y, x-y) = 0, x^{1} \gtrless y^{1}, \xi = \lambdaup (x-y), \lambdaup \in \R \big\}, 
			\label{eq: WF_Green_kernel_retarded_advanced_Minkowski}
			\\ 
			\WFPrime \fG 
			\! & = & \! 
			\left\{ \xxiyeta \in \dotCoTan \R^{4} \times \dotCoTan \R^{4} \,|\, \etaup (\xi, \xi) = 0, x-y = \lambdaup \xi, \lambdaup \in \R \right\}, 
			\label{eq: WF_Green_kernel_causal_Minkowski}
			\\ 
            \WFPrime \fW_{0} 
            \! & = & \! 
            \left\{ \xxiyeta \in \dotCoTan \R^{4} \times \dotCoTan \R^{4} \,|\, \etaup (\xi, \xi) = 0, x-y = \lambdaup \xi, \lambdaup \in \R, \xi_{1} > 0 \right\} \! . \qquad   
            \label{eq: WF_Minkowski_vacuum} 
		\end{eqnarray}
	\end{subequations}
\end{example}
%
%
%

However, there is \textit{no} analog of a vacuum state or even any other preferred state in a generic globally hyperbolic spacetime due to the lack of any non-trivial spacetime symmetry.  
There were numerous attempts to get around this problem and all suffered from significant limitations; see, for instance,  the 
review~\cite{Fredenhagen_JMP_2016} 
for a summary. 
In fact, under a few natural hypotheses, it has been proven that it is impossible to find any preferred state of a non-trivial locally covariant quantum field theory in a globally hyperbolic spacetime~\cite[Thm. 6.13]{Fewster_AHP_SPAS_2012} 
(see also the review~\cite[Thm. 1]{Fewster_IJMP_2018}). 
On the other hand, a generic algebraic state can be too singular to perform all the desired operations of a quantum field theory, e.g., defining non-linear observables like the quantum stress-energy tensor. 
Therefore, we utilise the preceding example as our guiding principle to select a physical state.  
\newline 

Since the tangent space of every spacetime is locally Minkowski and the singularity of a distribution is ``locally'' determined, we embrace the condition 
\begin{equation}
	\WF \fW = \WF \fW_{0}
\end{equation}
as a plausible criterion to select physical quantum states in any curved spacetime. 

%
%
%
\begin{definition} \label{def: Hadamard_state} 
	As in the terminologies of Definition~\ref{def: algebraic_state}, a state $\omega$ of $\sA (\sM; \sE)$ is called a \textbf{Hadamard state} if its $2$-point Wightman distribution has the following wavefront set~\cite{Radzikowski_CMP_1996}  
	\begin{equation}
		\WFPrime \fW = \{ \xxiyeta \in \coLightBun \times \coLightBun \,|\, \exists ! s \in \R_{\gtrless 0} : \xxi = \varPhi_{s} \yeta, \xi \vartriangleright 0 \}, 
	\end{equation} 
    where $\coLightBun \to \sM$ is the lightcone bundle over $\sM$, $\varPhi_{s}$ is the geodesic flow on the cotangent bundle and $\xi \vartriangleright 0$ means that $\xi$ is 
	future-directed\footnote{That is, $\xi (X) \geq 0$ for any futue-directed $X \in \tangent_{x} \sM$.}.  
\end{definition}
%
%
%

We note that Hadamard states are \textit{necessary} and \textit{sufficient} to define covariant non-linear quantum observables in any globally hyperbolic spacetime~\cite{Brunetti_CMP_2000, Hollands_CMP_2001, Hollands_CMP_2002, Khavkine_CMP_2016}  
(see also the reviews~\cite{Hollands_PR_2015, Fredenhagen_JMP_2016}). 
Therefore, only Hadamard states have been considered in this thesis. 

%
%
%
\begin{remark} \label{rem: nonunique_Hadamard_state}
	Hadamard states are \textit{never unique} as we have the freedom to add appropriate smooth terms. 
\end{remark}

%
%
%
%
%
%
%
%
%
%
\section{Time-ordered product} 
\label{sec: timeordered_product}
In order to construct non-linear quantum observables perturbatively one requires the notion of the time-ordered product between classical observables. 
Geometrically speaking, a classical field $\varphiup$ is a smooth section of $\sE$. 
A time-ordered product $\tT$ is a linear map from local and covariant polynomial expressions in $\varphiup$ to $\sF (\sM; \sE)$-valued distribution satisfying a number of postulates.     
The precise characterisation of those properties requires several new terminologies and, more importantly, the details are not essential to understand the subject matter of the ongoing discussion. 
Therefore, we abstain from providing the formal definition and refer~\cite{Hollands_CMP_2001} 
(see also the reviews~\cite{Hollands_PR_2015, Fredenhagen_Springer_2015} 
for details with earlier developments). 
It is sufficient for our purpose to note that 
\begin{subequations}
	\begin{eqnarray}
        && 
        \tT \big( \varphiup (x) \otimes \bar{\varphiup}^{*} (y) \big) 
        = \left\{ \begin{array}{cc}
                            \Phi_{x} \, \bar{\Phi}^{*}_{y}, & x \notin J^{-} (y), 
                            \\ 
                            \pm \bar{\Phi}^{*}_{y} \, \Phi_{x}, & x \notin J^{+} (y)
	                  \end{array} 
            \right. , 
        \\ 
        && 
        \tT \big( L \varphiup (x) \otimes \bar{\varphiup}^{*} (y) \big) = \ri \hbar \delta_{x, y} = \tT \big(  \varphiup (x) \otimes \bar{L}^{*} \bar{\varphiup}^{*} (y) \big), 
        \\ 
        && 
        \fW^{\tT}_{x, y} := \omega \left( \tT \big( \varphiup (x) \otimes \bar{\varphiup}^{*} (y) \big) \right), 
        \\ 
        && 
        \WF' \fW^{\tT} = C^{+}, 
    \end{eqnarray}
\end{subequations}
where $\fW^{\tT}$ is called the \textbf{time-ordered $2$-point Wightman distribution}, $C^{+}$ is the forward geodesic relation (Definition~\ref{def: geodesic_relation}) and $\pm$ is for Bosonic (resp. Fermionic) observables. 
The preceding two equations therefore entail that~\cite{Radzikowski_CMP_1996}  
\begin{equation} \label{eq: def_Green_kernel_Feynman_timeordered_product}
	\GreenKernelFeyn (u \boxtimes \phi) := - \ri \fW^{\tT} (u \boxtimes \phi)
\end{equation}
is the Feynman propagator (Definitions~\ref{def: fundamental_solution} and~\ref{def: Feynman_parametrix}). 

%
%
%
\begin{remark} \label{rem: state_depence_Feynman_propagator}
	The Feynman propagator \textit{depends} on the chosen quantum \textit{state} in contrast to the retarded, advanced, and causal propagators. 
\end{remark}
%
%
%

To comprehend this new propagator better, let us look again at the simplest case: the Klein-Gordon quantum field theory in Minkowski spacetime. 

%
%
%
\begin{example} \label{ex: Feynman_propagator_Minkowski}
	In the terminologies and notations introduced in Example~\ref{ex: Minkowski_vacuum}:  there is a unique Feynman propagator $\GreenKernelFeyn_{0}$ corresponding to the Minkowski vacuum: 
	\begin{equation}
		\GreenKernelFeyn_{0, x - y} = - \ri \big( \thetaup (x^{1} - y^{1}) \, (\Fourier \mu_{+}) (x-y) + \thetaup (y^{1} - x^{1}) \, (\Fourier \mu_{-}) (x-y) \big),  
	\end{equation}
	which is presented in the standard physics textbooks by\footnote{See Footnote~\ref{foot: a}} 
	(see e.g.~\cite[Sec. 1-3-1]{Itzykson_Dover_2005}) 
	\begin{equation} \label{eq: Green_kernel_Feynman_Minkowski_Fourier}
		\GreenKernelFeyn_{0, x-y} = - \ri \int_{\R^{3}} \left( \thetaup (x^{1} - y^{1}) \, \re^{- \ri \tE (x^{1} - y^{1})} + \thetaup (y^{1} - x^{1}) \, \re^{\ri \tE (x^{1} - y^{1})} \right) \re^{\ri \xi' \cdot (x'-y')} \frac{\rd \xi'}{(2 \pi)^{3} 2 \tE}  
	\end{equation}
    and its singularity structure is given by  
	(see e.g.~\cite[p. 118]{Strohmaier_Springer_2009}) 
	\begin{eqnarray} \label{eq: WF_Green_kernel_Feynman_Minkowski}
        \WFPrime \GreenKernelFeyn_{0} 
        & = & 
        \big\{ (x, \xi; x, \xi) \in \dotCoTan \R^{4} \times \dotCoTan \R^{4} \big\} \bigcup \big\{ (x, \xi; y, \xi) \in \dot{\R}^{4} \times \dot{\R}^{4} \times \dot{\R}^{4} \times \dot{\R}^{4} \,|
        \nonumber \\ 
		&& 
		\upeta (x-y, x-y) = 0, \xi = \lambdaup (x-y), \lambdaup \gtrless 0 \big\}.  
    \end{eqnarray}
\end{example}
%
%
%

These propagators have the following physical interpretation: 
\begin{itemize}
   	\item 
   	The retarded propagator~\eqref{eq: Green_kernel_retarded_Minkowski_Fourier} propagates positive ($\re^{- \ri \tE (x^{1} - y^{1})}$) and negative ($\re^{\ri \tE (x^{1} - y^{1})}$) frequencies forward in time, 
   	\item 
   	The advanced propagator~\eqref{eq: Green_kernel_advanced_Minkowski_Fourier} propagates positive and negative frequencies backward in time, 
   	\item 
   	The Feynman propagator~\eqref{eq: Green_kernel_Feynman_Minkowski_Fourier} propagates positive frequencies forward and negative frequencies backward in time, 
\end{itemize}
This formal expression of $\GreenKernelFeyn_{0}$ by means of the Fourier transform of $\lim_{\varepsilon \to 0} 1 / \big( \etaup (\xi, \xi) - \rrm^{2} + \ri \varepsilon \big)$ is often used in physics textbooks as the starting point to construct the Feynman propagator. 
However, neither the Fourier transformation nor the notation of a particle are available in an arbitrary globally hyperbolic spacetime. 
Thus, one cannot use this approach in a generic situation. 
Nevertheless, time-ordered product is well-defined in any globally hyperbolic spacetime and hence this can be used to define the Feynman propagators satisfactorily. 
One then needs to prove positivity of $\GreenKernelFeyn_{\omega}$ (cf. Proposition~\ref{prop: positivity_Feynman_minus_adv_NHOp}) in order to construct a Hadamard state $\omega$. 
Note, Feynman propagators are now \textit{no longer unique} as one has the freedom to add suitable smoothing terms without changing $C^{+}$.   
In other words,  $\GreenOpFeyn_{\omega}$ is \textit{unique} only up to those \textit{smoothing terms}, that is, there exists the notion of a \textit{unique Feynman parametrix} $\parametrixFeyn$ independent of any chosen $\omega$. 
\newline 

We can also view the one-to-one correspondence between $\omega$ and $\GreenKernelFeyn_{\omega}$ the way around. 
Since all Hadamard states have the same wavefront set, one can characterise the Feynman propagators in terms of the wavefront set as have been done in Definition~\ref{def: Feynman_parametrix}. 
The microlocal viewpoint has been adapted because it offers an intrinsic definition of the Feynman parametrix which has applications beyond quantum field theories as mentioned in Section~\ref{sec: intro_Feynman_propagator}.   
%
%
%
%
%
%
%
%
%
%
\section{Literature}
\label{sec: ref_QFT}
The idea of using the category of \textit{Hadamard states} as an appropriate replacement for the Minkowski vacuum and a crucial structural ingredient of any quantum field theory in curved spacetime was proposed by  
Bernard \textsc{Kay}~\cite[Sec. 1.1]{Ashtekar_Springer_1984} 
and by 
Rudolf \textsc{Haag}, Heide \textsc{Narnhofer}, and Ulrich \textsc{Stein}~\cite{Haag_CMP_1984}. 
Bernard \textsc{Kay} and Robert \textsc{Wald}~\cite{Kay_PR_1991} 
were the first to characterise these states precisely in terms of the short-distance behaviour of the $2$-point Wightman distribution. 
They had actually defined a local and a global Hadamard conditions. 
Recently, Valter \textsc{Moretti}~\cite{Moretti_LMP_2021}
located and closed a technical gap in their global condition.  
The modern definition of Hadamard states in the microlocal framework is due to 
Maciej \textsc{Radzikowski}~\cite{Radzikowski_CMP_1996}.  
Most of the earlier literature on Hadamard states was for scalar quantum field theories and the vector bundle generalisation developed 
subsequently~\cite{Koehler_PhD, Kratzert_AnnPhys_2000, Hollands_adiabatic_CMP_2001, Sahlmann_RMP_2001, Hollands_RMP_2008, Sanders_RMP_2010}.  
We refer, for instance, the 
exposition~\cite{Khavkine_Springer_2015}
for details. 
\newline 

A mathematically precise formalism of the \textit{time-ordered product} in curved spacetime was initiated by 
Romeo \textsc{Brunetti} and Klaus \textsc{Fredenhagen}~\cite{Brunetti_CMP_2000}, 
but there were a few ambiguities in their work.  
Employing the idea of local covariance of quantum observables, 
Stefan \textsc{Hollands} and Robert \textsc{Wald}  
have fixed those issues~\cite{Hollands_CMP_2001} 
and proven the existence of time-ordered product for scalar quantum field theory~\cite{Hollands_CMP_2002}. 
Later, the seminal result was extended by 
Stefan \textsc{Hollands}~\cite{Hollands_RMP_2008} 
for quantum Yang-Mills theory.

%% file: Gutzwiller_trace__onirban_islam.bbl
\begin{thebibliography}{236}%
\makeatletter
\providecommand \@ifxundefined [1]{%
 \@ifx{#1\undefined}
}%
\providecommand \@ifnum [1]{%
 \ifnum #1\expandafter \@firstoftwo
 \else \expandafter \@secondoftwo
 \fi
}%
\providecommand \@ifx [1]{%
 \ifx #1\expandafter \@firstoftwo
 \else \expandafter \@secondoftwo
 \fi
}%
\providecommand \natexlab [1]{#1}%
\providecommand \enquote  [1]{``#1''}%
\providecommand \bibnamefont  [1]{#1}%
\providecommand \bibfnamefont [1]{#1}%
\providecommand \citenamefont [1]{#1}%
\providecommand \href@noop [0]{\@secondoftwo}%
\providecommand \href [0]{\begingroup \@sanitize@url \@href}%
\providecommand \@href[1]{\@@startlink{#1}\@@href}%
\providecommand \@@href[1]{\endgroup#1\@@endlink}%
\providecommand \@sanitize@url [0]{\catcode `\\12\catcode `\$12\catcode
  `\&12\catcode `\#12\catcode `\^12\catcode `\_12\catcode `\%12\relax}%
\providecommand \@@startlink[1]{}%
\providecommand \@@endlink[0]{}%
\providecommand \url  [0]{\begingroup\@sanitize@url \@url }%
\providecommand \@url [1]{\endgroup\@href {#1}{\urlprefix }}%
\providecommand \urlprefix  [0]{URL }%
\providecommand \Eprint [0]{\href }%
\providecommand \doibase [0]{https://doi.org/}%
\providecommand \selectlanguage [0]{\@gobble}%
\providecommand \bibinfo  [0]{\@secondoftwo}%
\providecommand \bibfield  [0]{\@secondoftwo}%
\providecommand \translation [1]{[#1]}%
\providecommand \BibitemOpen [0]{}%
\providecommand \bibitemStop [0]{}%
\providecommand \bibitemNoStop [0]{.\EOS\space}%
\providecommand \EOS [0]{\spacefactor3000\relax}%
\providecommand \BibitemShut  [1]{\csname bibitem#1\endcsname}%
\let\auto@bib@innerbib\@empty
\bibitem [{\citenamefont {Kac}(1966)}]{Kac_AMM_1966}%
  \BibitemOpen
  \bibfield  {author} {\bibinfo {author} {\bibfnamefont {M.}~\bibnamefont
  {Kac}},\ }\bibfield  {title} {\bibinfo {title} {Can one hear the shape of a
  drum?},\ }\href {http://www.jstor.org/stable/2313748} {\bibfield  {journal}
  {\bibinfo  {journal} {Amer. Math. Monthly}\ }\textbf {\bibinfo {volume}
  {73}},\ \bibinfo {pages} {1} (\bibinfo {year} {1966})}\BibitemShut {NoStop}%
\bibitem [{\citenamefont {Uribe}(2000)}]{Uribe_Cuernavaca_1998}%
  \BibitemOpen
  \bibfield  {author} {\bibinfo {author} {\bibfnamefont {A.}~\bibnamefont
  {Uribe}},\ }\bibfield  {title} {\bibinfo {title} {Trace formulae},\ }in\
  \href {https://doi.org/10.1090/conm/260} {\emph {\bibinfo {booktitle} {First
  Summer School in Analysis and Mathematical Physics: Quantization, the
  {S}egal-{B}argmann Transform and Semiclassical Analysis}}},\ \bibinfo
  {editor} {edited by\ \bibinfo {editor} {\bibfnamefont {S.}~\bibnamefont
  {P\'{e}rez-Esteva}}\ and\ \bibinfo {editor} {\bibfnamefont {C.}~\bibnamefont
  {Villegas-Blas}}}\ (\bibinfo  {publisher} {American Mathematical Society},\
  \bibinfo {address} {Cuernavaca Morelos, Mexico},\ \bibinfo {year} {2000})\
  pp.\ \bibinfo {pages} {61 -- 90}\BibitemShut {NoStop}%
\bibitem [{\citenamefont {Levitan}(1952)}]{Levitan_1952}%
  \BibitemOpen
  \bibfield  {author} {\bibinfo {author} {\bibfnamefont {B.~M.}\ \bibnamefont
  {Levitan}},\ }\bibfield  {title} {\bibinfo {title} {On the asymptotic
  behavior of the spectral function of a self-adjoint differential equation of
  the second order},\ }\href {http://mi.mathnet.ru/eng/izv/v16/i4/p325}
  {\bibfield  {journal} {\bibinfo  {journal} {Izv. Akad. Nauk SSSR Ser. Mat.}\
  }\textbf {\bibinfo {volume} {16}},\ \bibinfo {pages} {325 } (\bibinfo {year}
  {1952})}\BibitemShut {NoStop}%
\bibitem [{\citenamefont {Levitan}(1955)}]{Levitan_1955}%
  \BibitemOpen
  \bibfield  {author} {\bibinfo {author} {\bibfnamefont {B.~M.}\ \bibnamefont
  {Levitan}},\ }\bibfield  {title} {\bibinfo {title} {On the asymptotic
  behavior of a spectral function and on expansion in eigenfunctions of a
  self-adjoint differential equation of second order.~{I}{I}},\ }\href
  {http://mi.mathnet.ru/eng/izv/v19/i1/p33} {\bibfield  {journal} {\bibinfo
  {journal} {Izv. Akad. Nauk SSSR Ser. Mat.}\ }\textbf {\bibinfo {volume}
  {19}},\ \bibinfo {pages} {33 } (\bibinfo {year} {1955})}\BibitemShut
  {NoStop}%
\bibitem [{\citenamefont {Avakumovi\v{c}}(1956)}]{Avakumovic_MathZ_1956}%
  \BibitemOpen
  \bibfield  {author} {\bibinfo {author} {\bibfnamefont {V.~G.}\ \bibnamefont
  {Avakumovi\v{c}}},\ }\bibfield  {title} {\bibinfo {title} {\"{U}ber die
  eigenfunktionen auf geschlossenen {R}iemannschen mannigfaltigkeiten},\ }\href
  {https://doi.org/10.1007/BF01473886} {\bibfield  {journal} {\bibinfo
  {journal} {Math. Z.}\ }\textbf {\bibinfo {volume} {65}},\ \bibinfo {pages}
  {327 } (\bibinfo {year} {1956})}\BibitemShut {NoStop}%
\bibitem [{\citenamefont {H{\"o}rmander}(1968)}]{Hoermander_ActaMath_1968}%
  \BibitemOpen
  \bibfield  {author} {\bibinfo {author} {\bibfnamefont {L.}~\bibnamefont
  {H{\"o}rmander}},\ }\bibfield  {title} {\bibinfo {title} {The spectral
  function of an elliptic operator},\ }\href
  {https://doi.org/10.1007/BF02391913} {\bibfield  {journal} {\bibinfo
  {journal} {Acta Math.}\ }\textbf {\bibinfo {volume} {121}},\ \bibinfo {pages}
  {193 } (\bibinfo {year} {1968})}\BibitemShut {NoStop}%
\bibitem [{\citenamefont {Ivrii}(2019)}]{Ivrii_I_Springer_2019}%
  \BibitemOpen
  \bibfield  {author} {\bibinfo {author} {\bibfnamefont {V.}~\bibnamefont
  {Ivrii}},\ }\href {https://doi.org/https://doi.org/10.1007/978-3-030-30557-4}
  {\emph {\bibinfo {title} {Microlocal Analysis, Sharp Spectral Asymptotics and
  Applications {I}: Semiclassical Microlocal Analysis and Local and Microlocal
  Semiclassical Asymptotics}}}\ (\bibinfo  {publisher} {Springer},\ \bibinfo
  {address} {Switzerland},\ \bibinfo {year} {2019})\BibitemShut {NoStop}%
\bibitem [{\citenamefont {Fang}(2017)}]{Fang_PhD}%
  \BibitemOpen
  \bibfield  {author} {\bibinfo {author} {\bibfnamefont {Y.-L.}\ \bibnamefont
  {Fang}},\ }\emph {\bibinfo {title} {Analysis of first order systems on
  manifolds without boundary: A spectral theoretic approach}},\ \href
  {https://discovery.ucl.ac.uk/id/eprint/1560979} {\bibinfo {type} {Ph{D}
  thesis}},\ \bibinfo  {school} {University College London} (\bibinfo {year}
  {2017})\BibitemShut {NoStop}%
\bibitem [{\citenamefont {Shubin}(2001)}]{Shubin_Springer_2001}%
  \BibitemOpen
  \bibfield  {author} {\bibinfo {author} {\bibfnamefont {M.}~\bibnamefont
  {Shubin}},\ }\href {https://doi.org/10.1007/978-3-642-56579-3} {\emph
  {\bibinfo {title} {Pseudodifferential Operators and Spectral Theory}}},\
  \bibinfo {edition} {2nd}\ ed.\ (\bibinfo  {publisher} {Springer-Verlag},\
  \bibinfo {address} {Berlin Heidelberg},\ \bibinfo {year} {2001})\BibitemShut
  {NoStop}%
\bibitem [{\citenamefont {Scott}(2010)}]{Scott_OUP_2010}%
  \BibitemOpen
  \bibfield  {author} {\bibinfo {author} {\bibfnamefont {S.}~\bibnamefont
  {Scott}},\ }\href {https://doi.org/10.1093/acprof:oso/9780198568360.001.0001}
  {\emph {\bibinfo {title} {Traces and Determinants of Pseudodifferential
  Operators}}},\ Oxford Mathematical Monographs\ (\bibinfo  {publisher} {Oxford
  University Press},\ \bibinfo {address} {UK},\ \bibinfo {year}
  {2010})\BibitemShut {NoStop}%
\bibitem [{\citenamefont {Duistermaat}\ and\ \citenamefont
  {Guillemin}(1975)}]{Duistermaat_InventMath_1975}%
  \BibitemOpen
  \bibfield  {author} {\bibinfo {author} {\bibfnamefont {J.~J.}\ \bibnamefont
  {Duistermaat}}\ and\ \bibinfo {author} {\bibfnamefont {V.~W.}\ \bibnamefont
  {Guillemin}},\ }\bibfield  {title} {\bibinfo {title} {The spectrum of
  positive elliptic operators and periodic bicharacteristics},\ }\href
  {https://doi.org/10.1007/BF01405172} {\bibfield  {journal} {\bibinfo
  {journal} {Invent. math.}\ }\textbf {\bibinfo {volume} {29}},\ \bibinfo
  {pages} {39} (\bibinfo {year} {1975})}\BibitemShut {NoStop}%
\bibitem [{\citenamefont {Grigis}\ and\ \citenamefont
  {Sj\"{o}strand}(1994)}]{Grigis_CUP_1994}%
  \BibitemOpen
  \bibfield  {author} {\bibinfo {author} {\bibfnamefont {A.}~\bibnamefont
  {Grigis}}\ and\ \bibinfo {author} {\bibfnamefont {J.}~\bibnamefont
  {Sj\"{o}strand}},\ }\href {https://doi.org/10.1017/CBO9780511721441} {\emph
  {\bibinfo {title} {Microlocal Analysis for Differential Operators: An
  Introduction}}},\ LMS Lecture Note Series\ (\bibinfo  {publisher} {Cambridge
  University Press},\ \bibinfo {year} {1994})\BibitemShut {NoStop}%
\bibitem [{\citenamefont {Gutzwiller}(1971)}]{Gutzwiller_JMP_1971}%
  \BibitemOpen
  \bibfield  {author} {\bibinfo {author} {\bibfnamefont {M.~C.}\ \bibnamefont
  {Gutzwiller}},\ }\bibfield  {title} {\bibinfo {title} {Periodic orbits and
  classical quantization conditions},\ }\href
  {https://doi.org/10.1063/1.1665596} {\bibfield  {journal} {\bibinfo
  {journal} {J. Math. Phys.}\ }\textbf {\bibinfo {volume} {12}},\ \bibinfo
  {pages} {343 } (\bibinfo {year} {1971})}\BibitemShut {NoStop}%
\bibitem [{\citenamefont {Gutzwiller}(1990)}]{Gutzwiller_Springer_1990}%
  \BibitemOpen
  \bibfield  {author} {\bibinfo {author} {\bibfnamefont {M.~C.}\ \bibnamefont
  {Gutzwiller}},\ }\href {https://doi.org/10.1007/978-1-4612-0983-6} {\emph
  {\bibinfo {title} {Chaos in Classical and Quantum Mechanics}}},\ \bibinfo
  {series} {Interdisciplinary Applied Mathematics}, Vol.~\bibinfo {volume} {1}\
  (\bibinfo  {publisher} {Springer-Verlag},\ \bibinfo {address} {New York},\
  \bibinfo {year} {1990})\BibitemShut {NoStop}%
\bibitem [{\citenamefont
  {Muratore-Ginanneschi}(2003)}]{Muratore-Ginanneschi_PR_2003}%
  \BibitemOpen
  \bibfield  {author} {\bibinfo {author} {\bibfnamefont {P.}~\bibnamefont
  {Muratore-Ginanneschi}},\ }\bibfield  {title} {\bibinfo {title} {Path
  integration over closed loops and {G}utzwiller's trace formula},\ }\href
  {https://doi.org/10.1016/S0370-1573(03)00212-6} {\bibfield  {journal}
  {\bibinfo  {journal} {Phys. Rep.}\ }\textbf {\bibinfo {volume} {383}},\
  \bibinfo {pages} {299 } (\bibinfo {year} {2003})},\ \Eprint
  {https://arxiv.org/abs/0210047} {arXiv:0210047 [nlin-CD]}\BibitemShut
  {NoStop}%
\bibitem [{\citenamefont {Colin~de
  Verdi\`{e}re}(1973{\natexlab{a}})}]{Verdiere_CompositioMathematica_1973}%
  \BibitemOpen
  \bibfield  {author} {\bibinfo {author} {\bibfnamefont {Y.}~\bibnamefont
  {Colin~de Verdi\`{e}re}},\ }\bibfield  {title} {\bibinfo {title} {Spectre du
  {L}aplacien et longueurs des g\'{e}od\'{e}siques p\'{e}riodiques. {I}},\
  }\href@noop {} {\bibfield  {journal} {\bibinfo  {journal} {Compositio
  Mathematica}\ }\textbf {\bibinfo {volume} {27}},\ \bibinfo {pages} {83}
  (\bibinfo {year} {1973}{\natexlab{a}})}\BibitemShut {NoStop}%
\bibitem [{\citenamefont {Colin~de
  Verdi\`{e}re}(1973{\natexlab{b}})}]{Verdiere_CompositioMathematica_1973_II}%
  \BibitemOpen
  \bibfield  {author} {\bibinfo {author} {\bibfnamefont {Y.}~\bibnamefont
  {Colin~de Verdi\`{e}re}},\ }\bibfield  {title} {\bibinfo {title} {Spectre du
  {L}aplacien et longueurs des g\'{e}od\'{e}siques p\'{e}riodiques. {I}{I}},\
  }\href@noop {} {\bibfield  {journal} {\bibinfo  {journal} {Compositio
  Mathematica}\ }\textbf {\bibinfo {volume} {27}},\ \bibinfo {pages} {159}
  (\bibinfo {year} {1973}{\natexlab{b}})}\BibitemShut {NoStop}%
\bibitem [{\citenamefont {Chazarain}(1974)}]{Chazarain_InventMath_1974}%
  \BibitemOpen
  \bibfield  {author} {\bibinfo {author} {\bibfnamefont {J.}~\bibnamefont
  {Chazarain}},\ }\bibfield  {title} {\bibinfo {title} {Formule de {P}oisson
  pour les vari{\'e}t{\'e}s {R}iemanniennes},\ }\href
  {https://doi.org/10.1007/BF01418788} {\bibfield  {journal} {\bibinfo
  {journal} {Invent. Math.}\ }\textbf {\bibinfo {volume} {24}},\ \bibinfo
  {pages} {65} (\bibinfo {year} {1974})}\BibitemShut {NoStop}%
\bibitem [{\citenamefont {Chazarain}(1980)}]{Chazarain_CPDE_1980}%
  \BibitemOpen
  \bibfield  {author} {\bibinfo {author} {\bibfnamefont {J.}~\bibnamefont
  {Chazarain}},\ }\bibfield  {title} {\bibinfo {title} {Spectre d'un
  hamiltonien quantique et mecanique classique},\ }\href
  {https://doi.org/10.1080/0360530800882148} {\bibfield  {journal} {\bibinfo
  {journal} {Commun. Partial Differ. Equ.}\ }\textbf {\bibinfo {volume} {5}},\
  \bibinfo {pages} {595 } (\bibinfo {year} {1980})}\BibitemShut {NoStop}%
\bibitem [{\citenamefont {Albeverio}\ \emph {et~al.}(1982)\citenamefont
  {Albeverio}, \citenamefont {Blanchard},\ and\ \citenamefont
  {H\o{}egh-Krohn}}]{Albeverio_CMP_1982}%
  \BibitemOpen
  \bibfield  {author} {\bibinfo {author} {\bibfnamefont {S.}~\bibnamefont
  {Albeverio}}, \bibinfo {author} {\bibfnamefont {P.}~\bibnamefont
  {Blanchard}},\ and\ \bibinfo {author} {\bibfnamefont {R.}~\bibnamefont
  {H\o{}egh-Krohn}},\ }\bibfield  {title} {\bibinfo {title} {Feynman path
  integrals and the trace formula for the schr\"{o}dinger operators},\ }\href
  {https://doi.org/https://doi.org/10.1007/BF01947071} {\bibfield  {journal}
  {\bibinfo  {journal} {Commun. Math. Phys.}\ }\textbf {\bibinfo {volume}
  {83}},\ \bibinfo {pages} {49 } (\bibinfo {year} {1982})}\BibitemShut
  {NoStop}%
\bibitem [{\citenamefont {Brummelhuis}\ and\ \citenamefont
  {Uribe}(1991)}]{Brummelhuis_CMP_1991}%
  \BibitemOpen
  \bibfield  {author} {\bibinfo {author} {\bibfnamefont {R.}~\bibnamefont
  {Brummelhuis}}\ and\ \bibinfo {author} {\bibfnamefont {A.}~\bibnamefont
  {Uribe}},\ }\bibfield  {title} {\bibinfo {title} {A semi-classical trace
  formula for {S}chr{\"o}dinger operators},\ }\href
  {https://doi.org/10.1007/BF02099074} {\bibfield  {journal} {\bibinfo
  {journal} {Commun. Math. Phys.}\ }\textbf {\bibinfo {volume} {136}},\
  \bibinfo {pages} {567 } (\bibinfo {year} {1991})}\BibitemShut {NoStop}%
\bibitem [{\citenamefont {Meinrenken}(1992)}]{Meinrenken_ReptMathPhys_1992}%
  \BibitemOpen
  \bibfield  {author} {\bibinfo {author} {\bibfnamefont {E.}~\bibnamefont
  {Meinrenken}},\ }\bibfield  {title} {\bibinfo {title} {Semiclassical
  principal symbols and {G}utzwiller's trace formula},\ }\href
  {https://doi.org/https://doi.org/10.1016/0034-4877(92)90019-W} {\bibfield
  {journal} {\bibinfo  {journal} {Rep. Math. Phys.}\ }\textbf {\bibinfo
  {volume} {31}},\ \bibinfo {pages} {279 } (\bibinfo {year}
  {1992})}\BibitemShut {NoStop}%
\bibitem [{\citenamefont {Paul}\ and\ \citenamefont
  {Uribe}(1995)}]{Paul_JFA_1995}%
  \BibitemOpen
  \bibfield  {author} {\bibinfo {author} {\bibfnamefont {T.}~\bibnamefont
  {Paul}}\ and\ \bibinfo {author} {\bibfnamefont {A.}~\bibnamefont {Uribe}},\
  }\bibfield  {title} {\bibinfo {title} {The semi-classical trace formula and
  propagation of wave packets},\ }\href
  {https://doi.org/https://doi.org/10.1006/jfan.1995.1105} {\bibfield
  {journal} {\bibinfo  {journal} {J. Funct. Anal.}\ }\textbf {\bibinfo {volume}
  {132}},\ \bibinfo {pages} {192 } (\bibinfo {year} {1995})}\BibitemShut
  {NoStop}%
\bibitem [{\citenamefont {Combescure}\ \emph {et~al.}(1999)\citenamefont
  {Combescure}, \citenamefont {Ralston},\ and\ \citenamefont
  {Robert}}]{Combescure_CMP_1999}%
  \BibitemOpen
  \bibfield  {author} {\bibinfo {author} {\bibfnamefont {M.}~\bibnamefont
  {Combescure}}, \bibinfo {author} {\bibfnamefont {J.}~\bibnamefont
  {Ralston}},\ and\ \bibinfo {author} {\bibfnamefont {D.}~\bibnamefont
  {Robert}},\ }\bibfield  {title} {\bibinfo {title} {A proof of the
  {G}utzwiller semiclassical trace formula using coherent states
  decomposition},\ }\href {https://doi.org/10.1007/s002200050591} {\bibfield
  {journal} {\bibinfo  {journal} {Commun. Math. Phys.}\ }\textbf {\bibinfo
  {volume} {202}},\ \bibinfo {pages} {463 } (\bibinfo {year} {1999})},\ \Eprint
  {https://arxiv.org/abs/9807005} {arXiv:9807005 [math-ph]}\BibitemShut
  {NoStop}%
\bibitem [{\citenamefont {H\"{o}rmander}(2009)}]{Hoermander_Springer_2009}%
  \BibitemOpen
  \bibfield  {author} {\bibinfo {author} {\bibfnamefont {L.}~\bibnamefont
  {H\"{o}rmander}},\ }\href {https://doi.org/10.1007/978-3-642-00136-9} {\emph
  {\bibinfo {title} {The Analysis of Linear Partial Differential Operators
  {I}{V}: {F}ourier Integral Operators}}},\ Classics in Mathematics\ (\bibinfo
  {publisher} {Springer-Verlag},\ \bibinfo {address} {Berlin, Heidelberg},\
  \bibinfo {year} {2009})\BibitemShut {NoStop}%
\bibitem [{\citenamefont {Safarov}\ and\ \citenamefont
  {Vassiliev}(1997)}]{Safarov_AMS_1997}%
  \BibitemOpen
  \bibfield  {author} {\bibinfo {author} {\bibfnamefont {Y.}~\bibnamefont
  {Safarov}}\ and\ \bibinfo {author} {\bibfnamefont {D.}~\bibnamefont
  {Vassiliev}},\ }\href {https://bookstore.ams.org/mmono-155} {\emph {\bibinfo
  {title} {The Asymptotic Distribution of Eigenvalues of Partial Differential
  Operators}}},\ \bibinfo {series} {Translations of Mathematical Monographs},
  Vol.\ \bibinfo {volume} {155}\ (\bibinfo  {publisher} {American Mathematical
  Society},\ \bibinfo {address} {USA},\ \bibinfo {year} {1997})\BibitemShut
  {NoStop}%
\bibitem [{\citenamefont {Guillemin}\ and\ \citenamefont
  {Sternberg}(2013)}]{Guillemin_InternationalP_2013}%
  \BibitemOpen
  \bibfield  {author} {\bibinfo {author} {\bibfnamefont {V.}~\bibnamefont
  {Guillemin}}\ and\ \bibinfo {author} {\bibfnamefont {S.}~\bibnamefont
  {Sternberg}},\ }\href
  {https://intlpress.com/site/pub/pages/books/items/00000409/reviews/index.html}
  {\emph {\bibinfo {title} {Semi-Classical Analysis}}}\ (\bibinfo  {publisher}
  {International Press of Boston, Inc.},\ \bibinfo {address} {US},\ \bibinfo
  {year} {2013})\BibitemShut {NoStop}%
\bibitem [{\citenamefont {Sandoval}(1999)}]{Sandoval_CPDE_1999}%
  \BibitemOpen
  \bibfield  {author} {\bibinfo {author} {\bibfnamefont {M.~R.}\ \bibnamefont
  {Sandoval}},\ }\bibfield  {title} {\bibinfo {title} {Wavetrace asymptotics
  for operators of {D}irac type},\ }\href
  {https://doi.org/10.1080/03605309908821487} {\bibfield  {journal} {\bibinfo
  {journal} {Commun. Partial. Differ. Equ.}\ }\textbf {\bibinfo {volume}
  {24}},\ \bibinfo {pages} {1903} (\bibinfo {year} {1999})}\BibitemShut
  {NoStop}%
\bibitem [{\citenamefont {Strohmaier}\ and\ \citenamefont
  {Zelditch}(2021{\natexlab{a}})}]{Strohmaier_AdvMath_2021}%
  \BibitemOpen
  \bibfield  {author} {\bibinfo {author} {\bibfnamefont {A.}~\bibnamefont
  {Strohmaier}}\ and\ \bibinfo {author} {\bibfnamefont {S.}~\bibnamefont
  {Zelditch}},\ }\bibfield  {title} {\bibinfo {title} {A {G}utzwiller trace
  formula for stationary space-times},\ }\href
  {https://doi.org/https://doi.org/10.1016/j.aim.2020.107434} {\bibfield
  {journal} {\bibinfo  {journal} {Adv. Math.}\ }\textbf {\bibinfo {volume}
  {376}},\ \bibinfo {pages} {107434} (\bibinfo {year} {2021}{\natexlab{a}})},\
  \Eprint {https://arxiv.org/abs/1808.08425} {arXiv:1808.08425 [math.AP]}\BibitemShut 
  {NoStop}%
\bibitem [{\citenamefont {Ivrii}(2016)}]{Ivrii_BullMathSci_2016}%
  \BibitemOpen
  \bibfield  {author} {\bibinfo {author} {\bibfnamefont {V.}~\bibnamefont
  {Ivrii}},\ }\bibfield  {title} {\bibinfo {title} {100 years of {W}eyl’s
  law},\ }\href {https://doi.org/10.1007/s13373-016-0089-y} {\bibfield
  {journal} {\bibinfo  {journal} {Bull. Math. Sci.}\ }\textbf {\bibinfo
  {volume} {6}},\ \bibinfo {pages} {379} (\bibinfo {year} {2016})},\ \Eprint
  {https://arxiv.org/abs/1608.03963} {arXiv:1608.03963 [math.SP]}\BibitemShut
  {NoStop}%
\bibitem [{\citenamefont {Avetisyan}\ \emph {et~al.}(2016)\citenamefont
  {Avetisyan}, \citenamefont {Fang},\ and\ \citenamefont
  {Vassiliev}}]{Avetisyan_JST_2016}%
  \BibitemOpen
  \bibfield  {author} {\bibinfo {author} {\bibfnamefont {Z.}~\bibnamefont
  {Avetisyan}}, \bibinfo {author} {\bibfnamefont {Y.-L.}\ \bibnamefont
  {Fang}},\ and\ \bibinfo {author} {\bibfnamefont {D.}~\bibnamefont
  {Vassiliev}},\ }\bibfield  {title} {\bibinfo {title} {Spectral asymptotics
  for first order systems},\ }\href {https://doi.org/10.4171/JST/137}
  {\bibfield  {journal} {\bibinfo  {journal} {J. Spectr. Theory}\ }\textbf
  {\bibinfo {volume} {6}},\ \bibinfo {pages} {695 } (\bibinfo {year} {2016})},\
  \Eprint {https://arxiv.org/abs/1512.06281} {arXiv:1512.06281 [math.SP]}\BibitemShut 
  {NoStop}%
\bibitem [{\citenamefont {Li}\ and\ \citenamefont
  {Strohmaier}(2016)}]{Li_JGP_2016}%
  \BibitemOpen
  \bibfield  {author} {\bibinfo {author} {\bibfnamefont {L.}~\bibnamefont
  {Li}}\ and\ \bibinfo {author} {\bibfnamefont {A.}~\bibnamefont
  {Strohmaier}},\ }\bibfield  {title} {\bibinfo {title} {The local counting
  function of operators of {D}irac and {L}aplace type},\ }\href
  {https://doi.org/https://doi.org/10.1016/j.geomphys.2016.02.006} {\bibfield
  {journal} {\bibinfo  {journal} {J. Geom. Phys.}\ }\textbf {\bibinfo {volume}
  {104}},\ \bibinfo {pages} {204 } (\bibinfo {year} {2016})},\ \Eprint
  {https://arxiv.org/abs/1509.00198} {arXiv:1509.00198 [math.SP]} \BibitemShut
  {NoStop}%
\bibitem [{\citenamefont {Haefliger}(1956)}]{Haefliger_1956}%
  \BibitemOpen
  \bibfield  {author} {\bibinfo {author} {\bibfnamefont {A.}~\bibnamefont
  {Haefliger}},\ }\bibfield  {title} {\bibinfo {title} {Sur l'{e}xtension du
  groupe structural d'une espace fibr\'{e}},\ }\href
  {https://gallica.bnf.fr/ark:/12148/bpt6k3195v/f562.image} {\bibfield
  {journal} {\bibinfo  {journal} {Comptes rendus de l'{A}cad\'{e}mie des
  {S}ciences (C. R. Acad. Sci. Paris)}\ }\textbf {\bibinfo {volume} {243}},\
  \bibinfo {pages} {558 } (\bibinfo {year} {1956})}\BibitemShut {NoStop}%
\bibitem [{\citenamefont {Borel}\ and\ \citenamefont
  {Hirzebruch}(1959)}]{Borel_AJM_1959}%
  \BibitemOpen
  \bibfield  {author} {\bibinfo {author} {\bibfnamefont {A.}~\bibnamefont
  {Borel}}\ and\ \bibinfo {author} {\bibfnamefont {F.}~\bibnamefont
  {Hirzebruch}},\ }\bibfield  {title} {\bibinfo {title} {Characteristic classes
  and homogeneous spaces, {I}{I}},\ }\href
  {http://www.jstor.org/stable/2372747} {\bibfield  {journal} {\bibinfo
  {journal} {Amer. J. Math.}\ }\textbf {\bibinfo {volume} {81}},\ \bibinfo
  {pages} {315} (\bibinfo {year} {1959})}\BibitemShut {NoStop}%
\bibitem [{\citenamefont {Milnor}(1963)}]{Milnor_Enseign_1963}%
  \BibitemOpen
  \bibfield  {author} {\bibinfo {author} {\bibfnamefont {J.}~\bibnamefont
  {Milnor}},\ }\bibfield  {title} {\bibinfo {title} {Spin structures on
  manifolds},\ }\href
  {https://www.e-periodica.ch/digbib/view?pid=ens-001%3A1963%3A9#292}
  {\bibfield  {journal} {\bibinfo  {journal} {Enseign. math.}\ }\textbf
  {\bibinfo {volume} {9}},\ \bibinfo {pages} {198 } (\bibinfo {year}
  {1963})}\BibitemShut {NoStop}%
\bibitem [{\citenamefont {Karoubi}(1968)}]{Karoubi_1968}%
  \BibitemOpen
  \bibfield  {author} {\bibinfo {author} {\bibfnamefont {M.}~\bibnamefont
  {Karoubi}},\ }\bibfield  {title} {\bibinfo {title} {Alg\`{e}bres de clifford
  et $k$-th\'{e}orie},\ }\href {https://doi.org/10.24033/asens.1163} {\bibfield
   {journal} {\bibinfo  {journal} {Annales scientifiques de l'\'{E}cole Normale
  Sup\'{e}rieure}\ }\textbf {\bibinfo {volume} {4e s{\'e}rie, 1}},\ \bibinfo
  {pages} {161} (\bibinfo {year} {1968})}\BibitemShut {NoStop}%
\bibitem [{\citenamefont {McCormick}(2022)}]{McCormick}%
  \BibitemOpen
  \bibfield  {author} {\bibinfo {author} {\bibfnamefont {A.}~\bibnamefont
  {McCormick}},\ }\href {https://doi.org/10.48550/ARXIV.2203.16729} {\bibinfo
  {title} {A trace formula on stationary {K}aluza-{K}lein spacetimes}}
  (\bibinfo {year} {2022}),\ \Eprint {https://arxiv.org/abs/2203.16729}
  {arXiv:2203.16729 [math-ph]}\BibitemShut {NoStop}%
\bibitem [{\citenamefont {Strohmaier}\ and\ \citenamefont
  {Zelditch}(2021{\natexlab{b}})}]{Strohmaier_IndagMath_2021}%
  \BibitemOpen
  \bibfield  {author} {\bibinfo {author} {\bibfnamefont {A.}~\bibnamefont
  {Strohmaier}}\ and\ \bibinfo {author} {\bibfnamefont {S.}~\bibnamefont
  {Zelditch}},\ }\bibfield  {title} {\bibinfo {title} {Semi-classical mass
  asymptotics on stationary spacetimes},\ }\href
  {https://doi.org/https://doi.org/10.1016/j.indag.2020.08.010} {\bibfield
  {journal} {\bibinfo  {journal} {Indag. Math.}\ }\textbf {\bibinfo {volume}
  {32}},\ \bibinfo {pages} {323} (\bibinfo {year} {2021}{\natexlab{b}})},\
  \bibinfo {note} {special Issue in memory of Hans Duistermaat},\ \Eprint
  {https://arxiv.org/abs/2002.01055} {arXiv:2002.01055 [math-ph]}\BibitemShut
  {NoStop}%
\bibitem [{\citenamefont {Strohmaier}\ and\ \citenamefont
  {Zelditch}(2021{\natexlab{c}})}]{Strohmaier_RMP_2021}%
  \BibitemOpen
  \bibfield  {author} {\bibinfo {author} {\bibfnamefont {A.}~\bibnamefont
  {Strohmaier}}\ and\ \bibinfo {author} {\bibfnamefont {S.}~\bibnamefont
  {Zelditch}},\ }\bibfield  {title} {\bibinfo {title} {Spectral asymptotics on
  stationary space-times},\ }\href {https://doi.org/10.1142/S0129055X20600077}
  {\bibfield  {journal} {\bibinfo  {journal} {Rev. Math. Phys.}\ }\textbf
  {\bibinfo {volume} {33}},\ \bibinfo {pages} {2060007} (\bibinfo {year}
  {2021}{\natexlab{c}})}\BibitemShut {NoStop}%
\bibitem [{\citenamefont {Duistermaat}\ and\ \citenamefont
  {H\"{o}rmander}(1972)}]{Duistermaat_ActaMath_1972}%
  \BibitemOpen
  \bibfield  {author} {\bibinfo {author} {\bibfnamefont {J.~J.}\ \bibnamefont
  {Duistermaat}}\ and\ \bibinfo {author} {\bibfnamefont {L.}~\bibnamefont
  {H\"{o}rmander}},\ }\bibfield  {title} {\bibinfo {title} {Fourier integral
  operators. {I}{I}},\ }\href {https://doi.org/10.1007/BF02392165} {\bibfield
  {journal} {\bibinfo  {journal} {Acta Math.}\ }\textbf {\bibinfo {volume}
  {128}},\ \bibinfo {pages} {183 } (\bibinfo {year} {1972})}\BibitemShut
  {NoStop}%
\bibitem [{\citenamefont {Brunetti}\ and\ \citenamefont
  {Fredenhagen}(2000)}]{Brunetti_CMP_2000}%
  \BibitemOpen
  \bibfield  {author} {\bibinfo {author} {\bibfnamefont {R.}~\bibnamefont
  {Brunetti}}\ and\ \bibinfo {author} {\bibfnamefont {K.}~\bibnamefont
  {Fredenhagen}},\ }\bibfield  {title} {\bibinfo {title} {Microlocal analysis
  and interacting quantum field theories: Renormalization on physical
  backgrounds},\ }\href {https://doi.org/10.1007/s002200050004} {\bibfield
  {journal} {\bibinfo  {journal} {Commun. Math. Phys.}\ }\textbf {\bibinfo
  {volume} {208}},\ \bibinfo {pages} {623} (\bibinfo {year} {2000})},\ \Eprint
  {https://arxiv.org/abs/9903028} {arXiv:9903028 [math-ph]}\BibitemShut
  {NoStop}%
\bibitem [{\citenamefont {Hollands}\ and\ \citenamefont
  {Wald}(2001)}]{Hollands_CMP_2001}%
  \BibitemOpen
  \bibfield  {author} {\bibinfo {author} {\bibfnamefont {S.}~\bibnamefont
  {Hollands}}\ and\ \bibinfo {author} {\bibfnamefont {R.~M.}\ \bibnamefont
  {Wald}},\ }\bibfield  {title} {\bibinfo {title} {Local {W}ick polynomials and
  time ordered products of quantum fields in curved spacetime},\ }\href
  {https://doi.org/10.1007/s002200100540} {\bibfield  {journal} {\bibinfo
  {journal} {Commun. Math. Phys.}\ }\textbf {\bibinfo {volume} {223}},\
  \bibinfo {pages} {289} (\bibinfo {year} {2001})},\ \Eprint
  {https://arxiv.org/abs/0103074} {arXiv:0103074 [gr-qc]}\BibitemShut
  {NoStop}%
\bibitem [{\citenamefont {Hollands}\ and\ \citenamefont
  {Wald}(2002)}]{Hollands_CMP_2002}%
  \BibitemOpen
  \bibfield  {author} {\bibinfo {author} {\bibfnamefont {S.}~\bibnamefont
  {Hollands}}\ and\ \bibinfo {author} {\bibfnamefont {R.~M.}\ \bibnamefont
  {Wald}},\ }\bibfield  {title} {\bibinfo {title} {Existence of local covariant
  time ordered products of quantum fields in curved spacetime},\ }\href
  {https://doi.org/10.1007/s00220-002-0719-y} {\bibfield  {journal} {\bibinfo
  {journal} {Commun. Math. Phys.}\ }\textbf {\bibinfo {volume} {231}},\
  \bibinfo {pages} {309} (\bibinfo {year} {2002})},\ \Eprint
  {https://arxiv.org/abs/0111108} {arXiv:0111108 [gr-qc]}\BibitemShut
  {NoStop}%
\bibitem [{\citenamefont {Radzikowski}(1996)}]{Radzikowski_CMP_1996}%
  \BibitemOpen
  \bibfield  {author} {\bibinfo {author} {\bibfnamefont {M.~J.}\ \bibnamefont
  {Radzikowski}},\ }\bibfield  {title} {\bibinfo {title} {Micro-local approach
  to the {H}adamard condition in quantum field theory on curved space-time},\
  }\href {http://projecteuclid.org/euclid.cmp/1104287114} {\bibfield  {journal}
  {\bibinfo  {journal} {Commun. Math. Phys.}\ }\textbf {\bibinfo {volume}
  {179}},\ \bibinfo {pages} {529} (\bibinfo {year} {1996})}\BibitemShut
  {NoStop}%
\bibitem [{\citenamefont {Lewandowski}(2022)}]{Lewandowski_JMP_2022}%
  \BibitemOpen
  \bibfield  {author} {\bibinfo {author} {\bibfnamefont {M.}~\bibnamefont
  {Lewandowski}},\ }\bibfield  {title} {\bibinfo {title} {Hadamard states for
  bosonic quantum field theory on globally hyperbolic spacetimes},\ }\href
  {https://doi.org/10.1063/5.0055753} {\bibfield  {journal} {\bibinfo
  {journal} {J. Math. Phys.}\ }\textbf {\bibinfo {volume} {63}},\ \bibinfo
  {pages} {013501} (\bibinfo {year} {2022})},\ \Eprint
  {https://arxiv.org/abs/2008.13156} {arXiv:2008.13156 [math-ph]}\BibitemShut
  {NoStop}%
\bibitem [{\citenamefont {B\"{a}r}\ and\ \citenamefont
  {Strohmaier}(2019)}]{Baer_AJM_2019}%
  \BibitemOpen
  \bibfield  {author} {\bibinfo {author} {\bibfnamefont {C.}~\bibnamefont
  {B\"{a}r}}\ and\ \bibinfo {author} {\bibfnamefont {A.}~\bibnamefont
  {Strohmaier}},\ }\bibfield  {title} {\bibinfo {title} {An index theorem for
  {L}orentzian manifolds with compact spacelike {C}auchy boundary},\ }\href
  {https://doi.org/10.1353/ajm.2019.0037} {\bibfield  {journal} {\bibinfo
  {journal} {Amer. J. Math.}\ }\textbf {\bibinfo {volume} {141}},\ \bibinfo
  {pages} {1421 } (\bibinfo {year} {2019})},\ \Eprint
  {https://arxiv.org/abs/1506.00959} {arXiv:1506.00959 [math.DG]}\BibitemShut
  {NoStop}%
\bibitem [{\citenamefont {B\"{a}r}\ and\ \citenamefont
  {Strohmaier}(2020)}]{Baer}%
  \BibitemOpen
  \bibfield  {author} {\bibinfo {author} {\bibfnamefont {C.}~\bibnamefont
  {B\"{a}r}}\ and\ \bibinfo {author} {\bibfnamefont {A.}~\bibnamefont
  {Strohmaier}},\ }\href@noop {} {\bibinfo {title} {Local index theory for
  {L}orentzian manifolds}} (\bibinfo {year} {2020}),\ \Eprint
  {https://arxiv.org/abs/2012.01364} {arXiv:2012.01364 [math.DG]}\BibitemShut
  {NoStop}%
\bibitem [{\citenamefont {Shen}\ and\ \citenamefont {Wrochna}(2021)}]{Shen}%
  \BibitemOpen
  \bibfield  {author} {\bibinfo {author} {\bibfnamefont {D.}~\bibnamefont
  {Shen}}\ and\ \bibinfo {author} {\bibfnamefont {M.}~\bibnamefont {Wrochna}},\
  }\href@noop {} {\bibinfo {title} {An index theorem on asymptotically static
  spacetimes with compact cauchy surface}} (\bibinfo {year} {2021}),\ \Eprint
  {https://arxiv.org/abs/2104.02816} {arXiv:2104.02816 [math.DG]}\BibitemShut
  {NoStop}%
\bibitem [{\citenamefont {Gell-Redman}\ \emph {et~al.}(2016)\citenamefont
  {Gell-Redman}, \citenamefont {Haber},\ and\ \citenamefont
  {Vasy}}]{GellRedman_CMP_2016}%
  \BibitemOpen
  \bibfield  {author} {\bibinfo {author} {\bibfnamefont {J.}~\bibnamefont
  {Gell-Redman}}, \bibinfo {author} {\bibfnamefont {N.}~\bibnamefont {Haber}},\
  and\ \bibinfo {author} {\bibfnamefont {A.}~\bibnamefont {Vasy}},\ }\bibfield
  {title} {\bibinfo {title} {The {F}eynman propagator on perturbations of
  {M}inkowski space},\ }\href {https://doi.org/10.1007/s00220-015-2520-8}
  {\bibfield  {journal} {\bibinfo  {journal} {Commun. Math. Phys.}\ }\textbf
  {\bibinfo {volume} {342}},\ \bibinfo {pages} {333} (\bibinfo {year}
  {2016})},\ \Eprint {https://arxiv.org/abs/1410.7113} {arXiv:1410.7113
  [math.AP]}\BibitemShut {NoStop}%
\bibitem [{\citenamefont {Dencker}(1982)}]{Dencker_JFA_1982}%
  \BibitemOpen
  \bibfield  {author} {\bibinfo {author} {\bibfnamefont {N.}~\bibnamefont
  {Dencker}},\ }\bibfield  {title} {\bibinfo {title} {On the propagation of
  polarization sets for systems of real principal type},\ }\href
  {https://doi.org/https://doi.org/10.1016/0022-1236(82)90051-9} {\bibfield
  {journal} {\bibinfo  {journal} {J. Funct. Anal.}\ }\textbf {\bibinfo {volume}
  {46}},\ \bibinfo {pages} {351 } (\bibinfo {year} {1982})}\BibitemShut
  {NoStop}%
\bibitem [{\citenamefont {G\'{e}rard}\ and\ \citenamefont
  {Wrochna}(2019)}]{Gerard_CMP_2019}%
  \BibitemOpen
  \bibfield  {author} {\bibinfo {author} {\bibfnamefont {C.}~\bibnamefont
  {G\'{e}rard}}\ and\ \bibinfo {author} {\bibfnamefont {M.}~\bibnamefont
  {Wrochna}},\ }\bibfield  {title} {\bibinfo {title} {Analytic {H}adamard
  states, {C}alder\'{o}n projectors and {W}ick rotation near analytic {C}auchy
  surfaces},\ }\href {https://doi.org/10.1007/s00220-019-03349-z} {\bibfield
  {journal} {\bibinfo  {journal} {Commun. Math. Phys.}\ }\textbf {\bibinfo
  {volume} {366}},\ \bibinfo {pages} {29} (\bibinfo {year} {2019})},\ \Eprint
  {https://arxiv.org/abs/1706.08942} {arXiv:1706.08942 [math-ph]}\BibitemShut
  {NoStop}%
\bibitem [{\citenamefont {Fulling}\ \emph {et~al.}(1978)\citenamefont
  {Fulling}, \citenamefont {Sweeny},\ and\ \citenamefont
  {Wald}}]{Fulling_CMP_1978}%
  \BibitemOpen
  \bibfield  {author} {\bibinfo {author} {\bibfnamefont {S.~A.}\ \bibnamefont
  {Fulling}}, \bibinfo {author} {\bibfnamefont {M.}~\bibnamefont {Sweeny}},\
  and\ \bibinfo {author} {\bibfnamefont {R.~M.}\ \bibnamefont {Wald}},\
  }\bibfield  {title} {\bibinfo {title} {Singularity structure of the two-point
  function quantum field theory in curved spacetime},\ }\href
  {https://projecteuclid.org:443/euclid.cmp/1103904566} {\bibfield  {journal}
  {\bibinfo  {journal} {Commun. Math. Phys.}\ }\textbf {\bibinfo {volume}
  {63}},\ \bibinfo {pages} {257} (\bibinfo {year} {1978})}\BibitemShut
  {NoStop}%
\bibitem [{\citenamefont {Brown}(1984)}]{Brown_JMP_1984}%
  \BibitemOpen
  \bibfield  {author} {\bibinfo {author} {\bibfnamefont {M.~R.}\ \bibnamefont
  {Brown}},\ }\bibfield  {title} {\bibinfo {title} {Symmetric {H}adamard
  series},\ }\href {https://doi.org/10.1063/1.526008} {\bibfield  {journal}
  {\bibinfo  {journal} {J. Math. Phys.}\ }\textbf {\bibinfo {volume} {25}},\
  \bibinfo {pages} {136 } (\bibinfo {year} {1984})}\BibitemShut {NoStop}%
\bibitem [{\citenamefont {Sahlmann}\ and\ \citenamefont
  {Verch}(2001)}]{Sahlmann_RMP_2001}%
  \BibitemOpen
  \bibfield  {author} {\bibinfo {author} {\bibfnamefont {H.}~\bibnamefont
  {Sahlmann}}\ and\ \bibinfo {author} {\bibfnamefont {R.}~\bibnamefont
  {Verch}},\ }\bibfield  {title} {\bibinfo {title} {Microlocal spectrum
  condition and {H}adamard form for vector-valued quantum fields in curved
  spacetime},\ }\href {https://doi.org/10.1142/S0129055X01001010} {\bibfield
  {journal} {\bibinfo  {journal} {Rev. Math. Phys.}\ }\textbf {\bibinfo
  {volume} {13}},\ \bibinfo {pages} {1203} (\bibinfo {year} {2001})},\ \Eprint
  {https://arxiv.org/abs/0008029} {arXiv:0008029 [math-ph]}\BibitemShut
  {NoStop}%
\bibitem [{\citenamefont {Marecki}(2003)}]{Marecki_master}%
  \BibitemOpen
  \bibfield  {author} {\bibinfo {author} {\bibfnamefont {P.}~\bibnamefont
  {Marecki}},\ }\emph {\bibinfo {title} {Quantum electrodynamics on background
  external fields}},\ \href@noop {} {Master's thesis} (\bibinfo {year}
  {2003}),\ \Eprint {https://arxiv.org/abs/0312304} {arXiv:0312304 [hep-th]}\BibitemShut 
  {NoStop}%
\bibitem [{\citenamefont {Hollands}(2008)}]{Hollands_RMP_2008}%
  \BibitemOpen
  \bibfield  {author} {\bibinfo {author} {\bibfnamefont {S.}~\bibnamefont
  {Hollands}},\ }\bibfield  {title} {\bibinfo {title} {Renormalized quantum
  {Y}ang-{M}ills fields in curved spacetime},\ }\href
  {https://doi.org/10.1142/S0129055X08003420} {\bibfield  {journal} {\bibinfo
  {journal} {Rev. Math. Phys.}\ }\textbf {\bibinfo {volume} {20}},\ \bibinfo
  {pages} {1033} (\bibinfo {year} {2008})},\ \Eprint
  {https://arxiv.org/abs/0705.3340} {arXiv:0705.3340 [gr-qc]}\BibitemShut
  {NoStop}%
\bibitem [{\citenamefont {Dappiaggi}\ \emph
  {et~al.}(2009{\natexlab{a}})\citenamefont {Dappiaggi}, \citenamefont {Hack},\
  and\ \citenamefont {Pinamonti}}]{Dappiaggi_RMP_2009}%
  \BibitemOpen
  \bibfield  {author} {\bibinfo {author} {\bibfnamefont {C.}~\bibnamefont
  {Dappiaggi}}, \bibinfo {author} {\bibfnamefont {T.-P.}\ \bibnamefont
  {Hack}},\ and\ \bibinfo {author} {\bibfnamefont {N.}~\bibnamefont
  {Pinamonti}},\ }\bibfield  {title} {\bibinfo {title} {The extended algebra of
  observables for {D}irac fields and the trace anomaly of their stress-energy
  tensor},\ }\href {https://doi.org/10.1142/S0129055X09003864} {\bibfield
  {journal} {\bibinfo  {journal} {Rev. Math. Phys.}\ }\textbf {\bibinfo
  {volume} {21}},\ \bibinfo {pages} {1241} (\bibinfo {year}
  {2009}{\natexlab{a}})},\ \Eprint {https://arxiv.org/abs/0904.0612}
  {arXiv:0904.0612 [math-ph]}\BibitemShut {NoStop}%
\bibitem [{\citenamefont {Fulling}\ \emph {et~al.}(1981)\citenamefont
  {Fulling}, \citenamefont {Narcowich},\ and\ \citenamefont
  {Wald}}]{Fulling_AnnPhys_1981}%
  \BibitemOpen
  \bibfield  {author} {\bibinfo {author} {\bibfnamefont {S.~A.}\ \bibnamefont
  {Fulling}}, \bibinfo {author} {\bibfnamefont {F.}~\bibnamefont {Narcowich}},\
  and\ \bibinfo {author} {\bibfnamefont {R.~M.}\ \bibnamefont {Wald}},\
  }\bibfield  {title} {\bibinfo {title} {Singularity structure of the two-point
  function in quantum field theory in curved spacetime, {I}{I}},\ }\href
  {https://doi.org/https://doi.org/10.1016/0003-4916(81)90098-1} {\bibfield
  {journal} {\bibinfo  {journal} {Ann. Phys.}\ }\textbf {\bibinfo {volume}
  {136}},\ \bibinfo {pages} {243} (\bibinfo {year} {1981})}\BibitemShut
  {NoStop}%
\bibitem [{\citenamefont {Murro}\ and\ \citenamefont
  {Volpe}(2021)}]{Murro_AGAG_2021}%
  \BibitemOpen
  \bibfield  {author} {\bibinfo {author} {\bibfnamefont {S.}~\bibnamefont
  {Murro}}\ and\ \bibinfo {author} {\bibfnamefont {D.}~\bibnamefont {Volpe}},\
  }\bibfield  {title} {\bibinfo {title} {Intertwining operators for symmetric
  hyperbolic systems on globally hyperbolic manifolds},\ }\href
  {https://doi.org/https://doi.org/10.1007/s10455-020-09739-0} {\bibfield
  {journal} {\bibinfo  {journal} {Ann. Glob. Anal. Geom.}\ }\textbf {\bibinfo
  {volume} {59}},\ \bibinfo {pages} {1} (\bibinfo {year} {2021})},\ \Eprint
  {https://arxiv.org/abs/2004.03300} {arXiv:2004.03300 [math.DG]}\BibitemShut
  {NoStop}%
\bibitem [{\citenamefont {Junker}(1996)}]{Junker_RMP_1996}%
  \BibitemOpen
  \bibfield  {author} {\bibinfo {author} {\bibfnamefont {W.}~\bibnamefont
  {Junker}},\ }\bibfield  {title} {\bibinfo {title} {Hadamard states, adiabatic
  vacua and the construction of physical states for scalar quantum fields on
  curved space-time},\ }\href {https://doi.org/10.1142/S0129055X9600041X}
  {\bibfield  {journal} {\bibinfo  {journal} {Rev. Math. Phys.}\ }\textbf
  {\bibinfo {volume} {8}},\ \bibinfo {pages} {1091} (\bibinfo {year} {1996})},\
  \bibinfo {note} {erratum: Rev. Math. Phys. \textbf{14}, 511
  (2002)}\BibitemShut {NoStop}%
\bibitem [{\citenamefont {Hollands}(2001)}]{Hollands_adiabatic_CMP_2001}%
  \BibitemOpen
  \bibfield  {author} {\bibinfo {author} {\bibfnamefont {S.}~\bibnamefont
  {Hollands}},\ }\bibfield  {title} {\bibinfo {title} {The {H}adamard condition
  for {D}irac fields and adiabatic states on {R}obertson-{W}alker spacetimes},\
  }\href {https://doi.org/10.1007/s002200000350} {\bibfield  {journal}
  {\bibinfo  {journal} {Commun. Math. Phys.}\ }\textbf {\bibinfo {volume}
  {216}},\ \bibinfo {pages} {635} (\bibinfo {year} {2001})},\ \Eprint
  {https://arxiv.org/abs/0102035} {arXiv:0102035 [math-ph]}\BibitemShut
  {NoStop}%
\bibitem [{\citenamefont {G\'{e}rard}\ and\ \citenamefont
  {Wrochna}(2014)}]{Gerard_CMP_2014}%
  \BibitemOpen
  \bibfield  {author} {\bibinfo {author} {\bibfnamefont {C.}~\bibnamefont
  {G\'{e}rard}}\ and\ \bibinfo {author} {\bibfnamefont {M.}~\bibnamefont
  {Wrochna}},\ }\bibfield  {title} {\bibinfo {title} {Construction of
  {H}adamard states by pseudo-differential calculus},\ }\href
  {https://doi.org/10.1007/s00220-013-1824-9} {\bibfield  {journal} {\bibinfo
  {journal} {Commun. Math. Phys.}\ }\textbf {\bibinfo {volume} {325}},\
  \bibinfo {pages} {713} (\bibinfo {year} {2014})},\ \Eprint
  {https://arxiv.org/abs/1209.2604} {arXiv:1209.2604 [math-ph]}\BibitemShut
  {NoStop}%
\bibitem [{\citenamefont {G\'{e}rard}\ and\ \citenamefont
  {Wrochna}(2015)}]{Gerard_CMP_2015}%
  \BibitemOpen
  \bibfield  {author} {\bibinfo {author} {\bibfnamefont {C.}~\bibnamefont
  {G\'{e}rard}}\ and\ \bibinfo {author} {\bibfnamefont {M.}~\bibnamefont
  {Wrochna}},\ }\bibfield  {title} {\bibinfo {title} {Hadamard states for the
  linearized {Y}ang-{M}ills equation on curved spacetime},\ }\href
  {https://doi.org/https://doi.org/10.1007/s00220-015-2305-0} {\bibfield
  {journal} {\bibinfo  {journal} {Commun. Math. Phys.}\ }\textbf {\bibinfo
  {volume} {337}},\ \bibinfo {pages} {253} (\bibinfo {year} {2015})},\ \Eprint
  {https://arxiv.org/abs/1403.7153} {arXiv:1403.7153 [math-ph]}\BibitemShut
  {NoStop}%
\bibitem [{\citenamefont {G\'{e}rard}\ and\ \citenamefont
  {Stoskopf}(2022{\natexlab{a}})}]{Gerard}%
  \BibitemOpen
  \bibfield  {author} {\bibinfo {author} {\bibfnamefont {C.}~\bibnamefont
  {G\'{e}rard}}\ and\ \bibinfo {author} {\bibfnamefont {T.}~\bibnamefont
  {Stoskopf}},\ }\bibfield  {title} {\bibinfo {title} {Hadamard property of the
  in and out states for {D}irac fields on asymptotically static spacetimes},\
  }\href {https://doi.org/https://doi.org/10.1007/s11005-022-01556-9}
  {\bibfield  {journal} {\bibinfo  {journal} {Lett. Math. Phys.}\ }\textbf
  {\bibinfo {volume} {112}},\ \bibinfo {pages} {63} (\bibinfo {year}
  {2022}{\natexlab{a}})},\ \Eprint {https://arxiv.org/abs/arXiv:2108.11955
  [math.AP]} {arXiv:2108.11955 [math.AP]}\BibitemShut {NoStop}%
\bibitem [{\citenamefont {G\'{e}rard}\ and\ \citenamefont
  {Stoskopf}(2022{\natexlab{b}})}]{Gerard_bounded_geometry}%
  \BibitemOpen
  \bibfield  {author} {\bibinfo {author} {\bibfnamefont {C.}~\bibnamefont
  {G\'{e}rard}}\ and\ \bibinfo {author} {\bibfnamefont {T.}~\bibnamefont
  {Stoskopf}},\ }\bibfield  {title} {\bibinfo {title} {Hadamard states for
  quantized {D}irac fields on {L}orentzian manifolds of bounded geometry},\
  }\href {https://doi.org/10.1142/S0129055X22500088} {\bibfield  {journal}
  {\bibinfo  {journal} {Rev. in Math. Phys.}\ }\textbf {\bibinfo {volume}
  {34}},\ \bibinfo {pages} {2250008} (\bibinfo {year} {2022}{\natexlab{b}})},\
  \Eprint {https://arxiv.org/abs/arXiv:2108.11630 [math.AP]} {arXiv:2108.11630
  [math.AP]} \BibitemShut {NoStop}%
\bibitem [{\citenamefont {G\'{e}rard}(2019)}]{Gerard_EMS_2019}%
  \BibitemOpen
  \bibfield  {author} {\bibinfo {author} {\bibfnamefont {C.}~\bibnamefont
  {G\'{e}rard}},\ }\href {https://doi.org/10.4171/094} {\emph {\bibinfo {title}
  {Microlocal Analysis of Quantum Fields on Curved Spacetimes}}},\ ESI Lectures
  in Mathematics and Physics\ (\bibinfo  {publisher} {European Mathematical
  Society},\ \bibinfo {address} {Germany},\ \bibinfo {year} {2019})\BibitemShut
  {NoStop}%
\bibitem [{\citenamefont {Moretti}(2008)}]{Moretti_CMP_2008}%
  \BibitemOpen
  \bibfield  {author} {\bibinfo {author} {\bibfnamefont {V.}~\bibnamefont
  {Moretti}},\ }\bibfield  {title} {\bibinfo {title} {Quantum ground states
  holographically induced by asymptotic flatness: Invariance under spacetime
  symmetries, energy positivity and {H}adamard property},\ }\href
  {https://doi.org/10.1007/s00220-008-0415-7} {\bibfield  {journal} {\bibinfo
  {journal} {Commun. Math. Phys.}\ }\textbf {\bibinfo {volume} {279}},\
  \bibinfo {pages} {31} (\bibinfo {year} {2008})},\ \Eprint
  {https://arxiv.org/abs/0610143} {arXiv:0610143 [gr-qc]}\BibitemShut
  {NoStop}%
\bibitem [{\citenamefont {Dappiaggi}\ \emph
  {et~al.}(2009{\natexlab{b}})\citenamefont {Dappiaggi}, \citenamefont
  {Moretti},\ and\ \citenamefont {Pinamonti}}]{Dappiaggi_JMP_2009}%
  \BibitemOpen
  \bibfield  {author} {\bibinfo {author} {\bibfnamefont {C.}~\bibnamefont
  {Dappiaggi}}, \bibinfo {author} {\bibfnamefont {V.}~\bibnamefont {Moretti}},\
  and\ \bibinfo {author} {\bibfnamefont {N.}~\bibnamefont {Pinamonti}},\
  }\bibfield  {title} {\bibinfo {title} {Distinguished quantum states in a
  class of cosmological spacetimes and their {H}adamard property},\ }\href
  {https://doi.org/10.1063/1.3122770} {\bibfield  {journal} {\bibinfo
  {journal} {J. Math. Phys.}\ }\textbf {\bibinfo {volume} {50}},\ \bibinfo
  {pages} {062304} (\bibinfo {year} {2009}{\natexlab{b}})},\ \Eprint
  {https://arxiv.org/abs/0812.4033} {arXiv:0812.4033 [gr-qc]}\BibitemShut
  {NoStop}%
\bibitem [{\citenamefont {Dappiaggi}\ \emph {et~al.}(2011)\citenamefont
  {Dappiaggi}, \citenamefont {Moretti},\ and\ \citenamefont
  {Pinamonti}}]{Dappiaggi_ATMP_2011}%
  \BibitemOpen
  \bibfield  {author} {\bibinfo {author} {\bibfnamefont {C.}~\bibnamefont
  {Dappiaggi}}, \bibinfo {author} {\bibfnamefont {V.}~\bibnamefont {Moretti}},\
  and\ \bibinfo {author} {\bibfnamefont {N.}~\bibnamefont {Pinamonti}},\
  }\bibfield  {title} {\bibinfo {title} {Rigorous construction and {H}adamard
  property of the {U}nruh state in {S}chwarzschild spacetime},\ }\href
  {https://doi.org/10.4310/ATMP.2011.v15.n2.a4} {\bibfield  {journal} {\bibinfo
   {journal} {Adv. Theor. Math. Phys.}\ }\textbf {\bibinfo {volume} {15}},\
  \bibinfo {pages} {355} (\bibinfo {year} {2011})},\ \Eprint
  {https://arxiv.org/abs/0907.1034} {arXiv:0907.1034 [gr-qc]}\BibitemShut
  {NoStop}%
\bibitem [{\citenamefont {G\'{e}rard}\ and\ \citenamefont
  {Wrochna}(2016)}]{Gerard_AnalPDE_2016}%
  \BibitemOpen
  \bibfield  {author} {\bibinfo {author} {\bibfnamefont {C.}~\bibnamefont
  {G\'{e}rard}}\ and\ \bibinfo {author} {\bibfnamefont {M.}~\bibnamefont
  {Wrochna}},\ }\bibfield  {title} {\bibinfo {title} {Construction of
  {H}adamard states by characteristic {C}auchy problem},\ }\href
  {https://doi.org/10.2140/apde.2016.9.111} {\bibfield  {journal} {\bibinfo
  {journal} {Anal. PDE}\ }\textbf {\bibinfo {volume} {9}},\ \bibinfo {pages}
  {111} (\bibinfo {year} {2016})},\ \Eprint {https://arxiv.org/abs/1409.6691}
  {arXiv:1409.6691 [math-ph]}\BibitemShut {NoStop}%
\bibitem [{\citenamefont {Dappiaggi}\ \emph {et~al.}(2017)\citenamefont
  {Dappiaggi}, \citenamefont {Moretti},\ and\ \citenamefont
  {Pinamonti}}]{Dappiaggi_Springer_2017}%
  \BibitemOpen
  \bibfield  {author} {\bibinfo {author} {\bibfnamefont {C.}~\bibnamefont
  {Dappiaggi}}, \bibinfo {author} {\bibfnamefont {V.}~\bibnamefont {Moretti}},\
  and\ \bibinfo {author} {\bibfnamefont {N.}~\bibnamefont {Pinamonti}},\ }\href
  {https://doi.org/10.1007/978-3-319-64343-4} {\emph {\bibinfo {title}
  {Hadamard States From Light-like Hypersurfaces}}},\ \bibinfo {series}
  {SpringerBriefs in Mathematical Physics}, Vol.~\bibinfo {volume} {25}\
  (\bibinfo {year} {2017})\ \Eprint {https://arxiv.org/abs/1706.09666}
  {arXiv:1706.09666 [math-ph]}\BibitemShut {NoStop}%
\bibitem [{\citenamefont {Vasy}(2017)}]{Vasy_AHP_2017}%
  \BibitemOpen
  \bibfield  {author} {\bibinfo {author} {\bibfnamefont {A.}~\bibnamefont
  {Vasy}},\ }\bibfield  {title} {\bibinfo {title} {On the positivity of
  propagator differences},\ }\href
  {https://doi.org/https://doi.org/10.1007/s00023-016-0527-0} {\bibfield
  {journal} {\bibinfo  {journal} {Ann. Henri Poincar\'{e}}\ }\textbf {\bibinfo
  {volume} {18}},\ \bibinfo {pages} {983 } (\bibinfo {year} {2017})},\ \Eprint
  {https://arxiv.org/abs/1411.7242} {arXiv:1411.7242 [math.AP]}\BibitemShut
  {NoStop}%
\bibitem [{\citenamefont {Derezi\'{n}ski}\ and\ \citenamefont
  {Siemssen}(2018)}]{Derezinski_RMP_2018}%
  \BibitemOpen
  \bibfield  {author} {\bibinfo {author} {\bibfnamefont {J.}~\bibnamefont
  {Derezi\'{n}ski}}\ and\ \bibinfo {author} {\bibfnamefont {D.}~\bibnamefont
  {Siemssen}},\ }\bibfield  {title} {\bibinfo {title} {Feynman propagators on
  static spacetimes},\ }\href {https://doi.org/10.1142/S0129055X1850006X}
  {\bibfield  {journal} {\bibinfo  {journal} {Rev. Math. Phys.}\ }\textbf
  {\bibinfo {volume} {30}},\ \bibinfo {pages} {1850006} (\bibinfo {year}
  {2018})},\ \Eprint {https://arxiv.org/abs/1608.06441} {arXiv:1608.06441
  [math-ph]}\BibitemShut {NoStop}%
\bibitem [{\citenamefont {Capoferri}\ \emph {et~al.}(2020)\citenamefont
  {Capoferri}, \citenamefont {Dappiaggi},\ and\ \citenamefont
  {Drago}}]{Capoferri_JMAA_2020}%
  \BibitemOpen
  \bibfield  {author} {\bibinfo {author} {\bibfnamefont {M.}~\bibnamefont
  {Capoferri}}, \bibinfo {author} {\bibfnamefont {C.}~\bibnamefont
  {Dappiaggi}},\ and\ \bibinfo {author} {\bibfnamefont {N.}~\bibnamefont
  {Drago}},\ }\bibfield  {title} {\bibinfo {title} {Global wave parametrices on
  globally hyperbolic spacetimes},\ }\href
  {https://doi.org/https://doi.org/10.1016/j.jmaa.2020.124316} {\bibfield
  {journal} {\bibinfo  {journal} {J. Math. Anal. Appl.}\ }\textbf {\bibinfo
  {volume} {490}},\ \bibinfo {pages} {124316} (\bibinfo {year} {2020})},\
  \Eprint {https://arxiv.org/abs/2001.04164} {arXiv:2001.04164 [math.AP]}\BibitemShut 
  {NoStop}%
\bibitem [{\citenamefont {Avetisyan}\ and\ \citenamefont
  {Capoferri}(2021)}]{Avetisyan_Math_2021}%
  \BibitemOpen
  \bibfield  {author} {\bibinfo {author} {\bibfnamefont {Z.}~\bibnamefont
  {Avetisyan}}\ and\ \bibinfo {author} {\bibfnamefont {M.}~\bibnamefont
  {Capoferri}},\ }\bibfield  {title} {\bibinfo {title} {Partial differential
  equations and quantum states in curved spacetimes},\ }\href
  {https://doi.org/10.3390/math9161936} {\bibfield  {journal} {\bibinfo
  {journal} {Mathematics}\ }\textbf {\bibinfo {volume} {9}},\ \bibinfo {pages}
  {16} (\bibinfo {year} {2021})}\BibitemShut {NoStop}%
\bibitem [{\citenamefont {Guillemin}\ and\ \citenamefont
  {Sternberg}(1977)}]{Guillemin_AMS_1977}%
  \BibitemOpen
  \bibfield  {author} {\bibinfo {author} {\bibfnamefont {V.}~\bibnamefont
  {Guillemin}}\ and\ \bibinfo {author} {\bibfnamefont {S.}~\bibnamefont
  {Sternberg}},\ }\href {https://bookstore.ams.org/surv-14} {\emph {\bibinfo
  {title} {Geometric Asymptotics}}},\ \bibinfo {series} {Mathematical Surveys
  and Monographs}, Vol.~\bibinfo {volume} {14}\ (\bibinfo  {publisher}
  {American Mathemtical Society},\ \bibinfo {address} {USA},\ \bibinfo {year}
  {1977})\BibitemShut {NoStop}%
\bibitem [{\citenamefont {G\"{u}nther}(1988)}]{Guenther_AP_1988}%
  \BibitemOpen
  \bibfield  {author} {\bibinfo {author} {\bibfnamefont {P.}~\bibnamefont
  {G\"{u}nther}},\ }\href
  {https://doi.org/https://doi.org/10.1016/C2013-0-10776-3} {\emph {\bibinfo
  {title} {Huygens' Principle and Hyperbolic Equations}}},\ \bibinfo {series}
  {Prespective in Mathematics}, Vol.~\bibinfo {volume} {5}\ (\bibinfo
  {publisher} {Academic Press},\ \bibinfo {address} {USA},\ \bibinfo {year}
  {1988})\BibitemShut {NoStop}%
\bibitem [{\citenamefont {H\"{o}rmander}(1965)}]{Hoermander_CPAM_1965}%
  \BibitemOpen
  \bibfield  {author} {\bibinfo {author} {\bibfnamefont {L.}~\bibnamefont
  {H\"{o}rmander}},\ }\bibfield  {title} {\bibinfo {title} {Pseudo-differential
  operators},\ }\href {https://doi.org/10.1002/cpa.3160180307} {\bibfield
  {journal} {\bibinfo  {journal} {Comm. Pure Appl. Math.}\ }\textbf {\bibinfo
  {volume} {18}},\ \bibinfo {pages} {501 } (\bibinfo {year}
  {1965})}\BibitemShut {NoStop}%
\bibitem [{\citenamefont {Weinstein}(1976)}]{Weinstein_BullAMS_1976}%
  \BibitemOpen
  \bibfield  {author} {\bibinfo {author} {\bibfnamefont {A.}~\bibnamefont
  {Weinstein}},\ }\bibfield  {title} {\bibinfo {title} {The principal symbol of
  a distribution},\ }\href
  {https://projecteuclid.org:443/euclid.bams/1183538121} {\bibfield  {journal}
  {\bibinfo  {journal} {Bull. Amer. Math. Soc.}\ }\textbf {\bibinfo {volume}
  {82}},\ \bibinfo {pages} {548} (\bibinfo {year} {1976})}\BibitemShut
  {NoStop}%
\bibitem [{\citenamefont {Weinstein}(1978)}]{Weinstein_TransAMS_1978}%
  \BibitemOpen
  \bibfield  {author} {\bibinfo {author} {\bibfnamefont {A.}~\bibnamefont
  {Weinstein}},\ }\bibfield  {title} {\bibinfo {title} {The order and symbol of
  a distribution},\ }\href
  {https://doi.org/https://doi.org/10.1090/S0002-9947-1978-0492288-9}
  {\bibfield  {journal} {\bibinfo  {journal} {Trans. Amer. Math. Soc.}\
  }\textbf {\bibinfo {volume} {241}},\ \bibinfo {pages} {1} (\bibinfo {year}
  {1978})}\BibitemShut {NoStop}%
\bibitem [{\citenamefont {H\"{o}rmander}(2003)}]{Hoermander_Springer_2003}%
  \BibitemOpen
  \bibfield  {author} {\bibinfo {author} {\bibfnamefont {L.}~\bibnamefont
  {H\"{o}rmander}},\ }\href {https://doi.org/10.1007/978-3-642-61497-2} {\emph
  {\bibinfo {title} {The Analysis of Linear Partial Differential Operators {I}:
  Distribution Theory and Fourier Analysis}}},\ \bibinfo {edition} {2nd}\ ed.,\
  Classics in Mathematics\ (\bibinfo  {publisher} {Springer-Verlag},\ \bibinfo
  {address} {Berlin, Heidelberg},\ \bibinfo {year} {2003})\BibitemShut
  {NoStop}%
\bibitem [{\citenamefont {H\"{o}rmander}(2007)}]{Hoermander_Springer_2007}%
  \BibitemOpen
  \bibfield  {author} {\bibinfo {author} {\bibfnamefont {L.}~\bibnamefont
  {H\"{o}rmander}},\ }\href {https://doi.org/10.1007/978-3-540-49938-1} {\emph
  {\bibinfo {title} {The Analysis of Linear Partial Differential Operators
  {I}{I}{I}: Pseudo-Differential Operators}}},\ Classics in Mathematics\
  (\bibinfo  {publisher} {Springer-Verlag},\ \bibinfo {address} {Berlin,
  Heidelberg},\ \bibinfo {year} {2007})\BibitemShut {NoStop}%
\bibitem [{\citenamefont {H\"{o}rmander}(1971)}]{Hoermander_ActaMath_1971}%
  \BibitemOpen
  \bibfield  {author} {\bibinfo {author} {\bibfnamefont {L.}~\bibnamefont
  {H\"{o}rmander}},\ }\bibfield  {title} {\bibinfo {title} {Fourier integral
  operators. {I}},\ }\href {https://doi.org/10.1007/BF02392052} {\bibfield
  {journal} {\bibinfo  {journal} {Acta Math.}\ }\textbf {\bibinfo {volume}
  {127}},\ \bibinfo {pages} {79 } (\bibinfo {year} {1971})}\BibitemShut
  {NoStop}%
\bibitem [{\citenamefont {H\"{o}rmander}(1970)}]{Hoermander_Nice_1970}%
  \BibitemOpen
  \bibfield  {author} {\bibinfo {author} {\bibfnamefont {L.}~\bibnamefont
  {H\"{o}rmander}},\ }\bibfield  {title} {\bibinfo {title} {Linear differential
  operators},\ }in\ \href {https://www.mathunion.org/icm/proceedings} {\emph
  {\bibinfo {booktitle} {Actes, Congr\'{e}s intern. math,}}}\ (\bibinfo
  {address} {Nice, France},\ \bibinfo {year} {1970})\ pp.\ \bibinfo {pages}
  {121 -- 133}\BibitemShut {NoStop}%
\bibitem [{\citenamefont {Strohmaier}(2009)}]{Strohmaier_Springer_2009}%
  \BibitemOpen
  \bibfield  {author} {\bibinfo {author} {\bibfnamefont {A.}~\bibnamefont
  {Strohmaier}},\ }\bibinfo {title} {Microlocal analysis},\ in\ \href
  {https://doi.org/10.1007/978-3-642-02780-2} {\emph {\bibinfo {booktitle}
  {Quantum Field Theory on Curved Spacetimes}}},\ \bibinfo {series} {Lect.
  Notes Phys.}, Vol.\ \bibinfo {volume} {786},\ \bibinfo {editor} {edited by\
  \bibinfo {editor} {\bibfnamefont {C.}~\bibnamefont {B\"{a}r}}\ and\ \bibinfo
  {editor} {\bibfnamefont {K.}~\bibnamefont {Fredenhagen}}}\ (\bibinfo
  {publisher} {Springer-Verlag},\ \bibinfo {address} {Berlin Heidelberg},\
  \bibinfo {year} {2009})\ pp.\ \bibinfo {pages} {85 -- 127}\BibitemShut
  {NoStop}%
\bibitem [{\citenamefont {Brouder}\ \emph {et~al.}(2014)\citenamefont
  {Brouder}, \citenamefont {Dang},\ and\ \citenamefont
  {H{\'{e}}lein}}]{Brouder_JPA_2014}%
  \BibitemOpen
  \bibfield  {author} {\bibinfo {author} {\bibfnamefont {C.}~\bibnamefont
  {Brouder}}, \bibinfo {author} {\bibfnamefont {N.~V.}\ \bibnamefont {Dang}},\
  and\ \bibinfo {author} {\bibfnamefont {F.}~\bibnamefont {H{\'{e}}lein}},\
  }\bibfield  {title} {\bibinfo {title} {A smooth introduction to the wavefront
  set},\ }\href {https://doi.org/10.1088/1751-8113/47/44/443001} {\bibfield
  {journal} {\bibinfo  {journal} {J. Phys. A: Math. Theor.}\ }\textbf {\bibinfo
  {volume} {47}},\ \bibinfo {pages} {443001} (\bibinfo {year} {2014})},\
  \Eprint {https://arxiv.org/abs/1404.1778} {arXiv:1404.1778 [math-ph]}\BibitemShut {NoStop}%
\bibitem [{\citenamefont {Peetre}(1960)}]{Peetre_MathScand_1960}%
  \BibitemOpen
  \bibfield  {author} {\bibinfo {author} {\bibfnamefont {J.}~\bibnamefont
  {Peetre}},\ }\bibfield  {title} {\bibinfo {title} {R\'{e}ctification \`{a}
  l'article "une caract\'{e}risation abstraite des op\'{e}rateurs
  diff\'{e}rentiels"},\ }\href {https://doi.org/10.7146/math.scand.a-10598}
  {\bibfield  {journal} {\bibinfo  {journal} {Math Scand}\ }\textbf {\bibinfo
  {volume} {8}},\ \bibinfo {pages} {116} (\bibinfo {year} {1960})}\BibitemShut
  {NoStop}%
\bibitem [{\citenamefont {Duistermaat}(1974)}]{Duistermaat_CPAM_1974}%
  \BibitemOpen
  \bibfield  {author} {\bibinfo {author} {\bibfnamefont {J.~J.}\ \bibnamefont
  {Duistermaat}},\ }\bibfield  {title} {\bibinfo {title} {Oscillatory
  integrals, {L}agrange immersions and unfolding of singularities},\ }\href
  {https://doi.org/10.1002/cpa.3160270205} {\bibfield  {journal} {\bibinfo
  {journal} {Comm. Pure Appl. Math.}\ }\textbf {\bibinfo {volume} {27}},\
  \bibinfo {pages} {207} (\bibinfo {year} {1974})}\BibitemShut {NoStop}%
\bibitem [{\citenamefont {Duistermaat}(2011)}]{Duistermaat_Birkhaeuser_2011}%
  \BibitemOpen
  \bibfield  {author} {\bibinfo {author} {\bibfnamefont {J.~J.}\ \bibnamefont
  {Duistermaat}},\ }\href {https://doi.org/10.1007/978-0-8176-8108-1} {\emph
  {\bibinfo {title} {Fourier Integral Operators}}},\ Modern Birkh\"{a}user
  Classics\ (\bibinfo  {publisher} {Birkh\"{a}user},\ \bibinfo {address} {New
  York},\ \bibinfo {year} {2011})\BibitemShut {NoStop}%
\bibitem [{\citenamefont {Treves}(1982)}]{Treves_Plenum_1980}%
  \BibitemOpen
  \bibfield  {author} {\bibinfo {author} {\bibfnamefont {J.-F.}\ \bibnamefont
  {Treves}},\ }\href@noop {} {\emph {\bibinfo {title} {Introduction to
  Pseudodifferential and Fourier Integral Operators: Fourier Integral
  Operators}}},\ \bibinfo {series} {University Series in Mathematics},
  Vol.~\bibinfo {volume} {2}\ (\bibinfo  {publisher} {Plenum Press},\ \bibinfo
  {address} {New York},\ \bibinfo {year} {1980; Second Printing
  1982})\BibitemShut {NoStop}%
\bibitem [{\citenamefont {Bates}\ and\ \citenamefont
  {Weinstein}(1997)}]{Bates_AMS_1997}%
  \BibitemOpen
  \bibfield  {author} {\bibinfo {author} {\bibfnamefont {S.}~\bibnamefont
  {Bates}}\ and\ \bibinfo {author} {\bibfnamefont {A.}~\bibnamefont
  {Weinstein}},\ }\href {https://bookstore.ams.org/bmln-8} {\emph {\bibinfo
  {title} {Lectures on the Geometry of Quantization}}},\ \bibinfo {series}
  {Berkeley Mathematics Lecture Notes}, Vol.~\bibinfo {volume} {8}\ (\bibinfo
  {publisher} {American Mathematical Society and Berkeley Center for Pure and
  Applied Mathematics},\ \bibinfo {address} {USA},\ \bibinfo {year}
  {1997})\BibitemShut {NoStop}%
\bibitem [{\citenamefont {Melrose}(2007)}]{Melrose_2007}%
  \BibitemOpen
  \bibfield  {author} {\bibinfo {author} {\bibfnamefont {R.}~\bibnamefont
  {Melrose}},\ }\href {http://math.mit.edu/~rbm/18.157-F07.html} {\bibinfo
  {title} {Introduction to microlocal analysis}} (\bibinfo {year} {2007}),\
  \bibinfo {note} {unpublished lecture notes taught at {M}{I}{T}}\BibitemShut
  {NoStop}%
\bibitem [{\citenamefont {Lee}(2013)}]{Lee_Springer_2013}%
  \BibitemOpen
  \bibfield  {author} {\bibinfo {author} {\bibfnamefont {J.~M.}\ \bibnamefont
  {Lee}},\ }\href {https://doi.org/https://doi.org/10.1007/978-1-4419-9982-5}
  {\emph {\bibinfo {title} {Introduction to Smooth Manifolds}}},\ \bibinfo
  {edition} {2nd}\ ed.,\ \bibinfo {series} {Graduate Texts in Mathematics},
  Vol.\ \bibinfo {volume} {218}\ (\bibinfo  {publisher} {Springer International
  Publishing},\ \bibinfo {year} {2013})\BibitemShut {NoStop}%
\bibitem [{\citenamefont {O'Niell}(1983)}]{ONeill_Academic_1983}%
  \BibitemOpen
  \bibfield  {author} {\bibinfo {author} {\bibfnamefont {B.}~\bibnamefont
  {O'Niell}},\ }\href
  {https://www.elsevier.com/books/semi-riemannian-geometry-with-applications-to-relativity/oneill/978-0-12-526740-3}
  {\emph {\bibinfo {title} {Semi-{R}iemannian Geometry}}}\ (\bibinfo
  {publisher} {Academic Press},\ \bibinfo {address} {USA},\ \bibinfo {year}
  {1983})\BibitemShut {NoStop}%
\bibitem [{\citenamefont {van~den Ban}\ and\ \citenamefont
  {Crainic}(2017)}]{vandenBan_2017}%
  \BibitemOpen
  \bibfield  {author} {\bibinfo {author} {\bibfnamefont {E.~P.}\ \bibnamefont
  {van~den Ban}}\ and\ \bibinfo {author} {\bibfnamefont {M.}~\bibnamefont
  {Crainic}},\ }\href
  {https://www.staff.science.uu.nl/~ban00101/geoman2017/AS-2017rev.pdf}
  {\bibinfo {title} {Analysis on manifolds}} (\bibinfo {year} {2017}),\
  \bibinfo {note} {unpublished lecture notes taught at Utrecht
  University}\BibitemShut {NoStop}%
\bibitem [{\citenamefont {Toth}\ and\ \citenamefont
  {Zelditch}(2013)}]{Toth_GFA_2013}%
  \BibitemOpen
  \bibfield  {author} {\bibinfo {author} {\bibfnamefont {J.~A.}\ \bibnamefont
  {Toth}}\ and\ \bibinfo {author} {\bibfnamefont {S.}~\bibnamefont
  {Zelditch}},\ }\bibfield  {title} {\bibinfo {title} {Quantum ergodic
  restriction theorems: Manifolds without boundary},\ }\href
  {https://doi.org/https://doi.org/10.1007/s00039-013-0220-0} {\bibfield
  {journal} {\bibinfo  {journal} {Geom. Funct. Anal.}\ }\textbf {\bibinfo
  {volume} {23}},\ \bibinfo {pages} {715 } (\bibinfo {year} {2013})},\ \Eprint
  {https://arxiv.org/abs/1104.4531} {arXiv:1104.4531 [math.SP]}\BibitemShut
  {NoStop}%
\bibitem [{\citenamefont {Low}(1989)}]{Low_JMP_1989}%
  \BibitemOpen
  \bibfield  {author} {\bibinfo {author} {\bibfnamefont {R.~J.}\ \bibnamefont
  {Low}},\ }\bibfield  {title} {\bibinfo {title} {The geometry of the space of
  null geodesics},\ }\href {https://doi.org/10.1063/1.528401} {\bibfield
  {journal} {\bibinfo  {journal} {J. Math. Phys.}\ }\textbf {\bibinfo {volume}
  {30}},\ \bibinfo {pages} {809} (\bibinfo {year} {1989})}\BibitemShut
  {NoStop}%
\bibitem [{\citenamefont {Low}(2001)}]{Low_NonlinearAnal_2001}%
  \BibitemOpen
  \bibfield  {author} {\bibinfo {author} {\bibfnamefont {R.}~\bibnamefont
  {Low}},\ }\bibfield  {title} {\bibinfo {title} {The space of null
  geodesics},\ }\href
  {https://doi.org/https://doi.org/10.1016/S0362-546X(01)00421-7} {\bibfield
  {journal} {\bibinfo  {journal} {Nonlinear Anal. Theory Methods. Appl.}\
  }\textbf {\bibinfo {volume} {47}},\ \bibinfo {pages} {3005 } (\bibinfo {year}
  {2001})},\ \bibinfo {note} {{P}roc. Third World Congress of Nonlinear
  Analysts}\BibitemShut {NoStop}%
\bibitem [{\citenamefont {Weinstein}(1975)}]{Weinstein_Nice_1975}%
  \BibitemOpen
  \bibfield  {author} {\bibinfo {author} {\bibfnamefont {A.}~\bibnamefont
  {Weinstein}},\ }\bibfield  {title} {\bibinfo {title} {On {M}aslov's
  quantization condition},\ }in\ \href {https://doi.org/10.1007/BFb0074200}
  {\emph {\bibinfo {booktitle} {Fourier Integral Operators and Partial
  Differential Equations}}},\ \bibinfo {editor} {edited by\ \bibinfo {editor}
  {\bibfnamefont {J.}~\bibnamefont {Chazarain}}}\ (\bibinfo  {publisher}
  {Springer Berlin Heidelberg},\ \bibinfo {address} {Berlin, Heidelberg},\
  \bibinfo {year} {1975})\ pp.\ \bibinfo {pages} {341--372}\BibitemShut
  {NoStop}%
\bibitem [{\citenamefont {Guillemin}(1993)}]{Guillemin_JFA_1993}%
  \BibitemOpen
  \bibfield  {author} {\bibinfo {author} {\bibfnamefont {V.}~\bibnamefont
  {Guillemin}},\ }\bibfield  {title} {\bibinfo {title} {Residue traces for
  certain algebras of {F}ourier integral operators},\ }\href
  {https://doi.org/https://doi.org/10.1006/jfan.1993.1096} {\bibfield
  {journal} {\bibinfo  {journal} {J. Funct. Anal.}\ }\textbf {\bibinfo {volume}
  {115}},\ \bibinfo {pages} {391 } (\bibinfo {year} {1993})}\BibitemShut
  {NoStop}%
\bibitem [{\citenamefont {Abraham}\ and\ \citenamefont
  {Marsden}(2008)}]{Abraham_AMS_1978}%
  \BibitemOpen
  \bibfield  {author} {\bibinfo {author} {\bibfnamefont {R.}~\bibnamefont
  {Abraham}}\ and\ \bibinfo {author} {\bibfnamefont {J.~E.}\ \bibnamefont
  {Marsden}},\ }\href {https://doi.org/https://doi.org/10.1090/chel/364} {\emph
  {\bibinfo {title} {Foundations of Mechanics}}},\ \bibinfo {edition} {2nd}\
  ed.,\ Vol.\ \bibinfo {volume} {364}\ (\bibinfo  {publisher} {AMS Chelsea
  Publishing},\ \bibinfo {address} {USA},\ \bibinfo {year} {1978; AMS reprint
  2008})\BibitemShut {NoStop}%
\bibitem [{\citenamefont {Eskin}(2011)}]{Eskin_AMS_2011}%
  \BibitemOpen
  \bibfield  {author} {\bibinfo {author} {\bibfnamefont {G.}~\bibnamefont
  {Eskin}},\ }\href {https://bookstore.ams.org/gsm-123/} {\emph {\bibinfo
  {title} {Lectures on Linear Partial Differential Equations}}},\ \bibinfo
  {series} {Graduate Studies in Mathematics}, Vol.\ \bibinfo {volume} {123}\
  (\bibinfo  {publisher} {American Mathematical Society},\ \bibinfo {address}
  {USA},\ \bibinfo {year} {2011})\BibitemShut {NoStop}%
\bibitem [{\citenamefont {Hintz}(2019)}]{Hintz_2019}%
  \BibitemOpen
  \bibfield  {author} {\bibinfo {author} {\bibfnamefont {P.}~\bibnamefont
  {Hintz}},\ }\href {https://math.mit.edu/~phintz/18.157-S19/} {\bibinfo
  {title} {Introduction to microlocal analysis}} (\bibinfo {year} {2019}),\
  \bibinfo {note} {unpublished lecture notes taught at {M}{I}{T}}\BibitemShut
  {NoStop}%
\bibitem [{\citenamefont {Kohn}\ and\ \citenamefont
  {Nirenberg}(1965)}]{Kohn_CPAM_1965}%
  \BibitemOpen
  \bibfield  {author} {\bibinfo {author} {\bibfnamefont {J.~J.}\ \bibnamefont
  {Kohn}}\ and\ \bibinfo {author} {\bibfnamefont {L.}~\bibnamefont
  {Nirenberg}},\ }\bibfield  {title} {\bibinfo {title} {An algebra of
  pseudo-differential operators},\ }\href
  {https://doi.org/10.1002/cpa.3160180121} {\bibfield  {journal} {\bibinfo
  {journal} {Comm. Pure Appl. Math.}\ }\textbf {\bibinfo {volume} {18}},\
  \bibinfo {pages} {269 } (\bibinfo {year} {1965})}\BibitemShut {NoStop}%
\bibitem [{\citenamefont {Lax}(1957)}]{Lax_DukeMathJ_1957}%
  \BibitemOpen
  \bibfield  {author} {\bibinfo {author} {\bibfnamefont {P.~D.}\ \bibnamefont
  {Lax}},\ }\bibfield  {title} {\bibinfo {title} {Asymptotic solutions of
  oscillatory initial value problems},\ }\href
  {https://doi.org/10.1215/S0012-7094-57-02471-7} {\bibfield  {journal}
  {\bibinfo  {journal} {Duke Math. J.}\ }\textbf {\bibinfo {volume} {24}},\
  \bibinfo {pages} {627} (\bibinfo {year} {1957})}\BibitemShut {NoStop}%
\bibitem [{\citenamefont {Maslov}\ and\ \citenamefont
  {Fedoriuk}(1981)}]{Maslov_Springer_1981}%
  \BibitemOpen
  \bibfield  {author} {\bibinfo {author} {\bibfnamefont {V.~P.}\ \bibnamefont
  {Maslov}}\ and\ \bibinfo {author} {\bibfnamefont {M.}~\bibnamefont
  {Fedoriuk}},\ }\href {https://link.springer.com/book/9789027712196} {\emph
  {\bibinfo {title} {Semi-Classical Approximation in Quantum Mechanics}}},\
  Mathematical Physics and Applied Mathematics\ (\bibinfo  {publisher}
  {Springer},\ \bibinfo {address} {Netherlands},\ \bibinfo {year}
  {1981})\BibitemShut {NoStop}%
\bibitem [{\citenamefont {Egorov}(1969)}]{Egorov_UMN_1969}%
  \BibitemOpen
  \bibfield  {author} {\bibinfo {author} {\bibfnamefont {Y.~V.}\ \bibnamefont
  {Egorov}},\ }\bibfield  {title} {\bibinfo {title} {The canonical
  transformations of pseudodifferential operators},\ }\href
  {http://www.mathnet.ru/php/archive.phtml?wshow=paper&jrnid=rm&paperid=5554&option_lang=eng}
  {\bibfield  {journal} {\bibinfo  {journal} {Uspekhi Mat. Nauk}\ }\textbf
  {\bibinfo {volume} {24}},\ \bibinfo {pages} {235} (\bibinfo {year}
  {1969})}\BibitemShut {NoStop}%
\bibitem [{\citenamefont {Bolte}\ and\ \citenamefont
  {Glaser}(2004{\natexlab{a}})}]{Bolte_CMP_2004}%
  \BibitemOpen
  \bibfield  {author} {\bibinfo {author} {\bibfnamefont {J.}~\bibnamefont
  {Bolte}}\ and\ \bibinfo {author} {\bibfnamefont {R.}~\bibnamefont {Glaser}},\
  }\bibfield  {title} {\bibinfo {title} {A semiclassical {E}gorov theorem and
  quantum ergodicity for matrix valued operators},\ }\href
  {https://doi.org/10.1007/s00220-004-1064-0} {\bibfield  {journal} {\bibinfo
  {journal} {Commun. Math. Phys.}\ }\textbf {\bibinfo {volume} {247}},\
  \bibinfo {pages} {391} (\bibinfo {year} {2004}{\natexlab{a}})},\ \Eprint
  {https://arxiv.org/abs/0204018} {arXiv:0204018 [quant-ph]}\BibitemShut
  {NoStop}%
\bibitem [{\citenamefont {Kordyukov}(2005)}]{Kordyukov_MPAG_2005}%
  \BibitemOpen
  \bibfield  {author} {\bibinfo {author} {\bibfnamefont {Y.~A.}\ \bibnamefont
  {Kordyukov}},\ }\bibfield  {title} {\bibinfo {title} {Egorov's theorem for
  transversally elliptic operators on foliated manifolds and noncommutative
  geodesic flow},\ }\href {https://doi.org/10.1007/s11040-004-6495-5}
  {\bibfield  {journal} {\bibinfo  {journal} {Math. Phys. Anal. Geom.}\
  }\textbf {\bibinfo {volume} {8}},\ \bibinfo {pages} {97} (\bibinfo {year}
  {2005})},\ \Eprint {https://arxiv.org/abs/0407435} {arXiv:0407435 [math.DG]}\BibitemShut {NoStop}%
\bibitem [{\citenamefont {Kordyukov}(2007)}]{Kordyukov_JGP_2007}%
  \BibitemOpen
  \bibfield  {author} {\bibinfo {author} {\bibfnamefont {Y.~A.}\ \bibnamefont
  {Kordyukov}},\ }\bibfield  {title} {\bibinfo {title} {The {E}gorov theorem
  for transverse {D}irac-type operators on foliated manifolds},\ }\href
  {https://doi.org/https://doi.org/10.1016/j.geomphys.2007.08.002} {\bibfield
  {journal} {\bibinfo  {journal} {J. Geom. Phys.}\ }\textbf {\bibinfo {volume}
  {57}},\ \bibinfo {pages} {2345 } (\bibinfo {year} {2007})},\ \Eprint
  {https://arxiv.org/abs/0708.1660} {arXiv:0708.1660 [math.DG]}\BibitemShut
  {NoStop}%
\bibitem [{\citenamefont {Jakobson}\ and\ \citenamefont
  {Strohmaier}(2007)}]{Jakobson_CMP_2007}%
  \BibitemOpen
  \bibfield  {author} {\bibinfo {author} {\bibfnamefont {D.}~\bibnamefont
  {Jakobson}}\ and\ \bibinfo {author} {\bibfnamefont {A.}~\bibnamefont
  {Strohmaier}},\ }\bibfield  {title} {\bibinfo {title} {High energy limits of
  {L}aplace-type and {D}irac-type eigenfunctions and frame flows},\ }\href
  {https://doi.org/10.1007/s00220-006-0176-0} {\bibfield  {journal} {\bibinfo
  {journal} {Commun. Math. Phys.}\ }\textbf {\bibinfo {volume} {270}},\
  \bibinfo {pages} {813} (\bibinfo {year} {2007})},\ \Eprint
  {https://arxiv.org/abs/0607616} {arXiv:0607616 [math.SP]}\BibitemShut
  {NoStop}%
\bibitem [{\citenamefont {Taylor}(1981)}]{Taylor_PUP_1981}%
  \BibitemOpen
  \bibfield  {author} {\bibinfo {author} {\bibfnamefont {M.~E.}\ \bibnamefont
  {Taylor}},\ }\href
  {https://press.princeton.edu/books/hardcover/9780691629865/pseudodifferential-operators-pms-34}
  {\emph {\bibinfo {title} {Pseudodifferential Operators}}},\ Princeton
  Mathematical Series\ (\bibinfo  {publisher} {Princeton University Press},\
  \bibinfo {year} {1981})\BibitemShut {NoStop}%
\bibitem [{\citenamefont {Minguzzi}\ and\ \citenamefont
  {S\'{a}nchez}(2008)}]{Minguzzi_EMS_2008}%
  \BibitemOpen
  \bibfield  {author} {\bibinfo {author} {\bibfnamefont {E.}~\bibnamefont
  {Minguzzi}}\ and\ \bibinfo {author} {\bibfnamefont {M.}~\bibnamefont
  {S\'{a}nchez}},\ }\bibinfo {title} {The causal hierarchy of spacetimes},\ in\
  \href {https://doi.org/10.4171/051-1/9} {\emph {\bibinfo {booktitle} {Recent
  Developments in Pseudo-{R}iemannian Geometry}}},\ \bibinfo {series and
  number} {ESI Lectures in Mathematics and Physics},\ \bibinfo {editor} {edited
  by\ \bibinfo {editor} {\bibfnamefont {D.~V.}\ \bibnamefont {Alekseevsky}}\
  and\ \bibinfo {editor} {\bibfnamefont {H.}~\bibnamefont {Baum}}}\ (\bibinfo
  {publisher} {European Mathematical Society},\ \bibinfo {address} {Germany},\
  \bibinfo {year} {2008})\ pp.\ \bibinfo {pages} {299 -- 358},\ \Eprint
  {https://arxiv.org/abs/0609119} {arXiv:0609119 [gr-qc]}\BibitemShut
  {NoStop}%
\bibitem [{\citenamefont {Geroch}(1970{\natexlab{a}})}]{Geroch_JMP_1970}%
  \BibitemOpen
  \bibfield  {author} {\bibinfo {author} {\bibfnamefont {R.}~\bibnamefont
  {Geroch}},\ }\bibfield  {title} {\bibinfo {title} {Domain of dependence},\
  }\href {https://doi.org/10.1063/1.1665157} {\bibfield  {journal} {\bibinfo
  {journal} {J. Math. Phys.}\ }\textbf {\bibinfo {volume} {11}},\ \bibinfo
  {pages} {437 } (\bibinfo {year} {1970}{\natexlab{a}})}\BibitemShut {NoStop}%
\bibitem [{\citenamefont {Bernal}\ and\ \citenamefont
  {S{\'a}nchez}(2003)}]{Bernal_CMP_2003}%
  \BibitemOpen
  \bibfield  {author} {\bibinfo {author} {\bibfnamefont {A.~N.}\ \bibnamefont
  {Bernal}}\ and\ \bibinfo {author} {\bibfnamefont {M.}~\bibnamefont
  {S{\'a}nchez}},\ }\bibfield  {title} {\bibinfo {title} {On smooth {C}auchy
  hypersurfaces and {G}eroch's splitting theorem},\ }\href
  {https://doi.org/10.1007/s00220-003-0982-6} {\bibfield  {journal} {\bibinfo
  {journal} {Commun. Math. Phys.}\ }\textbf {\bibinfo {volume} {243}},\
  \bibinfo {pages} {461} (\bibinfo {year} {2003})},\ \Eprint
  {https://arxiv.org/abs/0306108} {arXiv:0306108 [gr-qc]}\BibitemShut
  {NoStop}%
\bibitem [{\citenamefont {Bernal}\ and\ \citenamefont
  {S{\'a}nchez}(2005)}]{Bernal_CMP_2005}%
  \BibitemOpen
  \bibfield  {author} {\bibinfo {author} {\bibfnamefont {A.~N.}\ \bibnamefont
  {Bernal}}\ and\ \bibinfo {author} {\bibfnamefont {M.}~\bibnamefont
  {S{\'a}nchez}},\ }\bibfield  {title} {\bibinfo {title} {Smoothness of time
  functions and the metric splitting of globally hyperbolic spacetimes},\
  }\href {https://doi.org/10.1007/s00220-005-1346-1} {\bibfield  {journal}
  {\bibinfo  {journal} {Commun. Math. Phys.}\ }\textbf {\bibinfo {volume}
  {257}},\ \bibinfo {pages} {43 } (\bibinfo {year} {2005})},\ \Eprint
  {https://arxiv.org/abs/0401112} {arXiv:0401112 [gr-qc]}\BibitemShut
  {NoStop}%
\bibitem [{\citenamefont {Bernal}\ and\ \citenamefont
  {S{\'a}nchez}(2006)}]{Bernal_LMP_2006}%
  \BibitemOpen
  \bibfield  {author} {\bibinfo {author} {\bibfnamefont {A.~N.}\ \bibnamefont
  {Bernal}}\ and\ \bibinfo {author} {\bibfnamefont {M.}~\bibnamefont
  {S{\'a}nchez}},\ }\bibfield  {title} {\bibinfo {title} {Further results on
  the smoothability of {C}auchy hypersurfaces and {C}auchy time functions},\
  }\href {https://doi.org/10.1007/s11005-006-0091-5} {\bibfield  {journal}
  {\bibinfo  {journal} {Lett. Math. Phys.}\ }\textbf {\bibinfo {volume} {77}},\
  \bibinfo {pages} {183 } (\bibinfo {year} {2006})},\ \Eprint
  {https://arxiv.org/abs/0512095} {arXiv:0512095 [gr-qc]}\BibitemShut
  {NoStop}%
\bibitem [{\citenamefont {Bernal}\ and\ \citenamefont
  {S{\'a}nchez}(2007)}]{Bernal_CQG_2007}%
  \BibitemOpen
  \bibfield  {author} {\bibinfo {author} {\bibfnamefont {A.~N.}\ \bibnamefont
  {Bernal}}\ and\ \bibinfo {author} {\bibfnamefont {M.}~\bibnamefont
  {S{\'a}nchez}},\ }\bibfield  {title} {\bibinfo {title} {Globally hyperbolic
  spacetimes can be defined as ‘causal’ instead of ‘strongly causal’},\
  }\href {https://doi.org/10.1088/0264-9381/24/3/n01} {\bibfield  {journal}
  {\bibinfo  {journal} {Class. Quant. Grav.}\ }\textbf {\bibinfo {volume}
  {24}},\ \bibinfo {pages} {745} (\bibinfo {year} {2007})},\ \Eprint
  {https://arxiv.org/abs/0611138} {arXiv:0611138 [gr-qc]}\BibitemShut
  {NoStop}%
\bibitem [{\citenamefont {B\"{a}r}\ \emph {et~al.}(2007)\citenamefont
  {B\"{a}r}, \citenamefont {Ginoux},\ and\ \citenamefont
  {Pf\"{a}ffle}}]{Baer_EMS_2007}%
  \BibitemOpen
  \bibfield  {author} {\bibinfo {author} {\bibfnamefont {C.}~\bibnamefont
  {B\"{a}r}}, \bibinfo {author} {\bibfnamefont {N.}~\bibnamefont {Ginoux}},\
  and\ \bibinfo {author} {\bibfnamefont {F.}~\bibnamefont {Pf\"{a}ffle}},\
  }\href {https://doi.org/10.4171/037} {\emph {\bibinfo {title} {Wave Equations
  on {L}orentzian Manifolds and Quantization}}},\ E{S}{I} Lectures in
  Mathematics and Physics\ (\bibinfo  {publisher} {European Mathematical
  Society},\ \bibinfo {address} {Germany},\ \bibinfo {year} {2007})\ \Eprint
  {https://arxiv.org/abs/0806.1036} {arXiv:0806.1036 [math.DG]}\BibitemShut
  {NoStop}%
\bibitem [{\citenamefont {Hawking}\ and\ \citenamefont
  {Ellis}(1973)}]{Hawking_CUP_1973}%
  \BibitemOpen
  \bibfield  {author} {\bibinfo {author} {\bibfnamefont {S.~W.}\ \bibnamefont
  {Hawking}}\ and\ \bibinfo {author} {\bibfnamefont {G.~F.~R.}\ \bibnamefont
  {Ellis}},\ }\href {https://doi.org/10.1017/CBO9780511524646} {\emph {\bibinfo
  {title} {The Large Scale Structure of Space-Time}}},\ Cambridge Monographs on
  Mathematical Physics\ (\bibinfo  {publisher} {Cambridge University Press},\
  \bibinfo {address} {Cambridge},\ \bibinfo {year} {1973})\BibitemShut
  {NoStop}%
\bibitem [{\citenamefont {Green}\ and\ \citenamefont
  {Wald}(2014)}]{Green_CQG_2014}%
  \BibitemOpen
  \bibfield  {author} {\bibinfo {author} {\bibfnamefont {S.~R.}\ \bibnamefont
  {Green}}\ and\ \bibinfo {author} {\bibfnamefont {R.~M.}\ \bibnamefont
  {Wald}},\ }\bibfield  {title} {\bibinfo {title} {How well is our universe
  described by an {FLRW} model?},\ }\href
  {https://doi.org/10.1088/0264-9381/31/23/234003} {\bibfield  {journal}
  {\bibinfo  {journal} {Class. Quantum Grav.}\ }\textbf {\bibinfo {volume}
  {31}},\ \bibinfo {pages} {234003} (\bibinfo {year} {2014})},\ \Eprint
  {https://arxiv.org/abs/1407.8084} {arXiv:1407.8084 [gr-qc]}\BibitemShut
  {NoStop}%
\bibitem [{\citenamefont {Baum}\ and\ \citenamefont
  {Kath}(1996)}]{Baum_AGAG_1996}%
  \BibitemOpen
  \bibfield  {author} {\bibinfo {author} {\bibfnamefont {H.}~\bibnamefont
  {Baum}}\ and\ \bibinfo {author} {\bibfnamefont {I.}~\bibnamefont {Kath}},\
  }\bibfield  {title} {\bibinfo {title} {Normally hyperbolic operators, the
  {H}uygens property and conformal geometry},\ }\href
  {https://doi.org/10.1007/BF00129896} {\bibfield  {journal} {\bibinfo
  {journal} {Ann. Glob. Anal. Geom.}\ }\textbf {\bibinfo {volume} {14}},\
  \bibinfo {pages} {315} (\bibinfo {year} {1996})}\BibitemShut {NoStop}%
\bibitem [{\citenamefont {Berline}\ \emph {et~al.}(2004)\citenamefont
  {Berline}, \citenamefont {Getzler},\ and\ \citenamefont
  {Vergne}}]{Berline_Springer_2004}%
  \BibitemOpen
  \bibfield  {author} {\bibinfo {author} {\bibfnamefont {N.}~\bibnamefont
  {Berline}}, \bibinfo {author} {\bibfnamefont {E.}~\bibnamefont {Getzler}},\
  and\ \bibinfo {author} {\bibfnamefont {M.}~\bibnamefont {Vergne}},\ }\href
  {https://www.springer.com/gp/book/9783540200628} {\emph {\bibinfo {title}
  {Heat Kernels and {D}irac Operators}}},\ \bibinfo {edition} {2nd}\ ed.,\
  Grundlehren Text Editions\ (\bibinfo  {publisher} {Springer-Verlag},\
  \bibinfo {address} {Berlin Heidelberg},\ \bibinfo {year} {2004})\BibitemShut
  {NoStop}%
\bibitem [{\citenamefont {Wald}(2007)}]{Wald_Chicago_1984}%
  \BibitemOpen
  \bibfield  {author} {\bibinfo {author} {\bibfnamefont {R.~M.}\ \bibnamefont
  {Wald}},\ }\href
  {https://press.uchicago.edu/ucp/books/book/chicago/G/bo5952261.html} {\emph
  {\bibinfo {title} {General Relativity}}}\ (\bibinfo  {publisher} {University
  Chicago Press},\ \bibinfo {address} {USA},\ \bibinfo {year} {1984, {I}ndian
  reprint 2007})\BibitemShut {NoStop}%
\bibitem [{\citenamefont {Laptev}\ \emph {et~al.}(1994)\citenamefont {Laptev},
  \citenamefont {Safarov},\ and\ \citenamefont {Vassiliev}}]{Laptev_CPAM_1994}%
  \BibitemOpen
  \bibfield  {author} {\bibinfo {author} {\bibfnamefont {A.}~\bibnamefont
  {Laptev}}, \bibinfo {author} {\bibfnamefont {Y.}~\bibnamefont {Safarov}},\
  and\ \bibinfo {author} {\bibfnamefont {D.}~\bibnamefont {Vassiliev}},\
  }\bibfield  {title} {\bibinfo {title} {On global representation of
  {L}agrangian distributions and solutions of hyperbolic equations},\ }\href
  {https://doi.org/10.1002/cpa.3160471102} {\bibfield  {journal} {\bibinfo
  {journal} {Commun. Pure Appl. Math.}\ }\textbf {\bibinfo {volume} {47}},\
  \bibinfo {pages} {1411} (\bibinfo {year} {1994})}\BibitemShut {NoStop}%
\bibitem [{\citenamefont {Kamber}\ and\ \citenamefont
  {Tondeur}(1975)}]{Kamber_Springer_1975}%
  \BibitemOpen
  \bibfield  {author} {\bibinfo {author} {\bibfnamefont {F.~W.}\ \bibnamefont
  {Kamber}}\ and\ \bibinfo {author} {\bibfnamefont {P.}~\bibnamefont
  {Tondeur}},\ }\href {https://doi.org/10.1007/BFb0081558} {\emph {\bibinfo
  {title} {Foliated Bundles and Characteristic Classes}}},\ \bibinfo {series}
  {Lecture Notes in Mathematics}, Vol.\ \bibinfo {volume} {493}\ (\bibinfo
  {publisher} {Springer-Verlag},\ \bibinfo {address} {Berlin Heidelberg},\
  \bibinfo {year} {1975})\BibitemShut {NoStop}%
\bibitem [{\citenamefont {Hintz}(2017)}]{Hintz_JST_2017}%
  \BibitemOpen
  \bibfield  {author} {\bibinfo {author} {\bibfnamefont {P.}~\bibnamefont
  {Hintz}},\ }\bibfield  {title} {\bibinfo {title} {Resonance expansions for
  tensor-valued waves on asymptotically {K}err–de {S}itter spaces},\ }\href
  {https://doi.org/10.4171/jst/171} {\bibfield  {journal} {\bibinfo  {journal}
  {J. Spectr. Theory}\ }\textbf {\bibinfo {volume} {7}},\ \bibinfo {pages}
  {519–557} (\bibinfo {year} {2017})},\ \Eprint
  {https://arxiv.org/abs/1502.03183} {arXiv:1502.03183 [math.AP]}\BibitemShut
  {NoStop}%
\bibitem [{\citenamefont {Silva}(2007)}]{Silva_PublMat_2007}%
  \BibitemOpen
  \bibfield  {author} {\bibinfo {author} {\bibfnamefont {J.~D.}\ \bibnamefont
  {Silva}},\ }\bibfield  {title} {\bibinfo {title} {An accuracy improvement in
  {E}gorov’s theorem},\ }\href {https://doi.org/10.5565/PUBLMAT_51107_05}
  {\bibfield  {journal} {\bibinfo  {journal} {Publ. Mat.}\ }\textbf {\bibinfo
  {volume} {51}},\ \bibinfo {pages} {77 } (\bibinfo {year} {2007})}\BibitemShut
  {NoStop}%
\bibitem [{\citenamefont {Junker}\ and\ \citenamefont
  {Schrohe}(2002)}]{Junker_AHP_2002}%
  \BibitemOpen
  \bibfield  {author} {\bibinfo {author} {\bibfnamefont {W.}~\bibnamefont
  {Junker}}\ and\ \bibinfo {author} {\bibfnamefont {E.}~\bibnamefont
  {Schrohe}},\ }\bibfield  {title} {\bibinfo {title} {Adiabatic vacuum states
  on general space-time manifolds: Definition, construction, and physical
  properties},\ }\href {https://doi.org/10.1007/s000230200001} {\bibfield
  {journal} {\bibinfo  {journal} {Ann. Henri Poincar\'{e}}\ }\textbf {\bibinfo
  {volume} {3}},\ \bibinfo {pages} {1113} (\bibinfo {year} {2002})},\ \Eprint
  {https://arxiv.org/abs/0109010} {arXiv:0109010 [math-ph]}\BibitemShut
  {NoStop}%
\bibitem [{\citenamefont {Jubin}\ and\ \citenamefont
  {Schapira}(2016)}]{Jubin_LMP_2016}%
  \BibitemOpen
  \bibfield  {author} {\bibinfo {author} {\bibfnamefont {B.}~\bibnamefont
  {Jubin}}\ and\ \bibinfo {author} {\bibfnamefont {P.}~\bibnamefont
  {Schapira}},\ }\bibfield  {title} {\bibinfo {title} {Sheaves and {D}-modules
  on {L}orentzian manifolds},\ }\href
  {https://doi.org/10.1007/s11005-016-0832-z} {\bibfield  {journal} {\bibinfo
  {journal} {Lett. Math. Phys.}\ }\textbf {\bibinfo {volume} {106}},\ \bibinfo
  {pages} {607} (\bibinfo {year} {2016})},\ \Eprint
  {https://arxiv.org/abs/1510.01499} {arXiv:1510.01499 [math.AG]}\BibitemShut
  {NoStop}%
\bibitem [{\citenamefont {Penrose}(1972)}]{Penrose_1972}%
  \BibitemOpen
  \bibfield  {author} {\bibinfo {author} {\bibfnamefont {R.}~\bibnamefont
  {Penrose}},\ }\bibinfo {title} {On the nature of quantum geometry},\ in\
  \href@noop {} {\emph {\bibinfo {booktitle} {Magic Without Magic}}},\ \bibinfo
  {editor} {edited by\ \bibinfo {editor} {\bibfnamefont {J.~R.}\ \bibnamefont
  {Klauder}}}\ (\bibinfo  {publisher} {W. H. Freeman},\ \bibinfo {address} {San
  Francisco},\ \bibinfo {year} {1972})\BibitemShut {NoStop}%
\bibitem [{\citenamefont {Khesin}\ and\ \citenamefont
  {Tabachnikov}(2009)}]{Khesin_AdvMath_2009}%
  \BibitemOpen
  \bibfield  {author} {\bibinfo {author} {\bibfnamefont {B.}~\bibnamefont
  {Khesin}}\ and\ \bibinfo {author} {\bibfnamefont {S.}~\bibnamefont
  {Tabachnikov}},\ }\bibfield  {title} {\bibinfo {title} {Pseudo-{R}iemannian
  geodesics and billiards},\ }\href
  {https://doi.org/https://doi.org/10.1016/j.aim.2009.02.010} {\bibfield
  {journal} {\bibinfo  {journal} {Adv. Math.}\ }\textbf {\bibinfo {volume}
  {221}},\ \bibinfo {pages} {1364 } (\bibinfo {year} {2009})},\ \Eprint
  {https://arxiv.org/abs/0608620} {arXiv:0608620 [math.DG]}\BibitemShut
  {NoStop}%
\bibitem [{\citenamefont {B{\"a}r}\ and\ \citenamefont
  {Ginoux}(2012)}]{Baer_Springer_2012}%
  \BibitemOpen
  \bibfield  {author} {\bibinfo {author} {\bibfnamefont {C.}~\bibnamefont
  {B{\"a}r}}\ and\ \bibinfo {author} {\bibfnamefont {N.}~\bibnamefont
  {Ginoux}},\ }\bibinfo {title} {Classical and quantum fields on {L}orentzian
  manifolds},\ in\ \href {https://doi.org/10.1007/978-3-642-22842-1_12} {\emph
  {\bibinfo {booktitle} {Global Differential Geometry}}},\ \bibinfo {editor}
  {edited by\ \bibinfo {editor} {\bibfnamefont {C.}~\bibnamefont {B{\"a}r}},
  \bibinfo {editor} {\bibfnamefont {J.}~\bibnamefont {Lohkamp}},\ and\ \bibinfo
  {editor} {\bibfnamefont {M.}~\bibnamefont {Schwarz}}}\ (\bibinfo  {publisher}
  {Springer},\ \bibinfo {address} {Berlin, Heidelberg},\ \bibinfo {year}
  {2012})\ pp.\ \bibinfo {pages} {359--400},\ \Eprint
  {https://arxiv.org/abs/1104.1158} {arXiv:1104.1158 [math-ph]}\BibitemShut
  {NoStop}%
\bibitem [{\citenamefont {Baum}(1981)}]{Baum_1981}%
  \BibitemOpen
  \bibfield  {author} {\bibinfo {author} {\bibfnamefont {H.}~\bibnamefont
  {Baum}},\ }\href@noop {} {\emph {\bibinfo {title} {Spin-Strukturen und
  {D}irac-Operatoren über Pseudo-{R}iemannschen Mannigfaltigkeiten}}},\
  Teubner Texte zur Mathematik\ (\bibinfo {address} {Leipzig},\ \bibinfo {year}
  {1981})\BibitemShut {NoStop}%
\bibitem [{\citenamefont {B{\"a}r}\ \emph {et~al.}(2005)\citenamefont
  {B{\"a}r}, \citenamefont {Gauduchon},\ and\ \citenamefont
  {Moroianu}}]{Baer_MathZ_2005}%
  \BibitemOpen
  \bibfield  {author} {\bibinfo {author} {\bibfnamefont {C.}~\bibnamefont
  {B{\"a}r}}, \bibinfo {author} {\bibfnamefont {P.}~\bibnamefont {Gauduchon}},\
  and\ \bibinfo {author} {\bibfnamefont {A.}~\bibnamefont {Moroianu}},\
  }\bibfield  {title} {\bibinfo {title} {Generalized cylinders in
  semi-{R}iemannian and spin geometry},\ }\href
  {https://doi.org/10.1007/s00209-004-0718-0} {\bibfield  {journal} {\bibinfo
  {journal} {Math. Z.}\ }\textbf {\bibinfo {volume} {249}},\ \bibinfo {pages}
  {545} (\bibinfo {year} {2005})},\ \Eprint {https://arxiv.org/abs/0303095}
  {arXiv:0303095 [math.DG]}\BibitemShut {NoStop}%
\bibitem [{\citenamefont {Geroch}(1968)}]{Geroch_JMP_1968}%
  \BibitemOpen
  \bibfield  {author} {\bibinfo {author} {\bibfnamefont {R.}~\bibnamefont
  {Geroch}},\ }\bibfield  {title} {\bibinfo {title} {Spinor structure of
  space-times in general relativity. {I}},\ }\href
  {https://doi.org/10.1063/1.1664507} {\bibfield  {journal} {\bibinfo
  {journal} {J. Math. Phys.}\ }\textbf {\bibinfo {volume} {9}},\ \bibinfo
  {pages} {1739 } (\bibinfo {year} {1968})}\BibitemShut {NoStop}%
\bibitem [{\citenamefont {Geroch}(1970{\natexlab{b}})}]{Geroch_spin_JMP_1970}%
  \BibitemOpen
  \bibfield  {author} {\bibinfo {author} {\bibfnamefont {R.}~\bibnamefont
  {Geroch}},\ }\bibfield  {title} {\bibinfo {title} {Spinor structure of
  space-times in general relativity. {I}{I}},\ }\href
  {https://doi.org/10.1063/1.1665067} {\bibfield  {journal} {\bibinfo
  {journal} {J. Math. Phys.}\ }\textbf {\bibinfo {volume} {11}},\ \bibinfo
  {pages} {343 } (\bibinfo {year} {1970}{\natexlab{b}})}\BibitemShut {NoStop}%
\bibitem [{\citenamefont {Clarke}(1971)}]{Clarke_GRG_1971}%
  \BibitemOpen
  \bibfield  {author} {\bibinfo {author} {\bibfnamefont {C.~J.~S.}\
  \bibnamefont {Clarke}},\ }\bibfield  {title} {\bibinfo {title} {Magnetic
  charge, holonomy and characteristic classes: Illustrations of the methods of
  topology in relativity},\ }\href {https://doi.org/10.1007/BF02450517}
  {\bibfield  {journal} {\bibinfo  {journal} {Gen. Relat. Gravit.}\ }\textbf
  {\bibinfo {volume} {2}},\ \bibinfo {pages} {43 } (\bibinfo {year}
  {1971})}\BibitemShut {NoStop}%
\bibitem [{\citenamefont {Alagia}\ and\ \citenamefont
  {S\'{a}nchez}(1985)}]{Alagia_RUMA_1985}%
  \BibitemOpen
  \bibfield  {author} {\bibinfo {author} {\bibfnamefont {H.~R.}\ \bibnamefont
  {Alagia}}\ and\ \bibinfo {author} {\bibfnamefont {C.~U.}\ \bibnamefont
  {S\'{a}nchez}},\ }\bibfield  {title} {\bibinfo {title} {Spin structures on
  pseudo-{R}iemannian manifolds},\ }\href
  {https://inmabb.criba.edu.ar/revuma/pdf/v32n1/p064-078.pdf} {\bibfield
  {journal} {\bibinfo  {journal} {Revista de la Uni\'{o}n Matem\'{a}tica
  Argentina}\ }\textbf {\bibinfo {volume} {32}},\ \bibinfo {pages} {64 }
  (\bibinfo {year} {1985})}\BibitemShut {NoStop}%
\bibitem [{\citenamefont {Isham}(1978)}]{Isham_spin_PRSL_1978}%
  \BibitemOpen
  \bibfield  {author} {\bibinfo {author} {\bibfnamefont {C.~J.}\ \bibnamefont
  {Isham}},\ }\bibfield  {title} {\bibinfo {title} {Spinor fields in four
  dimensional space-time},\ }\href {https://doi.org/10.1098/rspa.1978.0219}
  {\bibfield  {journal} {\bibinfo  {journal} {Proc. R. Soc. Lond. A.}\ }\textbf
  {\bibinfo {volume} {364}},\ \bibinfo {pages} {591 } (\bibinfo {year}
  {1978})}\BibitemShut {NoStop}%
\bibitem [{\citenamefont {M\"{u}hlhoff}(2011)}]{Muehlhoff_JMP_2011}%
  \BibitemOpen
  \bibfield  {author} {\bibinfo {author} {\bibfnamefont {R.}~\bibnamefont
  {M\"{u}hlhoff}},\ }\bibfield  {title} {\bibinfo {title} {Cauchy problem and
  {G}reen's functions for first order differential operators and algebraic
  quantization},\ }\href {https://doi.org/10.1063/1.3530846} {\bibfield
  {journal} {\bibinfo  {journal} {J. Math. Phys.}\ }\textbf {\bibinfo {volume}
  {52}},\ \bibinfo {pages} {022303} (\bibinfo {year} {2011})},\ \Eprint
  {https://arxiv.org/abs/1001.4091} {arXiv:1001.4091 [math-ph]}\BibitemShut
  {NoStop}%
\bibitem [{\citenamefont {Dimock}(1982)}]{Dimock_AMS_1982}%
  \BibitemOpen
  \bibfield  {author} {\bibinfo {author} {\bibfnamefont {J.}~\bibnamefont
  {Dimock}},\ }\bibfield  {title} {\bibinfo {title} {Dirac quantum fields on a
  manifold},\ }\href
  {https://doi.org/https://doi.org/10.1090/S0002-9947-1982-0637032-8}
  {\bibfield  {journal} {\bibinfo  {journal} {Trans. Am. Math. Soc.}\ }\textbf
  {\bibinfo {volume} {269}},\ \bibinfo {pages} {133} (\bibinfo {year}
  {1982})}\BibitemShut {NoStop}%
\bibitem [{\citenamefont {Taylor}(2011)}]{Taylor_I_Springer_2011}%
  \BibitemOpen
  \bibfield  {author} {\bibinfo {author} {\bibfnamefont {M.}~\bibnamefont
  {Taylor}},\ }\href {https://doi.org/10.1007/978-1-4419-7055-8} {\emph
  {\bibinfo {title} {Partial Differential Equations {I}: Basic Theory}}},\
  \bibinfo {edition} {2nd}\ ed.,\ \bibinfo {series} {Applied Mathematical
  Sciences}, Vol.\ \bibinfo {volume} {115}\ (\bibinfo  {publisher}
  {Springer-Verlag},\ \bibinfo {address} {New York},\ \bibinfo {year}
  {2011})\BibitemShut {NoStop}%
\bibitem [{\citenamefont {Schr{\"o}dinger}(2020)}]{Schroedinger_GRG_2020}%
  \BibitemOpen
  \bibfield  {author} {\bibinfo {author} {\bibfnamefont {E.}~\bibnamefont
  {Schr{\"o}dinger}},\ }\bibfield  {title} {\bibinfo {title} {Republication of:
  Dirac electron in the gravitational field {I}},\ }\href
  {https://doi.org/https://doi.org/10.1007/s10714-019-2626-y} {\bibfield
  {journal} {\bibinfo  {journal} {Gen. Relativ. Gravit.}\ }\textbf {\bibinfo
  {volume} {52}},\ \bibinfo {pages} {4} (\bibinfo {year} {2020})}\BibitemShut
  {NoStop}%
\bibitem [{\citenamefont {Lichnerowicz}(1963)}]{Lichnerowicz_1963}%
  \BibitemOpen
  \bibfield  {author} {\bibinfo {author} {\bibfnamefont {A.}~\bibnamefont
  {Lichnerowicz}},\ }\bibfield  {title} {\bibinfo {title} {Spineurs
  harmoniques},\ }\href {http://gallica.bnf.fr/ark:/12148/bpt6k4007z/f7.image}
  {\bibfield  {journal} {\bibinfo  {journal} {Comptes rendus de
  l'{A}cad\'{e}mie des {S}ciences (C. R. Acad. Sci. Paris)}\ }\textbf {\bibinfo
  {volume} {257}},\ \bibinfo {pages} {7 } (\bibinfo {year} {1963})}\BibitemShut
  {NoStop}%
\bibitem [{\citenamefont {Rarita}\ and\ \citenamefont
  {Schwinger}(1941)}]{Rarita_PR_1941}%
  \BibitemOpen
  \bibfield  {author} {\bibinfo {author} {\bibfnamefont {W.}~\bibnamefont
  {Rarita}}\ and\ \bibinfo {author} {\bibfnamefont {J.}~\bibnamefont
  {Schwinger}},\ }\bibfield  {title} {\bibinfo {title} {On a theory of
  particles with half-integral spin},\ }\href
  {https://doi.org/10.1103/PhysRev.60.61} {\bibfield  {journal} {\bibinfo
  {journal} {Phys. Rev.}\ }\textbf {\bibinfo {volume} {60}},\ \bibinfo {pages}
  {61} (\bibinfo {year} {1941})}\BibitemShut {NoStop}%
\bibitem [{\citenamefont {Homma}\ and\ \citenamefont
  {Semmelmann}(2019)}]{Homma_CMP_2019}%
  \BibitemOpen
  \bibfield  {author} {\bibinfo {author} {\bibfnamefont {Y.}~\bibnamefont
  {Homma}}\ and\ \bibinfo {author} {\bibfnamefont {U.}~\bibnamefont
  {Semmelmann}},\ }\bibfield  {title} {\bibinfo {title} {The kernel of the
  {R}arita–{S}chwinger operator on {R}iemannian spin manifolds},\ }\href
  {https://doi.org/10.1007/s00220-019-03324-8} {\bibfield  {journal} {\bibinfo
  {journal} {Commun. Math. Phys.}\ }\textbf {\bibinfo {volume} {370}},\
  \bibinfo {pages} {853–871} (\bibinfo {year} {2019})},\ \Eprint
  {https://arxiv.org/abs/1804.10602} {arXiv:1804.10602 [math.DG]}\BibitemShut
  {NoStop}%
\bibitem [{\citenamefont {Wang}(1991)}]{Wang_IndianaUMJ_1991}%
  \BibitemOpen
  \bibfield  {author} {\bibinfo {author} {\bibfnamefont {M.~Y.}\ \bibnamefont
  {Wang}},\ }\bibfield  {title} {\bibinfo {title} {Preserving parallel spinors
  under metric deformations},\ }\href {http://www.jstor.org/stable/24896310}
  {\bibfield  {journal} {\bibinfo  {journal} {Indiana Univ. Math. J.}\ }\textbf
  {\bibinfo {volume} {40}},\ \bibinfo {pages} {815} (\bibinfo {year}
  {1991})}\BibitemShut {NoStop}%
\bibitem [{\citenamefont {Sorokin}(2005)}]{Sorokin_AIP_2005}%
  \BibitemOpen
  \bibfield  {author} {\bibinfo {author} {\bibfnamefont {D.}~\bibnamefont
  {Sorokin}},\ }\bibfield  {title} {\bibinfo {title} {Introduction to the
  classical theory of higher spins},\ }\href
  {https://doi.org/10.1063/1.1923335} {\bibfield  {journal} {\bibinfo
  {journal} {AIP Conf. Proc.}\ }\textbf {\bibinfo {volume} {767}},\ \bibinfo
  {pages} {172} (\bibinfo {year} {2005})},\ \Eprint
  {https://arxiv.org/abs/0405069} {arXiv:0405069 [hep-th]}\BibitemShut
  {NoStop}%
\bibitem [{\citenamefont {Rahman}\ and\ \citenamefont
  {Taronna}(2015)}]{Rahman_2015}%
  \BibitemOpen
  \bibfield  {author} {\bibinfo {author} {\bibfnamefont {R.}~\bibnamefont
  {Rahman}}\ and\ \bibinfo {author} {\bibfnamefont {M.}~\bibnamefont
  {Taronna}},\ }\bibfield  {title} {\bibinfo {title} {From higher spins to
  strings: A primer},\ }\href@noop {} {\  (\bibinfo {year} {2015})},\ \Eprint
  {https://arxiv.org/abs/1512.07932} {arXiv:1512.07932 [hep-th]}\BibitemShut
  {NoStop}%
\bibitem [{\citenamefont {Gibbons}(1976)}]{Gibbons_JPA_1976}%
  \BibitemOpen
  \bibfield  {author} {\bibinfo {author} {\bibfnamefont {G.~W.}\ \bibnamefont
  {Gibbons}},\ }\bibfield  {title} {\bibinfo {title} {A note on the
  {R}arita-{S}chwinger equation in a gravitational background},\ }\href
  {https://doi.org/10.1088/0305-4470/9/1/019} {\bibfield  {journal} {\bibinfo
  {journal} {J. Phys. A: Math. Gen.}\ }\textbf {\bibinfo {volume} {9}},\
  \bibinfo {pages} {145} (\bibinfo {year} {1976})}\BibitemShut {NoStop}%
\bibitem [{\citenamefont {Hack}\ and\ \citenamefont
  {Makedonski}(2013)}]{Hack_PLB_2013}%
  \BibitemOpen
  \bibfield  {author} {\bibinfo {author} {\bibfnamefont {T.-P.}\ \bibnamefont
  {Hack}}\ and\ \bibinfo {author} {\bibfnamefont {M.}~\bibnamefont
  {Makedonski}},\ }\bibfield  {title} {\bibinfo {title} {A no-go theorem for
  the consistent quantization of spin-3/2 fields on general curved
  spacetimes},\ }\href
  {https://doi.org/https://doi.org/10.1016/j.physletb.2012.11.033} {\bibfield
  {journal} {\bibinfo  {journal} {Phys. Lett. B}\ }\textbf {\bibinfo {volume}
  {718}},\ \bibinfo {pages} {1465 } (\bibinfo {year} {2013})},\ \Eprint
  {https://arxiv.org/abs/1106.6327} {arXiv:1106.6327 [hep-th]}\BibitemShut
  {NoStop}%
\bibitem [{\citenamefont {Dirac}(1936)}]{Dirac_PRSA_1936}%
  \BibitemOpen
  \bibfield  {author} {\bibinfo {author} {\bibfnamefont {P.~A.~M.}\
  \bibnamefont {Dirac}},\ }\bibfield  {title} {\bibinfo {title} {Relativistic
  wave equations},\ }\href {https://doi.org/10.1098/rspa.1936.0111} {\bibfield
  {journal} {\bibinfo  {journal} {Proc. Royal Soc. A}\ }\textbf {\bibinfo
  {volume} {155}},\ \bibinfo {pages} {447} (\bibinfo {year}
  {1936})}\BibitemShut {NoStop}%
\bibitem [{\citenamefont {Illge}\ and\ \citenamefont
  {Schimming}(1999)}]{Illge_AnnPhys_1999}%
  \BibitemOpen
  \bibfield  {author} {\bibinfo {author} {\bibfnamefont {R.}~\bibnamefont
  {Illge}}\ and\ \bibinfo {author} {\bibfnamefont {R.}~\bibnamefont
  {Schimming}},\ }\bibfield  {title} {\bibinfo {title} {Consistent field
  equations for higher spin on curved spacetimes},\ }\href
  {https://doi.org/10.1002/(SICI)1521-3889(199904)8:4<319::AID-ANDP319>3.0.CO;2-3}
  {\bibfield  {journal} {\bibinfo  {journal} {Ann. Phys.}\ }\textbf {\bibinfo
  {volume} {8}},\ \bibinfo {pages} {319} (\bibinfo {year} {1999})}\BibitemShut
  {NoStop}%
\bibitem [{\citenamefont {Buchdahl}(1982)}]{Buchdahl_JPA_1982}%
  \BibitemOpen
  \bibfield  {author} {\bibinfo {author} {\bibfnamefont {H.~A.}\ \bibnamefont
  {Buchdahl}},\ }\bibfield  {title} {\bibinfo {title} {On the compatibility of
  relativistic wave equations in {R}iemann spaces. {II}},\ }\href
  {https://doi.org/10.1088/0305-4470/15/4/012} {\bibfield  {journal} {\bibinfo
  {journal} {J. Phys. A: Mathematical and General}\ }\textbf {\bibinfo {volume}
  {15}},\ \bibinfo {pages} {1057} (\bibinfo {year} {1982})}\BibitemShut
  {NoStop}%
\bibitem [{\citenamefont {W\"{u}nsch}(1985)}]{Wuensch_GRG_1985}%
  \BibitemOpen
  \bibfield  {author} {\bibinfo {author} {\bibfnamefont {V.}~\bibnamefont
  {W\"{u}nsch}},\ }\bibfield  {title} {\bibinfo {title} {Cauchy's problem and
  {H}uygens' principle for relativistic higher spin wave equations in an
  arbitrary curved space-time},\ }\href {https://doi.org/10.1007/BF00760104}
  {\bibfield  {journal} {\bibinfo  {journal} {Gen Relat Gravit}\ }\textbf
  {\bibinfo {volume} {17}},\ \bibinfo {pages} {15} (\bibinfo {year}
  {1985})}\BibitemShut {NoStop}%
\bibitem [{\citenamefont {Illge}(1992)}]{Illge_ZAA_1992}%
  \BibitemOpen
  \bibfield  {author} {\bibinfo {author} {\bibfnamefont {R.}~\bibnamefont
  {Illge}},\ }\bibfield  {title} {\bibinfo {title} {On massless fields with
  arbitrary spin},\ }\href {https://doi.org/10.4171/ZAA/629} {\bibfield
  {journal} {\bibinfo  {journal} {Z. Anal. Anwend.}\ }\textbf {\bibinfo
  {volume} {11}},\ \bibinfo {pages} {25 } (\bibinfo {year} {1992})}\BibitemShut
  {NoStop}%
\bibitem [{\citenamefont {Illge}(1993)}]{Illge_CMP_1993}%
  \BibitemOpen
  \bibfield  {author} {\bibinfo {author} {\bibfnamefont {R.}~\bibnamefont
  {Illge}},\ }\bibfield  {title} {\bibinfo {title} {Massive fields of arbitrary
  spin in curved space-times},\ }\href
  {https://projecteuclid.org:443/euclid.cmp/1104254357} {\bibfield  {journal}
  {\bibinfo  {journal} {Comm. Math. Phys.}\ }\textbf {\bibinfo {volume}
  {158}},\ \bibinfo {pages} {433} (\bibinfo {year} {1993})}\BibitemShut
  {NoStop}%
\bibitem [{\citenamefont {Frauendiener}\ and\ \citenamefont
  {Sparling}(1999)}]{Frauendiener_JGP_1999}%
  \BibitemOpen
  \bibfield  {author} {\bibinfo {author} {\bibfnamefont {J.}~\bibnamefont
  {Frauendiener}}\ and\ \bibinfo {author} {\bibfnamefont {G.~A.}\ \bibnamefont
  {Sparling}},\ }\bibfield  {title} {\bibinfo {title} {On a class of consistent
  linear higher spin equations on curved manifolds},\ }\href
  {https://doi.org/https://doi.org/10.1016/S0393-0440(98)00050-3} {\bibfield
  {journal} {\bibinfo  {journal} {J. Geom. Phys.}\ }\textbf {\bibinfo {volume}
  {30}},\ \bibinfo {pages} {54 } (\bibinfo {year} {1999})},\ \Eprint
  {https://arxiv.org/abs/9511036} {arXiv:9511036 [gr-qc]}\BibitemShut
  {NoStop}%
\bibitem [{\citenamefont {Beem}\ \emph {et~al.}(1996)\citenamefont {Beem},
  \citenamefont {Ehrlich},\ and\ \citenamefont {Easley}}]{Beem_CRC_1996}%
  \BibitemOpen
  \bibfield  {author} {\bibinfo {author} {\bibfnamefont {J.~K.}\ \bibnamefont
  {Beem}}, \bibinfo {author} {\bibfnamefont {P.~E.}\ \bibnamefont {Ehrlich}},\
  and\ \bibinfo {author} {\bibfnamefont {K.~L.}\ \bibnamefont {Easley}},\
  }\href {https://doi.org/https://doi.org/10.1201/9780203753125} {\emph
  {\bibinfo {title} {Global {L}orentzian Geometry}}},\ \bibinfo {edition}
  {2nd}\ ed.,\ Chapman \& Hall/CRC Pure and Applied Mathematics\ (\bibinfo
  {publisher} {CRC Press},\ \bibinfo {address} {New York},\ \bibinfo {year}
  {1996})\BibitemShut {NoStop}%
\bibitem [{\citenamefont {Choquet-Bruhat}(2008)}]{Choquet_Bruhat_OUP_2008}%
  \BibitemOpen
  \bibfield  {author} {\bibinfo {author} {\bibfnamefont {Y.}~\bibnamefont
  {Choquet-Bruhat}},\ }\href
  {https://doi.org/10.1093/acprof:oso/9780199230723.001.0001} {\emph {\bibinfo
  {title} {General Relativity and the {E}instein Equations}}},\ Oxford
  Mathematical Monographs\ (\bibinfo  {publisher} {Oxford University Press},\
  \bibinfo {address} {UK},\ \bibinfo {year} {2008})\BibitemShut {NoStop}%
\bibitem [{\citenamefont {Leray}(1953)}]{Leray_Princeton_1953}%
  \BibitemOpen
  \bibfield  {author} {\bibinfo {author} {\bibfnamefont {J.}~\bibnamefont
  {Leray}},\ }\href@noop {} {\emph {\bibinfo {title} {Hyperbolic Differential
  Equations}}}\ (\bibinfo  {publisher} {Mimeographed notes},\ \bibinfo
  {address} {Princeton},\ \bibinfo {year} {1953})\BibitemShut {NoStop}%
\bibitem [{\citenamefont {Penrose}(1969)}]{Penrose_RivNuovoCim_1969}%
  \BibitemOpen
  \bibfield  {author} {\bibinfo {author} {\bibfnamefont {R.}~\bibnamefont
  {Penrose}},\ }\bibfield  {title} {\bibinfo {title} {Gravitational collapse:
  The role of general relativity},\ }\href@noop {} {\bibfield  {journal}
  {\bibinfo  {journal} {Riv. Nuovo Cim.}\ }\textbf {\bibinfo {volume} {1}},\
  \bibinfo {pages} {252} (\bibinfo {year} {1969})}\BibitemShut {NoStop}%
\bibitem [{\citenamefont {S\'{a}nchez}(2011)}]{Sanchez_AMS_2011}%
  \BibitemOpen
  \bibfield  {author} {\bibinfo {author} {\bibfnamefont {M.}~\bibnamefont
  {S\'{a}nchez}},\ }\bibinfo {title} {Recent progress on the notion of global
  hyperbolicity},\ in\ \href {https://bookstore.ams.org/amsip-49} {\emph
  {\bibinfo {booktitle} {Advances in {L}orentzian Geometry: Proceedings of the
  {L}orentzian Geometry Conference in {B}erlin}}},\ \bibinfo {series} {AMS/IP
  Studies in Advanced Mathematics}, Vol.~\bibinfo {volume} {49},\ \bibinfo
  {editor} {edited by\ \bibinfo {editor} {\bibfnamefont {M.}~\bibnamefont
  {Plaue}}, \bibinfo {editor} {\bibfnamefont {A.}~\bibnamefont {Rendall}},\
  and\ \bibinfo {editor} {\bibfnamefont {M.}~\bibnamefont {Scherfner}}}\
  (\bibinfo  {publisher} {American Mathematical Society and International
  Press},\ \bibinfo {address} {USA},\ \bibinfo {year} {2011})\ pp.\ \bibinfo
  {pages} {105--124},\ \Eprint {https://arxiv.org/abs/0712.1933}
  {arXiv:0712.1933 [gr-qc]}\BibitemShut {NoStop}%
\bibitem [{\citenamefont {Hadamard}(1908)}]{Hadamard_ActaMath_1908}%
  \BibitemOpen
  \bibfield  {author} {\bibinfo {author} {\bibfnamefont {J.}~\bibnamefont
  {Hadamard}},\ }\bibfield  {title} {\bibinfo {title} {Th\'{e}orie des
  \'{e}quations aux d\'{e}riv\'{e}es partielles lin\'{e}aires hyperboliques et
  du probl\`{e}me de cauchy},\ }\href {https://doi.org/10.1007/BF02415449}
  {\bibfield  {journal} {\bibinfo  {journal} {Acta Math.}\ }\textbf {\bibinfo
  {volume} {31}},\ \bibinfo {pages} {333} (\bibinfo {year} {1908})}\BibitemShut
  {NoStop}%
\bibitem [{\citenamefont {Hadamard}(2003)}]{Hadamard_Dover_2003}%
  \BibitemOpen
  \bibfield  {author} {\bibinfo {author} {\bibfnamefont {J.}~\bibnamefont
  {Hadamard}},\ }\href@noop {} {\emph {\bibinfo {title} {Lectures on {C}auchy's
  problem in linear partial differential equations}}}\ (\bibinfo  {publisher}
  {Dover},\ \bibinfo {address} {NY},\ \bibinfo {year} {2003})\BibitemShut
  {NoStop}%
\bibitem [{\citenamefont {Riesz}(1949)}]{Riesz_ActaMath_1949}%
  \BibitemOpen
  \bibfield  {author} {\bibinfo {author} {\bibfnamefont {M.}~\bibnamefont
  {Riesz}},\ }\bibfield  {title} {\bibinfo {title} {L'int\'{e}grale de
  riemann-liouville et le probl\`{e}me de {C}auchy},\ }\href
  {https://doi.org/10.1007/BF02395016} {\bibfield  {journal} {\bibinfo
  {journal} {Acta Math.}\ }\textbf {\bibinfo {volume} {81}},\ \bibinfo {pages}
  {1} (\bibinfo {year} {1949})}\BibitemShut {NoStop}%
\bibitem [{\citenamefont {Riesz}(1960)}]{Riesz_CPAM_1960}%
  \BibitemOpen
  \bibfield  {author} {\bibinfo {author} {\bibfnamefont {M.}~\bibnamefont
  {Riesz}},\ }\bibfield  {title} {\bibinfo {title} {A geometric solution of the
  wave equation in space-time of even dimension},\ }\href
  {https://doi.org/10.1002/cpa.3160130302} {\bibfield  {journal} {\bibinfo
  {journal} {Comm. Pure Appl. Math.}\ }\textbf {\bibinfo {volume} {13}},\
  \bibinfo {pages} {329} (\bibinfo {year} {1960})}\BibitemShut {NoStop}%
\bibitem [{\citenamefont {Friedlander}(1975)}]{Friedlander_CUP_1975}%
  \BibitemOpen
  \bibfield  {author} {\bibinfo {author} {\bibfnamefont {F.}~\bibnamefont
  {Friedlander}},\ }\href@noop {} {\emph {\bibinfo {title} {The Wave Equation
  on a Curved Space-Time}}},\ Cambridge Monographs on Mathematical Physics\
  (\bibinfo  {publisher} {Cambridge University Press},\ \bibinfo {year}
  {1975})\BibitemShut {NoStop}%
\bibitem [{\citenamefont {Feynman}(1949)}]{Feynman_PR_1949}%
  \BibitemOpen
  \bibfield  {author} {\bibinfo {author} {\bibfnamefont {R.~P.}\ \bibnamefont
  {Feynman}},\ }\bibfield  {title} {\bibinfo {title} {The theory of
  positrons},\ }\href {https://doi.org/10.1103/PhysRev.76.749} {\bibfield
  {journal} {\bibinfo  {journal} {Phys. Rev.}\ }\textbf {\bibinfo {volume}
  {76}},\ \bibinfo {pages} {749} (\bibinfo {year} {1949})}\BibitemShut
  {NoStop}%
\bibitem [{\citenamefont {Stueckelberg}(1941)}]{Stueckelberg_HPA_1941}%
  \BibitemOpen
  \bibfield  {author} {\bibinfo {author} {\bibfnamefont {E.}~\bibnamefont
  {Stueckelberg}},\ }\bibfield  {title} {\bibinfo {title} {Remarque a propos de
  la cr\'{e}ation de paires de particules en th\'{e}orie de la
  relativit\'{e}},\ }\href {http://dx.doi.org/10.5169/seals-111201} {\bibfield
  {journal} {\bibinfo  {journal} {Helv. phys. acta}\ }\textbf {\bibinfo
  {volume} {14}},\ \bibinfo {pages} {588} (\bibinfo {year} {1941})}\BibitemShut
  {NoStop}%
\bibitem [{\citenamefont {Dirac}(1928{\natexlab{a}})}]{Dirac_PRSA_1928}%
  \BibitemOpen
  \bibfield  {author} {\bibinfo {author} {\bibfnamefont {P.~A.~M.}\
  \bibnamefont {Dirac}},\ }\bibfield  {title} {\bibinfo {title} {The quantum
  theory of the electron},\ }\href
  {http://rspa.royalsocietypublishing.org/content/117/778/610} {\bibfield
  {journal} {\bibinfo  {journal} {Proc. R. Soc. Lond. A}\ }\textbf {\bibinfo
  {volume} {117}},\ \bibinfo {pages} {610 } (\bibinfo {year}
  {1928}{\natexlab{a}})}\BibitemShut {NoStop}%
\bibitem [{\citenamefont {Dirac}(1928{\natexlab{b}})}]{Dirac_PRSA_1928_P2}%
  \BibitemOpen
  \bibfield  {author} {\bibinfo {author} {\bibfnamefont {P.~A.~M.}\
  \bibnamefont {Dirac}},\ }\bibfield  {title} {\bibinfo {title} {The quantum
  theory of the electron. part {I}{I}.},\ }\href
  {http://rspa.royalsocietypublishing.org/content/118/779/351} {\bibfield
  {journal} {\bibinfo  {journal} {Proc. R. Soc. Lond. A}\ }\textbf {\bibinfo
  {volume} {118}},\ \bibinfo {pages} {351 } (\bibinfo {year}
  {1928}{\natexlab{b}})}\BibitemShut {NoStop}%
\bibitem [{\citenamefont {Atiyah}\ \emph {et~al.}(1975)\citenamefont {Atiyah},
  \citenamefont {Patodi},\ and\ \citenamefont
  {Singer}}]{Atiyah_MathProcCam_1975}%
  \BibitemOpen
  \bibfield  {author} {\bibinfo {author} {\bibfnamefont {M.~F.}\ \bibnamefont
  {Atiyah}}, \bibinfo {author} {\bibfnamefont {V.~K.}\ \bibnamefont {Patodi}},\
  and\ \bibinfo {author} {\bibfnamefont {I.~M.}\ \bibnamefont {Singer}},\
  }\bibfield  {title} {\bibinfo {title} {Spectral asymmetry and {R}iemannian
  geometry. {I}},\ }\href {https://doi.org/10.1017/S0305004100049410}
  {\bibfield  {journal} {\bibinfo  {journal} {Math. Proc. Cambridge Philos.}\
  }\textbf {\bibinfo {volume} {77}},\ \bibinfo {pages} {43} (\bibinfo {year}
  {1975})}\BibitemShut {NoStop}%
\bibitem [{\citenamefont {Branson}\ and\ \citenamefont
  {Gilkey}(1992)}]{Branson_JFA_1992}%
  \BibitemOpen
  \bibfield  {author} {\bibinfo {author} {\bibfnamefont {T.~P.}\ \bibnamefont
  {Branson}}\ and\ \bibinfo {author} {\bibfnamefont {P.~B.}\ \bibnamefont
  {Gilkey}},\ }\bibfield  {title} {\bibinfo {title} {Residues of the eta
  function for an operator of {D}irac type},\ }\href
  {https://doi.org/https://doi.org/10.1016/0022-1236(92)90146-A} {\bibfield
  {journal} {\bibinfo  {journal} {J. Funct. Anal.}\ }\textbf {\bibinfo {volume}
  {108}},\ \bibinfo {pages} {47 } (\bibinfo {year} {1992})}\BibitemShut
  {NoStop}%
\bibitem [{\citenamefont {Geroch}(1971)}]{Geroch_JMP_1971}%
  \BibitemOpen
  \bibfield  {author} {\bibinfo {author} {\bibfnamefont {R.}~\bibnamefont
  {Geroch}},\ }\bibfield  {title} {\bibinfo {title} {A method for generating
  solutions of {E}instein's equations},\ }\href
  {https://doi.org/10.1063/1.1665681} {\bibfield  {journal} {\bibinfo
  {journal} {J. Math. Phys.}\ }\textbf {\bibinfo {volume} {12}},\ \bibinfo
  {pages} {918 } (\bibinfo {year} {1971})}\BibitemShut {NoStop}%
\bibitem [{\citenamefont {Harris}(1992)}]{Harris_CQG_1992}%
  \BibitemOpen
  \bibfield  {author} {\bibinfo {author} {\bibfnamefont {S.}~\bibnamefont
  {Harris}},\ }\bibfield  {title} {\bibinfo {title} {Conformally stationary
  spacetimes},\ }\href {https://doi.org/10.1088/0264-9381/9/7/013} {\bibfield
  {journal} {\bibinfo  {journal} {Class. Quantum Grav.}\ }\textbf {\bibinfo
  {volume} {9}},\ \bibinfo {pages} {1823 } (\bibinfo {year}
  {1992})}\BibitemShut {NoStop}%
\bibitem [{\citenamefont {Javaloyes}\ and\ \citenamefont
  {S{\'{a}}nchez}(2008)}]{Javaloyes_CQG_2008}%
  \BibitemOpen
  \bibfield  {author} {\bibinfo {author} {\bibfnamefont {M.~A.}\ \bibnamefont
  {Javaloyes}}\ and\ \bibinfo {author} {\bibfnamefont {M.}~\bibnamefont
  {S{\'{a}}nchez}},\ }\bibfield  {title} {\bibinfo {title} {A note on the
  existence of standard splittings for conformally stationary spacetimes},\
  }\href {https://doi.org/10.1088/0264-9381/25/16/168001} {\bibfield  {journal}
  {\bibinfo  {journal} {Class. Quantum Grav.}\ }\textbf {\bibinfo {volume}
  {25}},\ \bibinfo {pages} {168001} (\bibinfo {year} {2008})},\ \Eprint
  {https://arxiv.org/abs/0806.0812} {arXiv:0806.0812 [gr-qc]}\BibitemShut
  {NoStop}%
\bibitem [{\citenamefont {Sanders}(2013)}]{Sanders_IJMPA_2013}%
  \BibitemOpen
  \bibfield  {author} {\bibinfo {author} {\bibfnamefont {K.}~\bibnamefont
  {Sanders}},\ }\bibfield  {title} {\bibinfo {title} {Thermal equilibrium
  states of a linear scalar quantum field in stationary space-times},\ }\href
  {https://doi.org/10.1142/S0217751X1330010X} {\bibfield  {journal} {\bibinfo
  {journal} {Int. J. Mod. Phys. A}\ }\textbf {\bibinfo {volume} {28}},\
  \bibinfo {pages} {1330010} (\bibinfo {year} {2013})},\ \Eprint
  {https://arxiv.org/abs/1209.6068} {arXiv:1209.6068 [math-ph]}\BibitemShut
  {NoStop}%
\bibitem [{\citenamefont {Caponio}\ \emph {et~al.}(2011)\citenamefont
  {Caponio}, \citenamefont {Javaloyes},\ and\ \citenamefont
  {Sánchez}}]{Caponio_RMI_2011}%
  \BibitemOpen
  \bibfield  {author} {\bibinfo {author} {\bibfnamefont {E.}~\bibnamefont
  {Caponio}}, \bibinfo {author} {\bibfnamefont {M.~A.}\ \bibnamefont
  {Javaloyes}},\ and\ \bibinfo {author} {\bibfnamefont {M.}~\bibnamefont
  {Sánchez}},\ }\bibfield  {title} {\bibinfo {title} {On the interplay between
  {L}orentzian causality and {F}insler metrics of {R}anders type},\ }\href
  {https://doi.org/10.4171/rmi/658} {\bibfield  {journal} {\bibinfo  {journal}
  {Rev. Matem. Iberoam.}\ ,\ \bibinfo {pages} {919}} (\bibinfo {year}
  {2011})},\ \Eprint {https://arxiv.org/abs/0903.3501} {arXiv:0903.3501
  [math.DG]}\BibitemShut {NoStop}%
\bibitem [{\citenamefont {Candela}\ \emph {et~al.}(2008)\citenamefont
  {Candela}, \citenamefont {Flores},\ and\ \citenamefont
  {S\'{a}nchez}}]{Candela_AdvMath_2008}%
  \BibitemOpen
  \bibfield  {author} {\bibinfo {author} {\bibfnamefont {A.}~\bibnamefont
  {Candela}}, \bibinfo {author} {\bibfnamefont {J.}~\bibnamefont {Flores}},\
  and\ \bibinfo {author} {\bibfnamefont {M.}~\bibnamefont {S\'{a}nchez}},\
  }\bibfield  {title} {\bibinfo {title} {Global hyperbolicity and
  {P}alais–{S}male condition for action functionals in stationary
  spacetimes},\ }\href
  {https://doi.org/https://doi.org/10.1016/j.aim.2008.01.004} {\bibfield
  {journal} {\bibinfo  {journal} {Adv. Math.}\ }\textbf {\bibinfo {volume}
  {218}},\ \bibinfo {pages} {515} (\bibinfo {year} {2008})},\ \Eprint
  {https://arxiv.org/abs/0610175} {arXiv:0610175 [math.DG]}\BibitemShut
  {NoStop}%
\bibitem [{\citenamefont {Anderson}(2000)}]{Anderson_AHP_2000}%
  \BibitemOpen
  \bibfield  {author} {\bibinfo {author} {\bibfnamefont {M.}~\bibnamefont
  {Anderson}},\ }\bibfield  {title} {\bibinfo {title} {On stationary vacuum
  solutions to the {E}instein equations},\ }\href
  {https://doi.org/10.1007/pl00001021} {\bibfield  {journal} {\bibinfo
  {journal} {Ann. Henri Poincar\'{e}}\ }\textbf {\bibinfo {volume} {1}},\
  \bibinfo {pages} {977 } (\bibinfo {year} {2000})},\ \Eprint
  {https://arxiv.org/abs/0001091} {arXiv:0001091 [gr-qc]}\BibitemShut
  {NoStop}%
\bibitem [{\citenamefont {B\"{a}r}(2000)}]{Baer_Seminar_2000}%
  \BibitemOpen
  \bibfield  {author} {\bibinfo {author} {\bibfnamefont {C.}~\bibnamefont
  {B\"{a}r}},\ }\bibfield  {title} {\bibinfo {title} {Dependence of the
  spectrum of the {D}irac operator on the spinor structure},\ }\href
  {https://smf.emath.fr/sites/default/files/2018-11/smf_sem-cong_4_17-33__sample.pdf}
  {\bibfield  {journal} {\bibinfo  {journal} {S\'{e}minaires \& Congr\`{e}s}\
  }\textbf {\bibinfo {volume} {4}},\ \bibinfo {pages} {17} (\bibinfo {year}
  {2000})}\BibitemShut {NoStop}%
\bibitem [{\citenamefont {Meinrenken}(1994)}]{Meinrenken_JGP_1994}%
  \BibitemOpen
  \bibfield  {author} {\bibinfo {author} {\bibfnamefont {E.}~\bibnamefont
  {Meinrenken}},\ }\bibfield  {title} {\bibinfo {title} {Trace formulas and the
  {C}onley-{Z}ehnder index},\ }\href
  {https://doi.org/https://doi.org/10.1016/0393-0440(94)90058-2} {\bibfield
  {journal} {\bibinfo  {journal} {J. Geom. Phys.}\ }\textbf {\bibinfo {volume}
  {13}},\ \bibinfo {pages} {1 } (\bibinfo {year} {1994})}\BibitemShut {NoStop}%
\bibitem [{\citenamefont {Capoferri}\ and\ \citenamefont
  {Murro}(2022)}]{Capoferri}%
  \BibitemOpen
  \bibfield  {author} {\bibinfo {author} {\bibfnamefont {M.}~\bibnamefont
  {Capoferri}}\ and\ \bibinfo {author} {\bibfnamefont {S.}~\bibnamefont
  {Murro}},\ }\href {https://doi.org/10.48550/ARXIV.2201.12104} {\bibinfo
  {title} {Global and microlocal aspects of {D}irac operators: propagators and
  {H}adamard states}} (\bibinfo {year} {2022}),\ \Eprint
  {https://arxiv.org/abs/2201.12104} {arXiv:2201.12104 [math.AP]}\BibitemShut
  {NoStop}%
\bibitem [{\citenamefont {Capoferri}\ and\ \citenamefont
  {Vassiliev}(2022)}]{Capoferri_2020}%
  \BibitemOpen
  \bibfield  {author} {\bibinfo {author} {\bibfnamefont {M.}~\bibnamefont
  {Capoferri}}\ and\ \bibinfo {author} {\bibfnamefont {D.}~\bibnamefont
  {Vassiliev}},\ }\bibfield  {title} {\bibinfo {title} {Global propagator for
  the massless {D}irac operator and spectral asymptotics},\ }\href
  {https://doi.org/https://doi.org/10.1007/s00020-022-02708-1} {\bibfield
  {journal} {\bibinfo  {journal} {Integr. Equ. Oper. Theory}\ }\textbf
  {\bibinfo {volume} {94}},\ \bibinfo {pages} {30} (\bibinfo {year} {2022})},\
  \Eprint {https://arxiv.org/abs/2004.06351} {arXiv:2004.06351 [math.AP]}\BibitemShut 
  {NoStop}%
\bibitem [{\citenamefont {Duistermaat}\ and\ \citenamefont
  {Kolk}(2010)}]{Duistermaat_Birkhaeuser_2010}%
  \BibitemOpen
  \bibfield  {author} {\bibinfo {author} {\bibfnamefont {J.~J.}\ \bibnamefont
  {Duistermaat}}\ and\ \bibinfo {author} {\bibfnamefont {J.~A.}\ \bibnamefont
  {Kolk}},\ }\href {https://doi.org/https://doi.org/10.1007/978-0-8176-4675-2}
  {\emph {\bibinfo {title} {Distributions: Theory and Applications}}},\
  Cornerstones\ (\bibinfo  {publisher} {Birkh\"{a}user},\ \bibinfo {address}
  {Boston},\ \bibinfo {year} {2010})\BibitemShut {NoStop}%
\bibitem [{\citenamefont {Safarov}(2001)}]{Safarov_JFA_2001}%
  \BibitemOpen
  \bibfield  {author} {\bibinfo {author} {\bibfnamefont {Y.}~\bibnamefont
  {Safarov}},\ }\bibfield  {title} {\bibinfo {title} {Fourier tauberian
  theorems and applications},\ }\href
  {https://doi.org/https://doi.org/10.1006/jfan.2001.3764} {\bibfield
  {journal} {\bibinfo  {journal} {J. Funct. Anal.}\ }\textbf {\bibinfo {volume}
  {185}},\ \bibinfo {pages} {111 } (\bibinfo {year} {2001})},\ \Eprint
  {https://arxiv.org/abs/0003014} {arXiv:0003014 [math.SP]}\BibitemShut
  {NoStop}%
\bibitem [{\citenamefont {Masiello}(1992)}]{Masiello_NonlinearAnal_1992}%
  \BibitemOpen
  \bibfield  {author} {\bibinfo {author} {\bibfnamefont {A.}~\bibnamefont
  {Masiello}},\ }\bibfield  {title} {\bibinfo {title} {Time-like periodic
  trajectories in stationary {L}orentz manifolds},\ }\href
  {https://doi.org/https://doi.org/10.1016/0362-546X(92)90019-B} {\bibfield
  {journal} {\bibinfo  {journal} {Nonlinear Anal. Theory Methods. Appl.}\
  }\textbf {\bibinfo {volume} {19}},\ \bibinfo {pages} {531} (\bibinfo {year}
  {1992})}\BibitemShut {NoStop}%
\bibitem [{\citenamefont
  {S\'{a}nchez}(1999{\natexlab{a}})}]{Sanchez_NonlinearAnal_1999}%
  \BibitemOpen
  \bibfield  {author} {\bibinfo {author} {\bibfnamefont {M.}~\bibnamefont
  {S\'{a}nchez}},\ }\bibfield  {title} {\bibinfo {title} {Geodesics in static
  spacetimes and t-periodic trajectories},\ }\href
  {https://doi.org/https://doi.org/10.1016/S0362-546X(97)00683-4} {\bibfield
  {journal} {\bibinfo  {journal} {Nonlinear Anal. Theory, Methods Appl.}\
  }\textbf {\bibinfo {volume} {35}},\ \bibinfo {pages} {677} (\bibinfo {year}
  {1999}{\natexlab{a}})}\BibitemShut {NoStop}%
\bibitem [{\citenamefont {S\'{a}nchez}(1999{\natexlab{b}})}]{Sanchez_AMS_1999}%
  \BibitemOpen
  \bibfield  {author} {\bibinfo {author} {\bibfnamefont {M.}~\bibnamefont
  {S\'{a}nchez}},\ }\bibfield  {title} {\bibinfo {title} {Timelike periodic
  trajectories in spatially compact {L}orentz manifolds},\ }\href
  {http://www.jstor.org/stable/119967} {\bibfield  {journal} {\bibinfo
  {journal} {Proc. Am. Math. Soc.}\ }\textbf {\bibinfo {volume} {127}},\
  \bibinfo {pages} {3057} (\bibinfo {year} {1999}{\natexlab{b}})}\BibitemShut
  {NoStop}%
\bibitem [{\citenamefont {Bartolo}(2001)}]{Bartolo_NonlinearAnal_2001}%
  \BibitemOpen
  \bibfield  {author} {\bibinfo {author} {\bibfnamefont {R.}~\bibnamefont
  {Bartolo}},\ }\bibfield  {title} {\bibinfo {title} {Periodic trajectories on
  stationary {L}orentzian manifolds},\ }\href
  {https://doi.org/https://doi.org/10.1016/S0362-546X(99)00246-1} {\bibfield
  {journal} {\bibinfo  {journal} {Nonlinear Anal. Theory Methods. Appl.}\
  }\textbf {\bibinfo {volume} {43}},\ \bibinfo {pages} {883} (\bibinfo {year}
  {2001})}\BibitemShut {NoStop}%
\bibitem [{\citenamefont {Reed}\ and\ \citenamefont {Simon}(1980)}]{Reed_I}%
  \BibitemOpen
  \bibfield  {author} {\bibinfo {author} {\bibfnamefont {M.}~\bibnamefont
  {Reed}}\ and\ \bibinfo {author} {\bibfnamefont {B.}~\bibnamefont {Simon}},\
  }\href@noop {} {\emph {\bibinfo {title} {Functional Analysis}}},\ \bibinfo
  {series} {Methods of Modern Mathematical Physics}, Vol.~\bibinfo {volume}
  {I}\ (\bibinfo  {publisher} {Acadmic Press},\ \bibinfo {address} {USA},\
  \bibinfo {year} {1980})\BibitemShut {NoStop}%
\bibitem [{\citenamefont {Guillemin}(1985)}]{Guillemin_AdvMath_1985}%
  \BibitemOpen
  \bibfield  {author} {\bibinfo {author} {\bibfnamefont {V.}~\bibnamefont
  {Guillemin}},\ }\bibfield  {title} {\bibinfo {title} {A new proof of {W}eyl's
  formula on the asymptotic distribution of eigenvalues},\ }\href
  {https://doi.org/https://doi.org/10.1016/0001-8708(85)90018-0} {\bibfield
  {journal} {\bibinfo  {journal} {Adv. Math.}\ }\textbf {\bibinfo {volume}
  {55}},\ \bibinfo {pages} {131 } (\bibinfo {year} {1985})}\BibitemShut
  {NoStop}%
\bibitem [{\citenamefont {Wunsch}(2013)}]{Wunsch_Zuerich_2008}%
  \BibitemOpen
  \bibfield  {author} {\bibinfo {author} {\bibfnamefont {J.}~\bibnamefont
  {Wunsch}},\ }\bibfield  {title} {\bibinfo {title} {Microlocal analysis and
  evolution equations: Lecture notes from the 2008 {C}{M}{I}/{E}{T}{H} summer
  school},\ }in\ \href {https://bookstore.ams.org/cmip-17} {\emph {\bibinfo
  {booktitle} {Evolution Equations}}},\ \bibinfo {series} {Clay Mathematics
  Proceedings}, Vol.~\bibinfo {volume} {17},\ \bibinfo {editor} {edited by\
  \bibinfo {editor} {\bibfnamefont {D.}~\bibnamefont {Ellwood}}, \bibinfo
  {editor} {\bibfnamefont {I.}~\bibnamefont {Rodnianski}}, \bibinfo {editor}
  {\bibfnamefont {G.}~\bibnamefont {Staffilani}},\ and\ \bibinfo {editor}
  {\bibfnamefont {J.}~\bibnamefont {Wunsch}}}\ (\bibinfo  {publisher} {American
  Mathematical Society},\ \bibinfo {year} {2013})\ pp.\ \bibinfo {pages}
  {1--72},\ \Eprint {https://arxiv.org/abs/0812.3181} {arXiv:0812.3181
  [math.AP]}\BibitemShut {NoStop}%
\bibitem [{\citenamefont {Colin~de Verdi\`{e}re}(2007)}]{Verdiere_AIF_2007}%
  \BibitemOpen
  \bibfield  {author} {\bibinfo {author} {\bibfnamefont {Y.}~\bibnamefont
  {Colin~de Verdi\`{e}re}},\ }\bibfield  {title} {\bibinfo {title} {Spectrum of
  the {L}aplace operator and periodic geodesics: thirty years after},\ }\href
  {https://doi.org/10.5802/aif.2339} {\bibfield  {journal} {\bibinfo  {journal}
  {Ann. Inst. Fourier}\ }\textbf {\bibinfo {volume} {57}},\ \bibinfo {pages}
  {2429} (\bibinfo {year} {2007})}\BibitemShut {NoStop}%
\bibitem [{\citenamefont {Bolte}\ and\ \citenamefont
  {Keppeler}(1998)}]{Bolte_PRL_1998}%
  \BibitemOpen
  \bibfield  {author} {\bibinfo {author} {\bibfnamefont {J.}~\bibnamefont
  {Bolte}}\ and\ \bibinfo {author} {\bibfnamefont {S.}~\bibnamefont
  {Keppeler}},\ }\bibfield  {title} {\bibinfo {title} {Semiclassical time
  evolution and trace formula for relativistic spin-1/2 particles},\ }\href
  {https://doi.org/10.1103/PhysRevLett.81.1987} {\bibfield  {journal} {\bibinfo
   {journal} {Phys. Rev. Lett.}\ }\textbf {\bibinfo {volume} {81}},\ \bibinfo
  {pages} {1987 } (\bibinfo {year} {1998})},\ \Eprint
  {https://arxiv.org/abs/9805041} {arXiv:9805041 [quant-ph]}\BibitemShut
  {NoStop}%
\bibitem [{\citenamefont {Bolte}\ and\ \citenamefont
  {Keppeler}(1999)}]{Bolte_AnnPhys_1999}%
  \BibitemOpen
  \bibfield  {author} {\bibinfo {author} {\bibfnamefont {J.}~\bibnamefont
  {Bolte}}\ and\ \bibinfo {author} {\bibfnamefont {S.}~\bibnamefont
  {Keppeler}},\ }\bibfield  {title} {\bibinfo {title} {A semiclassical approach
  to the {D}irac equation},\ }\href
  {https://doi.org/https://doi.org/10.1006/aphy.1999.5912} {\bibfield
  {journal} {\bibinfo  {journal} {Ann. Phys.}\ }\textbf {\bibinfo {volume}
  {274}},\ \bibinfo {pages} {125 } (\bibinfo {year} {1999})},\ \Eprint
  {https://arxiv.org/abs/9811025} {arXiv:9811025 [quant-ph]}\BibitemShut
  {NoStop}%
\bibitem [{\citenamefont {Bolte}(2001)}]{Bolte_FoundPhys_2001}%
  \BibitemOpen
  \bibfield  {author} {\bibinfo {author} {\bibfnamefont {J.}~\bibnamefont
  {Bolte}},\ }\bibfield  {title} {\bibinfo {title} {Semiclassical expectation
  values for relativistic particles with spin 1/2},\ }\href
  {https://doi.org/10.1023/A:1017502906292} {\bibfield  {journal} {\bibinfo
  {journal} {Found. Phys.}\ }\textbf {\bibinfo {volume} {31}},\ \bibinfo
  {pages} {423 } (\bibinfo {year} {2001})},\ \Eprint
  {https://arxiv.org/abs/0009052} {arXiv:0009052 [nlin.CD]}\BibitemShut
  {NoStop}%
\bibitem [{\citenamefont {Bolte}\ and\ \citenamefont
  {Glaser}(2004{\natexlab{b}})}]{Bolte_JPA_2004}%
  \BibitemOpen
  \bibfield  {author} {\bibinfo {author} {\bibfnamefont {J.}~\bibnamefont
  {Bolte}}\ and\ \bibinfo {author} {\bibfnamefont {R.}~\bibnamefont {Glaser}},\
  }\bibfield  {title} {\bibinfo {title} {Zitterbewegung and semiclassical
  observables for the {D}irac equation},\ }\href
  {https://doi.org/10.1088/0305-4470/37/24/012} {\bibfield  {journal} {\bibinfo
   {journal} {J. Phys. A: Math. Gen.}\ }\textbf {\bibinfo {volume} {37}},\
  \bibinfo {pages} {6359 } (\bibinfo {year} {2004}{\natexlab{b}})},\ \Eprint
  {https://arxiv.org/abs/0402154} {arXiv:0402154 [quant-ph]}\BibitemShut
  {NoStop}%
\bibitem [{\citenamefont {Chervova}\ \emph {et~al.}(2013)\citenamefont
  {Chervova}, \citenamefont {Downes},\ and\ \citenamefont
  {Vassiliev}}]{Chervova_JST_2013}%
  \BibitemOpen
  \bibfield  {author} {\bibinfo {author} {\bibfnamefont {O.}~\bibnamefont
  {Chervova}}, \bibinfo {author} {\bibfnamefont {R.}~\bibnamefont {Downes}},\
  and\ \bibinfo {author} {\bibfnamefont {D.}~\bibnamefont {Vassiliev}},\
  }\bibfield  {title} {\bibinfo {title} {The spectral function of a first order
  elliptic system},\ }\href {https://doi.org/10.4171/JST/47} {\bibfield
  {journal} {\bibinfo  {journal} {J. Spectr. Theory}\ }\textbf {\bibinfo
  {volume} {3}},\ \bibinfo {pages} {317 } (\bibinfo {year} {2013})},\ \Eprint
  {https://arxiv.org/abs/1208.6015} {arXiv:1208.6015 [math.SP]}\BibitemShut
  {NoStop}%
\bibitem [{\citenamefont {Rudolph}\ and\ \citenamefont
  {Schmidt}(2013)}]{Rudolph_Springer_2013}%
  \BibitemOpen
  \bibfield  {author} {\bibinfo {author} {\bibfnamefont {G.}~\bibnamefont
  {Rudolph}}\ and\ \bibinfo {author} {\bibfnamefont {M.}~\bibnamefont
  {Schmidt}},\ }\href {https://doi.org/10.1007/978-94-007-5345-7} {\emph
  {\bibinfo {title} {Differential Geometry and Mathematical Physics - {I}:
  Manifolds, Lie Groups and Hamiltonian Systems}}},\ Theoretical and
  Mathematical Physics\ (\bibinfo  {publisher} {Springer},\ \bibinfo {address}
  {Netherlands},\ \bibinfo {year} {2013})\BibitemShut {NoStop}%
\bibitem [{\citenamefont {Weinstein}(1977)}]{Weinstein_AMS_1977}%
  \BibitemOpen
  \bibfield  {author} {\bibinfo {author} {\bibfnamefont {A.}~\bibnamefont
  {Weinstein}},\ }\href {https://bookstore.ams.org/cbms-29} {\emph {\bibinfo
  {title} {Lectures on Symplectic Manifolds}}},\ \bibinfo {series} {CBMS
  Regional Conference Series in Mathematics}, Vol.~\bibinfo {volume} {29}\
  (\bibinfo  {publisher} {American Mathemtical Society},\ \bibinfo {address}
  {USA},\ \bibinfo {year} {1977})\BibitemShut {NoStop}%
\bibitem [{\citenamefont {Latour}(1991)}]{Latour_ASENS_1991}%
  \BibitemOpen
  \bibfield  {author} {\bibinfo {author} {\bibfnamefont {F.}~\bibnamefont
  {Latour}},\ }\bibfield  {title} {\bibinfo {title} {Lagrangian transversals,
  periodicity of {Bott} and generating forms for a {L}agrangian immersion in a
  cotangent},\ }\href {https://doi.org/10.24033/asens.1619} {\bibfield
  {journal} {\bibinfo  {journal} {Ann. Sci. \'{E}c. Norm. Sup\'{e}r.}\
  }\bibinfo {series} {4},\ \textbf {\bibinfo {volume} {24}},\ \bibinfo {pages}
  {3} (\bibinfo {year} {1991})}\BibitemShut {NoStop}%
\bibitem [{\citenamefont {Giroux}(1990)}]{Giroux_Springer_1990}%
  \BibitemOpen
  \bibfield  {author} {\bibinfo {author} {\bibfnamefont {E.}~\bibnamefont
  {Giroux}},\ }\bibinfo {title} {Formes generatrices d'immersions lagrangiennes
  dans un espace cotangent},\ in\ \href {https://doi.org/10.1007/BFb0097468}
  {\emph {\bibinfo {booktitle} {G{\'e}om{\'e}trie Symplectique et
  M{\'e}canique: Colloque International La Grande Motte, France, 23--28 Mai,
  1988}}},\ \bibinfo {editor} {edited by\ \bibinfo {editor} {\bibfnamefont
  {C.}~\bibnamefont {Albert}}}\ (\bibinfo  {publisher} {Springer},\ \bibinfo
  {address} {Berlin, Heidelberg},\ \bibinfo {year} {1990})\ pp.\ \bibinfo
  {pages} {139--145}\BibitemShut {NoStop}%
\bibitem [{\citenamefont {Arnol'd}(1967)}]{Arnold_FAA_1967}%
  \BibitemOpen
  \bibfield  {author} {\bibinfo {author} {\bibfnamefont {V.~I.}\ \bibnamefont
  {Arnol'd}},\ }\bibfield  {title} {\bibinfo {title} {Characteristic class
  entering in quantization conditions},\ }\href
  {https://doi.org/https://doi.org/10.1007/BF01075861} {\bibfield  {journal}
  {\bibinfo  {journal} {Funct. Anal. Its Appl.}\ }\textbf {\bibinfo {volume}
  {1}},\ \bibinfo {pages} {1} (\bibinfo {year} {1967})}\BibitemShut {NoStop}%
\bibitem [{\citenamefont {Ann{\'e}}(2006)}]{Anne_Cubo_2006}%
  \BibitemOpen
  \bibfield  {author} {\bibinfo {author} {\bibfnamefont {C.}~\bibnamefont
  {Ann{\'e}}},\ }\bibfield  {title} {\bibinfo {title} {A topological definition
  of the {M}aslov bundle},\ }\href
  {https://hal.archives-ouvertes.fr/hal-00017676} {\bibfield  {journal}
  {\bibinfo  {journal} {Cubo, a Mathematical Journal}\ }\textbf {\bibinfo
  {volume} {8}},\ \bibinfo {pages} {N. 1, p. 1} (\bibinfo {year} {2006})},\
  \Eprint {https://arxiv.org/abs/0601586} {arXiv:0601586 [math.DG]}\BibitemShut 
  {NoStop}%
\bibitem [{\citenamefont {Keller}(1958)}]{Keller_AnnPhys_1958}%
  \BibitemOpen
  \bibfield  {author} {\bibinfo {author} {\bibfnamefont {J.~B.}\ \bibnamefont
  {Keller}},\ }\bibfield  {title} {\bibinfo {title} {Corrected
  {B}ohr-{S}ommerfeld quantum conditions for nonseparable systems},\ }\href
  {https://doi.org/https://doi.org/10.1016/0003-4916(58)90032-0} {\bibfield
  {journal} {\bibinfo  {journal} {Ann. Phys.}\ }\textbf {\bibinfo {volume}
  {4}},\ \bibinfo {pages} {180 } (\bibinfo {year} {1958})}\BibitemShut
  {NoStop}%
\bibitem [{\citenamefont {Weinstein}(1971)}]{Weinstein_AdvMath_1971}%
  \BibitemOpen
  \bibfield  {author} {\bibinfo {author} {\bibfnamefont {A.}~\bibnamefont
  {Weinstein}},\ }\bibfield  {title} {\bibinfo {title} {Symplectic manifolds
  and their {L}agrangian submanifolds},\ }\href
  {https://doi.org/https://doi.org/10.1016/0001-8708(71)90020-X} {\bibfield
  {journal} {\bibinfo  {journal} {Adv. Math.}\ }\textbf {\bibinfo {volume}
  {6}},\ \bibinfo {pages} {329 } (\bibinfo {year} {1971})}\BibitemShut
  {NoStop}%
\bibitem [{\citenamefont {H\"{o}rmander}(1967)}]{Hoermander_AMS_1967}%
  \BibitemOpen
  \bibfield  {author} {\bibinfo {author} {\bibfnamefont {L.}~\bibnamefont
  {H\"{o}rmander}},\ }\bibfield  {title} {\bibinfo {title} {Pseudo-differential
  operators and hypoelliptic equations},\ }in\ \href
  {https://doi.org/https://doi.org/10.1090/pspum/010} {\emph {\bibinfo
  {booktitle} {Singular Integrals: Proc. Symposia Pure Math.}}},\ Vol.~\bibinfo
  {volume} {10},\ \bibinfo {editor} {edited by\ \bibinfo {editor}
  {\bibfnamefont {A.~P.}\ \bibnamefont {Calder\'{o}n}}}\ (\bibinfo  {publisher}
  {American Mathematical Society},\ \bibinfo {address} {USA},\ \bibinfo {year}
  {1967})\ pp.\ \bibinfo {pages} {138 -- 183}\BibitemShut {NoStop}%
\bibitem [{\citenamefont {Hollands}\ and\ \citenamefont
  {Wald}(2015)}]{Hollands_PR_2015}%
  \BibitemOpen
  \bibfield  {author} {\bibinfo {author} {\bibfnamefont {S.}~\bibnamefont
  {Hollands}}\ and\ \bibinfo {author} {\bibfnamefont {R.~M.}\ \bibnamefont
  {Wald}},\ }\bibfield  {title} {\bibinfo {title} {Quantum fields in curved
  spacetime},\ }\href
  {https://doi.org/http://dx.doi.org/10.1016/j.physrep.2015.02.001} {\bibfield
  {journal} {\bibinfo  {journal} {Phys. Rep.}\ }\textbf {\bibinfo {volume}
  {574}},\ \bibinfo {pages} {1 } (\bibinfo {year} {2015})},\ \Eprint
  {https://arxiv.org/abs/1401.2026} {arXiv:1401.2026 [gr-qc]}\BibitemShut
  {NoStop}%
\bibitem [{\citenamefont {Ali}\ and\ \citenamefont
  {Engli\v{s}}(2005)}]{Ali_RMP_2005}%
  \BibitemOpen
  \bibfield  {author} {\bibinfo {author} {\bibfnamefont {S.~T.}\ \bibnamefont
  {Ali}}\ and\ \bibinfo {author} {\bibfnamefont {M.}~\bibnamefont
  {Engli\v{s}}},\ }\bibfield  {title} {\bibinfo {title} {Quantization methods:
  A guide for physicists and analysts},\ }\href
  {https://doi.org/10.1142/s0129055x05002376} {\bibfield  {journal} {\bibinfo
  {journal} {Rev. Math. Phys.}\ }\textbf {\bibinfo {volume} {17}},\ \bibinfo
  {pages} {391 } (\bibinfo {year} {2005})},\ \Eprint
  {https://arxiv.org/abs/0405065} {arXiv:0405065 [math-ph]}\BibitemShut
  {NoStop}%
\bibitem [{\citenamefont {Fredenhagen}\ and\ \citenamefont
  {Rejzner}(2016)}]{Fredenhagen_JMP_2016}%
  \BibitemOpen
  \bibfield  {author} {\bibinfo {author} {\bibfnamefont {K.}~\bibnamefont
  {Fredenhagen}}\ and\ \bibinfo {author} {\bibfnamefont {K.}~\bibnamefont
  {Rejzner}},\ }\bibfield  {title} {\bibinfo {title} {Quantum field theory on
  curved spacetimes: Axiomatic framework and examples},\ }\href
  {https://doi.org/10.1063/1.4939955} {\bibfield  {journal} {\bibinfo
  {journal} {J. Math. Phys.}\ }\textbf {\bibinfo {volume} {57}},\ \bibinfo
  {pages} {031101} (\bibinfo {year} {2016})},\ \Eprint
  {https://arxiv.org/abs/1412.5125} {arXiv:1412.5125 [math-ph]}\BibitemShut
  {NoStop}%
\bibitem [{\citenamefont {Brunetti}\ \emph {et~al.}(2003)\citenamefont
  {Brunetti}, \citenamefont {Fredenhagen},\ and\ \citenamefont
  {Verch}}]{Brunetti_CMP_2003}%
  \BibitemOpen
  \bibfield  {author} {\bibinfo {author} {\bibfnamefont {R.}~\bibnamefont
  {Brunetti}}, \bibinfo {author} {\bibfnamefont {K.}~\bibnamefont
  {Fredenhagen}},\ and\ \bibinfo {author} {\bibfnamefont {R.}~\bibnamefont
  {Verch}},\ }\bibfield  {title} {\bibinfo {title} {The generally covariant
  locality principle -- a new paradigm for local quantum field theory},\ }\href
  {https://doi.org/10.1007/s00220-003-0815-7} {\bibfield  {journal} {\bibinfo
  {journal} {Commun. Math. Phys.}\ }\textbf {\bibinfo {volume} {237}},\
  \bibinfo {pages} {31} (\bibinfo {year} {2003})},\ \Eprint
  {https://arxiv.org/abs/0112041} {arXiv:0112041 [math-ph]}\BibitemShut
  {NoStop}%
\bibitem [{\citenamefont {Fewster}\ and\ \citenamefont
  {Verch}(2015)}]{Fewster_Springer_2015}%
  \BibitemOpen
  \bibfield  {author} {\bibinfo {author} {\bibfnamefont {C.~J.}\ \bibnamefont
  {Fewster}}\ and\ \bibinfo {author} {\bibfnamefont {R.}~\bibnamefont
  {Verch}},\ }\bibinfo {title} {Algebraic quantum field theory in curved
  spacetimes},\ in\ \href {https://doi.org/10.1007/978-3-319-21353-8_4} {\emph
  {\bibinfo {booktitle} {Advances in Algebraic Quantum Field Theory}}},\
  \bibinfo {editor} {edited by\ \bibinfo {editor} {\bibfnamefont
  {R.}~\bibnamefont {Brunetti}}, \bibinfo {editor} {\bibfnamefont
  {C.}~\bibnamefont {Dappiaggi}}, \bibinfo {editor} {\bibfnamefont
  {K.}~\bibnamefont {Fredenhagen}},\ and\ \bibinfo {editor} {\bibfnamefont
  {J.}~\bibnamefont {Yngvason}}}\ (\bibinfo  {publisher} {Springer
  International Publishing},\ \bibinfo {address} {Switzerland},\ \bibinfo
  {year} {2015})\ pp.\ \bibinfo {pages} {125 -- 189},\ \Eprint
  {https://arxiv.org/abs/1504.00586} {arXiv:1504.00586 [math-ph]}\BibitemShut
  {NoStop}%
\bibitem [{\citenamefont {Dimock}(1980)}]{Dimock_CMP_1980}%
  \BibitemOpen
  \bibfield  {author} {\bibinfo {author} {\bibfnamefont {J.}~\bibnamefont
  {Dimock}},\ }\bibfield  {title} {\bibinfo {title} {Algebras of local
  observables on a manifold},\ }\href {https://doi.org/10.1007/BF01269921}
  {\bibfield  {journal} {\bibinfo  {journal} {Commun. Math. Phys.}\ }\textbf
  {\bibinfo {volume} {77}},\ \bibinfo {pages} {219 } (\bibinfo {year}
  {1980})}\BibitemShut {NoStop}%
\bibitem [{\citenamefont {Kay}\ and\ \citenamefont {Wald}(1991)}]{Kay_PR_1991}%
  \BibitemOpen
  \bibfield  {author} {\bibinfo {author} {\bibfnamefont {B.~S.}\ \bibnamefont
  {Kay}}\ and\ \bibinfo {author} {\bibfnamefont {R.~M.}\ \bibnamefont {Wald}},\
  }\bibfield  {title} {\bibinfo {title} {Theorems on the uniqueness and thermal
  properties of stationary, nonsingular, quasifree states on spacetimes with a
  bifurcate {K}illing horizon},\ }\href
  {https://doi.org/http://dx.doi.org/10.1016/0370-1573(91)90015-E} {\bibfield
  {journal} {\bibinfo  {journal} {Phys. Rep.}\ }\textbf {\bibinfo {volume}
  {207}},\ \bibinfo {pages} {49} (\bibinfo {year} {1991})}\BibitemShut
  {NoStop}%
\bibitem [{\citenamefont {Sanders}(2010)}]{Sanders_RMP_2010}%
  \BibitemOpen
  \bibfield  {author} {\bibinfo {author} {\bibfnamefont {K.}~\bibnamefont
  {Sanders}},\ }\bibfield  {title} {\bibinfo {title} {The locally covariant
  {D}irac field},\ }\href {https://doi.org/10.1142/S0129055X10003990}
  {\bibfield  {journal} {\bibinfo  {journal} {Rev. Math. Phys.}\ }\textbf
  {\bibinfo {volume} {22}},\ \bibinfo {pages} {381} (\bibinfo {year} {2010})},\
  \Eprint {https://arxiv.org/abs/0911.1304} {arXiv:0911.1304 [math-ph]}\BibitemShut 
  {NoStop}%
\bibitem [{\citenamefont {Schambach}\ and\ \citenamefont
  {Sanders}(2018)}]{Schambach_ReptMathPhys_2018}%
  \BibitemOpen
  \bibfield  {author} {\bibinfo {author} {\bibfnamefont {M.}~\bibnamefont
  {Schambach}}\ and\ \bibinfo {author} {\bibfnamefont {K.}~\bibnamefont
  {Sanders}},\ }\bibfield  {title} {\bibinfo {title} {The {P}roca field in
  curved spacetimes and its zero mass limit},\ }\href
  {https://doi.org/https://doi.org/10.1016/S0034-4877(18)30086-7} {\bibfield
  {journal} {\bibinfo  {journal} {Rep. Math. Phys.}\ }\textbf {\bibinfo
  {volume} {82}},\ \bibinfo {pages} {203} (\bibinfo {year} {2018})},\ \Eprint
  {https://arxiv.org/abs/1709.01911} {arXiv:1709.01911 [math-ph]}\BibitemShut
  {NoStop}%
\bibitem [{\citenamefont {Pfenning}(2009)}]{Pfenning_CQG_2009}%
  \BibitemOpen
  \bibfield  {author} {\bibinfo {author} {\bibfnamefont {M.~J.}\ \bibnamefont
  {Pfenning}},\ }\bibfield  {title} {\bibinfo {title} {Quantization of the
  {M}axwell field in curved spacetimes of arbitrary dimension},\ }\href
  {https://doi.org/10.1088/0264-9381/26/13/135017} {\bibfield  {journal}
  {\bibinfo  {journal} {Class. Quantum Grav.}\ }\textbf {\bibinfo {volume}
  {26}},\ \bibinfo {pages} {135017} (\bibinfo {year} {2009})},\ \Eprint
  {https://arxiv.org/abs/0902.4887} {arXiv:0902.4887 [math-ph]}\BibitemShut
  {NoStop}%
\bibitem [{\citenamefont {Sanders}\ \emph {et~al.}(2014)\citenamefont
  {Sanders}, \citenamefont {Dappiaggi},\ and\ \citenamefont
  {Hack}}]{Sanders_CMP_2014}%
  \BibitemOpen
  \bibfield  {author} {\bibinfo {author} {\bibfnamefont {K.}~\bibnamefont
  {Sanders}}, \bibinfo {author} {\bibfnamefont {C.}~\bibnamefont {Dappiaggi}},\
  and\ \bibinfo {author} {\bibfnamefont {T.-P.}\ \bibnamefont {Hack}},\
  }\bibfield  {title} {\bibinfo {title} {Electromagnetism, local covariance,
  the aharonov--bohm effect and gauss' law},\ }\href
  {https://doi.org/10.1007/s00220-014-1989-x} {\bibfield  {journal} {\bibinfo
  {journal} {Commun. Math. Phys.}\ }\textbf {\bibinfo {volume} {328}},\
  \bibinfo {pages} {625} (\bibinfo {year} {2014})},\ \Eprint
  {https://arxiv.org/abs/1211.6420} {arXiv:1211.6420 [math-ph]}\BibitemShut
  {NoStop}%
\bibitem [{\citenamefont {Finster}\ and\ \citenamefont
  {Strohmaier}(2015)}]{Finster_AHP_2015}%
  \BibitemOpen
  \bibfield  {author} {\bibinfo {author} {\bibfnamefont {F.}~\bibnamefont
  {Finster}}\ and\ \bibinfo {author} {\bibfnamefont {A.}~\bibnamefont
  {Strohmaier}},\ }\bibfield  {title} {\bibinfo {title}
  {Gupta{\textendash}bleuler quantization of the {M}axwell field in globally
  hyperbolic space-times},\ }\href {https://doi.org/10.1007/s00023-014-0363-z}
  {\bibfield  {journal} {\bibinfo  {journal} {Ann. Henri Poincar{\'{e}}}\
  }\textbf {\bibinfo {volume} {16}},\ \bibinfo {pages} {1837} (\bibinfo {year}
  {2015})},\ \Eprint {https://arxiv.org/abs/1307.1632} {arXiv:1307.1632
  [math-ph]}\BibitemShut {NoStop}%
\bibitem [{\citenamefont {Fredenhagen}\ and\ \citenamefont
  {Rejzner}(2013)}]{Fredenhagen_CMP_2013}%
  \BibitemOpen
  \bibfield  {author} {\bibinfo {author} {\bibfnamefont {K.}~\bibnamefont
  {Fredenhagen}}\ and\ \bibinfo {author} {\bibfnamefont {K.}~\bibnamefont
  {Rejzner}},\ }\bibfield  {title} {\bibinfo {title} {Batalin-{V}ilkovisky
  formalism in perturbative algebraic quantum field theory},\ }\href
  {https://doi.org/10.1007/s00220-012-1601-1} {\bibfield  {journal} {\bibinfo
  {journal} {Commun. Math. Phys.}\ }\textbf {\bibinfo {volume} {317}},\
  \bibinfo {pages} {697} (\bibinfo {year} {2013})},\ \Eprint
  {https://arxiv.org/abs/1110.5232} {arXiv:1110.5232 [math-ph]}\BibitemShut
  {NoStop}%
\bibitem [{\citenamefont {Drago}\ and\ \citenamefont
  {Moretti}(2020)}]{Drago_LMP_2020}%
  \BibitemOpen
  \bibfield  {author} {\bibinfo {author} {\bibfnamefont {N.}~\bibnamefont
  {Drago}}\ and\ \bibinfo {author} {\bibfnamefont {V.}~\bibnamefont
  {Moretti}},\ }\bibfield  {title} {\bibinfo {title} {The notion of observable
  and the moment problem for $\ast$ -algebras and their {G}{N}{S}
  representations},\ }\href {https://doi.org/10.1007/s11005-020-01277-x}
  {\bibfield  {journal} {\bibinfo  {journal} {Lett. Math. Phys.}\ }\textbf
  {\bibinfo {volume} {110}},\ \bibinfo {pages} {1711} (\bibinfo {year}
  {2020})},\ \Eprint {https://arxiv.org/abs/1903.07496} {arXiv:1903.07496
  [math-ph]}\BibitemShut {NoStop}%
\bibitem [{\citenamefont {Reed}\ and\ \citenamefont {Simon}(1975)}]{Reed_II}%
  \BibitemOpen
  \bibfield  {author} {\bibinfo {author} {\bibfnamefont {M.}~\bibnamefont
  {Reed}}\ and\ \bibinfo {author} {\bibfnamefont {B.}~\bibnamefont {Simon}},\
  }\href@noop {} {\emph {\bibinfo {title} {Fourier Analysis}}},\ \bibinfo
  {series} {Methods of Modern Mathematical Physics}, Vol.~\bibinfo {volume}
  {II}\ (\bibinfo  {publisher} {Acadmic Press},\ \bibinfo {address} {USA},\
  \bibinfo {year} {1975})\BibitemShut {NoStop}%
\bibitem [{\citenamefont {Itzykson}\ and\ \citenamefont
  {Zuber}(2005)}]{Itzykson_Dover_2005}%
  \BibitemOpen
  \bibfield  {author} {\bibinfo {author} {\bibfnamefont {C.}~\bibnamefont
  {Itzykson}}\ and\ \bibinfo {author} {\bibfnamefont {J.-B.}\ \bibnamefont
  {Zuber}},\ }\href@noop {} {\emph {\bibinfo {title} {Quantum Field Theory}}}\
  (\bibinfo  {publisher} {Dover},\ \bibinfo {address} {New York},\ \bibinfo
  {year} {2005})\BibitemShut {NoStop}%
\bibitem [{\citenamefont {Fewster}\ and\ \citenamefont
  {Verch}(2012)}]{Fewster_AHP_SPAS_2012}%
  \BibitemOpen
  \bibfield  {author} {\bibinfo {author} {\bibfnamefont {C.~J.}\ \bibnamefont
  {Fewster}}\ and\ \bibinfo {author} {\bibfnamefont {R.}~\bibnamefont
  {Verch}},\ }\bibfield  {title} {\bibinfo {title} {Dynamical locality and
  covariance: What makes a physical theory the same in all spacetimes?},\
  }\href {https://doi.org/https://doi.org/10.1007/s00023-012-0165-0} {\bibfield
   {journal} {\bibinfo  {journal} {Ann. Henri Poincar\'{e}}\ }\textbf {\bibinfo
  {volume} {13}},\ \bibinfo {pages} {1613} (\bibinfo {year} {2012})},\ \Eprint
  {https://arxiv.org/abs/1106.4785} {arXiv:1106.4785 [math-ph]}\BibitemShut
  {NoStop}%
\bibitem [{\citenamefont {Fewster}(2018)}]{Fewster_IJMP_2018}%
  \BibitemOpen
  \bibfield  {author} {\bibinfo {author} {\bibfnamefont {C.~J.}\ \bibnamefont
  {Fewster}},\ }\bibfield  {title} {\bibinfo {title} {The art of the state},\
  }\href {https://doi.org/10.1142/S0218271818430071} {\bibfield  {journal}
  {\bibinfo  {journal} {Int. J. Mod. Phys. D}\ }\textbf {\bibinfo {volume}
  {27}},\ \bibinfo {pages} {1843007} (\bibinfo {year} {2018})},\ \Eprint
  {https://arxiv.org/abs/1803.06836} {arXiv:1803.06836 [gr-qc]}\BibitemShut
  {NoStop}%
\bibitem [{\citenamefont {Khavkine}\ and\ \citenamefont
  {Moretti}(2016)}]{Khavkine_CMP_2016}%
  \BibitemOpen
  \bibfield  {author} {\bibinfo {author} {\bibfnamefont {I.}~\bibnamefont
  {Khavkine}}\ and\ \bibinfo {author} {\bibfnamefont {V.}~\bibnamefont
  {Moretti}},\ }\bibfield  {title} {\bibinfo {title} {Analytic dependence is an
  unnecessary requirement in renormalization of locally covariant {Q}{F}{T}},\
  }\href {https://doi.org/10.1007/s00220-016-2618-7} {\bibfield  {journal}
  {\bibinfo  {journal} {Commun. Math. Phys.}\ }\textbf {\bibinfo {volume}
  {344}},\ \bibinfo {pages} {581} (\bibinfo {year} {2016})},\ \Eprint
  {https://arxiv.org/abs/1411.1302} {arXiv:1411.1302 [gr-qc]}\BibitemShut
  {NoStop}%
\bibitem [{\citenamefont {Fredenhagen}\ and\ \citenamefont
  {Rejzner}(2015)}]{Fredenhagen_Springer_2015}%
  \BibitemOpen
  \bibfield  {author} {\bibinfo {author} {\bibfnamefont {K.}~\bibnamefont
  {Fredenhagen}}\ and\ \bibinfo {author} {\bibfnamefont {K.}~\bibnamefont
  {Rejzner}},\ }\bibinfo {title} {Perturbative construction of models of
  algebraic quantum field theory},\ in\ \href
  {https://doi.org/10.1007/978-3-319-21353-8_2} {\emph {\bibinfo {booktitle}
  {Advances in Algebraic Quantum Field Theory}}},\ \bibinfo {editor} {edited
  by\ \bibinfo {editor} {\bibfnamefont {R.}~\bibnamefont {Brunetti}}, \bibinfo
  {editor} {\bibfnamefont {C.}~\bibnamefont {Dappiaggi}}, \bibinfo {editor}
  {\bibfnamefont {K.}~\bibnamefont {Fredenhagen}},\ and\ \bibinfo {editor}
  {\bibfnamefont {J.}~\bibnamefont {Yngvason}}}\ (\bibinfo  {publisher}
  {Springer International Publishing},\ \bibinfo {address} {Switzerland},\
  \bibinfo {year} {2015})\ pp.\ \bibinfo {pages} {31 -- 74},\ \Eprint
  {https://arxiv.org/abs/1503.07814} {arXiv:1503.07814 [math-ph]}\BibitemShut
  {NoStop}%
\bibitem [{\citenamefont {Ashtekar}(1984)}]{Ashtekar_Springer_1984}%
  \BibitemOpen
  \bibfield  {author} {\bibinfo {author} {\bibfnamefont {A.}~\bibnamefont
  {Ashtekar}},\ }\bibinfo {title} {General relativity and gravitation: Invited
  papers and discussion reports of the 10th international conference on general
  relativity and gravitation, {P}adua, {J}uly 3–8, 1983}\ (\bibinfo
  {publisher} {Springer},\ \bibinfo {address} {Netherlands},\ \bibinfo {year}
  {1984})\ Chap.\ \bibinfo {chapter} {Quantum Gravity and Quantum Field Theory
  in a Curved Space Report of Workshop D1}\BibitemShut {NoStop}%
\bibitem [{\citenamefont {Haag}\ \emph {et~al.}(1984)\citenamefont {Haag},
  \citenamefont {Narnhofer},\ and\ \citenamefont {Stein}}]{Haag_CMP_1984}%
  \BibitemOpen
  \bibfield  {author} {\bibinfo {author} {\bibfnamefont {R.}~\bibnamefont
  {Haag}}, \bibinfo {author} {\bibfnamefont {H.}~\bibnamefont {Narnhofer}},\
  and\ \bibinfo {author} {\bibfnamefont {U.}~\bibnamefont {Stein}},\ }\bibfield
   {title} {\bibinfo {title} {On quantum field theory in gravitational
  background},\ }\href {https://link.springer.com/article/10.1007%2FBF01209302}
  {\bibfield  {journal} {\bibinfo  {journal} {Commun. Math. Phys.}\ }\textbf
  {\bibinfo {volume} {94}},\ \bibinfo {pages} {219 } (\bibinfo {year}
  {1984})}\BibitemShut {NoStop}%
\bibitem [{\citenamefont {Moretti}(2021)}]{Moretti_LMP_2021}%
  \BibitemOpen
  \bibfield  {author} {\bibinfo {author} {\bibfnamefont {V.}~\bibnamefont
  {Moretti}},\ }\bibfield  {title} {\bibinfo {title} {On the global {H}adamard
  parametrix in {Q}{F}{T} and the signed squared geodesic distance defined in
  domains larger than convex normal neighbourhoods},\ }\href
  {https://doi.org/10.1007/s11005-021-01464-4} {\bibfield  {journal} {\bibinfo
  {journal} {Lett. Math. Phys.}\ }\textbf {\bibinfo {volume} {111}},\ \bibinfo
  {pages} {130} (\bibinfo {year} {2021})},\ \Eprint
  {https://arxiv.org/abs/2107.04903} {arXiv:2107.04903 [gr-qc]}\BibitemShut
  {NoStop}%
\bibitem [{\citenamefont {Koehler}(1995)}]{Koehler_PhD}%
  \BibitemOpen
  \bibfield  {author} {\bibinfo {author} {\bibfnamefont {M.}~\bibnamefont
  {Koehler}},\ }\emph {\bibinfo {title} {The stress energy tensor of a locally
  supersymmetric quantum field on a curved spacetime}},\ \href
  {https://doi.org/10.48550/ARXIV.GR-QC/9505014} {\bibinfo {type} {Ph{D}
  thesis}} (\bibinfo {year} {1995}),\ \Eprint {https://arxiv.org/abs/9505014}
  {arXiv:9505014 [gr-qc]}\BibitemShut {NoStop}%
\bibitem [{\citenamefont {Kratzert}(2000)}]{Kratzert_AnnPhys_2000}%
  \BibitemOpen
  \bibfield  {author} {\bibinfo {author} {\bibfnamefont {K.}~\bibnamefont
  {Kratzert}},\ }\bibfield  {title} {\bibinfo {title} {Singularity structure of
  the two point function of the free {D}irac field on a globally hyperbolic
  spacetime},\ }\href
  {https://doi.org/10.1002/1521-3889(200006)9:6<475::AID-ANDP475>3.0.CO;2-S}
  {\bibfield  {journal} {\bibinfo  {journal} {Ann. Phys.}\ }\textbf {\bibinfo
  {volume} {9}},\ \bibinfo {pages} {475} (\bibinfo {year} {2000})},\ \Eprint
  {https://arxiv.org/abs/0003015} {arXiv:0003015 [math-ph]}\BibitemShut
  {NoStop}%
\bibitem [{\citenamefont {Khavkine}\ and\ \citenamefont
  {Moretti}(2015)}]{Khavkine_Springer_2015}%
  \BibitemOpen
  \bibfield  {author} {\bibinfo {author} {\bibfnamefont {I.}~\bibnamefont
  {Khavkine}}\ and\ \bibinfo {author} {\bibfnamefont {V.}~\bibnamefont
  {Moretti}},\ }\bibinfo {title} {Algebraic {Q}{F}{T} in curved spacetime and
  quasifree {H}adamard states: An introduction},\ in\ \href
  {https://doi.org/10.1007/978-3-319-21353-8_5} {\emph {\bibinfo {booktitle}
  {Advances in Algebraic Quantum Field Theory}}},\ \bibinfo {editor} {edited
  by\ \bibinfo {editor} {\bibfnamefont {R.}~\bibnamefont {Brunetti}}, \bibinfo
  {editor} {\bibfnamefont {C.}~\bibnamefont {Dappiaggi}}, \bibinfo {editor}
  {\bibfnamefont {K.}~\bibnamefont {Fredenhagen}},\ and\ \bibinfo {editor}
  {\bibfnamefont {J.}~\bibnamefont {Yngvason}}}\ (\bibinfo  {publisher}
  {Springer International Publishing},\ \bibinfo {address} {Switzerland},\
  \bibinfo {year} {2015})\ pp.\ \bibinfo {pages} {191 -- 251},\ \Eprint
  {https://arxiv.org/abs/1412.5945} {arXiv:1412.5945 [math-ph]}\BibitemShut
  {NoStop}%
\end{thebibliography}
